\documentclass[11pt]{amsbook}
\usepackage{float, graphicx}
\usepackage[]{epsfig}
\usepackage{amsmath, amsthm, amssymb,}
\usepackage{epsfig}
\usepackage{verbatim}
\usepackage{multicol}
\usepackage{url}
\usepackage{adjustbox}
\usepackage{latexsym}
\usepackage{mathrsfs}
\usepackage[colorlinks, bookmarks=true]{hyperref}
\usepackage{graphicx}
\usepackage{enumerate}
\usepackage[normalem]{ulem}
\usepackage{bm}
\usepackage{stmaryrd}
\usepackage{enumitem}
\usepackage{tikz}
\usetikzlibrary{calc}
\usepackage{mathtools}
\usepackage{subcaption}

\usepackage{color}
\usepackage{xcolor}
\usepackage{geometry} 
\usepackage{chngcntr}

\counterwithout{section}{chapter}
\counterwithout{figure}{chapter}
\counterwithout{table}{chapter}
\counterwithout{footnote}{chapter}



\numberwithin{equation}{section}

\newcommand{\cC}{{\mathcal C}}

\newcommand{\cE}{{\mathcal E}}

\newcommand{\cG}{{\mathcal G}}

\newcommand{\cT}{{\mathcal T}}

\newcommand{\cN}{{\mathcal N}}

\newcommand{\cI}{{\mathcal I}}

\newcommand{ \sfx }{{\mathsf x}}

\newcommand{ \sft }{{\mathsf t}}
\newcommand{\sfa}{{\mathsf a}}
\newcommand{ \sfb }{{\mathsf b}}

\newcommand{ \sfh}{{\mathsf h}}

\newcommand{\sfC}{{\mathsf C}}
\newcommand{\sfD}{{\mathsf D}}

\newcommand{\sfP}{{\mathsf P}}
\newcommand{\sfH}{{\mathsf H}}
\newcommand{\sfS}{{\mathsf S}}
\newcommand{ \sfs }{{\mathsf s}}

\newcommand{\fa}{{\mathfrak a}}

\newcommand{\fA}{{\mathfrak A}}
\newcommand{\fU}{{\mathfrak U}}
\newcommand{\fN}{{\mathfrak N}}
\newcommand{\fb}{{\mathfrak b}}
\newcommand{\fc}{{\mathfrak c}}
\newcommand{\fC}{{\mathfrak C}}
\newcommand{\fL}{{\mathfrak L}}
\newcommand{\fF}{{\mathfrak F}}

\newcommand{\ft}{{\mathfrak t}}

\newcommand{\fT}{{\mathfrak T}}

\newcommand{\fP}{{\mathfrak P}}
\newcommand{\fO}{{\mathfrak O}}
\newcommand{\fW}{{\mathfrak W}}
\newcommand{\fR}{{\mathfrak R}}
\newcommand{\fB}{{\mathfrak B}}
\newcommand{\fV}{{\mathfrak V}}

\newcommand{\bmv}{{\bm{v}}}

\newcommand{\rd}{{\rm d}}

\newcommand{\ri}{\mathrm{i}}

\newcommand{\rb}{{\mathrm b}} 
\newcommand{\rw}{{\mathrm w}}

\newcommand{\bB}{{\mathbb B}}
\newcommand{\bC}{{\mathbb C}}
\newcommand{\bD}{{\mathbb D}}
\newcommand{\bX}{{\mathbb X}}

\newcommand{\bE}{\mathbb{E}}

\newcommand{\bP}{\mathbb{P}}

\newcommand{\bR}{{\mathbb R}}

\newcommand{\bT}{\mathbb T}

\newcommand{\bZ}{\mathbb{Z}}
\newcommand{\bW}{\mathbb{W}}

\newcommand{\al}{\alpha}

\DeclareMathOperator{\dist}{dist}

\DeclareMathOperator{\sgn}{sgn}
\DeclareMathOperator{\OO}{O}
\DeclareMathOperator{\oo}{o}
\DeclareMathOperator{\argmax}{argmax}

\DeclareMathOperator{\Res}{Res}

\renewcommand{\Re}{\mathop{\mathrm{Re}}}
\renewcommand{\Im}{\mathop{\mathrm{Im}}}

\renewcommand{\leq}{\leqslant}
\renewcommand{\geq}{\geqslant}
\renewcommand{\le}{\leqslant}
\renewcommand{\ge}{\geqslant}

\newcommand{\del}{\partial}

\newcommand{\wh}{\widehat}
\newcommand{\wt}{\widetilde}

\newcommand{\beq}{\begin{equation}}
\newcommand{\eeq}{\end{equation}}
\newcommand{\Adm}{\mathrm{Adm}}
\usepackage{tikz}
\usetikzlibrary{arrows}
\usepackage{cleveref}

\newtheorem{theorem}{Theorem}[section]
\newtheorem{proposition}[theorem]{Proposition}
\newtheorem{lemma}[theorem]{Lemma}
\newtheorem{cor}[theorem]{Corollary}
\newtheorem{conjecture}[theorem]{Conjecture}

\newtheorem{remark}[theorem]{Remark}
\newtheorem{definition}[theorem]{Definition}
\newtheorem{assumption}[theorem]{Assumption}

\title{Height fluctuation for Lozenge Tilings of Polygons}

    \author{Jiaoyang Huang}
    \address{University of Philadelphia, PA}
    \email{huangjy@wharton.upenn.edu}

\begin{document}

\title{Height Fluctuations of Lozenge Tilings of Polygons}

\begin{abstract}
We establish Gaussian free field fluctuations for uniformly random lozenge
tilings of simply connected polygonal domains with \(3d\) sides whose
directions cycle through the three lattice directions. More precisely,
assuming that the liquid region is connected and that the boundary data do
not force the height at any interior point, we prove that the fluctuations of centered height function converge to the Gaussian free field in the liquid region, confirming a prediction of Kenyon and Okounkov from 2007.

We introduce a tiling action function that encodes the geometry of
the limit shape through its critical points. The action function has a
complex conjugate pair of critical points in the liquid region, repeated
real critical points on the arctic boundary, and distinct real critical
points in the frozen region. Using this tiling action function, we construct
an approximation to the inverse Kasteleyn matrix in terms of explicit
single-contour and double-contour integrals and prove that the approximation is
uniform throughout the polygonal domain. The convergence to the Gaussian
free field then follows from standard kernel computations.
\end{abstract}
\maketitle
{
  \hypersetup{linkcolor=black}
  \setcounter{tocdepth}{1}
  \tableofcontents
}

\chapter{Results and Preliminaries}

\section{Introduction}

A perfect matching, or dimer configuration, on a simply connected planar
bipartite graph can be encoded by a discrete height function on the faces
of the graph, defined up to an additive constant
\cite{thurston1990conway}. In tiling models such as lozenge tilings, the
graph of the height function can be naturally realized as a stepped
surface in \(\bR^3\). More generally, after a suitable normalization,
dimer height functions belong to a class of discrete Lipschitz
functions. The dimer model may therefore be viewed as a model of random
Lipschitz surfaces. It is one of the few two-dimensional random-interface
models for which many microscopic observables can be computed exactly,
while the macroscopic geometry exhibits highly nontrivial phenomena,
including limit shapes and phase separation into frozen and liquid
regions. It consequently provides an important testing ground for
universality in random surface theory: one would like to determine which
features of the large-scale fluctuations depend on the microscopic
lattice structure and which are governed only by the macroscopic geometry
of the surface.

A central universality prediction is that, in the liquid region, the
centered height fluctuations are asymptotically described by the Gaussian
free field. The Gaussian free field is the canonical conformally
invariant, log-correlated Gaussian random distribution in two dimensions.
Establishing this convergence gives a field-level description of the
fluctuations and, in particular, reveals an emergent conformal invariance
that is not apparent in the underlying lattice model.

\begin{figure}
\begin{center}
 \includegraphics[scale=0.4,trim={0cm 8cm 0 7cm},clip]{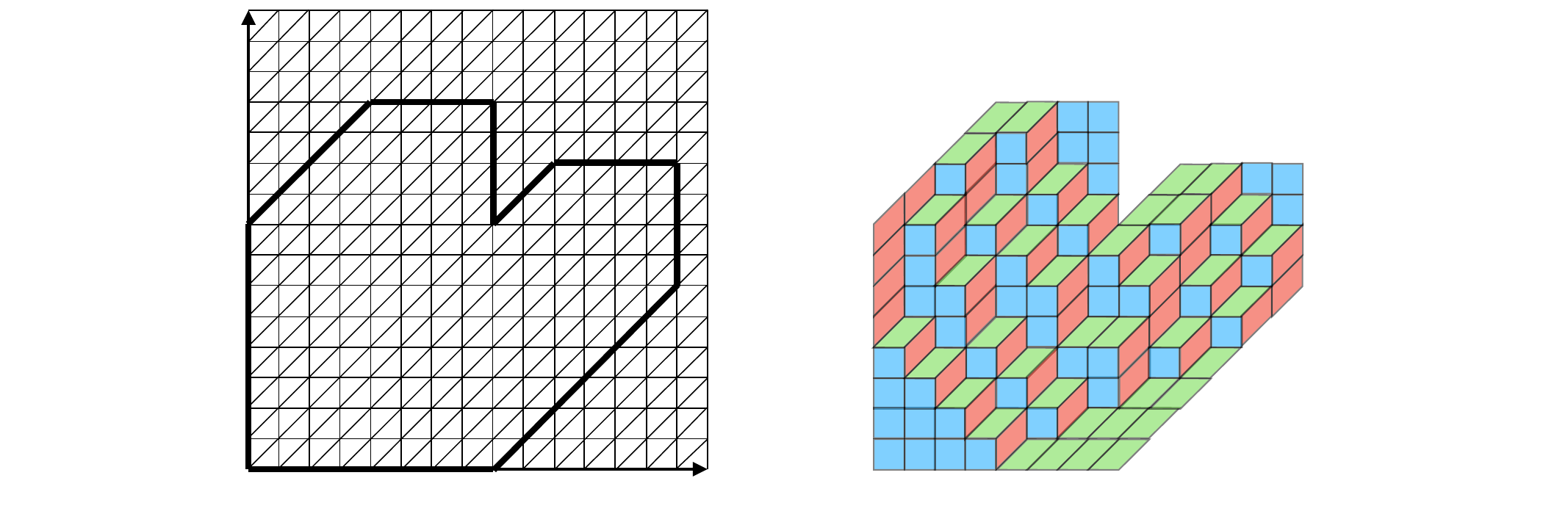}
 \caption{Lozenge tiling of polygonal domains on the standard square grid.}
 \label{f:lattice}
 \end{center}
 \end{figure}

In this paper, we study height functions of uniformly random lozenge
tilings of simply connected polygonal domains, or equivalently the
uniform dimer model on the honeycomb lattice. For notational convenience, we apply an affine transformation and use
coordinates in which the three lattice directions are parallel to
\((1,0)\), \((0,1)\), and \((1,1)\); see \Cref{f:lattice}.

A first fundamental question concerns the law of large numbers for the
height function. For a sequence of discrete domains approximating a
fixed macroscopic region $\fR$, the rescaled random height function
converges to a deterministic limit shape. This limit shape is
characterized as the minimizer of the surface-tension functional
\[
\int_{\fR} \sigma\bigl(\nabla h(x,s)\bigr)\,{\rm d}x\,{\rm d}s
\]
among all admissible Lipschitz functions with the prescribed boundary
values
\cite{cohn2001variational,kenyon2006dimers,kenyon2007limit}.
Equivalently, it maximizes the macroscopic entropy of the tiling. In the
liquid region, the Euler--Lagrange equation for this variational problem
can be rewritten as a complex Burgers equation, which endows the liquid
region with a natural complex structure
\cite{kenyon2007limit}. The regularity of the minimizer and the geometry
of the associated free boundary have been studied in considerable
generality. In particular, for natural classes of polygonal domains, the
frozen boundary is an algebraic curve and obeys the
Pokrovsky--Talapov law at its generic points
\cite{kenyon2007limit,astala2026dimer}.

One of the most striking features of these limit shapes is the
coexistence of frozen and liquid phases. In a frozen region, one lozenge
type has asymptotic density one. The limit shape is affine there, and the
microscopic tiling is essentially deterministic. In the liquid region,
all three lozenge types have positive densities, and the tiling remains
genuinely random at every scale. The curve separating these two phases
is called the \emph{arctic boundary}. Thus, even though the original
polygon is deterministic, the effective domain in which nontrivial
fluctuations survive is selected by the variational problem itself.

The fluctuation problem is therefore concentrated in the liquid region.
For dimer models, the conformal structure governing the limiting field
is generally not the Euclidean structure of the original coordinates.
Instead, it is determined by the complex slope of the limit shape
\cite{kenyon2006dimers,kenyon2007limit}. The prediction of Kenyon and
Okounkov is that, after passing to the corresponding natural complex
coordinate, the centered height fluctuations converge, as a random
distribution, to the Gaussian free field with zero boundary conditions
on the liquid region. The main result of this paper verifies this
prediction for uniformly random lozenge tilings of simply connected
polygonal domains, under the assumption that the liquid region is
connected.

\subsection{Approaches to Gaussian free field fluctuations}

Several complementary methods have been developed to prove convergence
of dimer height fluctuations to the Gaussian free field; see
\cite{gorin2021lectures} for a broader overview. We briefly review the approaches
most closely related to the present work.

\emph{Inverse Kasteleyn matrices and discrete complex analysis.}
For a bipartite planar dimer model, edge correlation functions are
determinants whose entries are given by the inverse Kasteleyn matrix.
The Kasteleyn operator can be interpreted as a discrete Dirac
operator, and entries of its inverse,
satisfy discrete holomorphic or discrete harmonic boundary-value
problems. Joint moments of height differences can therefore be
expressed in terms of products of the inverse Kasteleyn matrix. Once
they are shown to converge to their continuum
counterparts, these moment formulas identify the Green function and
the Wick contractions of the limiting field. This strategy was
developed by Kenyon for domino and lozenge tilings
\cite{kenyon2000conformal,kenyon2001dominos,kenyon2002laplacian,kenyon2008height}
and was subsequently extended to several classes of domains and
discretizations
\cite{de2007scaling,li2017conformal,
russkikh2018dimers,russkikh2018dominos}.

This approach is especially effective when the underlying discrete
complex structure remains uniformly nondegenerate. The presence of an
arctic boundary makes the analysis more delicate: the leading
asymptotics of the inverse Kasteleyn matrix change from oscillatory in
the liquid region to exponentially decaying in the frozen region, with
additional transition regimes near the arctic boundary.

\emph{\(t\)-embeddings.}
The theory of \(t\)-embeddings provides a more intrinsic version of
discrete complex analysis for general bipartite planar dimer models. In
this framework, the geometry needed for discrete complex analysis is
encoded by the embedding and its associated origami map, while
\(t\)-holomorphic functions replace the usual discrete holomorphic
observables. Compactness and regularity estimates then give general
criteria that reduce Gaussian free field convergence to the
construction and uniform control of suitable perfect \(t\)-embeddings
\cite{chelkak2020dimer,chelkak2021bipartite}.

The main difficulty in applying this framework is to construct the
relevant embeddings and prove the required nondegeneracy and convergence
of their associated surfaces. This program has been carried out for
uniformly weighted Aztec diamonds and regular hexagonal lozenge-tiling
domains
\cite{berggren2024perfect,berggren2024perfectb}.

\emph{Spanning trees and imaginary geometry.}
For dimer models admitting a Temperley-type correspondence, the height
function can be related to the winding of branches of an associated
uniform spanning tree. Berestycki, Laslier, and Ray proved that, under
an invariance principle for random walk and a Russo--Seymour--Welsh-type
crossing estimate, these winding fluctuations converge to a Gaussian
free field
\cite{berestycki2020dimers}. This approach is robust and uses the
exact solvability of the dimer model only to a limited extent. Among
its applications are lozenge tilings with boundary data lying in a
plane and Temperleyan domains in isoradial graphs. Related ideas have
also been used to study lozenge tilings with curved limit shapes
\cite{laslier2021central}.

\emph{Exact correlation kernels and steepest descent.}
For several integrable families of domains, the inverse Kasteleyn
matrix or an equivalent determinantal correlation kernel admits an
explicit double-contour representation. The critical points of the
integrand encode the complex slope of the limit shape. A
steepest-descent analysis of the contour integrals then produces the
continuum Cauchy-type kernels that determine the limiting covariance,
while the determinantal moment or cumulant formulas establish
Gaussianity.

This strategy has been used for anisotropic growth models and
interlacing particle systems
\cite{borodin2014anisotropic,duits2013gaussian,
kuan2014gaussian}. For lozenge tilings, Petrov first obtained an exact
double-contour formula for the correlation kernel through the
Gelfand--Tsetlin representation and then used it to prove Gaussian free
field fluctuations for trapezoid domains
\cite{petrov2014asymptotics,petrov2015asymptotics}.

\emph{Orthogonal-polynomial and recurrence methods.}
In determinantal models whose correlation kernels are built from
biorthogonal families, the cumulants of global linear statistics can
be expressed in terms of the recurrence matrices of the corresponding
polynomials. Their asymptotic behavior is governed by right limits of
these recurrence matrices. In particular, when the relevant right
limit is a Laurent matrix, the higher-order cumulants vanish and a
central limit theorem follows
\cite{breuer2017central}.

This method extends to multi-time and multi-level noncolliding
processes and gives another route to Gaussian free field fluctuations
of the associated random surfaces. In particular, it has been applied
to lozenge tilings of a hexagon and to their \(q\)-Racah deformations
\cite{duits2018global,duits2024lozenge}. The method is especially
effective when an appropriate finite-term recurrence relation is
available.

\emph{Schur generating functions and asymptotic representation
theory.}
Many random tiling models can be encoded by Schur processes or by
probability measures on signatures. Differential operators
diagonalized by Schur functions make it possible to extract moments
and cumulants of global height observables from derivatives of the
corresponding Schur generating functions. Asymptotic analysis of these
generating functions then yields laws of large numbers and Gaussian
fluctuations
\cite{bufetov2018fluctuations}.

This approach has been applied to several families of lozenge and
domino tilings and establishs the Kenyon--Okounkov prediction for certain
multiply connected polygonal domains
\cite{bufetov2019fourier}. Related Schur-process methods have also been
used for rectangular Aztec diamonds and for lozenge tilings of a
cylinder
\cite{bufetov2018asymptotics,ahn2022lozenge}.

\emph{Loop equations and discrete log-gases.}
Another approach is available when one-dimensional sections of a
random tiling can be identified with discrete log-gases. For
continuous \(\beta\)-ensembles, loop equations, also called
Dyson--Schwinger equations, were used by Johansson to prove
macroscopic central limit theorems
\cite{MR1487983}; see also
\cite{MR3010191,borot-guionnet2,KrSh}.
Borodin, Gorin, and Guionnet developed discrete analogues of these
identities, usually called Nekrasov equations, and used them to prove
Gaussian fluctuations for general discrete \(\beta\)-ensembles
\cite{borodin2017gaussian}. These equations are related to identities originating
in the work of Nekrasov and his collaborators
\cite{Nekrasov,Nek_PS,Nek_Pes}. Dimitrov and Knizel extended this
method to log-gases on quadratic lattices and to multilevel discrete
\(\beta\)-corners processes
\cite{dimitrov2022asymptotics,dimitrov2025multi}.

For random tilings, the conditional distributions of successive
parallel slices lead naturally to dynamical or multilevel loop
equations. Dynamical loop equations were used to prove Gaussian free
field fluctuations for polygons with one horizontal upper boundary
edge and, more generally, for several two-dimensional interacting
particle systems and deformed lozenge-tiling measures
\cite{huang2020height,gorin2024dynamical}.

More recently, Borot, Gorin, and Guionnet developed a general
asymptotic theory of discrete \(\beta\)-ensembles with several groups
of particles and applied it to random tilings
\cite{borot2026macroscopic}. Their framework allows fixed,
free, or constrained filling fractions and includes topologically
nontrivial tiling domains. Under off-criticality and appropriate
regularity assumptions, fixed filling fractions lead to Gaussian free
field fluctuations in orientable liquid regions. When the filling
fractions are allowed to vary, an additional discrete Gaussian
component appears. Their methods are complementary to the
inverse-Kasteleyn analysis developed here.

\emph{Our approach.}
Despite the numerous previous approaches and results, our setting is different. The main novelty of our \Cref{thm:moment-GFF} is that it applies when the tiled domain is an arbitrary simply connected polygonal domain.

Our method combines the inverse-Kasteleyn viewpoint with
contour-integral steepest descent. In the coordinates used here,
constructing the inverse Kasteleyn matrix amounts to solving a discrete
heat-type boundary-value problem, with homogeneous boundary conditions
determined by the polygonal domain. Rather than first deriving a single
exact double-contour representation for the inverse Kasteleyn matrix,
or for an equivalent particle correlation kernel, as in
\cite{petrov2014asymptotics,petrov2015asymptotics}, we construct local
double-contour approximations that serve as parametrices for the
inverse Kasteleyn matrix.

The form of these local parametrices depends on the geometry of the
relevant vertices, including the liquid, frozen, generic arctic,
tangency, cusp, and other transition regimes. We prove that the local
approximations are compatible in their overlap regions and patch them
together to obtain a uniform asymptotic approximation of the inverse
Kasteleyn matrix throughout the polygonal domain. In particular, our
estimates are not restricted to compact subsets of the liquid region:
they also control the frozen region and the transition neighborhoods
associated with the arctic and polygonal boundaries, with the
appropriate local asymptotic forms in each geometric regime.

This global control allows us to analyze joint moments of the height
function, identify the limiting covariance with the Green function of
the liquid region, and prove that all higher-order cumulants vanish.
We refer to \Cref{s:outline} for a more detailed overview of the
argument.

	\section{Main Results} 
	
	\label{Walks}

	\subsection{Lozenge tilings and Dimer models}

	\label{Tiling}

\begin{figure}
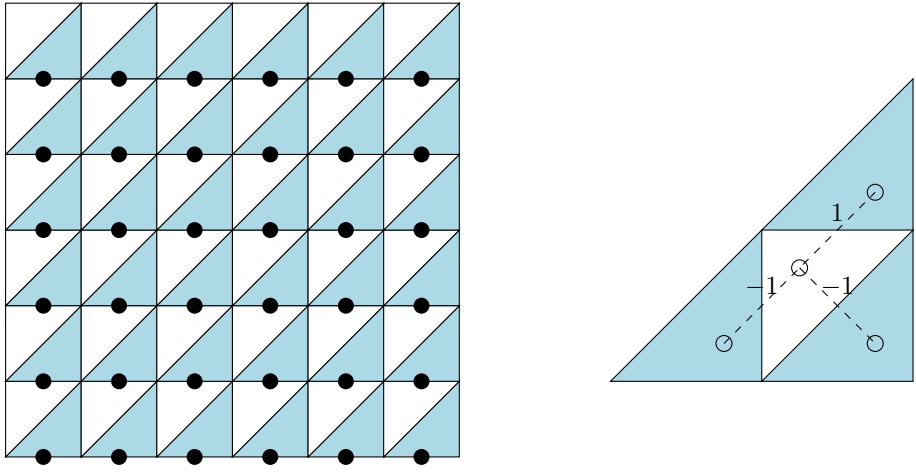

    \centering

    \caption{Triangular lattice $\mathbb{T}$ with blue and white triangle coloring, and its dual graph.}
    \label{fig:triangular_lattice}
\end{figure}

	We denote by $\bZ'=\bZ+1/2$ the set of half-integers. Let $\mathbb{T}$ be the \emph{triangular lattice}, defined as the graph with vertex set $\bZ'\times \mathbb{Z}$ and edge set consisting of pairs of vertices $(\sfx, \sfs), (\sfx', \sfs') \in \bZ'\times \mathbb{Z}$ such that
$
(\sfx' - \sfx, \sfs' - \sfs) \in \{ \pm(1, 0), \pm(0, 1), \pm(1, 1) \}
$.
The faces of \(\mathbb T\) are triangles of two types. A triangle with vertices
\[
\bigl\{(\sfx,\sfs),(\sfx+1,\sfs),(\sfx+1,\sfs+1)\bigr\}
\]
is called a \emph{blue triangle}, and we represent it by the integer point
\((\sfx+1/2,\sfs)\in\bZ^2\). A triangle with vertices
\[
\bigl\{(\sfx,\sfs),(\sfx,\sfs+1),(\sfx+1,\sfs+1)\bigr\}
\]
is called a \emph{white triangle}, and we represent it by the integer point
\((\sfx+1/2,\sfs+1)\in\bZ^2\); see \Cref{fig:triangular_lattice}. Here
\(\sfx\in\bZ'\) and \(\sfs\in\bZ\).

This coloring ensures that adjacent triangles always have opposite colors. The \emph{dual graph} of $\mathbb{T}$ is the \emph{hexagonal (honeycomb) lattice}, whose vertices correspond to the triangles of $\mathbb{T}$, and edges connect adjacent blue and white triangles.  We emphasize that \(\mathbb T\) is built on
the shifted lattice \(\bZ'\times\bZ\), rather than on the standard integer lattice
\(\bZ^2\). Nevertheless, with the above convention, triangular faces are naturally
labeled by integer points in \(\bZ^2\), together with their colors.

\begin{figure}
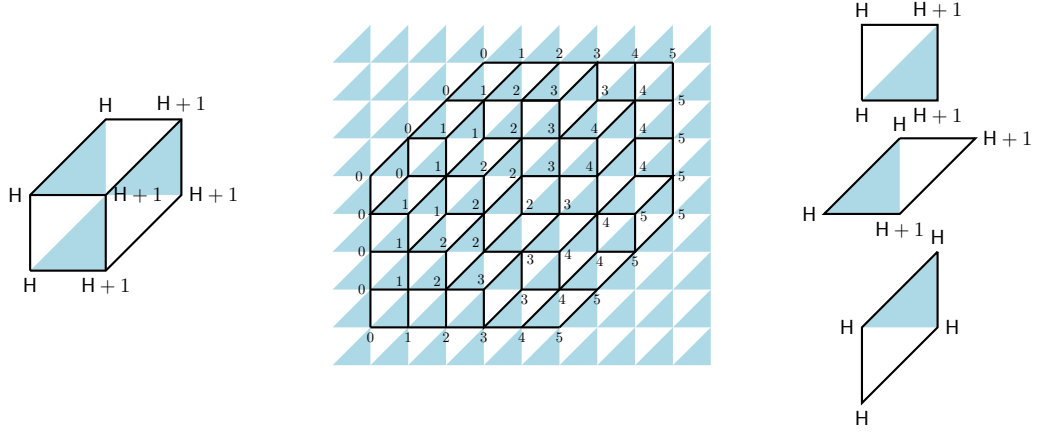

		
		\begin{center}		
			

			
		\end{center}
		
		\caption{\label{tilinghexagon} Depicted to the right are the three types of lozenges. Depicted in the middle is a lozenge tiling of a hexagon. One may view this tiling as a collection of nonintersecting Bernoulli paths which gives rise to a height function (shown in the middle). }
		
	\end{figure}

	A \emph{domain} $\mathsf{R} \subseteq \mathbb{T}$ is an induced subgraph of the triangular lattice $\mathbb{T}$ whose associated triangular-face realization is simply connected. By abuse of notation, we also write $\mathsf{R} \subseteq \mathbb{R}^2$ for the corresponding planar region, obtained as the union of all triangular faces of $\mathbb{T}$ whose vertices belong to $\mathsf{R}$. When $\mathsf{R}$ is viewed as this union of triangular faces, $\partial \mathsf{R}$ is the union of its boundary edges. A \emph{dimer covering} of a domain $\mathsf{R}$ is a perfect matching on the dual graph of $\mathsf{R}$. In this setting, each edge of the dual graph connects a pair of adjacent triangular faces in $\mathsf{R}$---specifically, one black and one white triangle. Every such matched pair forms a rhombus, or parallelogram, which we will refer to as a \emph{lozenge} or \emph{tile}. Lozenges can be oriented in one of three ways; see the right side of \Cref{tilinghexagon} for all three orientations. We refer to the topmost lozenge there, that is, one with vertices of the form
\[
\big\{ (\sfx, \sfs), (\sfx, \sfs + 1), (\sfx + 1, \sfs + 1), (\sfx + 1, \sfs) \big\},
\]
as a \emph{type I} lozenge. Similarly, we refer to the middle lozenge, with vertices of the form
\[
\big\{ (\sfx, \sfs), (\sfx + 1, \sfs), (\sfx + 2, \sfs + 1), (\sfx + 1, \sfs + 1) \big\},
\]
and the bottom lozenge, with vertices of the form
\[
\big\{ (\sfx, \sfs), (\sfx, \sfs + 1), (\sfx + 1, \sfs + 2), (\sfx + 1, \sfs + 1) \big\},
\]
as \emph{type II} and \emph{type III} lozenges, respectively. A dimer covering of $\mathsf{R}$ can equivalently be interpreted as a tiling of $\mathsf{R}$ by lozenges of types I, II, and III. Therefore, we will also refer to a dimer covering of $\mathsf{R}$ as a \emph{lozenge tiling}. We call $\mathsf{R}$ \emph{tileable} if it admits a tiling.

Associated with any tiling of $\mathsf{R}$ is a \emph{height function}, namely, a function $\mathsf{H}$ 
on the vertices of $\mathsf{R}$ that satisfies
\begin{align*}
    \mathsf{H}(\sfx',\sfs') - \mathsf{H}(\sfx,\sfs) \in \{0,1\},
    \quad
    \text{whenever } (\sfx',\sfs') \in \{(\sfx+1,\sfs), (\sfx,\sfs-1), (\sfx+1,\sfs+1)\},
\end{align*}
provided that both $(\sfx,\sfs)$ and $(\sfx',\sfs')$ are vertices of $\mathsf{R}$.

We refer to the restriction $\mathsf{h} = \mathsf{H}|_{\partial \mathsf{R}}$ as a \emph{boundary height function}, where $\mathsf{H}|_{\partial \mathsf{R}}$ means the restriction of $\mathsf{H}$ to the vertices lying on $\partial \mathsf{R}$. For any boundary height function $\mathsf{h} : \partial \mathsf{R} \rightarrow \mathbb{Z}$, let $\mathscr{G}(\mathsf{h})$ denote the set of all height functions $\mathsf{H}$ with $\mathsf{H}|_{\partial \mathsf{R}} = \mathsf{h}$. In what follows, any height function on a domain $\mathsf{R} \subseteq \mathbb{R}^2$ will always be implicitly extended by linearity to the faces of $\mathsf{R}$, so that it may be viewed as a piecewise linear function on $\mathsf{R}$.

For a fixed vertex $(\sfx_0,\sfs_0)$ of $\mathsf{R}$ and an integer $h_0 \in \mathbb{Z}$, one can associate with any tiling of $\mathsf{R}$ a height function $\mathsf{H}$ as follows. First, set $\mathsf{H}(\sfx_0,\sfs_0) = h_0$, and then define $\mathsf{H}$ at the remaining vertices of $\mathsf{R}$ in such a way that the height values on the four vertices of any lozenge in the tiling are of the form depicted on the right side of \Cref{tilinghexagon}. In particular, we require that
\begin{align}\label{e:height1}
\mathsf{H}(\sfx,\sfs) - \mathsf{H}(\sfx,\sfs+1) = 1
\end{align}
if and only if $(\sfx,\sfs)$ and $(\sfx,\sfs+1)$ are vertices of the same type II lozenge, and that
\begin{align}\label{e:height2}
\mathsf{H}(\sfx+1,\sfs) = \mathsf{H}(\sfx,\sfs)
\end{align}
if and only if $(\sfx,\sfs)$ and $(\sfx+1,\sfs)$ are vertices of the same type III lozenge. Since $\mathsf{R}$ is simply connected, a height function on $\mathsf{R}$ is uniquely determined by these conditions, together with the value $\mathsf{H}(\sfx_0,\sfs_0) = h_0$.

We refer to the middle panel of \Cref{tilinghexagon} for an example. As depicted there, we can also view a lozenge tiling of $\mathsf{R}$ as a family of nonintersecting Bernoulli paths after ignoring type III lozenges. In this case, the value $\mathsf{H}(\sfx,\sfs)$ of the height function at a vertex $(\sfx,\sfs) \in \mathsf{R}$ denotes the number of Bernoulli paths to the left of $(\sfx,\sfs)$. Observe in particular that, if a tiling $\mathscr{M}$ of $\mathsf{R}$ is associated with a height function $\mathsf{H}$, then the boundary height function $\mathsf{h} = \mathsf{H}|_{\partial \mathsf{R}}$ is independent of $\mathscr{M}$ and is uniquely determined by $\mathsf{R}$, up to a global additive constant. This constant is fixed by the normalization $\mathsf{H}(\sfx_0,\sfs_0) = h_0$.

		\subsection{Limit shapes and liquid region}

	\label{HeightLimit} 
	
	To analyze the limits of height functions of random tilings, it will be useful to introduce continuum analogs of the notions considered in \Cref{Tiling}. So, set 
	\begin{flalign} \label{t}
	 \mathcal{T} = \big\{ (s, t) \in (0,1) \times (-1,0): s+t>0 \big\} \subset \mathbb{R}^2, 
	\end{flalign} 

	\noindent and its closure $\overline{\mathcal{T}} = \big\{ (s, t) \in [0,1] \times [-1,0]: s+t\geq 0 \big\}$. We interpret $\overline{\mathcal{T}}$ as the set of possible gradients, also called \emph{slopes}, for a continuum height function; $\mathcal{T}$ is then the set of ``non-frozen'' or ``liquid'' slopes, whose associated tilings contain tiles of all types. For any simply-connected subset $\mathfrak{R} \subset \mathbb{R}^2$, we say that a function $H : \mathfrak{R} \rightarrow \mathbb{R}$ is \emph{admissible} if $H$ is $1$-Lipschitz and $\nabla H(u) \in \overline{\mathcal{T}}$ for almost all $u \in \mathfrak{R}$. We further say a function $h: \partial \mathfrak{R} \rightarrow \mathbb{R}$ \emph{admits an admissible extension to $\mathfrak{R}$} if $\Adm (\mathfrak{R}; h)$, the set of admissible functions $H: \mathfrak{R} \rightarrow \mathbb{R}$ with $H |_{\partial \mathfrak{R}} = h$, is not empty.
	
	We say that a sequence of domains $\mathsf{R}_1, \mathsf{R}_2, \ldots \subset \mathbb{T}$ \emph{converges} to a simply-connected subset $\mathfrak{R} \subset \mathbb{R}^2$ if $n^{-1} \mathsf{R}_n \subseteq \mathfrak{R}$ for each $n \geq 1$ and  $\lim_{n \rightarrow \infty} \dist (n^{-1} \del\mathsf{R}_n, \del\mathfrak{R}) = 0$. We further say that a sequence $\mathsf{h}_1, \mathsf{h}_2, \ldots $ of boundary height functions on $\mathsf{R}_1, \mathsf{R}_2, \ldots $, respectively, \emph{converges} to a boundary height function $h : \partial \mathfrak{R} \rightarrow \mathbb{R}$ if $\lim_{n \rightarrow \infty} n^{-1} \mathsf{h}_n (n v_n) = h (v)$ if $v_n$ is any point in $n^{-1} \del\mathsf{R}_n$ with limit $v\in \mathfrak \del \fR$.  
	
	To state results on the limiting height function of random tilings, for any $x \in \mathbb{R}_{\geq 0}$ and $(s, t) \in \overline{\mathcal{T}}$ we denote the \emph{Lobachevsky function} $L: \mathbb{R}_{\geq 0} \rightarrow \mathbb{R}$ and the \emph{surface tension} $\sigma : \overline{\mathcal{T}} \rightarrow \mathbb{R}$ by 
	\begin{flalign}
		\label{sigmal} 
		L(x) = - \displaystyle\int_0^x  \ln |2 \sin z| \mathrm{d} z; \qquad \sigma (s, t) = \displaystyle\frac{1}{\pi} \Big( L(\pi (1-s)) + L (- \pi t) + L \big( \pi ( s + t) \big) \Big).
	\end{flalign}
	
	\noindent For any $H \in \Adm (\mathfrak{R})$, we further denote the \emph{entropy functional}
	\begin{flalign}
		\label{efunctionh} 
		\mathcal{E} (H) = \displaystyle\int_{\mathfrak{R}} \sigma \big( \nabla H (z) \big) \mathrm{d}z.
	\end{flalign}
	
	The following variational principle of \cite{cohn2001variational}  states that the height function associated with a uniformly random tiling of a sequence of domains corresponding to $\mathfrak{R}$ converges to the maximizer of $\mathcal{E}$ with high probability.

	\begin{theorem}[{\cite[Theorem 1.1]{cohn2001variational}}]
		
		\label{hzh} 
		
		Let $\mathsf{R}_1, \mathsf{R}_2, \ldots \subset \mathbb{T}$ denote a sequence of tileable domains, with associated boundary height functions $\mathsf{h}_1, \mathsf{h}_2, \ldots $, respectively. Assume that they converge to a simply-connected subset $\mathfrak{R} \subset \mathbb{R}^2$ with piecewise smooth boundary, and a boundary height function $h : \partial \mathfrak{R} \rightarrow \mathbb{R}$, respectively. Denoting the height function associated with a uniformly random tiling of $\mathsf{R}_n$ with boundary height function $\mathsf h_n$ by $\mathsf{H}_n$, we have
		\begin{flalign*}
			\displaystyle\lim_{n \rightarrow \infty} \mathbb{P} \bigg( \displaystyle\max_{\mathsf{v} \in \mathsf{R}_n} \big| n^{-1} \mathsf{H}_n (\mathsf{v}) - H^* (n^{-1} \mathsf{v}) \big| > \varepsilon \bigg) = 0,
		\end{flalign*}  
	
		\noindent where $H^*$ is the unique maximizer of $\mathcal{E}$ on $\mathfrak{R}$ with boundary data $h$,
		\begin{flalign}
			\label{hmaximum}
			H^* = \displaystyle\argmax_{H \in \Adm (\mathfrak{R}; h)} \mathcal{E} (H).
		\end{flalign}
	\end{theorem} 

	\noindent The fact that there is a unique maximizer described as in \eqref{hmaximum} follows from \cite[Section 2]{cohn2001variational} and \cite[Proposition 4.5]{de2010minimizers}. 
	Under a suitable change of coordinates, this maximizer $H^*$ solves a complex variant of the Burgers equation \cite{kenyon2007limit}, which makes it amenable to further analysis; we will discuss this point in more detail in \Cref{Slopeft} below.
	 %and the \emph{arctic boundary} $\mathfrak{A} = \mathfrak{A} (\mathfrak{R}; h) \subset \overline{\mathfrak{R}}$ as the set of points $ u=(x,t) \in \del \fL$, such that for any sequence of points $u_n \in \fL$ converging to $u$,
%\begin{align}\begin{split}
%		\label{e:arctic} 
%		   \big( \partial_x H^* (u_n), \partial_t H^* (u_n) \big) \rightarrow \del\mathcal{T}, 
%\end{split}	\end{align}
%	where $H^*$ is as in \eqref{hmaximum}. The complement of the liquid region $\fR\setminus \fL$ is called the \emph{frozen region}. 

	\begin{figure}
	\begin{tikzpicture}[scale=3]
  % Axes
 \draw[->, thick, >=stealth, line width=1pt] (-1.5,0) -- (0.5,0); % x-axis
  \draw[->, thick, >=stealth, line width=1pt] (0,-1.2) -- (0,0.2); % y-axis

  % Triangle
  \coordinate (A) at (-1,0);
  \coordinate (B) at (0,0);
  \coordinate (C) at (-0.5,-1);
  \draw[thick] (A) -- (B) -- (C) -- cycle;

  % Labels
  \filldraw (A) circle (1pt) node[above left] {$-1$};
  \filldraw (B) circle (1pt) node[above right] {$0$};
  \filldraw (C) circle (1pt) node[below] {$f$};

     \definecolor{bluetri}{rgb}{0.68, 0.85, 0.9}
      \fill[bluetri] (-0.66,-0.23) -- (-0.66,-0.03) -- (-0.86,-0.23) -- cycle;
      \fill[bluetri](-0.14,-0.23) -- (-0.14,-0.03) -- (-0.34,-0.23) -- cycle;
       \fill[bluetri](-0.4,-0.73) -- (-0.4,-0.53) --(-0.6,-0.73) -- cycle;
     
  % Shaded parallelograms and squares
   \draw[thick] (-0.46,-0.03) -- (-0.66,-0.03) -- (-0.86,-0.23)--(-0.66,-0.23) -- cycle;

  \draw[thick] (-0.14,-0.23) -- (-0.14,-0.03) -- (-0.34,-0.23)--(-0.34,-0.43) -- cycle;

  \draw[thick] (-0.4,-0.73) -- (-0.4,-0.53) -- (-0.6,-0.53)--(-0.6,-0.73) -- cycle;

\end{tikzpicture}

%	\begin{center}
%	 \includegraphics[scale=0.3,trim={0cm 5cm 0 7cm},clip]{complex_slope.pdf}
	 \caption{Shown above the complex slope $f = f (u)$.}
	 \label{slope1}
	% \end{center}
	 \end{figure}
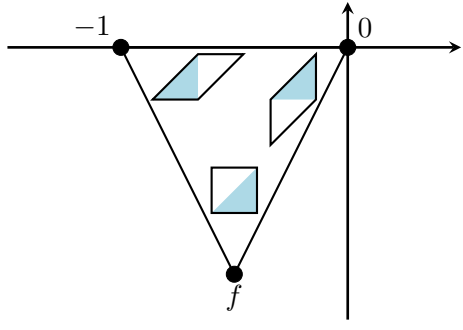

	\subsection{Kenyon–Okounkov conformal structure }
	\label{Slopeft}

	 For any simply-connected open subset $\mathfrak{R} \subset \mathbb{R}^2$ with Lipschitz boundary, and boundary height function $h: \partial \mathfrak{R} \rightarrow \mathbb{R}^2$ admitting an admissible extension to $\mathfrak{R}$, define the \emph{liquid region} $\mathfrak{L} = \mathfrak{L} (\mathfrak{R}; h) \subset \mathfrak{R}$ 
	\begin{align}\begin{split}
		\label{al} 
		&\mathfrak{L} = \big\{ u=(x,s) \in \mathfrak{R}: \big( \partial_x H^* (u), \partial_s H^* (u) \big) \in \mathcal{T} \big\}.
\end{split}	\end{align}

In this section we recall from \cite{kenyon2007limit,astala2026dimer} various complex analytic aspects of the tiling limit shapes. We define the \emph{complex slope} $f: \mathfrak{L} \rightarrow \bC_-$ by, for any $u \in \mathfrak{L}$, setting $f(u) \in \bC_-$ to be the unique complex number satisfying 
	\begin{flalign}
		\label{fh}
		\arg f(u) = - \pi \partial_x H^* (u); \qquad \arg \big( f(u) + 1 \big) = \pi \partial_t H^* (u),
	\end{flalign}

	\noindent where for any $z \in \bC\setminus \{0\}$ we have set $\arg z \in (-\pi,\pi]$; see \Cref{slope1} for a depiction, where there we interpret $1 - \partial_x H^* (u)$ and $-\partial_t H^* (u)$ as the approximate proportions of tiles of types III and II around $u \in \fR$, respectively (which follows from the definition of the height function from \Cref{Tiling}).   
	We note that our notation differs from that of \cite{kenyon2007limit}: their height function and complex slopes $(h, w, z)$ correspond, in our notation, to $(x - H^*,\, -f,\, 1 + f)$.

		 The following result from \cite{kenyon2007limit} indicates that the complex slope $f$ satisfies the complex Burgers equation on the liquid region. Moreover, it can be solved by complex characteristics. 
	  \begin{proposition}[{\cite[Theorem 1, Corollary 1]{kenyon2007limit}}]
	 	
	 	\label{fequation}
	 	
	 	For any $(x, s) \in \mathfrak{L}$,  $f_s (x) = f (x, s)$ is analytic on the liquid region $\fL$ as defined in \eqref{al}. And it satisfies the following complex Burgers equation 
	 	\begin{flalign}
	 		\label{ftx}
	 		\partial_s f_s (x) + \partial_x f_s (x) \displaystyle\frac{f_s (x)}{f_s (x) + 1} = 0.
	 	\end{flalign} 
Moreover, there exists an analytic function $Q$ of two variables such that
\begin{align}\label{e:analyticQ}
Q\left(f_s (x), x-s \frac{f_s (x)}{f_s (x)+1}\right)=0.
\end{align}	
	 	\end{proposition} 

The complex slope $f$ endows the liquid region $\fL$ with a natural complex structure, often referred to as the \emph{Kenyon--Okounkov conformal structure}.  
It was conjectured in \cite{kenyon2007limit} that, for random lozenge tilings, the fluctuations of the height function converge to the Gaussian Free Field (GFF) on $\fL$ with respect to this complex structure.  
We state the conjecture for uniform random lozenge tilings, and refer to \cite{gorin2021lectures} for further discussion and heuristics.

\begin{conjecture}[{\cite{kenyon2007limit,kenyon2008height}}]
Adopt the setting of \Cref{hzh}.  
In the liquid region $(x,s) \in \fL$, the centered height function satisfies
\begin{equation}
 \mathsf{H}_n(nx, ns) - \mathbb{E}[\mathsf{H}_n(nx, ns)]
 \xrightarrow[n \to \infty]{\;d\;}
 \mathrm{GFF},
\end{equation}
where $\mathrm{GFF}$ denotes the Gaussian Free Field on $\fL$ with Dirichlet boundary conditions, taken with respect to the complex structure induced by $f$.
\end{conjecture}

One way to describe the GFF on $\fL$ is via the analytic function $Q$ from \eqref{e:analyticQ}.  
In the liquid region, where $\Im[f_s (x)] < 0$, the map
\[
(x,s) \longmapsto \left( f_s (x),\, x - s\frac{\, f_s (x)}{f_s (x) + 1} \right)
\]
is a bijection between $\fL$ and a subset of the algebraic curve $Q = 0$ in $\mathbb{C}^2$.  
This subset, viewed as a complex curve, has a natural local coordinate system and hence supports a well-defined Gaussian Free Field.  
Pulling back this field via the above map yields the GFF on $\fL$ with respect to the Kenyon--Okounkov complex structure. We refer to \Cref{s:surface} for a more detailed discussion.

 \subsection{Polygonal domains}
 	
In this article, we focus on the case where $\mathfrak{R}$ is a polygonal domain, described as follows.  
In general, any region can be approximated by a polygonal domain; in this setting, it was shown in \cite{kenyon2007limit,astala2026dimer} that the associated analytic function $Q$ from \eqref{e:analyticQ} is algebraic and enjoys many benign properties.

 	\begin{definition} 		
 	\label{p}  	
	A subset \( \mathfrak{P} \subset \mathbb{R}^2 \) is called \emph{polygonal} if it is a simply connected polygon formed by $3d$ segments
with slopes $0, 1, \infty$ cyclically repeated as we follow the boundary in the
counterclockwise direction (see \Cref{f:lattice} for a 9-gon and \Cref{f:12-gon} for a 12-gon). We assume the domain\footnote{Throughout, we assume that all vertices of \(n \mathfrak{P}\) lie in \(\bZ'\times \mathbb{Z}\).}  
\[
\mathsf{P}_n := n \overline{\mathfrak{P}} \cap \mathbb{T}
\]
is tileable, and thus admits an associated boundary height function \(\sfh=\sfh_n\).  
Along vertical ($\infty$-slope) or unit-slope edges $\sfh$ is constant, while along horizontal edges (\(0\)-slope) it increases with unit slope in the \(x\)-direction; in particular, moving one unit to the right increases \(\sfh\) by one.  
The function \(\sfh\) is unique up to a global shift. 

With this normalization, we define a rescaled boundary height function
\[
h: \partial \mathfrak{P} \to \mathbb{R}, \quad h(u) := n^{-1} \mathsf{h}(nu), \quad \text{for each } u \in \partial \mathfrak{P}.
\]
We abbreviate the admissible height function class as \( \Adm(\mathfrak{P}) := \Adm(\mathfrak{P}; h) \) and the associated liquid region as \( \mathfrak{L}(\mathfrak{P}) := \mathfrak{L}(\mathfrak{P}; h) \).
The \emph{arctic boundary} is defined as the boundary of the liquid region:
\[
\mathfrak{A}(\mathfrak{P}) := \partial \mathfrak{L}(\mathfrak{P}).
\]
The complement \( \overline{\mathfrak{P}} \setminus \mathfrak{L}(\mathfrak{P}) \) is referred to as the \emph{frozen region}.
Both the liquid and frozen regions are independent of the particular global shift chosen to fix \( h \). Finally, we define the maximizer \( H^* \in \Adm(\mathfrak{P}; h) \) as in equation~\eqref{hmaximum}.

 	\end{definition}

	We will make use of the following results from \cite{kenyon2007limit,astala2026dimer} on the behavior of the limit shape $H^*$ and arctic boundary $\mathfrak{A}$ when $\mathfrak{R}=\fP$ is polygonal. The first statement in the below lemma is given by \cite[Theorem 7.5]{astala2026dimer} and the second by \cite[Theorem 6.13]{astala2026dimer} (see also \cite[Theorem 2]{kenyon2007limit}).
 	
 	\begin{lemma}[{\cite{kenyon2007limit,astala2026dimer}}]
	
	\label{pla}
	
	Adopt the notation of \Cref{p}, and assume that the domain $\mathfrak{R} = \mathfrak{P}$ is polygonal with at least $6$ sides. Then  the following two statements hold.
	
	\begin{enumerate}
		\item The arctic boundary $\fA$ is frozen in the sense that for any sequence $\{\zeta_i\}_{i\geq 1}\in \fL$ approaching $\zeta\in \fA$, we have $\nabla H^*(\zeta_i)\rightarrow \del \cT$.
		\item Each connected component of the liquid region $\fL$ is simply connected.

	\end{enumerate}
 
	\end{lemma}

\subsection{Main results}
Given the polygonal domain $\fP$ as in \Cref{p}, we recall the obstacle functions introduced in \cite{de2010minimizers,astala2026dimer}.  
There exist two functions $m, M \in \Adm(\fP; h)$ such that
\begin{flalign*} 
    m(u) \leq H(u) \leq M(u), \qquad \text{for all $H \in \Adm(\fP)$ and $u\in \fP$}.
\end{flalign*} 
Both $m$ and $M$ are piecewise linear, with gradients taking values in the set of vertices of $\mathcal{T}$, namely
$
\{ (0,0),\ (1,0),\ (1,-1) \}$, see \cite[Theorem 8.2]{astala2026dimer}.
The \emph{trivial set} is the set of points where the lower and upper obstacles coincide:
\begin{align}\label{e:trivialset}
\{ z \in \fP : m(z) = M(z) \}.
\end{align}

\begin{assumption} \label{a:asump}
Under the notation of \Cref{p}, we make the following assumptions:
\begin{enumerate}
    \item The liquid region is non-empty and connected.
  %  \item Any intersection point between $\mathfrak{A}$ and $\partial \mathfrak{P}$ must be a tangency location of $\mathfrak{A}$. %Moreover, $\nabla H^*(x,t)$ is continuous at any point on $\mathfrak{A}$ that is not a tangency location.
    \item The trivial set \eqref{e:trivialset} associated with the polygonal domain $\fP$ is empty. In other words, for every $v \in \fP$, one has $m(v) < M(v)$.
\end{enumerate}
\end{assumption}

\begin{remark}
In general, the liquid region may have several connected components. In this
case, we expect the fluctuations of the height function on different components
to converge to independent Gaussian Free Fields. We postpone a detailed study
of this situation to future work.

Furthermore, for any $H \in \Adm(\fP)$, one has $m(v)=M(v)=H(v)$
at each point \(v\) of the trivial set. Thus, the height is completely
determined there; in this sense, the lozenge tiling is frozen on the trivial
set. The tiling problem may therefore be reduced to the connected components
of its complement. These components are polygonal domains with sides parallel
to the three lattice directions, but they need not belong to the same
nondegenerate cyclic \(3d\)-gon class: some of the directions may occur only
through zero-length, degenerate sides.
\end{remark}

Given a polygonal domain $\fP$ as in \Cref{p}.		
We denote the height function associated with a uniformly random tiling of $\mathsf{P}_n:=n \overline{\mathfrak{P}} \cap \mathbb{T}$ with boundary height function $\mathsf h_n$ by $\mathsf{H}_n$.
We introduce the rescale height function, and the centered version as
\begin{align}\label{e:rescaled_height}
H_n(u)=n^{-1}\sfH_n(nu),\quad H_n^\circ(u)=H_n(u)-\bE[H_n(u)],\quad u\in \fP
\end{align}

To state our result, we first recall the Dirichlet Green function on the upper
half-plane. For $z,w\in \bC_+$, define
\begin{equation}\label{eq:Green-def}
    \cG(z,w)
    :=
    \frac{1}{2\pi}
     \ln\left|\frac{z-\overline w}{z-w}\right|.
\end{equation}
This is the covariance kernel of the Gaussian free field on $\bC_+$ with
Dirichlet boundary conditions. More precisely, if $\mathrm{GFF}_{\bC_+}$
denotes this field, then, in the usual formal notation, for distinct
$z,w\in\bC_+$,
\[
    \mathrm{Cov}\bigl(
        \mathrm{GFF}_{\bC_+}(z),
        \mathrm{GFF}_{\bC_+}(w)
    \bigr)
    =
    \cG(z,w).
\]
Thus, for any $k\geq 1$ and pairwise distinct points
$w_1,\ldots,w_k\in \bC_+$, Wick's formula gives
\[
    \bE\!\left[
        \mathrm{GFF}_{\bC_+}(w_1)\cdots
        \mathrm{GFF}_{\bC_+}(w_k)
    \right]
    =
    \begin{cases}
    \displaystyle
    \sum_{\pi\in\mathcal P_k}
    \prod_{\{a,b\}\in\pi}
        \cG(w_a,w_b),
    & k \text{ even},\\[2ex]
    0,
    & k \text{ odd},
    \end{cases}
\]
where $\mathcal P_k$ denotes the set of pairings of $\{1,\ldots,k\}$.

We prove the following convergence of moments of the height fluctuations to the
corresponding Wick moments of the Gaussian free field.

\begin{theorem}\label{thm:moment-GFF}
Let $\fP$ be a polygonal domain as in \Cref{p}, and suppose that
\Cref{a:asump} holds. Let $\fL$ be the liquid region, and let
\[
    \phi : \fL \to \bC_+
\]
be a uniformizing conformal map, where conformality is with respect to the
complex structure induced by the complex slope $f$ (see \Cref{s:cf_map} for more details). Fix $k\geq 1$ and
pairwise distinct points $
    (x_1,s_1),\ldots,(x_k,s_k)\in \fL$.
Set
\[
    w_j := \phi(x_j,s_j)\in \bC_+,
    \qquad j=1,\ldots,k.
\]
Then the centered height function $H_n^\circ$, defined in
\eqref{e:rescaled_height}, satisfies
\begin{equation}\label{eq:moment-conv}
    \lim_{n\to\infty}
    \pi^{k/2}\,
    \bE\!\left[
        H_n^\circ(x_1,s_1)\cdots H_n^\circ(x_k,s_k)
    \right]
    =
    \bE\!\left[
        \mathrm{GFF}_{\bC_+}(w_1)\cdots \mathrm{GFF}_{\bC_+}(w_k)
    \right].
\end{equation}
\end{theorem}

%\begin{remark}
%\Cref{thm:moment-GFF} proves joint moment convergence at pairwise distinct
%macroscopic bulk points. Convergence of the height field as a random
%distribution additionally requires a separate near-diagonal estimate:
%\begin{align}
%\bE\!\left[
%        |H_n^\circ(x,s)|^k
%    \right]\leq n^{\delta k}
%\end{align}
%for some small $\delta>0$.
% This type of 
%estimate can potentially be obtained from the kernel estimates in
%\Cref{subsec:edge-facet} by following the arguments in
%\cite{petrov2015asymptotics,borodin2014anisotropic}.
%\end{remark}

\section{Proof Strategy}
\label{s:outline} 
\subsection{Kasteleyn theory}\label{s:Kasteleyn}

In this section, we recall Kasteleyn theory for the uniform dimer model on the
hexagonal lattice. For a more general treatment, we refer the reader to
\cite[Section 3]{kenyon2009lectures}.

A \emph{Kasteleyn weighting} of a planar bipartite graph is an assignment of a
sign to each edge such that every bounded face whose number of edges is
congruent to \(0\pmod 4\) has an odd number of negative signs, whereas every
bounded face whose number of edges is congruent to \(2\pmod 4\) has an even
number of negative signs.

Recall that the dual graph of \(\bT\) is the hexagonal lattice: each vertex of
the dual graph corresponds to a triangle of \(\bT\), and two vertices are
adjacent if the corresponding triangles share an edge. Every face of the dual
graph is a hexagon and hence has six edges. We fix the following Kasteleyn
weighting.

Consider a white triangle of the form
\[
\big\{(\sfx-\tfrac12,\sfs),\,
      (\sfx-\tfrac12,\sfs+1),\,
      (\sfx+\tfrac12,\sfs+1)\big\},
\qquad \sfx,\sfs\in\bZ.
\]
It is adjacent to the following three blue triangles:
\begin{align*}
&\big\{(\sfx-\tfrac32,\sfs),\,
       (\sfx-\tfrac12,\sfs),\,
       (\sfx-\tfrac12,\sfs+1)\big\}, \\
&\big\{(\sfx-\tfrac12,\sfs),\,
       (\sfx+\tfrac12,\sfs),\,
       (\sfx+\tfrac12,\sfs+1)\big\}, \\
&\big\{(\sfx-\tfrac12,\sfs+1),\,
       (\sfx+\tfrac12,\sfs+1),\,
       (\sfx+\tfrac12,\sfs+2)\big\}.
\end{align*}
We assign Kasteleyn weights to the three dual edges incident to this white
triangle as follows; see the right panel of
\Cref{fig:triangular_lattice}:
\begin{itemize}
    \item the two edges connecting it to
    \[
    \big\{(\sfx-\tfrac32,\sfs),\,
          (\sfx-\tfrac12,\sfs),\,
          (\sfx-\tfrac12,\sfs+1)\big\}
 ,\quad
    \big\{(\sfx-\tfrac12,\sfs),\,
          (\sfx+\tfrac12,\sfs),\,
          (\sfx+\tfrac12,\sfs+1)\big\}
    \]
    have weight \(-1\);

    \item the edge connecting it to
    \[
    \big\{(\sfx-\tfrac12,\sfs+1),\,
          (\sfx+\tfrac12,\sfs+1),\,
          (\sfx+\tfrac12,\sfs+2)\big\}
    \]
    has weight \(1\).
\end{itemize}
Each hexagonal face has four edges of weight \(-1\) and two edges of weight
\(1\), so this is indeed a Kasteleyn weighting.

For any finite, simply connected, tileable region
\(\mathsf R\subseteq\bT\), the \emph{Kasteleyn matrix} is the signed adjacency
matrix of the dual graph restricted to \(\mathsf R\). Let \(\bW\) and \(\bB\)
denote the sets of white and blue triangles in \(\mathsf R\), respectively.
Using the Kasteleyn weighting defined above, we define
\[
K\in\bR^{|\bW|\times|\bB|}
\]
by
\[
K({\rw},{\rb})=
\begin{cases}
0, & \text{if \(\rw\) and \(\rb\) are not adjacent},\\
\text{the Kasteleyn weight of the edge \(\rw\rb\)},
   & \text{if \(\rw\) and \(\rb\) are adjacent},
\end{cases}
\]
for \(\rw\in\bW\) and \(\rb\in\bB\).

Since \(\mathsf R\) is tileable, \(K\) is square. Moreover, Kasteleyn theory
implies that \(\lvert\det K\rvert\) equals the number of dimer coverings of
\(\mathsf R\). In particular, \(K\) is invertible. Its inverse is characterized
by
\begin{align}
\delta_{\rw',\rw}
&=\sum_{\rb\in\bB}K(\rw',\rb)K^{-1}(\rb,\rw).
\end{align}

Kasteleyn theory also implies that the edges in a uniformly random dimer
covering form a determinantal point process whose correlation kernel is
expressed in terms of \(K^{-1}\).

\begin{theorem}\label{t:kasteleyn_kernel}
Let
\[
\bX=\{\rw_1\rb_1,\rw_2\rb_2,\ldots,\rw_k\rb_k\}
\]
be a collection of distinct edges in the dual graph of \(\mathsf R\), where
\(\rw_i\in\bW\) and \(\rb_i\in\bB\). Then the probability that all the edges in
\(\bX\) occur in a uniformly random dimer covering of \(\mathsf R\) is
\begin{align}
\left(\prod_{i=1}^k K(\rw_i,\rb_i)\right)
\det\!\left(K^{-1}(\rb_i,\rw_j)\right)_{1\le i,j\le k}.
\end{align}
\end{theorem}

We now specialize to \(\mathsf R=\sfP_n\). We extend \(K^{-1}(\rb,\rw)\) by zero
whenever the blue triangle \(\rb\) is not contained in \(\sfP_n\). We use the
following terminology for the exterior blue triangles adjacent to the boundary.

\begin{figure}
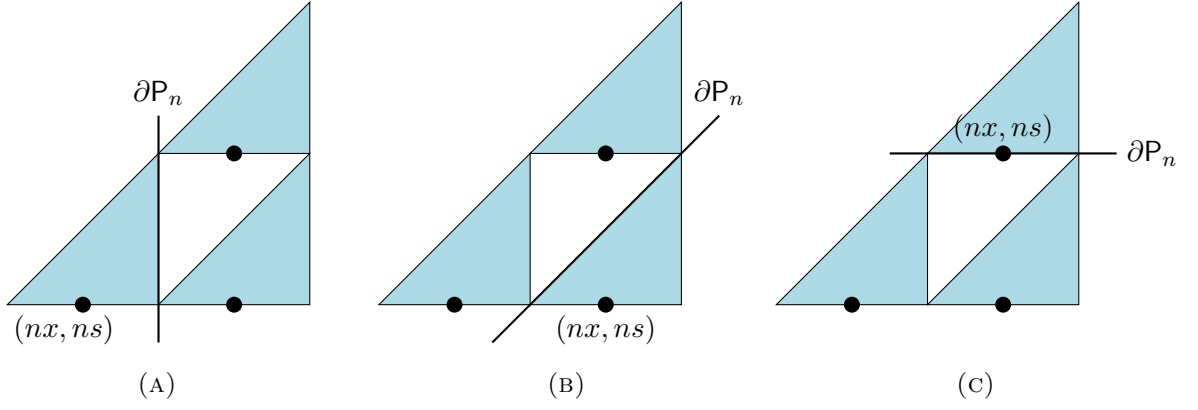


\begin{subfigure}[t]{0.32\textwidth}
\centering

\caption{}
\end{subfigure}
\begin{subfigure}[t]{0.32\textwidth}
\centering
%
\caption{}
\end{subfigure}
\begin{subfigure}[t]{0.32\textwidth}
\centering
%
\caption{}
\end{subfigure}
\caption{Boundary blue triangle}
\label{f:bb_triangle}
\end{figure}

\begin{definition}[Blue boundary triangles]\label{def:boundary_triangle}
A blue triangle \(\rb\) not contained in \(\sfP_n\) is called a
\emph{blue boundary triangle} if one of the following conditions holds:
\begin{itemize}
    \item its vertical edge lies on a vertical side of \(\sfP_n\); see
    panel (A) of \Cref{f:bb_triangle};

    \item its unit-slope edge lies on a side of \(\sfP_n\) of slope \(1\);
    see panel (B) of \Cref{f:bb_triangle};

    \item its horizontal edge lies on a horizontal side of \(\sfP_n\); see
    panel (C) of \Cref{f:bb_triangle}.
\end{itemize}
\end{definition}

We recall our convention for encoding triangular faces by lattice points.
Before rescaling, a blue triangle
\[
\big\{(\sfx-\tfrac12,\sfs),\,
      (\sfx+\tfrac12,\sfs),\,
      (\sfx+\tfrac12,\sfs+1)\big\}
\]
is encoded by the midpoint of its bottom edge, namely, by the point
\((\sfx,\sfs)\in\bZ^2\). Similarly, a white triangle
\[
\big\{(\sfx-\tfrac12,\sfs-1),\,
      (\sfx-\tfrac12,\sfs),\,
      (\sfx+\tfrac12,\sfs)\big\}
\]
is encoded by the midpoint of its top edge, again by the point
\((\sfx,\sfs)\in\bZ^2\); see \Cref{fig:triangular_lattice}.

We then rescale both coordinates by \(1/n\). Thus, the unscaled label
\((\sfx,\sfs)\in\bZ^2\) is identified with the rescaled point
$
(\sfx/n,\sfs/n)\in\bZ^2/n.
$
In particular, a unit change in either unscaled coordinate becomes a change of
size \(1/n\) in the corresponding rescaled coordinate.

With this convention, if a white triangle has rescaled coordinate
\(\rw'=(x,s)\in\bZ^2/n\), then its three adjacent blue triangles have rescaled
coordinates
\[
(x,s),\qquad
\left(x-\frac1n,s-\frac1n\right),\qquad
\left(x,s-\frac1n\right).
\]
The corresponding Kasteleyn weights are \(1,-1,-1\), respectively. Therefore,
for every \(w\in\bW\),
\begin{align}\label{e:discrete_kernel}
\delta_{\rw',\rw}
&=\sum_{b\in\bB}K(\rw',\rb)K^{-1}(\rb,\rw)\notag\\
&=K^{-1}((x,s),\rw)
  -K^{-1}\left(\left(x-\frac1n,s-\frac1n\right),\rw\right)
  -K^{-1}\left(\left(x,s-\frac1n\right),\rw\right),
\end{align}
where the first arguments on the right-hand side denote blue triangles in
rescaled coordinates. If one of these triangles is not contained in
\(\sfP_n\), the corresponding term is understood to be zero.

\begin{proposition}\label{p:inverse_kasteleyn_recursion}
Let \(\fP\) be a polygonal domain as in \Cref{p}, and let
\(\sfP_n=n\fP\). Assume that \(\sfP_n\) is tileable, and let \(K\) be the
Kasteleyn matrix defined above. We use rescaled coordinates, so that an
unscaled lattice point \((nx,ns)\in\bZ^2\) is written as
\((x,s)\in\bZ^2/n\).

Then \(K^{-1}\), together with its zero extension to blue boundary triangles,
is uniquely characterized by the following conditions.
\begin{enumerate}
\item If \((x,s)\) corresponds to a blue boundary triangle as in \eqref{def:boundary_triangle}, then
\begin{align}\label{e:zero_bb}
K^{-1}((x,s),(y,t))=0
\end{align}
for every white triangle in \(\sfP_n\) with rescaled coordinate \((y,t)\).

\item For any two white triangles in \(\sfP_n\) with rescaled coordinates
$
(x,s),\ (y,t)\in\bZ^2/n,
$
we have
\begin{align}\label{e:discrete_kernel2}
\delta_{(x,s),(y,t)}
&=
K^{-1}\bigl((x,s),(y,t)\bigr)
-K^{-1}\left(
    \left(x-\frac1n,s-\frac1n\right),(y,t)
  \right)
-K^{-1}\left(
    \left(x,s-\frac1n\right),(y,t)
  \right).
\end{align}
Here each first argument of \(K^{-1}\) denotes a blue triangle, whereas the
second argument denotes a white triangle.
\end{enumerate}
\end{proposition}

\begin{proof}
The first condition is precisely the zero-extension convention, and
\eqref{e:discrete_kernel2} follows from \eqref{e:discrete_kernel}.

To prove uniqueness, fix a white triangle \(\rw\in\bW\) and set
\[
\bmv(\rb)=K^{-1}(\rb,\rw),\qquad \rb\in\bB.
\]
After setting the values at exterior blue triangles equal to zero, the equations
in \eqref{e:discrete_kernel2} are exactly
$
K \bmv=\delta_\rw.
$
Since \(K\) is invertible, this system has the unique solution
\(\bmv=K^{-1}\delta_\rw\). Hence the stated conditions uniquely determine
\(K^{-1}\).
\end{proof}

\begin{remark}\label{r:walk}
If we ignore the boundary conditions, a causal fundamental solution of
\eqref{e:discrete_kernel2} is given by the following unnormalized discrete heat
kernel. For
$
(x,s),(y,t)\in\bZ^2/n,
$
define
\begin{align}\label{e:discrete_heat_kernel}
p\bigl((x,s),(y,t)\bigr)
&:=
{n(s-t)\choose n(x-y)}
\end{align}
By Pascal's identity,
\begin{align*}
\delta_{(x,s),(y,t)}
&=
p\bigl((x,s),(y,t)\bigr)
-p\left(
    \left(x-\frac1n,s-\frac1n\right),(y,t)
  \right)
-p\left(
    \left(x,s-\frac1n\right),(y,t)
  \right).
\end{align*}
Thus \(p\) satisfies the same inhomogeneous recursion \eqref{e:discrete_kernel2} as \(K^{-1}\).
\end{remark}

\subsection{Proof ideas and outline}

Our method combines the inverse-Kasteleyn viewpoint with
contour-integral steepest descent. The proof can be summarized in four
main steps. First, we use the complex Burgers equation to identify a tiling 
action function that governs the asymptotics of the model. Second, we
construct local contour integral approximations, or parametrices, for
the inverse Kasteleyn matrix in each relevant geometric regime. Third,
we prove that these local parametrices are compatible on overlapping
regions and glue them into a global approximation. Finally, we use the
global kernel estimates to compute the joint cumulants of the height
function and prove convergence to the Gaussian free field.

\emph{The inverse Kasteleyn equation and the tiling action.}
In the coordinates used in this article, constructing the inverse
Kasteleyn matrix amounts to solving a discrete heat-type boundary-value
problem; see \Cref{p:inverse_kasteleyn_recursion}. This observation
suggests that the inverse Kasteleyn matrix should admit a representation
in terms of discrete heat kernels.

For lozenge tilings of trapezoidal domains, the inverse Kasteleyn
matrix has an exact double-contour integral representation, as derived
in \cite{petrov2014asymptotics,petrov2015asymptotics}. After applying
Stirling's formula, the main part of this representation takes the form
\begin{align}
\frac{n}{(2\pi \ri)^2}
\oint\oint
(\cdots)
\frac{e^{n(S(w;x,s)-S(z;y,t))}}{w-z}
\,\rd z\,\rd w,
\end{align}
where the integrals are taken over suitable contours. The asymptotic
behavior of this integral is governed by an explicit action function
\(S(z;x,s)\), whose precise form depends on the shape of the
trapezoidal domain. A natural question is therefore
\[
\boxed{
\text{What is the action function }S(z;x,s)
\text{ for a general domain?}
}
\]

We next give a heuristic derivation of this action function, which we
call the \emph{tiling action}. Recall from \Cref{Slopeft} that the
complex slope \(f=f(x,s)\), regarded as a map from the liquid region
\(\fL\) into \(\bC_-\), satisfies the complex Burgers-type equation
\begin{align}\label{e:ftxcopy0}
\del_s\ln f(x,s)
+
\del_x\ln\!\bigl(f(x,s)+1\bigr)
=0,
\end{align}
with initial condition
$
f_0(x):=f(x,0).
$

Motivated by the Cole--Hopf transform, we formally introduce a
function \(\psi\) by
\begin{align}\label{e:colehopf0}
\ln f(x,s)
=
-\frac{1}{n}\,\partial_x\ln\psi(x,s)
=
-\frac{1}{n}\,
\frac{\partial_x\psi(x,s)}{\psi(x,s)},
\qquad
\psi(x,0)
:=
e^{-n\int_0^x\ln f_0(u)\,\rd u}.
\end{align}
Substituting \eqref{e:colehopf0} into \eqref{e:ftxcopy0} and formally
expanding at the lattice scale \(1/n\), we obtain
\begin{align}\label{e:discreteheat}
\psi(x,s)
\approx
\psi\left(x,s-\frac{1}{n}\right)
+
\psi\left(x-\frac{1}{n},s-\frac{1}{n}\right),
\qquad
x,s\in\frac{1}{n}\bZ.
\end{align}
Thus, at the level of this formal approximation, \(\psi\) satisfies a
discrete heat equation.

For \(ns\in\bZ_{\geq 0}\), its solution is formally given by convolution
with the corresponding binomial heat kernel:
\begin{align}\begin{split}\label{e:kernel_represent0}
\psi(x,s)
&\approx
\sum_{z\in\bZ/n}
{ns\choose n(x-z)}
\psi(z,0)
\approx
n\int_{\bR}
{ns\choose n(x-z)}
e^{-n\int_0^z\ln f_0(u)\,\rd u}
\,\rd z\\
&\approx
n\int_{\bR}
\sqrt{
\frac{s}
{2\pi n(x-z)(s-x+z)}
}
e^{nS(z;x,s)}
\,\rd z.
\end{split}\end{align}
where the tiling action $S(z;x,s)$ is given by  
\begin{align}\label{e:tiling_action}
S(z;x,s)
&:=
s\ln s
-(x-z)\ln(x-z)
-(s-x+z)\ln(s-x+z)
-\int_0^z\ln f_0(u)\,\rd u.
\end{align}

The logarithms and the primitive appearing in
\eqref{e:tiling_action} are naturally multivalued. Consequently, the
tiling action should not be viewed merely as a function on the complex
\(z\)-plane. Instead, it lives on a Riemann surface \(\cC\) determined
by the tiling model.

More precisely, the characteristic map
\begin{align}
(x,s)
\in\fL
\longmapsto
\left(
f(x,s),
x-\frac{s f(x,s)}{f(x,s)+1}
\right)
=(f,z)
\end{align}
maps the liquid region bijectively onto one half of the Riemann surface
\(\cC\). It extends to the arctic boundary, which is mapped to the real
locus of \(\cC\). In this way, the tiling action encodes both the global
geometry of the limit shape and the local asymptotic behavior of the
inverse Kasteleyn matrix.

The phase of the tiling model is reflected in the critical points of
the tiling action. When \((x,s)\) lies in the liquid region, the
relevant critical points form a nonreal complex-conjugate pair. At a
generic point of the arctic boundary, this pair coalesces at a real
degenerate critical point. Upon entering the frozen region, the
relevant critical points become real. Further degeneracies of these
critical points describe cusps and the various tangent regimes.

The limit shape and the associated Riemann surface are reviewed in
\Cref{s:limit_shape}. The tiling action is formally introduced in
\Cref{s:tilingaction}, and its basic critical-point properties are
collected in \Cref{s:prop_critical}.

\emph{Ansatz for the inverse Kasteleyn matrix.}
We next explain the construction in the simplest case, when both
\((x,s)\) and \((y,t)\) lie in the liquid region. Introduce
\begin{align}
\begin{split}
P_{ns}(nw,nx)
&:=
{ns\choose n(x-w)}
=
\frac{\Gamma(ns+1)}
{\Gamma(n(x-w)+1)
 \Gamma(n(w-(x-s))+1)},\\
Q_{nt}(nz,ny)
&:=
\frac{
\Gamma(n(y-z))
\Gamma(n(z-(y-t)))
}
{\Gamma(nt)}.
\end{split}
\end{align}
Here \(P\) is the binomial heat kernel appearing in
\eqref{e:kernel_represent0}, while \(Q\) is its backward counterpart.
Both \(P\) and \(Q\) satisfy exact forward and backward discrete
heat-type recursions; see \Cref{c:recursion}.

The local kernel ansatz is a sum of a single-contour integral and
double-contour integrals. For brevity, we display only the
double-contour term:
\begin{align}\begin{split}\label{e:bulk_illustrate}
&\frac{n}{(2\pi\ri)^2}
\oint\rd z
\oint\rd w\,
P_{ns}(nw,nx)
e^{-n\int_0^w\ln f(u)\,\rd u}
Q_{nt}(nz,ny)
e^{n\int_0^z\ln f(u)\,\rd u}
\frac{
\sqrt{\phi'(w)}\sqrt{\phi'(z)}
}{
\phi(w)-\phi(z)
}\\
&\approx
\frac{n}{(2\pi\ri)^2}
\oint\rd z
\oint\rd w\,
(\cdots)
e^{n(S(w;x,s)-S(z;y,t))}
\frac{
\sqrt{\phi'(w)}\sqrt{\phi'(z)}
}{
\phi(w)-\phi(z)
}.
\end{split}\end{align}
Here \(\phi\) is a uniformizing conformal map from one half of the
Riemann surface to the upper half-plane. Thus,
\[
\frac{\sqrt{\phi'(w)}\sqrt{\phi'(z)}}
{\phi(w)-\phi(z)}
\]
plays the role of a Cauchy kernel on the Riemann surface.

The \(w\)- and \(z\)-contours are chosen in small local charts around
the relevant critical points. The \(w\)-contour follows a local
steepest-descent path for \( S(w;x,s)\), whereas the \(z\)-contour
follows a local steepest-ascent path for \(S(z;y,t)\). The integral
therefore localizes around the corresponding critical points. The
detailed construction of the liquid-region ansatz is given in
\Cref{s:liquid_kernel}.

The functions \(P\) and \(Q\) are important not only for the
asymptotic analysis, but also for the algebraic structure of the
kernel. Their forward and backward heat-type recursions allow the
Kasteleyn recursion to be passed through the contour integrals.
Together with the accompanying single-contour term, the resulting
ansatz satisfies the same inhomogeneous recursion
\eqref{e:discrete_kernel2} as the inverse Kasteleyn matrix. The remaining factors in the integrand are chosen so
that the appropriate boundary conditions are also satisfied.

\emph{The different geometric regimes.}
The liquid-region formula is not uniform near the arctic or polygonal boundary. In these regions, critical points may coalesce, approach the real locus of the Riemann surface, or become spurious. The local charts,
contours, and factors in the integrand must therefore be
adapted to the local geometry.

Altogether, we distinguish $13$ geometric regimes:
\[
\begin{gathered}
\text{liquid},\qquad
\text{regular frozen},\qquad
\text{arctic},\qquad
\text{cusp},\\
\text{vertical tangent},\qquad
\text{vertical cusp},\qquad
\text{vertical frozen},\\
\text{unit-slope tangent},\qquad
\text{unit-slope cusp},\qquad
\text{unit-slope frozen},\\
\text{horizontal tangent},\qquad
\text{horizontal cusp},\qquad
\text{horizontal frozen}.
\end{gathered}
\]

According to the geometry of the relevant vertices, we construct a
phase-adapted open cover of the polygonal domain:
\begin{align}\label{e:covering}
\fP=\bigcup_{\al}\fN_\al.
\end{align}
For example, a liquid neighborhood is a small neighborhood of a point in
the liquid region, an arctic neighborhood is a small neighborhood of a
regular point on the arctic boundary, and a cusp neighborhood is a small
neighborhood of a cusp point.

Fix a lattice point $(y,t)\in\fP$. For each neighborhood $\fN_\al$ in the
covering \eqref{e:covering}, we choose suitable local charts on the Riemann
surface and construct a local analogue of \eqref{e:bulk_illustrate} of the
form
\begin{align}\label{e:localA}
A_\al((x,s),(y,t))
=
J^{(1)}((x,s),(y,t))
+
J^{(2)}((x,s),(y,t)),
\qquad
(x,s)\in\fN_\al,
\end{align}
where $J^{(1)}$ is a single-contour integral and $J^{(2)}$ is a
double-contour integral. The contours are localized near the relevant
critical points and follow the corresponding local steepest-descent and
steepest-ascent paths. Near the polygonal boundary, we modify the
integrands so as to enforce the zero boundary condition
\eqref{e:zero_bb}.

Thus, in every geometric regime, the local ansatz
$A_\al((x,s),(y,t))$ has the same two basic properties: it satisfies the
appropriate zero boundary condition and, as a function of $(x,s)$, obeys
the same inhomogeneous recursion \eqref{e:discrete_kernel2} as the inverse
Kasteleyn matrix. At the same time, its contours and local coordinates are
adapted to the critical-point geometry of the corresponding regime.

The construction of the local charts and the contours for the
double-contour integrals is given in \Cref{s:critical_point}, while
the corresponding single-contour integrals are constructed in
\Cref{s:single_integral}. The phase-adapted covering of the polygonal
domain is introduced in \Cref{s:assign_cover}, and the general kernel
ansatz is presented in \Cref{s:general_kernel_ansatz}. Refined
properties of the critical points and estimates for the contour
integrands are collected in
\Cref{s:revisit_critical,s:Integrand_est}.

\emph{Standard forms and  compatibility.}
The local contour formulas \eqref{e:localA} provide local parametrices for
the inverse Kasteleyn matrix. Their construction alone, however, does not
yet yield a global approximation, since the neighborhoods in the
phase-adapted cover overlap and different geometric descriptions may apply
to the same pair of vertices.

For example, a point $(x,s)$ may belong to the intersection of a liquid
neighborhood $\fN_\al$ and an arctic neighborhood $\fN_\beta$:
$
(x,s)\in\fN_\al\cap\fN_\beta.
$
On this overlap, one may use either the liquid-region ansatz
$A_\al((x,s),(y,t))$ or the arctic-region ansatz
$A_\beta((x,s),(y,t))$. Although the corresponding contour formulas may
look different, they approximate the same entry of the inverse Kasteleyn
matrix and must therefore agree to the required accuracy.

To establish this compatibility, we apply steepest-descent analysis and
rewrite the local parametrices in suitable standard forms. These standard
forms separate the universal contributions arising from the relevant
critical points from the geometry-dependent prefactors and the error
terms. Once the parametrices have been reduced to standard form, their
comparison on overlapping neighborhoods becomes transparent: their
leading contributions agree, while their difference is absorbed into the
required error. In particular, for
$(x,s)\in\fN_\al\cap\fN_\beta$, one has
\begin{align}\label{e:compatibility} 
\left|
A_\al((x,s),(y,t))
-
A_\beta((x,s),(y,t))
\right|
\leq \text{the required error}.
\end{align}

Estimates for the single-contour and double-contour integrals appearing in these local parametrices are provided in \Cref{s:cintegral}.
The standard forms in the liquid and non-liquid regimes are derived in \Cref{s:liquid_standard_form,s:non-liquid_standard_form}.
Their properties and the compatibility of the local approximations \eqref{e:compatibility} are proved in Section \Cref{s:compatibility_proof}. 

\emph{Uniform kernel asymptotics and the Gaussian free field.}
The compatibility property \eqref{e:compatibility} allows us to glue the
local parametrices into a single global approximation to the inverse
Kasteleyn matrix. Let $\{\rho_\al\}_\al$ be a partition of unity subordinate
to the covering \eqref{e:covering}, and define
\begin{align}\label{e:global_parametrix}
A((x,s),(y,t))
:=
\sum_\al
\rho_\al(x,s)A_\al((x,s),(y,t)).
\end{align}
Since each $A_\al$ satisfies the appropriate zero boundary condition, the
global parametrix $A$ satisfies the same boundary condition.

Applying $K$ in the first variable gives
\begin{align}\label{e:global_error}
KA=I+\cE,
\end{align}
where $\cE$ is an error matrix. Indeed, the terms in which $K$ acts on the
local parametrices reproduce the inhomogeneous recursion
\eqref{e:discrete_kernel2}. The remaining terms arise from the discrete
commutators between $K$ and the cutoff functions $\rho_\al$ and can be
rewritten, on overlapping neighborhoods, in terms of differences of the
form $A_\al-A_\beta$. They are therefore small by the compatibility
estimate \eqref{e:compatibility}. Together with the local error estimates,
this yields the desired bound on $\cE$. Moreover, for all sufficiently
large $n$, the Neumann series converges, and
\begin{align}
(I+\cE)^{-1}
=
I-\cE+\cE^2-\cE^3+\cdots=I+\OO(\cE).
\end{align}
It follows from \eqref{e:global_error} that
\begin{align}\label{e:global_approximation}
K^{-1}
&=
A(I+\cE)^{-1}
=
A-A\cE+A\cE^2-\cdots
=
A+\OO(A\cE).
\end{align}

The preceding construction therefore gives a uniform asymptotic
approximation to the inverse Kasteleyn matrix throughout the polygonal
domain. In particular, the resulting estimates are not restricted to
compact subsets of the liquid region: they also apply in the frozen
region and in the transition neighborhoods associated with the arctic
boundary.

This global control is the main input for the analysis of height
fluctuations. By the determinantal structure of the dimer model, the
joint cumulants of centered height functions can be expressed in
terms of sums of cyclic products of the Kasteleyn matrix and its
inverse \cite{petrov2015asymptotics,kenyon2008height}. Substituting the global kernel approximation into these
expressions, we show that the second cumulant converges to the Green
function of the liquid region, while all cumulants of order at least
three vanish. Consequently, the limiting joint moments are those of
the Gaussian free field.

The construction of the global approximation
$A((x,s),(y,t))$ to the inverse Kasteleyn matrix and the proof of
\eqref{e:global_approximation} are given in
\Cref{s:global_kernel_approximate}. The resulting estimates for the inverse
Kasteleyn matrix are then derived from \eqref{e:global_approximation} in
\Cref{s:final_kernel}.
The proof of convergence of the height function to the Gaussian free
field is given in \Cref{s:gff_moment_convergence}. Finally, additional
geometric properties of the frozen region are collected in
\Cref{s:frozen_structure}, and the detailed steepest-descent analysis of the
contour integrals is carried out in \Cref{s:path_analysis}.

\subsection{Related work}
The Gaussian free field describes the centered fluctuations of the
height function on the macroscopic scale, after the height function is
interpreted as a random distribution. This field-level limit should be
distinguished from local tiling statistics, which probe individual
lozenges either on the lattice scale in the liquid region or in
anisotropically shrinking neighborhoods of the arctic boundary.

At a fixed point in the interior of the liquid region, the microscopic
tiling converges to the infinite-volume, translation-invariant extremal
Gibbs measure whose slope is determined by the gradient of the limit
shape at that point. Under the standard particle encoding, this Gibbs
measure is determinantal, with correlation kernel given by the
incomplete beta kernel, a two-dimensional extension of the discrete
sine kernel
\cite{okounkov2003correlation,sheffield2005random,
kenyon2006dimers}. This bulk universality was proved for uniformly
random lozenge tilings of general simply connected domains in
\cite{aggarwal2023universality}; see also
\cite{johansson2002non,gorin2008nonintersecting,borodin2010q,petrov2014asymptotics,gorin2017bulk,gorin2019universality,laslier2019local}
for earlier results in more restrictive settings.

At a regular point of the arctic boundary that is neither a cusp nor a
boundary tangency, the natural transverse and tangential scales are of
orders \(n^{1/3}\) and \(n^{2/3}\), respectively, in lattice
coordinates. For simply connected polygonal domains satisfying a
technical genericity assumption on the limit shape, the rescaled
nonintersecting path ensemble converges to the Airy line ensemble
\cite{huang2024edge,aggarwal2025edge}. Earlier Airy limits for special
integrable domains were obtained in
\cite{johansson2002non,ferrari2003step,okounkov2007random,baik2007discrete,petrov2014asymptotics,duse2018universal}.

At a tangency between the arctic boundary and a straight side of the
domain, \cite{aggarwal2022gaussian} proved that the local statistics are
governed by the GUE--corners process. At an ordinary cusp, the natural
anisotropic scales are of orders \(n^{1/2}\) and \(n^{1/4}\), and, under
the genericity assumptions of \cite{huang2024pearcey}, the local point
process converges to the Pearcey process. These three limits correspond
to the regular, boundary-tangency, and ordinary-cusp geometries expected
for generic polygonal arctic boundaries. For the distinct Airy-cusp
geometry, the corresponding universality was established in
\cite{kumar2026airycusp}.

Nongeneric singularities can produce further universality classes. At
cusp--turning points, Cusp--Airy statistics have been proved in
certain integrable lozenge-tiling models
\cite{duse2016cusp}. Tacnode geometries give rise to
several distinct tacnode-type processes, including continuous,
GUE-minor, discrete, and hard-edge tacnode processes
\cite{adler2015coupled,
adler2015tacnode,
adler2018lozenge,
adler2018tilings,
ferrari2021fluctuations,
adler2022singular}.

Beyond uniform tilings, for periodically weighted bipartite dimer
models, the spectral curve may produce frozen, liquid, and gaseous
phases \cite{kenyon2006dimers}. The liquid region can be multiply
connected, with gaseous islands enclosed by liquid components. This
geometry has been analyzed in detail for periodic Aztec diamonds
\cite{chhita2016domino,berggren2021domino,
berggren2025geometry}.
Recent work shows that the
centered height function is described by a Gaussian free field on the
multiply connected liquid region together with an independent random
harmonic component whose boundary values form a discrete Gaussian
vector \cite{berggren2025gaussian}.
Thus, periodic weights can produce both nontrivial liquid topology and
additional finite-dimensional facet-height fluctuations beyond the
Dirichlet Gaussian free field.

The six-vertex model provides an important non-dimer counterpart.
At the free-fermion point, it can be related directly to a dimer
model. Perturbative renormalization-group methods show that sufficiently
small non-free-fermionic perturbations still produce logarithmically
correlated Gaussian height fluctuations, although the fluctuation
amplitude is generally nonuniversal
\cite{giuliani2017height}. More recently, convergence of the
full-plane six-vertex height function to a suitably normalized
Gaussian free field was proved throughout the parameter range
$
-1\leq \Delta\leq-1/2,
$
including anisotropic weights after an appropriate embedding
\cite{duminil2026gaussian}.

\subsection{Notation}

We denote by \(\bZ\) the set of integers and by
\[
\bZ':=\bZ+\frac12
=
\left\{k+\frac12:k\in\bZ\right\}
\]
the set of half-integers.

For a quantity \(X\) and a nonnegative quantity \(Y\), we write
$
X=\OO(Y)$ or $X\lesssim Y
$
if
$
|X|\leq \fC Y
$
for some constant \(\fC>0\) independent of the relevant asymptotic
parameters. The implicit constant may change from line to line. For
nonnegative quantities \(X\) and \(Y\), we write \(X\gtrsim Y\) if
\(Y\lesssim X\), and \(X\asymp Y\) if both \(X\lesssim Y\) and
\(X\gtrsim Y\). We write
$
X=\oo(Y)
$ or $
X\ll Y
$
if \(|X|/Y\to0\), and \(X\gg Y\) if \(Y/X\to0\).

For any set \(\mathfrak R\subseteq\bR^2\), we denote its boundary and
closure by \(\partial\mathfrak R\) and \(\overline{\mathfrak R}\),
respectively. For \(c\in\bR\) and \(u\in\bR^2\), we define
$
c\mathfrak R
:=
\{cr:r\in\mathfrak R\},
$ and $
\mathfrak R+u
:=
\{r+u:r\in\mathfrak R\}.
$
For \(x\in\bR^2\) and \(r>0\), let
\[
B_r(x)
:=
\{y\in\bR^2:\|x-y\|_2<r\}.
\]
For a set \(E\subseteq\bR^2\), define its \(r\)-neighborhood by
\[
B_r(E)
:=
\bigcup_{x\in E}B_r(x)
=
\{y\in\bR^2:\dist(y,E)<r\}.
\]

For \(z\in\bC\setminus\{0\}\), we denote its principal argument by
$
\arg z\in(-\pi,\pi].
$
Throughout this article, we use the principal branch of the logarithm,
\[
\ln:\bC\setminus(-\infty,0]\longrightarrow\bC,
\qquad
\ln z=\ln|z|+\ri\arg z,
\]
and the corresponding principal branch of the square root,
$
\sqrt z:=
\exp(\ln z)/2).
$

\subsection*{Acknowledgements}
The research of J.H. was supported in part by NSF Grant DMS-2246664, DMS-2331096  and a Sloan Research Fellowship.
The author would like to thank Amol Aggarwal, Vadim Gorin, Matthew Nicoletti, and Leonid Petrov for enlightening discussions, and Vadim Gorin for helpful comments on the manuscript.
The mathematical content of this paper was developed by the author. Large language model tools were used primarily for language editing and literature searches.

\chapter{Riemann Surface and Tiling Action}
\section{Limit Shape and Complex Slope}
\label{s:limit_shape}
	
	In this section we collect some results on the arctic boundary of lozenge tilings of polygonal domains, and decomposition of the frozen region.

\subsection{Frozen region structure}
Our first result records several properties of the arctic boundary
$\fA=\partial\fL$. These properties are essentially due to
\cite{astala2026dimer,de2010minimizers}; for completeness, we give the proof in
\Cref{s:limits}.

  \begin{proposition}[{\cite{astala2026dimer,de2010minimizers}}]\label{pa1}
With the notation of \Cref{p}, suppose that \Cref{a:asump} holds. Then the
following statements hold.

\begin{enumerate}
    \item The complex slope
  $
        f:\mathfrak L \to \mathbb C_-$
    extends continuously to
    $
        \overline{\mathfrak L}\to \mathbb C_-\cup \mathbb R\cup\{\infty\}.
   $
    Moreover, the slope of the arctic boundary at $(x,s)\in\mathfrak A$ is
    given by
    \begin{equation}\label{e:slope}
        \frac{1}{\chi(x,s)}
        :=
        \frac{f(x,s)+1}{f(x,s)} .
    \end{equation}

    \item 
    Define
    \begin{equation}\label{e:defg}
        g(x,s)
        =
        \frac{
            1-\ri (f(x,s)+1)^{-1}f(x,s)
        }{
            1+\ri (f(x,s)+1)^{-1}f(x,s)
        },
        \qquad (x,s)\in \fL .
    \end{equation}
    Since $f(x,s)\in \mathbb C_-$ for $(x,s)\in\fL$, one has
    $g(x,s)\in\mathbb D$. Let $d$ denote the degree of the map
    $g:\fL\to\mathbb D$. Then, as one traverses the arctic boundary
    counterclockwise, the boundary value of the complex slope
  $
        f:\fA\to \mathbb R\cup\{\infty\}
   $
    winds around $\mathbb R\cup\{\infty\}$ exactly $d$ times in the negative
    orientation.

    \item \label{i:no_tacnode} The arctic curve $\mathfrak A$ has $d-2$ inward-pointing cusps and
    no tacnodes. Moreover, $\partial\fL$ is locally smooth and strictly convex
    except at these $d-2$ cusps. More precisely, for every
    $\zeta\in\partial\fL$ which is not a cusp, there exists
    $\varepsilon>0$ such that
$
        B_\varepsilon(\zeta)\cap\fL
$
    is strictly convex.
\end{enumerate}
\end{proposition}

For a point $(x,s)\in \fA$, we call it a tangency location of $\fA$, if the tangent line to $\fA$ at that point has slope in $\{0, 1, \infty\}$. By \eqref{e:slope}, at such tangency locations, the value of $f(x,s)$ lies in $\{-1,\infty, 0\}$.

The following theorem describes the structure of the frozen region
$\overline{\fP}\setminus\fL$. The result is essentially due to
\cite{kenyon2007limit}, whose proof uses algebraic methods. For completeness,
we provide an alternative analytic proof in \Cref{s:frozen}, based on results
from \cite{de2010minimizers,astala2026dimer}.

\begin{figure}
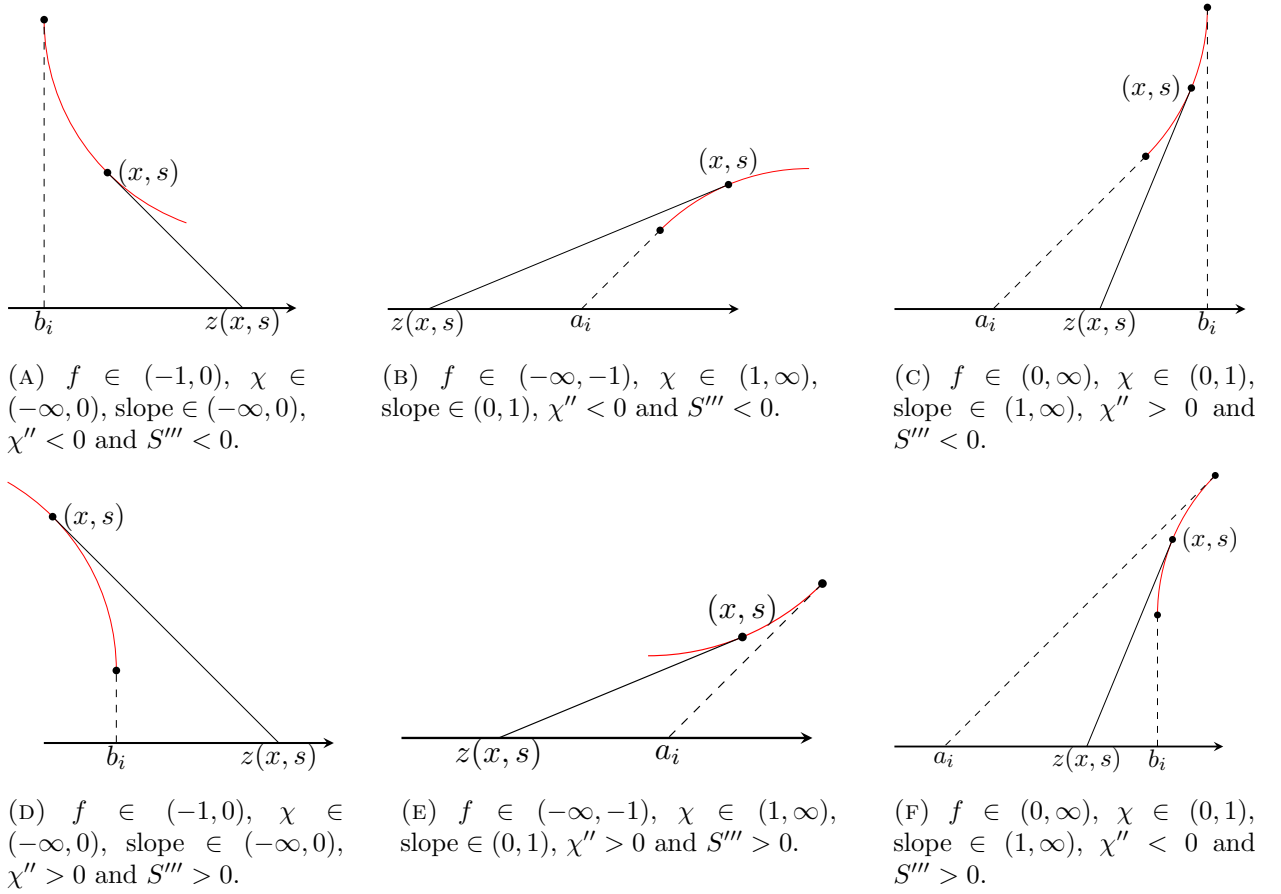


%================= Row 1 =================
\begin{subfigure}[t]{0.24\textwidth}
  \centering
  \resizebox{\linewidth}{!}{%

  }
   \caption{$f\in(-1,0)$, $\chi\in (-\infty,0)$, $\text{slope}\in (-\infty,0)$,  $\chi''<0$ and $S'''<0$.}
 \end{subfigure}
 \hfill
\begin{subfigure}[t]{0.35\textwidth}
 \centering
  \resizebox{\linewidth}{!}{%
    %
  }
  \caption{$f\in(-\infty,-1)$, $\chi\in (1,\infty)$, $\text{slope}\in (0,1)$, $\chi''<0$ and $S'''<0$.}
\end{subfigure}
\hfill
\begin{subfigure}[t]{0.29\textwidth}
  \centering
  \resizebox{\linewidth}{!}{%
    %
  }
   \caption{$f\in(0,\infty)$, $\chi\in (0,1)$, $\text{slope}\in (1,\infty)$, $\chi''>0$ and $S'''<0$.}
\end{subfigure}

\medskip

%================= Row 2 =================
\begin{subfigure}[t]{0.27\textwidth}
  \centering
  \resizebox{\linewidth}{!}{%
    %
  }
   \caption{$f\in(-1,0)$, $\chi\in (-\infty,0)$, $\text{slope}\in (-\infty,0)$, $\chi''>0$  and $S'''>0$.}
\end{subfigure}
\hfill
\begin{subfigure}[t]{0.35\textwidth}
  \centering
  \resizebox{\linewidth}{!}{%
    %
  }
   \caption{$f\in(-\infty,-1)$, $\chi\in (1,\infty)$, $\text{slope}\in (0,1)$, $\chi''>0$ and $S'''>0$.}
\end{subfigure}
\hfill
\begin{subfigure}[t]{0.29\textwidth}
  \centering
  \resizebox{\linewidth}{!}{%
    %
  }
   \caption{$f\in(0,\infty)$, $\chi\in (0,1)$, $\text{slope}\in (1,\infty)$, $\chi''<0$ and $S'''>0$.}
\end{subfigure}

\caption{For $(x,s)\in \fA$, $z(x,s)=x-s \chi(x,s)$ is the intersection of the tangent line to the arctic curve at $(x,s)$ with the $x$ axis.}
\label{f:arctic_boundary}
\end{figure}

\begin{theorem}[Frozen region structure]\label{t:frozen_structure}
With the notation of \Cref{p}, suppose that \Cref{a:asump} holds. The arctic curve
\(\fA\) is tangent to the sides of the polygon \(\fP\), or to their linear extensions.
The corresponding tangency points occur in the natural cyclic order of the sides of
\(\fP\).

For each vertex \(\zeta\) of \(\fP\), let \(\zeta_1,\zeta_2\in\fA\) be the two
tangency points corresponding to the two sides of \(\fP\) incident to \(\zeta\),
or to their linear extensions. Then the three boundary pieces
\[
[\zeta,\zeta_1],\qquad
[\zeta,\zeta_2],\qquad
\text{and the arc of }\fA\text{ between }\zeta_1\text{ and }\zeta_2
\]
bound a closed region. We call this region a \emph{curvilinear triangle}.

As \(\zeta\) ranges over the vertices of \(\fP\), these curvilinear triangles have
pairwise disjoint interiors and cover the entire frozen region. On each curvilinear triangle \(\fT\), the height function has constant gradient
\[
\nabla H^* \in \{(1,0),(1,-1),(0,0)\}.
\]
More precisely, the possible cases are as follows (see \Cref{f:arctic_boundary}):
\[
\begin{array}{c|c|c|c}
\nabla H^*\text{ on }\fT
&
\text{slope of the two linear sides}
&
\text{slope of }\fA\text{ along the arc}
&
f\text{-range}
\\
\hline
(1,0)
&
\infty,\ 0
&
(-\infty,0)
&
(-1,0)
\\
(1,-1)
&
0,\ 1
&
(0,1)
&
(-\infty,-1)
\\
(0,0)
&
1,\ \infty
&
(1,\infty)
&
(0,\infty)
\end{array}
\]
Here the third column refers to the slopes of the tangent lines to the arctic arc from \(\zeta_1\) to \(\zeta_2\).

Finally, let \(\angle \zeta_1\zeta\zeta_2\) denote the sector containing
\(B_\delta(\zeta)\cap\fP\) for all sufficiently small \(\delta>0\). This sector has
opening angle greater than \(180^\circ\) if and only if \(\zeta\) is a concave
corner of \(\fP\). The number of cusp points on the arctic arc from \(\zeta_1\) to
\(\zeta_2\) is
\begin{align}
\label{e:cusp_num}
\#\{\text{cusps on the arc from }\zeta_1\text{ to }\zeta_2\}
&=
\bm{1}\bigl(\text{\(\zeta\) is a concave corner of \(\fP\)}\bigr)
\nonumber\\
&
+\bm{1}\bigl(\text{\(\fA\) is tangent at \(\zeta_1\) on the exterior side of }
        \angle \zeta_1\zeta\zeta_2\bigr)
\nonumber\\
&
+\bm{1}\bigl(\text{\(\fA\) is tangent at \(\zeta_2\) on the exterior side of }
        \angle \zeta_1\zeta\zeta_2\bigr).
\end{align}
\end{theorem}

\begin{figure}
  \begin{subfigure}{0.25\textwidth}
    \centering
      \begin{tikzpicture}[scale=0.6]
        \draw[red] (0,4) arc[start angle=180, end angle=270, radius=4];
        \draw[thick] (-2,0)--(4,0);
        \draw[-stealth, thick] (0,-2)--(0,4);
      \end{tikzpicture}
  \end{subfigure}
   \begin{subfigure}{0.35\textwidth}
    \centering
      \begin{tikzpicture}[scale=0.6]
        \draw[red] (4,0) arc[start angle=90, end angle=135, radius={4*tan(67.5)}];
        \draw[thick] (-4,-4)--(2,2);
        \draw[-stealth, thick] (-2,0)--(4,0);
      \end{tikzpicture}
  \end{subfigure}
    \begin{subfigure}{0.25\textwidth}
    \centering
      \begin{tikzpicture}[scale=0.6]
        \draw[red] (0,4) arc[start angle=0, end angle=-45, radius={4*tan(67.5)}];
        \draw[-stealth,thick] (2,2)--(-4,-4);
        \draw[ thick] (0,-2)--(0,4);
      \end{tikzpicture}
  \end{subfigure}

\caption{Symmetry among the three types of curvilinear triangles. The left,
middle, and right panels correspond to
\(\nabla H^*=(1,0)\), \(\nabla H^*=(1,-1)\), and
\(\nabla H^*=(0,0)\), respectively.}
  \label{f:symmetry}
\end{figure}
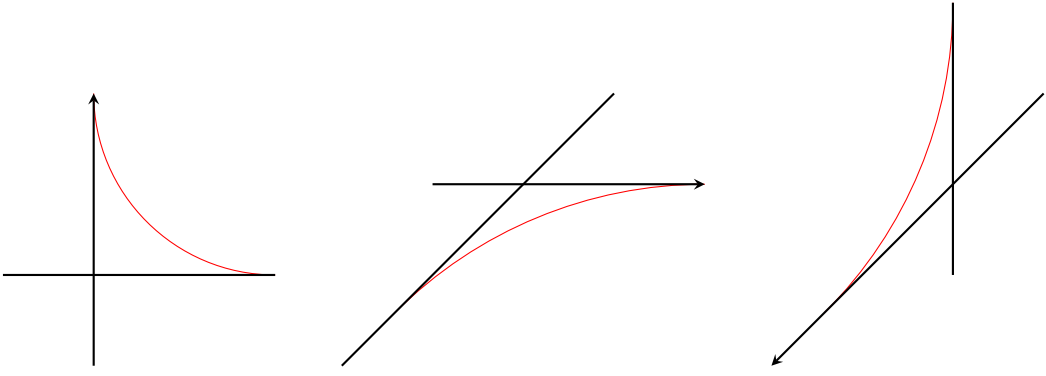

\begin{figure}
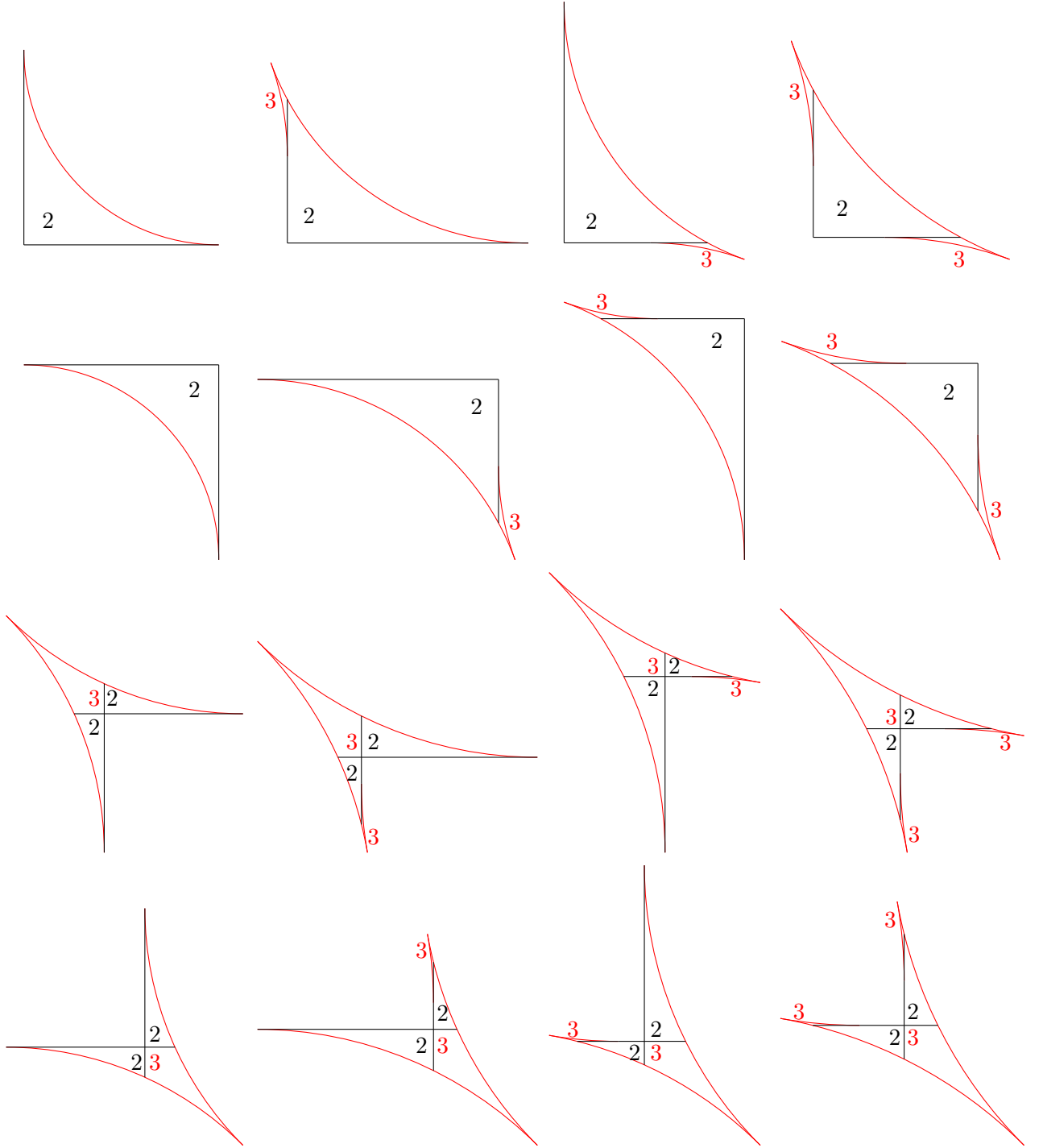

  \centering

  %================= Row 1 =================
  \begin{subfigure}{0.23\textwidth}
    \centering

  \end{subfigure}

  \caption{Curved triangle frozen region with arctic curves shown in red and extended sides in black.}
  \label{f:curvilinear_triangle}
\end{figure}

\begin{remark}\label{r:symmetry}
All possible configurations of curvilinear triangles with
\(\nabla H^*=(1,0)\) are shown in \Cref{f:curvilinear_triangle}. The three
gradient types are related by affine symmetries together with the corresponding
relabeling of the lozenge orientations.

In \Cref{f:symmetry}, the left, middle, and right panels correspond to
\[
\nabla H^*=(1,0),\qquad
\nabla H^*=(1,-1),\qquad
\nabla H^*=(0,0),
\]
respectively. Starting from the left panel, rotate \(90^\circ\) clockwise and
then apply the shear
\[
(u,v)\longmapsto(u+v,v).
\]
This sends the side slopes \((\infty,0)\) to \((0,1)\) and the tangent-slope
range \((-\infty,0)\) to \((0,1)\), giving the middle configuration. Rotating
the middle configuration \(135^\circ\) clockwise sends the side slopes
\((0,1)\) to \((1,\infty)\) and the tangent-slope range \((0,1)\) to
\((1,\infty)\), giving the right configuration. The arrows track the corresponding oriented side under these transformations.

\end{remark}

\begin{cor}\label{c:intersection}
With the notation of \Cref{p}, suppose that \Cref{a:asump} holds.
Let \(\fT\) be a  curvilinear triangle. If \(\ell\) is tangent to the
arctic boundary at a point on the portion of the arctic boundary contained in
\(\fT\), then \(\fT\cap \ell\) is a connected line segment.
\end{cor}

\begin{proof}
The claim follows by checking each case in the classification of curvilinear
triangles given in \Cref{f:curvilinear_triangle}.
\end{proof}

\begin{definition}[Extended sides]\label{d:extended_sides}
An \emph{extended side} is a maximal straight segment obtained by extending
a side of the polygon \(\fP\) along its supporting line while remaining in
the frozen region; see \Cref{f:curvilinear_triangle}. In particular, near
each vertex of \(\fP\), the two straight segments joining the vertex to the
corresponding tangency points on the arctic boundary are contained in
extended sides. When the context is clear, we also refer to the supporting
line of an extended side as an extended side.
\end{definition}

\begin{figure}
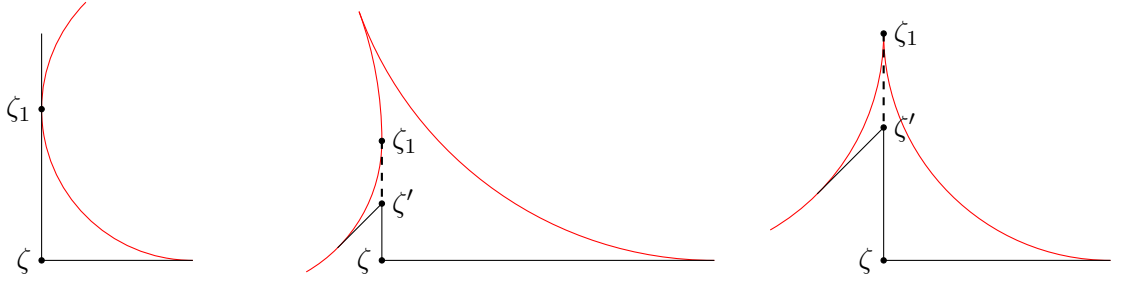

  \centering
  \begin{subfigure}{0.25\textwidth}
    \centering

  \end{subfigure}
  \caption{Two adjacent curvelinear triangles share a tangency location $\zeta_1$, and (possibly) a segment $[\zeta_1, \zeta']$ (the dashed segment).}
  \label{f:adjacent_curvilinear_triangle}
\end{figure}

\begin{remark}
Two adjacent curvilinear triangles share a tangency location $\zeta_1$ and, possibly, a segment $[\zeta_1,\zeta']$ that is a linear extension of a side of the polygon $\fP$.
In the left panel of \Cref{f:adjacent_curvilinear_triangle}, they share only the tangency location $\zeta_1$.
In the middle panel, they share a segment $[\zeta_1,\zeta']$, where $\zeta_1$ is a tangency location and $\zeta'$ is a vertex of $\fP$.
In the right panel, they share a segment $[\zeta_1,\zeta']$, where $\zeta_1$ is a cusp-turning location and $\zeta'$ is a vertex of $\fP$.
\end{remark}

\subsection{Complex slope and Riemann surface}
\label{s:surface}

	We consider the following map from the closure of the liquid region $\overline{\fL}=\fL\cup\fA$ to $\mathbb{CP}\times \mathbb{CP}$
\begin{align}\label{e:emb}
\Psi: (x,s)\in \overline\fL\mapsto (f(x,s), z(x,s)):=\left(f(x,s), x-s\chi(x,s)\right)\in \mathbb{CP}\times \bC\bP,
\end{align}
where 
\begin{align}\label{e:chi}
\chi(x,s)=\frac{f(x,s)}{(f(x,s)+1)}
\end{align}
 is as defined in \eqref{e:slope}. For $(x,s)\in \fA$, the map $(x,s)\mapsto z(x,s)$ has a simple geometric meaning. Recall from the first statement in \Cref{pa1} that the tangent vector to $(x,s)\in  \mathfrak{A}$ has slope $1/\chi(x,s)$. Thus $z(x,s)=x-s \chi(x,s)$ is the intersection of the tangent line to the arctic curve at $(x,s)$ with the $x$ axis; See Figure \ref{f:arctic_boundary}. 

By the third statement of \Cref{pa1}, the arctic boundary $\fA$ is locally strictly convex. As \( (x, s) \) traverses the arctic curve \( \mathfrak{A} \) in counterclockwise order, the imagine $z(x,s)$ increases monotonically until it diverges to 
$+\infty$, then wraps around to $-\infty$ and continues increasing. Over the entire curve \( \mathfrak{A} \), \( z(x, s) \) sweeps through \( \mathbb{R} \cup \{\infty\} \) exactly \( d \) times.

With the notations in \eqref{e:emb}, we can reformulate \Cref{fequation} as the following that there exists an analytic function $Q$ of two variables such that $Q(f(x,s), z(x,s))=0$. For polygonal domains as in \Cref{p}, the following \Cref{p:surface} states that the analytic function $Q$ is in fact a rational function. Then $Q(f,z)=0$ defines a Riemann surface, which can be identified as gluing the image of the map \eqref{e:emb} restricted to the liquid region $\fL$, and its complex conjugate along the image of the arctic boundary.  This result essentially follows from
\cite{kenyon2007limit, astala2026dimer}; for completeness, we give the proof in
\Cref{s:Rsurface}.

\begin{proposition}\label{p:surface}
Assume \Cref{a:asump} holds. 
%Consider the following map from the closure of \( \fL \), i.e., \( \overline{\fL} = \fL \cup \fA \), into \( \mathbb{CP}^2 \):
%\begin{align}\label{e:emb}
%(x, t) \in \overline{\fL} \mapsto (f(x,t), z(x,t)) := \left(f_t(x),\, x - t f_t(x)\right) \in \mathbb{CP}^2.
%\end{align}
Then the following statements hold for the map \eqref{e:emb}:
\begin{enumerate}
    \item Recall $\chi(x,s)$ from \eqref{e:chi}.
    Locally around a point \( (x_0, s_0) \in \fL \), if \( \partial_x \chi(x_0, s_0) \neq 0 \), then there exists an analytic function \( Q_0 \) defined near \( f(x_0, s_0) \) such that
    \begin{align}\label{e:introduceQ0}
        Q_0(f(x,s)) = z(x,s).
    \end{align}
    Locally around a point \( (x_0, s_0) \in \fL \), if \( \partial_x \chi(x_0, s_0) \neq 1/s_0 \), then there exists an analytic function \( f \) defined near \(z_0= z(x_0, s_0) \)  such that
    \begin{align}\label{e:introducef0}
        f(x,s) = f(z(x,s)).
    \end{align}
    \item 
    
    Locally around a point \( (x_0, s_0) \in \fA \) on the arctic boundary such that $f(x_0, s_0)\neq -1$, there exists a real analytic function \( \chi \) (i.e. $\chi(\overline{z})=\overline{\chi(z)}$) defined near \(z_0= z(x_0, s_0) \) such that
 \begin{align}\label{e:introducef0_boundary}
        \chi(x,s) = \chi(z(x,s)) \text{ for } (x,s)\in B_\delta (x_0,s_0)\cap \overline{\fL},
    \end{align}
Moreover we have the following cases
\begin{enumerate}
   \item  If \( (x_0, s_0) \in \fA \) is not a cusp location, then we have $\chi''(z_0)\neq 0$, and the arctic boundary can be parametrized as \( (\gamma(s), s) \) for \( s \in [s_0 - \varepsilon, s_0 + \varepsilon] \) with \( \varepsilon > 0 \) sufficiently small, and the following identities hold:
    \begin{align}\label{e:arctic_chider}
        -\frac{1}{s} = \chi'(z(\gamma(s), s)), \quad 
        \gamma'(s) = \chi(z(\gamma(s), s)), \quad 
        \gamma''(s) = -\frac{1}{s^3 \chi''(z(\gamma(s), s))}.
    \end{align}
The sign of \(\chi''\bigl(z(\gamma(s),s)\bigr)\) is determined by the convexity of
\(\gamma(s)\), as classified in \Cref{f:arctic_boundary}.

    \item 
    If  \( (x_0, s_0) \in \fA \) is a cusp location, then we have    \begin{align}\label{e:introducef0_cusp}
           \chi''(z_0) = 0, \quad 
        \chi'''(z_0) \neq 0.
    \end{align}
    Furthermore, the cusp points downward if \( \chi'''(z_0) > 0 \), and upward if \( \chi'''(z_0) < 0 \); see \Cref{f:cusp}.
\end{enumerate}

    \item 
    Locally around a \emph{horizontal tangent location} \( (x_0, s_0) \in \fA \) on the arctic boundary such that $f(x_0, s_0)=-1$, there exists an analytic function \(\chi \) defined near $\infty$ such that as \( z \to \infty \),
    \begin{align}\label{e:close_to_t}
        \chi(z) = \frac{-z+x_0}{s_0} + \sum_{i\geq 1}\frac{c_i}{(x_0-z)^i}= \wt \chi({\wt z})=\frac{1}{s_0{\wt z}}+\sum_{i\geq 1}c_i {\wt z}^i,\quad {\wt z}=\frac{1}{x_0-z}.
    \end{align}
    
    If $c_1\neq 0$, then \( (x_0, s_0)  \) is not a cusp location. Moreover, near $(x_0, s_0)$ the arctic boundary is locally convex in the $s$ direction when $c_1>0$, and locally concave when $c_1<0$.     
    If $c_1=0$, then \( (x_0, s_0)  \) is a cusp location, and $c_2\neq 0$. Furthermore, the cusp points leftward if \( c_2 > 0 \), and rightward if \( c_2< 0 \).   
     \item 
    There exists a nonzero rational function \( Q \) such that for all \( (x,s) \in \overline{\fL} \),
    \begin{align}\label{e:qfh}
        Q\big(f(x,s),\, z(x,s)\big) = 0.
    \end{align}
    \end{enumerate}
\end{proposition}

Let $\cC\subset\bC\bP\times\bC\bP$ be the projective curve defined by $Q(f,z)=0$,
where $Q$ has real coefficients; thus $\cC$ is real curve. The set $\cC$ is an immersed Riemann surface that can be realized by
gluing the image of \eqref{e:emb} on the liquid region $\fL$ to its complex conjugate
along the arctic boundary:
\[
\cC
=\bigl\{(f(x,s),\,z(x,s)):\,(x,s)\in\overline{\fL}\bigr\}
 \;\cup\;
 \bigl\{(\overline{f(x,s)},\,\overline{z(x,s)}):\,(x,s)\in\overline{\fL}\bigr\}.
\]
The real locus $\cC(\bR):=\cC\cap(\bR\bP\times\bR\bP)$ is homeomorphic to a circle and
separates $\cC$ into the two (conjugate) components
\[
\cC\cap(\bC_+\times\bC_-)
\qquad\text{and}\qquad
\cC\cap(\bC_-\times\bC_+).
\]

The curve $\cC$ is rational (its normalization has genus $0$), so topologically it is a Riemann sphere after resolving singularities. There is a single singularity at $(-1, \infty)$. The projection $(f,z)\in \cC\mapsto z\in \bC\bP$ has degree $d$ and the entire fiber over $z=\infty$ is supported at the single ambient point $(-1, \infty)$;  locally each of the $d$ branches admits a smooth chart. Everywhere else $\cC$ is smooth. In what follows we identify \(\cC\) with its normalization and denote the fiber over \(z=\infty\) by
\begin{align}\label{e:inftypoint}
(-1,\infty_1),\;(-1,\infty_2),\;\dots,\;(-1,\infty_d).
\end{align}

From the first statement in \Cref{p:surface}, the ramification points of the projection $(f,z)\in \cC\mapsto z\in \bC\bP$
correspond to points $\Psi(x_0,s_0), \overline{\Psi(x_0,s_0)}$ (recall
\eqref{e:emb}) with $(x_0,s_0)\in \fL$ and
$\del_x \chi(x_0,s_0)=1/s_0$. From \Cref{p:surface}, locally around
non-ramification points $(f,z)\in \cC$, we can parametrize $\cC$ as a function
$f(z)$ of $z$. When the context is clear, we will represent these points
\((f,z)\in\cC\) by their \(z\)-coordinate. With a slight abuse of notation,
define the meromorphic function $\chi:\cC\to\mathbb{CP}^1$ by
\begin{equation}\label{e:def-chi-on-C}
    \chi(f,z):=\frac{f}{f+1}.
\end{equation}
On a local branch where $f=f(z)$, we write simply
\begin{align}\label{e:def-chi-on-C2}
    \chi(z)=\frac{f(z)}{f(z)+1}.
\end{align}
This notation is consistent with \eqref{e:introducef0_boundary} and
\eqref{e:close_to_t}. In particular, on the liquid region we have
\begin{align}
\chi(z(x,s))=\chi(x,s), \quad (x,s)\in \fL\cup \fA.
\end{align}

With the above notation, we can rewrite the second coordinate of the map
$\Psi$ from \eqref{e:emb} as
\begin{align}\label{e:bulk_der}
x=z(x,s)+s\chi(z(x,s)).
\end{align}
By taking the derivative with respect to $x$ on both sides of
\eqref{e:bulk_der}, we obtain
\begin{align}
1
=
(\del_x z(x,s))\bigl(1+s\chi'(z(x,s))\bigr)
=
\bigl(1-s\del_x \chi(x,s)\bigr)
\bigl(1+s\chi'(z(x,s))\bigr).
\end{align}
So for $(x,s)\in\fL$, we have
\begin{align}\label{e:derchi_liquid}
\chi'(z(x,s))\neq -1/s.
\end{align}
This holds even when $(x,s)$ approaches a point $(x_0,s_0)$ with
$\del_x \chi(x_0,s_0)=1/s_0$, and in this case $\chi'(z(x,s))$ approaches
$\infty$. Moreover, by \eqref{e:arctic_chider}, the condition
$\chi'(z(x,s))=-1/s$ characterizes the arctic boundary.

We collect some estimates on $\chi(z)$ from \eqref{e:def-chi-on-C2}. The proof follows from a Taylor expansion, and is postponed to \Cref{s:Rsurface}.
\begin{lemma}\label{l:derivechi}
Fix \((x_0,s_0)\in\fA\), and let \((x,s)\in\fL\) be close to \((x_0,s_0)\).
Recall the embedding relation \(x=z(x,s)+s\chi(z(x,s))\) (see \eqref{e:bulk_der}).

Assume \((x_0,s_0)\) is bounded away from horizontal tangent locations, and write
\[
z_0:=z(x_0,s_0),\quad z(x,s)=z_0+a+\ri b, \quad a\in \bR,\quad  b>0.
\]
Then:
\begin{enumerate}
\item If \((x_0,s_0)\) is \emph{not} a cusp location, then $\chi''(z_0)\neq 0$, and 
\begin{equation}\label{e:curve_reg}
\frac1s+\chi'(z(x,s))=\ri\,\chi''(z_0)\,b+\OO\!\bigl(b(|a|+b)\bigr).
\end{equation}
\item If \((x_0,s_0)\) \emph{is} a cusp location, then \(\chi''(z_0)=0\) and \(\chi'''(z_0)\neq 0\), and
\begin{equation}\label{e:curve_cusp}
\frac1s+\chi'(z(x,s))=\chi'''(z_0)\,b\Bigl(\ri a-\frac{b}{3}\Bigr)+\OO\!\bigl(b(|a|+b)^2\bigr).
\end{equation}
\end{enumerate}

Assume \((x_0,s_0)\) is a horizontal tangent location, and set 
\[
\wt z(x,s):=\frac{1}{x_0-z(x,s)}=a+\ri b,\qquad b=\Im \wt w(x,s)>0,\quad a\in \bR,\quad b>0.
\]
Let \(\widetilde\chi\) be the analytic function near \(0\) such that
\(\chi(z(x,s))=\widetilde\chi(\wt z(x,s))\).
Then:
\begin{enumerate}
\item If \((x_0,s_0)\) is \emph{not} a cusp-turning location (equivalently \(\widetilde\chi'(0)\neq 0\)), then
\begin{equation}\label{e:curve_ht_reg}
\frac1s+{\wt z}(x,s)^2\,\widetilde\chi'({\wt z}(x,s))
=\del_{\wt z} ({\wt z}\widetilde\chi({\wt z}))|_{{\wt z}=0} \ri b\,(a+\ri b)+\OO\!\bigl(b(|a|+b)^2\bigr).
\end{equation}
\item If \((x_0,s_0)\) \emph{is} a cusp-turning location (equivalently \(\widetilde\chi'(0)=0\) and \(\widetilde\chi''(0)\neq 0\)), then
\begin{equation}\label{e:curve_ht_cusp}
\frac1s+{\wt z}(x,s)^2\,\widetilde\chi'({\wt z}(x,s))
=\del_{\wt z}^2 ({\wt z}\widetilde\chi({\wt z}))|_{{\wt z}=0} b\,(a+\ri b)\left(a\ri - \frac{b}{3}\right)+\OO\!\bigl(b(|a|+b)^3\bigr).
\end{equation}
\end{enumerate}
\end{lemma}

\subsection{Ramification points}\label{s:ramification}
From the first statement in \Cref{p:surface}, the ramification points of the projection $(f,z)\in \cC\mapsto z\in \bC\bP$ correspond to points $\Psi(x_0,s_0), \overline{\Psi(x_0,s_0)}$ (recall from \eqref{e:emb}) with $(x_0,s_0)\in \fL$ and $\del_x \chi(x_0,s_0)=1/s_0$.  To parametrize a neighborhood of a ramification point 
\((f,z)=\Psi(x_0,s_0)\in\cC\), we perform a small time shift. 
Fix \(\ft>0\) and define the shifted polygon and liquid region
\begin{align}\label{e:shifted_polygon}
\fP_{\ft}:=\{(x,s+\ft):(x,s)\in\fP\},\qquad
\fL_{\ft}:=\{(x,s+\ft):(x,s)\in\fL\}.
\end{align}
Since the complex slope is transported by the shift, for 
\((x,s+\ft)\in\overline{\fL}_{\ft}\) 
\begin{align}\label{e:chit}
f_\ft(x, s+\ft)=f(x,s),\quad \chi_\ft(x,s+\ft)=\chi(x,s)
\end{align}
and
the map \eqref{e:emb} induces a map
\begin{align}\begin{split}\label{e:emb_ft}
\Psi_\ft: (x,s)\in \overline \fL\mapsto (x,s+\ft)\in \overline\fL_\ft&\mapsto \left(f_\ft(x,s+\ft), x-(s+\ft)\chi_\sft(x,s+\ft)\right)
\\
&=(f(x,s), z(x,s)-\ft \chi(x,s))\in \mathbb{CP}\times \bC\bP,
\end{split}\end{align}
By the same argument as in \Cref{p:surface}, the above map also induces a Riemann surface $\cC_\ft$ by
\begin{align}\label{e:shifted_RS}
\cC_\ft=\{(f,w)\in\bC\bP\times\bC\bP: Q(f, w+\ft f/(f+1))=0\}. 
\end{align}
The map from $\cC$ to $\cC_\ft$
\begin{align}\label{e:CtoCt}
(f,z)\in \cC\mapsto (f, z-\ft f/(f+1))\in \cC_\ft. 
\end{align} 
is a biholomorphism with inverse \((f,w)\mapsto (f,\,w+\ft f/(f+1))\), 
and it yields the commutative diagram in \Cref{f:commute}.

At a ramification point for \eqref{e:emb}, one has 
\(\partial_x\chi(x_0,s_0)=1/s_0\). After the shift by time $\ft$, the criticality condition
would be \(\partial_x\chi_\ft(x_0,s_0+\ft)=\partial_x \chi(x_0, s_0)=1/(s_0+\ft)\), which fails for \(\ft\neq0\);
hence \(\Psi_\ft(x_0, s_0)=(f(x_0,s_0), z(x_0,s_0)-\ft \chi(x_0,s_0))\) is a regular point of the projection \((f,w)\in \cC_{\ft}\mapsto w\in \bC\bP\).
Consequently, in a neighborhood of \(\Psi_\ft(x_0, s_0)\) we can solve for \(f\) as a
holomorphic function of \(w\), i.e.
\begin{align}\label{e:def_ft}
f=f_{\ft}(w).
\end{align}

	\begin{figure}
	\begin{tikzpicture}[scale=3]
  % Axes
  \node[anchor=base east](A) at (0,0.5) {$\overline \fL\ni(x,s)$};
  \node[anchor=base west] (B) at (0.5,0.5) {$(f,z)=\left(f(x,s), x-s\frac{f(x,s)}{f(x,s)+1}\right)\in \cC$};
  \node[anchor=base east] (C) at (0,0) {$\overline\fL_\ft\ni(x,s+\ft)$};
  \node[anchor=base west] (D) at (0.5,0)  {$(f,z-\ft f/(f+1))=\left(f(x,s), x-(s+\ft) \frac{f(x,s)}{f(x,s)+1}\right)\in \cC_\ft$};
     \node[anchor=base east] (E) at (3.4,0.5) {$\bC\bP$};
  % arrows
  \draw[->] (A) -- node[above]{$\Psi$} (B);
  \draw[->] (-0.17,0.43) -- (-0.17,0.13);
  \draw[->] (0.67,0.43) -- (0.67,0.13);
  \draw[->] (C) --  (D);
  \draw[->] (A) -- node[above]{$\Psi_\ft$}  (0.53,0.13);
  
    \draw[->] (B) -- node[above]{$\phi$} (E);
      \draw[->] (3.28,0.13) -- node[right]{$\phi_\ft$} (3.28,0.43);
       
\end{tikzpicture}

%	\begin{center}
%	 \includegraphics[scale=0.3,trim={0cm 5cm 0 7cm},clip]{complex_slope.pdf}
	 \caption{Shown above is the commutative diagram corresponding to shifting the polygon and the map \eqref{e:CtoCt}.}
	 \label{f:commute}
	 \end{figure}

\subsection{The uniformizing conformal map}\label{s:cf_map}
In this section, we explicitly construct the uniformizing conformal map
\(\phi\) used in \Cref{thm:moment-GFF} and collect some of its properties.

Recall from \eqref{e:emb} that the embedding
\begin{align}\label{e:emb_copy}
\Psi:
(x,s)\in\fL
\longmapsto
(f(x,s),z(x,s))
:=
\bigl(f(x,s),\,x-s\chi(x,s)\bigr)
\in
\cC:=\{(f,z):Q(f,z)=0\}
\end{align}
is a bijection from the liquid region onto one half of the Riemann surface
\(\cC\); see the fourth statement of \Cref{p:surface}.

By the second statement of \Cref{pla}, the liquid region \(\fL\) is simply
connected. Let \(\phi(f,z)\) be the Riemann map from
\[
\Psi(\fL)=\cC\cap(\bC_-\times\bC_+)
\]
onto the upper half-plane \(\bC_+\). By the Schwarz reflection principle,
\(\phi\) extends to \(\cC\) and satisfies
\[
\phi(\bar f,\bar z)=\overline{\phi(f,z)}.
\]

Recall from \eqref{e:inftypoint} that \(\cC\) contains \(d\) points at
infinity,
$
\{(-1,\infty_i)\}_{i=1}^d.
$
After composing \(\phi\) with a real Möbius transformation if necessary, we
may assume that
\begin{align}
\phi(-1,\infty_i)=p_i\in\bR\cup\{\infty\},
\qquad
p_1=\infty.
\end{align}
Locally near \((-1,\infty_1)\), the map \(\phi\) has the expansion
\begin{align}\label{e:inftoinf}
\phi(f,z)
=
\fa_1z+\OO(1),
\qquad
\fa_1\neq0.
\end{align}
Locally near \((-1,\infty_i)\), for \(2\leq i\leq d\), the map \(\phi\)
has the expansion
\begin{align}\label{e:inftofinite}
\phi(f,z)
=
p_i+\frac{\fa_i}{z}
+\OO\left(\frac{1}{z^2}\right),
\qquad
\fa_i\neq0.
\end{align}

The map $\phi$ induces a diffeomorphism from $\fL$ to $\bC_+$; with a slight abuse of notation, we denote this induced map by the same symbol:
\begin{align}\label{e:defphi}
\phi:\ (x,t)\in\fL \longmapsto (f,z)=\bigl(f(x,t),z(x,t)\bigr)
\longmapsto \phi(f,z),\qquad
\phi(x,t):=\phi\bigl(f(x,t),z(x,t)\bigr).
\end{align}

We record the following lemma, which will be used later to specify the
orientations of the contour integrals.

\begin{lemma}\label{c:defsqrtphi}
There exists a globally defined \(C^1\) choice of the square root
\[
\sqrt{\partial_x\phi(x,t)}
\]
on \(\fL\). This choice is unique up to an overall sign.
\end{lemma}

\begin{proof}
Since $\phi$ is a diffeomorphism, its Jacobian is invertible at every point; in particular, the complex derivative $\partial_x\phi(x,s)$ never vanishes on $\fL$, since otherwise the first column of the Jacobian would vanish.
Because $\fL$ is simply connected, $\partial_x\phi$ admits a continuous branch of the logarithm,
$ \ln(\partial_x\phi):\fL\to\bC$.
Define
\[
\sqrt{\partial_x\phi(x,s)}:=\exp\!\left(\tfrac12\, \ln(\partial_x\phi(x,s))\right).
\]
This yields a continuous square root on $\fL$, determined uniquely up to a global choice of sign.
\end{proof}

When the context is clear, in particular when $(f,z)$ is bounded away from
ramification points, we will write $\phi(z):=\phi(f,z)$. In this regime,
differentiating \eqref{e:defphi} with respect to $x$ gives
\begin{align}\label{e:derphi}
\partial_x\phi(x,s)
=
\partial_x z(x,s)\,\phi'\bigl(z(x,s)\bigr)
=
\frac{\phi'\bigl(z(x,s)\bigr)}
     {1+s\,\chi'\bigl(z(x,s)\bigr)}.
\end{align}

In the following, we discuss the case near ramification points. Suppose
$(f,z)$ is close to a ramification point corresponding to
$\Psi(x_0,s_0)$ or $\overline{\Psi(x_0,s_0)}$, where
$(x_0,s_0)\in\fL$ and
$
    \partial_x\chi(x_0,s_0)=1/s_0$.
Recall the shifted polygon $\fP_{\ft}$ and liquid region $\fL_\ft$ from \eqref{e:shifted_polygon}, and the
associated Riemann surface $\cC_{\ft}$ from \eqref{e:shifted_RS}. Define the
induced map
\begin{align}\label{e:defphi_t}
\phi_{\ft}:\left(f,w=z-\ft f/(f+1)\right)\in\cC_{\ft}
\longmapsto \phi(f,z),\qquad (f,z)\in\cC,
\end{align}
so that the diagram in \Cref{f:commute} commutes. In a neighborhood of
$\Psi_{\ft}(x_0,s_0)$, we may solve $f=f_{\ft}(w)$ as a holomorphic function
of $w$; see \eqref{e:def_ft}. We then write
$\phi_{\ft}(w):=\phi_{\ft}(f_{\ft}(w),w)$.
Similarly to \eqref{e:defphi}, $\phi_\ft$ also induces a diffeomorphism from the shifted liquid region $\fL_\ft$ to $\bC_+$. We also denote it by the same symbol $\phi_\ft$, and
\begin{align}\label{e:defphi_t2}
\phi_\ft(x,s+\ft)=\phi(x,s)
\end{align}

\subsection{Tangent locations}
\label{s:tangent_location}
\begin{figure}
			
			\begin{center}		
				
				\begin{tikzpicture}[
					>=stealth,
					auto,
					style={
						scale = .475
					}
					]
					
					\node[] at (0,-0.2){\includegraphics[scale=0.6]{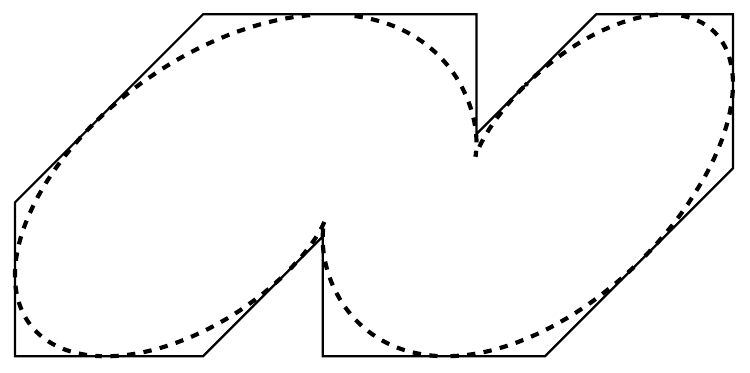}};
					%\draw[black, thick,->] (-7.65,-3) -- (8.4, -3);
					%\draw[black, thick,->] (-7.65,-3) -- (-7.65, 3.4);
					
					\filldraw[fill = black] (-7.7,-2) circle (0.1);
					\draw[] (-7.7,-2)  node[left, scale=0.6]{$(x_1^{(\infty)}, s_1^{(\infty)})$} ;
					
					\filldraw[fill = black] (-5.8,-3.82) circle (0.1);
					\draw[] (-5.8,-3.82) node[below, scale=0.6]{$(x_1^{(0)}, s_1^{(0)})$} ;

					\filldraw[fill = black] (-6,1.2) circle (0.1);
					\draw[] (-6,1.2) node[left, scale=0.6]{$(x_d^{(1)}, s_d^{(1)})$} ;

						\filldraw[fill = black] (-2,-2.15) circle (0.1);
						\draw[] (-2,-2.15) node[left, scale=0.6]{$(x_1^{(1)}, s_1^{(1)})$} ;
						
					\filldraw[fill = black] (-1,3.5) circle (0.1);
					\draw[] (-1,3.5) node[above, scale=0.6]{$(x_d^{(0)}, s_d^{(0)})$} ;

					\filldraw[fill = black] (2.15,1) circle (0.1);
					\draw[] (2.15,1) node[left, scale=0.6]{$(x_d^{(\infty)}, s_d^{(\infty)})$} ;

					 \draw (15,0) ellipse (5 and 2.5);	
					 \draw[] (15,0)  node[]{$\cC(\bR)$};

					 \filldraw[fill = black] (20,0) circle (0.05);
					\draw[] (20,0)  node[right, scale = .6]{$(\infty, a_1)$};
					
					\filldraw[fill = black] (10,0) circle (0.05);
					\draw[] (10,0)  node[left, scale = .6]{$(\infty, a_d)$} ;

					\filldraw[fill = black] (16,-2.45) circle (0.05);
					\draw[] (16,-2.45)  node[below, scale = .6]{$(0,b_1)$} ;
				
					\draw[] (16,2.45)  node[above, scale=0.8]{$\cdots$} ;

					\filldraw[fill = black] (18,-2.01) circle (0.05);
					\draw[] (18,-2.01)  node[below, scale = .6]{$(-1,\infty_1)$} ;

					\filldraw[fill = black] (12,2) circle (0.05);
					\draw[] (12,2)  node[above, scale = .6]{$(0,b_d)$} ;
					
					\filldraw[fill = black] (11,1.5) circle (0.05);
					\draw[] (11,1.5) node[left, scale = .6]{$(-1, \infty_d)$} ;

				\end{tikzpicture}
				
			\end{center}
			
			\caption{In the left panel, the arctic curve is tangent to the boundary of \(\mathfrak P\), in counterclockwise order, at $(x_1^{(\infty)}, s_1^{(\infty)}), (x_1^{(0)}, s_1^{(0)}), (x_1^{(1)}, s_1^{(1)}),\cdots, (x_d^{(\infty)}, s_d^{(\infty)}), (x_d^{(0)}, s_d^{(0)}), (x_d^{(1)}, s_d^{(1)})$; In the second subplot, the embedding map \eqref{e:emb} sends these tangency points to $(0, b_1), (-1, \infty_1), (\infty, a_1), \cdots, (0,b_d), (-1,\infty_d), (\infty, a_d)\in \cC(\bR)$.
			}\label{f:12-gon}
			
		\end{figure}
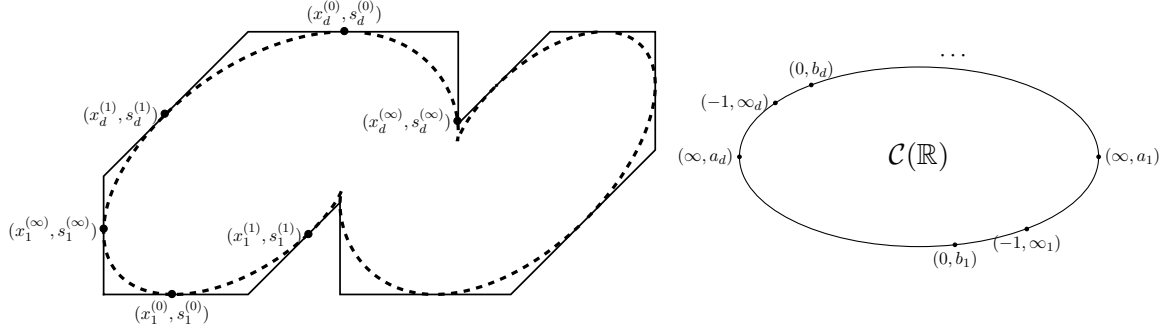

As we move counterclockwise along the boundary of the polygon
\(\mathfrak P\), the arctic boundary \(\fA\) is tangent successively to each
side, or its extension, of \(\mathfrak P\) at the points
\[
(x_1^{(\infty)},s_1^{(\infty)}),\,
(x_1^{(0)},s_1^{(0)}),\,
(x_1^{(1)},s_1^{(1)}),\,
\ldots,\,
(x_d^{(\infty)},s_d^{(\infty)}),\,
(x_d^{(0)},s_d^{(0)}),\,
(x_d^{(1)},s_d^{(1)}).
\]
The superscript indicates the slope; more precisely,
\((x_i^{(0)},s_i^{(0)})\), \((x_i^{(1)},s_i^{(1)})\), and
\((x_i^{(\infty)},s_i^{(\infty)})\) are the tangency points associated with
the sides of \(\mathfrak P\) of slopes \(0\), \(1\), and \(\infty\),
respectively, as shown in \Cref{f:12-gon}.
For $1\leq i\leq d$, using the relation \eqref{e:slope}
\begin{align}\label{e:bcond}
 f(x_i^{(\infty)},s_i^{(\infty)})=0,\quad
 f(x_i^{(0)},s_i^{(0)})=-1,\quad
 f(x_i^{(1)},s_i^{(1)})=\infty.
\end{align}
 Using the map \eqref{e:emb}, those tangency points correspond to points on the real locus $\cC(\bR)$, for $1\leq i\leq d$
  \begin{align}\begin{split}\label{e:defab}
      &( f(x_i^{(\infty)},s_i^{(\infty)}), x_i^{(\infty)}-s_i^{(\infty)} f(x_i^{(\infty)},s_i^{(\infty)})/( f(x_i^{(\infty)},s_i^{(\infty)})+1))
 =(0, x_i^{(\infty)})=:(0,b_i),\\
&( f(x_i^{(0)},s_i^{(0)}), x_i^{(0)}-s_i^{(0)} f(x_i^{(0)},s_i^{(0)})/( f(x_i^{(0)},s_i^{(0)})+1))
 =
 (-1, \infty_i), \\
  &( f(x_i^{(1)},s_i^{(1)}), x_i^{(1)}-s_i^{(1)} f(x_i^{(1)},s_i^{(1)})/( f(x_i^{(1)},s_i^{(1)})+1))
 =(\infty, x_i^{(1)}-s_i^{(1)})=:(\infty, a_i),
 \end{split}\end{align}
 These points lie on the real locus \(\cC(\bR)\) in counterclockwise order; See Figure~\ref{f:12-gon} for an example of $12$-gon.
When the context is clear, we will represent these points by their \(z\)-coordinate.

The point \((x_i^{(\infty)},s_i^{(\infty)})\) lies on the supporting line of a
vertical side of \(\fP\). Hence
$
b_i=x_i^{(\infty)}
$
is the constant \(x\)-coordinate of that supporting line and depends only on
the corresponding side. Similarly,
\((x_i^{(1)},s_i^{(1)})\) lies on the supporting line of a side of slope \(1\).
Consequently,
$
a_i=x_i^{(1)}-s_i^{(1)}
$
is the constant value of \(x-s\) along that supporting line and depends only
on the corresponding side. Therefore, the marked points
\[
(0,b_1),\ (-1,\infty_1),\ (\infty,a_1),\ \ldots,\ 
(0,b_d),\ (-1,\infty_d),\ (\infty,a_d)
\]
can be determined directly from the polygonal domain \(\fP\), without knowing
the arctic curve. With a slight abuse of notation, we associate the vertical,
horizontal, and unit-slope sides with \(b_i\), \(\infty_i\), and \(a_i\),
respectively.

Moreover, from our \Cref{p} of $\fP$, we have
\begin{align}\label{e:aibi}
na_i=n(x_i^{(1)}-s_i^{(1)}), \quad n b_i=n x_i^{(\infty)}\in \bZ'=\bZ+\frac{1}{2}, \quad 1\leq i\leq d.
\end{align}

Finally, under a vertical translation
$
(x,s)\mapsto(x,s+C),
$
each \(a_i\) is replaced by \(a_i-C\), whereas each \(b_j\) remains unchanged.
By translating \(\fP\) upward if necessary, we may assume that
\begin{align}\label{e:ab_order}
s>0 \quad\text{for all }(x,s)\in\fP,
\qquad
a_i=x_i^{(1)}-t_i^{(1)}
<b_j=x_j^{(\infty)},
\qquad 1\leq i,j\leq d.
\end{align}

For \(1\leq i\leq d\), with indices taken cyclically, the marked points
$
\cdots, b_i, \infty_i, a_i, b_{i+1},\cdots
$
occur in this order along the real locus \(\cC(\bR)\). They divide the real
locus into three open arcs. Under the map \eqref{e:emb}, each of these
arcs corresponds to a portion of the arctic boundary contained in a
curvilinear triangle. The range of \(f\) and the gradient of the limiting height function on the
corresponding curvilinear triangle are given by
\begin{equation}\label{e:arcCR}
\begin{aligned}
\text{on }(b_i,\infty_i):
&\qquad -1<f<0,
&\qquad \nabla H^*&=(1,0),\\
\text{on }(\infty_i,a_i):
&\qquad -\infty<f<-1,
&\qquad \nabla H^*&=(1,-1),\\
\text{on }(a_i,b_{i+1}):
&\qquad 0<f<\infty,
&\qquad \nabla H^*&=(0,0).
\end{aligned}
\end{equation}
Here and below, the intervals are understood as intervals in the real
locus $\cC(\bR)$.

\section{Tiling Action and Critical Points}
\label{s:tilingaction}

In this section, we introduce the \emph{tiling action} and study its properties and critical points.

\subsection{Tiling action function and complex Burgers equation}

Throughout this article we fix the branch of the logarithm as
\(\ln:\bC\setminus(-\infty, 0]\to\bC\). Then $\ln f$ is well defined when $f\in \bC\setminus (-\infty,0]$. In particularly, $\ln f(z)=\ln |f|+\ri \arg(f)$ is well defined on $\cC\setminus \{(f,z): f\in [-\infty,0]\}$. Recall from \eqref{e:arcCR}, this region is obtained from $\cC$ by removing the cuts $[b_i, a_i]\subset \cC(\bR)$ for $1\leq i\leq d$. 

\begin{lemma}\label{c:logf_dz_holo}
The $1$-form $\ln f\,\rd z$ is holomorphic on $\cC\setminus\bigl(\cup_{i=1}^d [b_i,a_i]\bigr)$,
where $\ln$ denotes the fixed branch $\ln:\bC\setminus(-\infty,0]\to\bC$ and we view $\ln f$
as the composition $(f,z)\mapsto f \mapsto \ln f$.
\end{lemma}

\begin{proof}
By the definition of the cuts (see \eqref{e:arcCR}), the projection $(f,z)\mapsto f$ maps
$\cC\setminus\bigl(\cup_{i=1}^d [b_i,a_i]\bigr)$ into $\bC\setminus(-\infty,0]$.
Hence $\ln f$ is a single-valued holomorphic function on
$\cC\setminus\bigl(\cup_{i=1}^d [b_i,a_i]\bigr)$, as the composition of the holomorphic map $f$
with the chosen holomorphic branch of $\ln$.

Thanks to \eqref{e:introducef0}, away from ramification points we may view $f=f(z)$
as a holomorphic function of $z$, and in particular $\rd z$ is holomorphic there.
In a neighborhood of a ramification point, by \eqref{e:introduceQ0} one can instead use $f$ as a local coordinate: writing
$z=Q_0(f)$, we have $\rd z = Q_0'(f)\,\rd f$, which is again holomorphic.
Therefore $\rd z$ is a holomorphic $1$-form on
$\cC\setminus\bigl(\cup_{i=1}^d [b_i,a_i]\bigr)$.

Since both $\ln f$ and $\rd z$ are holomorphic on this domain, their product $\ln f\,\rd z$ is a
holomorphic $1$-form on $\cC\setminus\bigl(\cup_{i=1}^d [b_i,a_i]\bigr)$.
\end{proof}

\begin{definition}[Tiling action]
We fix a basepoint $(0,b_1)\in\cC(\bR)$ (recall \eqref{e:defab}); after translating the $x$--coordinate, we may assume $b_1=0$.
Given any $(x,s)\in\fP$, we define an action on $\cC\setminus \cC(\bR)$ by
\begin{align}\label{e:def_action_C}
S\bigl(f,z;x,s\bigr)
:= s\ln s-(x-z)\ln (x-z)-(s-x+z)\ln (s-x+z)
-\int_{\sfC} \ln f(u)\,\rd u,
\end{align}
where the integral is taken along any path $\sfC$ parametrized by
$u\mapsto (f(u),u)\in\cC$ from $(0,0)$ to $(f,z)$ that lies in
$\cC\setminus \cC(\bR)$ except the starting point $(0,0)$ (so that $u$ denotes the $z$--coordinate along the path).
\end{definition}

Since $\cC\setminus \cC(\bR)$ has two simply connected components and the $1$-form $\ln f\,\rd z$
(equivalently $\ln(f(u))\,\rd u$ under the above parametrization) is holomorphic on each component,
the integral is well-defined and independent of the choice of $\sfC$.
We also note that the first three logarithmic terms are well-defined provided
\[
z\in \bC\setminus\bigl((-\infty, x-s]\cup [x,\infty)\bigr).
\]

When $(f,z)\in\cC\setminus \cC(\bR)$ is bounded away from ramification points, we may work on a fixed local sheet and identify $(f,z)$ with its $z$-coordinate.
With this identification we write
\begin{align}\label{e:def_action} 
S(z;x,s):=S\bigl(f,z;x,s\bigr)
= s\ln s-(x-z)\ln (x-z)-(s-x+z)\ln (s-x+z)
-\int_{0}^{z} \ln f(u)\,\rd u,
\end{align}
where the integral is taken along any path in the $z$-plane avoiding the cuts.
Differentiating \eqref{e:def_action}  yields
\begin{align}\label{e:critical}
S'(z;x,s)=\ln\!\left(\frac{x-z}{z-(x-s)}\right)-\ln f(z).
\end{align}
We defer the discussion of the case that $(f,z)$ is in a neighborhood of ramification points to \Cref{s:changetime}.

Next we give a heuristic motivation for the action \eqref{e:def_action}.
Specifically, we obtain \eqref{e:def_action} by formally solving the complex Burgers-type
equation
\begin{align}\label{e:ftxcopy}
\del_s \ln f_s+\del_x \ln\!\bigl(f_s+1\bigr)=0,
\end{align}
with initial condition \(f_0(x)=f(x,0)\). The discussion in this section is purely
heuristic and is included only for intuition; none of the conclusions will be used later.

We introduce the following Cole--Hopf type transform:
\begin{align}\label{e:colehopf}
\ln f(x,s)=-\frac{1}{n}\,\partial_x \ln \psi(x,s)=-\frac{1}{n}\,\frac{\partial_x\psi(x,s)}{\psi(x,s)},\quad \psi(x,0):=e^{-n\int_0^x \ln f(u,0)\rd u}.
\end{align}

Substituting \eqref{e:colehopf} into \eqref{e:ftxcopy} and integrating once in \(x\)
gives
\begin{align}\label{e:complexeq}
-\frac{1}{n}\del_s \ln \psi(x,s)
+\ln\left(1+e^{-\frac{1}{n}\del_x \ln \psi(x,s)}\right)
=C(s).
\end{align}
The function \(C(s)\) can be removed by multiplying \(\psi\) by a scalar function of
\(t\). Since such a multiplication does not change
\(-n^{-1}\del_x \ln\psi\), it does not change \(f\). Therefore, without loss of
generality, we set \(C(s)=0\).

With this normalization, \eqref{e:complexeq} becomes
\[
e^{-\frac{1}{n}\del_s \ln \psi(x,s)}
\left(
1+e^{-\frac{1}{n}\del_x \ln \psi(x,s)}
\right)=1.
\]
Multiplying by \(\psi(x,s)=e^{ \ln\psi(x,s)}\), we obtain
\begin{align}\label{e:formalshift}
\psi(x,s)=e^{ \ln\psi(x,s)-\frac{1}{n}\del_s \ln\psi(x,s)} +e^{ \ln\psi(x,s)
-\frac{1}{n}\del_s \ln\psi(x,s)
-\frac{1}{n}\del_x \ln\psi(x,s)}.
\end{align}
Interpreting the exponentials in \eqref{e:formalshift} as first-order Taylor
approximations to lattice shifts, we get
\begin{align}\label{e:discreteheat}
\psi(x,s)\approx \psi\left(x,s-\frac{1}{n}\right)
+\psi\left(x-\frac{1}{n},s-\frac{1}{n}\right),
\qquad x,s\in \frac{1}{n}\bZ .
\end{align}
Thus, at the level of this formal approximation, \(\psi\) satisfies the discrete
heat equation.

For \(nt\in \bZ_{\ge 0}\), the discrete heat kernel gives
\begin{align}\label{e:kernel_represent}
\psi(x,s)
&\approx\sum_{z\in \bZ/n}
{ns\choose n(x-z)} \psi(z,0)
\approx n \int_\bR
{ns\choose n(x-z)} e^{-n\int_0^z f_0(u)\rd u}\rd z.
\end{align}
We now apply Stirling's formula. Recall from \cite{spira1971calculation} that, uniformly for $z\in \bC\setminus (-\infty, 0]$,
\begin{align}\begin{split}\label{e:logGamma}
&\ln \Gamma(z)=(z-1/2)\ln z -z+\frac{1}{2}\ln (2\pi) +\OO\left(\frac{1}{\dist(z, (-\infty, 0])}\right).
\end{split}\end{align}

Hence, for \(w\in \bC\) such that $\delta:=\dist(w, (-\infty,x-s]\cup[x,\infty))\gg n^{-1}\),  we have
\begin{align}\begin{split}\label{e:binomial}
\ln {ns\choose n(x-w)}
&=
\ln\Gamma(ns+1)
-\ln\Gamma(n(x-w)+1)
-\ln\Gamma(n(s-x+w)+1) \\
&=
n\Bigl[
s\ln s
-(x-w)\ln(x-w)
-(s-x+w)\ln(s-x+w)
\Bigr]  \\
&
+\frac12\Bigl[
\ln s-\ln(x-w)-\ln(s-x+w)
\Bigr]
-\frac12\ln(2\pi n)
+\OO\left(\frac1{\delta n}\right).
\end{split}\end{align}
Substituting \eqref{e:binomial} into \eqref{e:kernel_represent}, we obtain
\begin{align}\label{e:psi_action_integral}
\psi(x,s)
&\approx
n\int_{x-s}^{x}
\sqrt{\frac{s}{2\pi n(x-z)(s-x+z)}}
\exp\{nS(z;x,s)\}\,\rd z,
\end{align}
where the phase is
\begin{align}\label{e:heuristic_action}
S(z;x,s)
&:=
s\ln s
-(x-z)\ln(x-z)
-(s-x+z)\ln(s-x+z)
-\int_0^z \ln f_0(u)\,\rd u .
\end{align}
This is the action as defined in \eqref{e:def_action}.

\subsection{Properties of Tiling Action Function}
In this section, we collect several elementary properties of the tiling action
function \eqref{e:def_action}. We recall from \eqref{e:arcCR} that the set
on which \(f\in[-\infty,0]\) is given by the union of cuts
$\cup_{i=1}^d [b_i,a_i]\subset \cC(\bR)$.

The first lemma says that the imaginary part of the primitive of \(\ln f\), after
multiplication by \(n\), is single-valued modulo \(2\pi\) away from the cuts, and
has an explicit jump across each cut.

\begin{lemma}\label{c:deflnf}
Fix \(ns,nx\in \bZ\). On $\cC\setminus\bigl(\cup_{i=1}^d [b_i,a_i]\bigr)$ the quantity
\begin{align}\label{e:defint}
\Im\left[n\int_0^z \ln f(u)\rd u\right]
\end{align}
is well defined modulo \(2\pi\), that is, it is independent modulo \(2\pi\) of the
path in $\cC\setminus\bigl(\cup_{i=1}^d [b_i,a_i]\bigr)$.
Moreover, for every finite $E\in [b_i,\infty_i)\cup(\infty_i,a_i]$,
where \(f(E)<0\) on the real locus, the non-tangential boundary values from the
upper and lower half-planes satisfy
\begin{align}\label{e:intdiff}
\Im\left[n\int_0^{E+0\ri}\ln f(u)\rd u\right]
-
\Im\left[n\int_0^{E-0\ri}\ln f(u)\rd u\right]
=
(1-2En)\pi
\pmod{2\pi}.
\end{align}
Here the signs \(+\) and \(-\) indicate the half-plane from which \(E\) is
approached in the \(z\)-plane. Since \(z\in\bC_+\) implies
\(f(z)\in\bC_-\), the boundary value \(f(E+\ri0)\) approaches the negative
real axis from below. Hence
$
\arg f(E+\ri0)=-\pi,
$ and $
\arg f(E-\ri0)=\pi.
$
\end{lemma}

The following lemmas record the corresponding jump of the logarithmic terms coming from
the binomial coefficient when \(w\) crosses one of its two real cuts.
\begin{lemma}\label{c:b-w}
We recall from \eqref{e:aibi} that $na_i, nb_i\in \bZ'=\bZ+1/2$. 
For any $E\in [b_i, \infty_i)\cup (\infty_i, a_i]$, we also have
\begin{align}\label{e:Imab}
\left.\Im\left[-n(b_i-w)\ln (b_i-w)-n(w-a_i)\ln (w-a_i)\right]\right|^{w=E+0\ri}_{w=E-0\ri }
=(1-2En)\pi \mod 2\pi.
\end{align}
\end{lemma}

%
%{\color{red}[not sure if we need]
%\begin{lemma}\label{c:walk-term}
%Fix \(ns,nx\in \bZ\). For every finite
%$
%E\in (-\infty,x-s]\cup[x,\infty)$,
%we have
%\begin{align}
%\left.
%\Im\left[
%-\left(n(x-w)+\frac12\right)\ln(x-w)
%-\left(n(s-x+w)+\frac12\right)\ln(s-x+w)
%\right]
%\right|_{w=E-0\ri}^{w=E+0\ri}
%=
%(1-2En)\pi
%\pmod{2\pi}.
%\end{align}
%\end{lemma}
%
%
%\begin{proposition}\label{p:integrand_analyticity}
%The integrand in \eqref{e:psi_action_integral}
%\begin{align}
%\exp\{nS(w;x,s)-\frac{1}{2}(\ln(x-w)+\ln(s-x+w))\}
%\end{align}
%extends analytically to 
%\begin{align}
%\cC\setminus\bigl(\cup_{i=1}^d [\min\{b_i, x\}, \max\{b_i, x\}]\cup \min[\{a_i, x-s\}, \max\{a_i, x-s\}]\bigr)
%\end{align}
%\end{proposition}}

\begin{proof}[Proof of \Cref{c:deflnf}]
We first prove that \eqref{e:defint} is well defined modulo \(2\pi\). It is
enough to check that, for any contour \(\omega_+\subset \cC\) surrounding one cut
\([b_i,a_i]\), we have
\begin{align}\label{e:period_lnf}
\Im\left[
n\int_{\omega_+}\ln f(u)\rd u
\right]
=0
\mod 2\pi .
\end{align}

By \eqref{e:defab} and \eqref{e:arcCR}, the arc \([b_i,a_i]\) contains
\(\infty_i\). Moreover, on \((b_i,\infty_i)\) we have \(-1<f<0\), while on
\((\infty_i,a_i)\) we have \(-\infty<f<-1\). The point
\((-1,\infty_i)\in\cC\) corresponds to a tangent location \((x,s)\in\fA\) with
\(f(x,s)=-1\). By the third statement of \Cref{p:surface}, locally near
\((-1,\infty_i)\) we have
\begin{align}\label{e:fzexp}
f(z)=\frac{\chi(z)}{1-\chi(z)},\qquad
\chi(z)=\frac{-z+x}{s}+\OO\left(\frac1z\right),
\qquad
f(z)=-1+\frac{s}{z}+\OO\left(\frac1{z^2}\right).
\end{align}

Let \(\omega\subset\cC\) be a contour surrounding \((-1,\infty_i)\), whose projection
to the \(z\)-coordinate meets the real locus at \(z=\pm R\), where \(R>0\) is large.
Deforming \(\omega_+\) to \(\omega\), together with the two sides of the cut, gives
\begin{align}\label{e:intlnf1}
\begin{split}
\int_{\omega_+}\ln f(u)\rd u
&=
\int_{\omega}\ln f(u)\rd u
+\int_{b_i-0\ri}^{R-0\ri}\ln f(u)\rd u
+\int_{R+0\ri}^{b_i+0\ri}\ln f(u)\rd u  \\
&
+\int_{-R-0\ri}^{a_i-0\ri}\ln f(u)\rd u
+\int_{a_i+0\ri}^{-R+0\ri}\ln f(u)\rd u  =
\int_{\omega}\ln f(u)\rd u
+2\pi\ri(R-b_i)+2\pi\ri(a_i+R).
\end{split}
\end{align}
Here we used the boundary values of \(\ln f\) across the cut. Namely, for
\(E\in [b_i,a_i]\) on the real locus,
$
\arg f(E+0\ri)=-\pi,
\arg f(E-0\ri)=\pi$,
and therefore
\[
\ln f(E-0\ri)-\ln f(E+0\ri)=2\pi\ri .
\]

It remains to evaluate the contribution from the small loop around \(\infty_i\).
\begin{align}\label{e:intlnf2}
\begin{split}
\int_{\omega}\ln f(u)\rd u
&=
\int_{\omega}
\ln\left(-1+\frac{s}{u}+\OO\left(\frac1{u^2}\right)\right)\rd u \\
&=
\int_{\omega}
\ln\left(1-\frac{s}{u}+\OO\left(\frac1{u^2}\right)\right)\rd u
-\int_{\omega\cap\bC_+}\pi\ri\,\rd u
+\int_{\omega\cap\bC_-}\pi\ri\,\rd u =
-4R\pi\ri,
\end{split}
\end{align}
where the first statement follows from \eqref{e:fzexp}; the second statement is from our choice for the branch of $\ln(\cdot)$; and the third statement follows from performing the integral. Combining \eqref{e:intlnf1} and \eqref{e:intlnf2}, we get
\[
\int_{\omega_+}\ln f(u)\rd u
=
2\pi\ri(a_i-b_i).
\]
Therefore
\[
\Im\left[
n\int_{\omega_+}\ln f(u)\rd u
\right]
=
2\pi n(a_i-b_i)
=0
\mod 2\pi,
\]
where the last equality follows from \eqref{e:aibi}. This proves the first statement in \Cref{c:deflnf}, that
\eqref{e:defint} is independent of the path modulo \(2\pi\).

We now prove the jump formula \eqref{e:intdiff}. First assume that \(E\in [b_i,\infty_i)\). Up to multiples of \(2\pi\), the
parts of the two paths from \(0\) to \(b_i\) cancel, so the left-hand side of
\eqref{e:intdiff} equals
\begin{align}
\Im\left[
n\int_{b_i}^{E+0\ri}\ln f(u)\rd u
\right]
-
\Im\left[
n\int_{b_i}^{E-0\ri}\ln f(u)\rd u
\right].
\end{align}
Using that for $u\in [b_i,E]$,
$
\arg f(u+0\ri)=-\pi$,
and $
\arg f(u-0\ri)=\pi$,
we obtain
\begin{align}
&\phantom{{}={}}\Im\left[
n\int_{b_i}^{E+0\ri}\ln f(u)\rd u
\right]
-
\Im\left[
n\int_{b_i}^{E-0\ri}\ln f(u)\rd u
\right]\\
&=
-n(E-b_i)\pi-n(E-b_i)\pi =
-2n(E-b_i)\pi =
(1-2En)\pi
\mod 2\pi,
\end{align}
where in the last step we used \(nb_i\in \bZ+1/2\) from \eqref{e:aibi}.
The case \(E\in(\infty_i,a_i]\) is analogous, so we omit. 
\end{proof}

\begin{proof}[Proof of \Cref{c:b-w}]
First assume that \(E\in [b_i,\infty_i)\). Then
\begin{align}
\Im[\ln(b_i-(E\pm0\ri))]=\mp\pi,
\qquad
\Im[\ln(E\pm0\ri-a_i)]=0.
\end{align}
Hence
\begin{align}
\left.\Im\left[-n(b_i-w)\ln (b_i-w)-n(w-a_i)\ln (w-a_i)\right]\right|^{w=E+0\ri}_{w=E-0\ri }=
2\pi n(b_i-E) =
(1-2En)\pi
\mod 2\pi,
\end{align}
where we used \(nb_i\in\bZ+1/2\).
The case that \(E\in(\infty_i,a_i]\) is analogous, so we omit. 
\end{proof}

\begin{figure}
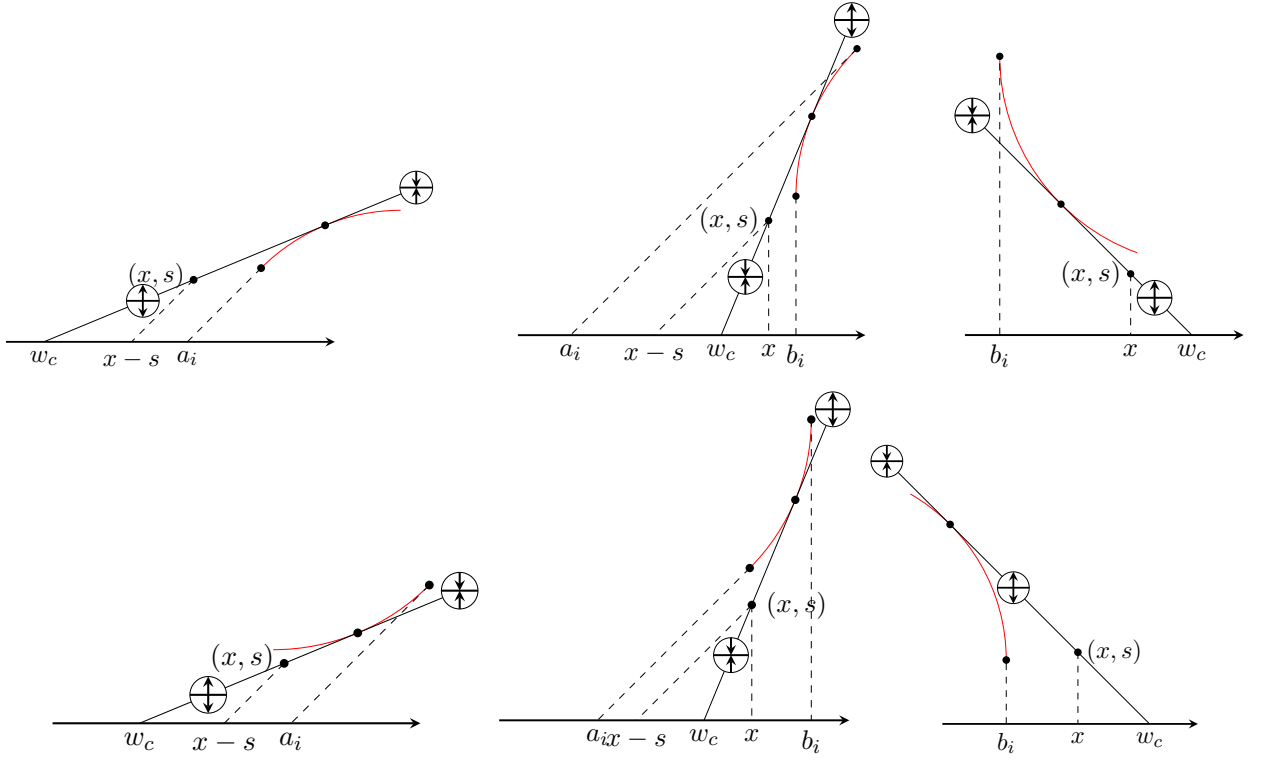

\centering

%================= Row 1 =================
\begin{subfigure}[t]{0.35\textwidth}
  \centering
  \resizebox{\linewidth}{!}{%

  }
\end{subfigure}

\medskip

%================= Row 2 =================
\tikzset{
  marked point/.style={circle, fill, inner sep=1.2pt},
  circled mark/.style={
    circle,
    draw,
    fill=white,
    inner sep=5pt
  },
  circledOut/.pic={
    \node[circled mark] {};
    \draw[stealth-, thick] (0,-0.25) -- (0,0);
    \draw[-stealth, thick] (0,0) -- (0,0.25);
    \draw[thick] (-0.25,0) -- (0.25,0);
  },
  circledIn/.pic={
    \node[circled mark] {};
    \draw[-stealth, thick] (0,-0.25) -- (0,0);
    \draw[stealth-, thick] (0,0) -- (0,0.25);
    \draw[thick] (-0.25,0) -- (0.25,0);
  }
}
\begin{subfigure}[t]{0.35\textwidth}
  \centering
  \resizebox{\linewidth}{!}{%

  }
\end{subfigure}

\caption{Critical points}
\label{f:critical}
\end{figure}

\subsection{Critical points and tangent lines}
\label{s:critical_points}

We recall from \eqref{e:critical} that, for any \((x,s)\in\fP\), the
formal critical points \(w_c\) of \(S(\,\cdot\,;x,s)\) satisfy
\begin{align}\label{e:critical_original}
f(w_c)=\frac{w_c-x}{x-s-w_c}.
\end{align}
Equivalently,
\begin{align}\label{e:critical_point}
w_c
=
x-s\frac{f(w_c)}{f(w_c)+1}
=
x-s\chi(w_c),
\end{align}
where \(\chi=f/(f+1)\), as in \eqref{e:def-chi-on-C2}. We refer to
\eqref{e:critical_point} as the \emph{critical-point equation}. In the
degenerate case \(\chi(w_c)=\infty\), we interpret
\eqref{e:critical_point} projectively and allow \(w_c=\infty\).

For each \((x,s)\in\fL\), let \(z(x,s)\) be the point defined in
\eqref{e:emb}. By \eqref{e:bulk_der}, it satisfies the critical-point
equation. We therefore associate with \((x,s)\) the pair of
complex-conjugate critical points
\begin{align}\label{e:two_points}
w_c:=z(x,s)=x-s\chi(z(x,s)),
\qquad
\overline{w_c}=\overline{z(x,s)}.
\end{align}

We next describe the real solutions of \eqref{e:critical_point}
geometrically. Let \((x',s')\in\fA\), and suppose that the line through
\((x,s)\) and \((x',s')\) is tangent to the arctic boundary \(\fA\) at
\((x',s')\). Let \(w_c\) be the intersection of this tangent line with the
axis \(\{s=0\}\). If the tangent line is horizontal, we use the projective
convention \(w_c=\infty\).

Suppose first that \(w_c\) is finite. By \Cref{pa1}, the slope of the
tangent line is \(1/\chi(w_c)\). Since the line passes through
\((w_c,0)\) and \((x',s')\), we have
\[
\frac{s'}{x'-w_c}=\frac{1}{\chi(w_c)}.
\]
Consequently,
\begin{align}\label{e:tangent_wc}
w_c
=
x'-s'\chi(w_c)
=
x-s\chi(w_c),
\end{align}
where the second equality follows because the same line passes through
\((x,s)\). Thus, \(w_c\) satisfies the critical-point equation
\eqref{e:critical_point}. The horizontal case is understood by the same
projective convention.

\begin{definition}\label{d:formal_genuine_critical}
A point \(w_c\) obtained either as one of the complex-conjugate solutions
of the critical-point equation in the liquid region or from a tangent line
to \(\fA\) is called a \emph{formal critical point}. A formal critical
point is called \emph{genuine} if it is a zero of
\(S'(\,\cdot\,;x,s)\), and \emph{spurious} otherwise.
\end{definition}

The distinction in \Cref{d:formal_genuine_critical} is relevant only for
tangent lines of slopes \(0\), \(1\), and \(\infty\). At these exceptional
slopes, logarithmic singularities in the expression for \(S'\) may cancel,
so a formal solution of \eqref{e:critical_point} need not be a genuine
critical point.

The following proposition summarizes the multiplicities of the formal and
genuine critical points arising from tangent lines.

\begin{proposition}\label{p:critical_multiplicity}
Suppose first that the tangent line has slope outside
\(\{0,1,\infty\}\). Then the formal and analytic multiplicities agree:
\begin{enumerate}
\item
At a regular arctic point, two tangent lines coalesce, and the corresponding
critical point has multiplicity two.

\item
At a cusp point, three tangent lines coalesce, and the corresponding
critical point has multiplicity three.

\item
Any other tangent line gives a simple genuine critical point.
\end{enumerate}

Suppose next that the tangent line has slope in
\(\{0,1,\infty\}\), so that it is the supporting line of an extended side
of \(\fP\). Then logarithmic cancellation reduces the analytic
multiplicity by one:
\begin{align}\label{e:formal_analytic_multiplicity}
\begin{array}{c|c|c}
\text{Location}
&
\text{Formal multiplicity}
&
\text{Analytic outcome}
\\
\hline
\text{Regular point of an extended side}
&1&\text{spurious point}
\\
\text{Tangency point}
&2&\text{simple genuine critical point}
\\
\text{Cusp-turning point}
&3&\text{genuine critical point of multiplicity two}.
\end{array}
\end{align}
\end{proposition}

\begin{proof}
Away from the exceptional slopes, the assertions follow from the usual
coalescence of tangent lines and the local expansions of the tiling action
at regular arctic and cusp points. We refer to \Cref{subsec:action-derivatives}.

For slopes \(0\), \(1\), and \(\infty\), the tangent line coincides with
the supporting line of an extended side. The corresponding local
expressions for the tiling action contain two logarithmic terms with
canceling singular parts. The cancellation removes one formal factor from
\(S'\), giving precisely the three outcomes in
\eqref{e:formal_analytic_multiplicity}; see
\Cref{s:vertical_tangent,s:unit_slope_tangent,s:horizontal_tangent,s:vertical_frozen_neighborhood,s:unit_slope_frozen_neighborhood,s:horizontal_frozen_neighborhood} for more details.
\end{proof}

\subsection{Derivatives of the tiling action at critical points}\label{subsec:action-derivatives}

In this subsection we record useful formulas for derivatives of the tiling action
$S(z;x,s)$ at a critical point $w_c$.

\subsubsection{Higher derivatives of the tiling action}
Fix \((x,s)\in\fP\), and assume that \((x,s)\) is neither a tangent location nor a
ramification point. Let \(w_c\) be a critical point of \(S(\,\cdot\,;x,s)\). By
\eqref{e:bcond}, the excluded cases imply
\[
f(w_c)\neq 0,-1,\infty,
\qquad
\chi(w_c)\neq 0,1,\infty .
\]
In particular, all denominators below are nonzero.

Differentiating the expression for \(S'\) and evaluating at \(w=w_c\), we obtain
\begin{align}\label{e:derSsecond}
S''(w_c;x,s)
&=
\frac{1}{w_c-x}
-\frac{1}{w_c-(x-s)}
-\frac{f'(w_c)}{f(w_c)} .
\end{align}
Using \(\chi=f/(f+1)\), we have
\[
\frac{f'}{f}
=
\frac{\chi'}{\chi(1-\chi)}.
\]
Moreover, the critical-point equation \eqref{e:critical_point} gives
\[
w_c-x=-s\,\chi(w_c),
\qquad
w_c-(x-s)=s\bigl(1-\chi(w_c)\bigr).
\]
Substituting these identities into \eqref{e:derSsecond} yields
\begin{align}\begin{split}\label{e:derSsecond-chi}
S''(w_c;x,s)
&=
\frac{1}{-s\,\chi(w_c)}
-\frac{1}{s(1-\chi(w_c))}
-\frac{\chi'(w_c)}
{\chi(w_c)(1-\chi(w_c))}  \\
&=
-\left(\chi'(w_c)+\frac{1}{s}\right)
\frac{1}{\chi(w_c)\bigl(1-\chi(w_c)\bigr)} .
\end{split}\end{align}

Further differentiating gives the higher derivatives of \(S\) at \(w_c\). For the
third derivative, one obtains
\begin{align}\label{e:derSthird0}
S'''(w_c;x,s)
&=
\frac{
\bigl((\chi'(w_c))^2-1/s^2\bigr)
\bigl(1-2\chi(w_c)\bigr)}
{\chi(w_c)^2\bigl(1-\chi(w_c)\bigr)^2}
-
\frac{\chi''(w_c)}
{\chi(w_c)\bigl(1-\chi(w_c)\bigr)} .
\end{align}
Similarly, the fourth derivative is
\begin{align}\label{e:derSfourth0}
S''''(w_c;x,s)
&=
-\frac{
2\bigl(3\chi(w_c)^2-3\chi(w_c)+1\bigr)
\bigl(\chi'(w_c)^3+1/s^3\bigr)}
{\chi(w_c)^3\bigl(1-\chi(w_c)\bigr)^3} \notag \\
&\quad
-\frac{
3\bigl(2\chi(w_c)-1\bigr)\chi'(w_c)\chi''(w_c)}
{\chi(w_c)^2\bigl(1-\chi(w_c)\bigr)^2}
-
\frac{\chi'''(w_c)}
{\chi(w_c)\bigl(1-\chi(w_c)\bigr)} .
\end{align}

Assume now that $(x,s)\in\fA$. By \eqref{e:tangent_wc}, the critical point $w_c$ is given by the
intersection of the tangent line through $(x,s)$ with the $x$-axis. Moreover, by the second statement of
\Cref{p:surface},
\begin{equation}\label{e:chi-prime-arctic}
\chi'(w_c)=-\frac{1}{s}.
\end{equation}
If $(x,s)$ is not a cusp location, then $\chi''(w_c)\neq 0$, and \eqref{e:derSthird0} simplifies to
\begin{align}\begin{split}\label{e:derSthird}
S'''(w_c;x,s)=-\frac{\chi''(w_c)}{\chi(w_c)\bigl(1-\chi(w_c)\bigr)}\neq 0.
\end{split}\end{align}
If $(x,s)$ is a cusp location, then $\chi''(w_c)=0$ but $\chi'''(w_c)\neq 0$, and
\eqref{e:derSthird0}--\eqref{e:derSfourth0} reduce to
\begin{align}\begin{split}\label{e:derSfourth}
S'''(w_c;x,s)=0,
\qquad
S''''(w_c;x,s)=
-\frac{\chi'''(w_c)}{\chi(w_c)\bigl(1-\chi(w_c)\bigr)}.
\end{split}\end{align}

\subsubsection{Derivative of the tiling action at critical point}

When the critical value \(S(w_c(x,s);x,s)\) is differentiated, the derivative
of \(w_c\) disappears because \(S'(w_c(x,s);x,s)=0\).  Direct differentiation of the tiling 
action gives
\begin{align}\begin{split}\label{e:critical_value_derivatives_gff}
 &\partial_xS(w_c;x,s)=-\ln\frac{x-w_c}{s-x+w_c}=-\ln\left(\frac{\chi(w_c)}{1-\chi(w_c}\right)=-\ln f(w_c),
\\
 &\partial_sS(w_c;x,s)=-\ln\left(\frac{s-x+w_c}{s}\right)=-\ln(1-\chi(w_c)),
 \quad
 (\partial_x+\partial_s)S(w_c;x,s)=-\ln\chi(w_c).
\end{split}\end{align}

\subsubsection{Horizontal tangency: change of coordinates}

If $(x,s)\in\fA$ is close to a horizontal tangency location $(x_0,s_0)\in\fA$, it is convenient to work
in the coordinate
\[
{\wt w}=\frac{1}{x_0-w},
\qquad\text{so that}\qquad w=x_0-\frac{1}{{\wt w}}.
\]
By the third statement in \Cref{p:surface}, we may expand
\begin{align}\label{e:chi-expand-horizontal}
\chi(w)=\frac{-w+x_0}{s_0}+\sum_{i\ge1}\frac{c_i}{(x_0-w)^i}
=: \wt\chi({\wt w})
=\frac{1}{s_0\wt w}+\sum_{i\ge1}c_i {\wt w}^i.
\end{align}
(Equivalently, $\wt\chi({\wt w})=\chi(x_0-1/{\wt w})$.) In this coordinate, the critical point of
$S'(w;x_0,s_0)$ is at $w=\infty$, i.e.\ ${\wt w}=0$.

Define the transformed tiling action
\begin{align}\label{e:change_coordinate}
\wt S({\wt w};x,s):&=S(w;x,s)=S\!\left(x_0-\frac{1}{{\wt w}};\,x,s\right)
=s\ln s-(x-x_0+1/{\wt w})\ln (x-x_0+1/{\wt w})\\
&\quad-(s-x+x_0-1/{\wt w})\ln (s-x+x_0-1/{\wt w})
-\int_0^{x_0-1/{\wt w}} \ln f(u)\,\rd u.
\end{align}
A critical point ${\wt w}_c=1/(x_0-w_c)$, where $w_c$ is a critical point of $S(\cdot;x,s)$, is characterized by
\begin{align}\label{e:critical_point2}
x_0-\frac{1}{{\wt w}_c}=x-s\,\wt \chi({\wt w}_c)
\qquad\Longleftrightarrow\qquad
s {\wt w}_c\,\wt \chi({\wt w}_c)-(x-x_0){\wt w}_c-1=0.
\end{align}

We will later use the identities 
\begin{align}\begin{split}\label{e:tchi_der}
\left.\del_{\wt w}\bigl({\wt w}\wt \chi({\wt w})\bigr)\right|_{{\wt w}={\wt w}_c}
&={\wt w}_c\wt \chi'({\wt w}_c)+\wt\chi({\wt w}_c)
=\frac{1}{{\wt w}_c}\Bigl(\chi'(w_c)+\frac{1}{s}\Bigr)+\frac{x-x_0}{s},\\
\left.\del_{\wt w}^2\bigl({\wt w}\wt \chi({\wt w})\bigr)\right|_{{\wt w}={\wt w}_c}
&={\wt w}_c\wt \chi''({\wt w}_c)+2\wt\chi'({\wt w}_c)
=\frac{1}{{\wt w}_c^3}\,\chi''(w_c),
\end{split}\end{align}
where we used $\chi'(w)={\wt w}^2\wt\chi'({\wt w})$ and $\chi''(w)=2{\wt w}^3\wt\chi'({\wt w})+{\wt w}^4\wt\chi''({\wt w})$.

\subsubsection{Derivatives of $\wt S$ at ${\wt w}_c$}

A direct chain-rule computation gives
\begin{align}\begin{split}\label{e:wS_der}
\wt S'({\wt w};x,s)
=S'(w;x,s)\,\del_{\wt w} w
=\frac{1}{{\wt w}^2}S'(w;x,s)=\frac{1}{{\wt w}^2}\left(\ln \frac{w-x}{x-s-w}-\ln \frac{\chi(w)}{1-\chi(w)}\right).
\end{split}\end{align}
Evaluating at a critical point ${\wt w}_c=1/(x_0-w_c)$, where $w_c$ is a critical point of $S(\cdot;x,s)$ and using
\eqref{e:derSsecond-chi} yields
\begin{align}\begin{split}\label{e:dertS2}
\wt S''({\wt w}_c;x,s)
&=\left.\del_{\wt w}\!\left(\frac{1}{{\wt w}^2}S'(w;x,s)\right)\right|_{{\wt w}={\wt w}_c}
=\frac{1}{{\wt w}_c^4}\,S''(w_c;x,s)\\
&=-\frac{1}{{\wt w}_c^4}\Bigl(\chi'(w_c)+\frac{1}{s}\Bigr)\frac{1}{\chi(w_c)\bigl(1-\chi(w_c)\bigr)}=-\frac{1}{{\wt w}_c^3}\,
\frac{\del_{\wt w}\bigl({\wt w}\wt \chi({\wt w})\bigr)\big|_{{\wt w}={\wt w}_c}-(x-x_0)/s}{\wt \chi({\wt w}_c)\bigl(1-\wt \chi({\wt w}_c)\bigr)}.
\end{split}\end{align}
Similarly,
\begin{align}\begin{split}\label{e:dertS3}
\wt S'''({\wt w}_c;x,s)
&=\left.\del_{\wt w}\!\left(-\frac{2}{{\wt w}^3}S'(w;x,s)+\frac{1}{{\wt w}^4}S''(w;x,s)\right)\right|_{{\wt w}={\wt w}_c}=-\frac{6}{{\wt w}_c^5}S''(w_c;x,s)+\frac{1}{{\wt w}_c^6}S'''(w_c;x,s)\\
&=\frac{6}{{\wt w}_c^5}\Bigl(\chi'(w_c)+\frac{1}{s}\Bigr)\frac{1}{\chi(w_c)\bigl(1-\chi(w_c)\bigr)}
-\frac{1}{{\wt w}_c^6}\frac{\chi''(w_c)}{\chi(w_c)\bigl(1-\chi(w_c)\bigr)}\\
&=\frac{6}{{\wt w}_c^4}\,
\frac{\del_{\wt w}\bigl({\wt w}\wt \chi({\wt w})\bigr)\big|_{{\wt w}={\wt w}_c}-(x-x_0)/s}{\wt \chi({\wt w}_c)\bigl(1-\wt \chi({\wt w}_c)\bigr)}
-\frac{1}{{\wt w}_c^3}\,
\frac{\del_{\wt w}^2\bigl({\wt w}\wt \chi({\wt w})\bigr)\big|_{{\wt w}={\wt w}_c}}{\wt \chi({\wt w}_c)\bigl(1-\wt \chi({\wt w}_c)\bigr)}.
\end{split}\end{align}

If $(x,s)\in\fA$ is not a cusp location, then $\chi'(w_c)=-1/s$. Hence, by \eqref{e:tchi_der},
\[
\del_{\wt w}\bigl({\wt w}\wt\chi({\wt w})\bigr)\big|_{{\wt w}={\wt w}_c}=\frac{x-x_0}{s}.
\]
Substituting this into \eqref{e:dertS2} gives $\wt S''({\wt w}_c;x,s)=0$, and the first term in \eqref{e:dertS3} vanishes as well. Therefore,
\begin{align}\begin{split}\label{e:dertS31}
\wt S''({\wt w}_c;x,s)=0,
\qquad
\wt S'''({\wt w}_c;x,s)
=-\frac{1}{{\wt w}_c^3}\,
\frac{\del_{\wt w}^2\bigl({\wt w}\wt \chi({\wt w})\bigr)\big|_{{\wt w}={\wt w}_c}}{\wt \chi({\wt w}_c)\bigl(1-\wt \chi({\wt w}_c)\bigr)}.
\end{split}\end{align}

\subsection{Change of time}
\label{s:changetime}
In this section we show that the critical value of the tiling action
\begin{align}\label{e:def_action_copy}
S(z;x,s)
&=
s\ln s
-(x-z)\ln (x-z)
-(s-x+z)\ln (s-x+z)
-\int_{0}^{z}\ln f(u)\,\rd u
\end{align}
is invariant under a vertical shift of time.

Fix \(\ft>0\), and shift the polygon upward by \(\ft\). Recall the shifted polygon
\(\fP_{\ft}\) from \eqref{e:shifted_polygon} and the associated Riemann surface
\(\cC_{\ft}\) from \eqref{e:shifted_RS}. Also recall from \eqref{e:emb_ft} the
embedding
\begin{align}\label{e:emb_ft_copy}
\Psi_\ft:\ (x,s)\in\overline{\fL}
\longmapsto
\bigl(f(x,s),\,x-(s+\ft)\chi(x,s)\bigr)\in\cC_\ft .
\end{align}
We also recall from \eqref{e:CtoCt}, the original curve \(\cC\) and the shifted curve \(\cC_\ft\) are biholomorphic via
the transport map
\begin{align}\label{e:maptoCt}
(f,z)\in\cC
\longmapsto
\biggl(f,\,z-\ft\frac{f}{f+1}\biggr)
=
\bigl(f,\,z-\ft\chi\bigr)
\in\cC_\ft .
\end{align}
The basepoint \((f,z)=(0,b_1)=(0,0)\in\cC(\bR)\) is mapped to the
basepoint \((0,0)\in\cC_\ft\). When no confusion can arise, we identify a point
\((f,w)\in\cC_\ft\) with its \(w\)-coordinate and write \(f=f_\ft(w)\).

Analogously to \eqref{e:def_action_copy}, define the shifted tiling action on \(\cC_\ft\)
by
\begin{align}\label{e:def_action_shift}
S_\ft(w;x,s+\ft)
&:=
(s+\ft)\ln (s+\ft)
-(x-w)\ln (x-w) \notag \\
&\quad
-(s+\ft-x+w)\ln (s+\ft-x+w)
-\int_0^w \ln f_\ft(u)\,\rd u,
\end{align}
where \(0\) denotes the basepoint \((0,0)\in\cC_\ft\).

The following lemma states that the critical value of the tiling action
\eqref{e:def_action_copy} is invariant under a vertical shift of time.

\begin{lemma}[Time-translation invariance]\label{c:change_time}
For any point \((x,s)\in\fL\), set
\[
w_c:=x-s\chi(x,s)\in\cC,
\qquad
w_{c,\ft}:=x-(s+\ft)\chi(x,s)=w_c-\ft\chi(x,s)\in\cC_\ft .
\]
Then \(w_c\) is a critical point of \(S(\,\cdot\,;x,s)\), \(w_{c,\ft}\) is a critical
point of \(S_\ft(\,\cdot\,;x,s+\ft)\), and
\[
S(w_c;x,s)
=
S_\ft(w_{c,\ft};x,s+\ft).
\]
\end{lemma}

\begin{proof}[Proof of \Cref{c:change_time}]
We identify a point \((f,z)\in\cC\)
with its \(z\)-coordinate and write \(f=f(z)\), and a point \((f,w)\in\cC_\ft\)
with its \(w\)-coordinate and write \(f=f_\ft(w)\). 
We denote \(\chi(z)=f(z)/(f(z)+1)\), and the transport map \eqref{e:maptoCt} is given by
\begin{align}\label{e:maptoCt_copy}
z\in \cC\mapsto z-\ft \chi(z) \in \cC_\ft, 
\qquad 
f(z)=f_\ft(z-\ft \chi(z)).
\end{align}

Then we have \(x-s\chi(w_c)=w_c\) and \(x-(s+\ft)\chi(w_c)=w_{c,\ft}\), and we can simplify 
\(S(w_c;x,s)\) as
\begin{align}\begin{split}\label{e:Sexp}
S(w_c;x,s)
&=s\ln s-(x-w_c)\ln (x-w_c)-(s-x+w_c)\ln (s-x+w_c)-\int_0^{w_c}\ln f(u)\rd u\\
&=s\ln s-(s\chi(w_c))\ln (s\chi(w_c))-(s(1-\chi(w_c)))\ln (s(1-\chi(w_c)))-\int_0^{w_c}\ln f(u)\rd u\\
&=-s \left(\chi(w_c)\ln \chi(w_c)+(1-\chi(w_c))\ln (1-\chi(w_c))\right)-\int_0^{w_c}\ln f(u)\rd u,
\end{split}\end{align}
where the second line follows from \(x-w_c=s\chi(w_c)\), and the third line follows from cancelling terms. 
By the same argument, we have
\begin{align}\begin{split}\label{e:Binterm}
&\phantom{{}={}}(s+\ft)\ln (s+\ft)-(x-w_{c,\ft})\ln (x-w_{c,\ft})-(s+\ft-x+w_{c,\ft})\ln (s+\ft-x+w_{c,\ft})\\
&=(s+\ft)\ln (s+\ft)-((s+\ft)\chi(w_c))\ln ((s+\ft)\chi(w_c))
-((s+\ft)(1-\chi(w_c)))\ln ((s+\ft)(1-\chi(w_c)))\\
&=-(s+\ft) \left(\chi(w_c)\ln \chi(w_c)+(1-\chi(w_c))\ln (1-\chi(w_c))\right).
\end{split}\end{align}

By \eqref{e:maptoCt_copy}, \(u\mapsto u-\ft\chi(u)\) maps the path from \(0\) to \(w_c\) in \(\cC\) onto a path from \(0\) to \(w_{c,\ft}\) in \(\cC_\ft\), and \(f(u)=f_\ft(u-\ft\chi(u))\) along this path. We can perform a change of variable:
\begin{align}\begin{split}\label{e:intterm1}
\int_0^{w_{c,\ft}} \ln f_\ft(u)\rd u
&=\int_0^{w_c} \ln f_\ft(u-\ft \chi(u))\rd (u-\ft \chi(u))
=\int_0^{w_c} \ln f(u)\rd (u-\ft \chi(u))\\
&=\int_0^{w_c} \ln f(u)\rd u-\ft \int_0^{w_c} \ln f(u)\rd \chi(u).
\end{split}\end{align}

For the second term, using integration by parts and \(f(u)=\chi(u)/(1-\chi(u))\),
\begin{align}\begin{split}\label{e:intterm2}
\int_0^{w_c} \ln f(u)\,\rd \chi(u)
&=\Bigl[\chi(u)\,\ln f(u)\Bigr]_{0}^{w_c}\;-\;\int_0^{w_c} \chi(u)\,\rd\!\bigl(\ln f(u)\bigr) \\
&=\Bigl[\chi(u)\,\ln f(u)\Bigr]_{0}^{w_c}\;-\;\int_0^{w_c} \frac{f(u)}{1+f(u)}\,\frac{\rd f(u)}{f(u)} \\
&=\Bigl[\chi(u)\,\ln f(u)\Bigr]_{0}^{w_c}\;-\;\Bigl[\ln\bigl(1+f(u)\bigr)\Bigr]_{0}^{w_c}.
\end{split}\end{align}
Since \(f=\chi/(1-\chi)\), the boundary expression can be written as
\[
  \chi\,\ln f - \ln(1+f)
  \;=\; \chi\ln\!\frac{\chi}{1-\chi}\;-\;\ln\!\frac{1}{1-\chi}
  \;=\; \chi\ln\chi + (1-\chi)\ln(1-\chi).
\]
With the choice of the basepoint \(u=0\) so that \(\chi(0)=0\) and \(f(0)=0\), hence
\[
\left.\chi\ln\chi + (1-\chi)\ln(1-\chi)\right|_{u=0}=0,
\]
\eqref{e:intterm2} yields
\begin{equation}\label{e:intterm2_final}
  \int_0^{w_c} \ln f(u)\,\rd \chi(u)
  =
  \chi(w_c)\ln \chi(w_c)+(1-\chi(w_c))\ln (1-\chi(w_c)).
\end{equation}

Combining \eqref{e:Binterm}, \eqref{e:intterm1}, and \eqref{e:intterm2_final}, we conclude that
\begin{align*}
S_\ft(w_{c,\ft};x,s+\ft)
&=(s+\ft)\ln (s+\ft)-(x-w_{c,\ft})\ln (x-w_{c,\ft})\\
&\quad -(s+\ft-x+w_{c,\ft})\ln (s+\ft-x+w_{c,\ft}) -\int_0^{w_{c,\ft}} \ln f_\ft(u)\rd u\\
&=-(s+\ft)\left(\chi(w_c)\ln \chi(w_c)+(1-\chi(w_c))\ln (1-\chi(w_c))\right)\\
&\quad -\left[
\int_0^{w_c}\ln f(u)\rd u
-\ft\left(\chi(w_c)\ln \chi(w_c)+(1-\chi(w_c))\ln (1-\chi(w_c))\right)
\right]\\
&=-s \left(\chi(w_c)\ln \chi(w_c)+(1-\chi(w_c))\ln (1-\chi(w_c))\right)
-\int_0^{w_c}\ln f(u)\rd u=S(w_c;x,s).
\end{align*}
\end{proof}

\section{Descent and Ascent Critical Points}
\label{s:prop_critical}

In this section, we assign to each point $(x,s)\in \fP$ a subset of its associated critical points. We then introduce the notions of descent and ascent critical points, together with several related concepts.

\subsection{Associated critical points}
\label{s:dacritical}

Allowing tangencies to all of \(\fA\) would generally produce more formal
critical points than are needed later. We therefore specify the relevant
subset of tangent lines.

Recall from \Cref{t:frozen_structure} that  the
frozen region is decomposed into curvilinear triangles. Each such triangle
is bounded by two segments joining a vertex of \(\fP\) to two tangency
points of \(\fA\), together with a piece of the arctic curve \(\fA\). The
curvilinear triangles have pairwise disjoint interiors and cover
\(\fP\setminus\fL\).

\begin{proposition}\label{p:associate_critical_points}
For each point \((x,s)\in\fP\), we associate formal critical points as
follows.

\begin{enumerate}
\item
\textbf{Liquid region.}
If \((x,s)\in\fL\) we
associate with \((x,s)\) the two complex-conjugate critical points in
\eqref{e:two_points}.

\item
\textbf{A unique curvilinear triangle.}
Suppose that \((x,s)\in\fP\setminus\fL\) belongs to a unique
curvilinear triangle \(\fT\). We associate with \((x,s)\) the formal
critical points corresponding to tangent lines from \((x,s)\) whose
points of tangency lie on the portion of the arctic boundary contained in
\(\fT\).

\item
\textbf{A boundary shared by two curvilinear triangles.}
Suppose that \((x,s)\) lies on a boundary shared by two adjacent
curvilinear triangles \(\fT_A\) and \(\fT_B\). We associate with
\((x,s)\) the formal critical points corresponding to tangent lines from
\((x,s)\) whose points of tangency lie on the portion of the arctic
boundary contained in
$
\fT_A\cup\fT_B.
$
\end{enumerate}

In each of the last two cases, the total number of associated formal
critical points, counted with formal multiplicity, is two or three.
\end{proposition}

\begin{figure}
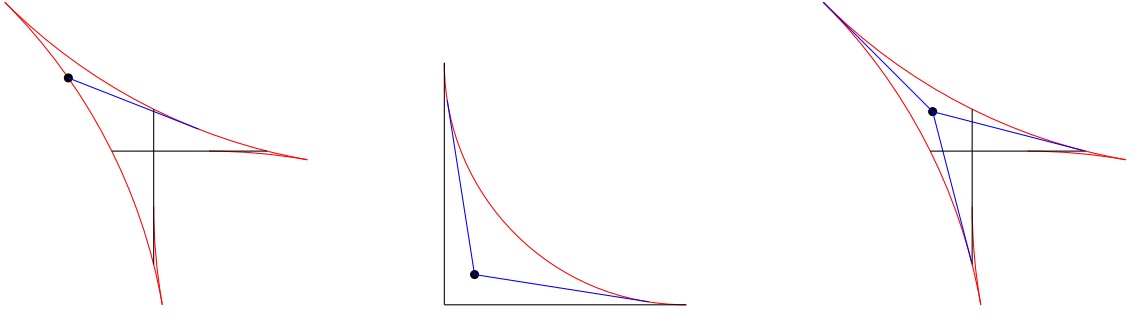

  \centering

 \begin{subfigure}{0.32\textwidth}
    \centering

  \end{subfigure}

  \caption{Examples of critical points and tangent lines.}
  \label{f:example_tangent_line}
\end{figure}

The possible multiplicities in the last statement of
\Cref{p:associate_critical_points} can be summarized as follows:
\begin{enumerate}
\item
At a cusp point that is not a cusp-turning point, the distinguished tangent
line gives one formal critical point of multiplicity three.

\item
At a regular arctic point that is not a tangent location, the tangent line
at the point gives one formal critical point of multiplicity two. There may
also be one additional simple formal critical point; see the left panel of
\Cref{f:example_tangent_line}.

\item
At a point in the interior of the frozen region that does not lie on a
shared boundary, there are two or three simple formal critical points; see
the middle and right panels of \Cref{f:example_tangent_line}.

\item
If two adjacent curvilinear triangles meet only at a tangency point (see the left panel of
\Cref{f:adjacent_curvilinear_triangle}), then
the tangent line at that point gives a formal critical point of multiplicity
two.

\item
At the tangency endpoint of a nontrivial shared segment (see the middle panel of
\Cref{f:adjacent_curvilinear_triangle}), the tangent line at
the endpoint gives a formal critical point of multiplicity two, and there is
one additional simple formal critical point.

\item
At a cusp-turning endpoint of a shared segment (see the right panel of
\Cref{f:adjacent_curvilinear_triangle}), the distinguished tangent
line gives one formal critical point of multiplicity three.

\item
At an interior point of a nontrivial shared segment (see the middle and right panel of
\Cref{f:adjacent_curvilinear_triangle}), there are three
distinct formal critical points. One of them corresponds to the supporting
line of the shared segment and is spurious, while the other two correspond
to tangencies on the two adjacent portions of the arctic boundary.
\end{enumerate}
\begin{proof}[Proof of \Cref{p:associate_critical_points}]
The assertions for the liquid region follow from \eqref{e:two_points}. It remains to count the tangent lines in the frozen region. By symmetry, it
suffices to consider a curvilinear triangle \(\fT\) on which
\[
\nabla H^*=(1,0).
\]
We refer to \Cref{f:curvilinear_triangle} for the complete list of such
curvilinear triangles.

As \((x,s)\) varies within any of the open regions determined by the
extended sides, the number of tangent lines from \((x,s)\) to the portion
of \(\fA\) contained in \(\fT\), counted with multiplicity, is locally
constant. When \((x,s)\) crosses an extended side, exactly one tangent line
enters or leaves the portion of the arctic boundary contained in \(\fT\).
Inspecting the four configurations in \Cref{f:curvilinear_triangle} gives
\begin{align}\label{e:formal_critical_count}
\begin{array}{c|c|c|c|c}
\text{Open quadrant}
&\text{First}&\text{Second}&\text{Third}&\text{Fourth}\\
\hline
\text{Number of formal critical points}
&2&3&2&3.
\end{array}
\end{align}
On the arctic boundary, the corresponding tangent lines coalesce with the
multiplicities described in \Cref{p:critical_multiplicity}.

The same argument, applied to the union of two adjacent curvilinear
triangles, gives the shared-boundary cases. In the relative interior of a
shared extended side, the supporting line contributes one spurious formal
critical point, while each of the two adjacent portions of the arctic
boundary contributes one additional tangent line. Thus, there are three
formal critical points. At a tangency or cusp-turning endpoint, some of
these tangent lines coalesce, giving the stated multiplicities.
\end{proof}

\begin{remark}\label{r:associated_critical_arc}
Suppose that \((x,s)\) is contained in a curvilinear triangle \(\fT\).
By \eqref{e:arcCR}, every formal critical point associated with \((x,s)\)
through \(\fT\) lies in
\[
\begin{cases}
[b_i,\infty_i],
& \nabla H^*=(0,0)\text{ on }\fT,\\
[\infty_i,a_i],
& \nabla H^*=(1,0)\text{ on }\fT,\\
[a_i,b_{i+1}],
& \nabla H^*=(1,-1)\text{ on }\fT.
\end{cases}
\]
The arcs are closed because a formal critical point may coincide with a
tangent location. If \((x,s)\) lies on a boundary shared by
\(\fT_A\) and \(\fT_B\), then the associated formal critical points lie
in the union of the corresponding two arcs.
\end{remark}

\subsection{Descent and ascent critical points}
\label{s:descent_ascent_critical_points}

We now classify the genuine critical points according to their local
steepest-descent and steepest-ascent geometry.

\begin{definition}[Descent and ascent critical points]
\label{d:descent_ascent_critical}
A genuine critical point \(w_c\) associated with \((x,s)\) is called a
\emph{descent critical point} if at least one nonreal local
steepest-descent path emanates from \(w_c\). It is called an
\emph{ascent critical point} if at least one nonreal local steepest-ascent
path emanates from \(w_c\).

These two classes are not disjoint: a genuine critical point may be both a
descent and an ascent critical point. In particular, each of the two
complex-conjugate critical points associated with a liquid point is both a
descent and an ascent critical point.
\end{definition}

Suppose that \((x,s)\notin\fA\), and let \(w_c\in\cC(\bR)\) be a simple
genuine critical point associated with a tangent line from \((x,s)\) to
\(\fA\) at \((x',s')\). Assume that the slope of the tangent line does
not belong to \(\{0,1,\infty\}\). The tangent line is
\begin{align}\label{e:tangent_line0}
L(w_c)
:=
\{(u,v):u-v\chi(w_c)=w_c\},
\end{align}
and
\[
w_c=x'-s'\chi(w_c)=x-s\chi(w_c).
\]
We recall the second derivative of the tiling action from \eqref{e:derSsecond-chi}, and using
\(\chi'(w_c)=-1/s'\), we obtain
\begin{align}\label{e:critical_second_derivative}
S''(w_c;x,s)
=
-\left(\chi'(w_c)+\frac1s\right)
\frac{1}{\chi(w_c)(1-\chi(w_c))}
=
\frac{s-s'}
{ss'\chi(w_c)(1-\chi(w_c))}.
\end{align}

For a simple real critical point, the local quadratic expansion gives
\begin{align}\label{e:descent_ascent_sign}
\begin{split}
w_c\text{ is a descent critical point}
&\quad\Longleftrightarrow\quad
S''(w_c;x,s)>0,\\
w_c\text{ is an ascent critical point}
&\quad\Longleftrightarrow\quad
S''(w_c;x,s)<0.
\end{split}
\end{align}
Indeed, if \(S''(w_c;x,s)>0\), then two nonreal steepest-descent paths and
two real steepest-ascent paths emanate from \(w_c\). If
\(S''(w_c;x,s)<0\), the roles of ascent and descent are reversed.

Since \(s,s'>0\), the sign in
\eqref{e:critical_second_derivative} is determined by \(s-s'\) and
\(\chi(w_c)(1-\chi(w_c))\). The tangency point \((x',s')\) divides
\(L(w_c)\) into two components, and the type of \(w_c\) remains constant
as \((x,s)\) varies within either component. The resulting geometric
criterion is
\begin{align}\label{e:geometric_descent_ascent}
\begin{array}{c|c|c|c}
\nabla H^*\text{ on }\fT
&
\operatorname{slope}L(w_c)
&
\text{Descent side}
&
\text{Ascent side}
\\
\hline
(1,0)
&(-\infty,0)
&\text{below/right}
&\text{above/left}
\\
(1,-1)
&(0,1)
&\text{below/left}
&\text{above/right}
\\
(0,0)
&(1,\infty)
&\text{above/right}
&\text{below/left}.
\end{array}
\end{align}
Here ``above,'' ``below,'' ``left,'' and ``right'' describe the position
of \((x,s)\) relative to the tangency point \((x',s')\) along
\(L(w_c)\). The two descriptions in each entry are equivalent because
\((x,s)\) and \((x',s')\) lie on the same tangent line; see \Cref{f:critical}.

The degenerate and exceptional cases are summarized in the following
proposition.

\begin{proposition}\label{p:descent_ascent_degenerate}
The following statements hold.

\begin{enumerate}
\item
At a regular arctic point that is not a tangent location, the corresponding
critical point has multiplicity two and is both a descent and an ascent
critical point.

\item
At a cusp point that is not a cusp-turning point, the corresponding
critical point has multiplicity three and is both a descent and an ascent
critical point.

\item
At a regular point of an extended side that is neither a tangency point nor
a cusp-turning point, the corresponding formal critical point is spurious
and is therefore neither a descent nor an ascent critical point.

\item
At a tangency point on an extended side, logarithmic cancellation leaves a
simple genuine critical point. It is a descent critical point in each of
the following cases:
\begin{enumerate}
\item
the extended side is vertical and the arctic boundary is locally tangent to
it from the left;

\item
the extended side has unit slope and the arctic boundary is locally tangent
to it from the right;

\item
the extended side is horizontal and the arctic boundary is locally tangent
to it from above.
\end{enumerate}
In the three opposite configurations, it is an ascent critical point.

\item
At a cusp-turning point, logarithmic cancellation leaves a genuine critical
point of multiplicity two. This point is both a descent and an ascent
critical point.
\end{enumerate}
\end{proposition}

\begin{proof}
The first two statements follow from the local cubic and quartic expansions
of the tiling action at regular arctic and cusp points. The last three
statements follow from \Cref{p:critical_multiplicity} and the local expansions of the tiling action in 
\Cref{s:vertical_tangent,s:unit_slope_tangent,s:horizontal_tangent,s:vertical_frozen_neighborhood,s:unit_slope_frozen_neighborhood,s:horizontal_frozen_neighborhood}.
\end{proof}

Suppose that all associated real critical points are distinct and
nondegenerate. When ordered along the relevant real arc, the signs of
\(S''\) alternate. Consequently, two associated critical points consist of
one descent and one ascent critical point, whereas three associated
critical points have one of the two patterns
\[
\text{ascent--descent--ascent}
\qquad\text{or}\qquad
\text{descent--ascent--descent}.
\]

\begin{figure}
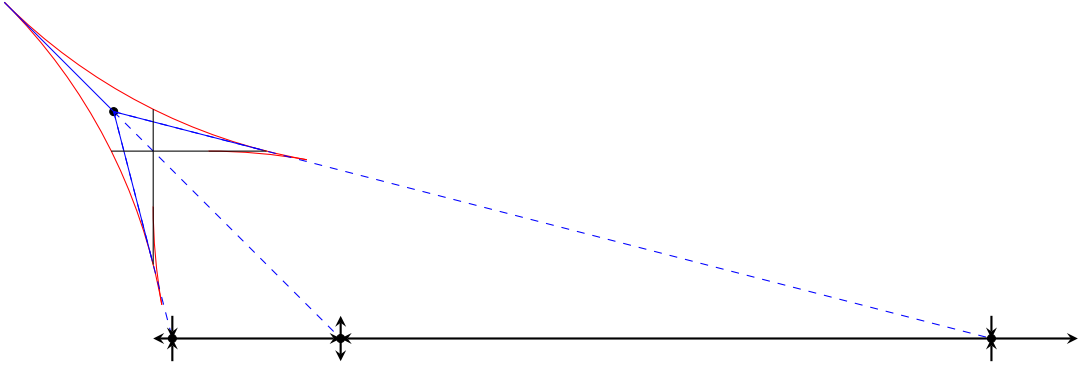

    \centering
      %

  \caption{Critical point pattern}
  \label{f:critical_point_pattern}
\end{figure}

For example, in \Cref{f:critical_point_pattern}, the three critical points
\(w_{c,1},w_{c,2},w_{c,3}\), ordered from left to right, satisfy
\[
S''(w_{c,1};x,s)<0,\qquad
S''(w_{c,2};x,s)>0,\qquad
S''(w_{c,3};x,s)<0.
\]
They therefore have the ascent--descent--ascent pattern.

For \(\nabla H^*=(1,0)\), the numbers of descent and ascent critical
points at a point in each open quadrant relative to the
corresponding corner are
\begin{align}\label{e:descent_ascent_count}
\begin{array}{c|c|c|c}
\text{Open quadrant}
&
\text{Total}
&
\text{Descent}
&
\text{Ascent}
\\
\hline
\text{First}
&2&1&1\\
\text{Second}
&3&1&2\\
\text{Third}
&2&1&1\\
\text{Fourth}
&3&2&1
\end{array}
\end{align}
The cases \(\nabla H^*=(1,-1)\) and
\(\nabla H^*=(0,0)\) follow from the symmetries in
\Cref{f:symmetry}.

Finally, we introduce notation for the sets of descent and ascent critical
points.

\begin{definition}\label{def:ascent_descent_critical}
For each point \((x,s)\in\fP\), let
$
\operatorname{Crit}^{\rm d}(x,s)
$ and $
\operatorname{Crit}^{\rm a}(x,s)
$
denote, respectively, the sets of descent and ascent critical points
associated with \((x,s)\) as in \Cref{p:associate_critical_points}.
\end{definition}

\subsection{Descent and ascent cuts}
\label{s:descent_ascent_cuts}

Fix a curvilinear triangle \(\fT\), let
\((x',s')\in\fA\cap\fT\), and let \(w_c\) be the formal critical point
corresponding to the tangent line to \(\fA\) at \((x',s')\). Thus,
$
w_c=x'-s'\chi(w_c),
$
and the tangent line is
\begin{align}\label{e:tangent_line}
L(w_c)
:=
\{(x,s):w_c=x-s\chi(w_c)\}.
\end{align}
When \(\chi(w_c)=\infty\), this equation is understood projectively.

By \Cref{c:intersection}, \(L(w_c)\cap\fT\) is a single closed interval.
We first define the cuts when the slope of \(L(w_c)\) does not belong to
\(\{0,1,\infty\}\).

\begin{definition}\label{d:nonexceptional_cuts}
Suppose that the slope of \(L(w_c)\) in \eqref{e:tangent_line} does not
belong to \(\{0,1,\infty\}\). We define the \emph{descent cut} and the
\emph{ascent cut} by
\begin{align}\label{e:l0}
\begin{split}
\ell_-(w_c;\fT)
&:=
\left\{
(x,s)\in L(w_c)\cap\fT:
S''(w_c;x,s)>0
\right\}
\cup\{(x',s')\},
\\
\ell_+(w_c;\fT)
&:=
\left\{
(x,s)\in L(w_c)\cap\fT:
S''(w_c;x,s)<0
\right\}
\cup\{(x',s')\}.
\end{split}
\end{align}
Thus, if \((x,s)\in\ell_-(w_c;\fT)\), then \(w_c\) is a descent
critical point associated with \((x,s)\); if
\((x,s)\in\ell_+(w_c;\fT)\), then \(w_c\) is an ascent critical point
associated with \((x,s)\). The tangency point
\((x',s')\) belongs to both cuts, and $w_c$ is both a descent and an ascent critical point.  
\end{definition}

At a regular arctic point, the removal of \((x',s')\) separates
\(L(w_c)\cap\fT\) into two components, one forming the descent cut and the
other forming the ascent cut; see Panel (A) of \Cref{f:l}. Their positions are given by
\eqref{e:geometric_descent_ascent}. 

At a cusp point, \((x',s')\) is an endpoint of
\(L(w_c)\cap\fT\), and the sign of \(S''(w_c;x,s)\) is constant away
from \((x',s')\). Hence one of the two cuts is reduced to the cusp point,
while the other is all of \(L(w_c)\cap\fT\); see Panel (B) of \Cref{f:l}.  The same geometric
description \eqref{e:geometric_descent_ascent} remains valid.

We next define the cuts at the exceptional slopes $\{0,1,\infty\}$. In this case,
\((x',s')\) is a tangency point or a cusp-turning point, and
\(L(w_c)\) is the supporting line of an extended side. The line is a
common boundary of two adjacent curvilinear triangles \(\fT_A\) and
\(\fT_B\). We fix one of these triangles,
\[
\fT\in\{\fT_A,\fT_B\}.
\]

\begin{definition}\label{d:exceptional_cuts}
Suppose that the slope of \(L(w_c)\) belongs to
\(\{0,1,\infty\}\). Let
\[
(x,s)\in
\bigl(L(w_c)\cap\fT\bigr)\setminus\{(x',s')\}.
\]
Choose a sufficiently small perturbation
\((\widehat x,\widehat s)\) of \((x,s)\) into the side of \(L(w_c)\)
occupied by \(\fT\), so that the corresponding nearby tangency point
\((\widehat x',\widehat s')\) lies on the portion of the arctic boundary
contained in \(\fT\). Let \(\widehat w_c\) be the corresponding nearby
genuine critical point. Then from \eqref{e:critical_second_derivative}
\begin{align}\label{e:frozen_S''3}
S''(\widehat w_c;\widehat x,\widehat s)
=
\frac{\widehat s-\widehat s'}
{\widehat s\,\widehat s'\,
\chi(\widehat w_c)(1-\chi(\widehat w_c))}
\end{align}
is well defined.

We place \((x,s)\) in \(\ell_-(w_c;\fT)\) if the quantity in
\eqref{e:frozen_S''3} is positive for all sufficiently small admissible
perturbations, and in \(\ell_+(w_c;\fT)\) if it is negative. By convention,
the tangency point \((x',s')\) belongs to both cuts. See Panel~(C) of
\Cref{f:l} for the tangent-location case and Panel~(D) for the cusp-turning
case.
\end{definition}

The following lemma collects the basic properties of these definitions.

\begin{lemma}\label{l:descent_ascent_cuts}
The definition in \Cref{d:exceptional_cuts} is independent of the choice of
all sufficiently small admissible perturbations. Moreover, in both the
exceptional and nonexceptional cases,
\begin{align}\label{e:cut_partition}
\ell_-(w_c;\fT)\cup\ell_+(w_c;\fT)
&=
L(w_c)\cap\fT,
\nonumber\\
\ell_-(w_c;\fT)\cap\ell_+(w_c;\fT)
&=
\{(x',s')\}.
\end{align}
In particular, \(\ell_-(w_c;\fT)\) and
\(\ell_+(w_c;\fT)\) are closed subintervals of
\(L(w_c)\cap\fT\), possibly reduced to the tangency point.

The side descriptions in \eqref{e:geometric_descent_ascent} remain valid at
tangency and cusp-turning points by one-sided continuation within \(\fT\).

Finally, suppose that \(L(w_c)\) contains a common boundary portion of two
adjacent curvilinear triangles \(\fT_A\) and \(\fT_B\), and set
\[
B:=L(w_c)\cap\fT_A\cap\fT_B.
\]
Then
\begin{align}\begin{split}\label{e:cut_swap}
\ell_-(w_c;\fT_A)\cap B
&=
\ell_+(w_c;\fT_B)\cap B,
\\
\ell_+(w_c;\fT_A)\cap B
&=
\ell_-(w_c;\fT_B)\cap B.
\end{split}\end{align}
\end{lemma}

\begin{proof}
For a nonexceptional slope, the assertions follow directly from
\eqref{e:critical_second_derivative}. Indeed, the tangency point
\((x',s')\) is the only point at which \(S''(w_c;x,s)=0\), and the sign
is constant on each component of
\[
\bigl(L(w_c)\cap\fT\bigr)\setminus\{(x',s')\}.
\]

At an exceptional slope, an admissible perturbation remains on a fixed side
of \(L(w_c)\), and the nearby tangency point remains on the fixed portion
of the arctic boundary contained in \(\fT\). Consequently,
\(\widehat s-\widehat s'\) has a fixed sign for all sufficiently small
admissible perturbations. The sign of
\[
\chi(\widehat w_c)(1-\chi(\widehat w_c))
\]
is also fixed by the curvilinear triangle \(\fT\). Hence the sign in
\eqref{e:frozen_S''3} is independent of the perturbation. This proves that
\Cref{d:exceptional_cuts} is well defined.

The partition properties in \eqref{e:cut_partition} follow because every
point other than \((x',s')\) has exactly one of the two signs. Closedness
follows from the inclusion of the tangency point in both cuts and from the
one-sided continuity of the endpoints.

If the tangent line is shared by \(\fT_A\) and \(\fT_B\), the admissible
perturbations associated with the two triangles approach \(w_c\) from
opposite sides. The corresponding signs in
\eqref{e:frozen_S''3} are therefore opposite. This proves
\eqref{e:cut_swap}.
\end{proof}

\begin{figure}
	\begin{subfigure}{0.2\textwidth}
	
		\centering

			\caption{}
		\end{subfigure}		
		\begin{subfigure}{0.2\textwidth}
    \begin{center}
      %
      \end{center}
      	\caption{}
  \end{subfigure}
%	\begin{subfigure}{0.24\textwidth}
%    \begin{center}
%      %
  
      	\caption{}    
  \end{subfigure}
    \begin{subfigure}{0.24\textwidth}
    \centering
      %
      	\caption{}
  \end{subfigure}

\caption{Four cases for the descent cut associated with a curvilinear
triangle satisfying \(\nabla H^*=(1,0)\). In each case, the descent cut is
the possibly degenerate subinterval of \(L(w_c)\cap\fT\) lying below, or
equivalently to the right of, the tangency point \((x',s')\).}
\label{f:l}
	
\end{figure}

\subsection{Cell assignments}
\label{s:assign_cell}
As recalled in \Cref{s:Kasteleyn}, a lattice point
$
(x,s)\in\bZ^2/n\cap\fP
$
represents either a blue or a white triangle. We associate with each such
lattice point a \emph{cell}, which is either the liquid region \(\fL\) or one of
the curvilinear triangles. We denote the cell assignments for blue and white
triangles by
$
\operatorname{Cell}^{\rb}(x,s)
$ and $
\operatorname{Cell}^{\rw}(x,s),
$
respectively.

If the lattice point lies in the liquid region, we set
\begin{align}\label{e:cell_liquid}
\operatorname{Cell}^{\rb}(x,s)
=
\operatorname{Cell}^{\rw}(x,s)
:=
\fL.
\end{align}
If a lattice point in \(\fP\setminus\fL\) belongs to a unique curvilinear
triangle \(\fT\), we set
\[
\operatorname{Cell}^{\rb}(x,s)
=
\operatorname{Cell}^{\rw}(x,s)
:=
\fT.
\]
By \Cref{p}, every vertical extended side is supported on a line whose
\(x\)-coordinate lies in
$
\bZ'/n=(\bZ+1/2)/n,
$
whereas every unit-slope extended side is supported on a line whose
\(x-s\)-coordinate lies in \(\bZ'/n\). On the other hand, a lattice point
\((x,s)\in\bZ^2/n\) satisfies
$
x,\ x-s\in\bZ/n.
$
Thus, no lattice point can lie on the supporting line of a vertical or
unit-slope extended side. Consequently, an ambiguity can arise only when a
lattice point lies on a horizontal boundary shared by two adjacent
curvilinear triangles.

\begin{definition}\label{def:Cell}
Suppose that \((x,s)\in\bZ^2/n\cap\fP\) lies on a horizontal boundary shared
by two adjacent curvilinear triangles, and let \(L(w_0)\) be its supporting
line. The line \(L(w_0)\) is tangent to the arctic boundary at a point
\(\zeta\). There are two cases.

\begin{enumerate}
\item\label{i:point}
The two curvilinear triangles meet only at the point \(\zeta\); see Panels
(A) and (B) of \Cref{f:c_triangle}. In this case, \((x,s)=\zeta\), and we set
\begin{align}\label{e:cell_point}
\operatorname{Cell}^{\rb}(\zeta)
=
\operatorname{Cell}^{\rw}(\zeta)
:=
\text{the curvilinear triangle on which }
\nabla H^*=(1,0).
\end{align}

\item\label{i:segment}
The shared boundary is a nondegenerate horizontal segment. In this case, let
\(\fT^\uparrow\) and \(\fT^\downarrow\) denote the adjacent curvilinear
triangles lying above and below the segment, respectively; see Panels (C)
and (D) of \Cref{f:c_triangle}. We set
\begin{align}\label{e:cell_segment}
\operatorname{Cell}^{\rb}(x,s)
:=
\fT^\uparrow,
\qquad
\operatorname{Cell}^{\rw}(x,s)
:=
\fT^\downarrow.
\end{align}
\end{enumerate}
\end{definition}

With these conventions,
\(\operatorname{Cell}^{\rb}(x,s)\) and
\(\operatorname{Cell}^{\rw}(x,s)\) are uniquely defined for every blue and
white lattice point in \(\fP\), respectively. Moreover, if the lattice point
lies in \(\fP\setminus\fL\), its assigned cell is a curvilinear triangle
containing the lattice point and a portion of the corresponding microscopic
triangle.

We recall from \Cref{p:associate_critical_points} that every point of
\(\fP\) is associated with one or more critical points. We first describe the
critical points associated with a shared horizontal boundary.

\begin{lemma}\label{l:describe_critical_point}
Suppose that \((x,s)\in\fP\) lies on a horizontal boundary shared by two
adjacent curvilinear triangles, and let \(L(w_0)\) be its supporting line.

\begin{enumerate}
\item
In the setting of \Cref{i:point} in \Cref{def:Cell}, one has
\((x,s)=\zeta\). The point \(w_0\) is a descent critical point if the arctic
boundary is tangent to \(L(w_0)\) from above, and an ascent critical point if
the arctic boundary is tangent to \(L(w_0)\) from below.

\item
In the setting of \Cref{i:segment} in \Cref{def:Cell}, suppose that
\((x,s)\) lies in the relative interior of the shared horizontal segment.
Then \(w_0\) is a spurious critical point, and each curvilinear triangle
contributes one additional critical point. We denote these points by
\[
\xi_c^{\rm d}
\quad\text{for }\fT^\uparrow,
\qquad
\xi_c^{\rm a}
\quad\text{for }\fT^\downarrow.
\]
At a tangency point or a cusp-turning endpoint, some of the points
\(\xi_c^{\rm d},w_0,\xi_c^{\rm a}\) may coalesce.
\end{enumerate}
\end{lemma}

\begin{proof}
The descent/ascent classification follows from
\Cref{p:descent_ascent_degenerate}. In the relative interior of a horizontal
extended side, \(w_0\) is spurious, and the types of the two remaining
critical points follow from \eqref{e:geometric_descent_ascent}. The endpoint
cases follow by continuity as the corresponding tangent lines coalesce.
\end{proof}

\begin{lemma}\label{l:number_critical_point}
The following statements hold.

\begin{enumerate}
\item
If \((x,s)\in\fL\cap\bZ^2/n\), then
\[
\operatorname{Cell}^{\rb}(x,s)
=
\operatorname{Cell}^{\rw}(x,s)
=
\fL,
\]
and \((x,s)\) is associated with a pair of complex-conjugate critical
points.

\item
Let \((x,s)\in\fP\setminus\fL\) represent a blue triangle. Then
\((x,s)\) is associated with one or two descent critical points, and the
point of tangency corresponding to each such critical point lies on the
portion of the arctic boundary contained in
\(\operatorname{Cell}^{\rb}(x,s)\).

\item
Let \((x,s)\in\fP\setminus\fL\) represent a white triangle. Then
\((x,s)\) is associated with one or two ascent critical points, and the
point of tangency corresponding to each such critical point lies on the
portion of the arctic boundary contained in
\(\operatorname{Cell}^{\rw}(x,s)\).
\end{enumerate}
\end{lemma}

\begin{proof}[Proof of \Cref{l:number_critical_point}]
The statement in the liquid region follows immediately from
\eqref{e:cell_liquid} and \Cref{p:associate_critical_points}.

We prove the statement for a white triangle represented by \((x,s)\). The
statement for blue triangles follows from the same argument, with ascent and
descent interchanged.

Suppose that \((x,s)\in\fP\setminus\fL\). By
\Cref{p:associate_critical_points} and
\eqref{e:descent_ascent_count}, \((x,s)\) is associated with one or two
ascent critical points. It remains to show that their points of tangency lie
on the portion of the arctic boundary contained in
\[
\fT:=\operatorname{Cell}^{\rw}(x,s).
\]

If \((x,s)\) belongs to a unique curvilinear triangle, this follows directly
from \Cref{p:associate_critical_points}. It remains to consider the two
ambiguous cases in \Cref{def:Cell}. Let \(L(w_0)\) be the line supporting the
shared horizontal boundary. The corresponding critical points are described
in \Cref{l:describe_critical_point}.

In the setting of \Cref{i:point} in \Cref{def:Cell}, since \((x,s)\)
represents a white triangle, we have
\[
(x,s-\varepsilon)\in\fP
\]
for every sufficiently small \(\varepsilon>0\). Thus, the polygon \(\fP\)
lies locally below the horizontal supporting line, and Panel (B) of
\Cref{f:c_triangle} cannot occur.

In Panel (A) of \Cref{f:c_triangle}, the unique ascent critical point is
\(w_0\), and its point of tangency is \(\zeta\). By
\eqref{e:cell_point}, this point lies on the portion of the arctic boundary
contained in \(\fT\).

In the setting of \Cref{i:segment} in \Cref{def:Cell}, the assignment
convention gives
\[
\fT=\fT^\downarrow.
\]
If \((x,s)\) lies in the relative interior of the shared segment, the unique
ascent critical point is \(\xi_c^{\rm a}\), whose point of tangency lies on
the portion of the arctic boundary contained in \(\fT^\downarrow\). At a
tangency point or a cusp-turning endpoint, the same conclusion follows by
continuity as the critical points coalesce.

The statement for blue triangles follows analogously, so we omit the proof.
\end{proof}

\subsection{Perturbed critical point sets}
In this section we introduce the limits of critical point sets
under downward perturbations and collect some related estimates.

Let
$
(x,s)\in\bZ^2/n
$
represent a white triangle contained in \(\fP\).
Suppose first that \((x,s)\) lies on a horizontal extended side. Let
\(L(w_0)\) be its supporting line, so that \(w_0=\infty\) is the
corresponding formal critical point, and set
\[
\fT:=\operatorname{Cell}^{\rw}(x,s).
\]
Since the white triangle represented by \((x,s)\) is contained in \(\fP\),
the assigned cell \(\fT\) contains points \((x',s')\) arbitrarily close to
\((x,s)\) with \(s'<s\). We define
\begin{align}\label{e:defCdown}
(\operatorname{Crit}^{\rm d}(x,s^-), \operatorname{Crit}^{\rm a}(x,s^-))
:=
\lim_{\substack{(x',s')\to(x,s)\\
s'<s,\ (x',s')\in\operatorname{int}\fT}}
(\operatorname{Crit}^{\rm d}(x',s'),\operatorname{Crit}^{\rm a}(x',s')),
\end{align}
where the limit is taken through generic points. As we will see in
\Cref{l:descent_critical_below}, the downward perturbation may turn a
spurious critical point back into a descent critical point, but it does not
alter the ascent critical points.

If \((x,s)\) does not lie on a horizontal extended side, we set
\[
(\operatorname{Crit}^{\rm d}(x,s^-),\operatorname{Crit}^{\rm a}(x,s^-))
:=
(\operatorname{Crit}^{\rm d}(x,s), \operatorname{Crit}^{\rm a}(x,s)).
\]

The following lemma gives an explicit description of the one-sided limits
of the critical point sets and, in particular, shows that the limit in
\eqref{e:defCdown} is well defined.

\begin{lemma}\label{l:descent_critical_below}
Let \((x,s)\) represent a white triangle contained in \(\fP\), and suppose
that \((x,s)\) lies on a horizontal extended side with supporting line
\(L(w_0)\). Then
\(\operatorname{Crit}^{\rm a}(x,s^-)=\operatorname{Crit}^{\rm a}(x,s)\).

\begin{enumerate}
\item
If \((x,s)\) lies on a horizontal boundary shared by two adjacent
curvilinear triangles, then
\begin{align}\label{e:critical_below_shared}
\operatorname{Crit}^{\rm d}(x,s^-)=\{w_0\}.
\end{align}

\item
Suppose that \((x,s)\) belongs to a unique curvilinear triangle \(\fT\).
Then
\begin{align}\label{e:critical_below_unique}
\operatorname{Crit}^{\rm d}(x,s^-)
=
\begin{cases}
\operatorname{Crit}^{\rm d}(x,s)\cup\{w_0\},
& \text{if }(x,s)\in\ell_-(w_0;\fT),\\[1mm]
\operatorname{Crit}^{\rm d}(x,s),
& \text{if }(x,s)\notin\ell_-(w_0;\fT).
\end{cases}
\end{align}
\end{enumerate}
\end{lemma}

\begin{remark}\label{r:change_critical}
It follows from the proof of \Cref{l:descent_critical_below} that the
following statements hold. Suppose that \((x,s)\in\fT\) lies on a
horizontal extended side. If \((x',s')\in\fT\) approaches \((x,s)\) from
below, then all ascent critical points persist. By symmetry, if
\((x',s')\in\fT\) approaches \((x,s)\) from above, then all descent
critical points persist.

The corresponding statements for vertical and unit-slope extended sides
follow by symmetry. For a vertical extended side, ``from above'' and ``from
below'' are replaced by ``from the left'' and ``from the right,''
respectively. For a unit-slope extended side, they are replaced by ``from
the right'' and ``from the left,'' respectively.
\end{remark}

\begin{figure}
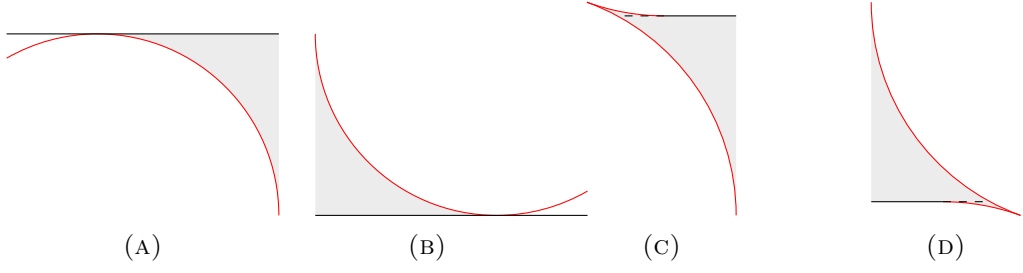


   \begin{subfigure}{0.26\textwidth}
    \centering

          \caption{}
  \end{subfigure}
  \begin{subfigure}{0.18\textwidth}
    \centering
      %
      \caption{}
  \end{subfigure}
    \begin{subfigure}{0.18\textwidth}
    \centering
      %
          \caption{}
  \end{subfigure}
\begin{subfigure}{0.26\textwidth}
    \centering
      %
          \caption{}
  \end{subfigure}
 
  \caption{Curvilinear triangles with a horizontal side}
  \label{f:c_triangle}
\end{figure}

\begin{proof}
Suppose first that the two curvilinear triangles meet only at
\(\zeta=(x,s)\), as in \Cref{i:point} of \Cref{def:Cell}. For a white triangle contained in
\(\fP\), only Panel (A) of \Cref{f:c_triangle} can occur. In this case, by
\Cref{p:descent_ascent_degenerate}, \(w_0\) is an ascent critical point at
\((x,s)\), so
\[
\operatorname{Crit}^{\rm d}(x,s)=\emptyset,
\quad
\operatorname{Crit}^{\rm a}(x,s)=\{w_0\}.
\]
Under a downward perturbation into
\(\operatorname{Cell}^{\rw}(x,s)\), there are nearby ascent and descent
critical points that converge to \(w_0\), and the slopes of their
corresponding tangent lines tend to \(0\). Hence
\[
\operatorname{Crit}^{\rm d}(x,s^-)
=
\operatorname{Crit}^{\rm a}(x,s^-)
=
\{w_0\}.
\]

Suppose next that the shared boundary is a nondegenerate horizontal segment.
Then
\[
\operatorname{Cell}^{\rw}(x,s)=\fT^\downarrow.
\]
If \((x,s)\) lies in the relative interior of the shared segment,
\Cref{l:describe_critical_point} gives
\[
\operatorname{Crit}^{\rm a}(x,s)=\{\xi_c^{\rm a}\},
\qquad
\operatorname{Crit}^{\rm d}(x,s)=\{\xi_c^{\rm d}\},
\]
where \(\xi_c^{\rm a}\) is associated with \(\fT^\downarrow\) and
\(\xi_c^{\rm d}\) is associated with \(\fT^\uparrow\). Under a downward
perturbation into \(\fT^\downarrow\), the branch converging to
\(\xi_c^{\rm a}\) remains an ascent branch, whereas a descent critical point
converges to the formal point \(w_0\). The critical point
\(\xi_c^{\rm d}\), which belongs to the upper triangle, is not seen by this
downward perturbation. Therefore,
\[
\operatorname{Crit}^{\rm d}(x,s^-)=\{w_0\},
\quad
\operatorname{Crit}^{\rm a}(x,s^-)=\{\xi_c^{\rm a}\}.
\]
At a tangency point or a cusp-turning endpoint, the same conclusion follows
by continuity as some of the points
\(\xi_c^{\rm d},w_0,\xi_c^{\rm a}\) coalesce. This proves
\eqref{e:critical_below_shared}.

Finally, suppose that \((x,s)\) belongs to a unique curvilinear triangle
\(\fT\). By the geometric description
\eqref{e:geometric_descent_ascent}, a downward perturbation does not alter
the ascent critical points, and a descent critical point converges to
\(w_0\) precisely when
\[
(x,s)\in\ell_-(w_0;\fT).
\]
This gives \eqref{e:critical_below_unique} and completes the proof.
\end{proof}

\begin{proposition}\label{p:critical_compare}
For any lattice point \((u,v)\in\bZ^2/n\) representing a white triangle
contained in \(\fP\), and any
\[
\xi_c^{\rm d}
\in
\operatorname{Crit}^{\rm d}(u,v^-),
\qquad
\xi_c^{\rm a}
\in
\operatorname{Crit}^{\rm a}(u,v),
\]
one has
\begin{align}\label{e:critical_compare}
\Re[S(\xi_c^{\rm a};u,v)]
\geq
\Re[S(\xi_c^{\rm d};u,v)].
\end{align}
If we further assume that \((u,v)\in\fP\setminus\fL\) and
\(\xi_c^{\rm a}\) and \(\xi_c^{\rm d}\) are bounded away from each other,
then there exists \(\fc>0\) such that
\begin{align}\label{e:critical_compare2}
\Re [S(\xi_c^{\rm a};u,v)]
\geq
\Re [S(\xi_c^{\rm d};u,v)]+\fc.
\end{align}
\end{proposition}

\begin{proof}[Proof of \Cref{p:critical_compare}]
Suppose first that \((u,v)\in\fL\). Then
\[
\operatorname{Crit}^{\rm d}(u,v^-)
=
\operatorname{Crit}^{\rm d}(u,v)
=
\operatorname{Crit}^{\rm a}(u,v)
\]
consists of a pair of complex-conjugate critical points. Since the real part
of the action is invariant under complex conjugation,
\eqref{e:critical_compare} holds with equality.

Suppose next that \((u,v)\in\fP\setminus\fL\) belongs to a single
curvilinear triangle \(\fT\) and lies neither on the arctic boundary nor on
a horizontal extended side. Then
\[
\operatorname{Crit}^{\rm d}(u,v^-)
=
\operatorname{Crit}^{\rm d}(u,v).
\]
By \Cref{p:associate_critical_points}, the point \((u,v)\) has two or three
distinct nondegenerate associated critical points. Order them along the real
arc of \(\cC(\bR)\) corresponding to \(\fT\).

Along this arc, descent critical points are local minima of
\(\Re[ S(\,\cdot\,;u,v)]\), whereas ascent critical points are local maxima,
and the two types alternate. Thus, in the two-critical-point case, there is
one critical point of each type, while in the three-critical-point case the
pattern is
\[
\text{ascent--descent--ascent}
\qquad\text{or}\qquad
\text{descent--ascent--descent};
\]
see \Cref{f:critical_point_pattern}. Since
\(\Re [S(\,\cdot\,;u,v)]\) is strictly monotone between consecutive critical
points, every ascent critical value is larger than every descent critical
value. This proves \eqref{e:critical_compare} in the generic case. If
\(\xi_c^{\rm a}\) and \(\xi_c^{\rm d}\) are bounded away from each other,
continuity and compactness give a uniform \(\fc>0\), and hence
\eqref{e:critical_compare2}.

Now suppose that \((u,v)\) belongs to a single curvilinear triangle
\(\fT\), but lies on the arctic boundary or on a horizontal extended side.
We perturb \((u,v)\) within \(\fT\) to generic points
\[
(u',v')\in\operatorname{int}\fT,
\qquad
v'<v,
\qquad
(u',v')\longrightarrow(u,v),
\]
in the same manner as in the definition of
\(\operatorname{Crit}^{\rm d}(u,v^-)\) in \eqref{e:defCdown}. In this way,
\(\operatorname{Crit}^{\rm d}(u',v')\) converges to
\(\operatorname{Crit}^{\rm d}(u,v^-)\), while
\(\operatorname{Crit}^{\rm a}(u',v')\) converges to
\(\operatorname{Crit}^{\rm a}(u,v)\) by \Cref{l:descent_critical_below}.

Under this perturbation, the corresponding critical values vary
continuously. Applying the generic inequality before taking the limit and
using continuity of the critical values gives
\eqref{e:critical_compare}. If \(\xi_c^{\rm a}\) and
\(\xi_c^{\rm d}\) are bounded away from each other, they do not arise from
a coalescing pair. The strict monotonicity on the limiting real arc,
together with continuity and compactness, then gives
\eqref{e:critical_compare2}.

It remains to consider a point \((u,v)\) on a horizontal boundary shared by
two adjacent curvilinear triangles. Let \(L(w_0)\) be the common supporting
line, and let \(w_0=\infty\) be its corresponding formal critical point.
In this case, \Cref{l:describe_critical_point,l:descent_critical_below}
classify the corresponding critical points and give
\[
\operatorname{Crit}^{\rm d}(u,v^-)=\{w_0\},\quad \xi_c^{\rm d}=w_0.
\]

If the two curvilinear triangles meet only at the point
\(\zeta=(u,v)\), then \(w_0\) is the unique ascent critical point, and the
claim is immediate. Otherwise, their common boundary is a nondegenerate
horizontal segment. Let \(\fT^\uparrow\) and \(\fT^\downarrow\) denote the
adjacent curvilinear triangles lying above and below this segment,
respectively.

Suppose first that \((u,v)\) lies in the relative interior of the shared
segment. Besides \(w_0\), each triangle contributes one additional critical
point. Denote these points by
\[
\widetilde{\xi}_c^{\rm d}
\quad\text{for }\fT^\uparrow,
\qquad
\xi_c^{\rm a}
\quad\text{for }\fT^\downarrow.
\]
Along the real arc corresponding to \(\fT^\downarrow\), the real part of the
action increases from \(w_0\) to the ascent point \(\xi_c^{\rm a}\). Along
the real arc corresponding to \(\fT^\uparrow\), it decreases from \(w_0\)
to the descent point \(\widetilde{\xi}_c^{\rm d}\). Hence
\[
\Re S(\xi_c^{\rm a};u,v)
\geq
\Re S(w_0;u,v)
\geq
\Re S(\widetilde{\xi}_c^{\rm d};u,v).
\]
Since
$
\operatorname{Crit}^{\rm d}(u,v^-)=\{w_0\}
$ and $\operatorname{Crit}^{\rm a}(u,v)=\{\xi_c^{\rm a}\}$
this gives the required comparison. Moreover, if
\(\xi_c^{\rm a}\) and \(w_0\) are bounded away from each other, the first
inequality is strict, and continuity and compactness give
\eqref{e:critical_compare2}.

At a tangency point or a cusp-turning endpoint, some of these critical
points coalesce. The same inequalities follow by continuity from the
relative interior of the shared segment. If \(\xi_c^{\rm a}\) and
\(\xi_c^{\rm d}\) are bounded away from each other, they do not belong to a
coalescing pair, so the inequality remains strict and
\eqref{e:critical_compare2} follows by continuity and compactness. This
completes the proof.
\end{proof}

\chapter{Charts and Contours}
	\section{Ansatz for Liquid-Region Kernel}\label{s:liquid_kernel}
In this section, we construct the liquid-region ansatz
\eqref{e:bulk_ansatz} for the inverse Kasteleyn matrix. The ansatz is given
by the sum of a single-contour integral and a double-contour integral, both
built from the discrete heat kernel \eqref{e:discrete_heat_kernel}. We then
introduce several modified factors that will later be used to construct the
ansatz for the general kernel.

\subsection{Liquid neighborhood}\label{s:critical_bulk}

For any  $(x_0, s_0)\in \fL$ bounded away from the ramification point and the arctic boundary, let
\begin{align}
w_0=x_0-s_0 \chi(x_0, s_0)=x_0-s_0 \chi(w_0)\in \bC_+.
\end{align}

We recall from \eqref{e:two_points}, $w_0$ is a critical point of $S(\cdot; x_0,s_0)$ and $S'(w_0;x_0,s_0)=0$. Moreover, from \eqref{e:derSsecond-chi}, locally around $w_0$, $S(w;x_0,s_0)$ is holomorphic, and we have the Taylor expansion
\begin{align}\begin{split}\label{e:liquid_neighborhood}
&S(w;x_0,s_0)=S(w_0;x_0,s_0)+S''(w_0; x_0, s_0)\frac{(w-w_0)^2}{2} +\OO(|w-w_0|^3),\\
 &S''(w_0; x_0, s_0)=-\left(\chi'(w_0)+\frac{1}{s_0}\right)\frac{1}{\chi(w_0)(1-\chi(w_0))}
\end{split}\end{align}
In this case $S''(w_0; x_0,s_0)$ is bounded away from $0$, and the steepest descent direction is determined by
\begin{align}\label{e:deftheta0}
2\theta_0 := \pi -\arg S''(w_0; x_0, s_0),\quad
e^{2\ri\theta_0}:=-\frac{|S''(w_0; x_0, s_0)|}{S''(w_0; x_0, s_0)}.
\end{align}

The following lemma introduces the liquid chart and records several of its properties. Its proof, based on a Taylor expansion of the tiling action, is deferred to \Cref{s:liquid_chart_proof}.
\begin{lemma}[$\fc$-liquid chart]\label{c:bulk}
Fix $\delta_0>0$ and let $(x_0,s_0)\in\fL$ be a point whose distance is at least $\delta_0$
from both the arctic boundary and the set of ramification points. Then for every sufficiently
small $\fc>0$ (depending only on $\delta_0$) there exists $\delta=\delta(\fc)>0$, also sufficiently
small and depending only on $\delta_0$, such that the following holds.

Let $w_0=x_0-s_0\,\chi(x_0,s_0)$ be the critical point corresponding to $(x_0,s_0)$, and set
\[
\fU:=\{w\in\bC:\ |w-w_0|\le \fc\}.
\]
On $\fU$ the Riemann surface $\cC$ can be parametrized as $(f(w),w)$. We will therefore identify
$\fU$ with its image in $\cC$ under the map $w\mapsto (f(w),w)$.  Moreover, for all $w\in\fU$,
\begin{align}\label{e:bulk_bounds}
\Im[w], -\Im[\chi(w)]\asymp 1,
\end{align}
with implicit constants depending only on $\delta_0$.

Now let $(x,s)\in\fP$ satisfy $\|(x,s)-(x_0,s_0)\|_2\le \delta$. Then the tiling action $S(\,\cdot\,;x,s)$ has exactly one critical point $w_c$ inside $\fU$, satisfying
\begin{align}\label{e:wcbb}
|w_c-w_0|
\lesssim \|(x,s)-(x_0,s_0)\|_2
\leq \delta.
\end{align}
Moreover, for every $w\in\fU$,
\begin{align}\label{e:bulk_Taylor}
S(w;x,s)-S(w_c;x,s)=d\,(w-w_c)^2+\cE(w),
\qquad d:=\frac12\,S''(w_0;x_0,s_0),\quad  |d|\asymp 1.
\end{align}
The error term satisfies
\begin{align}\label{e:bulk_err}
|\cE(w)|\leq C\bigl(\delta\ln(1/\delta)+|w-w_c|^3\bigr)\le \frac{|d|\,\fc^2}{100}.
\end{align}

In this situation, we call $\fU$ a $\fc$-liquid chart (centered at $w_0$), and we say that the
point $(x,s)$ is \emph{adapted} to $\fU$. 
\end{lemma}

\begin{remark}
Since the Riemann surface \(\cC\) and the critical-point equation
\eqref{e:critical_point} are invariant under complex conjugation, the point
\(\overline{w_0}\) is the other critical point corresponding to \((x_0,s_0)\).
We also call the conjugate disk
\[
\overline{\fU}:=\{w\in\bC:\ |w-\overline{w_0}|\le \fc\}
\]
a \(\fc\)-liquid chart centered at \(\overline{w_0}\). On \(\overline{\fU}\),
the Riemann surface \(\cC\) can again be parametrized as \((f(w),w)\). Moreover,
if \(w_c\in\fU\) is the critical point associated with \((x,s)\), then the
corresponding critical point in \(\overline{\fU}\) is \(\overline{w_c}\).
\end{remark}

\begin{figure}
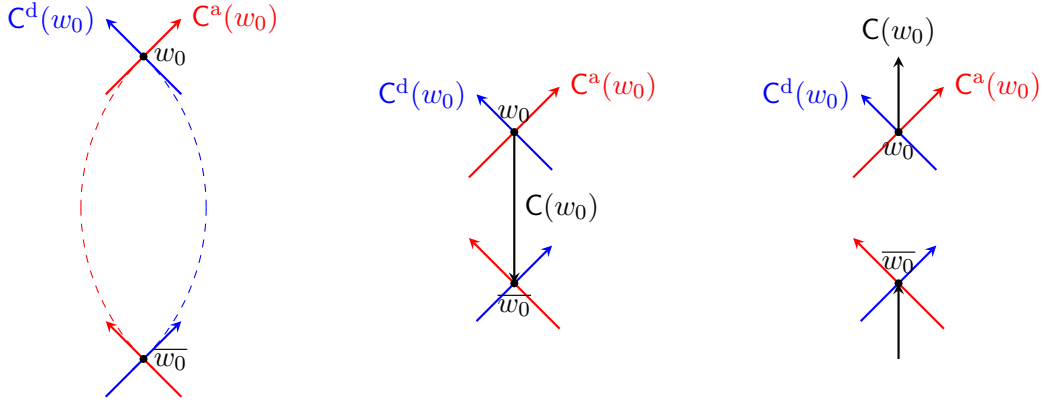

		\begin{subfigure}[t]{0.3\textwidth}
		\begin{center}


				
			\end{center}
			\end{subfigure}

	\caption{
	Local descent/ascent segments used in the kernel ansatz \eqref{e:bulk_ansatz2}; Middle panel for $s<t$; Right panel for $s\geq t$. }
	\label{f:bulk_path}
\end{figure}

We now introduce the local descent and ascent segments through \(w_0\); see \Cref{f:bulk_path}. Let
\(\pm e^{\ri\theta_0}\) be the steepest-descent directions at \(w_0\), where
\(\theta_0\) is defined in \eqref{e:deftheta0}. For \(\fc>0\) sufficiently small,
define
\begin{equation}\label{e:gamma_path}
\sfC^{\rm d}(w_0)
:=
\bigl\{\,w_0+r e^{\ri\theta_0}:\ -\fc\le r\le \fc\,\bigr\},
\qquad
\sfC^{\rm a}(w_0)
:=
\bigl\{\,w_0-\ri r e^{\ri\theta_0}:\ -\fc\le r\le \fc\,\bigr\}.
\end{equation}
Thus \(\sfC^{\rm d}(w_0)\) is the local descent segment, while
\(\sfC^{\rm a}(w_0)\) is obtained from it by a rotation by \(-\pi/2\) and is the
corresponding local ascent segment.

At the conjugate critical point \(\overline{w_0}\), the steepest-descent directions
are \(\pm e^{-\ri\theta_0}\). We define the local descent and ascent segments through
\(\overline{w_0}\) by
\begin{equation}\label{e:conj_gamma_path}
\sfC^{\rm d}(\overline{w_0})
:=
\bigl\{\,\overline{w_0}+r e^{-\ri\theta_0}:\ -\fc\le r\le \fc\,\bigr\}
=
\overline{\sfC^{\rm d}(w_0)},
\qquad
\sfC^{\rm a}(\overline{w_0})
:=
\bigl\{\,\overline{w_0}+\ri r e^{-\ri\theta_0}:\ -\fc\le r\le \fc\,\bigr\}
=
\overline{\sfC^{\rm a}(w_0)}.
\end{equation}
We notice  that
\[
\sfC^{\rm d}(w_0),\ \sfC^{\rm a}(w_0)\subset \fU,
\qquad
\sfC^{\rm d}(\overline{w_0}),\ \sfC^{\rm a}(\overline{w_0})\subset \overline{\fU}.
\]

\begin{figure}
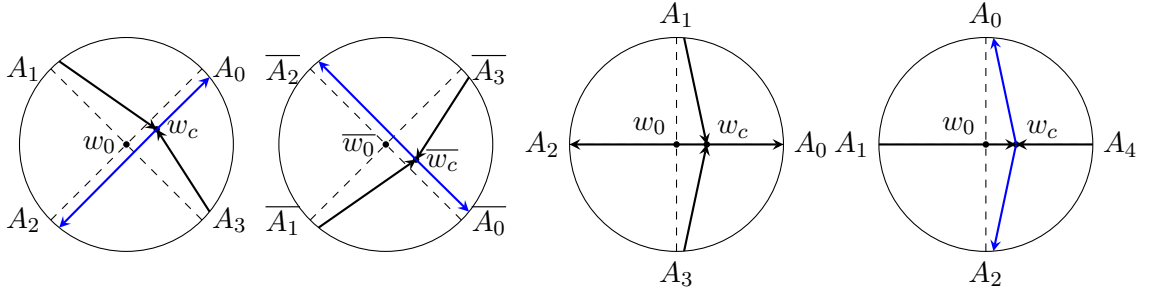

		\begin{subfigure}[t]{0.2\textwidth}
		\begin{center}


				
			\end{center}
			\end{subfigure}

	\caption{
Gradient flow of \(S(\cdot;x,s)\) in a liquid chart (left) and a regular frozen chart (right).
}
	\label{f:bulk_critical}
\end{figure}

The following lemma shows that the local descent and ascent paths can be deformed into steepest-descent and steepest-ascent paths with negligible error. Its proof is deferred to \Cref{s:liquid_chart_proof}.
\begin{lemma}\label{c:bulk_steepest}
Adopt the assumptions and notation of \Cref{c:bulk}. Then there exists a
constant $\fc'>0$ such that the following statements hold.
\begin{itemize}
\item \noindent\emph{Replacing \(\mathsf C^{\rm d}(w_0)\) by \(\mathsf D^{\rm d}(w_c)\).}
We define the local steepest-descent set $\mathsf D^{\rm d}(w_c)$ at \(w_c\)
to be the portions of the steepest-descent trajectories of
$S(\,\cdot\,;x,s)$ issuing from \(w_c\) up to their first exit from the disk
\[
\{\,w: |w-w_0|\le \fc\,\};
\]
see \Cref{f:bulk_critical}. Then $\mathsf D^{\rm d}(w_c)$ has total length
$\OO(1)$. Moreover, $\sfC^{\rm d}(w_0)$ can be deformed to
$\mathsf D^{\rm d}(w_c)$, together with finitely many arcs of total length
$\OO(1)$, on which
\begin{align}\label{e:bulk_subarc}
n\,\Re[ S(w;x,s)-S(w_c;x,s)]\le -\,n\fc'.
\end{align}

\item \noindent\emph{Replacing \(\mathsf D^{\rm d}(w_c)\) by \(\mathsf S^{\rm d}(w_c)\).}
Let
$
r_n:={\ln n}/{\sqrt n}$.
We define the truncated steepest-descent set $\mathsf S^{\rm d}(w_c)$ at
\(w_c\) to be the portions of the steepest-descent trajectories issuing from
\(w_c\) up to their first exit from the disk
\[
\{\,w: |w-w_c|\le r_n\,\}.
\]
Then, $\sfS^{\rm d}(w_c)$ has total length $\OO(r_n)$, and  for all \(w\in \mathsf D^{\rm d}(w_c)\setminus
\mathsf S^{\rm d}(w_c)\),
\[
e^{n\Re[S(w;x,s)]}
\le
e^{n\Re[S(w_c;x,s)]}e^{-\fc'(\ln n)^2}.
\]
\end{itemize}
We define the analogous ascent sets
\(\mathsf D^{\rm a}(w_c)\) and \(\mathsf S^{\rm a}(w_c)\). The same
statements hold with \(S\) replaced by \(-S\) and with the superscript
\({\rm d}\) replaced by \({\rm a}\).
\end{lemma}

\subsection{Ansatz for liquid region kernel}
\label{s:bulk_kernel}
Fix a point \((x_0,s_0)\in\fL\) bounded away from both the arctic boundary and the
ramification points, and let
$
w_0=x_0-s_0\chi(w_0)
$
be the corresponding critical point. Recall \(\delta>0\) from \Cref{c:bulk}, and set
\begin{align}\label{e:defN}
\fN=\bigl\{(x,s): \|(x,s)-(x_0,s_0)\|_2\le \delta\bigr\}.
\end{align}
In this section we explain the kernel ansatz for rescaled lattice points
$
(x,s),(y,t)\in \fN\cap \bZ^2/n $.

We introduce the following two 
functions \(P\) and \(Q\), where $P$ is the binomial coefficient \eqref{e:binomial} from the solution of complex Burgers equation, while $Q$ is its backward counterpart.
\begin{align}\label{e:defPQ}
\begin{split}
P_{ns}(nw,nx)
&:=
{ns\choose n(x-w)}
=
\frac{\Gamma(ns+1)}
{\Gamma(n(x-w)+1)\Gamma(n(w-(x-s))+1)}, \\
Q_{nt}(nz,ny)
&:=
\frac{\Gamma(n(y-z))\Gamma(n(z-(y-t)))}{\Gamma(nt)} .
\end{split}
\end{align}
Both $P$ and $Q$ satisfy heat-type recursions, which follows from elementary recursions for the two Gamma-function
factors in \eqref{e:defPQ}.

\begin{lemma}\label{c:recursion}
The functions \(P\) and \(Q\) in \eqref{e:defPQ} satisfy the discrete forward and
backward heat-type recursions
\begin{align}\label{e:PQ}
\begin{split}
P_{ns}(nw,nx)
&=
P_{ns-1}(nw,nx)+P_{ns-1}(nw,nx-1), \\
Q_{nt}(nz,ny)
&=
Q_{nt+1}(nz,ny)+Q_{nt+1}(nz,ny+1).
\end{split}
\end{align}
\end{lemma}

\begin{proof}[Proof of \Cref{c:recursion}]
The first identity in \eqref{e:PQ} is Pascal's identity. For \(Q\), using
\(\Gamma(u+1)=u\Gamma(u)\), we obtain
\begin{align*}
&\phantom{{}={}}
Q_{nt+1}(nz,ny)+Q_{nt+1}(nz,ny+1) \\
&=
\frac{
\Gamma(n(y-z))\Gamma(n(z-(y-t))+1)
+
\Gamma(n(y-z)+1)\Gamma(n(z-(y-t)))
}
{\Gamma(nt+1)} \\
&=
\frac{
n(z-(y-t))\Gamma(n(y-z))\Gamma(n(z-(y-t)))
+
n(y-z)\Gamma(n(y-z))\Gamma(n(z-(y-t)))
}
{nt\,\Gamma(nt)} \\
&=
\frac{\Gamma(n(y-z))\Gamma(n(z-(y-t)))}{\Gamma(nt)}
=
Q_{nt}(nz,ny).
\end{align*}
\end{proof}

We recall the uniformizing conformal map $\phi$ from \Cref{s:cf_map},  and introduce the following kernel ansatz:
\begin{align}\label{e:bulk_ansatz}
\begin{split}
&\phantom{{}={}}A((x,s),(y,t))
:=
\frac{n}{2\pi \ri}
\int_{\sfC}
P_{ns}(nz,nx)\,Q_{nt}(nz,ny)\,\rd z \\
&\quad
+\frac{n}{(2\pi \ri)^2}
\int_{\sfC^{\rm a}}\rd z
\int_{\sfC^{\rm d}}\rd w\,
P_{ns}(nw,nx)
e^{-n\int_0^w \ln f(u)\rd u} 
\times
Q_{nt}(nz,ny)
e^{n\int_0^z \ln f(u)\rd u}
\frac{\sqrt{\phi'(w)}\sqrt{\phi'(z)}}{\phi(w)-\phi(z)} .
\end{split}
\end{align}
Here the contours
\[
\sfC\subset\bC,
\qquad
\sfC^{\rm d},\sfC^{\rm a}\subset \fU\cup\overline{\fU},
\]
as well as the branches of \(\sqrt{\phi'(z)}\) and \(\sqrt{\phi'(w)}\) appearing in
\eqref{e:bulk_ansatz}, will be specified in \Cref{s:liquid_integral_contours}. 
We remark that the integrand in the single-contour integral in
\eqref{e:bulk_ansatz} is the residue at \(w=z\) of the integrand in the
double-contour integral.

Asymptotically, the kernel ansatz \eqref{e:bulk_ansatz} reduces to single- and
double-contour integrals whose exponential parts are governed by the tiling action
\eqref{e:def_action}.  Stirling's formula may be applied uniformly
to the Gamma functions appearing in \(P\) and \(Q\).

\begin{lemma}\label{l:PIQI_bound}
Given any complex number \(w\in \bC\) and assume
\begin{align}\label{e:d1}
\delta_w:=\operatorname{dist}\left(w,(-\infty, x-s]\cup [x, \infty)\right)\gg n^{-1}.
\end{align}
Then
\begin{align}
P_{ns}(nw,nx) e^{-n\int_0^w \ln f(u)\rd u} 
=
\frac{\sqrt{s}}{\sqrt{2\pi n}\sqrt{x-w}\sqrt{w-(x-s)}} e^{n S(w;x,s)+\OO(1/(\delta_w n))}
\end{align}
Here \(S(\,\cdot\,;x,s)\) is the tiling action function
defined in \eqref{e:def_action}. 

Similarly, given any complex number \(z\in \bC\) and assume
\begin{align}\label{e:d2}
\delta_z:=\operatorname{dist}\left(z,(-\infty, y-t]\cup [y, \infty)\right)\gg n^{-1}.
\end{align}
Then
\begin{align}
Q_{nt}(nz,ny) e^{n\int_0^z \ln f(u)\rd u} = \frac{\sqrt{2\pi t}}{\sqrt{ n}\sqrt{y-z}\sqrt{z-(y-t)}} e^{-n S(z;y,t)+\OO(1/(\delta_z n))}.
\end{align}
\end{lemma}
\begin{proof}[Proof of \Cref{l:PIQI_bound}]

Recalling \eqref{e:binomial}, under assumption \eqref{e:d1} we have
\begin{align}\label{e:binomial1}
\begin{split}
\ln P_{ns}(nw,nx)
%&=
%\ln \Gamma(ns+1)
%-\ln \Gamma(n(x-w)+1)
%-\ln\Gamma(n(s-x+w)+1) \\
&=
n\Bigl[
s\ln s
-(x-w)\ln (x-w)
-(s-x+w)\ln (s-x+w)
\Bigr] \\
&\quad
+\frac12\Bigl[
\ln s-\ln (x-w)-\ln (s-x+w)
\Bigr]
-\frac12\ln(2\pi n)
+\OO\left(\frac1{\delta_w n}\right).
\end{split}
\end{align}
The claim follows from the definition of the tiling action \eqref{e:def_action}.
Similarly, using \eqref{e:logGamma}, under assumption \eqref{e:d2} we have
\begin{align}\label{e:binomial2}
\begin{split}
\ln Q_{nt}(nz,ny)
&=
\ln \Gamma(n(y-z))
+\ln\Gamma(n(z-(y-t)))
-\ln \Gamma(nt) \\
&=
-n\Bigl[
t\ln t
-(y-z)\ln(y-z)
-(z-(y-t))\ln(z-(y-t))
\Bigr] \\
&\quad
+\frac12\Bigl[
\ln t-\ln(y-z)-\ln(z-(y-t))
\Bigr]
+\frac12\ln\left(\frac{2\pi}{n}\right)
+\OO\left(\frac1{\delta_z n}\right).
\end{split}
\end{align}
\end{proof}

Because the contours in \(\fU\cup\overline{\fU}\) remain a fixed
positive distance from the real axis, \Cref{l:PIQI_bound} applies along
these contours. After extracting the leading prefactors, the ansatz
\eqref{e:bulk_ansatz} takes the following asymptotic form:
\begin{align}\label{e:bulk_ansatz2}
\begin{split}
&\phantom{{}={}}A((x,s),(y,t))
=
\frac{\sqrt{st}}{2\pi \ri}
\int_{\sfC}
\frac{
e^{n\left(S(z;x,s)-S(z;y,t)+\OO(1/n))\right)}
}
{\sqrt{x-w}\sqrt{w-(x-s)}\sqrt{y-z}\sqrt{z-(y-t)}}
\rd z \\
&\quad
+\frac{\sqrt{st}}{(2\pi \ri)^2}
\int_{\sfC^{\rm a}}\rd z
\int_{\sfC^{\rm d}}\rd w\,
\frac{\sqrt{\phi'(w)}\sqrt{\phi'(z)}}
{\sqrt{x-w}\sqrt{w-(x-s)}\sqrt{y-z}\sqrt{z-(y-t)}}
\frac{
e^{n\left(S(w;x,s)-S(z;y,t)+\OO(1/n)\right)}
}
{\phi(w)-\phi(z)}
.
\end{split}
\end{align}
Here \(S(\,\cdot\,;x,s)\) and \(S(\,\cdot\,;y,t)\) denote the tiling
action functions defined in \eqref{e:def_action}.

\subsection{Integral contours}
\label{s:liquid_integral_contours}

We first fix a sign convention for transverse intersections of oriented
curves.

\begin{definition}[Positive/negative intersection]\label{d:positive_negative} 
Let $\sfC_1,\sfC_2\subset\bC\simeq\bR^2$ be oriented $C^1$ curves that intersect transversely at a point $\xi$. Let $\tau_1(\xi)$ and $\tau_2(\xi)$ denote tangent vectors to $\sfC_1$ and $\sfC_2$ at $\xi$ compatible with their orientations. We say that $\sfC_1$ and $\sfC_2$ \emph{intersect positively at $\xi$} if rotating $\tau_1(\xi)$ counterclockwise by an angle in $(0,\pi)$ points toward $\tau_2(\xi)$ and \emph{intersect negatively at $\xi$} otherwise. \end{definition}

For the single-contour term in \eqref{e:bulk_ansatz}, we set
\[
\sfC=\sfC(w_0;(x,s),(y,t)),
\]
where the contour on the right-hand side is defined as follows; see
\Cref{f:bulk_path}.

\begin{definition}
\label{d:defCw}
Let \(w_0\in\bC_+\).

\begin{enumerate}
\item
If \(s\geq t\), then
\(\sfC(w_0;(x,s),(y,t))\) is the union of two oriented paths: a path
from \(w_0\) to \(+\infty\) contained in the upper half-plane, and a
path from \(-\infty\) to \(\overline{w_0}\) contained in the lower
half-plane.

\item
If \(s<t\), then
\(\sfC(w_0;(x,s),(y,t))\) is an oriented path from \(w_0\) to
\(\overline{w_0}\) that crosses the real axis at a point of
\[
\bigl(\max\{x-s,y-t\},\,\min\{x,y\}\bigr).
\]
\end{enumerate}
\end{definition}

As we will show later in \Cref{s:pof_single_integral}, these topological choices
allow \(\sfC(w_0;(x,s),(y,t))\) to be deformed, in the corresponding
case, to a steepest-descent contour for the phase
$
S(z;x,s)-S(z;y,t)
$
appearing in the single-contour term of \eqref{e:bulk_ansatz2}.

We next specify the contours in the double integral in \eqref{e:bulk_ansatz}. Recall the local
steepest-descent and steepest-ascent segments
\(\sfC^{\rm d}(w_0)\) and \(\sfC^{\rm a}(w_0)\), together with their
lower-half-plane counterparts, from \eqref{e:gamma_path} and
\eqref{e:conj_gamma_path}. We set
\[
\sfC^{\rm d}
=
\sfC^{\rm d}(w_0)\cup\sfC^{\rm d}(\overline{w_0}),
\qquad
\sfC^{\rm a}
=
\sfC^{\rm a}(w_0)\cup\sfC^{\rm a}(\overline{w_0}).
\]

The steepest descent condition determines an unoriented tangent line at
\(w_0\), but not an orientation of that line. Equivalently, a normalized
descent vector \(v^{\rm d}(w_0)\) satisfies
\[
v^{\rm d}(w_0)^2
=
-\frac{1}{S''(w_0;(x_0,s_0))}
\]
and is therefore determined only up to sign. We couple this sign choice
to the local branch of \(\sqrt{\phi'(w_0)}\).

At the level of squared quantities, \eqref{e:derSsecond-chi} and
\eqref{e:derphi} give
\begin{align}\label{e:Sdirection}
-\frac{1}{S''(w_0;(x_0,s_0))}
=
\frac{s_0\chi(w_0)\bigl(1-\chi(w_0)\bigr)}
     {1+s_0\chi'(w_0)}
=
s_0\chi(w_0)\bigl(1-\chi(w_0)\bigr)
\frac{\partial_x\phi(x_0,s_0)}{\phi'(w_0)}.
\end{align}

By \Cref{c:defsqrtphi}, the function
\(\sqrt{\partial_x\phi(x,s)}\) is defined globally on the liquid region
\(\fL\), up to one global sign. We fix this sign once and for all.
Because \(\chi(w_0)\in\bC_-\), we use the principal branches of
\(\sqrt{\chi(w_0)}\) and \(\sqrt{1-\chi(w_0)}\). For a chosen local
branch of \(\sqrt{\phi'(w_0)}\), define
\begin{equation}
\label{eq:factor}
v^{\rm d}(w_0)
:=
\sqrt{s_0}
\sqrt{\chi(w_0)}
\sqrt{1-\chi(w_0)}
\frac{\sqrt{\partial_x\phi(x_0,s_0)}}{\sqrt{\phi'(w_0)}}.
\end{equation}
We orient \(\sfC^{\rm d}(w_0)\) so that at \(w_0\) it points in
the direction \(v_{\rm d}(w_0)\).

The branch-dependent part of \eqref{eq:factor} is
\begin{equation}
\label{eq:factor2}
\frac{
\sqrt{\chi(w_0)}\sqrt{1-\chi(w_0)}
}{
\sqrt{\phi'(w_0)}
}.
\end{equation}
This factor also appears reciprocally in the prefactor of the double-contour integral \eqref{e:bulk_ansatz2}. As \(w\to w_0\) and \((x,s)\to(x_0,s_0)\),
we obtain
\begin{equation}\label{eq:integrandlocal}
\frac{\sqrt{\phi'(w)}}{\sqrt{x-w}\sqrt{w-(x-s)}}
=
(1+\oo(1))\frac{\sqrt{\phi'(w_0)}}{\sqrt{x_0-w_0}\sqrt{w_0-(x_0-s_0)}}
=
\frac{1}{s}\frac{\sqrt{\phi'(w_0)}}{\sqrt{\chi(w_0)}\sqrt{1-\chi(w_0)}},
\end{equation}
where we used the critical-point relation $w_0=x_0-s_0\,\chi(w_0)$ in the last step. 
Changing the local branch of \(\sqrt{\phi'(w_0)}\) changes the signs of
\eqref{eq:integrandlocal}. At the same time,
\eqref{eq:factor} shows that it reverses the tangent vector
\(v^{\rm d}(w_0)\), and hence reverses the orientation of
\(\sfC^{\rm d}(w_0)\). These two sign changes compensate in the contour
integral.

We now fix the corresponding choices at the conjugate critical point.
Using the Schwarz-symmetry identities
\[
\overline{\phi(w)}=\phi(\overline w),
\qquad
\overline{S(w;x,s)}=S(\overline w;x,s),
\]
we orient the descent segment in the lower half plane by setting
\begin{align}\label{e:vdlow}
v_{\rm d}(\overline{w_0})
:=
-\overline{v_{\rm d}(w_0)}.
\end{align}
Thus, if the oriented unit tangent at \(w_0\) is \(e^{\ri\theta_0}\),
then the oriented unit tangent at \(\overline{w_0}\) is
\(-e^{-\ri\theta_0}\). Consequently, complex conjugation maps each
local descent segment to the other with the opposite orientation; see
\Cref{f:bulk_path}.

We choose the lower-half-plane square-root branch by conjugation:
\[
\frac{\sqrt{\phi'(\overline w)}}{\sqrt{x-\overline  w}\sqrt{\overline  w-(x-s)}}
=
\overline{\frac{\sqrt{\phi'(w)}}{\sqrt{x-w}\sqrt{w-(x-s)}}}
\]
for \(w\) near \(w_0\). 

Finally, define the oriented ascent tangent vectors by
\begin{align}\label{e:va}
v^{\rm a}(w_0):=-\ri\,v^{\rm d}(w_0),
\qquad
v^{\rm a}(\overline{w_0})
:=\ri\,v^{\rm d}(\overline{w_0}),
\end{align}
and orient the corresponding ascent segments accordingly. Since
\[
v^{\rm a}(w_0)^2
=
\frac{1}{S''(w_0;(x_0,s_0))},
\qquad
v^{\rm a}(\overline{w_0})^2
=
\frac{1}{S''(\overline{w_0};(x_0,s_0))},
\]
these are indeed steepest-ascent directions. By
\Cref{d:positive_negative},
\(\sfC^{\rm d}(w_0)\) and \(\sfC^{\rm a}(w_0)\) intersect negatively at
\(w_0\), whereas
\(\sfC^{\rm d}(\overline{w_0})\) and
\(\sfC^{\rm a}(\overline{w_0})\) intersect positively at
\(\overline{w_0}\).

In summary, choosing a local branch of $\sqrt{\phi'(w_0)}$ fixes the orientation of the steepest-descent direction $v^{\rm d}(w_0)$ through \eqref{eq:factor}. Reversing the branch reverses the contour orientation, and the resulting sign changes cancel in the double-contour integral \eqref{e:bulk_ansatz2}, provided that the prefactor $\sqrt{\phi'(w)}$ is defined using the same branch convention.

\subsection{Inverse Kasteleyn matrix relation}
In the rest of this section, we check that the ansatz
\(A((x,s),(y,t))\), defined in \eqref{e:bulk_ansatz}, satisfies the
Kasteleyn kernel relation \eqref{e:discrete_kernel2}.

\begin{proposition}\label{l:Aeq}
For \((x,s),(y,t)\in \fN\cap \bZ^2/n\) (recall from \eqref{e:defN}), the kernel
\(A(\cdot,\cdot)\) defined in \eqref{e:bulk_ansatz} satisfies
\begin{align}\label{e:Aeq0}
\begin{split}
&\phantom{{}={}}
A((x,s),(y,t))
-
A\left(\left(x,s-\frac1n\right),(y,t)\right)
-
A\left(\left(x-\frac1n,s-\frac1n\right),(y,t)\right) =
\delta_{(x,s)=(y,t)} .
\end{split}
\end{align}
\end{proposition}

Before proving \Cref{l:Aeq}, we collect some basic properties of the integrand
$P_{ns}Q_{nt}$ in \Cref{l:integrand_prop}, together with an identity for the
single-contour integral \Cref{c:single_int}.

\begin{figure}
  \centering

  %================= Row 1 =================
  \begin{subfigure}{0.24\textwidth}
    \centering
    \resizebox{\linewidth}{!}{%
      \begin{tikzpicture}
        \draw[] (-2.5,0)--(-0.5,0);
        \foreach \x in {-2,-1.8,-1.6,-1.4,-1.2,-1} {
    \fill (\x,0) circle (1pt);
  }
 % \foreach \x in {1, 1.2, 1.4, 1.6, 1.8, 2} {
 %   \fill (\x,0) circle (1pt);
 %   \draw[fill=white] (0,0) circle (1pt);
 % }
  \node at (-2, -0.2) {\tiny $x-s$};
  \node at (-1, -0.2) {\tiny $y-t$};
  %\node at (1, -0.2) {$y$};
  %\node at (2, -0.2) {$x$};
      \end{tikzpicture}
    }
  \end{subfigure}
 \begin{subfigure}{0.24\textwidth}
    \centering
    \resizebox{\linewidth}{!}{%
      \begin{tikzpicture}
        \draw[] (-2.5,0)--(-0.5,0);
        \foreach \x in {-1.8,-1.6,-1.4,-1.2} {
    \fill[red] (\x,0) circle (1pt);
  }
   \draw[fill=white] (-2,0) circle (1pt);
      \draw[fill=white] (-1,0) circle (1pt);
 % \foreach \x in {1, 1.2, 1.4, 1.6, 1.8, 2} {
 %   \fill (\x,0) circle (1pt);
 %   \draw[fill=white] (0,0) circle (1pt);
 % }
  \node at (-2, -0.2) {\tiny $y-t$};
  \node at (-1, -0.2) {\tiny $x-s$};
  %\node at (1, -0.2) {$y$};
  %\node at (2, -0.2) {$x$};
      \end{tikzpicture}
    }
  \end{subfigure}
  \begin{subfigure}{0.24\textwidth}
    \centering
    \resizebox{\linewidth}{!}{%
      \begin{tikzpicture}
        \draw[] (0.5,0)--(2.5,0);
    \foreach \x in {1, 1.2, 1.4, 1.6, 1.8, 2} {
    \fill (\x,0) circle (1pt);
  }
  \node at (1, -0.2) {\tiny $y$};
  \node at (2, -0.2) {\tiny $x$};
      \end{tikzpicture}
    }
  \end{subfigure}
  \begin{subfigure}{0.24\textwidth}
    \centering
    \resizebox{\linewidth}{!}{%
      \begin{tikzpicture}
        \draw[] (0.5,0)--(2.5,0);
    \foreach \x in {1.2, 1.4, 1.6, 1.8} {
    \fill[red] (\x,0) circle (1pt);
  }
    \draw[fill=white] (2,0) circle (1pt);
      \draw[fill=white] (1,0) circle (1pt);
  \node at (1, -0.2) {\tiny $x$};
  \node at (2, -0.2) {\tiny $y$};
      \end{tikzpicture}
    }
  \end{subfigure}
  \caption{Black dots indicate poles, and red dots indicate zeros.}
  \label{f:pole_zero}
\end{figure}
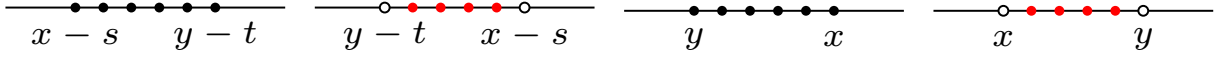

\begin{lemma}\label{l:integrand_prop}
The integrand $P_{ns}(nz, nx) Q_{nt}(nz,ny)$ from \eqref{e:bulk_ansatz} is given by
\begin{align}\label{e:single_integrand}
P_{ns}(nz, nx) Q_{nt}(nz,ny) =
\frac{\Gamma(ns+1)}{\Gamma(nt)}
\frac{\Gamma(n(y-z))\Gamma(n(z-(y-t)))}
{\Gamma(n(x-z)+1)\Gamma(n(z-(x-s))+1)}.
\end{align}
There are several cases for the relative locations of $s,t$, $x,y$ and $x-s, y-t$; see \Cref{f:pole_zero}.
\begin{enumerate}
\item 
There exists a constant $C\neq 0$, as $|z|\rightarrow \infty$, we have the following asymptotic behavior:
\begin{align}\label{e:asymptotic}
P_{ns}(nz, nx) Q_{nt}(nz,ny) \sim (C+\oo(1))\frac{z^{n(t-s)}}{z^2}.
\end{align}
Thus if $s\geq t$, $z^{n(t-s)}/z^2=\OO(z^{-2})$, and the differential $P_{ns}(nz, nx) Q_{nt}(nz,ny) \rd z$
does not have residual at $\infty$;  if $t>s$ it has a pole at $\infty$.

\item If $x\geq y$ then $\Gamma(n(y-z))/\Gamma(n(x-z)+1)$ has poles at 
\begin{align}
x, x-1/n, x-2/n,\cdots, y;
\end{align}
 If $x<y$ then it has zeros at 
 \begin{align}
 y-1/n, y-2/n, \cdots, x+1/n.
 \end{align}
\item If $y-t\geq x-s $ then $\Gamma(n(z-(y-t)))/\Gamma(n(z-(x-s))+1)$ has poles at 
\begin{align}
y-t, y-t-1/n, y-t-2/n,\cdots, x-s;
\end{align}
 If $y-t<x-s$ then it has zeros at 
 \begin{align}
 x-s-1/n, x-s-2/n, \cdots, y-t+1/n.
 \end{align}
\end{enumerate}
In particular, away from $\infty$ and from the two closed intervals
\[
[x-s, y-t],\quad [y,x],
\]
the integrand $P_{ns}(nz, nx) Q_{nt}(nz,ny)$ is holomorphic.
\end{lemma}

\begin{lemma}\label{c:single_int}
Assume \(s\ge t\). Then
\begin{align}\label{e:single_int}
\frac{n}{2\pi \ri}
\int_{-\ri\infty}^{\ri\infty}
P_{ns}(nz,nx)\,Q_{nt}(nz,ny)\,\rd z
=
\bm 1(x\ge y)
{n(s-t)\choose n(x-y)} .
\end{align}
Here the contour is oriented from \(-\ri\infty\) to \(\ri\infty\) and crosses the
real axis at a point in
\[
\bigl(\max\{x-s,y-t\},\,\min\{x,y\}\bigr).
\]
\end{lemma}
\begin{proof}[Proof of \Cref{l:Aeq}]
By \Cref{c:recursion}, the double-contour integral in \eqref{e:bulk_ansatz}
satisfies the homogeneous recursion corresponding to \eqref{e:Aeq0}. Indeed, the
only dependence on \((x,s)\) in that term through the Kasteleyn operator is through
\(P_{ns}(nw,nx)\), and the first identity in \eqref{e:PQ} gives cancellation.

It remains to identify the contribution of the single-contour integral. For
\(s'\ge t\), \Cref{c:single_int} gives
\begin{align}\label{e:single_int2}
\begin{split}
&\phantom{{}={}}
\frac{n}{2\pi \ri}
\left(
\int_{w_0}^{\ri\infty}
+
\int_{-\ri\infty}^{\overline{w_0}}
\right)
P_{ns'}(nz,nx')\,Q_{nt}(nz,ny)\,\rd z \\
&=
\frac{n}{2\pi \ri}
\int_{-\ri\infty}^{\ri\infty}
P_{ns'}(nz,nx')\,Q_{nt}(nz,ny)\,\rd z
+
\frac{n}{2\pi \ri}
\int_{w_0}^{\overline{w_0}}
P_{ns'}(nz,nx')\,Q_{nt}(nz,ny)\,\rd z \\
&=
\bm 1(x'\ge y)
{n(s'-t)\choose n(x'-y)}
+
\frac{n}{2\pi \ri}
\int_{w_0}^{\overline{w_0}}
P_{ns'}(nz,nx')\,Q_{nt}(nz,ny)\,\rd z .
\end{split}
\end{align}
The second term on the right-hand side of \eqref{e:single_int2} satisfies the
homogeneous recursion by \Cref{c:recursion}, and therefore cancels after substitution
into \eqref{e:Aeq0}. Hence the only remaining contribution is the discrete heat
kernel, and the claim \eqref{e:Aeq0} follows from the following relation
\begin{align*}
&\phantom{{}={}}
\bm 1(s\ge t,\ x\ge y)
{n(s-t)\choose n(x-y)}
-
\bm 1\left(s-\frac1n\ge t,\ x\ge y\right)
{n(s-t)-1\choose n(x-y)} \\
&\quad
-
\bm 1\left(s-\frac1n\ge t,\ x-\frac1n\ge y\right)
{n(s-t)-1\choose n(x-y)-1}
=
\bm 1(s=t,\ x=y).
\end{align*}
Therefore \eqref{e:Aeq0} follows.
\end{proof}

\begin{proof}[Proof of \Cref{l:integrand_prop}]
We use the elementary identity
\begin{equation}\label{e:Gamma_shift_basic}
\frac{\Gamma(u)}{\Gamma(u+m)}
=
\begin{cases}
\displaystyle \frac{1}{u(u+1)\cdots(u+m-1)}, & m\in\mathbb Z_{\geq 0},\\[6pt]
\displaystyle (u+m)(u+m+1)\cdots(u-1), & m\in\mathbb Z_{<0}.
\end{cases}
\end{equation}
From the above identity, we have
 \(\Gamma(u)/\Gamma(u+m)\sim u^{-m}\) and
\begin{align}\label{e:gamma_ration}
\frac{\Gamma(n(y-z))}{\Gamma(n(x-z)+1)}
\sim C_1z^{n(y-x)-1}
,\quad 
\frac{\Gamma(n(z-(y-t)))}{\Gamma(n(z-(x-s))+1)}
\sim C_2 z^{n(x-y+t-s)-1}.
\end{align}
Multiplying these estimates and the constant factor
\(\Gamma(ns+1)/\Gamma(nt)\) gives
\[
P_{ns}(nz, nx) Q_{nt}(nz,ny)\sim  Cz^{n(t-s)-2},
\]
proving \eqref{e:asymptotic}.

For the first ratio in \eqref{e:gamma_ration}, set \(u=n(y-z)\). Then
$
n(x-z)+1=u+n(x-y)+1$.
If \(x\geq y\), then \(n(x-y)+1>0\), so \eqref{e:Gamma_shift_basic} gives
poles at \(u=0,-1,\cdots,-n(x-y)\), namely at
$
z=y,y+1/n,\cdots, x$.
If \(x<y\), the same identity gives zeros at
\(u=1,2,\cdots,n(y-x)-1\), namely at
$
z=y-1/n,y-2/n,\cdots, x+1/n$,
when this set is nonempty. This proves the second item.

For the second ratio, set \(v=n(z-(y-t))\). Then
$
n(z-(x-s))+1=v+n\bigl((y-t)-(x-s)\bigr)+1$.
Applying \eqref{e:Gamma_shift_basic} again gives the pole and zero sets listed
in the third item.

The final holomorphicity
statement follows directly from the listed pole and zero sets.
\end{proof}

\begin{proof}[Proof of \Cref{c:single_int}]

If $x<y$, by the second statement of \Cref{l:integrand_prop}, the integrand is holomorphic on $[y-t/2, \infty)$. We can deform the integral contour to $\{-\pi/2\leq \theta\leq \pi/2: Re^{\ri\theta}\}$ for sufficiently large $R$, and using \eqref{e:asymptotic} the integral is zero.

Now assume \(x\ge y\),  by the second statement of  \Cref{l:integrand_prop}, the integrand has poles at $x, x-1/n, x-2/n,\cdots, y$. %, and at $y-t, y-t-1/n, y-t-2/n,\cdots, x-s$. 
Thus we can deform the contour from $-\ri \infty$ to $\ri\infty$, to a contour enclosing $x, x-1/n, x-2/n,\cdots, y$ (but not $(-\infty, y-t/2]$) in clockwise.  Then we have
\begin{align}
&\phantom{{}={}}\frac{n}{2\pi \ri}\int_{-\ri\infty}^{\ri\infty}\frac{\Gamma(ns+1)}{\Gamma(nt)}\frac{\Gamma(n(y-z))\Gamma(n(z-(y-t)))}{\Gamma(n(x-z)+1)\Gamma(n(z-(x-s))+1)}  \rd z\\
&=\frac{n}{2\pi \ri}\int_{-\ri\infty}^{\ri\infty} \frac{\Gamma(ns+1)}{\Gamma(nt)}
\frac{1}{n(y-z) n(y+1/n-z)\cdots n(x-z)}
\frac{\Gamma(n(z-(y-t)))}{\Gamma(n(z-(x-s))+1)}
\rd z\\
&=\sum_{j=0}^{n(x-y)}\frac{\Gamma(ns+1)}{\Gamma(nt)}
\frac{(-1)^j}{j! (n(x-y)-j)!}\frac{\Gamma(nt+j)}{\Gamma(n(y-(x-s))+j+1)}\\
&=\sum_{j=0}^{n(x-y)}
\frac{\Gamma(ns+1)}{ (n(x-y)-j)!\Gamma(n(y-(x-s))+j+1)}\frac{(-1)^j\Gamma(nt+j)}{j!\Gamma(nt)}\\
&=\sum_{j=0}^{n(x-y)}{ns\choose n(x-y)-j}{-nt \choose j}={n(s-t)\choose n(x-y)},
% \frac{1}{n(z-(y-t)) n(z-(y-t-1/n)) \cdots n(z-(x-s))}  \rd z
\end{align}
where in the first statement we compute the ratio of two Gamma functions using $\Gamma(z+1)=z\Gamma(z)$; in the second statement, we compute the residual at $z=y+j/n$, namely the $j$-th term corresponds to the residual at $z=y+j/n$; in the third statement we rearrange the expression; the last statement follows from Chu–Vandermonde identity (see \cite[Section 1.1]{kuznetsov2006orthogonal}).
This proves \eqref{e:single_int}.
\end{proof}

\subsection{Modified factors}\label{s:factor}
When \((x,s)\) and \((y,t)\) both lie in the liquid region and stay a fixed
distance away from the arctic boundary and the ramification points, the
expression \eqref{e:bulk_ansatz} provides a natural ansatz for the
correlation kernel \(K^{-1}((x,s),(y,t))\). To extend this ansatz to
$
(x,s),(y,t)\in\fP\cap(\bZ^2/n),
$
we introduce some additional notation on the Riemann surface \(\cC\). We
recall the cuts \([b_i,a_i]\subset\cC(\bR)\), \(1\leq i\leq d\), from
\eqref{e:arcCR}.

\begin{enumerate}
\item\noindent\textbf{A basic factor.}
Away from the ramification points, \(\cC\) can be locally parametrized as
\((f(w),w)\). We define
\begin{equation}\label{e:defI}
I(w):=\exp\Bigl(-n\int_0^w \ln f(u)\,\rd u\Bigr).
\end{equation}
By \Cref{c:deflnf}, the integral in \eqref{e:defI}, and hence \(I(w)\), is
well defined on
$
\cC\setminus\cup_{i=1}^d[b_i,a_i].
$

\medskip

\item\noindent\textbf{Modification near the real locus.}
When \(w\) approaches the segment \([b_i,a_i]\subset\cC(\bR)\), it is
convenient to replace \(I(w)\) by a modified factor \(I_i(w)\) that
incorporates the local Gamma-type behavior. We set
\begin{align}
\begin{split}\label{e:defI2}
&I_i(w)
:=\exp\Bigl(-n\int_0^w \ln f(u)\,\rd u\Bigr)\,
\Gamma\!\bigl(n(b_i-w)+\tfrac12\bigr)\,
\Gamma\!\bigl(n(w-a_i)+\tfrac12\bigr)\\
&\quad\times
\exp\Bigl(-n\bigl((b_i-w)\ln(b_i-w)+(w-a_i)\ln(w-a_i)\bigr)
-n(b_i-a_i)\ln(n/e)-\ln(2\pi)\Bigr).
\end{split}
\end{align}

\item\noindent\textbf{Modification near ramification points.}
The shifted map \(\Psi_{\ft}\) from \eqref{e:emb_ft} sends a ramification
point \((x_0,t_0)\) to a point at which the projection
$
(f,w)\in\cC_{\ft}\longmapsto w\in\bC\bP
$
is regular. In a neighborhood of \(\Psi_{\ft}(x_0,t_0)\), we can
parametrize \(\cC_\ft\) as \((f_\ft(w),w)\); recall \eqref{e:def_ft}. Define
\begin{align}
I_\ft(w):=\exp\left(-n\int_0^w\ln f_\ft(u)\,\rd u\right).
\end{align}
Here \((0,0)\) is a fixed reference point, and the integral is taken along
any path from \((0,0)\) to \((f_\ft(w),w)\) in \(\cC_\ft\).
\end{enumerate}

Similarly to the bulk kernel ansatz \eqref{e:bulk_ansatz}, the general
kernel ansatz is given by the sum of a single-contour integral, which may
vanish, and a finite sum of double-contour integrals.

The single-contour integral term, if nonvanishing, is of one of the following
forms:
\begin{align}\label{e:single_term0}
\frac{n}{2\pi\ri}
\int_{\sfC}
P_{ns}(nz,nx)Q_{nt}(nz,ny)\,\rd z,
\qquad\text{or}\qquad
\frac{n}{2\pi\ri}
\int_{\sfC}
P_{n(s+\ft)}(nz,nx)Q_{n(t+\ft)}(nz,ny)\,\rd z,
\end{align}
where the contour \(\sfC\subset\bC\) depends on the location of \((x,s)\)
relative to \((y,t)\). We give the precise definition in
\Cref{s:single_integral}.

For the double-contour integrals, there are
then several cases, depending on whether \((x,s)\) or \((y,t)\) is close to or far from the
ramification points.
\begin{enumerate}
\item
Suppose that both \((x,s)\) and \((y,t)\) are bounded away from the
ramification points. The double-contour integral is a finite sum of terms of
the following form:
\begin{align}
\begin{split}\label{e:all_term0}
\frac{n}{(2\pi\ri)^2}
\int_{\sfC^{\rm a}}\!\!\int_{\sfC^{\rm d}}
P_{ns}(nw,nx)I_+(w)\,
Q_{nt}(nz,ny)I_-^{-1}(z)
\frac{\sqrt{\phi'(w)}\sqrt{\phi'(z)}}{\phi(w)-\phi(z)}
\,\rd w\,\rd z,
\end{split}
\end{align}
where each of \(I_+\) and \(I_-\) is either \(I\) or \(I_i\) for some
\(1\leq i\leq d\).

\item
Suppose that \((x,s)\) lies in a small neighborhood of a ramification
point. In the expression above, we replace
\[
P_{ns}(nw,nx), I_+(w), \phi(w), \phi'(w)\quad
\text{by}\quad 
P_{n(s+\ft)}(nw,nx),I_\ft(w),
\phi_\ft(w), \phi_\ft'(w),
\]
respectively. Similarly, if \((y,t)\) lies in a small neighborhood of a
ramification point, we replace
\[
Q_{nt}(nz,ny), I_-(z), \phi(z), \phi'(z)
\quad
\text{by}\quad 
Q_{n(t+\ft)}(nz,ny), I_\ft(z),
\phi_\ft(z), \phi_\ft'(z),
\]
respectively.

For example, if \((x,s)\) lies in a small neighborhood of a ramification
point but \((y,t)\) does not, the double-contour integral \(J^{(2)}\) is a
finite sum of terms of the following form:
\begin{equation}\label{e:J2_critical_bulk}
\frac{n}{(2\pi\ri)^2}
\int_{\sfC^{\rm a}}\!\!\int_{\sfC^{\rm d}}
P_{n(s+\ft)}(nw,nx)I_\ft(w)\,
Q_{nt}(nz,ny)I_-^{-1}(z)
\frac{\sqrt{\phi_\ft'(w)}\sqrt{\phi'(z)}}
{\phi_\ft(w)-\phi(z)}
\,\rd w\,\rd z,
\end{equation}
where \(I_-\) is either \(I\) or \(I_i\) for some \(1\leq i\leq d\).
\end{enumerate}
The integral contours \(\sfC^{\rm a}\) and \(\sfC^{\rm d}\), which depend
on the local geometry near \((x,s)\) and \((y,t)\), will be specified in
\Cref{s:critical_point}.

\subsection{Zeros and poles of the integrand}

In this section we collect some basic properties on the zeros and poles of the integrand in \eqref{e:all_term0}.

\begin{lemma}\label{l:properties_Ii}
Recall $I_i(w)$ from \eqref{e:defI2}.
\begin{enumerate}

\item The product
\begin{equation}\label{e:diffterm}
\exp\Bigl(-n\int_0^w \ln f(u)\,\rd u\Bigr)\,
\exp\Bigl(-n\bigl((b_i-w)\ln(b_i-w)+(w-a_i)\ln(w-a_i)\bigr)\Bigr)
\end{equation}
admits a single-valued analytic continuation in a sufficiently small neighborhood of $[b_i,a_i]$, except at $\infty_i$. 

\item Let $(-1,\infty_i)\in\cC$ correspond to the horizontal tangency point
$(x_i^{(0)},t_i^{(0)})\in\fA$. Then, as $w\to \infty_i$,
\[
e^{-n\int_0^w \ln f(u)\,\rd u -n\bigl((b_i-w)\ln(b_i-w)+(w-a_i)\ln(w-a_i)\bigr)}
=(C+\oo(1)) w^{-n(b_i-a_i-t_i^{(0)})},
\]
for a constant $C$ independent of $w$.
\end{enumerate}
\end{lemma}

\begin{lemma}\label{c:Iiproperty}
Recall that $na_i, n b_i\in \bZ'=\bZ+1/2$ from \eqref{e:aibi}.
\begin{enumerate}
\item
The following quantity
\begin{align}\begin{split}\label{e:Pterm}
&\phantom{{}={}}P_{ns}(nw,nx)  \Gamma(n(b_i-w)+1/2)\Gamma(n(w-a_i)+1/2)\\
&=\Gamma(ns+1)\frac{\Gamma(n(b_i-w)+1/2)\Gamma(n(w-a_i)+1/2)}{\Gamma(n(x-w)+1)\Gamma(n(w+s-x)+1)},
\end{split}\end{align}
is meromorphic in $w$, and its poles are given by
\begin{align}
\{b_i+1/(2n), b_i+1/(2n)+1/n,\cdots, x\}
\cup
\{x-s, x-s+1/n,\cdots, a_i-1/(2n)\}.
\end{align}
Here and below, an arithmetic progression is understood to be empty if its endpoint is not reached. Moreover, as $w\to\infty$,
\begin{align}\label{e:Pterm_limit}
\eqref{e:Pterm}=( C+\oo(1))w^{n(b_i-a_i-s)-1},
\end{align}
for a constant $C\neq 0$ independent of $w$.

\item
The following quantity
\begin{align}\begin{split}\label{e:Qterm}
&\phantom{{}={}}\frac{Q_{nt}(nz,ny)}{\Gamma(n(b_i-z)+1/2)\Gamma(n(z-a_i)+1/2)}\\
&=\frac{1}{\Gamma(nt)}\frac{\Gamma(n(y-z))\Gamma(n(z-(y-t)))}{\Gamma(n(b_i-z)+1/2)\Gamma(n(z-a_i)+1/2)},
\end{split}\end{align}
is meromorphic in $z$, and its poles are given by
\begin{align}
\{y, y+1/n, \cdots, b_i-1/(2n)\}
\cup
\{a_i+1/(2n), a_i+1/(2n)+1/n,\cdots, y-t\}.
\end{align}
Moreover, as $z\to\infty$ 
\begin{align}
\eqref{e:Qterm}=(C+\oo(1))z^{n(a_i-b_i+t)-1}
\end{align}
for a constant $C\neq 0$ independent of $z$.
\end{enumerate}

\end{lemma}

\begin{lemma}\label{c:PIproperty}
Recall from \eqref{e:aibi} that $na_i,nb_i\in \bZ'=\bZ+1/2$, and recall the
definition of $I_i(w)$ from \eqref{e:defI2}.

\begin{enumerate}
\item
In a sufficiently small neighborhood of $[b_i,a_i]$, the function
$P_{ns}(nw,nx) I_i(w)$ is meromorphic in $w$, and its poles are given by
\begin{align}\label{e:PI_poles}
\{b_i+1/(2n),\, b_i+1/(2n)+1/n,\,\cdots,\, x\}
\cup
\{x-s,\, x-s+1/n,\,\cdots,\, a_i-1/(2n)\}.
\end{align}
Here and below, an arithmetic progression is understood to be empty if its
endpoint is not reached. Moreover, let $(-1,\infty_i)\in\cC$ correspond to the horizontal tangency point
$(x_i^{(0)},t_i^{(0)})\in\fA$, as $w\to\infty_i$,
\begin{align}\label{e:limitPI}
P_{ns}(nw,nx) I_i(w)\frac{\sqrt{\phi'(w)}}{\phi(w)-u}
=
(C+\oo(1))w^{n(t_i^{(0)}-s)-2}, \quad u\in \{\phi(z), \phi_\ft(z)\}
\end{align}
for some constant $C\neq 0$ independent of $w$.

\item
In a sufficiently small neighborhood of $[b_i,a_i]$, the function
$Q_{nt}(nz,ny) I_i^{-1}(z)$ is meromorphic in $z$, and its poles are given by
\begin{align}\label{e:QI_poles}
\{y,\, y+1/n,\, \cdots,\, b_i-1/(2n)\}
\cup
\{a_i+1/(2n),\, a_i+1/(2n)+1/n,\,\cdots,\, y-t\}.
\end{align}
Moreover, let $(-1,\infty_i)\in\cC$ correspond to the horizontal tangency point
$(x_i^{(0)},t_i^{(0)})\in\fA$, as $z\to\infty_i$,
\begin{align}\label{e:limitQI}
Q_{nt}(nz,ny) I_i^{-1}(z)\frac{\sqrt{\phi'(z)}}{u-\phi(z)}
=
(C+\oo(1))z^{n(t-t_i^{(0)})-2}, \quad u\in \{\phi(w), \phi_\ft(w)\}
\end{align}
for some constant $C\neq 0$ independent of $z$.
\end{enumerate}
\end{lemma}

\begin{proof}[Proof of Lemma \ref{l:properties_Ii}]

By combining \eqref{e:intdiff} and \eqref{e:Imab}, the jump of the exponent across $[b_i,\infty_i)\cup (\infty_i, a_i]$ is zero modulo $2\pi$. More precisely, for $E\in [b_i,\infty_i)\cup (\infty_i, a_i]$,
\begin{align}
\left.\Im\left[-n \int_0^w \ln f(u)\,\rd u-n\bigl((b_i-w)\ln (b_i-w)+(w-a_i)\ln (w-a_i)\bigr)\right]\right|_{E-0\ri }^{E+0\ri}=0\mod 2\pi.
\end{align}
It follows that \eqref{e:diffterm} is analytic in a neighborhood of $[b_i,a_i]$, except at $\infty_i$.

For $u$ in a neighborhood of $\infty_i$, by \eqref{e:fzexp} we have $f(u)=-1+t_i^{(0)}/u+\OO(1/u^2)$. It follows that for $w$ in a neighborhood of $\infty_i$, with $\Im[w]>0$,
\begin{align}\begin{split}\label{e:intf_term}
&\phantom{{}={}}\int_0^w \ln f(u)\,\rd u
=\int_0^w \ln(-1+t_i^{(0)}/u+\OO(1/u^2))\,\rd u\\
&=\int_0^w \ln(1-t_i^{(0)}/u+\OO(1/u^2))\,\rd u -w\pi \ri
=-t_i^{(0)}\ln w-w\pi\ri+\OO(1),
\end{split}\end{align}
and
\begin{align}\label{e:bwterm}
(b_i-w)\ln (b_i-w)+(w-a_i)\ln (w-a_i)
=(b_i-a_i)\ln w -(b_i-w)\pi \ri +\OO(1).
\end{align}
Thus, combining \eqref{e:intf_term} and \eqref{e:bwterm}, we obtain
\begin{align}
e^{-n \int_0^w \ln f(u)\,\rd u }e^{-n((b_i-w)\ln (b_i-w)+(w-a_i)\ln (w-a_i) )}
=(C+\oo(1)) w^{-n(b_i-a_i-t_i^{(0)})}.
\end{align}

\end{proof}

\begin{proof}[Proof of \Cref{c:Iiproperty}]
The proof follows the same argument as that of \Cref{l:integrand_prop}, so we omit.
\end{proof}

\begin{proof}[Proof of \Cref{c:PIproperty}]
The first claim \eqref{e:PI_poles} for $P_{ns}(nw,nx) I_i(w)$ follows from
\eqref{e:diffterm}, the first statement of \Cref{c:Iiproperty}, and the decomposition
\begin{align}\label{e:PI_decomp}
P_{ns}(nw,nx) I_i(w)
=
 \eqref{e:diffterm}\times
 \eqref{e:Pterm}.
\end{align}
For the second claim \eqref{e:limitPI}, we note that when $w\rightarrow \infty$,
either $\phi(w)\rightarrow \infty$, in which case
$\phi(w)=\fa_i w+\OO(1)$ with some $\fa_i\neq 0$; or $\phi(w)\rightarrow p_i\neq u$,
in which case
$\phi(w)=p_i+\fa_i/w+\OO(1/w^2)$ with some $\fa_i\neq 0$. In both cases, when $w\rightarrow \infty_i$, by \eqref{e:inftoinf} and \eqref{e:inftofinite}
\begin{align}
\frac{\sqrt{\phi'(w)}}{\phi(w)-u}=(C+\oo(1))\frac{1}{w}.
\end{align}
Thus the claim \eqref{e:limitPI} follows from \eqref{e:Pterm_limit}, the first statement of \Cref{c:Iiproperty}, and the decomposition \eqref{e:PI_decomp}.

The claims for $Q_{nt}(nz,ny) I_i^{-1}(z)$ follow from the analogous
decomposition, so we omit.
\end{proof}

\section{Charts and Contours}\label{s:critical_point}
In \Cref{s:critical_bulk}, we constructed the liquid chart and the
corresponding local descent/ascent paths, which were then used to formulate
the kernel ansatz \eqref{e:bulk_ansatz} in the liquid region. In this section,
we construct analogous charts for the remaining types of distinguished points:
ramification points, arctic boundary points, cusp points, frozen points, tangent points,
cusp-turning points, and tangent frozen points. We also associate to each of these charts suitable local
descent/ascent paths.

\subsection{Ramification chart}\label{s:ramification}
We recall from the first statement in \Cref{p:surface} that there are only finitely
many ramification points \((x_0,s_0)\in \fL\), characterized by
$\del_x \chi(x_0,s_0)={1}/{s_0}$.
Fix any small \(\ft>0\). Recall the shifted polygon
\(\fP_{\ft}\) from \eqref{e:shifted_polygon}, the function $\chi_\ft(x,s+\ft)=\chi(x,s)$ from \eqref{e:chit}, the associated Riemann surface
\(\cC_{\ft}\) from \eqref{e:shifted_RS},  and the shifted tiling action
\(S_\ft(w;x_0,s_0+\ft)\) from \eqref{e:def_action_shift}. Let
\begin{align}
w_0
:=
x_0-(s_0+\ft)\chi_\ft(x_0,s_0+\ft)
=
x_0-(s_0+\ft)\chi(x_0,s_0).
\end{align}
Then \(w_0\) is a critical point of \(S_\ft(\cdot;x_0,s_0+\ft)\), and
$S_\ft'(w_0;x_0,s_0+\ft)=0$.
From \eqref{e:def_ft}, the point \((f_0,w_0)\in \cC_\ft\), where
\(f_0:=f(x_0,s_0)\), is bounded away from the ramification points of
\(\cC_\ft\), and we can parametrize \(\cC_\ft\) locally as \(f=f_\ft(w)\).

In this case \(S_\ft''(w_0;x_0,s_0+\ft)\) is bounded away from \(0\), and the
steepest-descent direction is determined by
\begin{align}
2\theta_0
:=
\pi-\arg S_\ft''(w_0;x_0,s_0+\ft),
\qquad
e^{2\ri\theta_0}
:=
-\frac{|S_\ft''(w_0;x_0,s_0+\ft)|}
{S_\ft''(w_0;x_0,s_0+\ft)}.
\end{align}

After this shift of time, \Cref{c:bulk} still holds after replacing \(S\) by
\(S_\ft\) and \(w_c=x-s\chi(x,s)\) by
$
w_{c,\ft}=x-(s+\ft)\chi_\ft(x,s+\ft)=x-(s+\ft)\chi(x,s)$.

\begin{lemma}[$\fc$-ramification chart]\label{c:ramification}
Let $(x_0,s_0)\in\fL$ be a ramification point. Then for every sufficiently
small $\fc>0$ there exists $\delta=\delta(\fc)>0$, also sufficiently
small, such that the following holds.

Let $w_0=x_0-(s_0+\ft)\chi_\ft(x_0,s_0+\ft)$ be the critical point of $S_\ft(\cdot;x_0,s_0+\ft)$ corresponding to $(x_0,s_0)$, and set
\[
\fU_\ft:=\{w\in\bC:\ |w-w_0|\le \fc\}.
\]
On $\fU_\ft$ the Riemann surface $\cC_\ft$ can be parametrized as $(f_\ft(w),w)$. We will therefore identify
$\fU_\ft$ with its image in $\cC_\ft$ under the map $w\mapsto (f_\ft(w),w)$.  Moreover, for all $w\in\fU_\ft$,
\begin{align}
\Im[w], -\Im[\chi_\ft(w)]\asymp 1.
\end{align}

Now let $(x,s)\in\fP$ satisfy $\|(x,s)-(x_0,s_0)\|_2\le \delta$. Then the tiling action $S_\ft(\,\cdot\,;x,s+\ft)$ has exactly one critical point $w_{c,\ft}$ inside $\fU_\ft$, satisfying
\begin{align}
|w_{c,\ft}-w_0|
\lesssim \|(x,s)-(x_0,s_0)\|_2
\leq \delta.
\end{align}
Moreover, for every $w\in\fU_\ft$,
\begin{align}
S_\ft(w;x,s+\ft)-S_\ft(w_{c,\ft};x,s+\ft)=d\,(w-w_{c,\ft})^2+\cE(w),
\qquad d:=\frac12\,S_\ft''(w_0;x_0,s_0+\ft),\quad  |d|\asymp 1.
\end{align}
The error term satisfies
\begin{align}
|\cE(w)|\leq C\bigl(\delta\ln(1/\delta)+|w-w_{c,\ft}|^3\bigr)\le \frac{|d|\,\fc^2}{100},\quad |d|\asymp 1.
\end{align}

In this situation, we call $\fU_\ft$ a $\fc$-ramification chart (centered at $w_0$), and we say that the point
 $(x,s)$ is \emph{adapted} to $\fU_\ft$. 
We also introduce the corresponding chart in \(\cC\) by pulling back
\(\fU_\ft\subset\cC_\ft\) under the transport map \eqref{e:maptoCt}:
\[
\fU
:=
\left\{
\left(f,w+\frac{\ft f}{f+1}\right):(f,w)\in\fU_\ft
\right\}.
\]

\end{lemma}

Similarly to \eqref{e:gamma_path} and \eqref{e:conj_gamma_path}, we introduce the
local paths through \(w_0\) and \(\overline{w_0}\) as
\begin{align}\label{e:critical_gamma_path}
\begin{split}
\sfC^{\rm d}(w_0)
&:=
\bigl\{\,w_0+r e^{\ri\theta_0}:\ -\fc\le r\le \fc\,\bigr\},
\qquad
\sfC^{\rm a}(w_0)
:=
\bigl\{\,w_0-\ri r e^{\ri\theta_0}:\ -\fc\le r\le \fc\,\bigr\},\\
\sfC^{\rm d}(\overline{w_0})
&:=
\bigl\{\,\overline{w_0}+r e^{-\ri\theta_0}:\ -\fc\le r\le \fc\,\bigr\},
\qquad
\sfC^{\rm a}(\overline{w_0})
:=
\bigl\{\,\overline{w_0}+\ri r e^{-\ri\theta_0}:\ -\fc\le r\le \fc\,\bigr\},
\end{split}
\end{align}
for a sufficiently small \(\fc>0\).

In the following we orient the local paths \eqref{e:critical_gamma_path}, using the
same convention as that of the contours in the liquid region introduced in
\Cref{s:critical_bulk}. We recall $\phi_\ft$ from \eqref{e:defphi_t} and \eqref{e:defphi_t2}, $\chi_\ft(w)=f_\ft(w)/(f_\ft(w)+1)$,
and recall the factorization for descent vector from \eqref{eq:factor}:
\begin{align}\label{e:sqrtSder}
\begin{split}
v^{\rm d}(w_0)=\sqrt{-\frac{1}{S_\ft''(w_0;x_0,s_0+\ft)}}
&=
(s_0+\ft)^{1/2}
\sqrt{\chi_\ft(w_0)}\sqrt{1-\chi_\ft(w_0)}
\cdot
\frac{\sqrt{\partial_x\phi_\ft(x_0,s_0+\ft)}}
{\sqrt{\phi_\ft'(w_0)}} \\
&=
(s_0+\ft)^{1/2}
\sqrt{\chi_\ft(w_0)}\sqrt{1-\chi_\ft(w_0)}
\cdot
\frac{\sqrt{\partial_x\phi(x_0,s_0)}}
{\sqrt{\phi_\ft'(w_0)}},
\end{split}
\end{align}
where in the second line we used $
\partial_x\phi_\ft(x_0,s_0+\ft)=\partial_x\phi(x_0,s_0)$.
After fixing the global branch of \(\sqrt{\partial_x\phi(x,s)}\), using the principal branches of $\sqrt{\chi_\ft(w_0)}$ and $\sqrt{1-\chi_\ft(w_0)}$  in
\eqref{e:sqrtSder}, the local choice of
$
\sqrt
{\phi_\ft'(w_0)}
$
 gives the orientation of \(\sfC^{\rm d}(w_0)\). Reversing the branch reverses the contour orientation, and the resulting sign changes cancel in the double-contour integral over the ramification chart (see \eqref{e:J2_critical_bulk}), provided that the prefactor $\sqrt{\phi_\ft'(w)}$ is defined using the same branch convention.

The direction of \(\sfC^{\rm d}(\overline{w_0})\) is chosen to be the complex
conjugate of the direction of \(\sfC^{\rm d}(w_0)\), with an additional minus sign: $v_{\rm d}(\overline{w_0})
:=
-\overline{v_{\rm d}(w_0)}$.
Then we orient \(\sfC^{\rm a}(w_0)\) by multiplying the direction of
\(\sfC^{\rm d}(w_0)\) by \(-\ri\), and we orient
\(\sfC^{\rm a}(\overline{w_0})\) by multiplying the direction of
\(\sfC^{\rm d}(\overline{w_0})\) by \(\ri\).

\begin{lemma}\label{c:ramification_steepest}
Adopt the assumptions and notation of \Cref{c:ramification}. Then there exists a
constant $\fc'>0$ such that the following statements hold.
\begin{itemize}
\item \noindent\emph{Replacing \(\mathsf C^{\rm d}(w_0)\) by \(\mathsf D^{\rm d}(w_{c,\ft})\).}
We define the local steepest-descent set $\mathsf D^{\rm d}(w_{c,\ft})$ at \(w_{c,\ft}\)
to be the portions of the steepest-descent trajectories of
$S_\ft(\,\cdot\,;x,s+\ft)$ issuing from \(w_{c,\ft}\) up to their first exit from the disk
\[
\{\,w: |w-w_0|\le \fc\,\}.
\]
Then $\mathsf D^{\rm d}(w_{c,\ft})$ has total length
$\OO(1)$. Moreover, $\sfC^{\rm d}(w_0)$ can be deformed to
$\mathsf D^{\rm d}(w_{c,\ft})$, together with finitely many arcs of total length
$\OO(1)$, on which
\begin{align}
n\,\Re[ S_\ft(w;x,s+\ft)-S_\ft(w_{c,\ft};x,s+\ft)]\le -\,n\fc'.
\end{align}

\item \noindent\emph{Replacing \(\mathsf D^{\rm d}(w_{c,\ft})\) by \(\mathsf S^{\rm d}(w_{c,\ft})\).}
Let
$
r_n:={\ln n}/{\sqrt n}$.
We define the truncated steepest-descent set $\mathsf S^{\rm d}(w_{c,\ft})$ at
\(w_{c,\ft}\) to be the portions of the steepest-descent trajectories issuing from
\(w_{c,\ft}\) up to their first exit from the disk
\[
\{\,w: |w-w_{c,\ft}|\le r_n\,\}.
\]
Then, $\sfS^{\rm d}(w_{c,\ft})$ has total length $\OO(r_n)$, and  for all \(w\in \mathsf D^{\rm d}(w_{c,\ft})\setminus
\mathsf S^{\rm d}(w_{c,\ft})\),
\[
e^{n\Re[S_\ft(w;x,s+\ft)]}
\le
e^{n\Re[S_\ft(w_{c,\ft};x,s+\ft)]}e^{-\fc'(\ln n)^2}.
\]
\end{itemize}
We define the analogous ascent sets
\(\mathsf D^{\rm a}(w_{c,\ft})\) and \(\mathsf S^{\rm a}(w_{c,\ft})\). The same
statements hold with \(S_\ft\) replaced by \(-S_\ft\) and with the superscript
\({\rm d}\) replaced by \({\rm a}\).
\end{lemma}

The proofs of \Cref{c:ramification} and \Cref{c:ramification_steepest}
are identical to those of \Cref{c:bulk} and \Cref{c:bulk_steepest},
respectively, after replacing \(S(\cdot;x,s)\) by
\(S_\ft(\cdot;x,s+\ft)\), so we omit them.

\subsection{Arctic chart}
\label{s:critical_arctic}

Fix any $(x_0,s_0)\in \fA$ bounded away from tangent locations and cusp locations. Define \(w_0\) by the characteristic relation
\begin{align}
w_0=x_0-s_0 \chi(x_0;s_0)=x_0-s_0 \chi(w_0)\in \bR.
\end{align}
Then $w_0$ is a critical point, and $S'(w_0;x_0,s_0)=0$.
By the second statement in \Cref{p:surface}, we have
\begin{align}\label{e:chi_property}
\chi'(w_0)=-1/s_0, \quad |\chi(w_0)|, |1-\chi(w_0)|,  |\chi''(w_0)|\asymp 1,
\end{align}
and from \eqref{e:derSthird}
\begin{align}\label{e:derSthird2}
S'(w_0;x_0,s_0)=S''(w_0;x_0,s_0)=0,\quad
S'''(w_0;x_0,s_0)=-\frac{\chi''(w_0)}{\chi(w_0)(1-\chi(w_0))}\neq 0.
\end{align}
The sign of  \(S'''(w_0;x_0,s_0)\) is determined by the signs of $\chi''(w_0)$ and $\chi(w_0)(1-\chi(w_0))$, as classified in \Cref{f:arctic_boundary}.
Locally around $w_0$,  $\chi(w)$ satisfies
 \begin{align}
 |\chi(w)|, |1-\chi(w)|\asymp 1, \quad |\chi(w)|, |\chi'(w)|, |\chi''(w)|\lesssim 1,
 \end{align} 
 and 
 $S(w;x_0,s_0)$ is holomorphic and we have the Taylor expansion
\begin{align}\label{e:arctic_neighborhood}
&S(w;x_0,s_0)=S(w_0;x_0,s_0)+S'''(w_0; x_0, s_0)\frac{(w-w_0)^3}{6} +\OO(|w-w_0|^4).
\end{align}

The following lemma introduces the arctic chart and records several of its properties. Its proof, based on a Taylor expansion of the tiling action, is deferred to \Cref{s:arctic_chart_proof}.
\begin{lemma}[$\fc$-arctic chart]\label{c:arctic_critical}
Fix $\delta_0>0$ and let $(x_0,s_0)\in\fL$ be a point whose distance is at least $\delta_0$
from all cusp and tangency points. Then for every sufficiently small $\fc>0$
(depending only on $\delta_0$) there exists $\delta=\delta(\fc)>0$, also sufficiently small
and depending only on $\delta_0$, such that the following holds.

Let $w_0=x_0-s_0\,\chi(x_0,s_0)$ be the critical point corresponding to $(x_0,s_0)$, and set
\[
\fU:=\{w\in\bC:\ |w-w_0|\le \fc\}.
\]
On $\fU$ the Riemann surface $\cC$ can be parametrized as $(f(w),w)$. We will therefore identify
$\fU$ with its image in $\cC$ under the map $w\mapsto (f(w),w)$. Moreover, for all $w\in\fU$,
\begin{align}\label{e:arctic_bounds}
|\chi(w)|,\ |1-\chi(w)|,\ |\chi''(w)|\asymp 1,
\end{align}
with implicit constants depending only on $\delta_0$.

Now let $(x,s)\in\fP$ satisfy $\|(x,s)-(x_0,s_0)\|_2\le \delta$.
Then the tiling action $S(\,\cdot\,;x,s)$ has exactly two critical points $w_c$ in $\fU$, satisfying
\begin{align}
|w_c-w_0|
\lesssim \|(x,s)-(x_0,s_0)\|^{1/2}_2
\leq \sqrt\delta.
\end{align}

Moreover, for each such critical point $w_c$ and every $w\in\fU$,
\begin{align}\label{e:arctic_cubic}
S(w;x,s)-S(w_c;x,s)=d\,(w-w_0)^3+\cE(w),
\quad d:=\frac{1}{6}\,S'''(w_0;x_0,s_0), \quad |d|\asymp 1.
\end{align}
The error term satisfies
\begin{align}\label{e:arctic_err}
|\cE(w)|\leq C(\delta\ln(1/\delta)+|w-w_0|^4)
\le \frac{|d|\,\fc^3}{100}.
\end{align}

In this situation, we call $\fU$ a $\fc$-arctic chart (centered at $w_0$), and we say that
the point $(x,s)$ is \emph{adapted} to $\fU$.
\end{lemma}

We recall that in Panel (A), (B) and (C) of \Cref{f:arctic_boundary}, $S'''(w_0;x_0, s_0)<0$;  in Panel (D), (E) and (F) of \Cref{f:arctic_boundary}, $S'''(w_0;x_0, s_0)>0$.
Depending on the sign of $S'''(w_0; x_0, s_0)$, there are two cases for the descent and ascent paths 
\begin{enumerate}
\item If $S'''(w_0;x_0,s_0)<0$, then the steepest-descent directions at $w_0$ are $1$ and $e^{\pm 2\ri \pi/3}$. We introduce the following local paths, consisting of the non-real descent and ascent rays
\begin{align}\label{e:arctic_contour2}
\sfC^{\rm d}(w_0)= \{w_0+r e^{\pm 2\ri \pi/3}: 0\leq r\leq \fc\},\quad  \sfC^{\rm a}(w_0)= \{w_0+r e^{\pm \ri \pi/3}: 0\leq r\leq \fc\}.
\end{align}
\item If $S'''(w_0; x_0, s_0)>0$, then the steepest-descent directions at $w_0$ are $-1$ and $e^{\pm \ri \pi/3}$. We introduce the following local paths, consisting of the non-real descent and ascent rays
\begin{align}\label{e:arctic_contour1}
\sfC^{\rm d}(w_0)= \{w_0+r e^{\pm \ri \pi/3}: 0\leq r\leq \fc\},\quad \sfC^{\rm a}(w_0)= \{w_0+r e^{\pm 2\ri \pi/3}: 0\leq r\leq \fc\}.
\end{align}
\end{enumerate}
Note that \(\sfC^{\rm d}(w_0)\) and  \(\sfC^{\rm a}(w_0)\)  are contained in \(\fU\).

In the following we orient the local paths \eqref{e:arctic_contour1} and \eqref{e:arctic_contour2}, so they are compatible with the contours in the liquid region as introduced in \Cref{s:critical_bulk}. For $(x,s)\in \fL$ close to $(x_0, s_0)$, let $w_c=x-s\chi(w_c)\in\bC_+$,  we recall the factorization for descent vector from \eqref{eq:factor}:
\begin{align}
v^{\rm d}(w_c)=\sqrt{-\frac{1}{S''(w_c;x,s)}}
=\frac{\sqrt{\chi(w_c)}\sqrt{1-\chi(w_c)}}{\sqrt{1/s+\chi'(w_c)}}
=s^{1/2}\sqrt{\chi(w_c)}\sqrt{1-\chi(w_c)}\cdot
\frac{\sqrt{\partial_x\phi(x,s)}}{\sqrt{\phi'(w_c)}}.
\end{align}
Thus, a local choice of the branch of $\sqrt{\phi'(w_c)}$ determines the
oriented descent direction $v^{\rm d}(w_c)$, and conversely. By analytic
continuation, this choice also determines the branch of the prefactor
$\sqrt{\phi'(w)}$ appearing in the double-contour integral
\eqref{e:all_term0}.

We next give a geometric description of this choice as $(x,s)$ approaches
$(x_0,s_0)$ and $w_c$ approaches $w_0$.  Let $w_c=w_0+a+\ri b$, then $b>0$ and \eqref{e:curve_reg} gives
\begin{align}
1/s+\chi'(w_c)=\chi''(w_0)b\ri +\OO(b(|a|+b)), \quad \chi(w_c)=\chi(w_0)+\OO(|a|+b)
\end{align}
Thus as $(x,s)\rightarrow (x_0, s_0)$ the direction of $v^{\rm d}(w_c)=\sqrt{-1/S''(w_c;x,s)}$ is given by
\begin{align}\label{e:arc_direction}
%\sqrt{\frac{\chi(w_c)\bigl(1-\chi(w_c)\bigr)}{\chi''(w_0) \ri}}\approx
(1+\oo(1))\frac{\sqrt{\chi(w_0)}\sqrt{1-\chi(w_0)}}{\sqrt{b\chi''(w_0) \ri}}= \frac{(1+\oo(1))}{\sqrt{b}}\sqrt{\frac{\ri}{S'''(w_0;x_0, s_0)}},
\end{align}
where the second equality follows from \eqref{e:derSthird2}. Since
$\sqrt{b}>0$, the limiting direction of $v^{\rm d}(w_c)$ is determined by
the square root in the last expression.

\begin{figure}
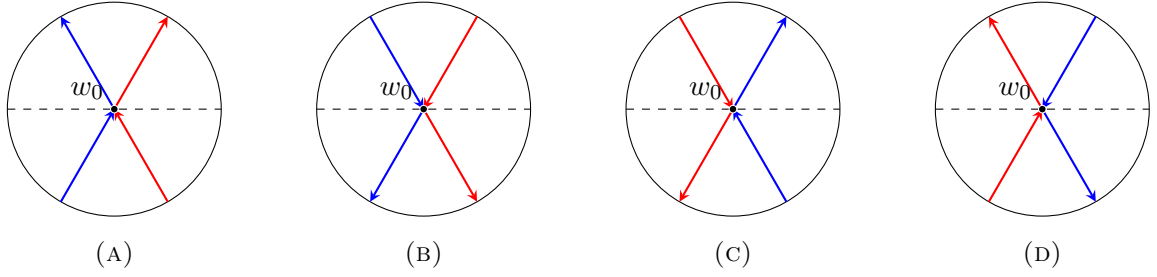

		
			\begin{subfigure}{0.24\textwidth}
	
		\begin{center}	

		\end{center}
		\caption{}
		\end{subfigure}	
		\begin{subfigure}{0.24\textwidth}
	
		\begin{center}	
		%
		\end{center}
	\caption{}
		\end{subfigure}			
			\begin{subfigure}{0.24\textwidth}
	
		\begin{center}	
		%
		\end{center}
	\caption{}
		\end{subfigure}		
			\begin{subfigure}{0.24\textwidth}
	
		\begin{center}	
		%
		\end{center}
	\caption{}
		\end{subfigure}

				\caption{\label{f:edge_path}  
Local paths associated with arctic charts.
}
\end{figure}

If $S'''(w_0;x_0,s_0)<0$, choosing the limiting direction
$e^{3\pi\ri/4}$ gives the blue paths in panel~(A) of
\Cref{f:edge_path}, whereas choosing the opposite direction
$-e^{3\pi\ri/4}$ gives the blue paths in panel~(B).

If $S'''(w_0;x_0,s_0)>0$, choosing the limiting direction
$e^{\pi\ri/4}$ gives the blue paths in panel~(C) of
\Cref{f:edge_path}, whereas choosing the opposite direction
$-e^{\pi\ri/4}$ gives the blue paths in panel~(D).

We orient the ascent contour $\sfC^{\rm a}(w_0)$ as indicated by the red paths in
\Cref{f:edge_path}. With this convention, in the upper half-plane the ordered pair
$
\bigl(\sfC^{\rm d}(w_0),\sfC^{\rm a}(w_0)\bigr)$
intersects negatively at \(w_0\), while in the lower half-plane they intersect positively at \(w_0\), in the sense of
\Cref{d:positive_negative}.

\begin{figure}
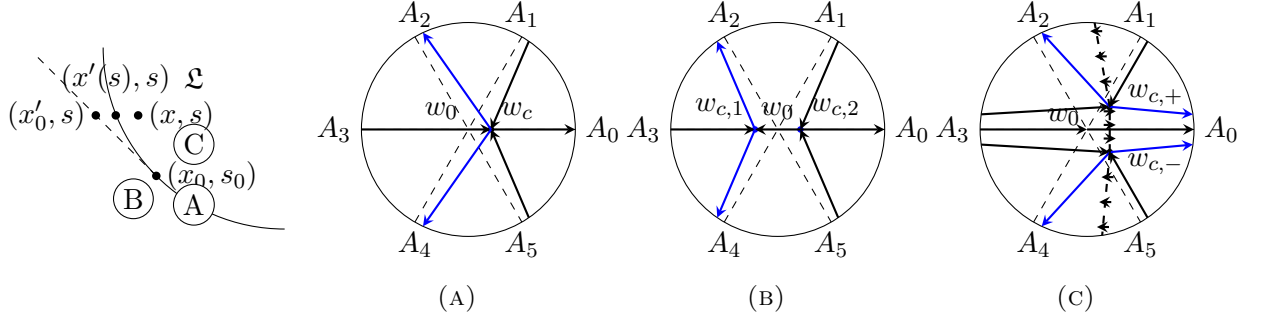

	\begin{subfigure}[t]{0.25\textwidth}
	
		\centering

		\caption{}
		\end{subfigure}		
			\begin{subfigure}[t]{0.24\textwidth}
	
			%
		
\caption{}
	\end{subfigure}
		\begin{subfigure}[t]{0.24\textwidth}
	
			%

	\caption{}
	\end{subfigure}

\caption{\label{f:cubic_saddle}
Gradient flow of \(S(\cdot;x,s)\) in an arctic chart, with
\(d=S'''(w_0;x_0,s_0)/6<0\).
}
\end{figure}

The following lemma show that the local descent and ascent paths can be deformed into steepest-descent and steepest-ascent paths with negligible error. Its proof is deferred to \Cref{s:arctic_chart_proof}.
\begin{lemma}\label{c:arctic_steepest}
Adopt the assumptions and notation of \Cref{c:arctic_critical}. Then there exists a
constant $\fc'>0$ such that the following statements hold.

\begin{itemize}
\item \noindent\emph{Replacing \(\mathsf C^{\rm d}(w_0)\) by \(\mathsf D^{\rm d}(w_c)\).}
For each critical point \(w_c\), we define the local steepest--descent path
\(\mathsf D^{\rm d}(w_c)\) at \(w_c\) to be the portions of the non-real
steepest--descent trajectories of \(S(\,\cdot\,;x,s)\) starting at \(w_c\),
up to their first exit from the disk
\[
\{w: |w-w_0|\le \fc\};
\]
see \Cref{f:cubic_saddle}. If \(w_c\) is not a descent critical point, we set
\(\mathsf D^{\rm d}(w_c)=\emptyset\). The length of
\(\mathsf D^{\rm d}(w_c)\) is \(\OO(1)\).

There are two cases:
\begin{enumerate}
  \item If \((x,s)\in\fL\), then there are two complex conjugate critical
  points. We can deform \(\mathsf C^{\rm d}(w_0)\) to the union of the
  steepest--descent paths \(\mathsf D^{\rm d}(w_c)\) associated with these two
  critical points, together with several arcs of total length \(\OO(1)\).

  \item If \((x,s)\in\fP\setminus \fL\), then there are two real critical
  points, counted with multiplicity. There is one descent critical point
  \(w_c\) and one ascent critical point \(w_c'\). We can deform
  \(\mathsf C^{\rm d}(w_0)\) to the steepest--descent path
  \(\mathsf D^{\rm d}(w_c)\), together with several arcs of total length
  \(\OO(1)\).
\end{enumerate}
In both cases, on these extra arcs,
\[
n\,\Re\!\bigl[S(w;x,s)-S(w_c;x,s)\bigr]\le -\,n\fc'.
\]

\item \noindent\emph{Replacing \(\mathsf D^{\rm d}(w_c)\) by \(\mathsf S^{\rm d}(w_c)\).}
Let \((x'(s),s)\in\fA\) denote the nearby point on the arctic curve with second coordinate \(s\), then $|x-x'(s)|\asymp \dist((x,s),\fA)$;  see the first
panel of \Cref{f:cubic_saddle}. We distinguish two regimes.

\begin{enumerate}
\item \emph{Close to the arctic boundary.}
If
$
|x-x'(s)|\leq {(\ln n)^2}/{n^{2/3}},
$
then, for each critical point \(w_c\), we set
$
r_n:={(\ln n)^2}/{n^{1/3}}.
$

\item \emph{Liquid and frozen regimes.}
If
$
|x-x'(s)|\geq {(\ln n)^2}/{n^{2/3}},
$
then, for each critical point \(w_c\), we set
$
r_n:={\ln n}/({n^{1/2}|x-x'(s)|^{1/4}})\ll\sqrt{|x-x'(s)|}.
$
\end{enumerate}

For each descent critical point \(w_c\), define the truncated local
steepest--descent path \(\mathsf S^{\rm d}(w_c)\) to be the portion of
\(\mathsf D^{\rm d}(w_c)\) starting at \(w_c\) and stopped at its first exit
from the disk
\[
\{w: |w-w_c|\leq r_n\}.
\]
Then, for all
$
w\in \mathsf D^{\rm d}(w_c)\setminus \mathsf S^{\rm d}(w_c),
$
we have
\[
e^{n\Re[S(w;x,s)]}
\leq
e^{n\Re[S(w_c;x,s)]}e^{-\fc'(\ln n)^2}.
\]
\end{itemize}

The analogous statements also hold for the local steepest--ascent paths and
their truncations.

Finally, we record the following estimate for later use. If
\((x,s)\in\fP\setminus \fL\), then there is one descent critical point \(\xi_c^{\rm d}\)
and one ascent critical point \(\xi_c^{\rm a}\), and
\begin{align}\label{e:Sdiff_at_critical}
\Re\!\bigl[S(\xi_{c}^{\rm a};x,s)-S(\xi^{\rm d}_c;x,s)\bigr]
\asymp |x-x'(s)|^{3/2}.
\end{align}

\end{lemma}

\begin{figure}
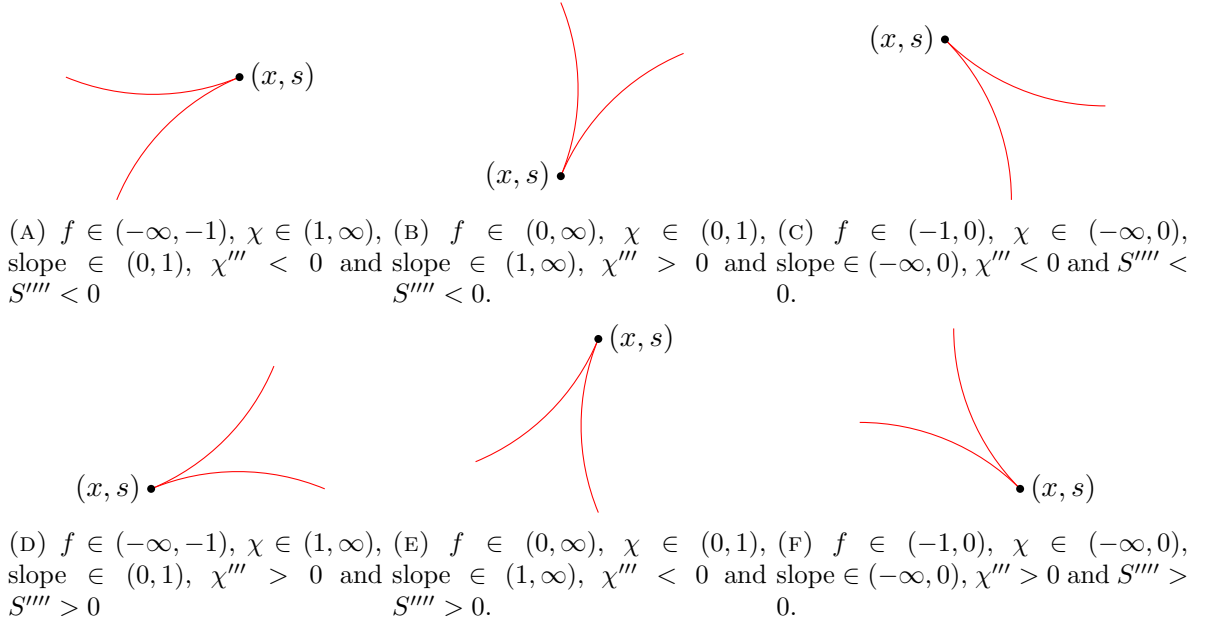
  
    \begin{subfigure}{0.3\textwidth}
    \centering
      %
      \caption{$f\in(-\infty,-1)$, $\chi\in (1,\infty)$, $\text{slope}\in (0,1)$, $\chi'''<0$ and $S''''<0$}
  \end{subfigure}
    \begin{subfigure}{0.3\textwidth}
    \centering
      %
       \caption{$f\in(0,\infty)$, $\chi\in (0,1)$, $\text{slope}\in (1,\infty)$, $\chi'''>0$ and $S''''<0$.}
  \end{subfigure}
 \begin{subfigure}{0.33\textwidth}
    \centering
      %
       \caption{$f\in(-1,0)$, $\chi\in (-\infty,0)$, $\text{slope}\in (-\infty,0)$, $\chi'''<0$ and $S''''<0$.}
  \end{subfigure}

 \begin{subfigure}{0.3\textwidth}
    \centering
      %
       \caption{$f\in(-\infty,-1)$, $\chi\in (1,\infty)$, $\text{slope}\in (0,1)$, $\chi'''>0$ and $S''''>0$}
  \end{subfigure}
  \begin{subfigure}{0.3\textwidth}
    \centering
      %
      \caption{$f\in(0,\infty)$, $\chi\in (0,1)$, $\text{slope}\in (1,\infty)$, $\chi'''<0$ and $S''''>0$.}
  \end{subfigure}
    \begin{subfigure}{0.33\textwidth}
    \centering
      %
      \caption{$f\in(-1,0)$, $\chi\in (-\infty,0)$, $\text{slope}\in (-\infty,0)$, $\chi'''>0$ and $S''''>0$.}
  \end{subfigure}

	 \caption{Cusp locations}
	 \label{f:cusp}
	 \end{figure}

\subsection{Cusp chart}
\label{s:critical_cusp}

For any cusp $(x_0,s_0)\in \fA$ which is not a cusp-turning point. Define \(w_0\) by the characteristic relation
\begin{align}
w_0=x_0-s_0 \chi(x_0;s_0)=x_0-s_0 \chi(w_0)\in \bR.
\end{align}
Then $S'(w_0;x_0,s_0)=0$, and $w_0$ is a critical point.
By the second statement in \Cref{p:surface}, we have
\begin{align}\label{e:cusp_expandchi}
\chi'(w_0)=-1/s_0,  \quad \chi''(w_0)=0,\quad |\chi(w_0)|, |1-\chi(w_0)|, |\chi'''(w_0)|\asymp 1,
\end{align}
and from \eqref{e:derSfourth}
\begin{align}\label{e:derSfourth2}
S'(w_0;x_0,s_0)=S''(w_0;x_0,s_0)=S'''(w_0;x_0,s_0)=0,\quad
S''''(w_0;x_0,s_0)=-\frac{\chi'''(w_0)}{\chi(w_0)(1-\chi(w_0))}\neq 0.
\end{align}
The sign of  \(S''''(w_0;x_0,s_0)\) is determined by the signs of $\chi'''(w_0)$ and $\chi(w_0)(1-\chi(w_0))$, as classified in \Cref{f:cusp}.
Locally around $w_0$, $S(w;x_0,s_0)$ is holomorphic and we have the Taylor expansion
\begin{align}\label{e:cusp_neighborhood}
&S(w;x_0,s_0)=S(w_0;x_0,s_0)+S''''(w_0; x_0, s_0)\frac{(w-w_0)^4}{24} +\OO(|w-w_0|^5),
\end{align}

The following lemma introduces the cusp chart and records several of its properties. Its proof, based on a Taylor expansion of the tiling action, is deferred to \Cref{s:cusp_chart_proof}.
\begin{lemma}[$\fc$-cusp chart]\label{c:cusp_critical}
Given a cusp location $(x_0, s_0)\in \fA$ which is not a tangent location. Then for every $\fc>0$ sufficiently small  there exists
$\delta=\delta(\fc)>0$ sufficiently small such that the following holds.

Let $w_0=x_0-s_0\,\chi(x_0,s_0)$ be the critical point corresponding to $(x_0,s_0)$, and set
\[
\fU:=\{w\in\bC:\ |w-w_0|\le \fc\}.
\]
On $\fU$ the Riemann surface $\cC$ can be parametrized as $(f(w),w)$. We will therefore identify
$\fU$ with its image in $\cC$ under the map $w\mapsto (f(w),w)$. Moreover, for all $w\in\fU$,
\begin{align}\label{e:cusp_bounds}
 |\chi(w)|, |1-\chi(w)|, |\chi'''(w)|\asymp 1,\quad   |\chi''(w)| \asymp |w-w_0|
\end{align}
with implicit constants depending only on $\delta_0$.

Now let $(x,s)\in\fL$ satisfy $\|(x,s)-(x_0,s_0)\|_2\le \delta$.
Then the tiling action $S(\,\cdot\,;x,s)$ has exactly three critical points $w_c$ inside $\fU$ satisfying
\begin{align}
|w_c-w_0|
\lesssim \|(x,s)-(x_0,s_0)\|^{1/3}_2
\leq \delta^{1/3}.
\end{align}

Moreover, for each such critical point $w_c$ and every $w\in\fU$,
\begin{align}\label{e:cusp_quartic}
S(w;x,s)-S(w_c;x,s)=d(w-w_0)^4 +\cE(w),\quad d:=S''''(w_0;x_0,s_0)/24,\quad  |d|\asymp 1.
\end{align}
The error term satisfies
\begin{align}\label{e:arctic_err}
|\cE(w)|\leq C(\delta\ln(1/\delta)+|w-w_0|^5)\leq\frac{|d|\fc^4}{100}.
\end{align}

In this situation, we call $\fU$ a $\fc$-cusp chart (centered at $w_0$), and we say that
the point $(x,s)$ is \emph{adapted} to $\fU$.
\end{lemma}

We recall that in Panel (A), (B) and (C) of \Cref{f:cusp}, $S''''(w_0;x_0, s_0)<0$;  in Panel (D), (E) and (F) of \Cref{f:cusp}, $S''''(w_0;x_0, s_0)>0$.
Depending on the sign of $S''''(w_0; x_0, s_0)$, there are two cases for the local descent and ascent paths:
\begin{enumerate}
\item If $S''''(w_0;( x_0, s_0)) < 0$, the steepest–descent directions at $w_0$ are 
$\pm 1$ and $\pm \ri$. 
We introduce the following local paths, consisting of the non-real descent and ascent rays
\begin{align}\label{e:cusp_contour2}
\mathsf{C}^{\rm d}(w_0)
 = \bigl\{\, w_0 + r \ri : -\fc \le r \le \fc \,\bigr\},\quad \mathsf{C}^{\rm a}(w_0)
 = \bigl\{\, w_0 + r e^{\pm \ri\pi/4}, w_0 + r e^{\pm 3\ri\pi/4} : 0 \le r \le \fc \,\bigr\}.
\end{align}
\item If $S''''(w_0; x_0, s_0) > 0$, the steepest–descent directions at $w_0$ are 
$e^{\pm i\pi/4}$ and $e^{\pm 3\ri\pi/4}$. 
We introduce the following local paths, consisting of the non-real descent and ascent rays
\begin{align}\label{e:cusp_contour1}
\sfC^{\rm d} (w_0)
 = \bigl\{\, w_0 + r e^{\pm \ri\pi/4}, w_0 + r e^{\pm 3\ri\pi/4} : 0 \le r \le \fc \,\bigr\},\quad \sfC^{\rm a}(w_0)
 = \bigl\{\, w_0 + r \ri : -\fc \le r \le \fc \,\bigr\}.
\end{align}
\end{enumerate}
Note that \(\sfC^{\rm d}(w_0)\) and  \(\sfC^{\rm a}(w_0)\) are contained in \(\fU\).

In the following we orient the local paths \eqref{e:cusp_contour2} and  \eqref{e:cusp_contour1}, so they are compatible with the contours in the liquid region as introduced in \Cref{s:critical_bulk}. 

For $(x,s)\in \fL$ close to $(x_0, s_0)$, let $w_c=x-s\chi(w_c)\in\bC_+$,  we recall the factorization for descent vector from \eqref{eq:factor}:
\begin{align}
v^{\rm d}(w_c)=\sqrt{-\frac{1}{S''(w_c;x,s)}}
=\frac{\sqrt{\chi(w_c)}\sqrt{1-\chi(w_c)}}{\sqrt{1/s+\chi'(w_c)}}
=s^{1/2}\sqrt{\chi(w_c)}\sqrt{1-\chi(w_c)}\cdot
\frac{\sqrt{\partial_x\phi(x,s)}}{\sqrt{\phi'(w_c)}}.
\end{align}
Thus, a local choice of the branch of $\sqrt{\phi'(w_c)}$ determines the
oriented descent direction $v^{\rm d}(w_c)$, and conversely. By analytic
continuation, this choice also determines the branch of the prefactor
$\sqrt{\phi'(w)}$ appearing in the double-contour integral
\eqref{e:all_term0}.

We next give a geometric description of this choice as $(x,s)$ approaches
$(x_0,s_0)$ and $w_c$ approaches $w_0$. Let $w_c=w_0+a+\ri b$, then $b>0$ and \eqref{e:curve_cusp} gives
\begin{align}
1/s+\chi'(w_c)=\chi'''(w_0)b(a\ri-b/3) +\OO(b(|a|+b)^2)
\end{align}
Thus, as $(x,s)\to(x_0,s_0)$, $v^{\rm d}(w_c)=\sqrt{-1/S''(w_c;x,s)}$ is given by
\begin{align}\label{e:cusp_direction}
(1+\oo(1))\frac{\sqrt{\chi(w_0)}\sqrt{1-\chi(w_0)}}{\sqrt{b(a\ri-b/3) \chi'''(w_0)}}=\frac{(1+\oo(1))}{\sqrt{(a^2+b^2/9)b}} \sqrt{\frac{b/3+a\ri}{S''''(w_0;x_0, s_0)}},
\end{align}
where we used \eqref{e:derSfourth2}. 

\begin{figure}
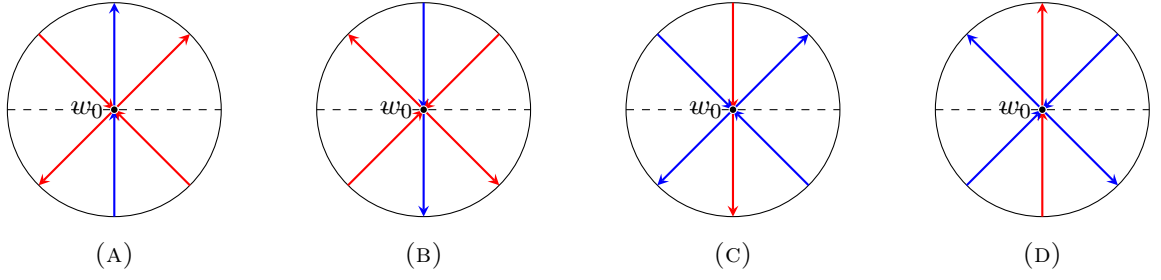

		\begin{subfigure}{0.24\textwidth}
			\begin{center}		

			\end{center}\caption{}
			\end{subfigure}
			\begin{subfigure}{0.24\textwidth}
			\begin{center}		
			%
			\end{center}\caption{}
			\end{subfigure}
			\begin{subfigure}{0.24\textwidth}
			\begin{center}		
			%
			\end{center}\caption{}
			\end{subfigure}
			\begin{subfigure}{0.24\textwidth}
			\begin{center}		
			%
			\end{center}\caption{}
			\end{subfigure}

				\caption{\label{f:cusp_path}  
Local paths associated with regular cusp region.
}
\end{figure}

If $S''''(w_0;x_0,s_0)<0$, then the direction of the square root in~\eqref{e:cusp_direction} lies in
$
\{\pm e^{\ri\theta}:\ \pi/4<\theta<3\pi/4\}.
$
If it lies in $\{e^{\ri\theta}:\ \pi/4<\theta<3\pi/4\}$, the corresponding direction is illustrated by the blue paths in panel~(A) of \Cref{f:cusp_path}. If it lies in $\{-e^{\ri\theta}:\ \pi/4<\theta<3\pi/4\}$, the corresponding direction is illustrated in panel~(B) of \Cref{f:cusp_path}.

If $S''''(w_0;x_0,s_0)>0$, then the direction of the limiting square root in~\eqref{e:cusp_direction} lies in
$
\{\pm e^{\ri\theta}:\ -\pi/4<\theta<\pi/4\}.
$
If it lies in $\{e^{\ri\theta}:\ -\pi/4<\theta<\pi/4\}$, the corresponding direction is illustrated by the blue curves in panel~(C) of \Cref{f:cusp_path}. If it lies in $\{-e^{\ri\theta}:\ -\pi/4<\theta<\pi/4\}$, the corresponding direction is illustrated in panel~(D) of \Cref{f:cusp_path}.

We orient the ascent contour $\sfC^{\rm a}(w_0)$ as indicated by the red paths in
\Cref{f:cusp_path}. With this convention, in the upper half-plane the ordered pair
$
\bigl(\sfC^{\rm d}(w_0),\sfC^{\rm a}(w_0)\bigr)$
intersects negatively at \(w_0\), while in the lower half-plane they intersect positively at \(w_0\), in the sense of
\Cref{d:positive_negative}.

\begin{figure}
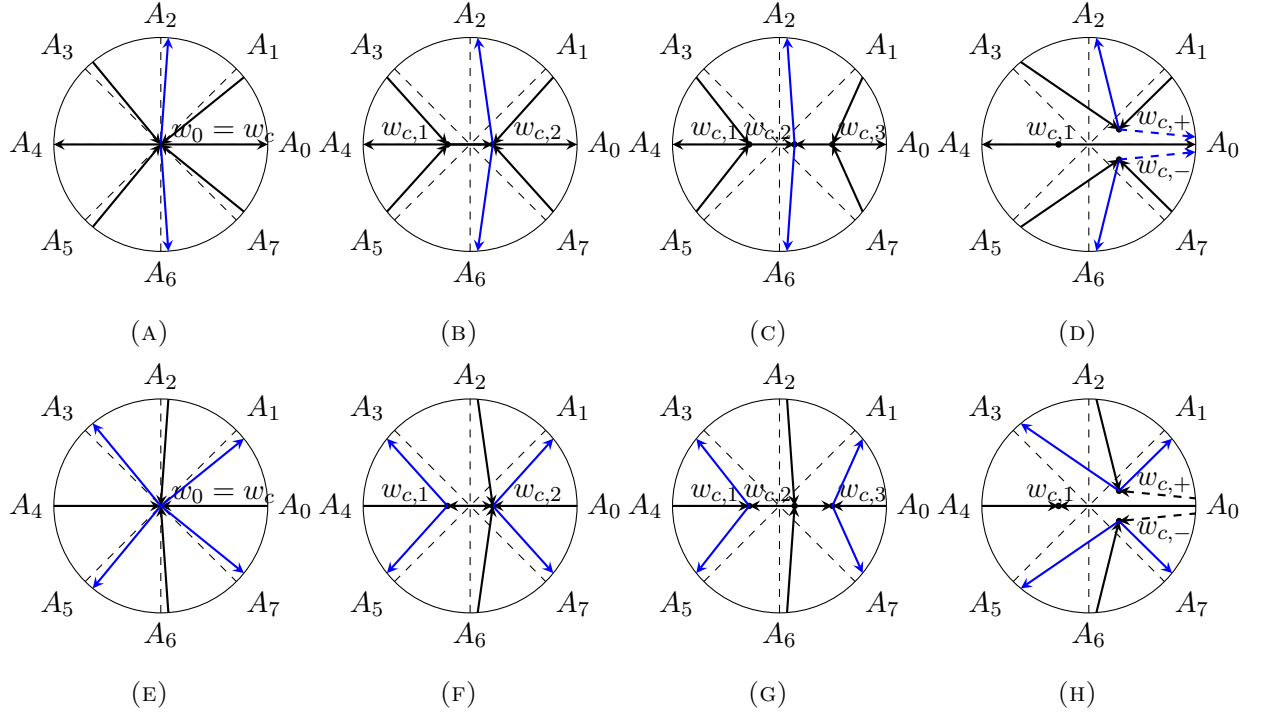


	\begin{subfigure}[t]{0.24\textwidth}
	\begin{center}		
		

		
			\end{center}
			\caption{}
			\end{subfigure}
				\begin{subfigure}[t]{0.24\textwidth}
	\begin{center}		
		
		%
		
			\end{center}
				\caption{}
			\end{subfigure}
			\begin{subfigure}[t]{0.24\textwidth}
			\begin{center}		
			%
			\end{center}
				\caption{}
			\end{subfigure}
			\begin{subfigure}[t]{0.24\textwidth}
			\begin{center}		
			%
			\end{center}
				\caption{}
			\end{subfigure}
					\caption{\label{f:quartic_saddle}  
Gradient flow of \(S(\cdot;x,s)\) in a cusp chart. Top row: \(S''''<0\). Bottom row: \(S''''>0\) (gradient directions reversed).
}
\end{figure}

The following lemma show that the local descent and ascent paths can be deformed into steepest-descent and steepest-ascent paths with negligible error. Its proof is deferred to \Cref{s:cusp_chart_proof}.
\begin{lemma}\label{l:cusp_steepest}
Adopt the assumptions and notation of \Cref{c:cusp_critical}. Then there exists a
constant $\fc'>0$ such that the following statements hold.

For each critical point \(w_c\), define the local steepest--descent path
\(\mathsf D^{\rm d}(w_c)\) at \(w_c\) to be the union of the portions of the
non-real steepest--descent trajectories of \(S(\,\cdot\,;x,s)\) starting at
\(w_c\), stopped at their first exit from the disk
\[
\{w: |w-w_0|\le \fc\};
\]
see \Cref{f:quartic_saddle}. If \(w_c\) is not a descent critical point, we set
\(\mathsf D^{\rm d}(w_c)=\emptyset\). The length of
\(\mathsf D^{\rm d}(w_c)\) is \(\OO(1)\).

There are two cases.
\begin{enumerate}
  \item If \((x,s)\in\fL\), then there are two complex conjugate critical
  points and one real critical point. We can deform \(\mathsf C^{\rm d}(w_0)\)
  to the union of the paths \(\mathsf D^{\rm d}(w_c)\) associated with the
  complex critical points, together with several arcs of total length
  \(\OO(1)\).

  \item If \((x,s)\in\fP\setminus \fL\), then there are three real critical
  points, counted with multiplicity. There are one or two descent critical
  points. We can deform \(\mathsf C^{\rm d}(w_0)\) to the union of the paths
  \(\mathsf D^{\rm d}(w_c)\) associated with all descent critical points,
  together with several arcs of total length \(\OO(1)\).
\end{enumerate}
In both cases, on these extra arcs,
\[
n\,\Re \bigl[S(w;x,s)-S(w_c;x,s)\bigr]\le -\,n\fc',
\]
where \(w_c\) denotes the corresponding descent critical point.

The analogous statements also hold for the local steepest--ascent paths.
\end{lemma}

%\begin{lemma}\label{l:cusp_S}
%Adopt the notation from \Cref{c:cusp_critical}. Let $\fU$ be a cusp chart, and let
%$w,z\in \fU$. Then, uniformly for $w,z\in \fU$, \eqref{e:PI+_bound} and \eqref{e:QI-_bound} hold.
%\end{lemma}

\begin{figure}
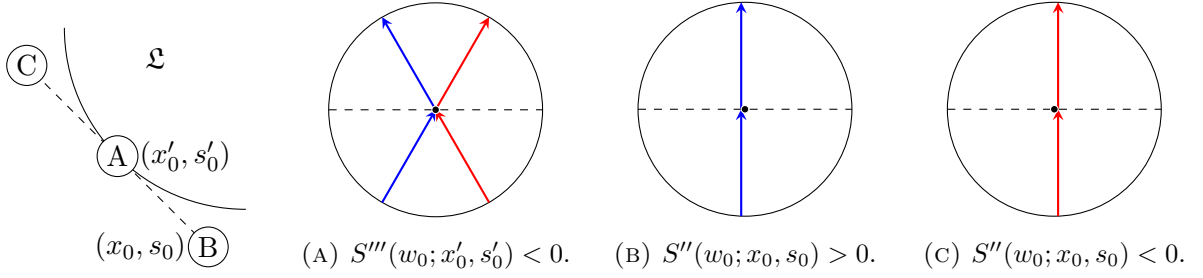

	\begin{subfigure}{0.24\textwidth}
	
		\centering
		%
		\end{center}
		\caption{$S'''(w_0;x'_0,s'_0)<0$.}
		\end{subfigure}	
		\begin{subfigure}{0.24\textwidth}
			\begin{center}		
			%
			\end{center}
			\caption{$S''(w_0;x_0,s_0)>0$.}
			\end{subfigure}
			\begin{subfigure}{0.24\textwidth}
			\begin{center}		
			%
				
			\end{center}
			\caption{$S''(w_0;x_0,s_0)<0$.}
			\end{subfigure}

				\caption{\label{f:frozen_path1}  
Local path associated with regular frozen chart, where the tangent location is a regular arctic point.
}
\end{figure}

\subsection{Regular frozen chart and cusp frozen chart}
\label{s:frozen_neighborhood}

For any $(x_0,s_0)\in \fP$ and a line through \((x_0,s_0)\) that is tangent to the arctic boundary at \((x_0',s_0')\). Assume \((x_0,s_0)\) is bounded away from \((x_0',s_0')\), and \((x_0',s_0')\) is bounded away from any tangent locations, see \Cref{f:frozen_path1}. We remark that it is possible that $(x_0, s_0)$ is itself on the arctic boundary.
Define \(w_0\) by the characteristic relation
\[
w_0 \;=\; x_0' - s_0'\,\chi(x_0',s_0') \;=\; x'_0 - s'_0\,\chi(w_0)\in\bR.
\]
By the second statement of \Cref{p:surface}, we have \(w_0\in\mathbb{R}\) and
\[
\chi'(w_0) = -\frac{1}{s_0'} , 
\qquad 
|\chi(w_0)|, |1-\chi(w_0)| \asymp 1 .
\]
Moreover, the slope of the tangent line is \(\chi(w_0)\), and it follows that
\begin{align}
w_0=x_0'-s_0'\chi(w_0)=x_0-s_0\chi(w_0).
\end{align}
Consequently $S'(w_0;x_0,s_0)=0$, so $w_0$ is a critical point.

\begin{figure}

  \begin{subfigure}{0.3\textwidth}
    \begin{center}
      \begin{tikzpicture}[scale=1.1]
      \draw[] (0,0) arc (45:0:3);
       \draw[] (0,0) arc (225:270:3);
         \fill (0,0) circle (1.5pt);
    \draw[](-0.2,0) node[left] { $(x'_0,  s'_0)$};
    \node[] at (1,0) {$\fL$};
    
	%\draw[white, fill=black]  ({-3/sqrt(2)+3*cos(30)},{-3/sqrt(2)+3*sin(30)}) circle (0.05);
	   %\draw[]({-3/sqrt(2)+3*cos(30)},{-3/sqrt(2)+3*sin(30)}) node[left] { $(x'_0, s'_0)$};

	%\draw[dashed] ({-3/sqrt(2)+3*cos(30)},{-3/sqrt(2)+3*sin(30)})--({-3/sqrt(2)+3*cos(30)+1.5*cos(-60)},{-3/sqrt(2)+3*sin(30)+1.5*sin(-60)});

\draw[dashed] (0,0)--(1.5,-1.5);

\node[circle, draw, fill=white, inner sep=1.5pt] at (0,0) {A};

	\node[circle, draw, fill=white, inner sep=1.5pt] at (1.5,-1.5) {B};	
	 
	 \draw[]({1.5+0.2},{-1.5}) node[right] { $(x_0, s_0)$};
             
\node[] at (0,-2.5) {\phantom{=}};
      \end{tikzpicture}
      \end{center}
  \end{subfigure}
\begin{subfigure}{0.24\textwidth}
			\begin{center}		
			\begin{tikzpicture}[scale=1]
			\draw (0,0) circle [radius=sqrt(2)];
			\draw[dashed] (-{sqrt(2)},0)--({sqrt(2)},0);
			
			\draw[-stealth, thick, blue] (0,0)--(0,{sqrt(2)});
			\draw[-stealth, thick, blue] (0,-{sqrt(2)})--(0,0);
			
			\draw[stealth-, thick, red]  (-1,-1)--(0,0);
			\draw[stealth-, thick, red]  (1,1)--(0,0);
			\draw[stealth-, thick, red]  (0,0)--(1,-1);
			\draw[stealth-, thick, red]  (0,0)--(-1,1);
			\draw[white, fill=black]  (0,0) circle (0.05);
			
			\end{tikzpicture}
			\end{center}\caption{$S''''(w_0;x_0, s_0)<0$}
			\end{subfigure}
			\begin{subfigure}{0.24\textwidth}
			\begin{center}		
			\begin{tikzpicture}[scale=1]
			\draw (0,0) circle [radius=sqrt(2)];
			\draw[dashed] (-{sqrt(2)},0)--({sqrt(2)},0);
			
			\draw[-stealth, thick, blue] (0,0)--(0,{sqrt(2)});
			\draw[-stealth, thick, blue] (0,-{sqrt(2)})--(0,0);
			
			\draw[white, fill=black]  (0,0) circle (0.05);
			
			\end{tikzpicture}
			\end{center}\caption{$S''(w_0;x_0,s_0)>0$.}
			\end{subfigure}  
  
	 \caption{\label{f:frozen_path2}  
Local path associated with regular frozen chart, where the tangent location is a regular cusp point.
}
	 \end{figure}
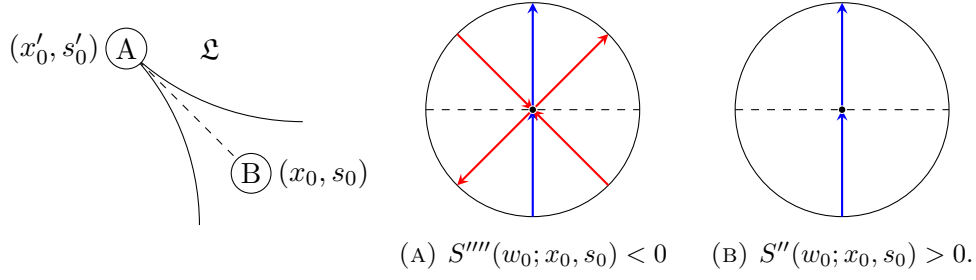

 Locally round $w_0$, $S(w;x_0,s_0)$ is holomorphic, and we have the Taylor expansion
\begin{align}\begin{split}\label{e:frozen_S''}
&S(w;x_0,s_0)=S(w_0;x_0,s_0)+S''(w_0; x_0,s_0)\frac{(w-w_0)^2}{2} +\OO(|w-w_0|^3),\\
 &S''(w_0; x_0, s_0)=-\left(\chi'(w_0)+\frac{1}{s_0}\right)\frac{1}{\chi(w_0)(1-\chi(w_0))}=\frac{s_0-s'_0}{s_0s'_0\chi(w_0)(1-\chi(w_0))}
\end{split}\end{align}
In this case, $S''(w_0;x_0,s_0)\neq 0$. Recall from \Cref{s:dacritical}
that if $S''(w_0;x_0,s_0)>0$, then $w_0$ is a descent critical point and
issues two non-real steepest-descent paths; if $S''(w_0;x_0,s_0)<0$, then
$w_0$ is an ascent critical point and issues two non-real steepest-ascent
paths; see \Cref{f:frozen_path1} and \Cref{f:frozen_path2}.

%{\color{red}
%If $(x_0',s_0')$ lies within distance $\delta$ of a tangent location, then $(x_0,s_0)$ lies within distance
%$\OO(\delta)$ of an extended side.
%}
The following lemmas introduce the regular frozen chart and cusp frozen chart,  and record several of their properties. Their proof, based on a Taylor expansion of the tiling action, are deferred to \Cref{s:frozen_chart_proof}.

\begin{lemma}[$\fc$-frozen chart (regular)]\label{c:regular_frozen_critical}
Fix $\delta_0>0$. Let $(x_0,s_0)\in\fB$, and let $(x_0',s_0')$ be a point on the arctic boundary such that the line
through $(x_0,s_0)$ is tangent to the arctic boundary at $(x_0',s_0')$. Assume that $(x_0', s_0')$ is at least $\delta_0$ distance from $(x_0, s_0)$ and all cusp and tangency locations.
Then for every sufficiently small $\fc>0$ (depending only on $\delta_0$) there exists
$\delta=\delta(\fc)>0$, also sufficiently small and depending on $\delta_0$, such that the following holds. 

Let $w_0=x_0-s_0\chi(x_0', s_0')$ be the critical point corresponding to $(x_0, s_0)$, and set
\begin{align}
\fU:=\{w\in \bC:|w-w_0|\leq \fc\}.
\end{align}
On $\fU$ the Riemann surface $\cC$ can be parametrized as $(f(w), w)$. We will therefore identify
$\fU$ with its image in $\cC$ under the map $w\mapsto (f(w), w)$.  Moreover, for all $w\in\fU$,
\begin{align}
 |\chi(w)|, |1-\chi(w)|, |\chi''(w)|\asymp 1,
\end{align}
with implicit constants depending only on $\delta_0$. 

Now let $(x,s)\in\fL$ satisfy $\|(x,s)-(x_0,s_0)\|_2\le \delta$.
Then the tiling action $S(\,\cdot\,;x,s)$ has exactly one critical point $w_c$ inside $\fU$, satisfying 
\begin{align}\label{e:frozen_wc1}
 |w_c-w_0|\lesssim \|(x,s)-(x_0,s_0)\|_2\leq {\delta}.
\end{align}
Moreover, for every $w\in\fU$, 
\begin{align}\label{e:frozen_wc2}
S(w;x,s)-S(w_c;x,s)=dw^2 +\cE(w),\quad d:=S''(0;x_0,s_0)/2,\quad |d|\asymp 1.
\end{align}
The error term satisfies
\begin{align}\label{e:frozen_wc3}
|\cE(w)|\leq C(\delta\ln(1/\delta)+|w-w_0|^3)\leq\frac{|d|\fc^2}{100}.
\end{align}

In this situation, we call $\fU$ a regular $\fc$-frozen chart (centered at $w_0$), and we say that
the point $(x,s)$ is \emph{adapted} to $\fU$.
\end{lemma}

\begin{lemma}[$\fc$-frozen chart (cusp)]\label{c:cusp_frozen_critical}
Fix $\delta_0>0$. Let $(x_0,s_0)\in\fB$, and let $(x_0',s_0')$ be a point on the arctic boundary such that the line
through $(x_0,s_0)$ is tangent to the arctic boundary at $(x_0',s_0')$. Assume that $(x_0', s_0')$ is a cusp location (but not a cusp-turning location). Assume that $(x_0', s_0')$ is at least distance $\delta_0$  from $(x_0, s_0)$. 
Then for every sufficiently small $\fc>0$ there exists
$\delta=\delta(\fc)>0$, also sufficiently small and depending on $\delta_0$,  such that the following holds. 

Let $w_0=x_0-s_0\chi(x_0', s_0')$ be the critical point corresponding to $(x_0, s_0)$, and set
\begin{align}
\fU:=\{w\in \bC:|w-w_0|\leq \fc\}.
\end{align}
On $\fU$ the Riemann surface $\cC$ can be parametrized as $(f(w), w)$. We will therefore identify
$\fU$ with its image in $\cC$ under the map $w\mapsto (f(w), w)$.  Moreover, for all $w\in\fU$,
\begin{align}
 |\chi(w)|, |1-\chi(w)|, |\chi'''(w)|\asymp 1,\quad |\chi''(w)|\asymp |w-w_0|.
\end{align}

Now let $(x,s)\in\fL$ satisfy $\|(x,s)-(x_0,s_0)\|_2\le \delta$.
Then the tiling action $S(\,\cdot\,;x,s)$ has exactly one critical point $w_c$ inside $\fU$, satisfying 
\eqref{e:frozen_wc1}, \eqref{e:frozen_wc2} and \eqref{e:frozen_wc3}.

In this situation, we call $\fU$ a cusp $\fc$-frozen chart (centered at $w_0$), and we say that
the point $(x,s)$ is \emph{adapted} to $\fU$.
\end{lemma}

If $S''(w_0;x_0,s_0)>0$; see Panel (B) of \Cref{f:frozen_path1}, then we
associate with $w_0$ the local descent contour
\begin{align}\label{e:frozen_direction1}
\sfC^{\rm d}(w_0)
= \{\, w_0+r\ri : -\fc \le r \le \fc \,\},
\qquad
\sfC^{\rm a}(w_0)=\emptyset .
\end{align}

If $S''(w_0;x_0,s_0)<0$; see Panel (C) of \Cref{f:frozen_path1}, then we
associate with $w_0$ the local ascent contour
\begin{align}\label{e:frozen_direction2}
\sfC^{\rm a}(w_0)
= \{\, w_0+r\ri : -\fc \le r \le \fc \,\},
\qquad
\sfC^{\rm d}(w_0)=\emptyset .
\end{align}

In what follows, we orient the local paths
\eqref{e:frozen_direction1} and \eqref{e:frozen_direction2}.

We first consider the case where the tangency point
$(x'_0,s'_0)\in\fA$ is not a cusp location. In \Cref{s:critical_arctic}, we
associated with $(x'_0,s'_0)$ a local descent path and a local ascent path;
see Panel (A) of \Cref{f:frozen_path1}. We orient
\eqref{e:frozen_direction1} in the same direction as the descent path at the
tangency point, and \eqref{e:frozen_direction2} in the same direction as the
ascent path at the tangency point. In this way, the orientation of the local
descent/ascent paths remains unchanged as we move from $(x'_0,s'_0)$ to
$(x_0,s_0)$ along the tangent line.

When $(x'_0,s'_0)$ is a cusp location, in \Cref{s:critical_cusp} we associated
with the cusp local descent and ascent paths; see Panel (A) of
\Cref{f:frozen_path2}. There are two cases:
\begin{enumerate}
\item If $\chi'''(w_0)>0$, then by the second statement of \Cref{p:surface}
      (see also \Cref{f:cusp}), the cusp points downward, and hence
      $s_0>s'_0$.
\item If $\chi'''(w_0)<0$, then by the second statement of \Cref{p:surface}
      (see also \Cref{f:cusp}), the cusp points upward, and hence
      $s_0<s'_0$.
\end{enumerate}
In both cases, the following two quantities have opposite signs:
\begin{align}
S{''''}(w_0;x'_0,s'_0)
=-\frac{\chi'''(w_0)}
{\chi(w_0)\bigl(1-\chi(w_0)\bigr)},\quad 
S''(w_0;x_0,s_0)
=\frac{s_0-s'_0}
{s_0s'_0\chi(w_0)\bigl(1-\chi(w_0)\bigr)} .
\end{align}
Thus, if $S''(w_0;x_0,s_0)>0$, then
$S''''(w_0;x'_0,s'_0)<0$, and we orient
\eqref{e:frozen_direction1} in the same direction as the descent paths
associated with the cusp; see Panels (A) and (B) of \Cref{f:cusp_path} and
\Cref{f:frozen_path2}. If $S''(w_0;x_0,s_0)<0$, then
$S''''(w_0;x'_0,s'_0)>0$, and we orient
\eqref{e:frozen_direction2} in the same direction as the ascent paths
associated with the cusp; see Panels (C) and (D) of \Cref{f:cusp_path}. Again,
with this convention, the orientation of the local descent/ascent paths remains
unchanged as we move from $(x'_0,s'_0)$ to $(x_0,s_0)$ along the tangent line.

The following lemma show that the local descent and ascent paths can be deformed into steepest-descent and steepest-ascent paths with negligible error. Its proof is deferred to \Cref{s:frozen_chart_proof}.
\begin{lemma}\label{c:frozen_steepest}
Under the assumptions and with the notation of either
\Cref{c:regular_frozen_critical} or \Cref{c:cusp_frozen_critical}, there exists
a constant $\fc'>0$ such that the following statements hold.

Suppose first that \(S''(w_0;x_0,s_0)>0\). Then \(w_c\) is a descent critical
point, and $\sfC^{\rm a}(w_0)=\emptyset$.

\begin{itemize}
\item \emph{Replacing \(\sfC^{\rm d}(w_0)\) by \(\sfD^{\rm d}(w_c)\).}
We define the local steepest-descent set $\sfD^{\rm d}(w_c)$ to be the union of
the portions of the steepest-descent trajectories of $S(\,\cdot\,;x,s)$ issuing
from \(w_c\), stopped at their first exit from the disk
\[
\{\,w: |w-w_0|\le \fc\,\};
\]
see \Cref{f:bulk_critical}. Then $\sfD^{\rm d}(w_c)$ has total length $\OO(1)$.
Moreover, $\sfC^{\rm d}(w_0)$ can be deformed into the union of
$\sfD^{\rm d}(w_c)$ and finitely many arcs of total length $\OO(1)$, along
which
\begin{align}\label{e:bulk_subarc}
n\,\Re\bigl[ S(w;x,s)-S(w_c;x,s)\bigr]\le -\,n\fc'.
\end{align}

\item \emph{Replacing \(\sfD^{\rm d}(w_c)\) by \(\sfS^{\rm d}(w_c)\).}
Let
$
r_n:={\ln n}/{\sqrt n}.
$
We define the truncated steepest-descent set $\sfS^{\rm d}(w_c)$ to be the union
of the portions of the steepest-descent trajectories issuing from \(w_c\),
stopped at their first exit from the disk
\[
\{\,w: |w-w_c|\le r_n\,\}.
\]
Then $\sfS^{\rm d}(w_c)$ has total length $\OO(r_n)$, and for all
\(w\in \sfD^{\rm d}(w_c)\setminus \sfS^{\rm d}(w_c)\),
\[
e^{n\Re[S(w;x,s)]}
\le
e^{n\Re[S(w_c;x,s)]}e^{-\fc'(\ln n)^2}.
\]
\end{itemize}

If \(S''(w_0;x_0,s_0)<0\), then \(w_c\) is an ascent critical point, and
$\sfC^{\rm d}(w_0)=\emptyset$. In this case, we define the analogous ascent
sets \(\sfD^{\rm a}(w_c)\) and \(\sfS^{\rm a}(w_c)\). The preceding statements
hold with \(S\) replaced by \(-S\) and with the superscript \({\rm d}\) replaced
by \({\rm a}\).
\end{lemma}

\begin{figure}
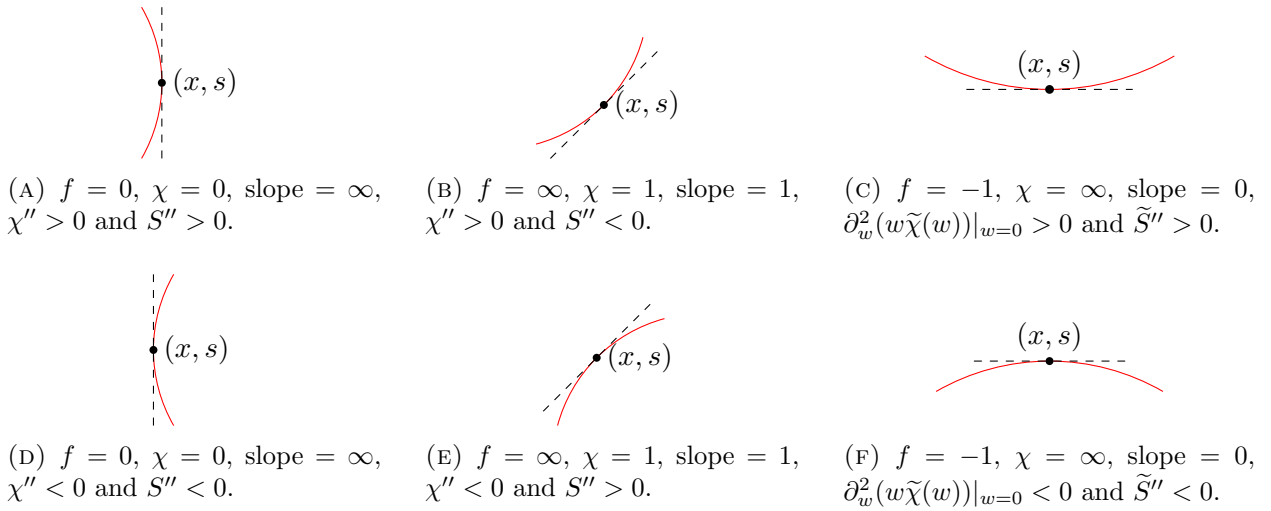

  \centering

  % Row 1
  \begin{subfigure}[t]{0.3\textwidth}
    \centering

    \caption{$f=0$, $\chi=0$, $\text{slope}=\infty$, $\chi''>0$ and $S''>0$.}
  \end{subfigure}
  \hfill
  \begin{subfigure}[t]{0.3\textwidth}
    \centering
    %
    \caption{$f=\infty$, $\chi=1$, $\text{slope}=1$,  $\chi''>0$ and $S''<0$.}
  \end{subfigure}
  \hfill
  \begin{subfigure}[t]{0.33\textwidth}
    \centering
    %
    \caption{$f=-1$, $\chi=\infty$, $\text{slope}=0$, %$(w\widetilde\chi(w))|_{w=0}=1/s$, 
    $\partial_w^2(w\widetilde\chi(w))|_{w=0}>0$ and $\wt S''>0$.}
  \end{subfigure}

  \vspace{1em}

  % Row 2
  \begin{subfigure}[t]{0.3\textwidth}
    \centering
    %
    \caption{$f=0$, $\chi=0$, $\text{slope}=\infty$, $\chi''<0$ and $S''<0$.}
  \end{subfigure}
  \hfill
  \begin{subfigure}[t]{0.3\textwidth}
    \centering
    %
    \caption{$f=\infty$, $\chi=1$, $\text{slope}=1$,  $\chi''<0$ and $S''>0$.}
  \end{subfigure}
  \hfill
  \begin{subfigure}[t]{0.33\textwidth}
    \centering
    %
    \caption{$f=-1$, $\chi=\infty$, $\text{slope}=0$, %$(w\widetilde\chi(w))|_{w=0}=1/s$, 
    $\partial^2_w(w\widetilde\chi(w))|_{w=0}<0$ and $\wt S''<0$.}
  \end{subfigure}

  \caption{Tangent locations.}
  \label{f:tangent}
\end{figure}

\begin{figure}
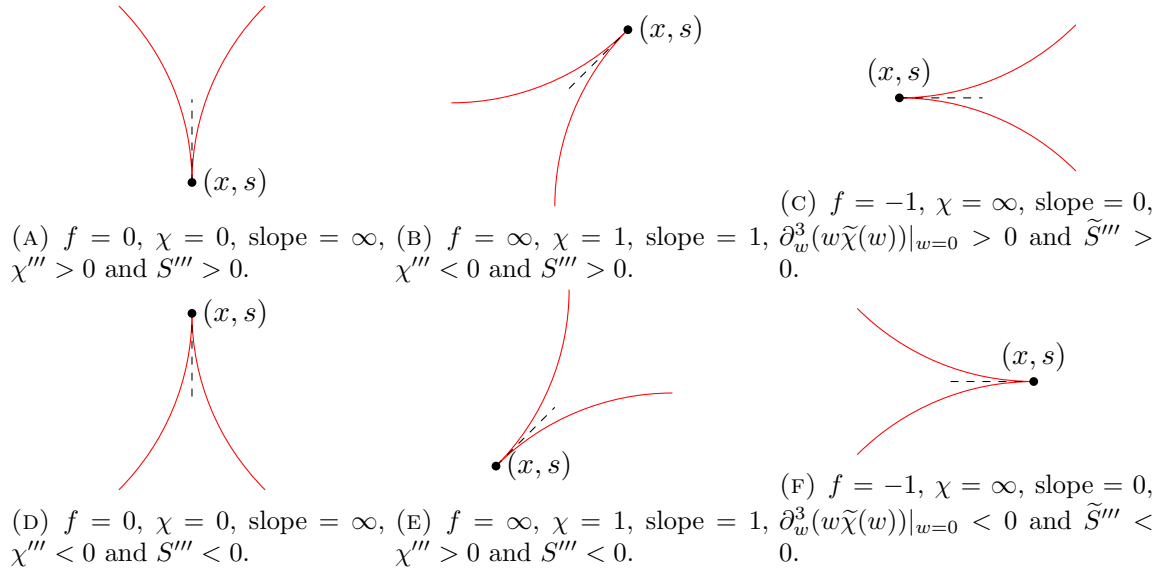

\begin{subfigure}{0.3\textwidth}
    \centering
      %
      \caption{$f=0$, $\chi=0$, $\text{slope}=\infty$, $\chi'''>0$ and $S'''>0$.}
  \end{subfigure}
  \begin{subfigure}{0.3\textwidth}
    \centering
      %
      \caption{$f=\infty$, $\chi=1$, $\text{slope}=1$,  $\chi'''<0$ and $S'''>0$.}
  \end{subfigure}
\begin{subfigure}{0.3\textwidth}
    \centering
      %
      \caption{$f=-1$, $\chi=\infty$, $\text{slope}=0$, $\del^3_w(w\wt \chi(w))|_{w=0}>0$ and $\wt S'''>0$.}
  \end{subfigure}

\begin{subfigure}{0.3\textwidth}
    \centering
      %
      \caption{$f=0$, $\chi=0$, $\text{slope}=\infty$, $\chi'''<0$ and $S'''<0$.}
  \end{subfigure}
\begin{subfigure}{0.3\textwidth}
    \centering
      %
      \caption{$f=\infty$, $\chi=1$, $\text{slope}=1$,  $\chi'''>0$ and $S'''<0$.}
  \end{subfigure}
\begin{subfigure}{0.3\textwidth}
    \centering
      %
      \caption{$f=-1$, $\chi=\infty$, $\text{slope}=0$, $\del^3_w(w\wt \chi(w))|_{w=0}<0$ and $\wt S'''<0$.}
  \end{subfigure}

	 \caption{cusp-turning locations.}
	 \label{f:cusp_turning}
	 \end{figure}

\subsection{Vertical tangent chart and cusp-turning chart}\label{s:vertical_tangent}
%There are three types of tangent locations: vertical tangent locations ($f=0$), unit-slope tangent locations ($f=\infty$), and horizontal tangent locations ($f=-1$).  
For any vertical tangency point \((x_0,s_0)\in\fA\), from \eqref{e:bcond} and
\eqref{e:arctic_chider}, we have
$
f(x_0,s_0)=\chi(x_0,s_0)=0$,
and
\begin{align}
w_0:=x_0-s_0\chi(x_0,s_0)=x_0,
\qquad
\chi(w_0)=0,
\qquad
\chi'(w_0)=-\frac1{s_0}.
\end{align}
Moreover, \((x_0,s_0)\) lies on the common boundary of two adjacent curvilinear
triangles \(\fT_A\) and \(\fT_B\); see \Cref{f:adjacent_curvilinear_triangle}.
On each of \(\fT_A\) and \(\fT_B\), the height function has constant slope.

In a neighborhood of \(w_0\), write
$
\chi(w)=(w_0-w)\chi_0(w)$.
Then, differentiating both sides, we get
\begin{align}
-\frac1{s_0}
=
\chi'(w_0)
=
\left[-\chi_0(w)+(w_0-w)\chi_0'(w)\right]_{w=w_0}
=
-\chi_0(w_0),
\end{align}
and hence
\[
\chi_0(w_0)=\frac1{s_0}>0.
\]
Thus, locally around \(w_0\), we can expand \(\chi(w)\) and \(\chi_0(w)\) as
\begin{align}\label{e:expandchi_b}
\begin{split}
\chi(w)
&=
-\frac{w-w_0}{s_0}
+c_1(w-w_0)^2
+c_2(w-w_0)^3+\cdots,\\
\chi_0(w)
&=
\frac1{s_0}
+c_1(w_0-w)
-c_2(w_0-w)^2+\cdots,
\end{split}
\end{align}
where $c_1=\chi''(w_0)/2$ and $c_2=\chi'''(w_0)/6$.
Moreover,
\begin{align}\label{e:vert_tangent}
\begin{split}
S'(w;x_0,s_0)
&=
\ln \frac{w-x_0}{x_0-s_0-w}
-\ln \frac{\chi(w)}{1-\chi(w)}\\
&=
\ln \frac{(1-\chi(w))(w-w_0)}
{\chi(w)(x_0-s_0-w)}
=
\ln \frac{1-(w_0-w)\chi_0(w)}
{\chi_0(w)(w-x_0+s_0)}\\
&=
\ln
\frac{
1+(w-w_0)/s_0-c_1(w-w_0)^2+\OO\!\bigl((w-w_0)^3\bigr)}
{
\bigl(1-s_0c_1(w-w_0)-s_0c_2(w-w_0)^2
+\OO\!\bigl((w-w_0)^3\bigr)\bigr)
\bigl(1+(w-w_0)/s_0\bigr)
}\\
&=
s_0c_1(w-w_0)
+\Bigl(s_0c_2+\frac12s_0^2c_1^2-c_1\Bigr)(w-w_0)^2
+\OO\!\bigl((w-w_0)^3\bigr).
\end{split}
\end{align}
In particular, \(S'(w;x_0,s_0)\) is analytic in a neighborhood of \(w_0\). If
\((x_0,s_0)\) is not a cusp location, then by \eqref{e:arctic_chider},
\begin{align}\label{e:c_1}
c_1={\chi''(w_0)}/2,\quad |c_1|\asymp 1.
\end{align}
If \((x_0,s_0)\) is a cusp-turning location, then by
\eqref{e:introducef0_cusp},
\begin{align}\label{e:c_2}
\chi''(w_0)=0,\quad
c_2={\chi'''(w_0)}/6,\quad |c_2|\asymp 1.
\end{align}
 In these cases, the signs of  \(S''(w_0;x_0,s_0)\) and \(S'''(w_0;x_0,s_0)\) are determined by the signs of $\chi''(w_0)$ and $\chi'''(w_0)$, as classified in \Cref{f:tangent} and \Cref{f:cusp_turning}.

For \((x,s)\) in a sufficiently small neighborhood of \((x_0,s_0)\), we have
\begin{align}\label{e:vert_tangent_diff}
\begin{split}
S'(w;x,s)-S'(w;x_0,s_0)
=
\ln \frac{w-x}{x-s-w}
-\ln \frac{w-x_0}{x_0-s_0-w}=
\ln \frac{w-x}{w-x_0}
+\ln \frac{x_0-s_0-w}{x-s-w}.
\end{split}
\end{align}
The first logarithm has a branch cut along the segment
\[
[\min\{x_0,x\},\,\max\{x_0,x\}],
\]
and locally around \(w_0=x_0\), the difference
\eqref{e:vert_tangent_diff} is analytic on the complement of this segment.

The following lemmas introduce the vertical tangent chart and vertical cusp-turning chart, and  record several of their properties. Their proof, based on a Taylor expansion of the tiling action, is deferred to \Cref{s:vertical_tangent_chart_proof}.
\begin{lemma}[$\fc$-tangent chart (vertical)]\label{c:tangent_critical1}
Given a vertical tangent location $(x_0, s_0)\in \fA$ which is not a cusp location. Then for every $\fc>0$ sufficiently small  there exists
$\delta=\delta(\fc)>0$ sufficiently small such that the following holds.

Let $w_0=x_0-s_0\,\chi(x_0,s_0)=x_0$ be the critical point corresponding to $(x_0,s_0)$, and set
\[
\fU:=\{w\in\bC:\ |w-w_0|\le \fc\}.
\]
On $\fU$ the Riemann surface $\cC$ can be parametrized as $(f(w),w)$. We will therefore identify
$\fU$ with its image in $\cC$ under the map $w\mapsto (f(w),w)$.  Moreover, for all $w\in\fU$,
\begin{align}\label{e:chi_v_tangent}
  |1-\chi(w)|, |\chi''(w)|\asymp 1,\quad  |\chi(w)|\asymp |w-w_0|.
\end{align}

Now let $(x,s)\in\fL$ satisfy $\|(x,s)-(x_0,s_0)\|_2\le \delta$.
Then the tiling action $S(\,\cdot\,;x,s)$ has exactly two (formal) critical points $w_c$ inside $\fU$, satisfying 
\begin{align}\begin{split}\label{e:wc_vertical_tangent}
&\delta^{-1/2}|x-x_0|\lesssim |w_c-w_0|\lesssim \|(x,s)-(x_0,s_0)\|^{1/2}_2\leq {\delta}^{1/2},\\
&\dist\bigl(w_c,[\min\{x_0,x\},\max\{x_0,x\}]\bigr)
\gtrsim |x-x_0|
\end{split}\end{align}

Moreover, for each such critical point $w_c$ and every $w\in\fU$, 
\begin{align}\label{e:tangent_S}
S(w;x,s)-S(w_c;x,s)=d(w-w_0)^2 +\cE(w),\quad d:=S''(w_0;x_0,s_0)/2, \quad |d|\asymp 1.
\end{align}
The error term satisfies
\begin{align}\label{e:tangent_err}
|\cE(w)|\leq C(\delta\ln(1/\delta)+|w-w_0|^3)\leq\frac{|d|\fc^2}{100}.
\end{align}

In this situation, we call $\fU$ a $\fc$-tangent chart (centered at $w_0$), and we say that
the point $(x,s)$ is \emph{adapted} to $\fU$.
\end{lemma}

\begin{lemma}[$\fc$-cusp-turning chart (vertical)]\label{c:cusp_turning_critical1}
Given a vertical cusp-turning location $(x_0, s_0)\in \fA$. Then for every $\fc>0$ sufficiently small  there exists
$\delta=\delta(\fc)>0$ sufficiently small such that the following holds.

Let $w_0=x_0-s_0\,\chi(x_0,s_0)=x_0$ be the critical point corresponding to $(x_0,s_0)$, and set
\[
\fU:=\{w\in\bC:\ |w-w_0|\le \fc\}.
\]
On $\fU$ the Riemann surface $\cC$ can be parametrized as $(f(w),w)$. We will therefore identify
$\fU$ with its image in $\cC$ under the map $w\mapsto (f(w),w)$. Moreover, for all $w\in\fU$,
\begin{align}\label{e:chi_v_cusp_turning}
 |1-\chi(w)|, |\chi'''(w)| \asymp 1,\quad  |\chi(w)|, |\chi''(w)|\asymp  |w-w_0|.
\end{align}

Now let $(x,s)\in\fL$ satisfy $\|(x,s)-(x_0,s_0)\|_2\le \delta$.
Then the tiling action $S(\,\cdot\,;x,s)$ has exactly three (formal) critical points $w_c$ inside $\fU$, satisfying 
\begin{align}\begin{split}\label{e:wc_vertical_cusp}
&\delta^{-2/3}|x-x_0|\lesssim |w_c-w_0|\lesssim \|(x,s)-(x_0,s_0)\|^{1/3}_2\leq {\delta}^{1/3},\\
&
\dist\bigl(w_c,[\min\{x_0,x\},\max\{x_0,x\}]\bigr)
\gtrsim |x-x_0|.
\end{split}\end{align}

Moreover, for each such critical point $w_c$ and every $w\in\fU$, 
\begin{align}\label{e:cusp_turning_S}
S(w;x,s)-S(w_c;x,s)=d(w-w_0)^3 +\cE(w),\quad d:=S'''(w_0;x_0,s_0)/6, \quad |d|\asymp 1.
\end{align}
The error term satisfies
\begin{align}\label{e:cusp_turning_err}
|\cE(w)|\leq C(\delta\ln(1/\delta)+|w-w_0|^4)\leq\frac{|d|\fc^3}{100}.
\end{align}

In this situation, we call $\fU$ a $\fc$-cusp-turning chart (centered at $w_0$), and we say that
the point $(x,s)$ is \emph{adapted} to $\fU$.
\end{lemma}

\begin{figure}					
	\begin{subfigure}[t]{0.30\textwidth}
			\centering

			\caption*{$S''>0$}
		\end{subfigure}	
		\begin{subfigure}[t]{0.3\textwidth}
			\centering
			%
		\caption{}

	\end{subfigure}
	\begin{subfigure}[t]{0.3\textwidth}
			\centering
			%
		\caption{}

	\end{subfigure}

	\caption{
\label{f:vertical_tangent1}
Local path associated with vertical tangent locations, $S''>0$.}
	\end{figure}

\begin{figure}			
		\begin{subfigure}[t]{0.30\textwidth}
			\centering
			%
					\caption*{$S''>0$}
					\end{subfigure}	
			\begin{subfigure}[t]{0.3\textwidth}
			\centering
			%
		\caption{}
	\end{subfigure}
\begin{subfigure}[t]{0.3\textwidth}
			\centering
			%
		\caption{}
	\end{subfigure}
\caption{
\label{f:vertical_tangent2}
Local path associated with vertical tangent locations, $S''>0$.}
	\end{figure}

\begin{figure}			
		\begin{subfigure}[t]{0.3\textwidth}

		\centering
			%
				
			\caption*{$S''<0$}
			\end{subfigure}
		\begin{subfigure}[t]{0.3\textwidth}
			\centering
			%
		\caption{}

	\end{subfigure}
		\begin{subfigure}[t]{0.3\textwidth}
			\centering
			%
		\caption{}

	\end{subfigure}
	
\caption{
Local path associated with vertical tangent locations, $S''<0$.}
	\label{f:vertical_tangent3}
	\end{figure}

\begin{figure}

	\begin{subfigure}[t]{0.3\textwidth}
		\centering
			%
			\caption*{$S''<0$}
			\end{subfigure}
			\begin{subfigure}[t]{0.3\textwidth}
			\centering
			%

		\caption{}
	\end{subfigure}
			
	\caption{
\label{f:vertical_tangent4}
Local path associated with vertical tangent locations, $S''<0$.}
	\end{figure}

The tangent location $(x_0,s_0)$ lies on the common boundary of two adjacent curvilinear triangles
$\fT_A$ and $\fT_B$, corresponding to regions $A$ (on which $\nabla H^*=(0,0)$) and $B$ (on which $\nabla H^*=(1,0)$) in
\Cref{f:vertical_tangent1,f:vertical_tangent2,f:vertical_tangent3,f:vertical_tangent4}.
For instance, in \Cref{f:vertical_tangent1}, locally near $(x_0,s_0)$, every point
$(x,s)\in\fT_A$ satisfies $x\le x_0$, whereas $\fT_B$ may also contain nearby points with $x<x_0$.
This distinction will matter for the double-integral representation \eqref{e:all_term0}, where we encounter
two types of integrals:
\begin{equation}\label{e:twoint}
\int_{\sfC^{\rm d}} P_{ns}(nw,nx)\, I_i(w)\,(\cdots)\,\rd w,
\qquad
\int_{\sfC^{\rm a}} Q_{ns}(nw,nx)\, I_i(w)^{-1}\,(\cdots)\,\rd w .
\end{equation}
By \Cref{c:Iiproperty} (with $b_i=x_0$), when $x>x_0$ the integrand in the first integral have poles
for $w\in[x_0,x]$, while when $x<x_0$ the integrand in the second integral have poles for
$w\in[x,x_0]$.

Depending on the sign of $S''(w_0;x_0,s_0)$ (equivalently, the sign of $\chi''(w_0)$) and the local geometry,
we distinguish several cases for the choice of steepest-descent paths:
\begin{enumerate}
\item Suppose $S''(w_0;x_0,s_0)>0$ and locally near $(x_0,s_0)$, every point
$(x,s)\in\fT_A$ satisfies $x\le x_0$, whereas $\fT_B$ may also contain nearby points with $x>x_0$
(see \Cref{f:vertical_tangent1}).

\begin{itemize}
\item For $\fT_A$, we take
\begin{align}\begin{split}\label{e:tangent_contour1}
\sfC^{\rm d}(w_0)&= \{\,w_0\pm r\fc \ri +(1-r)\fc/2: -1\le r\le 1\,\}
\cup \{\,w_0+\fc e^{\ri \theta}/4: 0\le \theta<2\pi\,\},\\
\sfC^{\rm a}(w_0)&= \{\,w_0+\fc e^{\ri \theta}/3: 0\le \theta<2\pi\,\}.
\end{split}\end{align}

\item For $\fT_B$, we take
\begin{align}\begin{split}\label{e:tangent_contour2}
\sfC^{\rm d}(w_0)&= \{\,w_0\pm r\fc \ri +(1-r)\fc/2: -1\le r\le 1\,\}
\cup \{\,w_0+\fc e^{\ri \theta}/3: 0\le \theta<2\pi\,\},\\
\sfC^{\rm d}(w_0)&= \{\,w_0+\fc e^{\ri \theta}/4: 0\le \theta<2\pi\,\}.
\end{split}\end{align}
\end{itemize}

In both subcases, $\sfC^{\rm a}(w_0)$ is a small circle around $w_0$, while
$\sfC^{\rm d}(w_0)$ consists of two pieces that can be deformed together into
a path from $w_0-\fc\ri$ to $w_0+\fc\ri$ passing to the right of $w_0$. This
right-passing condition is crucial: if $(x,s)\in\fT_B$ with $x>x_0$, then the
tangent line from $(x,s)$ to the portion of the arctic boundary contained in
$\fT_B$ produces a \emph{descent critical point} $w_c$ lying to the right of
$w_0$; recall from \Cref{s:descent_ascent_critical_points} that such a critical point issues two
non-real steepest-descent paths. Hence $\sfC^{\rm d}(w_0)$ can be further
deformed to the steepest-descent paths passing through $w_c$. On the other
hand, the integrand in \eqref{e:twoint} has poles confined to a small
neighborhood of $w_0$, so if $\sfC^{\rm d}(w_0)$ passed to the left of $w_0$,
then such a deformation to the steepest-descent paths through $w_c$ would not
be possible without crossing these poles.

The only difference between \eqref{e:tangent_contour1} and
\eqref{e:tangent_contour2} is the relative nesting of $\sfC^{\rm d}(w_0)$ and
$\sfC^{\rm a}(w_0)$. This matters for the double integral \eqref{e:all_term0},
whose integrand has poles when $w=z$. The nesting relation is determined as
follows. When $(x,s)\in\fT_A$ with $x\leq x_0$, by \Cref{c:PIproperty}, the integrand in the first term
of \eqref{e:twoint} has no poles in the relevant neighborhood. Therefore, the
circular piece of $\sfC^{\rm d}(w_0)$, namely the blue circle in Panel (A) of
\Cref{f:vertical_tangent1}, can be contracted to the empty contour without
crossing $\sfC^{\rm a}(w_0)$, namely the red circle. When $(x,s)\in\fT_B$ with
$x\leq x_0$, the tangent line from $(x,s)$ to the portion of the arctic
boundary contained in $\fT_B$ produces a descent critical point $w_c$ lying to
the right of $w_0$; again, recall from \Cref{s:descent_ascent_critical_points} that it issues two
non-real steepest-descent paths. Deforming $\sfC^{\rm d}(w_0)$ in Panel (B) of
\Cref{f:vertical_tangent1} so that it passes through $w_c$ does not cross
$\sfC^{\rm a}(w_0)$, the red circle.

\item If $S''(w_0;x_0,s_0)>0$ and locally near $(x_0,s_0)$, every point
$(x,s)\in\fT_B$ satisfies $x\le x_0$, whereas $\fT_A$ may also contain nearby points with $x>x_0$,
we choose the local paths for $\fT_A$ and $\fT_B$ as in \Cref{f:vertical_tangent2}.

\item If $S''(w_0;x_0,s_0)<0$ and locally near $(x_0,s_0)$, every point
$(x,s)\in\fT_A$ satisfies $x\ge x_0$, whereas $\fT_B$ may also contain nearby points with $x<x_0$,
we choose the local paths for $\fT_A$ and $\fT_B$ as in \Cref{f:vertical_tangent3}.

\item If $S''(w_0;x_0,s_0)<0$ and locally near $(x_0,s_0)$, every point
$(x,s)\in\fT_B$ satisfies $x\ge x_0$, whereas $\fT_A$ may also contain nearby points with $x<x_0$,
we choose the local paths for $\fT_A$ and $\fT_B$ as in \Cref{f:vertical_tangent4}.
\end{enumerate}

\begin{remark}
If $(x_0,s_0)$ lies on a vertical side of $\fP$, then the two cases shown in
\Cref{f:vertical_tangent1} and \Cref{f:vertical_tangent2} coincide, and the corresponding
contour prescriptions are equivalent. Indeed, in this case, locally around $(x_0,s_0)$,
every point $(x,s)\in\fT_A\cup\fT_B$ satisfies $x\le x_0$, so the integrand in the first
term of \eqref{e:twoint} has no poles in the relevant neighborhood. Consequently, the
contribution from the circular piece of $\sfC^{\rm d}(w_0)$, namely the blue circle in
Panel (A) of \Cref{f:vertical_tangent1} and Panel (B) of \Cref{f:vertical_tangent2},
vanishes by Cauchy's theorem.

Thus we can deform $\sfC^{\rm d}(w_0)$ in Panel (A) of \Cref{f:vertical_tangent1} to
$\sfC^{\rm d}(w_0)$ in Panel (A) of \Cref{f:vertical_tangent2}, and similarly deform
$\sfC^{\rm d}(w_0)$ in Panel (B) of \Cref{f:vertical_tangent1} to $\sfC^{\rm d}(w_0)$
in Panel (B) of \Cref{f:vertical_tangent2}. An analogous argument shows that, in the same
setting, the contours in \Cref{f:vertical_tangent3} and \Cref{f:vertical_tangent4} are also
equivalent.
\end{remark}

In the following we orient the local paths \eqref{e:tangent_contour1} and \eqref{e:tangent_contour2}, so they are compatible with the contours in the liquid region as introduced in \Cref{s:critical_bulk}. 
For $(x,s)\in \fL$ close to $(x_0, s_0)$, let $w_c=x-s\chi(w_c)\in\bC_+$,  we recall the factorization for descent vector from \eqref{eq:factor}:
\begin{align}\label{e:vertical_tangent_dir}
v^{\rm d}(w_c)=\sqrt{-\frac{1}{S''(w_c;x,s)}}
=\frac{\sqrt{\chi(w_c)}\sqrt{1-\chi(w_c)}}{\sqrt{1/s+\chi'(w_c)}}
=s^{1/2}\sqrt{\chi(w_c)}\sqrt{1-\chi(w_c)}\cdot
\frac{\sqrt{\partial_x\phi(x,s)}}{\sqrt{\phi'(w_c)}}.
\end{align}
Choosing a branch of $\sqrt{\phi'(w_c)}$ is equivalent to fixing the oriented descent direction $v^{\rm d}(w_c)$ and, by analytic continuation, determines the branch of $\sqrt{\phi'(w)}$ in \eqref{e:all_term0}. This choice 

We next give a geometric description of this choice as $(x,s)$ approaches
$(x_0,s_0)$ and $w_c$ approaches $w_0$.  Let $w_c=w_0+a+\ri b$, then \eqref{e:curve_reg} gives
\begin{align}
&1/s+\chi'(w_c)=\chi''(w_0)b\ri +\OO(b(|a|+b)), \\
& \chi(w_c)(1-\chi(w_c))=(1+\OO(|a|+b))\chi'(w_0)(w_c-w_0)=-(1+\OO(|a|+b))(a+\ri b)/s_0,
\end{align}
where in the second line we used $\chi(w_0)=0$ and $\chi'(w_0)=-1/s_0$.
Thus as $(x,s)\rightarrow (x_0, s_0)$,  $v^{\rm d}(w_c)=\sqrt{-1/S''(w_c;x,s))}$ is given by
\begin{align}\label{e:tangent_direction}
(1+\oo(1)) \sqrt{-\frac{(a+\ri b)}{ s_0b \chi''(w_0) \ri}}=\frac{1+\oo(1)}{\sqrt{2b}} \sqrt{\frac{a\ri -b}{S''(w_0;x_0, s_0)}}.
\end{align}
where we used \eqref{e:vert_tangent}.

For $S''(w_0;x_0,s_0)>0$,  then the direction of the square root in~\eqref{e:tangent_direction} lies in
$
\{\pm e^{\ri\theta}:\ \pi/4<\theta<3\pi/4\}.
$
If it lies in $\{e^{\ri\theta}:\ \pi/4<\theta<3\pi/4\}$, the corresponding direction is illustrated by the blue paths in \Cref{f:vertical_tangent1,f:vertical_tangent2};
if it lies in $\{-e^{\ri\theta}:\ \pi/4<\theta<3\pi/4\}$, the corresponding direction should be reversed.

For $S''(w_0;x_0,s_0)<0$,  then the direction of the square root in~\eqref{e:tangent_direction} lies in
$
\{\pm e^{\ri\theta}:\ -\pi/4<\theta<\pi/4\}.
$
If it lies in $\{-e^{\ri\theta}:\ -\pi/4<\theta<\pi/4\}$, the corresponding direction is illustrated by the blue paths in \Cref{f:vertical_tangent3,f:vertical_tangent4};
if  it lies in $\{e^{\ri\theta}:\ -\pi/4<\theta<\pi/4\}$, the corresponding direction should be reversed.

Finally, we orient the ascent contour $\sfC^{\rm a}(w_0)$ as indicated by the red paths in
\Cref{f:vertical_tangent1,f:vertical_tangent2,f:vertical_tangent3,f:vertical_tangent4}.
With this convention, if we deform the upward piece to a segment passing through $w_0$,
then in the upper half-plane the ordered pair
$
\bigl(\sfC^{\rm d}(w_0),\sfC^{\rm a}(w_0)\bigr)$
has a negative intersection at $w_0$, while in the lower half-plane it has a positive
intersection at $w_0$, in the sense of \Cref{d:positive_negative}.

\begin{figure}
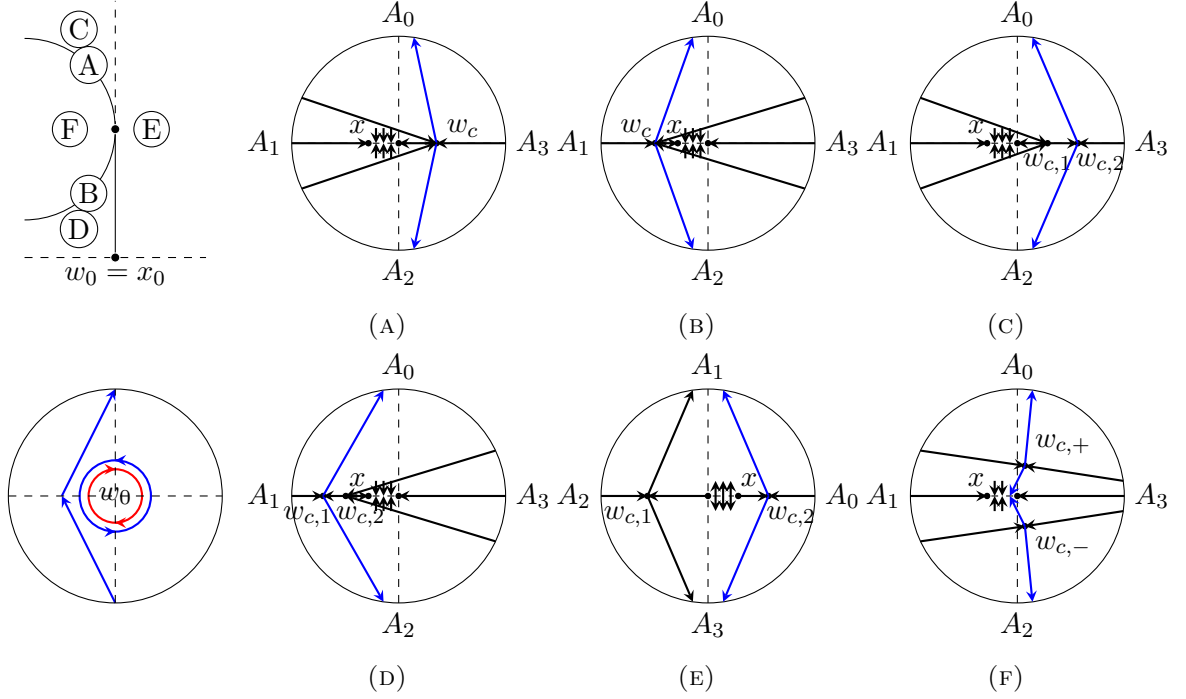

			\begin{subfigure}[t]{0.18\textwidth}

		\centering

				
			\caption{}
			\end{subfigure}
			\begin{subfigure}[t]{0.24\textwidth}
			
			%
				
			\caption{}
			\end{subfigure}
			\begin{subfigure}[t]{0.24\textwidth}
			
			%
				
			\caption{}
			\end{subfigure}

				\begin{subfigure}[t]{0.18\textwidth}
			\centering
			%
				
			\caption{}
			\end{subfigure}
	\begin{subfigure}[t]{0.24\textwidth}
			
			%
			\caption{}
			\end{subfigure}
					
				\caption{
				Gradient flow of \(S(\cdot;x,s)\) in a vertical tangent chart, with \(d=S''(w_0;x_0,s_0)/2<0\).}
	\label{f:tangent1}
\end{figure}

The following lemma shows that the local descent and ascent paths can be deformed into steepest-descent and steepest-ascent paths with negligible error. Its proof is deferred to \Cref{s:vertical_tangent_chart_proof}.
\begin{lemma}\label{l:vertical_tangent_steepest}
Adopt the assumptions and notation of \Cref{c:tangent_critical1}. For each critical
point \(w_c\), let \(\mathsf D^{\rm d}(w_c)\) denote the portions of the non-real
steepest--descent trajectories of \(S(\,\cdot\,;x,s)\) emanating from \(w_c\), stopped
upon their first exit from the disk
\[
\{w: |w-w_0|\le \fc\};
\]
see \Cref{f:tangent1} and \Cref{f:tangent2}. If \(w_c\) is not a descent critical
point, set
\[
\mathsf D^{\rm d}(w_c)=\emptyset .
\]
The length of each path \(\mathsf D^{\rm d}(w_c)\) is \(\OO(1)\).

The following alternatives hold.

\begin{enumerate}
  \item If \((x,s)\in\fL\), then there are two complex-conjugate critical
  points. In this case, the contour
  \(\mathsf C^{\rm d}(w_0)\) can be deformed to the union of the local
  steepest--descent paths \(\mathsf D^{\rm d}(w_c)\) associated with the two
  complex critical points, together with finitely many additional arcs of total
  length \(\OO(1)\).

  \item If \((x,s)\in\fP\setminus\fL\), then there are two real critical points,
  counted with multiplicity. Among them, there is at most one descent critical
  point \(w_c\) assigned to \((x,s)\) by
  \Cref{p:associate_critical_points}. If such a point exists, then
  \(\mathsf C^{\rm d}(w_0)\) can be deformed to
  \(\mathsf D^{\rm d}(w_c)\), together with finitely many additional arcs of
  total length \(\OO(1)\). If no such descent critical point exists, then
  \(\mathsf C^{\rm d}(w_0)\) can be deformed to the empty contour.
\end{enumerate}
In both cases, on these extra arcs,
\[
n\,\Re \bigl[S(w;x,s)-S(w_c;x,s)\bigr]\le -\,n\fc',
\]

The analogous statements hold for the local steepest--ascent paths.
\end{lemma}

\begin{figure}
			
			\begin{subfigure}[t]{0.18\textwidth}

			\centering
			%
			\caption{}	

			\end{subfigure}
%			
%\begin{subfigure}[t]{0.24\textwidth}
%			
%			%
			\caption{}
			\end{subfigure}

				\caption{
		Gradient flow of \(S(\cdot;x,s)\) in a vertical tangent chart, with \(d=S''(w_0;x_0,s_0)/2>0\).}
	\label{f:tangent2}
\end{figure}

\begin{figure}
	\begin{subfigure}[t]{0.3\textwidth}
			\centering
			%
				\caption*{$S'''>0$}
			\end{subfigure}%
			\begin{subfigure}[t]{0.3\textwidth}
			\centering
			%
		\caption{}

	\end{subfigure}
	\begin{subfigure}[t]{0.3\textwidth}
			\centering
			%
		\caption{}

	\end{subfigure}

	\caption{
\label{f:c_vertical_cusp1}
Local path associated with cusp-turning locations, $S'''>0$.}
	\end{figure}

\begin{figure}
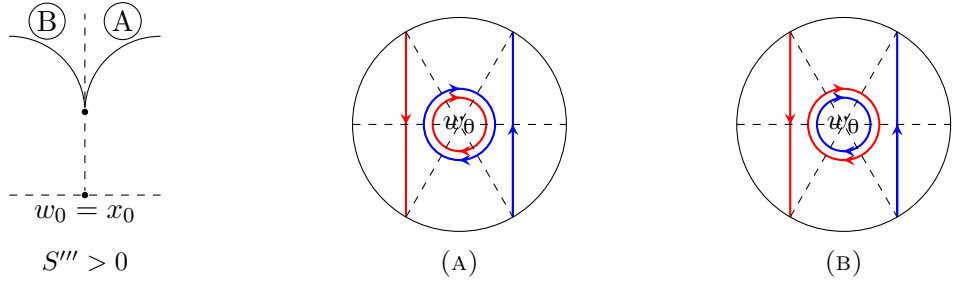
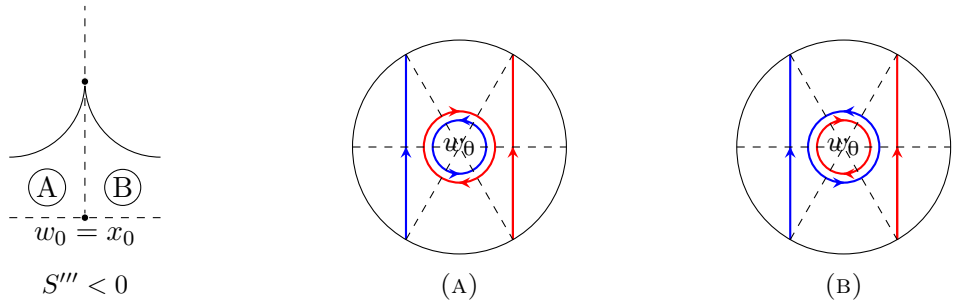
			
\begin{subfigure}[t]{0.3\textwidth}
			\centering
			%
				\caption*{$S'''<0$}
			\end{subfigure}%
			\begin{subfigure}[t]{0.3\textwidth}
			\centering
			%
		\caption{}

	\end{subfigure}
	\begin{subfigure}[t]{0.3\textwidth}
			\centering
			%
		\caption{}

	\end{subfigure}
		\caption{
\label{f:c_vertical_cusp2}
Local path associated with cusp-turning locations, $S'''<0$.}
	\end{figure}

If $(x_0,s_0)\in\fA$ is a cusp-turning location, then $c_1=0$ and $c_2\neq 0$, and
\begin{equation}\label{e:cusp_Sprime}
S'(w;x_0,s_0)= s_0c_2\,(w-x_0)^2+\OO\!\bigl((w-x_0)^3\bigr).
\end{equation}

The cusp-turning location $(x_0,s_0)$ lies on the common boundary of two adjacent curvilinear triangles
$\fT_A$ and $\fT_B$, corresponding to regions $A$ (where $\nabla H^*=(0,0)$) and $B$
(where $\nabla H^*=(1,0)$) in \Cref{f:c_vertical_cusp1,f:c_vertical_cusp2}.
Depending on the sign of $S'''(w_0;x_0,s_0)$ (equivalently, the sign of $\chi'''(w_0)$), we distinguish the
following cases for the choice of local steepest-descent paths.

\begin{enumerate}
\item \textbf{Case $S'''(w_0;x_0,s_0)>0$ (cusp pointing downward).}
See \Cref{f:c_vertical_cusp1}.

\begin{itemize}
\item For $\fT_A$, we take
\begin{align}\begin{split}\label{e:cusp_turning_contour1}
\sfC^{\rm d}(w_0)&=\{\,w_0+(1+\sqrt3\,r\ri)\tfrac{\fc}{2}:\,-1\le r\le 1\,\}
\;\cup\;\{\,w_0+\tfrac{\fc}{3}e^{\ri\theta}:\,0\le \theta<2\pi\,\},\\
\sfC^{\rm a}(w_0)&=\{\,w_0+(-1+\sqrt3\,r\ri)\tfrac{\fc}{2}:\,-1\le r\le 1\,\}
\;\cup\;\{\,w_0+\tfrac{\fc}{4}e^{\ri\theta}:\,0\le \theta<2\pi\,\}.
\end{split}\end{align}

\item For $\fT_B$, we take
\begin{align}\begin{split}\label{e:cusp_turning_contour2}
\sfC^{\rm d}(w_0)&=\{\,w_0+(1+\sqrt3\,r\ri)\tfrac{\fc}{2}:\,-1\le r\le 1\,\}
\;\cup\;\{\,w_0+\tfrac{\fc}{4}e^{\ri\theta}:\,0\le \theta<2\pi\,\},\\
\sfC^{\rm a}(w_0)&=\{\,w_0+(-1+\sqrt3\,r\ri)\tfrac{\fc}{2}:\,-1\le r\le 1\,\}
\;\cup\;\{\,w_0+\tfrac{\fc}{3}e^{\ri\theta}:\,0\le \theta<2\pi\,\}.
\end{split}\end{align}
\end{itemize}

In both subcases, $\sfC^{\rm d}(w_0)$ consists of two pieces that can be
deformed together into a path from $w_0+e^{-\pi\ri/3}\fc$ to
$w_0+e^{\pi\ri/3}\fc$ passing to the left of $w_0$, while
$\sfC^{\rm a}(w_0)$ consists of two pieces that can be deformed together into
a path from $w_0+e^{-2\pi\ri/3}\fc$ to $w_0+e^{2\pi\ri/3}\fc$ passing to the
right of $w_0$. Here we use the fact that all poles relevant to
\eqref{e:twoint} are confined to a small neighborhood of $w_0$.

Moreover, when $(x,s)\in\fT_A$, the tangent line from $(x,s)$ to the portion
of the arctic boundary contained in $\fT_A$ produces a \emph{descent critical
point} $w_c$ lying to the left of $w_0$; recall from \Cref{s:descent_ascent_critical_points} that
such a critical point issues two non-real steepest-descent paths. Hence
$\sfC^{\rm d}(w_0)$ can be further deformed to pass through $w_c$. Similarly,
when $(x,s)\in\fT_B$, the tangent line from $(x,s)$ to the portion of the
arctic boundary contained in $\fT_B$ produces an \emph{ascent critical point}
$w_c$ lying to the right of $w_0$, and hence $\sfC^{\rm a}(w_0)$ can be further
deformed to pass through $w_c$.

The only difference between \eqref{e:cusp_turning_contour1} and
\eqref{e:cusp_turning_contour2} is the relative nesting of
$\sfC^{\rm d}(w_0)$ and $\sfC^{\rm a}(w_0)$. This matters for the double
integral \eqref{e:all_term0}, whose integrand has a pole along the diagonal
$w=z$. The nesting is determined by the following observations. If
$(x,s)\in\fT_A$ with $x\ge x_0$, then the integrand in the second term of
\eqref{e:twoint} is holomorphic inside the circular piece of
$\sfC^{\rm a}(w_0)$, namely the red circle in Panel (A) of
\Cref{f:c_vertical_cusp1}. Hence this circle can be contracted to a point
without crossing $\sfC^{\rm d}(w_0)$, the blue contour. Likewise, if
$(x,s)\in\fT_B$ with $x\le x_0$, then the integrand in the first term of
\eqref{e:twoint} is holomorphic inside the circular piece of
$\sfC^{\rm d}(w_0)$, namely the blue circle in Panel (B) of
\Cref{f:c_vertical_cusp2}. Thus this circle can be contracted to a point
without crossing $\sfC^{\rm a}(w_0)$, the red contour.

\item \textbf{Case $S'''(w_0;x_0,s_0)<0$ (cusp pointing upward).}
See \Cref{f:c_vertical_cusp2}. We choose the local paths for $\fT_A$ and $\fT_B$ as indicated in
\Cref{f:c_vertical_cusp2}.
\end{enumerate}

In the following we orient the local paths \eqref{e:cusp_turning_contour1} and \eqref{e:cusp_turning_contour2}, so they are compatible with the contours in the liquid  region as introduced in \Cref{s:critical_bulk}.

For $(x,s)\in \fL$ close to $(x_0, s_0)$, let $w_c=x-s\chi(w_c)\in\bC_+$. The same as in \eqref{e:vertical_tangent_dir}, a local choice of the branch of $\sqrt{\phi'(w_c)}$ determines the
oriented descent direction $v^{\rm d}(w_c)$, and conversely. By analytic
continuation, this choice also determines the branch of the prefactor
$\sqrt{\phi'(w)}$ appearing in the double-contour integral
\eqref{e:all_term0}.

We next give a geometric description of this choice as $(x,s)$ approaches
$(x_0,s_0)$ and $w_c$ approaches $w_0$. Let $w_c=w_0+a+\ri b$ with $b>0$, then \eqref{e:curve_cusp} give
\begin{align}
&1/s+\chi'(w_c)=\chi'''(w_0)b(a\ri-b/3) +\OO(b(|a|+b)^2),\\
&\chi(w_c)(1-\chi(w_c))=(1+\OO(|a|+b))\chi'(w_0)(w_c-w_0)=-(1+\OO(|a|+b))(a+\ri b)/s_0
\end{align}
where in the second line we used $\chi(w_0)=0$ and $\chi'(w_0)=-1/s_0$.

Thus, as $(x,s)\to(x_0,s_0)$,   
$v^{\rm d}(w_c)=\sqrt{-1/S''(w_c;x,s))}$ is given by 
\begin{align}\label{e:cusp_turning_direction} 
(1+\oo(1))\sqrt{-\frac{(a+\ri b)}{s_0 b(a\ri-b/3)\chi'''(w_0) }} =\frac{1+\oo(1)}{\sqrt{b(b^2/9+a^2)}}\sqrt{\frac{-2ab/3+(a^2+b^2/3)\ri}{3S'''(w_0;x_0,s_0)}}.
\end{align}
where we used \eqref{e:cusp_Sprime}

For $S'''(w_0;x_0,s_0)>0$,  then the direction of the square root in~\eqref{e:cusp_turning_direction} lies in
$
\{\pm e^{\ri\theta}:0<\theta<\pi/2\}.
$
If it lies in $\{e^{\ri\theta}:0<\theta<\pi/2\}$, the corresponding direction is illustrated by the blue paths in \Cref{f:c_vertical_cusp1};
if it lies in $\{-e^{\ri\theta}:0<\theta<\pi/2\}$, the corresponding direction should be reversed.

For $S'''(w_0;x_0,s_0)<0$,  then the direction of the square root in~\eqref{e:cusp_turning_direction} lies in
$
\{\pm e^{\ri\theta}:\ \pi/2<\theta<\pi\}.
$
If it lies in $\{e^{\ri\theta}:\ \pi/2<\theta<\pi\}$, the corresponding direction is illustrated by the blue paths in \Cref{f:c_vertical_cusp2};
if  it lies in $\{-e^{\ri\theta}:\ \pi/2<\theta<\pi\}$, the corresponding direction should be reversed.

Finally, we orient the ascent contour $\sfC^{\rm a}(w_0)$ as indicated by the
red paths in \Cref{f:c_vertical_cusp1,f:c_vertical_cusp2}. With this
convention, if the descent and ascent contours are each deformed into a single
path so that they cross, then in the upper half-plane the ordered pair
$
\bigl(\sfC^{\rm d}(w_0),\sfC^{\rm a}(w_0)\bigr)
$
has a negative intersection, while in the lower half-plane the
corresponding ordered pair has a positive intersection, in the sense
of \Cref{d:positive_negative}.

The following lemma show that the local descent and ascent paths can be deformed into steepest-descent and steepest-ascent paths with negligible error. Its proof is deferred to \Cref{s:vertical_tangent_chart_proof}.
\begin{lemma}\label{c:vertical_cusp_steepest}
Adopt the assumptions and notation of \Cref{c:cusp_turning_critical1}. For each critical
point \(w_c\), let \(\mathsf D^{\rm d}(w_c)\) denote the portions of the non-real
steepest--descent trajectories of \(S(\,\cdot\,;x,s)\) emanating from \(w_c\), stopped
upon their first exit from the disk
\[
\{w: |w-w_0|\le \fc\};
\]
see \Cref{f:vertical_cusp1} and \Cref{f:vertical_cusp2}. If \(w_c\) is not a descent critical
point, set
\[
\mathsf D^{\rm d}(w_c)=\emptyset .
\]
The length of each path \(\mathsf D^{\rm d}(w_c)\) is \(\OO(1)\).

The following alternatives hold.

\begin{enumerate}
  \item If \((x,s)\in\fL\), then there are two complex-conjugate critical
  points and one real critical point. In this case, the contour
  \(\mathsf C^{\rm d}(w_0)\) can be deformed to the union of the local
  steepest--descent paths \(\mathsf D^{\rm d}(w_c)\) associated with the two
  complex critical points, together with finitely many additional arcs of total
  length \(\OO(1)\).

  \item If \((x,s)\in\fP\setminus\fL\), then there are three real critical points,
  counted with multiplicity. Among them, there is one descent critical
  point \(w_c\) assigned to \((x,s)\) by
  \Cref{p:associate_critical_points}. In this case, the contour
  \(\mathsf C^{\rm d}(w_0)\) can be deformed to
  \(\mathsf D^{\rm d}(w_c)\), together with finitely many additional arcs of
  total length \(\OO(1)\). 
\end{enumerate}
In both cases, on these extra arcs,
\[
n\,\Re \bigl[S(w;x,s)-S(w_c;x,s)\bigr]\le -\,n\fc',
\]

The analogous statements hold for the local steepest--ascent paths.

\end{lemma}

\begin{figure}			
\begin{subfigure}[t]{0.24\textwidth}
			\centering

				
			\caption{}
			\end{subfigure}
	\begin{subfigure}[t]{0.24\textwidth}
			
			%
				
			\caption{}
			\end{subfigure}
		\begin{subfigure}[t]{0.24\textwidth}
			
			%
				
			\caption{}
			\end{subfigure}
	
	\begin{subfigure}[t]{0.24\textwidth}
			\centering
			%
				
			\caption{}
			\end{subfigure}
\begin{subfigure}[t]{0.24\textwidth}
			
			%
				
			\caption{}
			\end{subfigure}

	\caption{
\label{f:vertical_cusp1}
Gradient flow of \(S(\cdot;x,s)\) in a vertical cusp chart, with \(d=S'''(w_0;x_0,s_0)/6<0\).}
	\end{figure}

\begin{figure}
			\begin{subfigure}[t]{0.24\textwidth}
			\centering
			%
				
			\caption{}
			\end{subfigure}
			
			\begin{subfigure}[t]{0.24\textwidth}
			\centering
			%
				
			\caption{}
			\end{subfigure}
	\caption{
\label{f:vertical_cusp2}
Gradient flow of \(S(\cdot;x,s)\) in a vertical cusp chart, with \(d=S'''(w_0;x_0,s_0)/6>0\).
}
	\end{figure}

\subsection{Vertical tangent frozen chart}\label{s:vertical_frozen_neighborhood}
Suppose that $(x_0,s_0)\in \fP\setminus \fL$ lies on a vertical extended side tangent to the arctic boundary at the vertical tangency point $(x_0,s_0')\in \fA$. From \eqref{e:bcond}  we have
 $f(x_0, s'_0)=\chi(x_0;s'_0)=0$ and 
\begin{align}
 \quad w_0:=x_0-s_0'\chi(w_0)=x_0\in \cC(\bR), \quad \chi(w_0)=0. 
\end{align}
In a neighborhood of $w_0$, we denote $\chi(w)=(w_0-w)\chi_0(w)$, and the same as in \eqref{e:expandchi_b}, we have
\begin{align}
\begin{split}
\chi(w)
&=
-\frac{w-w_0}{s'_0}
+c_1(w-w_0)^2
+c_2(w-w_0)^3+\cdots,\\
\chi_0(w)
&=
\frac1{s'_0}
+c_1(w_0-w)
-c_2(w_0-w)^2+\cdots,
\end{split}
\end{align}
where $c_1=\chi''(w_0)/2$ and $c_2=\chi'''(w_0)/6$.
 It follows
\begin{align}\begin{split}\label{e:S_extend}
S'(w;x_0, s_0)
&=\ln \frac{w-w_0}{w_0-s_0-w}-\ln \frac{\chi(w)}{1-\chi(w)}
=\ln \frac{(1-\chi(w))(w-w_0)}{\chi(w)(w_0-s_0-w)}
=\ln \frac{(1-(w_0-w)\chi_0(w))}{\chi_0(w)(w-w_0+s_0)}\\
&=\ln \frac{s'_0}{s_0}+\ln \frac{1+(w-w_0)/s'_0+\OO((w-w_0)^2)}{(1-s'_0a_1(w-w_0)+\OO((w-w_0)^2))(1+(w-w_0)/s_0)}\\
&=\ln \frac{s'_0}{s_0}+\left(\frac{1}{s'_0}-\frac{1}{s_0}+s'_0 a_1\right)(w-w_0)+\OO((w-w_0)^2).
\end{split}\end{align}

The following lemmas introduce the vertical frozen chart and  record several of its properties. Their proof, based on a Taylor expansion of the tiling action, is deferred to \Cref{s:vertical_frozen_chart_proof}.
\begin{lemma}[$\fc$-frozen chart (vertical tangent)]\label{c:vertical_frozen_critical1}
Fix $\delta_0>0$ and let $(x_0,s_0)\in \fP$.  Assume there exists a vertical tangent location $(x_0',s_0')\in \fA$ such that
the line through $(x_0,s_0)$ and $(x_0',s_0')$ is tangent to the arctic curve at $(x_0',s_0')$, and moreover $(x_0, s_0)$ is at least distance $\delta_0$ from $(x_0',s_0')$. Then for every sufficiently small $\fc>0$ there exists
$\delta=\delta(\fc)>0$, also sufficiently small and depending on $\delta_0$,  such that the following holds.

Let $w_0=x_0-s_0\,\chi(x_0,s_0)=x_0$ be the critical point corresponding to $(x_0,s_0)$, and set
\[
\fU:=\{w\in\bC:\ |w-w_0|\le \fc\}.
\]
On $\fU$ the Riemann surface $\cC$ can be parametrized as $(f(w),w)$. We will therefore identify
$\fU$ with its image in $\cC$ under the map $w\mapsto (f(w), w)$.   Moreover, for all $w\in\fU$,
\begin{align}\label{e:chi_v_tangent}
  |1-\chi(w)| \asymp 1,\quad  |\chi(w)|\asymp |w-w_0|.
\end{align}

Now let $(x,s)\in\fL$ satisfy $\|(x,s)-(x_0,s_0)\|_2\le \delta$.
Then the tiling action $S(\,\cdot\,;x,s)$ has exactly one (formal) critical point $w_c$ inside $\fU$, satisfying 
\begin{align}\label{e:wcest}
\dist(w_c, [\min\{x_0, x\}, \max\{x_0, x\}])|\asymp |x-x_0|,
\end{align}
Moreover, for every $w\in\fU$, 
\begin{align}\label{e:tangent_frozen_S}
S(w;x,s)-S(w_c;x,s)=d(w-w_0) +\cE(w),\quad d:=S'(w_0;x_0,s_0), 
\end{align}
The error term satisfies
\begin{align}\label{e:tangent_frozen_err}
|\cE(w)|\leq C(\delta\ln(1/\delta)+|w-w_0|^2)\leq\frac{|d|\fc}{100}.
\end{align}

In this situation, we call $\fU$ a (vertical tangent) $\fc$-frozen chart (centered at $w_0$), and we say that
the point $(x,s)$ is \emph{adapted} to $\fU$.
\end{lemma}

\begin{figure}
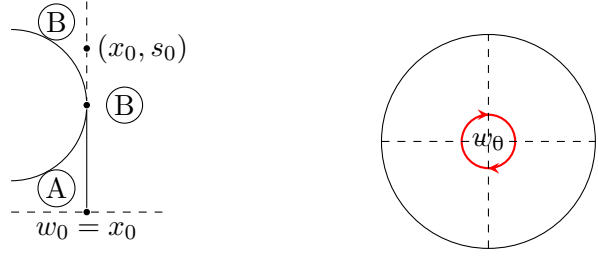
					
	\begin{subfigure}[t]{0.30\textwidth}
			\centering

		%\caption{}

	\end{subfigure}

	\caption{
\label{f:tangent_frozen3}
Local path associated with $(x_0, s_0)\in \fP\setminus \fL$ which lies in (or on the boundary of) a single curvilinear triangle $\fT=\fT_B$.}
	\end{figure}

\begin{figure}
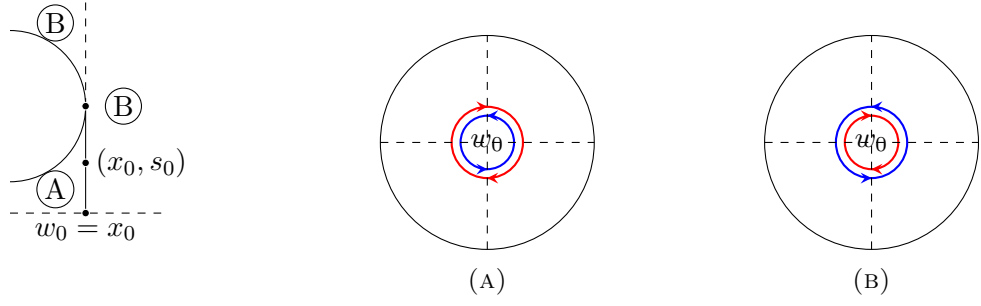
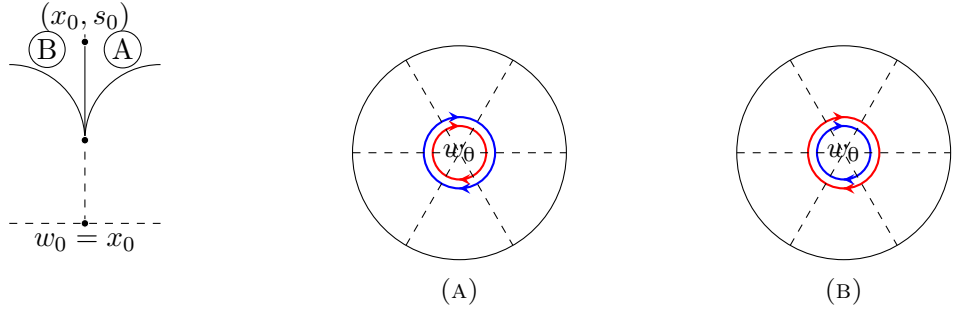

	\begin{subfigure}[t]{0.3\textwidth}
			\centering
			%
		\caption{}

	\end{subfigure}
	\begin{subfigure}[t]{0.3\textwidth}
			\centering
			%
		\caption{}

	\end{subfigure}

	\caption{
\label{f:tangent_frozen2}
Local path associated with $(x_0, s_0)\in \fP\setminus \fL$ which lies on the common boundary of two adjacent curvilinear triangles.}
	\end{figure}

Next, we describe the integration contours associated with vertical-tangent
$\fc$-frozen charts. The corresponding tangent or cusp-turning point
$(x_0,s_0')$ lies on the common boundary of two adjacent curvilinear
triangles $\fT_A$ and $\fT_B$; see
\Cref{f:adjacent_curvilinear_triangle}. There are two cases.
\begin{enumerate}
\item
Suppose that $(x_0,s_0)$ does not lie on the common boundary of
$\fT_A$ and $\fT_B$. Let $\fT\in\{\fT_A,\fT_B\}$ denote the region
containing $(x_0,s_0)$; see \Cref{f:tangent_frozen3}.

The point $(x_0,s_0')$ divides the tangent line through it into the
descent and ascent cuts $\ell_-(w_0;\fT)$ and $\ell_+(w_0;\fT)$
defined in \Cref{d:exceptional_cuts}. Suppose first that
\[
(x_0,s_0)\in\ell_-(w_0;\fT).
\]
In this case, if we perturb $(x_0,s_0)$ slightly to a nearby point
$(\widehat x_0,\widehat s_0)$ away from the extended side, chosen so
that the corresponding tangent line remains in $\fT$, then the
associated critical point $\widehat w_0$ is a descent critical point.

In \Cref{s:vertical_tangent}, we associated local descent and ascent
contours with the tangency or cusp-turning point $(x_0,s_0')\in\fA$,
each containing a small circular component surrounding $w_0$. In the
present frozen setting, we associate with $(x_0,s_0)$ the circular
component of the local descent contour (blue) associated with
$(x_0,s_0')$. 

If instead
\[
(x_0,s_0)\in\ell_+(w_0;\fT),
\]
we analogously associate with $(x_0,s_0)$ the circular component of
the local ascent contour (red); see \Cref{f:tangent_frozen3}.

\item
Suppose that $(x_0,s_0)$ lies on the common boundary of the two
adjacent curvilinear triangles:
\[
(x_0,s_0)\in\fT_A\cap\fT_B;
\]
see \Cref{f:tangent_frozen1,f:tangent_frozen2}. On this shared
boundary, the descent cut for $\fT_A$ is the ascent cut for $\fT_B$,
and vice versa. Thus $(x_0, s_0)$ belongs to both descent and ascent cuts. 

We therefore associate with $(x_0,s_0)$ both the circular component
of the local descent contour (blue) and that of the local ascent
contour (red) associated with $(x_0,s_0')$ in
\Cref{s:vertical_tangent}, for $\fT_A$ or $\fT_B$.

For example, in \Cref{f:tangent_frozen1}, the point $(x_0,s_0')$ is a
tangency point. The contours for $\fT_A$ and $\fT_B$ in Panels (A)
and (B), respectively, are obtained from the corresponding panels of
\Cref{f:vertical_tangent1} by retaining only the circular components.

Similarly, in \Cref{f:tangent_frozen2}, the corresponding point is a
cusp-turning point. The contours for $\fT_A$ and $\fT_B$ in Panels
(A) and (B), respectively, are obtained from the corresponding panels
of \Cref{f:c_vertical_cusp1} by retaining only the circular
components.
\end{enumerate}

\begin{figure}					
	\begin{subfigure}[t]{0.30\textwidth}
			\centering
			%
		\caption{}

	\end{subfigure}
	\begin{subfigure}[t]{0.3\textwidth}
			\centering
			%
		\caption{}

	\end{subfigure}

	\caption{
\label{f:tangent_frozen1}
Local path associated with $(x_0, s_0)\in \fP\setminus \fL$ which lies on the common boundary of two adjacent curvilinear triangles.}
	\end{figure}

\begin{figure}
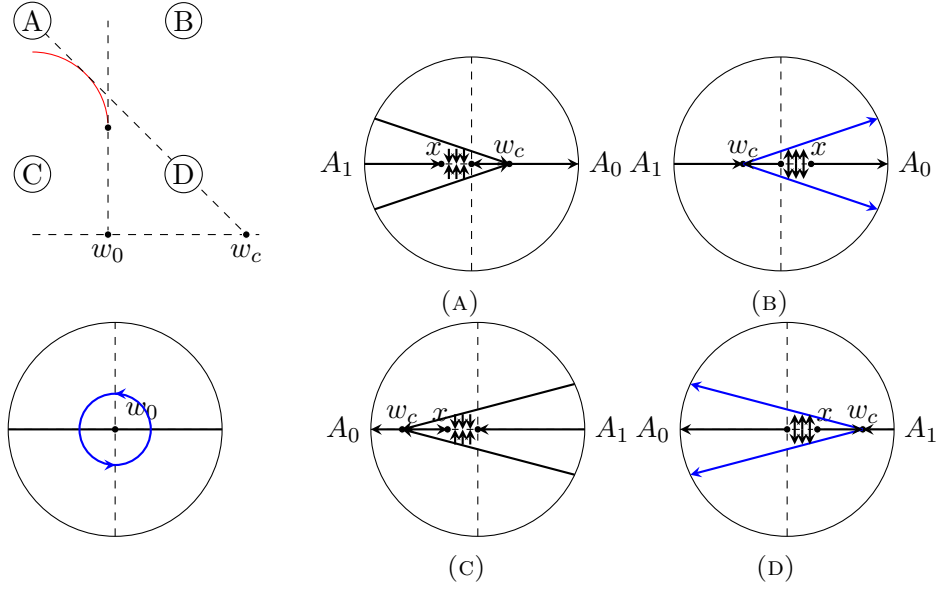

		\begin{subfigure}[t]{0.23\textwidth}

				
			\caption{}
			\end{subfigure}
		\begin{subfigure}[t]{0.24\textwidth}
			
			%
				
			\caption{}
			\end{subfigure}
			
	\begin{subfigure}[t]{0.24\textwidth}
			
			%
				
			\caption{}
			\end{subfigure}
		\begin{subfigure}[t]{0.24\textwidth}
			
			%
				
			\caption{}
			\end{subfigure}
	\caption{
	Gradient flow structure near a nondegenerate saddle (\(m=2\) in
\Cref{prop:exit}), which is correspond to a vertical tangent location.}
	\label{f:side}
\end{figure}

The following lemma shows that the local descent and ascent paths can be deformed into steepest-descent and steepest-ascent paths with negligible error. Its proof is deferred to \Cref{s:vertical_frozen_chart_proof}.
\begin{lemma}\label{l:local_descent_deformation}
Adopt the assumptions and notation of
\Cref{c:vertical_frozen_critical1}. There exists a constant $\fc'>0$ such that the following statements hold.
Suppose that the tangent frozen chart
contains local descent paths, that is,
\[
    \sfC^{\rm d}(w_0)\neq \emptyset.
\]

For the critical point \(w_c\), let \(\mathsf D^{\rm d}(w_c)\) denote the
union of the non-real steepest--descent trajectories of
\(S(\,\cdot\,;x,s)\) emanating from \(w_c\), each stopped upon its first exit
from the disk
\[
    \{w: |w-w_0|\le \fc\};
\]
see \Cref{f:side}.

If \(w_c\) is not a descent critical point, then
\(\mathsf D^{\rm d}(w_c)=\emptyset\), and the local contour
\(\sfC^{\rm d}(w_0)\) can be removed by deformation; equivalently, it can be
deformed to
\[
    \mathsf D^{\rm d}(w_c)=\emptyset.
\]

If \(w_c\) is a descent critical point, then
\(\mathsf D^{\rm d}(w_c)\neq \emptyset\) and has total length $\OO(1)$, and
\(w_c\) is assigned to
\((x,s)\) by \Cref{p:associate_critical_points}.
\begin{itemize} 
\item \emph{Replacing \(\sfC^{\rm d}(w_0)\) by \(\sfD^{\rm d}(w_c)\).}
The contour
\(\sfC^{\rm d}(w_0)\) can be deformed to
$\mathsf D^{\rm d}(w_c)$, together with several arcs of total length $\OO(1)$.
Along these auxiliary arcs,
\[
    n\,\Re [S(w;x,s)-S(w_c;x,s)]
    \le -\,n\fc'
\]
for some constant \(\fc'>0\).

\item \emph{Replacing \(\sfD^{\rm d}(w_c)\) by \(\sfS^{\rm d}(w_c)\).}
We distinguish two regimes:
\begin{enumerate}

\item \emph{Close to the tangent line.} If
$|x-x_0|\leq (\ln n)^3/n$, we define
$r_n=(\ln n)^5/n$.

\item \emph{Other frozen regime.} If
$|x-x_0|>(\ln n)^3/n$, we define
$r_n:=(\ln n)^2 (|x-x_0|/n)^{1/2}$.

\end{enumerate}
We define the truncated steepest--descent set
$\mathsf S^{\rm d}(w_c)$ at \(w_c\) to be the portions of the
steepest--descent trajectories starting at \(w_c\) up to their first exit from
the disk
\[
\{w: |w-w_c|\leq r_n\}.
\]

Then, for all
\[
w\in \mathsf D^{\rm d}(w_c)\setminus \mathsf S^{\rm d}(w_c),
\]
we have
\[
e^{n\Re[S(w;x,s)]}
\leq
e^{n\Re[S(w_c;x,s)]}e^{-\fc'(\ln n)^5}.
\]
\end{itemize}

The analogous statements hold for the local steepest--ascent paths and
their truncations.
\end{lemma}

\subsection{Unit-slope tangent chart and cusp-turning chart}\label{s:unit_slope_tangent}
\begin{figure}
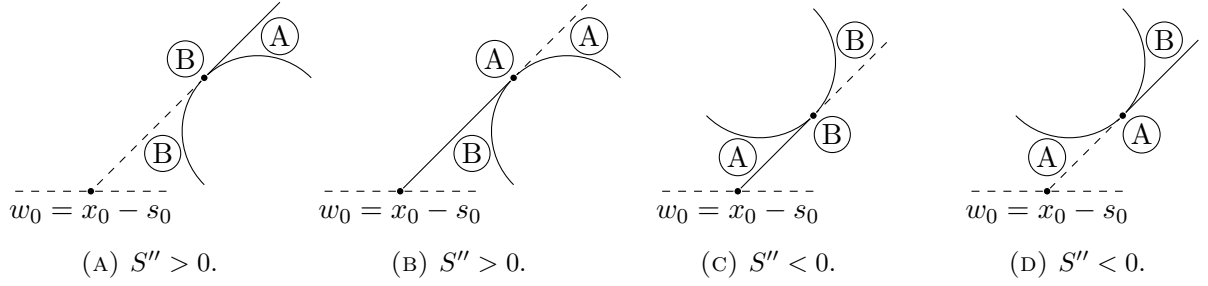
	
\begin{subfigure}[t]{0.24\textwidth}
			\centering

			\caption{$S''>0$.}
		\end{subfigure}			
\begin{subfigure}[t]{0.24\textwidth}
			\centering
			%
			\caption{$S''>0$.}
		\end{subfigure}		
\begin{subfigure}[t]{0.24\textwidth}
			\centering
			%
			\caption{$S''<0$.}
		\end{subfigure}			
	\begin{subfigure}[t]{0.24\textwidth}
			\centering
			%
			\caption{$S''<0$.}
		\end{subfigure}

	\caption{
\label{f:unit_slope_tangent}
Local path associated with unit-slope tangent locations}
	\end{figure}

Let $(x_0,s_0)\in \fA$ be a unit-slope tangent location. Then, by
\eqref{e:bcond} and \eqref{e:arctic_chider}, we have
$f(x_0,s_0)=\infty$ and
\begin{align}
w_0:=x_0-s_0\chi(x_0,s_0)=x_0-s_0,\quad \chi(w_0)=1,\quad \chi'(w_0)=-1/s_0.
\end{align}
In this case, $(x_0,s_0)$ lies on the common boundary of two adjacent
curvilinear triangles $\fT_A$ and $\fT_B$; see \Cref{f:unit_slope_tangent}.

In a neighborhood of $w_0$, write
$
1-\chi(w)=(w-w_0)\chi_0(w)$.
Taking derivatives gives
\begin{align}
-\frac{1}{s_0}
=\chi'(w_0)
=\left[-\chi_0(w)-(w-w_0)\chi_0'(w)\right]_{w=w_0}
=-\chi_0(w_0),
\end{align}
and hence $\chi_0(w_0)=1/s_0>0$. Thus, locally around $w_0$, we have the
expansions
\begin{align}\begin{split}\label{e:expandchi_a}
1-\chi(w)
&=\frac{w-w_0}{s_0}+c_1(w-w_0)^2+c_2(w-w_0)^3+\cdots,\\
\chi_0(w)
&=\frac{1}{s_0}+c_1(w-w_0)+c_2(w-w_0)^2+\cdots,
\end{split}\end{align}
where $c_1=-\chi''(w_0)/2$ and $c_2=-\chi'''(w_0)/6$. Moreover,
\begin{align}\begin{split}\label{e:unit_slope_tangent}
S'(w;x_0,s_0)
&=\ln \frac{w-x_0}{x_0-s_0-w}
 -\ln \frac{\chi(w)}{1-\chi(w)} =\ln \frac{(1-\chi(w))(w-w_0-s_0)}{\chi(w)(w_0-w)}
=\ln \frac{\chi_0(w)(s_0-(w-w_0))}
{1-(w-w_0)\chi_0(w)}\\
&=\ln \frac{
\bigl(1+s_0c_1(w-w_0)+s_0c_2(w-w_0)^2+\OO((w-w_0)^3)\bigr)
\bigl(1-(w-w_0)/s_0\bigr)}
{1-(w-w_0)/s_0-c_1(w-w_0)^2+\OO((w-w_0)^3)}\\
&=s_0c_1(w-w_0)
+\left(s_0c_2-\frac{1}{2}s_0^2c_1^2+c_1\right)(w-w_0)^2
+\OO((w-w_0)^3).
\end{split}\end{align}
In particular, \(S'(w;x_0,s_0)\) is analytic in a neighborhood of \(w_0\).
If $(x_0,s_0)$ is not a cusp-turning location, then
$c_1=-\chi''(w_0)/2\neq 0$, and
$
S''(w_0;x_0,s_0)=s_0c_1=-s_0\chi''(w_0)/2$.
If $(x_0,s_0)$ is a cusp-turning location, then $\chi''(w_0)=0$ while
$c_2=-\chi'''(w_0)/6\neq 0$, and
$
S'''(w_0;x_0,s_0)=2s_0c_2=-s_0\chi'''(w_0)/3$.
Thus the relevant signs are determined by the signs of $\chi''(w_0)$ and
$\chi'''(w_0)$, as classified in \Cref{f:tangent} and
\Cref{f:cusp_turning}.

For \((x,s)\) in a sufficiently small neighborhood of \((x_0,s_0)\), we have
\begin{align}\label{e:vert_tangent_diff2}
S'(w;x,s)-S'(w;x_0,s_0)
=\ln \frac{w-x}{x-s-w}
 -\ln \frac{w-x_0}{x_0-s_0-w} =\ln \frac{x_0-s_0-w}{x-s-w}
 +\ln \frac{w-x}{w-x_0}.
\end{align}
The first logarithm has a branch cut along the segment
\[
[\min\{x_0-s_0,x-s\},\,\max\{x_0-s_0,x-s\}].
\]
Since the second logarithm is analytic in a sufficiently small neighborhood of
\(w_0=x_0-s_0\), the difference \eqref{e:vert_tangent_diff2} is locally
analytic near \(w_0\) on the complement of this segment.

The following lemma introduces the unit-slope tangent chart and  record several of its properties.The statement can be proven in the same way as in the vertical tangent case \Cref{c:tangent_critical1}, so we omit.
\begin{lemma}[$\fc$-tangent chart (unit-slope)]\label{c:tangent_critical2}
Given a unit-slope tangent location $(x_0, s_0)\in \fA$ which is not a cusp location. Then for every $\fc>0$ sufficiently small  there exists
$\delta=\delta(\fc)>0$ sufficiently small such that the following holds.

Let $w_0=x_0-s_0\,\chi(x_0,s_0)=x_0-s_0$ be the critical point corresponding to $(x_0,s_0)$, and set
\[
\fU:=\{w\in\bC:\ |w-w_0|\le \fc\}.
\]
On $\fU$ the Riemann surface $\cC$ can be parametrized as $(f(w),w)$. We will therefore identify
$\fU$ with its image in $\cC$ under the map $w\mapsto (f(w),w)$. Moreover, for all $w\in\fU$,
\begin{align}
  |\chi(w)|, |\chi''(w)|\asymp 1,\quad  |1-\chi(w)|\asymp |w-w_0|.
\end{align}

Now let $(x,s)\in\fL$ satisfy $\|(x,s)-(x_0,s_0)\|_2\le \delta$.
Then the tiling action $S(\,\cdot\,;x,s)$ has exactly two (formal) critical points $w_c$ inside $\fU$, satisfying 
\begin{align}\begin{split}\label{e:wc_unit_slope_tangent}
&\delta^{-1/2}|(x-s)-(x_0-s_0)|\lesssim |w_c-w_0|\lesssim \|(x,s)-(x_0,s_0)\|^{1/2}_2\leq {\delta}^{1/2},\\
&\dist\bigl(w_c,[\min\{x_0-s_0,x-s\},\max\{x_0-s_0,x-s\}]\bigr)
\gtrsim |(x-s)-(x_0-s_0)|.
\end{split}\end{align}

Moreover, for each such critical point $w_c$ and every $w\in\fU$, \eqref{e:tangent_S} and \eqref{e:tangent_err} hold.

In this situation, we call $\fU$ a $\fc$-tangent chart (centered at $w_0$), and we say that
the point $(x,s)$ is \emph{adapted} to $\fU$.
\end{lemma}

As in the vertical tangency case, there are four choices of local descent and ascent paths, depending on the sign of
$S''(w_0;x_0,s_0)$ and the local geometry; see \Cref{f:unit_slope_tangent}. The four configurations in
\Cref{f:unit_slope_tangent} correspond to those in
\Cref{f:vertical_tangent1,f:vertical_tangent2,f:vertical_tangent3,f:vertical_tangent4}
after a $135^\circ$ counterclockwise rotation. By the symmetry among the three types of curvilinear triangles described in
\Cref{r:symmetry}, we may use the same local contour prescriptions as in the vertical tangency case. Namely, we use the contours from
\Cref{f:vertical_tangent1} for Panel~(A) of \Cref{f:unit_slope_tangent}, those from
\Cref{f:vertical_tangent2} for Panel~(B), and, similarly, those from
\Cref{f:vertical_tangent3,f:vertical_tangent4} for Panels~(C) and~(D), respectively.

The following lemma introduces the unit-slope cusp-turning chart and  record several of its properties. The statement can be proven in the same way as in the vertical cusp-turning case \Cref{c:cusp_turning_critical1}, so we omit.
\begin{lemma}[$\fc$-cusp-turning chart (unit-slope)]\label{c:cusp_turning_critical2}
Given a unit-slope cusp-turning location $(x_0, s_0)\in \fA$. Then for every $\fc>0$ sufficiently small  there exists
$\delta=\delta(\fc)>0$ sufficiently small such that the following holds.

Let $w_0=x_0-s_0\,\chi(x_0,s_0)=x_0-s_0$ be the critical point corresponding to $(x_0,s_0)$, and set
\[
\fU:=\{w\in\bC:\ |w-w_0|\le \fc\}.
\]
On $\fU$ the Riemann surface $\cC$ can be parametrized as $(f(w),w)$. Moreover, for all $w\in\fU$,
\begin{align}
 |\chi(w)|, |\chi'''(w)| \asymp 1,\quad  |1-\chi(w)|, |\chi''(w)|\asymp  |w-w_0|.
\end{align}

Now let $(x,s)\in\fL$ satisfy $\|(x,s)-(x_0,s_0)\|_2\le \delta$.
Then the tiling action $S(\,\cdot\,;x,s)$ has exactly three (formal) critical points $w_c$ inside $\fU$, satisfying 
\begin{align}\begin{split}\label{e:wc_unit_slope_cusp}
&\delta^{-2/3}|(x-s)-(x_0-s_0)|\lesssim |w_c-w_0|\lesssim \|(x,s)-(x_0,s_0)\|^{1/3}_2\leq {\delta}^{1/3},
\\
&\dist\bigl(w_c,[\min\{x_0-s_0,x-s\},\max\{x_0-s_0,x-s\}]\bigr)
\gtrsim |(x-s)-(x_0-s_0)|.
\end{split}\end{align}

Moreover, for each such critical point $w_c$ and every $w\in\fU$, \eqref{e:cusp_turning_S} and \eqref{e:cusp_turning_err} hold.

In this situation, we call $\fU$ a $\fc$-cusp-turning chart (centered at $w_0$), and we say that
the point $(x,s)$ is \emph{adapted} to $\fU$.
\end{lemma}

\begin{figure}	
\begin{subfigure}[t]{0.24\textwidth}
			\centering
			\begin{tikzpicture}
			\draw[] (0,0) arc (135:225:1);
			\draw (0,0) arc (-45:-135:1);
			\draw[dashed] (0,0)--(-1.5,-1.5);
			\draw[dashed] (-2.5,-1.5)--(-0.5,-1.5);

			\draw[white, fill=black]  (0,0) circle (0.05);
			%\draw[](0,0) node[right]{$(b,s_0)$};
			\draw[white, fill=black]  (-1.5,-1.5) circle (0.05);
			\draw[](-1.5,-1.5) node[below]{$w_0=x_0-s_0$};
			
			\node[circle, draw, fill=white, inner sep=1pt] at ({-0.55, -1}) {B};
			\node[circle, draw, fill=white, inner sep=1pt] at ({-1,-0.55}) {A};

			%\draw[white, fill=black]  ({1-2*sqrt(2)},{-sqrt(2)}) circle (0.05);
			%\draw[]({1-2*sqrt(2)},{-sqrt(2)}) node[below]{$w_c$};
			\end{tikzpicture}
			\caption{$S'''>0$.}
		\end{subfigure}				
	\begin{subfigure}[t]{0.24\textwidth}
			\centering
			\begin{tikzpicture}
			\draw[] (0,0) arc (135:45:1);
			\draw (0,0) arc (-45:45:1);
			\draw[dashed] (0,0)--(-0.5,-0.5);
			\draw[dashed] (0,0)--(1,1);
			\draw[dashed] (-1.5,-0.5)--(0.5,-0.5);

			\draw[white, fill=black]  (0,0) circle (0.05);
			%\draw[](0,0) node[right]{$(b,s_0)$};
			\draw[white, fill=black]  (-0.5,-0.5) circle (0.05);
			\draw[](-0.5,-0.5) node[below]{$w_0=x_0-s_0$};
			
			\node[circle, draw, fill=white, inner sep=1pt] at ({0.55, 1}) {B};
			\node[circle, draw, fill=white, inner sep=1pt] at ({1, 0.55}) {A};

			%\draw[white, fill=black]  ({1-2*sqrt(2)},{-sqrt(2)}) circle (0.05);
			%\draw[]({1-2*sqrt(2)},{-sqrt(2)}) node[below]{$w_c$};
			\end{tikzpicture}
			\caption{$S'''<0$.}
		\end{subfigure}

	\caption{
\label{f:unit_slope_cusp_turning}
Local path associated with unit-slope cusp-turning locations}
	\end{figure}
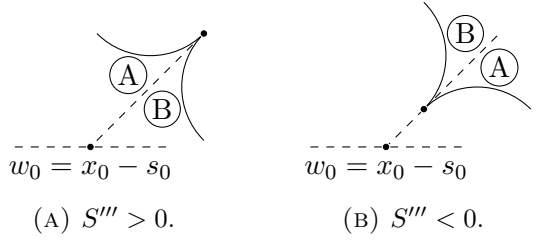

As in the vertical cusp-turning case, there are two choices of local descent and ascent paths, determined by the sign of
$S'''(w_0;x_0,s_0)$; see \Cref{f:unit_slope_cusp_turning}. The two configurations in
\Cref{f:unit_slope_cusp_turning} correspond to those in
\Cref{f:c_vertical_cusp1,f:c_vertical_cusp2} after a $135^\circ$ counterclockwise rotation. By the symmetry among the three types of curvilinear triangles described in
\Cref{r:symmetry}, we may adopt the same local contour prescriptions as in the vertical cusp-turning case. Namely, we use the contours from
\Cref{f:c_vertical_cusp1} for Panel~(A) of \Cref{f:unit_slope_cusp_turning} and those from
\Cref{f:c_vertical_cusp2} for Panel~(B).

The following lemma shows that the local descent and ascent paths can be deformed into steepest-descent and steepest-ascent paths with negligible error. The statements can be proven in the same way as in the vertical tangent/cusp-turning case \Cref{l:vertical_tangent_steepest} and \Cref{c:vertical_cusp_steepest}, so we omit. 
\begin{lemma}\label{l:unit_slope_tangent_steepest}
Under the assumptions and notation of \Cref{c:tangent_critical2},
the statement of \Cref{l:vertical_tangent_steepest} remains valid after
replacing the vertical-tangency setting with the unit-slope tangency setting.

Likewise, under the assumptions and notation of
\Cref{c:cusp_turning_critical2}, the statement of
\Cref{c:vertical_cusp_steepest} remains valid after replacing the vertical
cusp-turning setting with the unit-slope cusp-turning setting.
\end{lemma}

\subsection{Unit-slope tangent frozen chart}\label{s:unit_slope_frozen_neighborhood}
Suppose that $(x_0,s_0)\in \fP\setminus \fL$ lies on a unit-slope extended side tangent to the arctic boundary at the unit-slope tangency point $(x'_0,s_0')\in \fA$. From \eqref{e:bcond}  we have $f(x'_0, s'_0)=\infty$,  $\chi(x'_0;s'_0)=1$ and 
\begin{align}
w_0:=x'_0-s'_0 \chi(w_0)=x_0-s_0 \chi(w_0)=x_0-s_0\in \cC(\bR),\quad \chi(w_0)=1.
\end{align}
In a neighborhood of $w_0$, we denote $1-\chi(w)=(w-w_0)\chi_0(w)$, and the same as in  \eqref{e:expandchi_a}, we have
\begin{align}\begin{split}\label{e:expandchi_a2}
& 1-\chi(w)=-\frac{(w-w_0)}{s'_0}+a_1(w-w_0)^2+a_2(w-w_0)^3+\cdots,\\
&\chi_0(w)=\frac{1}{s'_0}+a_1(w-w_0)+a_2(w-w_0)^2+\cdot,
\end{split}\end{align}
where $c_1=-\chi''(w_0)/2$ and $c_2=-\chi'''(w_0)/6$. It follows
\begin{align}
S'(w;x_0, s_0)
&=\ln \frac{w-x_0}{x_0-s_0-w}-\ln \frac{\chi(w)}{1-\chi(w)}
=\ln \frac{(1-\chi(w))(w-w_0-s_0)}{\chi(w)(w_0-w)}
=\ln \frac{\chi_0(w)(s_0-(w-w_0))}{(1-(w-w_0)\chi_0(w))}\\
&=\ln \frac{s_0}{s_0'}+\ln \frac{(1+s'_0a_1(w-w_0)+\OO((w-w_0)^2))(1-(w-w_0)/s_0)}{1-(w-w_0)/s'_0 +\OO((w-w_0)^2)}\\
&=\ln \frac{s_0}{s_0'}+\left(\frac{1}{s_0'}-\frac{1}{s_0}+s_0 a_1\right)(w-w_0)+\OO((w-w_0)^2).
\end{align}

The following lemma introduces the unit-slope frozen chart and  record several of its properties. The statement can be proven in the same way as in the vertical case \Cref{c:vertical_frozen_critical1}, so we omit.
\begin{lemma}[$\fc$-frozen chart (unit-slope tangent)]\label{c:unit_slope_frozen_critical1}
Fix $\delta_0>0$ and let $(x_0,s_0)\in \fP$.  Assume there exists a unit-slope tangent location $(x_0',s_0')\in \fA$ such that
the line through $(x_0,s_0)$ and $(x_0',s_0')$ is tangent to the arctic curve at $(x_0',s_0')$, and moreover $(x_0, s_0)$ is at least distance $\delta_0$ from $(x_0',s_0')$. Then for every sufficiently small $\fc>0$ there exists
$\delta=\delta(\fc)>0$, also sufficiently small and depending on $\delta_0$,  such that the following holds.

Let $w_0=x_0-s_0\,\chi(x_0,s_0)=x_0-s_0$ be the critical point corresponding to $(x_0,s_0)$, and set
\[
\fU:=\{w\in\bC:\ |w-w_0|\le \fc\}.
\]
On $\fU$ the Riemann surface $\cC$ can be parametrized as $(f(w),w)$. We will therefore identify
$\fU$ with its image in $\cC$ under the map $w\mapsto (f(w), w)$.  Moreover, for all $w\in\fU$,
\begin{align}
  |\chi(w)| \asymp 1,\quad  |1-\chi(w)|\asymp |w-w_0|.
\end{align}

Now let $(x,s)\in\fL$ satisfy $\|(x,s)-(x_0,s_0)\|_2\le \delta$.
Then the tiling action $S(\,\cdot\,;x,s)$ has exactly one (formal) critical point $w_c$ inside $\fU$, satisfying 
\begin{align}\label{e:wc_unit_slope_frozen}
\dist(w_c, [\min\{x_0-s_0, x-s\}, \max\{x_0-s_0, x-s\}])|\asymp |(x-s)-(x_0-s_0)|,
\end{align}
Moreover, for every $w\in\fU$, \eqref{e:tangent_frozen_S} and \eqref{e:tangent_frozen_err} hold.

In this situation, we call $\fU$ a (unit-slope tangent) $\fc$-frozen chart (centered at $w_0$), and we say that
the point $(x,s)$ is \emph{adapted} to $\fU$.
\end{lemma}

As in the vertical tangency case, the unit-slope tangency point
$(x'_0,s'_0)\in\fA$ lies on the common boundary of two adjacent curvilinear
triangles $\fT_A$ and $\fT_B$. There are two cases for the integration
contours.

Suppose first that $(x_0,s_0)$ does not lie on the common boundary of
$\fT_A$ and $\fT_B$, and let $\fT\in\{\fT_A,\fT_B\}$ denote the region
containing $(x_0,s_0)$. Depending on whether
\[
(x_0,s_0)\in\ell_-(w_0;\fT)
\qquad\text{or}\qquad
(x_0,s_0)\in\ell_+(w_0;\fT),
\]
we associate with $(x_0,s_0)$ the circular component of the local descent
or ascent contour, respectively, associated with $(x'_0,s'_0)$ for $\fT$
in \Cref{s:unit_slope_tangent}.

If $(x_0,s_0)$ lies on the common boundary of the two adjacent curvilinear
triangles, we associate with $(x_0,s_0)$ both the circular component of the
local descent contour (blue) and that of the local ascent contour (red)
associated with $(x'_0,s'_0)$ for either $\fT_A$ or $\fT_B$ in
\Cref{s:unit_slope_tangent}.

The following lemma shows that the local descent and ascent paths can be deformed into steepest-descent and steepest-ascent paths with negligible error. The statement can be proven in the same way as in the vertical case \Cref{l:local_descent_deformation}, so we omit.
\begin{lemma}\label{l:local_descent_deformation2}
Adopt the assumptions and notation of
\Cref{c:unit_slope_frozen_critical1}. The statement of
\Cref{l:local_descent_deformation} remains valid after replacing the
vertical-tangency setting with the unit-slope tangency setting and replacing
the definition of $r_n$ by the following:
\begin{enumerate}

\item \emph{Close to the tangent line.} If
$|(x-s)-(x_0-s_0)|\leq (\ln n)^3/n$, we define
$r_n=(\ln n)^5/n$.

\item \emph{Other frozen regime.} If
$|(x-s)-(x_0-s_0)|>(\ln n)^3/n$, we define
$r_n:=(\ln n)^2 (|(x-s)-(x_0-s_0)|/n)^{1/2}$.

\end{enumerate}
\end{lemma}

\subsection{Horizontal tangent chart and cusp-turning chart}\label{s:horizontal_tangent}

Fix any horizontal tangent location $(x_0, s_0)\in \fA$. Then from \eqref{e:bcond} and
\eqref{e:arctic_chider}, we have $f(x_0, s_0)=-1$, $\chi(x_0;s_0)=\infty$ and by the third statement in \Cref{p:surface}
    \begin{align}\label{e:chi_exp_infinite}
        \chi(w) = \frac{-w+x_0}{s_0} + \sum_{i\geq 1}\frac{c_i}{(x_0-w)^i}.
    \end{align}
The critical points of $S'(w;x_0, s_0)$ is at $w=\infty$. We can change the coordinate, set
\begin{align}\label{e:change_coordinate}
\wt w=\frac{1}{x_0-w}, \quad \wt S({\wt w}; x_0, s_0)=S(w; x_0, s_0),\quad \chi(w)=\wt \chi({\wt w}) = \frac{1}{s_0 {\wt w}} + \sum_{i\geq 1}c_i{\wt w}^i,
\end{align}
Then
\begin{align}
\wt S'({\wt w}; x_0, s_0)\del_w {\wt w}=\wt S'({\wt w}; x_0, s_0){\wt w}^2=S'(w; x_0, s_0)
\end{align}
and we conclude that
\begin{align}\begin{split}\label{e:hor_tangent}
\wt S'({\wt w}; x_0, s_0)
&=\frac{S'(w; x_0, s_0)}{{\wt w}^2}
=\frac{1}{{\wt w}^2}\left(\ln \frac{w-x_0}{x_0-s_0-w}-\ln \frac{\chi(w)}{1-\chi(w)}\right)\\
&=\frac{1}{{\wt w}^2}\left(\ln \frac{1}{s_0{\wt w}-1}-\ln \frac{1+\sum_{i\geq 1} s_0 c_i {\wt w}^{i+1}}{s_0{\wt w} -1-\sum_{i\geq 1} c_i s_0 {\wt w}^{i+1}}\right)\\
&=\frac{1}{{\wt w}^2}\left(\ln \frac{1-s_0{\wt w} +\sum_{i\geq 1} c_i s_0 {\wt w}^{i+1}}{(1-s_0{\wt w})(1+\sum_{i\geq 1} s_0 c_i {\wt w}^{i+1})}\right)=s_0^2 c_1 {\wt w} + (c_1 s_0^3+c_2 s_0^2){\wt w}^2+\OO({\wt w}^3).
\end{split}\end{align}

In particular, \(\wt S'({\wt w};x_0,s_0)\) is analytic in a neighborhood of \(0\). If $(x_0, s_0)$ is not a cusp location, then 
\begin{align}\label{e:hor_c_1} 
c_1=(1/2)\del_{\wt w}^2({\wt w}\wt\chi({\wt w}))|_{{\wt w}=0},\quad |c_1|\asymp 1.
\end{align}
 If $(x_0, s_0)$ is a cusp-turning location, then 
 \begin{align}\label{e:hor_c_2}
 \del_{\wt w}^2({\wt w}\wt\chi({\wt w}))|_{{\wt w}=0}=0,\quad 
 c_2=(1/6)\del_{\wt w}^3({\wt w}\wt\chi({\wt w}))|_{{\wt w}=0},\quad  |c_2|\asymp 1.
 \end{align}
  In these cases, the signs of  \(\wt S''(0;x_0,s_0)\) and \(\wt S'''(0;x_0,s_0)\) are determined by the signs of $\del_{\wt w}^2({\wt w}\wt\chi({\wt w}))|_{{\wt w}=0}$ and $\del_{\wt w}^3({\wt w}\wt\chi({\wt w}))|_{{\wt w}=0}$, as classified in \Cref{f:tangent} and \Cref{f:cusp_turning}.

For \((x,s)\) in a sufficiently small neighborhood of \((x_0,s_0)\), as
\(\wt w\to0\), we have
\begin{align}\label{e:hor_tangent_diff}
\wt S'(\wt w;x,s)-\wt S'(\wt w;x_0,s_0)
&=
\frac{1}{\wt w^2}
\left(
\ln\frac{(x_0-x)\wt w-1}{(x-x_0-s)\wt w+1}
-\ln\frac{1}{s_0\wt w-1}
\right)=
\frac{s-s_0}{\wt w}+\OO(1).
\end{align}
Since \(\wt S'(\wt w;x_0,s_0)\) is bounded as \(\wt w\to0\), it follows
that
\[
\wt S'(\wt w;x,s)
=
\frac{s-s_0}{\wt w}+\OO(1).
\]
Thus, if \(s>s_0\), then \(\wt w=0\) is a sink for the steepest-descent
flow of \(\Re\wt S(\cdot;x,s)\); that is, the steepest-descent paths point
toward \(0\). If \(s<s_0\), then \(\wt w=0\) is a source; that is, the
steepest-descent paths point away from \(0\).

%Integrating \eqref{e:hor_tangent_diff} along a small path avoiding ${\wt w}=0$ yields
%\[
%\wt S({\wt w};x,s)-\wt S({\wt w};x_0,s_0)=(s-s_0)\ln {\wt w}+\OO({\wt w}),
%\quad {\wt w}\to 0,
%\]
%and therefore
%\[
%e^{n\wt S({\wt w};x,s)} = {\wt w}^{\,n(s-s_0)}\,e^{\OO(1)},\quad {\wt w}\to 0.
%\]
%In particular, if $s<s_0$ then $e^{n\wt S({\wt w};x,s)}$ has a pole at ${\wt w}=0$.

The following lemmas introduce the horizontal tangent chart and horizontal cusp-turning chart, and  record several of their properties. Their proofs, based on a Taylor expansion of the tiling action, are deferred to \Cref{s:horizontal_tangent_chart_proof}.
\begin{lemma}[$\fc$-tangent chart (horizontal)]\label{c:tangent_critical3}
Let \((x_0,s_0)\in\fA\) be a horizontal tangent location which is not a cusp location.
Then, for every sufficiently small \(\fc>0\), there exists a sufficiently small
\(\delta=\delta(\fc)>0\) such that the following holds.

Let
\[
\wt \fU:=\{\wt w\in\bC:\ |\wt w|\le \fc\}.
\]
On \(\wt \fU\), the Riemann surface \(\cC\) can be parametrized as
\[
\wt w\longmapsto \left(\frac{\wt\chi(\wt w)}{1-\wt\chi(\wt w)},\wt w\right).
\]
We identify \(\wt \fU\) with its image in \(\cC\) under this map. Moreover, for all
\(\wt w\in\wt \fU\),
\begin{align}
|\wt w \wt \chi(\wt w)|,  \left|\del_{\wt w}^2\bigl(\wt w\,\wt\chi(\wt w)\bigr)\right|\asymp 1.
\end{align}

Now let \((x,s)\in\fL\) satisfy
\[
\|(x,s)-(x_0,s_0)\|_2\le \delta .
\]
Then the tiling action \(\wt S(\,\cdot\,;x,s)\) has exactly two (formal) critical points
\(\wt w_c\) inside \(\wt \fU\), satisfying
\begin{align}\label{e:wc_horizontal_tangent}
\delta^{-1/2}|s-s_0|
\lesssim
|\wt w_c|
\lesssim
\|(x,s)-(x_0,s_0)\|_2
\le \delta^{1/2},
\quad |\wt w_c|
\gtrsim |s-s_0|
\end{align}
Moreover, for each such critical point \(\wt w_c\) and every \(\wt w\in\fU\),
\begin{align}
\wt S(\wt w;x,s)-\wt S(\wt w_c;x,s)
=
d\,\wt w^2+\cE(\wt w),
\qquad
d:=\frac{\wt S''(0;x_0,s_0)}{2},
\end{align}
where the error term satisfies
\begin{align}
|\cE(\wt w)|
\leq C\bigl(\delta\log(1/\delta)+|\wt w|^3\bigr)
\le
\frac{|d|\fc^2}{100}.
\end{align}

In this situation, we call $\wt \fU$ a $\fc$-tangent chart. Moreover, for each
critical point $\wt w_c$ of $\wt S(\,\cdot\,;x,s)$ in $\wt\fU$, we say that the point
$(x,s)$ is \emph{adapted} to $\wt \fU$.

We also introduce the corresponding chart in the $w$-plane, together with the
corresponding pulled-back critical points, by pulling back the $\wt w$-chart:
\[
\fU=\{w:\wt w\in \wt\fU\}
    =\left\{x_0-\frac{1}{\wt w}:\wt w\in \wt\fU\right\},
    \qquad
w_c=x_0-\frac{1}{\wt w_c}.
\]
We call $\fU$ a $\fc$-tangent chart as well, and say that the point
$(x,s)$ is \emph{adapted} to $\fU$.

\end{lemma}

\begin{figure}
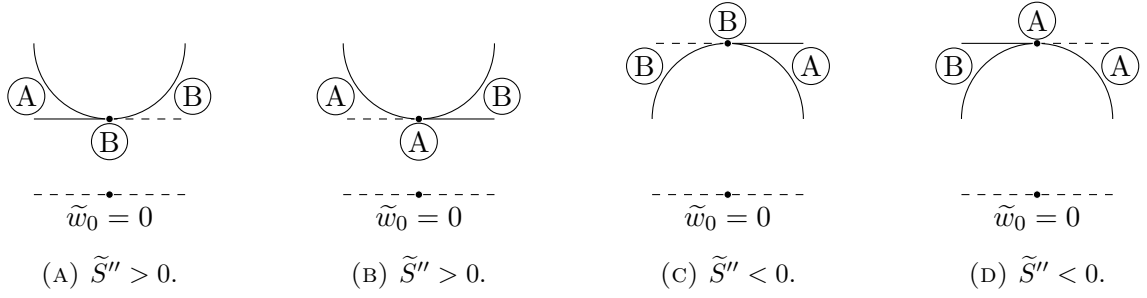
	
\begin{subfigure}[t]{0.24\textwidth}
			\centering

			\caption{$\wt S''>0$.}
\end{subfigure}	
\begin{subfigure}[t]{0.24\textwidth}
			\centering
			%
			\caption{$\wt S''>0$.}
		\end{subfigure}			
\begin{subfigure}[t]{0.24\textwidth}
			\centering
			%
			\caption{$\wt S''<0$.}
		\end{subfigure}	
\begin{subfigure}[t]{0.24\textwidth}
			\centering
			%
			\caption{$\wt S''<0$.}
		\end{subfigure}	
	\caption{
\label{f:horizontal_tangent}
Local path associated with horizontal tangent locations}
	\end{figure}

As in the vertical tangency case, there are four choices of local descent and ascent paths, depending on the sign of
$\wt S''(0;x_0,s_0)$ and the local geometry; see \Cref{f:horizontal_tangent}. The four configurations in
\Cref{f:horizontal_tangent} correspond to those in
\Cref{f:vertical_tangent1,f:vertical_tangent2,f:vertical_tangent3,f:vertical_tangent4} after a $90^\circ$
clockwise rotation. By the symmetry among the three types of curvilinear triangles described in
\Cref{r:symmetry}, we may use the same local contour prescriptions as in the vertical tangency case. Namely, we use the contours from \Cref{f:vertical_tangent1} for
panel (A) of \Cref{f:horizontal_tangent}, those from \Cref{f:vertical_tangent2} for panel (B), and similarly for
panels (C) and (D). We denote the corresponding local descent and ascent contours as $\wt \sfC^{\rm d}(0)$ and $\wt \sfC^{\rm a}(0)$.

In what follows, we orient the local paths $\wt \sfC^{\rm d}(0), \wt \sfC^{\rm a}(0)$  associated with
\Cref{f:horizontal_tangent} so that they are compatible with the contours in the
liquid region introduced in \Cref{s:critical_bulk}. The resulting orientations
are exactly the same as those in
\Cref{f:vertical_tangent1,f:vertical_tangent2,f:vertical_tangent3,f:vertical_tangent4}.

For $(x,s)\in \fL$ close to $(x_0, s_0)$, let $w_c=x-s\chi({w}_c)$, and ${\wt w}_c=1/(x_0-w_c)$. We recall the factorization for descent vector from \eqref{eq:factor} (composing with the map $w\mapsto 1/(x_0-w)$ it becomes)
\begin{align}\label{eq:forward_factor}
\wt v^{\rm d}(w_c)={\wt w}_c^2\sqrt{-\frac{1}{S''(w_c;x,s)}}=
{\wt w}_c^2 \frac{\sqrt{\chi(w_c)}\sqrt{1-\chi(w_c)}}{\sqrt{1/s+\chi'(w_c)}}
={\wt w}_c^2\frac{\sqrt{\wt \chi({\wt w}_c)}\sqrt{1-\wt\chi({\wt w}_c)}}{\sqrt{1/s+{\wt w}_c^2\wt \chi'({\wt w}_c)}},
\end{align}

We next give a geometric description of this choice as $(x,s)$ approaches
$(x_0,s_0)$ and $w_c$ approaches $w_0$.  
  As $(x,s)$ approaches to $(x_0, s_0)$, ${\wt w}_c$ approaches $0$. 
 Let ${\wt w}_c=a+\ri b$ with $b>0$, then \eqref{e:curve_ht_reg} and \eqref{e:change_coordinate} give
\begin{align}
&1/s+{\wt w}_c^2 \wt \chi'({\wt w}_c)=(\del_{\wt w} \wt \chi({\wt w}))|_{{\wt w}=0}\ri b(a+b\ri)+\OO(b(|a|+b)^2),\quad \wt \chi({\wt w}_c)=1/(s_0{\wt w}_c)+\OO({\wt w}_c).
\end{align}
Thus as $(x,s)\rightarrow (x_0, s_0)$, the direction of $v^{\rm d}(w_c)={\wt w}_c^2\sqrt{-1/S''(w_c;x,s))}$ is given by
\begin{align}\label{e:hor_tangent_direction}
(1+\oo(1))(a+\ri b)^2\sqrt{-\frac{1}{s_0^2 (a+\ri b)^2((\del_{\wt w} \wt \chi({\wt w}))|_{{\wt w}=0} \ri b(a+b\ri))} }=\frac{1+\oo(1)}{\sqrt{2b}}\sqrt{\frac{a\ri- b}{2b \wt S''(0;x_0, s_0)}}
\end{align}
where we used \eqref{e:hor_tangent}.

For $\wt S''(0;x_0,s_0)>0$, then the direction of the  square root in~\eqref{e:hor_tangent_direction} lies in
$
\{\pm e^{\ri\theta}:\pi/4<\theta<3\pi/4\}.
$
If it lies in $\{e^{\ri\theta}:\pi/4<\theta<3\pi/4\}$, the corresponding direction is illustrated by the blue paths in \Cref{f:vertical_tangent1} and \Cref{f:vertical_tangent2};
if  it lies in $\{-e^{\ri\theta}:\pi/4<\theta<3\pi/4\}$, the corresponding direction should be reversed.

For $\wt S''(0;x_0,s_0)<0$, then the direction of the  square root in~\eqref{e:hor_tangent_direction} lies in
$
\{\pm e^{\ri\theta}:-\pi/4<\theta<\pi/4\}.
$
If it lies in $\{-e^{\ri\theta}:-\pi/4<\theta<\pi/4\}$, the corresponding direction is illustrated by the blue paths in \Cref{f:vertical_tangent3} and \Cref{f:vertical_tangent4};
if  it lies in $\{e^{\ri\theta}:-\pi/4<\theta<\pi/4\}$, the corresponding direction should be reversed.

Finally, we orient the ascent contour $\wt\sfC^{\rm a}(0)$ as indicated by the red paths in
\Cref{f:vertical_tangent1,f:vertical_tangent2,f:vertical_tangent3,f:vertical_tangent4}.

If $(x_0, s_0)\in \fA$ is also a cusp location, then $c_1=0$ and $c_2\neq 0$, and 
\begin{align}\label{e:tS'''}
\wt S'({\wt w};x_0, s_0)=c_2s^2_0 {\wt w}^2+\OO({\wt w}^3).
\end{align}

\begin{lemma}[$\fc$-cusp-turning chart (horizontal)]\label{c:cusp_turning_critical3}
Given a horizontal cusp-turning location $(x_0, s_0)\in \fA$. Then for every $\fc>0$ sufficiently small  there exists
$\delta=\delta(\fc)>0$ sufficiently small such that the following holds.

Let
\[
\wt \fU:=\{{\wt w}\in\bC:\ |{\wt w}|\le \fc\}.
\]
On $\wt \fU$ the Riemann surface $\cC$ can be parametrized as $(\wt \chi({\wt w})/(1-\wt \chi({\wt w})),{\wt w})$. We will therefore identify
$\wt \fU$ with its image in $\cC$ under the map ${\wt w}\mapsto (\wt \chi({\wt w})/(1-\wt \chi({\wt w})),{\wt w})$.   Moreover, for all ${\wt w}\in\wt \fU$,
\begin{align}
|\wt \chi(\wt w)|,  |\del_{\wt w}^3 ({\wt w}\wt\chi({\wt w}))|\asymp 1,\quad |\del_{\wt w}^2 ({\wt w}\wt\chi({\wt w}))|\asymp |{\wt w}|
\end{align}

Now let $(x,s)\in\fL$ satisfy $\|(x,s)-(x_0,s_0)\|_2\le \delta$.
Then the tiling action \(\wt S(\,\cdot\,;x,s)\) has exactly three (formal) critical points
\(\wt w_c\) inside \(\wt \fU\), satisfying
\begin{align}\label{e:wc_horizontal_cusp}
\delta^{-2/3}|s-s_0|\lesssim |{\wt w}_c|\lesssim \|(x,s)-(x_0,s_0)\|^{1/3}_2\leq {\delta}^{1/3},\quad |\wt w_c|\gtrsim |s-s_0|
\end{align}

Moreover, for each such critical point ${\wt w}_c$ and every ${\wt w}\in\wt \fU$, 
\begin{align}
\wt S({\wt w};x,s)-\wt S({\wt w}_c;x,s)=d{\wt w}^3 +\cE({\wt w}),\quad d:=\wt S'''(0;x_0,s_0)/6, 
\end{align}
The error term satisfies
\begin{align}\label{e:hor_cusp_turning_err}
|\cE({\wt w})|=\OO(\delta\ln(1/\delta)+|{\wt w}|^4)\leq\frac{|d|\fc^3}{100}.
\end{align}

In this situation, we call $\wt \fU$ a $\fc$-cusp-turning chart (centered at ${\wt w}_0$), and we say that
the point $(x,s)$ is \emph{adapted} to $\wt \fU$.

We also introduce the corresponding chart in the $w$-plane, together with the
corresponding pulled-back critical points, by pulling back the $\wt w$-chart:
\[
\fU=\{w:\wt w\in \wt\fU\}
    =\left\{x_0-\frac{1}{\wt w}:\wt w\in \wt\fU\right\},
    \qquad
w_c=x_0-\frac{1}{\wt w_c}.
\]
We call $\fU$ a $\fc$-tangent chart as well, and say that the point
$(x,s)$ is \emph{adapted} to $\fU$.

\end{lemma}

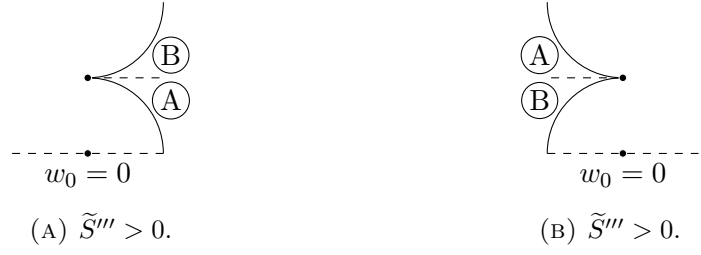
\begin{figure}	
\begin{subfigure}[t]{0.4\textwidth}
			\centering
			\begin{tikzpicture}
			\draw[] (0,0) arc (-90:0:1);
			\draw[] (0,0) arc (90:0:1);
			\draw[dashed] (0,0)--(1,0);
			\draw[dashed] (-1,-1)--(1,-1);

			\draw[white, fill=black]  (0,0) circle (0.05);
			%\draw[](0,0) node[right]{$(b,s_0)$};
			\draw[white, fill=black]  (0,-1) circle (0.05);
			\draw[](0,-1) node[below]{$w_0=0$};
			
			\node[circle, draw, fill=white, inner sep=1pt] at ({1.1,-0.3}) {A};
			\node[circle, draw, fill=white, inner sep=1pt] at ({1.1,0.3}) {B};

			%\draw[white, fill=black]  ({1-2*sqrt(2)},{-sqrt(2)}) circle (0.05);
			%\draw[]({1-2*sqrt(2)},{-sqrt(2)}) node[below]{$w_c$};
			\end{tikzpicture}
			\caption{$\wt S'''>0$.}
\end{subfigure}	
\begin{subfigure}[t]{0.4\textwidth}
			\centering
			\begin{tikzpicture}
			\draw[] (0,0) arc (90:180:1);
			\draw[] (0,0) arc (-90:-180:1);
			\draw[dashed] (0,0)--(-1,0);
			\draw[dashed] (-1,-1)--(1,-1);

			\draw[white, fill=black]  (0,0) circle (0.05);
			%\draw[](0,0) node[right]{$(b,s_0)$};
			\draw[white, fill=black]  (0,-1) circle (0.05);
			\draw[](0,-1) node[below]{$w_0=0$};
			
			\node[circle, draw, fill=white, inner sep=1pt] at ({-1.1,-0.3}) {B};
			\node[circle, draw, fill=white, inner sep=1pt] at ({-1.1,0.3}) {A};

			%\draw[white, fill=black]  ({1-2*sqrt(2)},{-sqrt(2)}) circle (0.05);
			%\draw[]({1-2*sqrt(2)},{-sqrt(2)}) node[below]{$w_c$};
			\end{tikzpicture}
			\caption{$\wt S'''>0$.}
		\end{subfigure}			
	\caption{
\label{f:horizontal_cusp_turning}
Local path associated with horizontal cusp-turning locations}
	\end{figure}

As in the vertical cusp-turning case, there are two choices of local steepest--descent paths, determined by the sign of
$S'''({\wt w}_0;x_0,s_0)$; see \Cref{f:horizontal_cusp_turning}. The two configurations in
\Cref{f:horizontal_cusp_turning} correspond to those in
\Cref{f:c_vertical_cusp1,f:c_vertical_cusp2} after a $90^\circ$ clockwise rotation. Therefore, by the same
argument as before, we may adopt the same local contour prescriptions as in the vertical cusp-turning case: we use the
contours from \Cref{f:c_vertical_cusp1} for panel (A) of \Cref{f:horizontal_cusp_turning}, and those from
\Cref{f:c_vertical_cusp2} for panel (B). We denote the corresponding local descent and ascent contours as $\wt \sfC^{\rm d}(0)$ and $\wt \sfC^{\rm a}(0)$.

In the following we orient the local paths $\wt \sfC^{\rm d}(0), \wt \sfC^{\rm a}(0)$ associated with  \Cref{f:horizontal_cusp_turning}, so they are compatible with the contours in the liquid  region as introduced in \Cref{s:critical_bulk}. Again, the resulting orientations
are exactly the same as those in
\Cref{f:c_vertical_cusp1,f:c_vertical_cusp2}.

For $(x,s)\in \fL$ close to $(x_0, s_0)$, let $w_c=x+s\chi({ w}_c)$, and ${\wt w}_c=1/(x_0-w_c)$. We recall the factorization for descent vector from \eqref{eq:forward_factor}.
We next give a geometric description of this choice as $(x,s)$ approaches
$(x_0,s_0)$ and $w_c$ approaches $w_0$.  
 Let ${\wt w}_c=a+\ri b$ with $b>0$, then \eqref{e:curve_ht_reg} and \eqref{e:change_coordinate}  gives
\begin{align}
&1/s+{\wt w}_c^2 \wt \chi'({\wt w}_c)
=\del_{\wt w}^2 ({\wt w}\widetilde\chi({\wt w}))|_{{\wt w}=0} b\,(a+\ri b)\left(a\ri - \frac{b}{3}\right)+\OO(b(|a|+b)^3), \\
& \wt \chi({\wt w}_c)=1/(s_0{\wt w}_c)+\OO({\wt w}).
\end{align}
Thus as $(x,s)\rightarrow (x_0, s_0)$, the direction of $\wt v^{\rm d}(w_c)={\wt w}_c^2\sqrt{-1/S''(w_c;x,s))}$ is given by
\begin{align}\label{e:hor_cusp_direction}
(1+\oo(1))(a+\ri b)^2\sqrt{-\frac{1}{s_0^2 (a+\ri b)^3((\del^2_{\wt w} \wt \chi({\wt w}))|_{{\wt w}=0} b(a\ri-b/3))} }=\frac{1+\oo(1)}{\sqrt{b(b^2+9a^2)}}\sqrt{\frac{((3a^2+b^2)\ri +2ab)}{\wt S'''(0;x_0, s_0)}}
\end{align}
where we used \eqref{e:tS'''}.

For $\wt S'''(0;x_0,s_0)>0$, then the direction of the square root in~\eqref{e:hor_cusp_direction} lies in
$
\{\pm e^{\ri\theta}:0<\theta<\pi/2\}.
$
If it lies in $\{e^{\ri\theta}:0<\theta<\pi/2\}$, the corresponding direction is illustrated by the blue paths in \Cref{f:c_vertical_cusp1};
if  it lies in $\{-e^{\ri\theta}:0<\theta<\pi/2\}$, the corresponding direction should be reversed.

For $\wt S'''(0;x_0,s_0)<0$, then the direction of the square root in~\eqref{e:hor_cusp_direction} lies in
$
\{\pm e^{\ri\theta}:\pi/2<\theta<\pi\}.
$
If it lies in $\{e^{\ri\theta}:\pi/2<\theta<\pi\}$, the corresponding direction is illustrated in \Cref{f:c_vertical_cusp2};
if  it lies in $\{-e^{\ri\theta}:\pi/2<\theta<\pi\}$, the corresponding direction should be reversed.

Finally, we also orient the ascent contour $\wt \sfC^{\rm a}(0)$ as indicated by the red paths in \Cref{f:c_vertical_cusp1} and \Cref{f:c_vertical_cusp2}.

In both the horizontal tangent and cusp-turning cases, we define the corresponding local descent and ascent contours in the \(w\)-plane by pulling back the contours in the \(\wt w\)-plane:
\begin{align}
\sfC^{\rm d}(\infty)
=\{w: \wt w\in \wt \sfC^{\rm d}(0)\}, \quad 
\sfC^{\rm a}(\infty)
=\{w: \wt w\in \wt \sfC^{\rm a}(0)\},\quad \wt w=\frac{1}{x_0-w}.
\end{align}

The following lemma shows that the local descent and ascent paths can be deformed into steepest-descent and steepest-ascent paths with negligible error. Its proof is deferred to \Cref{s:horizontal_tangent_chart_proof}.
\begin{lemma}\label{c:horizontal_tangent_steepest}
Under the assumptions and notation of \Cref{c:tangent_critical3}, the
statement of \Cref{l:vertical_tangent_steepest} remains valid after replacing
the vertical-tangency setting with the horizontal-tangency setting, setting
$w_0=0$, and replacing $S,\sfC,\sfD$ with
$\wt S,\wt\sfC,\wt\sfD$, respectively.

Likewise, under the assumptions and notation of
\Cref{c:cusp_turning_critical3}, the statement of
\Cref{c:vertical_cusp_steepest} remains valid under the same replacements.
\end{lemma}

\subsection{Horizontal tangent frozen chart}\label{s:horizontal_frozen_neighborhood}
Suppose that $(x_0,s_0)\in \fP\setminus \fL$ lies on a horizontal extended side tangent to the arctic boundary at the horizontal tangency point $(x'_0,s_0)\in \fA$. From \eqref{e:bcond} we have $f(x_0', s_0)=-1$ and $\chi(x_0',s_0)=\infty$. The same as in \eqref{e:chi_exp_infinite}, we have the expansion    \begin{align}\label{e:chi_exp_infinite2}
        \chi(w) = \frac{-w+x'_0}{s_0} + \sum_{i\geq 1}\frac{c_i}{(x'_0-w)^i}.
    \end{align}
We can change the coordinate, set
\begin{align}\label{e:frozen_change_coordinate}
\wt w=\frac{1}{x'_0-w}, \quad \wt S(\wt w; x_0, s_0)=S(w; x_0, s_0),\quad \wt \chi(\wt w)=\chi(w) = \frac{1}{s_0 \wt w} + \sum_{i\geq 1}c_i\wt w^i
\end{align}
Then
\begin{align}
\wt S'({\wt w}; x_0, s_0)\del_w {\wt w}=\wt S'({\wt w}; x_0, s_0){\wt w}^2=S'(w; x_0, s_0)
\end{align}
and we conclude that
\begin{align}\begin{split}\label{e:h_S_extend}
\wt S'({\wt w}; x_0, s_0)
&=\frac{S'(w; x_0, s_0)}{{\wt w}^2}
=\frac{1}{{\wt w}^2}\left(\ln \frac{w-x_0}{x_0-s_0-w}-\ln \frac{\chi(w)}{1-\chi(w)}\right)\\
&=\frac{1}{{\wt w}^2}\left(\ln \frac{1+(x_0-x_0'){\wt w}}{(s_0-x_0+x_0'){\wt w}-1}-\ln \frac{1+\sum_{i\geq 1} s_0 c_i {\wt w}^{i+1}}{s_0{\wt w} -1-\sum_{i\geq 1} c_i s_0 {\wt w}^{i+1}}\right)\\
&=\frac{1}{{\wt w}^2}\left(\ln \frac{(1-s_0 {\wt w})(1+(x_0-x_0'){\wt w})}{1-(s_0-x_0+x_0'){\wt w}}-\ln \frac{(1-s_0 {\wt w})\left(1+\sum_{i\geq 1} s_0 c_i {\wt w}^{i+1}\right)}{1-s_0{\wt w} +\sum_{i\geq 1} c_i s_0 {\wt w}^{i+1}}\right)\\
&=\frac{1}{{\wt w}^2}\left(-\frac{s_0(x_0-x_0'){\wt w}^2}{1-(s_0-x_0+x_0'){\wt w}}+s_0^2 c_1 {\wt w}^3+\OO({\wt w}^4)\right)\\
&=s_0(x_0'-x_0)+(s_0^2 c_1-s_0(x_0-x_0')(x_0-x_0'-s_0)) {\wt w} +\OO({\wt w}^2).
\end{split}\end{align}

The following lemma introduces the horizontal frozen chart and  record several of its properties. The proof is based on a Taylor expansion of the tiling action, and is deferred to \Cref{s:horizontal_frozen_chart_proof}.
\begin{lemma}[$\fc$-frozen chart (horizontal tangent)]\label{c:horizontal_frozen_critical1}
Fix $\delta_0>0$ and let $(x_0,s_0)\in \fP$.  Assume there exists a horizontal tangent location $(x_0',s_0)\in \fA$ such that
the line through $(x_0,s_0)$ and $(x_0',s_0)$ is tangent to the arctic curve at $(x_0',s_0)$, and moreover $(x_0, s_0)$ is at least distance $\delta_0$ from $(x_0',s_0)$. Then for every sufficiently small $\fc>0$ there exists
$\delta=\delta(\fc)>0$, also sufficiently small and depending on $\delta_0$,  such that the following holds.

Let
\[
\wt \fU:=\{{\wt w}\in\bC:\ |{\wt w}|\le \fc\}.
\]
On $\wt \fU$ the Riemann surface $\cC$ can be parametrized as $(\wt \chi({\wt w})/(1-\wt \chi({\wt w})),{\wt w})$. We will therefore identify
$\wt \fU$ with its image in $\cC$ under the map $w\mapsto (\wt \chi({\wt w})/(1-\wt \chi({\wt w})), {\wt w})$. Moreover, for all
\(\wt w\in\wt \fU\),
\begin{align}\label{e:wtchi}
|\wt \chi(\wt w)| \asymp 1.
\end{align}

Now let $(x,s)\in\fL$ satisfy $\|(x,s)-(x_0,s_0)\|_2\le \delta$.
The tiling action $\wt S(\,\cdot\,;x,s)$ has exactly one (formal) critical point $\wt w_c$ in $\wt \fU$ satisfying 
\begin{align}\label{e:twc_bound}
|{\wt w}_c|\asymp |s-s_0|,
\end{align}
Moreover, for each such critical point ${\wt w}_c$ and every ${\wt w}\in\wt \fU$, 
\begin{align}
\wt S({\wt w};x,s)-\wt S({\wt w}_c;x,s)=d{\wt w} +\cE({\wt w}),\quad d:=\wt S'(0;x_0,s_0),\quad |d|\asymp 1.
\end{align}
The error term satisfies
\begin{align}\label{e:hor_tangent_frozen_err}
|\cE({\wt w})|\leq C(\delta\ln(1/\delta)+|{\wt w}|^2)\leq\frac{|d|\fc}{100}.
\end{align}

In this situation, we call $\wt \fU$ a (horizontal tangent) $\fc$-frozen chart (centered at $0$), and we say that
the point $(x,s)$ is \emph{adapted} to $\wt \fU$.

We also introduce the corresponding chart in the $w$-plane, together with the
corresponding pulled-back critical points, by pulling back the $\wt w$-chart:
\[
\fU=\{w:\wt w\in \wt\fU\}
    =\left\{x'_0-\frac{1}{\wt w}:\wt w\in \wt\fU\right\},
    \qquad
w_c=x'_0-\frac{1}{\wt w_c}.
\]
We call $\fU$ a $\fc$-frozen chart (horizontal tangent) as well, and say that the point
$(x,s)$ is \emph{adapted} to $\fU$.

\end{lemma}

As in the vertical tangency case, the horizontal tangency point
$(x'_0,s_0)\in\fA$ lies on the common boundary of two adjacent curvilinear
triangles $\fT_A$ and $\fT_B$. There are two cases for the integration
contours.

Suppose first that $(x_0,s_0)$ does not lie on the common boundary of
$\fT_A$ and $\fT_B$, and let $\fT\in\{\fT_A,\fT_B\}$ denote the region
containing $(x_0,s_0)$. Depending on whether
\[
(x_0,s_0)\in\ell_-(\infty;\fT)
\qquad\text{or}\qquad
(x_0,s_0)\in\ell_+(\infty;\fT),
\]
we associate with $(x_0,s_0)$ the circular component of the local descent
or ascent contour, respectively, associated with $(x'_0,s_0)$ for $\fT$
in \Cref{s:horizontal_tangent}. 

If $(x_0,s_0)$ lies on the common boundary of the two adjacent curvilinear
triangles, we associate with $(x_0,s_0)$ both the circular component of the
local descent contour (blue) and that of the local ascent contour (red)
associated with $(x'_0,s_0)$ for either $\fT_A$ or $\fT_B$ in
\Cref{s:horizontal_tangent}.

 We define the corresponding local descent and ascent contours in the \(w\)-plane by pulling back the contours in the \(\wt w\)-plane:
\begin{align}
\sfC^{\rm d}(\infty)
=\{w: \wt w\in \wt \sfC^{\rm d}(0)\}, \quad 
\sfC^{\rm a}(\infty)
=\{w: \wt w\in \wt \sfC^{\rm a}(0)\},\quad \wt w=\frac{1}{x_0'-w}.
\end{align}

The following lemma shows that the local descent and ascent paths can be deformed into steepest-descent and steepest-ascent paths with negligible error. Its proof is deferred to \Cref{s:horizontal_frozen_chart_proof}.
\begin{lemma}\label{l:local_descent_deformation3}
Adopt the assumptions and notation of
\Cref{c:horizontal_frozen_critical1}. The statement of
\Cref{l:local_descent_deformation} remains valid after replacing the
vertical-tangency setting with the horizontal tangency, setting $w_0=0$ and replacing $S, \sfC, \sfD$ by $\wt S, \wt \sfC, \wt \sfD$, and the definition of $r_n$ by the following:
\begin{enumerate}

\item \emph{Close to the tangent line.} If
$|s-s_0|\leq (\ln n)^3/n$, we define
$r_n=(\ln n)^5/n$.

\item \emph{Other frozen regime.} If
$|s-s_0|>(\ln n)^3/n$, we define
$r_n:=(\ln n)^2 (|(s-s_0)|/n)^{1/2}$.

\end{enumerate}
\end{lemma}

%\begin{proof}
%The proof is identical to that of \Cref{l:local_descent_deformation} and is
%therefore omitted.
%\end{proof}

\section{Single-Contour Integral}
\label{s:single_integral}
In this section, we describe the single-contour integral
\eqref{e:single_term0}, its contours, the corresponding steepest-descent
contours, and related properties.
\subsection{Contours for the single-contour integral}

We recall the contour $\sfC(w_0;(x,s),(y,t))$ for $w_0\in\bC_+$ from \Cref{d:defCw}. We now extend this definition to the case
$w_0=E\in\bR$.
\begin{definition}\label{def:CE}
For any \(E\in\bR\), we define the oriented contour
$
\sfC(E;(x,s),(y,t))\subset\bC
$
as follows:
\begin{enumerate}
\item
If \(s\geq t\), then \(\sfC(E;(x,s),(y,t))\) consists of a path from
\(-\ri\infty\) to \(E-\ri0\) through the lower half-plane, followed by a
path from \(E+\ri0\) to \(+\ri\infty\) through the upper half-plane.

\item
If \(s<t\), then \(\sfC(E;(x,s),(y,t))\) is a loop based at \(E\). It
starts from \(E+\ri0\), travels through the upper half-plane to a point in
\[
\bigl(\max\{x-s,y-t\},\,\min\{x,y\}\bigr),
\]
and then returns through the lower half-plane to \(E-\ri0\).
\end{enumerate}
\end{definition}
For the general kernel ansatz, the single-contour integral has the same form
as the corresponding integral in \eqref{e:bulk_ansatz}, except that the
contour \(\sfC\) is replaced by
$
\sfC\bigl(\xi;(x,s),(y,t)\bigr),
$
where \(\xi\in\bC\) depends on the position of \((x,s)\) relative to
\((y,t)\). In the remainder of this section, we collect several useful
identities for analyzing such single-contour integral.

%
%In the following we give a geometric criterion for when a chart from the $(y,t)$-collection and a chart from the
%$(x,s)$-collection are concentric, with the $(y,t)$-chart contributing a local \emph{steepest--ascent} path and the
%$(x,s)$-chart contributing a local \emph{steepest--descent} path.
%
%\medskip
%\noindent\textbf{Liquid-type concentric charts.}
%If the concentric charts are of liquid type, then $(y,t),(x,s)\in\fL$ are bounded away from the arctic boundary
%and the ramification points, and each point is assigned a pair of complex-conjugate liquid charts (see
%\Cref{f:bulk_path}). In this case both ascent and descent contours are available, and no further discussion is
%needed.

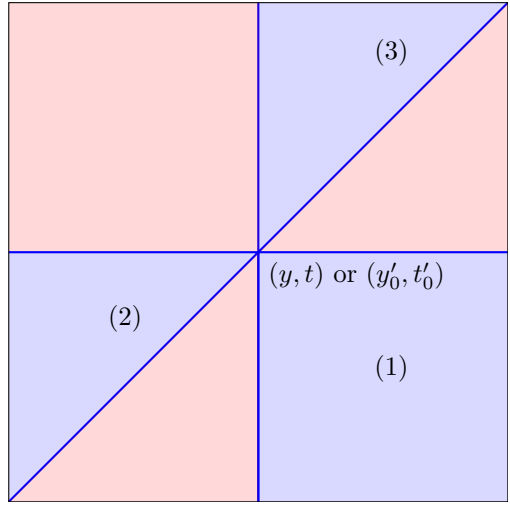
\begin{figure}
% Preamble:
% \usepackage{tikz}
% \usetikzlibrary{calc}

\begin{tikzpicture}[scale=1.1, font=\small]
  % ---- plotting window ----
  \def\xmin{-3}
  \def\xmax{ 3}
  \def\ymin{-3}
  \def\ymax{ 3}

  % ---- clip to a box ----
  \clip (\xmin,\ymin) rectangle (\xmax,\ymax);

  % ---- fill the 6 regions ----
  % i1: x>=0, s>=x, s>0  (triangle in first quadrant above diagonal)
  \fill[blue, opacity=0.15] (0,0) -- (\xmax,\xmax) -- (0,\ymax) -- cycle;
  \draw[blue,thick] (\xmax,\xmax) --(0,0)-- (0,\xmax);
 
  % i2: x>=0, 0<s<x  (triangle in first quadrant below diagonal)
  \fill[red, opacity=0.15] (0,0) -- (\xmax,0) -- (\xmax,\xmax) -- cycle;
  \draw[blue,thick] (0,\xmin) --(0,0);

  % i3: x>=0, s<=0, s<x  (fourth quadrant; effectively x>=0, s<=0)
  \fill[blue, opacity=0.15] (0,\ymin) rectangle (\xmax,0);
  \draw[blue,thick] (0,\xmin) --(0,0)--(\xmax, 0);
  
  % i4: x<0, s>0 (second quadrant)
  \fill[red, opacity=0.15] (\xmin,0) rectangle (0,\ymax);

  % i5: x<0, s<=0, s>=x (wedge between s=0 and s=x for x<0)
  \fill[blue, opacity=0.15] (\xmin,\xmin) -- (\xmin,0) -- (0,0) -- cycle;
  \draw[blue,thick] (\xmin,0) --(0,0)--(\xmin, \xmin);

  % i6: x<0, s<x (below diagonal in third quadrant for x<0)
  \fill[red, opacity=0.15] (\xmin,\ymin) -- (0,\ymin) -- (0,0) -- cycle;

    % ---- labels ----
  \node at (1.6,2.4) {$(3)$};
  %\node at (1.6,0.8) {$(2)$};
  \node at (1.6,-1.4) {$(1)$};

 % \node at (-1.6,1.4) {$(4)$};
  \node at (-1.6,-0.8) {$(2)$};
 % \node at (-1.6,-2.4) {$(6)$};
  \node [below right]at (0,0) {$(y,t) \text{ or } (y_0',t_0')$};

  % ---- frame (optional) ----
  \draw[thin] (\xmin,\ymin) rectangle (\xmax,\ymax);
\end{tikzpicture}

\caption{There are three cases for a point relative to the location of $(y,t) \text{ or } (y_0',t_0')$.}
\label{f:l0}
\end{figure}

\begin{lemma}
Fix $(y,t)$ and $z_0\in \bC$, and write
$\sfC(z_0;x,s)=\sfC(z_0;(x,s),(y,t))$. Then, for any $(x,s)$,
\begin{align}
\frac{n}{2\pi \ri}
\left(
\int_{\sfC(z_0; x,s)}
-\int_{\sfC(z_0; x-1/n,s-1/n)}
-\int_{\sfC(z_0; x,s-1/n)}
\right)
P_{ns}(nz,nx)\,Q_{nt}(nz,ny)\,\rd z
=\delta_{(x,s),(y,t)}.
\end{align}
\end{lemma}

\begin{proof}
The claim follows from exactly the same argument as in \Cref{l:Aeq}, so we omit.
\end{proof}

\begin{lemma}\label{c:single_contour}
For $\max\{x-s,y-t\}< E< \min\{x,y\}$
\begin{align}\label{e:gH00}
\frac{1}{2\pi \ri}\int_{\sfC(E;(x,s),(y,t))} P_{ns}(nw,nx)\,Q_{nt}(nz,ny)  \rd z=\bm1(s\geq t)\bm1(x\geq y) {n(s-t)\choose n(x-y)}
\end{align}

For $E> \max\{x,y\}$
\begin{align}\label{e:gH10}
\frac{1}{2\pi \ri}\int_{\sfC(E;(x,s),(y,t))} P_{ns}(nw,nx)\,Q_{nt}(nz,ny)  \rd z=-\bm1(t>s)\bm1(x\geq y) {n(s-t)\choose n(x-y)}
\end{align}

For $E<\min\{x-s,y-t\}$
\begin{align}\label{e:gH1-1}
\frac{1}{2\pi \ri}\int_{\sfC(E;(x,s),(y,t))} P_{ns}(nw,nx)\,Q_{nt}(nz,ny)  \rd z=-\bm1(t>s)\bm1(y-t\geq x-s) {n(s-t)\choose n((y-t)-(x-s))}
\end{align}

\end{lemma}

\begin{remark}\label{r:order_condition}
We remark that \eqref{e:gH10} is nonzero, only if 
\begin{align}
t>s, \quad x\geq y
\end{align}
equivalently $(x,s)\in \rm{Region} (1)$ (see \Cref{f:l0}) relative to $(y,t)$; 
\eqref{e:gH1-1} is nonzero, only if 
\begin{align}
t>s, \quad y-t\geq x-s
\end{align}
equivalently $(x,s)\in \rm{Region} (2)$ (see \Cref{f:l0}) relative to $(y,t)$; \eqref{e:gH00} is nonzero, only if 
\begin{align}
y-t\geq x-s, \quad x\geq y
\end{align}
equivalently $(x,s)\in \rm{Region} (3)$ (see \Cref{f:l0}) relative to $(y,t)$.
\end{remark}

\begin{proof}[Proof of \Cref{c:single_contour}]
The first statement \eqref{e:gH00} follows from \Cref{c:single_int}.

For \eqref{e:gH10}, if $s\ge t$, then by definition $\sfC(E;(x,s),(y,t))$ is a  contour from $-\ri\infty$ to $+\ri\infty$ crossing the real
axis at $E>\max\{x,y\}$. Since all poles of the integrand lie to the left of $E$, the contour may be deformed
to the far right without crossing any pole, and the integral vanishes. This gives the factor $\bm 1(t>s)$. Now assume $s<t$. if $x<y$, then the integrand has no poles enclosed by $\sfC(E;(x,s),(y,t))$, so the integral vanishes. Hence it remains to
consider the case $s<t$ and $x\ge y$. In this case, the integrand has poles at
\begin{align}\label{e:poles}
z=y,\ y+\frac1n,\ \cdots,\ x.
\end{align}
Thus we may deform $\sfC(E;(x,s),(y,t))$ to a counterclockwise contour enclosing precisely these poles.
\begin{align}
&\frac{1}{2\pi \ri}\int_{\sfC(E;(x,s),(y,t))}
\frac{\Gamma(ns+1)}{\Gamma(nt)}
\frac{\Gamma(n(y-z))\Gamma(n(z-(y-t)))}
{\Gamma(n(x-z)+1)\Gamma(n(z-(x-s))+1)}\,\rd z \notag\\
&=
\frac{1}{2\pi \ri}\int_{\sfC(E;(x,s),(y,t))}
\frac{\Gamma(ns+1)}{\Gamma(nt)}
\frac{(-1)^{n(x-y)+1}}{n(z-y)\,n(z-y-1/n)\cdots n(z-x)}
\frac{\Gamma(n(z-(y-t)))}{\Gamma(n(z-(x-s))+1)}\,\rd z,
\end{align}
where we compute the ratio of two Gamma functions using $\Gamma(z+1)=z\Gamma(z)$.
Taking residues at \eqref{e:poles}, we obtain
\begin{align}\begin{split}
&\sum_{j=0}^{n(x-y)}
\frac{\Gamma(ns+1)}{\Gamma(nt)}
\frac{(-1)^{j+1}}{j!(n(x-y)-j)!}
\frac{\Gamma(nt+j)}{\Gamma(ns-n(x-y)+j+1)} \\
&=
-\sum_{j=0}^{n(x-y)}
\frac{\Gamma(nt+j)}{\Gamma(nt)\,j!}
\frac{\Gamma(ns+1)}{(n(x-y)-j)!\,\Gamma(ns-n(x-y)+j+1)}\\
&=
-\sum_{j=0}^{n(x-y)}{-nt\choose j}{ns\choose n(x-y)-j}={n(s-t)\choose n(x-y)},
\end{split}\end{align}
where the first line follows from computing the residual at $z=y+j/n$, namely the $j$-th term corresponds to the residual at $z=y+j/n$; in the second statement we rearrange the expression; the last statement follows from Chu–Vandermonde identity (see \cite[Section 1.1]{kuznetsov2006orthogonal}).

For \eqref{e:gH1-1}, by the same argument as for \eqref{e:gH1-1} the integral vanishes unless $s<t$ and $y-t\ge x-s$. In this case the integrand has poles at $y-t, y-t-1/n, \cdots, x-s$. %, and at $y-t, y-t-1/n, y-t-2/n,\cdots, x-s$. 
Thus we may deform $\sfC(E;(x,s),(y,t))$ to a \emph{clockwise} contour enclosing precisely these poles. By the same argument as in \eqref{e:gH10}, we have
\begin{align}
&\phantom{{=}}\frac{1}{2\pi \ri}\int_{\sfC(E;(x,s),(y,t))} \frac{\Gamma(ns+1)}{\Gamma(nt)}\frac{\Gamma(n(y-z))\Gamma(n(z-(y-t)))}{\Gamma(n(x-z)+1)\Gamma(n(z-(x-s))+1)}  \rd z\\
&=\frac{1}{2\pi \ri}\int_{\sfC(E;(x,s),(y,t))} \frac{\Gamma(ns+1)}{\Gamma(nt)}
\frac{\Gamma(n(y-z))}{\Gamma(n(x-z)+1)}\frac{1}{n(z-(x-s)) n(z-(x-s)-1/n)\cdots n(z-(y-t))} 
\rd z\\
&=-\sum_{j=0}^{n((y-t)-(x-s))}\frac{\Gamma(ns+1)}{\Gamma(nt)}
\frac{(-1)^{n((y-t)-(x-s))-j}}{j! (n((y-t)-(x-s))-j)!}\frac{\Gamma(n(y-(x-s))-j)}{\Gamma(ns-j+1)}\\
&=-\sum_{j=0}^{n((y-t)-(x-s))}
\frac{(-1)^{n((y-t)-(x-s))-j}\Gamma(n(y-(x-s))-j)}{\Gamma(nt) (n((y-t)-(x-s))-j)!}\frac{\Gamma(ns+1)}{j!\Gamma(ns-j+1)}\\
&=-\sum_{j=0}^{n((y-t)-(x-s))}{-nt\choose n((y-t)-(x-s))-j}{ns \choose j}=-{n(s-t)\choose n((y-t)-(x-s))}.
% \frac{1}{n(z-(y-t)) n(z-(y-t-1/n)) \cdots n(z-(x-s))}  \rd z
\end{align}

\end{proof}

\begin{figure}
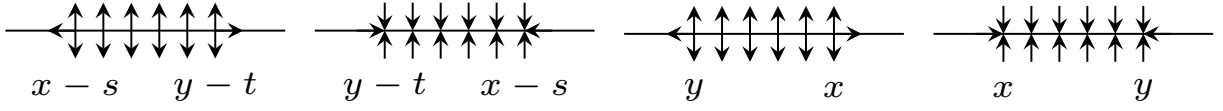

  \centering

  %================= Row 1 =================
  \begin{subfigure}{0.24\textwidth}
    \centering
    \resizebox{\linewidth}{!}{%

    }
  \end{subfigure}
  \caption{Steepest descent paths around $[x-s,y-t]$ and $[y,x]$.}
  \label{f:in-out}
\end{figure}

\subsection{Phase of the single-contour integral}

We recall the single-contour integral from \eqref{e:single_term0}:
\begin{equation}
\frac{n}{2\pi \ri}
\int_{\sfC(z_0;(x,s),(y,t))}
P_{ns}(nz,nx)\,Q_{nt}(nz,ny)\,\rd z,
\end{equation}
where \(z_0\in\bC\) and the contour
\(\sfC(z_0;(x,s),(y,t))\) is defined as in \Cref{d:defCw} and \Cref{def:CE}.

Using Stirling's formula \eqref{e:logGamma}, uniformly for \(z\) bounded
away from the cuts between \(x-s\) and \(y-t\), and between \(x\) and \(y\),
we have
\begin{align*}
P_{ns}(nz,nx)Q_{nt}(nz,ny)
&=
\frac{\Gamma(ns+1)}{\Gamma(nt)}
\frac{\Gamma(n(y-z))\Gamma(n(z-(y-t)))}
{\Gamma(n(x-z)+1)\Gamma(n(z-(x-s))+1)}
\notag\\
&=
\frac{1}{n}
\sqrt{
\frac{s}{(x-z)(z-(x-s))}
\frac{t}{(y-z)(z-(y-t))}
}
e^{n(S(z;x,s)-S(z;y,t))+\OO(1/n)}.
\end{align*}
Here
\begin{align*}
S(z;x,s)-S(z;y,t)
&=
s\ln s
-(x-z)\ln(x-z)
-(z-(x-s))\ln(z-(x-s))
\notag\\
&\quad
-\left(
t\ln t
-(y-z)\ln(y-z)
-(z-(y-t))\ln(z-(y-t))
\right).
\end{align*}

Observe that if \((x,s)=(y,t)\), then
\[
S(z;x,s)-S(z;y,t)\equiv0.
\]
Below, the ``nonconstant case'' refers to
\((x,s)\neq(y,t)\). In this case, the critical points of
\(S(z;x,s)-S(z;y,t)\) are determined by 
\begin{align}
S'(z;x,s)-S'(z;y,t)
=
\ln\frac{x-z}{z-(x-s)}
-
\ln\frac{y-z}{z-(y-t)}
=
0.
\end{align}
Equivalently,
\[
(x-z)(z-(y-t))
=
(y-z)(z-(x-s)).
\]
Thus, if \(x\neq y\), \(x-s\neq y-t\), and \(s\neq t\), the unique
candidate critical point is
\begin{align}\label{e:zc_loc}
z_c
=
\frac{xt-ys}{t-s}
\in\bR.
\end{align}
Moreover,
\begin{align}
S''(z_c;x,s)-S''(z_c;y,t)
=
\frac{(t-s)^3}
{st(x-y)(x-y+t-s)}.
\end{align}
Thus, \(z_c\) is a real critical point.

As \(|z|\to\infty\), we have
\begin{align}\label{e:Sder_exp}
S'(z;x,s)-S'(z;y,t)
=
\frac{t-s}{z}
+
\frac{-2xs+s^2+2yt-t^2}{2z^2}
+
\OO\left(\frac{1}{z^3}\right).
\end{align}
On the circle \(z=Re^{\ri\theta}\), with \(R\) sufficiently large, the
steepest-descent direction through \(z\) is
\begin{align}
-\overline{(S'(z;x,s)-S'(z;y,t))}
=
-\frac{(t-s)e^{\ri\theta}}{R}
+
\OO\left(\frac{1}{R^2}\right).
\end{align}
Hence:
\begin{enumerate}
\item
If \(t>s\), then \(\infty\) is a source: descent paths emanate from
infinity and point inward in the \(z\)-plane.

\item
If \(t<s\), then \(\infty\) is a sink: descent paths point outward in the
\(z\)-plane and terminate at infinity.

\item
If \(t=s\), then
\[
S'(z;x,s)-S'(z;y,t)
=
\frac{t(y-x)}{z^2}
+
\OO\left(\frac{1}{z^3}\right),
\qquad
|z|\to\infty.
\]
Hence, in the nonconstant case, \(\infty\) is a critical point that is
neither a sink nor a source.
\end{enumerate}

Moreover, after fixing the branches of the logarithms, the function
\(S'(z;x,s)-S'(z;y,t)\) is holomorphic on
\[
\bC\setminus
\left(
[\min\{x-s,y-t\},\max\{x-s,y-t\}]
\cup
[\min\{x,y\},\max\{x,y\}]
\right).
\]
On these intervals, the behavior of the steepest-descent paths depends on
the relative positions of \(x,y\) and \(x-s,y-t\); see
\Cref{f:in-out}.
\begin{enumerate}
\item
If \(x>y\), the interval \([y,x]\) is a source: the steepest-descent paths
point outward from it. If \(x<y\), the interval \([x,y]\) is a sink: the
steepest-descent paths point inward toward it. If \(x=y\), then, in the
nonconstant case, \(z_c=x\) is neither a sink nor a source and is not a critical
point.

\item
If \(y-t>x-s\), the interval \([x-s,y-t]\) is a source: the
steepest-descent paths point outward from it. If \(y-t<x-s\), the interval
\([y-t,x-s]\) is a sink: the steepest-descent paths point inward toward it.
If \(x-s=y-t\), then, in the nonconstant case, \(z_c=x-s\) is neither a sink
nor a source and is not a critical point.
\end{enumerate}

\subsection{Steepest descent paths for the single-contour integral}

\label{s:pof_single_integral}

For any \(z_0\in\bC_+\), suppose first that
$
(x,s)\neq(y,t).
$
Let
$
\mathsf{D}(z_0;(x,s),(y,t))
$
denote the steepest--descent path for the phase
$
S(z;x,s)-S(z;y,t)
$
starting at $z_0$.

The contour
$\mathsf{D}(z_0;(x,s),(y,t))$ is oriented so that it starts
from $z_0$. As will be shown below, it remains in the upper half-plane and
terminates either at the critical point $z_c$, when it exists, or at one of the following, provided it is a sink:
\begin{align}
\infty,\qquad
[\min\{x-s,y-t\},\max\{x-s,y-t\}],\qquad
[\min\{x,y\},\max\{x,y\}].
\end{align}
If
$
(x,s)=(y,t),
$
then
$
S(z;x,s)-S(z;y,t)\equiv0.
$
In this degenerate case, we choose
$
\sfD(z_0;(x,s),(y,t))
$
to be any path in the upper half-plane from \(z_0\) to \(\infty\).

We write
\[
\sfD(\overline{z_0};(x,s),(y,t))
\]
for the complex conjugate of
\(\sfD(z_0;(x,s),(y,t))\), with the orientation reversed, so that it ends
at \(\overline{z_0}\). We then define the full oriented contour by
\begin{align}\label{e:def_full_single_descent_contour}
\widehat{\sfD}(z_0;(x,s),(y,t))
:=
\sfD(z_0;(x,s),(y,t))
\cup
\sfD(\overline{z_0};(x,s),(y,t)).
\end{align}
We extend this definition to \(E\in\bR\) by taking the boundary values
obtained by perturbing \(E\) infinitesimally into the upper and lower
half-planes, denoted by \(E+\ri0\) and \(E-\ri0\), respectively. We then
define
\begin{align}
\widehat{\sfD}(E;(x,s),(y,t))
:=
\sfD(E+\ri0;(x,s),(y,t))
\cup
\sfD(E-\ri0;(x,s),(y,t)).
\end{align}

\begin{lemma}\label{l:deform_single_contour}
For any \(z_0\in\bC_+\cup\bR\), the contour
$
\sfC(z_0;(x,s),(y,t))
$
can be deformed to the steepest-descent contour
$
\widehat{\sfD}(z_0;(x,s),(y,t)).
$
In particular,
\begin{align}
&\frac{n}{2\pi\ri}
\int_{\sfC(z_0;(x,s),(y,t))}
P_{ns}(nz,nx)\,Q_{nt}(nz,ny)\,\rd z
=
\frac{n}{2\pi\ri}
\int_{\widehat{\sfD}(z_0;(x,s),(y,t))}
P_{ns}(nz,nx)\,Q_{nt}(nz,ny)\,\rd z.
\end{align}
\end{lemma}

\begin{proof}
For simplicity, write
\[
\sfC(z_0):=\sfC(z_0;(x,s),(y,t)),
\quad
\sfD(z_0):=\sfD(z_0;(x,s),(y,t)),
\quad
\widehat{\sfD}(z_0)
:=
\widehat{\sfD}(z_0;(x,s),(y,t)).
\]
It suffices to consider \(z_0\in\bC_+\); the result for \(z_0\in\bR\)
follows by taking the boundary value \(z_0+\ri0\).

Recall that the two intervals are disjoint:
\[
\max\{x-s,y-t\}<\min\{x,y\}.
\]
If \(s\geq t\), then, by \Cref{l:integrand_prop}, the integrand decays as
\(\OO(z^{-2})\) at infinity. In this case, \(\sfC(z_0)\) consists of an
upper-half-plane piece running from \(z_0\) to \(\infty\), followed by a
lower-half-plane piece running from \(\infty\) to \(\overline{z_0}\).

If \(s<t\), then \(\sfC(z_0)\) consists of an upper-half-plane piece
running from \(z_0\) to a point
\[
\xi_0\in
\bigl(\max\{x-s,y-t\},\min\{x,y\}\bigr),
\]
followed by a lower-half-plane piece running from \(\xi_0\) to
\(\overline{z_0}\). In both cases, we regard \(\sfC(z_0)\) as a single
oriented contour obtained by concatenating these two pieces.

If \((x,s)=(y,t)\), then both \(\sfC(z_0)\) and
\(\widehat{\sfD}(z_0)\) consist of an upper-half-plane path from \(z_0\)
to \(\infty\), together with a lower-half-plane path from \(\infty\) to
\(\overline{z_0}\). By \Cref{l:integrand_prop}, the integrand is analytic
on \(\bC\setminus\bR\) and decays as \(\OO(z^{-2})\) at infinity. Hence
Cauchy's theorem gives the desired equality.

Assume now that \((x,s)\neq(y,t)\). Recall from \eqref{e:zc_loc} that,
when the critical point exists, it is
\[
z_c=\frac{xt-ys}{t-s}.
\]
The following identities will be used to locate it:
\begin{align}
z_c-y
&=
\frac{t(x-y)}{t-s},
&
z_c-x
&=
\frac{s(x-y)}{t-s},
\notag\\
z_c-(y-t)
&=
\frac{t(x-y+t-s)}{t-s},
&
z_c-(x-s)
&=
\frac{s(x-y+t-s)}{t-s}.
\end{align}
Moreover,
\begin{align}
S''(z_c;x,s)-S''(z_c;y,t)
=
\frac{(t-s)^3}
{st(x-y)(x-y+t-s)}.
\end{align}
Thus, when
$
S''(z_c;x,s)-S''(z_c;y,t)>0,
$
the two nonreal trajectories through \(z_c\) are descent trajectories.
When
$
S''(z_c;x,s)-S''(z_c;y,t)<0,
$
they are ascent trajectories when oriented away from \(z_c\), and hence
steepest-descent trajectories may terminate at \(z_c\).

It remains to distinguish six cases.

\begin{enumerate}
\item\label{i1}
Suppose that
\[
x\geq y,
\qquad
y-t\geq x-s.
\]
Since \((x,s)\neq(y,t)\), this implies \(s>t\). Hence \(\infty\) is the
only sink. If \(z_c\) exists, then
$
y-t<z_c<y,
$
and the two nonreal trajectories through \(z_c\) are descent trajectories.
Therefore, the upper steepest-descent path \(\sfD(z_0)\) terminates at
\(\infty\). Taking also the conjugate lower-half-plane path, the full
contour \(\sfC(z_0)\) can be deformed to
\(\widehat{\sfD}(z_0)\).

\item\label{i2}
Suppose that
\[
x\geq y,
\qquad
y-t<x-s,
\qquad
t<s.
\]
Then both \([y-t,x-s]\) and \(\infty\) are sinks, while \([y,x]\) is the
only source. If \(z_c\) exists, then
$
z_c<y-t,
$
and the two nonreal trajectories through \(z_c\) are ascent trajectories.
Therefore, \(\sfD(z_0)\) terminates at \([y-t,x-s]\), at \(z_c\), or at
\(\infty\).

By \Cref{l:integrand_prop}, the integrand is \(\OO(z^{-2})\) as
\(|z|\to\infty\) and is analytic across \((-\infty,y)\). Every finite
terminal point of \(\sfD(z_0)\) lies in \((-\infty,x-s]\). Hence the
upper-half-plane piece of \(\sfC(z_0)\) can be deformed to
\(\sfD(z_0)\), and the conjugate deformation applies in the lower
half-plane. Therefore, \(\sfC(z_0)\) can be deformed to
\(\widehat{\sfD}(z_0)\).

\item\label{i3}
Suppose that
\[
x\geq y,
\qquad
y-t<x-s,
\qquad
t\geq s.
\]
Then \([y-t,x-s]\) is the only sink. If \(z_c\) exists, then necessarily
\(t>s\), and
$
z_c>x.
$
In this case, the two nonreal trajectories through \(z_c\) are descent
trajectories. Hence \(\sfD(z_0)\) terminates at \([y-t,x-s]\).

By \Cref{l:integrand_prop}, the integrand is analytic across
\((-\infty,y)\). If \(t>s\), the contour \(\sfC(z_0)\) passes through
\(\xi_0\). If \(t=s\), its endpoint at infinity may be moved using the
\(\OO(z^{-2})\) decay from \Cref{l:integrand_prop}. In either case, the
upper-half-plane piece of \(\sfC(z_0)\) can be deformed to
\(\sfD(z_0)\), and the conjugate deformation applies in the lower
half-plane. Thus, \(\sfC(z_0)\) can be deformed to
\(\widehat{\sfD}(z_0)\).

\item\label{i4}
Suppose that
\[
x<y,
\qquad
y-t\geq x-s,
\qquad
t<s.
\]
Then both \([x,y]\) and \(\infty\) are sinks, while
\([x-s,y-t]\) is the only source. If \(z_c\) exists, then
$
z_c>y,
$
and the two nonreal trajectories through \(z_c\) are ascent trajectories.
Therefore, \(\sfD(z_0)\) terminates at \([x,y]\), at \(z_c\), or at
\(\infty\).

By \Cref{l:integrand_prop}, the integrand is \(\OO(z^{-2})\) as
\(|z|\to\infty\) and is analytic across \((y-t,\infty)\). Every finite
terminal point of \(\sfD(z_0)\) lies in \([x,\infty)\). Hence the
upper-half-plane piece of \(\sfC(z_0)\) can be deformed to
\(\sfD(z_0)\), and the conjugate deformation applies in the lower
half-plane. Therefore, \(\sfC(z_0)\) can be deformed to
\(\widehat{\sfD}(z_0)\).

\item\label{i5}
Suppose that
\[
x<y,
\qquad
y-t\geq x-s,
\qquad
t\geq s.
\]
Then \([x,y]\) is the only sink. If \(z_c\) exists, then necessarily
\(t>s\), and
$
z_c<x-s.
$
In this case, the two nonreal trajectories through \(z_c\) are descent
trajectories. Hence \(\sfD(z_0)\) terminates at \([x,y]\).

By \Cref{l:integrand_prop}, the integrand is analytic across
\((y-t,\infty)\). If \(t>s\), the contour \(\sfC(z_0)\) passes through
\(\xi_0\). If \(t=s\), its endpoint at infinity may be moved using the
\(\OO(z^{-2})\) decay from \Cref{l:integrand_prop}. In either case, the
upper-half-plane piece of \(\sfC(z_0)\) can be deformed to
\(\sfD(z_0)\), and the conjugate deformation applies in the lower
half-plane. Thus, \(\sfC(z_0)\) can be deformed to
\(\widehat{\sfD}(z_0)\).

\item\label{i6}
Suppose that
\[
x<y,
\qquad
y-t<x-s.
\]
Then necessarily \(s<t\). Hence both \([y-t,x-s]\) and \([x,y]\) are
sinks, while \(\infty\) is the only source. If \(z_c\) exists, then
$
x-s<z_c<x,
$
and the two nonreal trajectories through \(z_c\) are ascent trajectories.
Therefore, \(\sfD(z_0)\) terminates at \([y-t,x-s]\), at \(z_c\), or at
\([x,y]\).

By \Cref{l:integrand_prop}, the integrand is analytic across \(\bR\).
Hence the upper-half-plane piece of \(\sfC(z_0)\) can be deformed to
\(\sfD(z_0)\), and the conjugate deformation applies in the lower
half-plane. Therefore, \(\sfC(z_0)\) can be deformed to
\(\widehat{\sfD}(z_0)\).
\end{enumerate}

This proves that the full oriented contour \(\sfC(z_0)\) can be deformed
to \(\widehat{\sfD}(z_0)\) in all cases. Cauchy's theorem therefore gives
the desired identity of contour integrals. The case \(z_0\in\bR\) follows
by taking the boundary value \(z_0+\ri0\).
\end{proof}

\begin{lemma}
Assume that \((x,s)\neq(y,t)\). Let \(z_0,w_0\in\bC_+\) satisfy
\begin{align}\label{e:ImSdiff}
\Im\bigl[S(z_0;x,s)-S(z_0;y,t)\bigr]
=
\Im\bigl[S(w_0;x,s)-S(w_0;y,t)\bigr].
\end{align}
Then \(z_0\) and \(w_0\) lie on the same steepest-descent path of
$
S(z;x,s)-S(z;y,t).
$
\end{lemma}

\begin{proof}
By the discussion above, in each parameter regime there is either a unique sink or a unique source. We treat the case where there is a unique sink, namely \Cref{i1,i3,i5}; the case of a unique source is analogous, replacing steepest descent by steepest ascent.

We recall that $ \sfD(z_0)= \sfD(z_0; (x,s),(y,t))$ and $ \sfD(w_0)= \sfD(w_0; (x,s),(y,t))$ are the steepest-descent paths starting from $z_0$ and $w_0$, respectively. Both curves lie in $\bC_+$ and terminate at the unique sink. Moreover, by \eqref{e:ImSdiff}, the quantity $\Im[S(z;x,s)-S(z;y,t)]$ takes the same value along both $ \sfD(z_0)$ and $ \sfD(w_0)$. 

In the setting of \Cref{i1}, both $ \sfD(z_0)$ and $ \sfD(w_0)$ terminate at $\infty$.
Fix $R>0$ large and consider the upper semicircle $\{Re^{\ri\theta}:0<\theta<\pi\}$. Using \eqref{e:Sder_exp}, we compute for $0<\theta<\pi$,
\begin{align}
&\phantom{{}={}}\del_\theta \Im[S(Re^{\ri \theta};x,s)-S(Re^{\ri \theta};y,t)]
=R \Re[e^{\ri\theta}(S'(Re^{\ri \theta};x,s)-S'(Re^{\ri \theta};y,t))]\\
&=
R\Re\left[
e^{\ri\theta}
\left(
\frac{t-s}{Re^{\ri\theta}}
+\OO\left(\frac{1}{R^2}\right)
\right)
\right]
=
(t-s)+\OO\left(\frac{1}{R}\right)<0.
\end{align}
So $\theta\mapsto \Im[S(Re^{\ri \theta};x,s)-S(Re^{\ri \theta};y,t)]$ is strictly decreasing on $(0,\pi)$. For \(R\) sufficiently large, both $ \sfD(z_0)$ and $ \sfD(w_0)$ intersect the upper semicircle $\{Re^{\ri\theta}:0<\theta<\pi\}$ exactly once, and by the monotonicity above they must intersect it at the same point. Therefore, the two paths share a common point; by uniqueness of the gradient flow, they must overlap.

In the setting of \Cref{i3}, both $ \sfD(z_0)$ and $ \sfD(w_0)$ terminate at the interval $[y-t,x-s]$. Fix a small $\varepsilon>0$, and consider the curve $\omega=\{y-t+\sqrt 2\varepsilon e^{\ri\theta}: 3\pi/4\leq \theta<\pi\}\cup \{u+\ri \varepsilon: y-t-\varepsilon\leq u\leq x-s+\varepsilon\}\cup \{x-s+\sqrt 2\varepsilon e^{\ri\theta}: 0<\theta\leq \pi/4\}$.
We compute, for $3\pi/4<\theta<\pi$,
\begin{align}
&\phantom{{}={}}
\del_\theta
\Im\left[
S(y-t+\sqrt 2\varepsilon e^{\ri\theta};x,s)
-
S(y-t+\sqrt 2\varepsilon e^{\ri\theta};y,t)
\right]
\notag\\
&=
\sqrt 2\varepsilon
\Re\left[
e^{\ri\theta}
\left(
S'(y-t+\sqrt 2\varepsilon e^{\ri\theta};x,s)
-
S'(y-t+\sqrt 2\varepsilon e^{\ri\theta};y,t)
\right)
\right]
\notag\\
&=
\sqrt 2\varepsilon
\Re\left[
e^{\ri\theta}
\left(
\ln
\frac{x-(y-t)-\sqrt 2\varepsilon e^{\ri\theta}}
{(y-t)+\sqrt 2\varepsilon e^{\ri\theta}-(x-s)}
-
\ln
\frac{t-\sqrt 2\varepsilon e^{\ri\theta}}
{\sqrt 2\varepsilon e^{\ri\theta}}
\right)
\right]
\notag\\
&=
\sqrt 2\varepsilon
\Re\left[
e^{\ri\theta}
\left(
\OO(1)+\ln(\sqrt 2\varepsilon e^{\ri\theta})
\right)
\right]
=
\sqrt 2\varepsilon
\left(
\cos\theta\ln\varepsilon+\OO(1)
\right)
>0.
\end{align}
and the same argument gives for $0<\theta\leq \pi/4$
\begin{align}
&\del_\theta
\Im\left[
S(x-s+\sqrt 2\varepsilon e^{\ri\theta};x,s)
-
S(x-s+\sqrt 2\varepsilon e^{\ri\theta};y,t)
\right]
=
-\sqrt 2\varepsilon
\left(
\cos\theta\ln\varepsilon+\OO(1)
\right)
>0.
\end{align}
Finally, for $y-t-\varepsilon\leq u\leq x-s+\varepsilon$
\begin{align}
&\phantom{{}={}}
\del_u
\Im\left[
S(u+\ri\varepsilon;x,s)-S(u+\ri\varepsilon;y,t)
\right]
\notag\\
&=
\Im\left[
\ln
\frac{x-u-\ri\varepsilon}{u+\ri\varepsilon-(x-s)}
-
\ln
\frac{y-u-\ri\varepsilon}{u+\ri\varepsilon-(y-t)}
\right]
\notag\\
&=
\arg(u+\ri\varepsilon-(y-t))
-
\arg(u+\ri\varepsilon-(x-s))
+
\OO(\varepsilon)
\leq
-\frac{\pi}{4}+\OO(\varepsilon)<0.
\end{align}
So $ \Im[S(z;x,s)-S(z;y,t)]$ is strictly decreasing as $z$ travels along $\omega$ from left to right. For \(\varepsilon\) sufficiently small, each of $ \sfD(z_0)$ and $ \sfD(w_0)$ intersects $\omega$ exactly once. Consequently, they intersect $\omega$ at the same point and hence must overlap.

In the setting of \Cref{i5}, both $ \sfD(z_0)$ and $ \sfD(w_0)$ terminate at the interval $[x,y]$. By essentially the same argument as in the last case, we have that $ \sfD(z_0)$ and $ \sfD(w_0)$ must overlap. 

\end{proof}

\chapter{Ansatz for General Kernels}

\section{Construction of a Phase-Adapted Open Cover of $\fP$}
\label{s:assign_cover}

In this section, we construct a phase-adapted open cover of the polygonal domain. The neighborhoods in the cover are classified according to their local geometry, and each is equipped with a collection of pairwise disjoint charts.

\subsection{Assigning neighborhoods for points in $\fP$}
\label{s:assign_neighborhood}
This following proposition constructs, uniformly in a neighborhood of each
\((x,s)\in\fP\), a finite collection of pairwise disjoint phase-adapted
charts containing all critical points associated with every nearby point.

\begin{proposition}\label{p:construct_neighborhood1}
Fix \((x,s)\in\fP\). Then there exist a sufficiently small \(\fc>0\) and an
open neighborhood \(\fN_{(x,s)}\subset\fP\) of \((x,s)\) such that the following
holds.

There is a finite collection of pairwise disjoint \(\fc\)-charts on the
Riemann surface \(\cC\), together with the corresponding
steepest-descent/ascent paths introduced in \Cref{s:critical_point}, such that,
for every \((x',s')\in\fN_{(x,s)}\),
\begin{enumerate}
\item \((x',s')\) is adapted to each chart in the collection; and
\item the union of the charts contains all critical points associated with
\((x',s')\) in \Cref{p:associate_critical_points}.
\end{enumerate}

Depending on the location of \((x,s)\), the collection of charts has one of
the following forms:
\begin{enumerate}
\item if \((x,s)\) lies in the liquid region, then the collection consists of a
pair of complex-conjugate liquid charts or ramification charts;

\item if \((x,s)\) is close to a cusp or cusp-turning point, then the
collection consists of a single cusp chart or cusp-turning chart;

\item if \((x,s)\) is close to the arctic curve, but bounded away from cusp and cusp-turning point, then the collection consists
of an arctic chart, or a tangent chart, possibly together
with one additional frozen chart of regular, cusp, or  tangency type;

\item if \((x,s)\) lies in the interior of the frozen region, then the
collection consists of two or three frozen charts of regular, cusp, or  tangency type.
\end{enumerate}

If \(\fN_{(x,s)}\) meets a boundary shared by two adjacent curvilinear triangles
\(\fT_A\) and \(\fT_B\), then the steepest-descent/ascent paths may be chosen
relative to either \(\fT_A\) or \(\fT_B\).
\end{proposition}

We choose $\fN_{(x,s)}$ small enough so that the local geometric configuration is stable on $\fN_{(x,s)}$, in the
following sense. The tangent lines from $(x',s')$ to the arctic boundary that are relevant near $(x,s)$, together
with their associated critical points, vary continuously as $(x',s')$ moves in $\fN_{(x,s)}$. However, the list of
critical points recorded in \Cref{p:associate_critical_points} need not be locally constant: by definition, that
proposition retains only those tangent lines whose points of tangency lie on the portion of the arctic boundary
contained in the given curvilinear triangle $\fT$. When $(x',s')$ crosses an extended side, one of the tangency
points may move out of $\fT$, and the corresponding critical point is then no longer included in
\Cref{p:associate_critical_points}, even though the underlying tangent line and critical point still persist
continuously.

For the present chart construction, we retain the charts associated with all such locally continuous tangent lines,
including those whose tangency points may leave $\fT$ under this selection rule. In this way, although the subset of
critical points singled out in \Cref{p:associate_critical_points} may change across an extended side, the total
number and type of charts remain constant on $\fN_{(x,s)}$.

\begin{proof}[Proof of \Cref{p:construct_neighborhood1}]
We argue by cases, according to the location of $(x,s)$. In each case we proceed in the same way: 
\begin{enumerate}
\item[(i)] choose $\fN_{(x,s)}$ small enough so that the local geometric configuration is stable on $\fN_{(x,s)}$.
\item[(ii)] choose $\fc$-charts in $\cC$ covering the associated critical points;
\item[(iii)]  shrink $\fN_{(x,s)}$ further, if needed, so that every $(x',s')\in\fN_{(x,s)}$ is adapted to each chosen chart.
And pairwise disjointness is ensured by taking $\fc$ sufficiently small.
\end{enumerate}
In the following discussion, for each point \((x_0,s_0)\in\fA\), we recall
the critical point \(z(x_0,s_0)\) from \eqref{e:emb}.

\medskip
\noindent\textbf{Liquid neighborhood.}
Assume $(x,s)\in\fL$ stays a fixed distance away from $\fA$ and from the ramification points. Choose $\fN_{(x,s)}$
with the same property. Then each $(x',s')\in\fN_{(x,s)}$ has two complex-conjugate critical points
$w_c$ and $\overline{w_c}$, and these vary continuously. We take a $\fc$-liquid chart  and its complex conjugate chart, with local paths as in \Cref{s:critical_bulk}, and shrink $\fN_{(x,s)}$ if needed
so that for all $(x',s')\in\fN_{(x,s)}$ the two critical points remain in the respective charts and $(x',s')$ is
adapted to both charts. We call $\fN_{(x,s)}$ a \emph{liquid neighborhood}.

\medskip
\noindent \textbf{Ramification neighborhood.}
Let \((x_0,s_0)\) be a ramification point, and set
\((x,s)=(x_0,s_0)\). As in \Cref{s:ramification}, associate with
\((x,s)\) a \(\fc\)-ramification chart in \(\cC_\ft\), its
complex-conjugate chart, and the corresponding local descent and
ascent paths.

Choose a sufficiently small neighborhood
\(\fN_{(x_0,s_0)}\) of \((x_0,s_0)\) such that, for every
\((x',s')\in\fN_{(x_0,s_0)}\), the tiling action
$
S_\ft(\,\cdot\,;x',s'+\ft)
$
has two complex-conjugate critical points \(w_{c,\ft}\) and
\(\overline{w_{c,\ft}}\). We further require that $(x',s')$ is adapted to the ramification chart and its complex-conjugate chart,
respectively. We call \(\fN_{(x_0,s_0)}\) a
\emph{ramification neighborhood}.

%$
%((x',s'+\ft),w_{c,\ft})
%$ and
%$
%((x',s'+\ft),\overline{w_{c,\ft}})
%$
%are adapted to the ramification chart and its complex-conjugate chart,
%respectively. We call \(\fN_{(x_0,s_0)}\) a
%\emph{ramification neighborhood}.

\medskip
\noindent\textbf{Cusp neighborhood.}
Let \((x_0,s_0)\) be a cusp point, and take \((x,s)=(x_0,s_0)\). Choose a small
neighborhood \(\fN_{(x_0, s_0)}\) of \((x_0,s_0)\). For
\((x',s')\in\fN_{(x_0, s_0)}\), if \((x',s')\in\fL\), then there are two
complex-conjugate critical points near \(z(x_0,s_0)\). In addition,
there is a tangent line from \((x',s')\) to the arctic boundary whose point of
tangency is close to \((x_0,s_0)\); this tangent line produces an additional
real critical point. If \((x',s')\in\fP\setminus\fL\), then there are three
real critical points, counted with multiplicity, arising from tangent lines
whose tangency points are close to \((x_0,s_0)\).

We choose a \(\fc\)-cusp chart with local paths as in
\Cref{s:critical_cusp}, and shrink \(\fN_{(x_0, s_0)}\) so that every
\((x',s')\in\fN_{(x_0, s_0)}\) is adapted to this chart and all the associated
critical points remain inside it. We call \(\fN_{(x_0, s_0)}\) a
\emph{cusp neighborhood}.

\medskip
\noindent\textbf{Arctic neighborhood.}
Let \((x_0,s_0)\in\fA\) be a regular arctic point, bounded away from tangency
points, cusp points, and cusp-turning points, and take \((x,s)=(x_0,s_0)\).
Choose a small neighborhood \(\fN_{(x_0, s_0)}\) of \((x_0,s_0)\) with the same
property. For \((x',s')\in\fN_{(x_0, s_0)}\), if \((x',s')\in\fL\), then there are
two complex-conjugate critical points near \(z(x_0,s_0)\). If
\((x',s')\in\fP\setminus\fL\), then there are two real critical points,
counted with multiplicity, near \(z(x_0,s_0)\), arising from tangencies near
\((x_0,s_0)\).
We choose a \(\fc\)-arctic chart with local paths as in
\Cref{s:critical_arctic}, and shrink \(\fN_{(x_0, s_0)}\) so that every
\((x',s')\in\fN_{(x_0, s_0)}\) is adapted to this chart and the above critical
points remain inside it.

In addition, there may be a further tangent line from \((x',s')\) to the
arctic boundary. In this case, we choose \(\fN_{(x_0, s_0)}\) small enough so that
this tangent line persists continuously for all
\((x',s')\in\fN_{(x_0, s_0)}\). This gives rise to an additional real critical
point, which remains bounded away from \(z(x_0,s_0)\). According to whether
this critical point is bounded away from cusp and tangency points, close to a
tangency point, or close to a cusp point, we associate to \((x_0,s_0)\) the
corresponding \(\fc\)-frozen chart of regular, tangency, or cusp type,
together with the local paths as in
\Cref{s:frozen_neighborhood,s:vertical_frozen_neighborhood,s:unit_slope_frozen_neighborhood,s:horizontal_frozen_neighborhood}.
For \(\fc\) sufficiently small, this frozen chart can be chosen disjoint from
the arctic chart. Shrinking \(\fN_{(x_0, s_0)}\) further if necessary, we may
assume that every point in \(\fN_{(x_0, s_0)}\) is adapted to both charts. We call
\(\fN_{(x_0, s_0)}\) an \emph{arctic neighborhood}.

\medskip

\noindent \textbf{Interior frozen neighborhood.}
Suppose $(x,s)\in\fT$ is bounded away from the arctic boundary and from all other curvilinear triangles. Choose a small neighborhood $\fN_{(x,s)}$ of $(x,s)$ with the same property.

For each $(x',s')\in\fN_{(x,s)}$, there are two or three tangent lines associated with $(x',s')$. We choose $\fN_{(x,s)}$ small enough so that these tangent line persists continuously for all
$(x',s')\in\fN_{(x,s)}$.

 Each such tangent line yields a critical point, and these critical points remain separated from one another. For each such tangent line, 
according to whether the critical point is bounded away from cusp and tangency points, close to a tangency point, or close to a cusp point, we associate to $(x,s)$ the corresponding $\fc$-frozen chart of regular, tangency, or cusp type, together with the local paths as in \Cref{s:frozen_neighborhood,s:vertical_frozen_neighborhood,s:unit_slope_frozen_neighborhood,s:horizontal_frozen_neighborhood}. Choosing $\fc$ sufficiently small, these charts may be taken pairwise disjoint. Shrinking $\fN_{(x,s)}$ if necessary, we ensure that every $(x',s')\in\fN_{(x,s)}$ is adapted to all of them. We call $\fN_{(x,s)}$ an \emph{interior frozen neighborhood}.

\medskip
\medskip
\noindent\textbf{Tangent, cusp-turning, and interface frozen neighborhoods.}
Now suppose that \((x,s)\) lies on a boundary shared by two adjacent
curvilinear triangles \(\fT_A\) and \(\fT_B\). We distinguish the following
subcases.

\begin{figure}
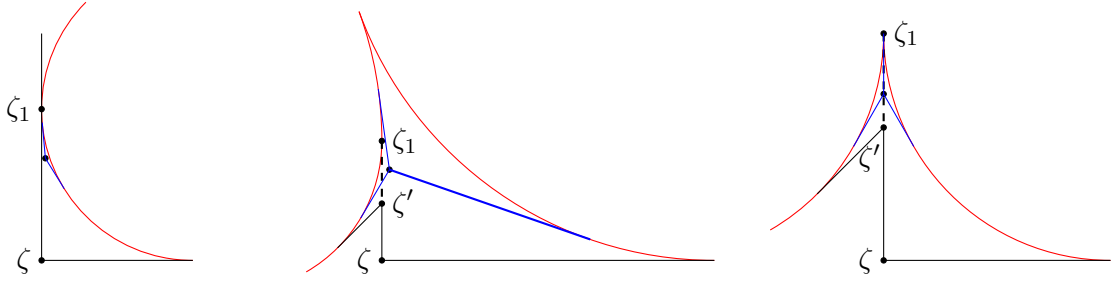

  \centering
  \begin{subfigure}{0.25\textwidth}
    \centering

  \end{subfigure}
\caption{Tangent lines from a point near a common boundary of two adjacent curvilinear triangles.}
  \label{f:adjacent_curvilinear_triangle2}
\end{figure}

\begin{enumerate}
\item \textbf{Shared tangency point.}
Suppose \(\fT_A\) and \(\fT_B\) share only a tangency point
\(\zeta_1=(x_0,s_0)\) lying on a side of \(\fP\); see the left panel of
\Cref{f:adjacent_curvilinear_triangle2}.

Take \((x,s)=(x_0,s_0)\), and choose a small neighborhood
\(\fN_{(x_0, s_0)}\) of \((x_0,s_0)\). For
\((x',s')\in\fN_{(x_0, s_0)}\), if \((x',s')\in\fL\), then there are two
complex-conjugate critical points. If \((x',s')\in\fP\setminus\fL\), then
there are two tangent lines from \((x',s')\) to the portion of the arctic
boundary contained in \(\fT_A\cup\fT_B\); see the left panel of
\Cref{f:adjacent_curvilinear_triangle2}.
 Both tangent lines are close to the
corresponding side of \(\fP\), and their tangency points are close to
\((x_0,s_0)\). They give two real critical points, counted with multiplicity.
In either case, the two critical points are close to \(z(x_0,s_0)\).

We choose a \(\fc\)-tangent chart with local descent/ascent paths chosen
relative to either \(\fT_A\) or \(\fT_B\), as in
\Cref{s:vertical_tangent,s:unit_slope_tangent,s:horizontal_tangent}. Shrinking
\(\fN_{(x_0, s_0)}\) if necessary, we may assume that every point in
\(\fN_{(x_0, s_0)}\) is adapted to this chart. We call \(\fN_{(x_0, s_0)}\) a
\emph{tangent neighborhood}.

\item \textbf{Shared segment from a tangency point to a vertex.}
\label{i:share_caseb}
Suppose \(\fT_A\) and \(\fT_B\) share a segment \([\zeta_1,\zeta']\), where
\(\zeta_1=(x_0,s_0)\) is a tangency point and \(\zeta'\) is a vertex of
\(\fP\); see the middle panel of \Cref{f:adjacent_curvilinear_triangle2}.

First take \((x,s)=(x_0,s_0)\), and choose a small neighborhood
\(\fN_{(x_0, s_0)}\) of \((x_0,s_0)\). For every
\((x',s')\in\fN_{(x_0, s_0)}\), there are two critical points close to
\(z(x_0,s_0)\), as in the previous subcase. We therefore associate to
\((x_0,s_0)\) the same \(\fc\)-tangent chart as above.

In addition, there exsits one further tangent line from \((x',s')\) to the
arctic boundary; see the thick blue tangent line in the middle panel of
\Cref{f:adjacent_curvilinear_triangle2}.
 The corresponding critical point remains bounded away from
\(z(x_0,s_0)\). We therefore also associate the corresponding
\(\fc\)-frozen chart of regular, tangency, or cusp type, together with
local paths as in \Cref{s:frozen_neighborhood,s:vertical_frozen_neighborhood,s:unit_slope_frozen_neighborhood,s:horizontal_frozen_neighborhood}.
For \(\fc\) sufficiently small, this frozen chart can be chosen disjoint from
the tangent chart. Shrinking \(\fN_{(x_0, s_0)}\) further if necessary, we may
assume that every point in \(\fN_{(x_0, s_0)}\) is adapted to both charts. We call
\(\fN_{(x_0, s_0)}\) a \emph{tangent neighborhood}.

Now suppose that \((x,s)\in(\zeta_1,\zeta']\) is bounded away from
\(\zeta_1\). Choose a small neighborhood \(\fN_{(x,s)}\) of \((x,s)\) that is
also bounded away from \(\zeta_1\). Then, for every
\((x',s')\in\fN_{(x,s)}\), there are three tangent lines from \((x',s')\) to the
arctic boundary. One of them lies in a small neighborhood of the shared segment
\([\zeta_1,\zeta']\), has tangency point close to \((x_0,s_0)\), and gives a
critical point close to \(z(x_0,s_0)\). We associate to \((x,s)\) the
corresponding tangency-type \(\fc\)-frozen chart, together with
local descent/ascent paths chosen relative to either \(\fT_A\) or
\(\fT_B\), as in \Cref{s:vertical_frozen_neighborhood,s:unit_slope_frozen_neighborhood,s:horizontal_frozen_neighborhood}.

The other two tangent lines lie in \(\fT_A\) and \(\fT_B\), respectively, and
their corresponding critical points are bounded away from \(z(x_0,s_0)\). We
associate to \((x,s)\) the corresponding \(\fc\)-frozen charts of regular,
tangency, or cusp type, together with local paths as in
\Cref{s:frozen_neighborhood,s:vertical_frozen_neighborhood,s:unit_slope_frozen_neighborhood,s:horizontal_frozen_neighborhood}. For \(\fc\)
sufficiently small, these charts can be chosen pairwise disjoint. Shrinking
\(\fN_{(x,s)}\) further if necessary, we may assume that every point in
\(\fN_{(x,s)}\) is adapted to all of them. We call \(\fN_{(x,s)}\) an
\emph{interface frozen neighborhood}.

\item \textbf{Shared segment from a cusp-turning point to a vertex.}
Suppose \(\fT_A\) and \(\fT_B\) share a segment \([\zeta_1,\zeta']\), where
\(\zeta_1=(x_0,s_0)\) is a cusp-turning point and \(\zeta'\) is a vertex of
\(\fP\); see the right panel of \Cref{f:adjacent_curvilinear_triangle2}.

First take \((x,s)=(x_0,s_0)\). Exactly as in the cusp case, we associate to
\((x_0,s_0)\) the corresponding \(\fc\)-cusp-turning chart, together with
local paths chosen relative to either \(\fT_A\) or \(\fT_B\), as in
\Cref{s:vertical_tangent,s:unit_slope_tangent,s:horizontal_tangent}. A
sufficiently small neighborhood \(\fN_{(x_0, s_0)}\) can be chosen so that every
point in it is adapted to this chart. We call \(\fN_{(x_0, s_0)}\) a
\emph{cusp-turning neighborhood}.

Now suppose that \((x,s)\in(\zeta_1,\zeta']\) is bounded away from
\(\zeta_1\). Choose a small neighborhood \(\fN_{(x,s)}\) of \((x,s)\) that is
also bounded away from \(\zeta_1\). Then, exactly as in the case of a shared
segment from a tangency point to a vertex, see \Cref{i:share_caseb}, we
associate to \((x,s)\) three pairwise disjoint \(\fc\)-frozen charts, so that
every point in \(\fN_{(x,s)}\) is adapted to all of them. We call
\(\fN_{(x,s)}\) an \emph{interface frozen neighborhood}.
\end{enumerate}

From the discussion above, the polygon $\fP$ can be covered by liquid, ramification, cusp, arctic, interior frozen, tangency, cusp-turning, and interface frozen neighborhoods. This proves \Cref{p:construct_neighborhood1}.
\end{proof}

\begin{remark}\label{r:shrink_c}
We remark that the collection of $\fc$-charts associated with $(y,t)$ admits a
natural refinement structure as $\fc$ decreases. If we shrink $\fc$ continuously, then
either the corresponding (shrunk) $\fc$-charts still satisfy the conclusion of
\Cref{p:construct_neighborhood1}, or else the conclusion fails and we must reconstruct the
collection.

In the latter case, the reconstruction amounts to a finite sequence of local refinements:
\begin{enumerate}
\item A cusp chart may split into an arctic chart together with one frozen chart, a pair of complex-conjugate liquid charts, three regular frozen charts, or one cusp-frozen chart together with two regular frozen charts.

\item A cusp-turning chart may be replaced by an arctic chart, or may split into a pair of complex-conjugate liquid charts, two regular frozen charts, or one tangent frozen chart together with two regular frozen charts.

\item A tangent chart may be replaced by an arctic chart, or may split into two smaller liquid charts or two frozen charts.

\item An arctic chart may split into two smaller frozen charts, or a pair of complex-conjugate liquid charts.

\item A ramification chart may be replaced by a smaller liquid chart.
\item A tangent frozen or cusp-frozen chart may be replaced by a smaller regular frozen chart.
\end{enumerate}

Such reconstructions may increase the number of charts and/or shift their centers.
By contrast, regular frozen and liquid  charts are stable under shrinking $\fc$ (up to taking
smaller charts with the same center).
\end{remark}

\begin{proposition}\label{p:replace_center}
Adopt the assumptions and notation of \Cref{p:construct_neighborhood1}.
After possibly decreasing \(\fc\) and shrinking \(\fN_{(x,s)}\), the
centers of the charts in \Cref{p:construct_neighborhood1} may be chosen
canonically as follows: an \(\fc\)-arctic chart may be centered at the
midpoint, in the local \(w\)-coordinate, of the two critical points of
\(S(\cdot;x,s)\) contained in the chart; an \(\fc\)-liquid chart or
\(\fc\)-regular frozen chart may be centered at the corresponding critical
point of \(S(\cdot;x,s)\).
\end{proposition}

The
liquid and regular frozen cases follow from the continuity of the corresponding nondegenerate
critical point, while the arctic case follows from the following \Cref{c:lclose}. This proposition allows us to shift the center within a
fixed local chart, while keeping all relevant critical points inside the chart. 

\begin{proposition}\label{c:lclose}
Let \(\fU\subset\cC\) be a chart of any type constructed in
\Cref{s:critical_point}, centered at \(w_0\in\cC(\bR)\). Then, for every
sufficiently small \(\fc>0\), there exists \(\fb=\fb(\fc)>0\), depending
also on \(w_0\), such that the following statements hold.
\begin{enumerate}
\item
If \(w_0\neq\infty\), then \(\cC\) can be parametrized on \(\fU\) as a
graph \((f(w),w)\). For any \((x,s)\in\fP\), if \(S(\cdot;x,s)\) has a
critical point \(w_c\in\fU\) satisfying
$
|w_c-w_0|<\fb,
$
then there exists an \(\fc\)-chart centered at \(w_0\) such that a
sufficiently small neighborhood of \((x,s)\) is adapted to this chart.

\item
If \(w_0=\infty\), corresponding to a horizontal tangency point or a
cusp-turning point \((x_0,s_0)\), recall the \(\wt w\)-chart
\begin{align}
\wt\fU
:=
\left\{\frac{1}{x_0-w}:w\in\fU\right\},
\end{align}
centered at \(0\). For any \((x,s)\in\fP\), if
\(\wt S(\cdot;x,s)\) has a critical point
\(\wt w_c\in\wt\fU\) satisfying
$
|\wt w_c|<\fb,
$
then there exists an \(\fc\)-chart centered at \(0\) such that a
sufficiently small neighborhood of \((x,s)\) is adapted to this chart.
\end{enumerate}
\end{proposition}

\begin{proof}[Proof of \Cref{p:replace_center}]
Start with the collection of charts given by
\Cref{p:construct_neighborhood1}. Since this collection is finite, it
suffices to replace its charts one at a time.

Suppose first that \((x,s)\) is adapted to an arctic chart, and let
\(w_1\) and \(w_2\) be the two critical points of \(S(\cdot;x,s)\)
contained in this chart. Set
\[
w_0:=\frac{w_1+w_2}{2}\in \bR,
\]
where the midpoint is taken in the local \(w\)-coordinate. If \(w_1\) and
\(w_2\) are sufficiently close to \(w_0\), then, by
\Cref{c:lclose}, after shrinking
\(\fN_{(x,s)}\) if necessary, every point in \(\fN_{(x,s)}\) is adapted
to a \(\fc\)-arctic chart centered at \(w_0\).

If \(w_1\) and \(w_2\) are not sufficiently close to \(w_0\), then they
are separated and nondegenerate. In this case, we replace the arctic chart
by two sufficiently small disjoint liquid charts or regular frozen charts,
centered at \(w_1\) and \(w_2\), respectively.

Suppose next that \((x,s)\) is adapted to a liquid chart or a regular
frozen chart. Let \(w_c\) be the unique critical point of
\(S(\cdot;x,s)\) contained in this chart. Since \(w_c\) is nondegenerate,
the corresponding critical point varies continuously with \((x,s)\). Therefore, after shrinking
\(\fN_{(x,s)}\) if necessary, every point in \(\fN_{(x,s)}\) is adapted
to a liquid chart or regular frozen chart, respectively, centered at
\(w_c\). This completes the proof.
\end{proof}

\begin{proof}[Proof of \Cref{c:lclose}]
We prove the statement when \(w_0\neq\infty\). The case \(w_0=\infty\)
follows by the same argument in the \(\wt w\)-coordinate, using the
corresponding critical point relation, so we omit it.

\medskip
\noindent\emph{Step 1: a compactness argument along the reference tangent
line.}
The point \(w_0\in\cC(\bR)\) corresponds to a point
\((x_0,s_0)\in\fA\). Consider the tangent line to the arctic boundary at
\((x_0,s_0)\),
\[
L_0:=\{(x',s'):\ w_0=x'-s'\chi(w_0)\},
\]
and set
$
\ell_0:=L_0\cap\fP.
$
Since \(\fP\) is compact, \(\ell_0\) is compact.

Fix \(\fc>0\) sufficiently small. For each \(\xi\in\ell_0\), the
constructions in \Cref{s:critical_point} provide an \(\fc\)-chart centered
at \(w_0\) and a number \(\delta(\xi)>0\) such that
\(B_{\delta(\xi)}(\xi)\cap\fP\) is adapted to this chart. Thus,
\[
\ell_0
\subset
\bigcup_{\xi\in\ell_0}B_{\delta(\xi)}(\xi).
\]
By compactness, there exist
\(\xi_1,\dots,\xi_k\in\ell_0\) such that, writing
\(\delta_i:=\delta(\xi_i)\),
\[
\ell_0
\subset
\bigcup_{i=1}^k B_{\delta_i}(\xi_i).
\]
Since the right-hand side is an open neighborhood of the compact set
\(\ell_0\), there exists \(\delta>0\) such that
\begin{equation}\label{e:inclusion_lclose}
B_\delta(\ell_0)
\subset
\bigcup_{i=1}^k B_{\delta_i}(\xi_i).
\end{equation}

\medskip
\noindent\emph{Step 2: closeness of \(w_c\) to \(w_0\) forces
\((x,s)\) to be close to \(\ell_0\).}
We claim that there exists \(\fb>0\) sufficiently small such that, if
\((x,s)\in\fP\) admits a critical point \(w_c\in\fU\) satisfying
$
|w_c-w_0|<\fb,
$
then
\begin{equation}\label{e:disL_lclose}
\dist\bigl((x,s),\ell_0\bigr)<\delta.
\end{equation}

Assuming \eqref{e:disL_lclose}, we have
\((x,s)\in B_\delta(\ell_0)\), and
\eqref{e:inclusion_lclose} implies that
\((x,s)\in B_{\delta_i}(\xi_i)\) for some \(i\in\{1,\dots,k\}\).
By construction, a sufficiently small neighborhood of \((x,s)\) is
adapted to an \(\fc\)-chart centered at \(w_0\). This proves the
proposition.

It remains to prove \eqref{e:disL_lclose}. We distinguish two cases.

\begin{enumerate}
\item
\emph{Suppose that \(w_c\notin\bR\).}
By complex conjugation, it suffices to consider \(w_c\in\bC_+\). Then
\((x,s)\in\fL\). By \eqref{e:emb}, the map
$
(x,s)\in\overline{\fL}\longmapsto z(x,s)
$
is a homeomorphism onto its image. Since
$
w_c=z(x,s),
$ and $
w_0=z(x_0,s_0),
$
the condition \(|w_c-w_0|<\fb\) implies
\[
\|(x,s)-(x_0,s_0)\|_2=\oo_{\fb}(1).
\]
Because \((x_0,s_0)\in\ell_0\), choosing \(\fb\) sufficiently small gives
\eqref{e:disL_lclose}.

\item
\emph{Suppose that \(w_c\in\bR\).}
Since \(w_c\) is a critical point of \(S(\,\cdot\,;x,s)\), it satisfies
\begin{equation}\label{e:crit_relation_lclose}
w_c=x-s\chi(w_c).
\end{equation}
By continuity of \(\chi\) on \(\fU\),
\[
|\chi(w_c)-\chi(w_0)|=\oo_{\fb}(1)
\]
whenever \(|w_c-w_0|<\fb\). 
Choose \(M>0\) such that \(|s'|\leq M\) for every
\((x',s')\in\fP\), and notice that $(w_0+s\chi(w_0), s)\in L_0$. Consequently,
\begin{align*}
\dist\bigl((x,s),L_0\bigr)
&\leq
\bigl|x-(w_0+s\chi(w_0))\bigr| =
\bigl|(w_c+s\chi(w_c))-(w_0+s\chi(w_0))\bigr| \notag\\
&\leq
|w_c-w_0|
+
M|\chi(w_c)-\chi(w_0)|
=
\oo_{\fb}(1).
\end{align*}

Since \(\fP\) is compact and
\(\ell_0=\fP\cap L_0\), closeness to \(L_0\) implies, uniformly for points
in \(\fP\), closeness to \(\ell_0\). Indeed, otherwise there would exist a
sequence \((x_j,s_j)\in\fP\) such that
\[
\dist\bigl((x_j,s_j),L_0\bigr)\longrightarrow0,
\qquad
\dist\bigl((x_j,s_j),\ell_0\bigr)\geq\delta.
\]
After passing to a subsequence,
\((x_j,s_j)\to(x_\infty,s_\infty)\in\fP\). The first relation implies
\((x_\infty,s_\infty)\in L_0\), and hence
\((x_\infty,s_\infty)\in\ell_0\), contradicting the second relation.
Therefore, \eqref{e:disL_lclose} holds for \(\fb\) sufficiently small.
\end{enumerate}

In either case, choosing \(\fb=\fb(\fc)>0\) sufficiently small yields
\eqref{e:disL_lclose} and completes the proof.
\end{proof}

\subsection{Assigning charts for a pair of points in $\fP$}
When we later construct the approximation of the inverse Kasteleyn matrix $K^{-1}((x,s),(y,t))$, we fix $(y,t)\in \fP$ and let $(x,s)\in \fP$ vary. Accordingly, we must associate charts in $\cC$ to both $(y,t)$ and $(x,s)$. The next proposition asserts that we can first fix the collection of charts associated with $(y,t)$. And then, for each $(x,s)\in \fP$, we can associate it with charts such that any two charts (from the combined collections) are either disjoint, or concentric (i.e. they have the same center). 

\begin{proposition}\label{p:construct_neighborhood2}
Fix $(y,t)\in \fP$. Then there
exists \(\fc_0>0\), sufficiently small, such that \((y,t)\) can be
associated with a finite collection of pairwise disjoint \(\fc_0\)-charts
satisfying the conclusions of \Cref{p:construct_neighborhood1}, with the
canonical centers specified in \Cref{p:replace_center}.

Moreover, for every \((x,s)\in\fP\), there exists \(\fc'>0\) such that the
following holds. After possibly shrinking the \(\fc_0\)-charts associated
with \((y,t)\), while keeping their centers fixed, the conclusions of
\Cref{p:construct_neighborhood1} continue to hold. In addition, \((x,s)\)
can be associated with a finite collection of pairwise disjoint
\(\fc'\)-charts satisfying the conclusions of
\Cref{p:construct_neighborhood1}.

Finally, these two collections can be chosen compatibly in the following
sense: any chart from the \((y,t)\)-collection and any chart from the
\((x,s)\)-collection are either disjoint or concentric, that is, they have
the same center.
\end{proposition}

%{\color{red}add a figure. In addition, the concentric overlaps between the two collections are either none, exactly one,
%or all; in the third case (all), $(x,s)$ is close to $(y,t)$.}

 Then there exist $\fc>0$ sufficiently small such that $(y,t)$ can be associated with a finite collection
of $\fc$-charts satisfying \Cref{p:construct_neighborhood1} and \Cref{p:replace_center}. The centers of
these charts belong to a finite list of possibilities, namely:
\begin{enumerate}
\item {Cusp/cusp-turning/tangent charts:} the center $w_0$ is the
corresponding cusp, cusp-turning, or tangency point in $\cC$.

\item {Arctic charts:} the center $w_0\in\bR$ is the midpoint of the two
critical points associated with $(y,t)$ contained in the chart; moreover, $w_0$ is
bounded away from the cusp, cusp-turning, and tangency points in $\cC$.

\item {Tangent frozen/cusp frozen charts:} the center $w_0$ is the
corresponding tangency or cusp point in $\cC$.

\item {Regular frozen/liquid/ramification charts:} the center $w_0$ is a critical
point associated with $(y,t)$.
\end{enumerate}

\begin{lemma}[Nested charts and a gap scale]\label{c:key_prop}
Fix \((y,t)\in\fP\) bounded away from the ramification points. Then there
exist constants
\[
0<\fb_0<\frac{\fc_0}{4}
\]
such that the following statements hold.

\begin{enumerate}
\item\label{i:exist}
\emph{Nested chart assignment.}
For every \(\fc\in[\fb_0,\fc_0]\), we can associate with \((y,t)\) a
finite collection of \(\fc\)-charts in \(\cC\) satisfying the conclusions
of \Cref{p:construct_neighborhood1}. Moreover, this assignment can be
chosen to be nested in \(\fc\): as \(\fc\) decreases, the charts shrink
while their centers remain fixed. Let \(P_0\) denote the finite set of
centers of these charts.

Furthermore, for each \(\fc\in[\fb_0,\fc_0]\), there exists
\[
0<\fb(\fc)<\fc
\]
such that, for every \(w_0\in P_0\cap\cC(\bR)\), the stability statement
of \Cref{c:lclose} holds with this common choice of \(\fb(\fc)\). 
\item\label{i:gap}
\emph{Existence of a gap scale.}
Fix \((x,s)\in\fP\), and let \(P(x,s)\) be the finite set of critical
points \(w_c\in\cC\) associated with \((x,s)\) by
\Cref{p:associate_critical_points}. Define the finite set of distances
\(D(x,s)\) as follows.
For each \(w_0\in P_0\):
\begin{enumerate}
\item
if \(w_0\neq\infty\), let
\(\fU(w_0)\subset\cC\) denote the \(\fc_0\)-chart centered at
\(w_0\) from part~\ref{i:exist}, and include in \(D(x,s)\) the distances
\begin{equation}\label{e:defD}
\Bigl\{
|w_c-w_0|:
w_c\in P(x,s)\cap\fU(w_0)
\Bigr\};
\end{equation}

\item
if \(w_0=\infty\), let
\(\wt\fU(0)\) denote the corresponding \(\fc_0\)-chart centered
at \(0\) in the \(\wt w\)-coordinate, and write \(\wt w_c\) for the
\(\wt w\)-coordinate of \(w_c\). Include in \(D(x,s)\) the distances
\begin{equation}\label{e:defD1}
\Bigl\{
|\wt w_c|:
w_c\in P(x,s),\ 
\wt w_c\in\wt\fU(0)
\Bigr\}.
\end{equation}
\end{enumerate}

Then there exists \(\fc'\in[\fb_0,\fc_0/4]\) such that
\begin{equation}\label{e:seperateD}
D(x,s)\cap\bigl(\fb(\fc'),4\fc'\bigr)=\emptyset.
\end{equation}
\end{enumerate}
\end{lemma}

\begin{proof}[Proof of \Cref{c:key_prop}]
We prove \eqref{i:exist} and \eqref{i:gap} in two steps.

\medskip
\noindent\emph{Step 1: Pigeonhole gap for distances.}
Fix $(x,s)\in\fP$. By the construction in \Cref{p:associate_critical_points}, the number of relevant points is uniformly bounded:
\[
|P_0|,| P(x,s)|\le 3=:M,
\]
and hence the set of pairwise distances $D(x,s)$ has cardinality at most $M^2$.

Given any sufficiently small $\fc>0$, define a decreasing sequence of scales by
\[
\fc_{1}:=\fc,\qquad \fc_{i+1}:=\fb(\fc_i/4)<\fc_i/4,\qquad i=1,2,\dots,2M^2.
\]
We can view $\fc_{2M^2}$ as a function of $\fc$:
\begin{equation}\label{e:intro_g}
\fc_{2M^2}=g(\fc).
\end{equation}

The $2M^2-1$ intervals $(\fc_{i+1},\fc_i)$ for $i=1,2,\dots,2M^2-1$ are disjoint and contained in $[g(\fc),\fc]$.
Since $D(x,s)$ has at most $M^2<2M^2-1$ elements, there exists an index $i\in\{1,2,\dots,2M^2-1\}$ such that
\[
D(x,s)\cap(\fc_{i+1},\fc_i)=\emptyset.
\]
Set $\fc':=\fc_i/4$. Then $\fc'\in[g(\fc),\fc]$, and because $\fc_{i+1}=\fb(\fc_i/4)=\fb(\fc')$, the emptiness above implies
\begin{equation}\label{e:Dgap}
D(x,s)\cap(\fb(\fc'),\,4\fc')=\emptyset,
\end{equation}
which is exactly the type of statement as in \Cref{i:gap}.

\medskip
\noindent\emph{Step 2: Existence of a nested interval and fixed centers.}
Start from any sufficiently small $\fc^{(0)}>0$. Define inductively a decreasing sequence
\begin{equation}\label{e:iter}
\fc^{(n+1)} := g(\fc^{(n)}),\qquad n\ge 0,
\end{equation}
where $g(\cdot)$ is the function produced in \eqref{e:intro_g} (so that $0<g(\fc)<\fc$ for $\fc$ small).
Consider the interval $[\fc^{(n+1)},\fc^{(n)}]$. By \Cref{r:shrink_c}, as $\fc$ decreases within this interval,
either the collection of $\fc$-charts associated with $(y,t)$ can be chosen nested (shrinking with fixed
centers) throughout the whole interval, or else a reconstruction event occurs.

Each reconstruction event is one of the local refinements listed in \Cref{r:shrink_c}. In particular, such an
event either increases the number of charts and/or changes their centers, but it does so in a controlled,
finite manner: cusp/tangent/arctic/ramification types can only refine a finite number of times before all
charts become stable types (regular frozen or bulk), which, by \Cref{r:shrink_c}, do not undergo further
reconstructions under shrinking.

Consequently, after finitely many steps, there exists an index $n_\ast$ such that no reconstruction occurs on the
entire interval $[\fc^{(n_\ast+1)},\fc^{(n_\ast)}]$. We then set
\[
\fb_0:=\fc^{(n_\ast+1)},\qquad \fc_0:=\fc^{(n_\ast)}.
\]
On this interval, the $\fc$-charts associated with $(y,t)$ may be chosen nested in $\fc$ (with fixed
centers). Denote by $P_0$ the resulting (finite) set of centers. This proves \Cref{i:exist}. The second statement
\Cref{i:gap} follows from \eqref{e:Dgap}. This completes the proof of \Cref{c:key_prop}.
\end{proof}

\begin{proof}[Proof of \Cref{p:construct_neighborhood2}]

If \((y,t)\) is sufficiently close to a ramification point
\((y_0,t_0)\), the statement is straightforward. It is associated with a
complex-conjugate pair of ramification charts. If \((x,s)\) is also
sufficiently close to \((y_0,t_0)\), then we associate it with the same
complex-conjugate pair of ramification charts. If \((x,s)\) is sufficiently
far from \((y_0,t_0)\), then, after possibly shrinking the ramification
charts associated with \((y,t)\), we can associate \((x,s)\) with charts
disjoint from them.

In the following we assume that $(y,t)\in\fP$ is bounded away from ramification points.
Fix \((y,t)\in\fP\), and choose
$
0<\fb_0<\fc_0/4
$
and the nested chart assignment from \Cref{c:key_prop}. At the scale
\(\fc_0\), this gives the first assertion of the proposition, after choosing
the canonical centers as in \Cref{p:replace_center}.

Now fix \((x,s)\in\fP\), and let
\[
\fc'\in[\fb_0,\fc_0/4]
\]
be the scale provided by part~\ref{i:gap} of \Cref{c:key_prop}. Since
\[
\fc',\,4\fc'\in[\fb_0,\fc_0],
\]
the nested assignment provides both the \(\fc'\)-charts and the
\(4\fc'\)-charts associated with \((y,t)\), with the same set of centers
\(P_0\). For each \(w_0\in P_0\), denote these charts by
\[
\fU_{\fc'}(w_0)
\qquad\text{and}\qquad
\fU_{4\fc'}(w_0),
\]
respectively. When \(w_0=\infty\), these charts and all distances below are
understood in the corresponding \(\wt w\)-coordinate.

Using \Cref{p:construct_neighborhood1}, and shrinking the resulting charts
if necessary, construct a finite collection of pairwise disjoint
\(\fc'\)-charts associated with \((x,s)\). We modify this collection as
follows.

Let \(\fV\) be a chart in the \((x,s)\)-collection that intersects
\(\fU_{\fc'}(w_0)\) for some \(w_0\in P_0\). Since both charts have scale
\(\fc'\), every critical point \(w_c\) contained in \(\fV\) lies in
\(\fU_{4\fc'}(w_0)\). By the nestedness of the charts,
\[
\fU_{4\fc'}(w_0)\subset\fU_{\fc_0}(w_0),
\]
so the corresponding distance from \(w_c\) to \(w_0\) belongs to
\(D(x,s)\). Moreover, this distance is strictly smaller than \(4\fc'\).
Therefore, by \eqref{e:seperateD}, it is at most \(\fb(\fc')\).

The number \(\fb(\fc')\) in \Cref{c:key_prop} is the stability radius furnished by \Cref{c:lclose}. Hence the
stability statement in part~\ref{i:exist} of \Cref{c:key_prop} implies that
a sufficiently small neighborhood of \((x,s)\) is adapted to the
\(\fc'\)-chart \(\fU_{\fc'}(w_0)\). We may therefore replace all charts in
the \((x,s)\)-collection that intersect \(\fU_{\fc'}(w_0)\) by the single
chart \(\fU_{\fc'}(w_0)\).

We perform this replacement for every \(w_0\in P_0\). The replacement is
unambiguous: the \(4\fc'\)-charts
\[
\bigl\{\fU_{4\fc'}(w_0):w_0\in P_0\bigr\}
\]
are pairwise disjoint, so a chart in the \((x,s)\)-collection cannot
intersect \(\fc'\)-charts centered at two distinct points of \(P_0\).
Moreover, the replacement charts are pairwise disjoint, and every chart in
the \((x,s)\)-collection that is not replaced is disjoint from all charts
in the \((y,t)\)-collection.

After shrinking the neighborhood of \((x,s)\) so that it is adapted to all
of the finitely many resulting charts, we obtain a pairwise disjoint
collection of \(\fc'\)-charts satisfying the conclusions of
\Cref{p:construct_neighborhood1}. By construction, every chart in this
collection is either disjoint from every chart in the \((y,t)\)-collection
or concentric with one of them.

Finally, the \(\fc'\)-charts associated with \((y,t)\) are obtained from
the original \(\fc_0\)-charts by shrinking while keeping their centers
fixed, and part~\ref{i:exist} of \Cref{c:key_prop} ensures that the
conclusions of \Cref{p:construct_neighborhood1} continue to hold. This
completes the proof.
\end{proof}

\subsection{An open cover of $\fP$}\label{s:finite_cover}

Fix \((y,t)\in\fP\). We use
\Cref{p:construct_neighborhood2} to construct an open cover of \(\fP\).

For each \((x,s)\in\fP\), there exist \(\fc_0>0\) and
\(\fc'=\fc'(x,s)>0\), together with a collection of \(\fc_0\)-charts
associated with \((y,t)\) and a collection of \(\fc'\)-charts associated
with \((x,s)\), such that:
\begin{enumerate}
\item the charts within each collection are pairwise disjoint;
\item each chart in the \((y,t)\)-collection and each chart in the
\((x,s)\)-collection are either disjoint or concentric.
\end{enumerate}

Moreover, as $(x,s)$ varies, the  chart collections associated with the point $(y,t)$ have the same centers and differ
only by shrinking. Thus, after decreasing \(\fc_0\) if necessary, we may
choose a single collection of \(\fc_0\)-charts associated with \((y,t)\).  We
then choose a sufficiently small open neighborhood \(\fN_{(y,t)}\) of
\((y,t)\), 
\begin{align}\label{e:ytchart}
(y,t)\in \fN_{(y,t)}
\end{align}
 of one of the types appearing in
\Cref{p:construct_neighborhood1}---liquid, ramification, cusp, arctic,
interior frozen, tangency, cusp-turning, or interface frozen---such that
every point in \(\fN_{(y,t)}\) is adapted to all charts in this fixed
\((y,t)\)-collection.

For each \((x,s)\in\fP\), let \(\fN_{(x,s)}\) be a small open neighborhood
of \((x,s)\), again of one of the types appearing in
\Cref{p:construct_neighborhood1}, such that every
\((x',s')\in\fN_{(x,s)}\) is adapted to all of the \(\fc'\)-charts in the
\((x,s)\)-collection and all critical points associated with \((x',s')\)
are contained in the union of these charts. Then
\(\{\fN_{(x,s)}\}_{(x,s)\in\fP}\) is an open cover of \(\fP\). Since
\(\fP\) is compact, we can extract a finite subcover. We denote it by
\begin{equation}\label{e:finite_cover}
\fP\subset\bigcup_{\alpha}\fN_{\alpha}.
\end{equation}

\section{Ansatz for the General Kernels}\label{s:general_kernel_ansatz}

We recall $\fN_{(y,t)}$ from \eqref{e:ytchart},  the open cover of $\fP$ from \eqref{e:finite_cover}, and consider
\begin{align}\label{e:ytxs}
(y,t)\in \fN_{(y,t)}, \quad (x,s)\in \fN_\alpha.
\end{align}
Then both
$(x,s)$ and $(y,t)$ are adapted to finite collections of charts. According to the type of neighborhood, the chart
collection associated with $(x,s)$ (and similarly with $(y,t)$) is one of the following:
\begin{enumerate}
\item two liquid charts (or two ramification charts), which are complex conjugates of each other;
\item a single cusp chart or cusp-turning chart;
\item an arctic chart, with or without an additional frozen chart;
\item a tangency chart, with or without an additional frozen chart;
\item two or three frozen charts, of regular, tangency, or cusp type.
\end{enumerate}

For each chart associated with $(x,s)$ and $(y,t)$, we denote its center by \(w_0\) and \(z_0\), respectively, and
recall the local descent and ascent contours \(\sfC^{\rm d}(w_0)\) and \(\sfC^{\rm a}(z_0)\) from \Cref{s:critical_bulk,s:critical_point}. Any two charts from the combined
collections are either disjoint or concentric (i.e., they have the same center).

In this section we construct ansatz for the inverse Kasteleyn matrix $K^{-1}((x,s),(y,t))$,
\begin{equation}\label{e:def_Aalpha}
A_\alpha((x,s),(y,t))
:=
J^{(1)}((x,s),(y,t))+J^{(2)}((x,s),(y,t)).
\end{equation}
where $J^{(1)}((x,s),(y,t))$ is a single-contour integral (which may vanish), and $J^{(2)}((x,s),(y,t))$ is a sum of double-contour integrals.

\subsection{Double-contour integral}

Let \((x,s),(y,t)\in\fP\) be as in \eqref{e:ytxs}. In this section, we
specify the double-contour integral \(J^{(2)}((x,s),(y,t))\) in
\eqref{e:def_Aalpha}. We recall the functions \(P_{ns}(nw,nx)\) and
\(Q_{nt}(nz,ny)\) from \eqref{e:defPQ}, together with the factors
\(I\), \(I_i\), and \(I_\ft\) from \Cref{s:factor}.

There are two cases, depending on the types of the neighborhoods
\(\fN_\alpha\) and \(\fN_{(y,t)}\):
\begin{enumerate}
\item
\emph{Neither \(\fN_\alpha\) nor \(\fN_{(y,t)}\) is a ramification
neighborhood.}
For each relevant pair of charts associated with \((x,s)\) and \((y,t)\),
centered at \(w_0\) and \(z_0\), respectively, we include in \(J^{(2)}\) a
double-contour integral of the form
\begin{equation}\label{e:all_term}
\frac{n}{(2\pi\ri)^2}
\int_{\sfC^{\rm a}}\!\!\int_{\sfC^{\rm d}}
P_{ns}(nw,nx)\,Q_{nt}(nz,ny)\,
\frac{I_+(w)}{I_-(z)}\,
\frac{\sqrt{\phi'(w)}\sqrt{\phi'(z)}}{\phi(w)-\phi(z)}\,
\rd w\,\rd z.
\end{equation}
If \(w_0\in[b_i,a_i]\) for some \(1\leq i\leq d\), we set
\(I_+(w)=I_i(w)\); otherwise, we set \(I_+(w)=I(w)\). Similarly, if
\(z_0\in[b_i,a_i]\) for some \(1\leq i\leq d\), we set
\(I_-(z)=I_i(z)\); otherwise, we set \(I_-(z)=I(z)\). We take
$
\sfC^{\rm d}=\sfC^{\rm d}(w_0)$,
and $
\sfC^{\rm a}=\sfC^{\rm a}(z_0),
$
as constructed in \Cref{s:critical_bulk,s:critical_point}.

\item
\emph{At least one of \(\fN_\alpha\) and \(\fN_{(y,t)}\) is a
ramification neighborhood.}
For each relevant pair of charts associated with \((x,s)\) and \((y,t)\),
centered at \(w_0\) and \(z_0\), respectively, we include in \(J^{(2)}\)
the double-contour integral obtained from \eqref{e:all_term} by making the
following replacements. If the chart centered at \(w_0\) is a ramification
chart, we replace
\[
P_{ns}(nw,nx), I_+(w), \phi(w), \phi'(w)\quad
\text{by}\quad 
P_{n(s+\ft)}(nw,nx),I_\ft(w),
\phi_\ft(w), \phi_\ft'(w),
\]
respectively. Similarly, if the chart centered at \(z_0\) is a
ramification chart, we replace
\[
Q_{nt}(nz,ny), I_-(z), \phi(z), \phi'(z)
\quad
\text{by}\quad 
Q_{n(t+\ft)}(nz,ny), I_\ft(z),
\phi_\ft(z), \phi_\ft'(z),
\]
respectively. The factors corresponding to nonramification charts are
chosen as in the first case.
\end{enumerate}

\begin{remark}\label{r:choose_set}
We recall from
\Cref{s:vertical_tangent,s:unit_slope_tangent,s:horizontal_tangent,s:vertical_frozen_neighborhood,s:unit_slope_frozen_neighborhood,s:horizontal_frozen_neighborhood}
that a neighborhood centered at a point on a boundary shared by two
curvilinear triangles \(\fT_A\) and \(\fT_B\) is of tangency,
cusp-turning, or tangent frozen type and comes with two possible choices of
local descent/ascent paths: one relative to \(\fT_A\) and one relative to
\(\fT_B\). In these cases, we must choose which set of local paths
\(\sfC^{\rm d}(w_0)\) and \(\sfC^{\rm a}(z_0)\) to use in
\eqref{e:all_term}.

In all such cases, the two contour prescriptions differ only in the relative
nesting of \(\sfC^{\rm d}(w_0)\) and \(\sfC^{\rm a}(z_0)\). In particular,
when considered separately, the paths \(\sfC^{\rm d}(w_0)\) and
\(\sfC^{\rm a}(z_0)\) are the same for the two choices. Hence this choice
does not affect the double-contour integral \eqref{e:all_term} unless the chart collections associated with
\(\fN_\alpha\) and \(\fN_{(y,t)}\) contain concentric charts centered at
\(w_0=z_0\).

In that situation, we may choose the local paths in \eqref{e:all_term}  relative to either \(\fT_A\) or \(\fT_B\), and
choose the contour in the single-contour integral in \eqref{e:def_Aalpha}
accordingly; see \Cref{s:single_integral}. With this consistent choice, the
approximate kernel \(A_\alpha\) in \eqref{e:def_Aalpha} is independent of
whether the paths are chosen relative to \(\fT_A\) or \(\fT_B\).
\end{remark}

\subsection{Concentric charts}

We recall the local descent and ascent paths associated with the charts
constructed in \Cref{s:critical_bulk} and \Cref{s:critical_point}. Except
for frozen charts, every chart carries both a local descent path and a local
ascent path. A frozen chart generally carries only one of these paths,
depending on the sign of the relevant second derivative. The only exceptions
are the tangent frozen charts corresponding to points on the common boundary
of two adjacent curvilinear triangles.
These charts carry both local paths.

Consequently, a chart associated with \((y,t)\) that does not carry a local
ascent path, or a chart associated with \((x,s)\) that does not carry a local
descent path, does not contribute to the terms in \eqref{e:all_term}. Equivalently, in the definition of \(J^{(2)}\) in
\eqref{e:def_Aalpha}, we retain only those charts associated with \((y,t)\)
that carry a local ascent path and those charts associated with \((x,s)\)
that carry a local descent path. The
following lemma gives necessary geometric conditions for retaining a
chart.

A chart is called a \emph{real chart} if it is centered at
\(w_0\in\cC(\bR)\). Given a curvilinear triangle \(\fT\), we recall the
descent and ascent cuts \(\ell_-(z_0;\fT)\) and \(\ell_+(z_0;\fT)\),
respectively, from \Cref{s:descent_ascent_cuts}.

\begin{lemma}\label{l:chart_criterion}
If \(\fN_{(y,t)}\) is associated with a real chart centered at \(z_0\in \cC(\bR)\) that
contains an ascent critical point, then, for some curvilinear triangle
\(\fT\),
\begin{align}\label{e:selectw0}
\fN_{(y,t)}\cap\ell_+(z_0;\fT)\neq\emptyset.
\end{align}

If \((x,s)\in\fN_{\al}\) is associated with a real chart centered at \(w_0\in \cC(\bR)\)
that contains a descent critical point, then, for some curvilinear triangle
\(\fT\),
\begin{align}\label{e:selectw02}
\fN_{\al}\cap\ell_-(w_0;\fT)\neq\emptyset.
\end{align}

Moreover, suppose that \(\fN_{\al}\) and \(\fN_{(y,t)}\) share a pair of
concentric real charts satisfying the preceding assumptions. Then one of
the following holds:
\begin{enumerate}
\item The conditions \eqref{e:selectw0} and \eqref{e:selectw02} hold with
\(w_0=z_0\) and the same curvilinear triangle \(\fT\). In this case, both
\(\fN_{\al}\) and \(\fN_{(y,t)}\) intersect \(\fT\).

\item The conditions \eqref{e:selectw0} and \eqref{e:selectw02} hold with
\(w_0=z_0\) and two different curvilinear triangles \(\fT_A\) and
\(\fT_B\), in either order. In this case, either both
\(\fN_{\al}\) and \(\fN_{(y,t)}\) intersect \(\fT_A\), or both intersect
\(\fT_B\).
\end{enumerate}
\end{lemma}

\begin{proof}[Proof of \Cref{l:chart_criterion}]
A real chart associated with \(\fN_{(y,t)}\) is centered at a point
\(z_0\in\cC(\bR)\) and has an associated base point
$
(y_0,t_0)\in\fN_{(y,t)}
$
satisfying
\[
z_0=y_0-t_0\chi(z_0).
\]
The center \(z_0\) determines the tangent line \(L(z_0)\) defined in
\eqref{e:tangent_line}, which is tangent to the arctic boundary at a point
\((y_0',t_0')\in\fA\).

Recall the descent and ascent cuts
\(\ell_-(z_0;\fT)\) and \(\ell_+(z_0;\fT)\) from
\eqref{e:l0}. If the chart centered at \(z_0\) contains an
ascent critical point, then, for one of the curvilinear triangles \(\fT\)
used in its construction,
\[
(y_0,t_0)\in
\fN_{(y,t)}\cap\ell_+(z_0;\fT).
\]
This proves \eqref{e:selectw0}. The same argument applied to a real chart
centered at \(w_0\) that contains a descent critical point proves
\eqref{e:selectw02}.

It remains to prove the final assertion. Suppose that
\(\fN_{\al}\) and \(\fN_{(y,t)}\) share a pair of concentric real charts.
Then \(w_0=z_0\). If \eqref{e:selectw0} and \eqref{e:selectw02} hold with
the same curvilinear triangle \(\fT\), then both neighborhoods intersect
\(\fT\), and the first case follows.

Suppose instead that \eqref{e:selectw0} and \eqref{e:selectw02} hold with
two different curvilinear triangles \(\fT_A\) and \(\fT_B\). These triangles
are adjacent, and \(L(z_0)\) separates them. In particular,
the slope of \(L(z_0)\) belongs to \(\{0,1,\infty\}\), and
\((y_0',t_0')\) is either a tangent location or a cusp-turning point.

If either \(\fN_{\al}\) or \(\fN_{(y,t)}\) intersects both
\(\fT_A\) and \(\fT_B\), then the second case follows. We may therefore
assume that neither neighborhood intersects both curvilinear triangles.

By symmetry, it suffices to consider the case in which \(L(z_0)\) is
vertical and
\[
\nabla H^*=(0,0)\quad\text{on }\fT_A,
\qquad
\nabla H^*=(1,0)\quad\text{on }\fT_B.
\]
There are two cases.

\begin{enumerate}
\item Suppose that \eqref{e:selectw0} holds with \(\fT=\fT_A\), while
\(\fN_{(y,t)}\) does not intersect \(\fT_B\), and that
\eqref{e:selectw02} holds with \(\fT=\fT_B\), while
\(\fN_{\al}\) does not intersect \(\fT_A\).

By \eqref{e:geometric_descent_ascent},
\(\ell_+(z_0;\fT_A)\) and \(\ell_-(z_0;\fT_B)\) lie at or below
\((y_0',t_0')\). The portions of \(L(z_0)\cap\fT_A\) and
\(L(z_0)\cap\fT_B\) below \((y_0',t_0')\) are intervals with the common
endpoint \((y_0',t_0')\), and hence one is contained in the other.
However, the assumptions imply both that the first interval contains a
point not belonging to the second and that the second interval contains a
point not belonging to the first. This is impossible.

\item Suppose that \eqref{e:selectw0} holds with \(\fT=\fT_B\), while
\(\fN_{(y,t)}\) does not intersect \(\fT_A\), and that
\eqref{e:selectw02} holds with \(\fT=\fT_A\), while
\(\fN_{\al}\) does not intersect \(\fT_B\).

By \eqref{e:geometric_descent_ascent},
\(\ell_+(z_0;\fT_B)\) and \(\ell_-(z_0;\fT_A)\) lie at or above
\((y_0',t_0')\). The portions of \(L(z_0)\cap\fT_A\) and
\(L(z_0)\cap\fT_B\) above \((y_0',t_0')\) are again intervals with the
common endpoint \((y_0',t_0')\), and hence are nested. As in the first
case, the assumptions would require each interval to contain a point not
belonging to the other, which is impossible.
\end{enumerate}

Thus, either both
\(\fN_{\al}\) and \(\fN_{(y,t)}\) intersect \(\fT_A\), or both intersect
\(\fT_B\).

\end{proof}

\begin{lemma}\label{l:number_of_cuts}
We say that \(\fN_{(y,t)}\) is associated with an admissible descent cut
\(\ell_-(z_0;\fT)\) if \(\fN_{(y,t)}\) is associated with a chart centered
at \(z_0\) and \eqref{e:selectw0} holds.

If \(\fN_{(y,t)}\) intersects exactly one curvilinear triangle \(\fT\),
then it is associated with either one or two admissible descent cuts. In the
latter case, the two cuts are separated from each other.

If \(\fN_{(y,t)}\) intersects two adjacent curvilinear triangles
\(\fT_A\) and \(\fT_B\) whose common boundary is contained in the tangent
line \(L(z_0)\), then it is associated with the combined descent cut
\[
\ell_-(z_0;\fT_A)\cup\ell_-(z_0;\fT_B)\subset L(z_0)
\]
and with at most one additional admissible descent cut
\(\ell_-(z_0';\fT)\), where \(z_0'\neq z_0\) and
\(\fT\in\{\fT_A,\fT_B\}\). If this additional cut exists, it is separated
from
\(\ell_-(z_0;\fT_A)\cup\ell_-(z_0;\fT_B)\).
\end{lemma}

\begin{proof}[Proof of \Cref{l:number_of_cuts}]
We first consider the case in which \(\fN_{(y,t)}\) intersects only one
curvilinear triangle \(\fT\) and
\[
\nabla H^*=(1,0)
\]
on \(\fT\). The other two gradient types follow by the symmetry in
\Cref{f:symmetry}.

By \eqref{e:selectw0}, a descent cut
\(\ell_-(z_0;\fT)\) is admissible when \(\fN_{(y,t)}\) intersects the
corresponding ascent cut \(\ell_+(z_0;\fT)\). By
\eqref{e:geometric_descent_ascent}, this ascent cut consists of the tangency
point \((y_0',t_0')\) and the portion of \(L(z_0)\cap\fT\) lying above, or
equivalently to the left of, \((y_0',t_0')\).

By the construction of the charts, the chart corresponding to each
admissible descent cut contains an ascent critical point of
\(S(\,\cdot\,;y,t)\). Distinct nonoverlapping admissible descent cuts have
disjoint associated charts and therefore correspond to distinct ascent
critical points counted in \eqref{e:descent_ascent_count}. Hence there are
at most two admissible descent cuts. Inspection of the configurations in
\Cref{f:curvilinear_triangle} shows that a second admissible descent cut can
occur only in the following cases: \(\fN_{(y,t)}\) is contained in the open
second quadrant; it intersects the vertical boundary at or above the
vertical tangency location; or it intersects the horizontal boundary at or
to the left of the horizontal tangency location.

Whenever two distinct nonoverlapping admissible descent cuts occur,
inspection of the same configurations shows that they are separated from
each other. If \(\fN_{(y,t)}\) is a cusp neighborhood, the two candidate
descent cuts overlap and are therefore counted as a single descent cut. In
all remaining configurations, exactly one admissible descent cut occurs.
This proves the first assertion.

Suppose now that \(\fN_{(y,t)}\) intersects two adjacent curvilinear
triangles \(\fT_A\) and \(\fT_B\). We assume that the tangent line
\(L(z_0)\) is vertical; the other two cases follow in the same way by
symmetry. Then
\[
\nabla H^*=(0,0)\quad\text{on }\fT_A,
\qquad
\nabla H^*=(1,0)\quad\text{on }\fT_B.
\]
The common tangent line \(L(z_0)\) gives the combined descent cut
\[
\ell_-(z_0;\fT_A)\cup\ell_-(z_0;\fT_B).
\]

By \Cref{p:construct_neighborhood1}, if \(\fN_{(y,t)}\) is a tangent or
cusp-turning neighborhood, then, in addition to \(L(z_0)\), it is associated
with at most one additional tangent line to the portion of the arctic
boundary contained in \(\fT_A\cup\fT_B\). This line produces at most one
additional admissible descent cut. If this additional cut exists, the
geometric configurations and the disjointness of the associated charts show
that it is separated from
\[
\ell_-(z_0;\fT_A)\cup\ell_-(z_0;\fT_B).
\]

Otherwise, \(\fN_{(y,t)}\) is a tangent frozen neighborhood. In this case,
in addition to \(L(z_0)\), it is associated with two tangent lines to the
portions of the arctic boundary contained in \(\fT_A\) and \(\fT_B\),
respectively. Moreover, \(\fN_{(y,t)}\) lies to the left of the tangency
point of one of these lines and to the right of the tangency point of the
other. If both additional tangent lines produced admissible descent cuts,
then \eqref{e:geometric_descent_ascent} and \eqref{e:selectw0} would require
\(\fN_{(y,t)}\) to lie to the left of both tangency points, which is
impossible. Therefore, at most one of the two additional tangent lines
produces an admissible descent cut. Whenever this additional cut exists, it
is separated from
\[
\ell_-(z_0;\fT_A)\cup\ell_-(z_0;\fT_B).
\]
This proves the second assertion.
\end{proof}

By \Cref{p:construct_neighborhood2}, any chart in the
\((y,t)\)-collection and any chart in the \((x,s)\)-collection are either
disjoint or concentric; in the latter case, they have the same center. The
following lemma classifies the possible concentric pairs. 

\begin{lemma}[Concentric charts]\label{c:concentric}
Suppose that a chart associated with
\((y,t)\in\fN_{(y,t)}\) and carrying a local ascent path is concentric with
a chart associated with \((x,s)\in\fN_\alpha\) and carrying a local descent
path. Then exactly one of the following alternatives holds:
\begin{enumerate}
\item Both \(\fN_{(y,t)}\) and \(\fN_\alpha\) are liquid neighborhoods.
In this case, the two collections share two pairs of concentric liquid
charts, and the two common centers form a complex-conjugate pair.

\item At least one of \(\fN_{(y,t)}\) and \(\fN_\alpha\) is not a liquid
neighborhood. In this case, the common center is a real point
\(z_0\in\cC(\bR)\) corresponding to a regular arctic point, a cusp point,
a tangent location, or a cusp-turning point. Moreover, the two collections
share at most one concentric pair of real charts.
\end{enumerate}
\end{lemma}

\begin{proof}[Proof of \Cref{c:concentric}]
Suppose first that both neighborhoods are liquid. Each collection contains
a pair of liquid charts with complex-conjugate centers. If one center is
common to the two collections, then its complex conjugate is also common.
Hence the two collections share two pairs of concentric liquid charts, and
the two common centers form a complex-conjugate pair. This proves the first
alternative.

Assume now that at least one of the two neighborhoods is not liquid. By
\Cref{p:construct_neighborhood2}, a non-real center can occur only for a
liquid chart. Therefore, the common center must be a real point
\(z_0\in\cC(\bR)\). The corresponding point on the arctic boundary is a
regular arctic point, a cusp point, a tangent location, or a cusp-turning
point.

It remains to prove that the two collections share at most one concentric
pair of real charts.

Suppose first that \(\fN_{(y,t)}\) intersects a curvilinear triangle
\(\fT\), but no other curvilinear triangle. By
\Cref{l:number_of_cuts}, it is associated with either one descent cut or
two descent cuts satisfying \eqref{e:selectw0}; in the latter case, the two
cuts are separated from each other. By \Cref{l:chart_criterion}, if a
chart in the \((x,s)\)-collection is concentric with one of the
corresponding charts in the \((y,t)\)-collection, then \(\fN_\alpha\) must
intersect the corresponding descent cut, as required by
\eqref{e:selectw02}. By the construction in
\Cref{p:construct_neighborhood2}, \(\fN_\alpha\) can intersect at most one
of these descent cuts. Therefore, the two collections share at most one
concentric pair of real charts in this case.

Suppose next that \(\fN_{(y,t)}\) intersects two curvilinear triangles
\(\fT_A\) and \(\fT_B\). By \Cref{l:number_of_cuts}, it is associated with
the cut
\begin{align}\label{e:vertical_cut}
\ell_-(z_0;\fT_A)\cup\ell_-(z_0;\fT_B),
\end{align}
and possibly with one additional descent cut. The two cuts appearing in
\eqref{e:vertical_cut} have the same center \(z_0\) and hence together
correspond to only one possible concentric pair. If the additional cut
exists, it has a different center and is separated from the cut in
\eqref{e:vertical_cut}.

By \Cref{l:chart_criterion}, a concentric pair can occur only if
\(\fN_\alpha\) intersects the cut in \eqref{e:vertical_cut} or the
additional descent cut, as required by \eqref{e:selectw02}. By the
construction in \Cref{p:construct_neighborhood2}, \(\fN_\alpha\) can
intersect at most one of these separated cuts. Consequently, the two
collections share at most one concentric pair of real charts in this case
as well. This proves the second alternative.
\end{proof}

\subsection{The region \(\fC(\fT;\fN_{(y,t)})\)}
\label{s:deffC}

Assume that \(\fN_{(y,t)}\) intersects a curvilinear triangle \(\fT\). In
this section, we define a subset
\[
\fC(\fT;\fN_{(y,t)})\subset\fT,
\]
which will be used to specify when \(J^{(1)}=0\) in \eqref{e:def_Aalpha}.

The set \(\fC(\fT;\fN_{(y,t)})\) may be viewed as the stability region of
the local contour configuration determined by \((y,t)\). We regard
\((x,s)\) as moving inside \(\fP\), starting from \((y,t)\). As \((x,s)\)
moves, the associated neighborhood \(\fN_\alpha\) and the corresponding
local descent paths vary. As long as \(\fN_\alpha\) intersects
\(\fC(\fT;\fN_{(y,t)})\), we say that it has not escaped from
\(\fC(\fT;\fN_{(y,t)})\). Under this condition, the local descent paths
associated with \((x,s)\) do not cross the local ascent paths associated
with \((y,t)\).

Recall from \Cref{t:frozen_structure} that the curvilinear triangle \(\fT\)
is enclosed by the three boundary pieces
\begin{align}\label{e:boundary_piece}
[\zeta,\zeta_1],\qquad
[\zeta,\zeta_2],\qquad
\text{and the arc of \(\fA\) between \(\zeta_1\) and \(\zeta_2\)},
\end{align}
where \(\zeta\) is a vertex of \(\fP\), and the two sides of \(\fP\)
incident to \(\zeta\), or their linear extensions, are tangent to the
arctic boundary at the corresponding tangency points
\(\zeta_1,\zeta_2\in\fA\).

We now define \(\fC(\fT;\fN_{(y,t)})\). Recall the descent and ascent cuts
\(\ell_-(z_0;\fT)\) and \(\ell_+(z_0;\fT)\), respectively, defined in
\Cref{s:descent_ascent_cuts}. The set
\(\fC(\fT;\fN_{(y,t)})\) is obtained by cutting \(\fT\) along the descent
cuts \(\ell_-(z_0;\fT)\) from \Cref{l:number_of_cuts}, together with
certain boundary cuts described below. Throughout the definition, we use
only descent cuts \(\ell_-(z_0;\fT)\) with centers \(z_0\) satisfying \eqref{e:selectw0}.

\begin{figure}
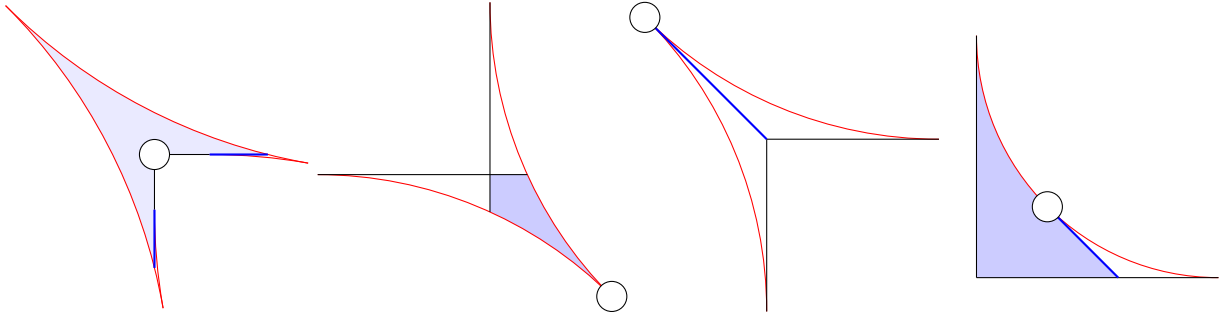

\centering

\begin{subfigure}[t]{0.23\textwidth}
  \centering
  \begin{minipage}[t][4cm][c]{\linewidth}
    \centering

  \end{minipage}

\end{subfigure}
\caption{The region $\fC(\fT;\fN_{(y,t)})$ in blue. In panel (A), $\fN_{(y,t)}$ is a frozen neighborhood. In panel (B), $\fN_{(y,t)}$ is a cusp neighborhood with $\ell_-(z_0;\fT)=\{(y_0',t_0')\}$. In panel (C), $\fN_{(y,t)}$ is a cusp neighborhood with $\ell_+(z_0;\fT)=\{(y_0',t_0')\}$. In panel (D), $\fN_{(y,t)}$ is an arctic neighborhood.} 
\label{f:Cregion}
\end{figure}

\begin{enumerate}
\item Suppose first that \(\nabla H^*=(1,0)\) on \(\fT\). The two linear
boundary pieces in \eqref{e:boundary_piece} determine a vertical tangent
line \(L_\infty\) and a horizontal tangent line \(L_0\), which meet at the
vertex \(\zeta\) of \(\fP\).

If \(\fN_{(y,t)}\) lies to the right of \(L_\infty\), then there is only one
admissible descent cut \(\ell_-(z_0;\fT)\) satisfying
\eqref{e:selectw0}, and this cut also lies to the right of \(L_\infty\).
In this case, we include
\[
L_\infty\cap\fT
\]
as an additional boundary cut. Similarly, if \(\fN_{(y,t)}\) lies below
\(L_0\), then there is only one admissible descent cut
\(\ell_-(z_0;\fT)\) satisfying \eqref{e:selectw0}, and this cut also lies
below \(L_0\). In this case, we include
\[
L_0\cap\fT
\]
as an additional boundary cut.

The descent cuts \(\ell_-(z_0;\fT)\), together with any additional
boundary cuts described above, divide \(\fT\) into closed pieces. If
\(\fN_{(y,t)}\) is a frozen neighborhood, we define
\(\fC(\fT;\fN_{(y,t)})\) to be the unique piece that intersects
\(\fN_{(y,t)}\); see Panel (A) of \Cref{f:Cregion}.

It remains to specify the limiting convention when \(\fN_{(y,t)}\) is not
a frozen neighborhood. If \(\fN_{(y,t)}\) is a cusp, cusp-turning, arctic,
or tangent neighborhood, then several closed pieces may intersect
\(\fN_{(y,t)}\). In this case, we define
\(\fC(\fT;\fN_{(y,t)})\) by approaching \(\fN_{(y,t)}\) through frozen
neighborhoods.

The line \(L(z_0)\) containing the  descent cut
\(\ell_-(z_0;\fT)\) is tangent to the arctic boundary at
\((y_0',t_0')\). The limiting convention selects the closed piece lying,
relative to each descent cut, on the side toward which the outward normal
to the arctic boundary at \((y_0',t_0')\) points. 
At a cusp or
cusp-turning point, this direction is understood as the corresponding
limit from nearby regular arctic points. More explicitly:
\begin{itemize}
\item If there is a single descent cut, centered at \(z_0\), and
\[
\ell_-(z_0;\fT)=\{(y_0',t_0')\},
\]
then this cut does not further divide \(\fT\). The limiting convention
selects the unique piece that intersects \(\fN_{(y,t)}\); see Panel (B) of
\Cref{f:Cregion}.

\item If there is a single descent cut, centered at \(z_0\), and
\[
\ell_+(z_0;\fT)=\{(y_0',t_0')\},
\]
then the limiting piece degenerates to
\[
\fC(\fT;\fN_{(y,t)})
=
\ell_-(z_0;\fT);
\]
see Panel (C) of \Cref{f:Cregion}.

\item In all remaining cases, neither \(\ell_-(z_0;\fT)\) nor
\(\ell_+(z_0;\fT)\) is a singleton for any admissible center \(z_0\). We
define \(\fC(\fT;\fN_{(y,t)})\) to be the unique closed piece that
intersects \(\fN_{(y,t)}\) and contains a nontrivial portion of the
corresponding ascent cut \(\ell_+(z_0;\fT)\) for every admissible center
\(z_0\). Thus, if more than one admissible descent cut is present, the
corresponding selection conditions are imposed simultaneously. See Panel (D) of \Cref{f:Cregion}.
\end{itemize}

\item The remaining two gradient types are treated by the same
construction, with the case-dependent data replaced according to the
following table:
\[
\begin{array}{c|c|c|c}
\hline
\nabla H^*
&
\text{boundary tangent lines}
&
\text{descent cut}
&
\text{additional boundary cuts}
\\
\hline
(1,-1)
&
L_1,\ L_0
&
\text{below/left}
&
\begin{gathered}
L_1\cap\fT \quad \text{if } \fN_{(y,t)}
   \text{ lies to the left of } L_1,\\
L_0\cap\fT \quad \text{if } \fN_{(y,t)}
   \text{ lies below } L_0
\end{gathered}
\\[2.5ex]
\hline
(0,0)
&
L_\infty,\ L_1
&
\text{above/right}
&
\begin{gathered}
L_\infty\cap\fT \quad \text{if } \fN_{(y,t)}
   \text{ lies to the right of } L_\infty,\\
L_1\cap\fT \quad \text{if } \fN_{(y,t)}
   \text{ lies to the left of } L_1
\end{gathered}
\\
\hline
\end{array}
\]

More precisely, if \(\nabla H^*=(1,-1)\), then the two linear boundary
pieces in \eqref{e:boundary_piece} determine the unit-slope tangent line
\(L_1\) and the horizontal tangent line \(L_0\). These lines meet at
\(\zeta\), and all tangent lines to the portion of \(\fA\) contained in
\(\fT\) have slopes in \([0,1]\).

If \(\nabla H^*=(0,0)\), then the two linear boundary pieces determine the
vertical tangent line \(L_\infty\) and the unit-slope tangent line \(L_1\).
These lines meet at \(\zeta\), and all tangent lines to the portion of
\(\fA\) contained in \(\fT\) have slopes in \([1,\infty]\).

The descent cuts \(\ell_-(z_0;\fT)\) satisfying \eqref{e:selectw0},
together with the additional boundary cuts listed in the last column
whenever the corresponding side conditions hold, divide \(\fT\) into
closed pieces. We then define \(\fC(\fT;\fN_{(y,t)})\) by the same
selection rule as in the case \(\nabla H^*=(1,0)\): if
\(\fN_{(y,t)}\) is a frozen neighborhood, we take the unique piece that
intersects \(\fN_{(y,t)}\); otherwise, we take the limiting piece obtained
by approaching \(\fN_{(y,t)}\) through frozen neighborhoods, using the
same conventions as above. If several admissible descent cuts are present,
the corresponding selection conditions are imposed simultaneously.
\end{enumerate}

\begin{lemma}\label{l:boundary_piece_in_C}
Let \(\fT\) be a curvilinear triangle, and let \(L\) be the tangent line
containing one of the two linear boundary pieces of \(\fT\) in
\eqref{e:boundary_piece}. If
\[
\fN_{(y,t)}\cap(L\cap\fT)\neq\emptyset,
\]
then
\[
L\cap\fT\subset\fC(\fT;\fN_{(y,t)}).
\]
\end{lemma}

\begin{proof}
The claim follows by inspecting the classification of curvilinear triangles in
\Cref{f:curvilinear_triangle}. In each possible configuration, once
\(\fN_{(y,t)}\) intersects the boundary piece \(L\cap\fT\), the region
\(\fC(\fT;\fN_{(y,t)})\) contains the entire boundary piece \(L\cap\fT\).
\end{proof}

\subsection{Single-contour integral}\label{s:single_contour}
In this section we specify $J^{(1)}$ in \eqref{e:def_Aalpha}. Assume that $(y,t)$ is bounded away from ramification points.
There are three cases:
\begin{enumerate}
\item $\fN_{(y,t)}$ is a liquid/ramification neighborhood;
\item $\fN_{(y,t)}$ intersects a curvilinear triangle $\fT$ and is bounded away from all other curvilinear triangles;
\item $\fN_{(y,t)}$ is on the shared boundary of two curvilinear triangles $\fT_A$ and $\fT_B$, in which case it is a tangency neighborhood, a cusp-turning neighborhood, or an interface frozen neighborhood.
\end{enumerate}

\subsubsection{Liquid/ramification neighborhood}

Assume that \(\fN_{(y,t)}\) is either a liquid neighborhood or a
ramification neighborhood. Then \((y,t)\) is associated with a
complex-conjugate pair of liquid charts or ramification charts,
respectively.

There are two possibilities for \((x,s)\). First, \((x,s)\) may also be
associated with a complex-conjugate pair of liquid or ramification charts
having the same centers as those associated with \((y,t)\). In this case,
we define \(J^{(1)}\) by one of the following single-contour integrals:
\begin{equation}\label{e:single_term}
\frac{n}{2\pi\ri}
\int_{\sfC}
P_{ns}(nz,nx)\,Q_{nt}(nz,ny)\,\rd z,
\qquad\text{or}\qquad
\frac{n}{2\pi\ri}
\int_{\sfC}
P_{n(s+\ft)}(nz,nx)\,Q_{n(t+\ft)}(nz,ny)\,\rd z.
\end{equation}
The first integral is used in the liquid case, with the contour \(\sfC\)
chosen as in \Cref{s:liquid_integral_contours}. The second integral is used
in the ramification case, with
$
\sfC=\sfC(w_0;(x,s),(y,t)),
$
where \(w_0\in\bC_+\) is the center of the shared ramification chart in
\(\cC_\ft\).

Otherwise, the charts associated with \((x,s)\) are disjoint from those
associated with \((y,t)\). In this case, we set
$
J^{(1)}=0
$
in \eqref{e:def_Aalpha}.

\subsubsection{\(\fN_{(y,t)}\) intersects \(\fT\)}

Assume that \(\fN_{(y,t)}\) intersects a curvilinear triangle \(\fT\) and is
bounded away from all other curvilinear triangles. Then the charts associated
with \((y,t)\) are centered on \(\cC(\bR)\). Recall the region
$
        \fC(\fT;\fN_{(y,t)})\subset\fT
$
from \Cref{s:deffC}.

If
$
        \fC(\fT;\fN_{(y,t)})\cap\fN_\alpha=\emptyset,
$
we set \(J^{(1)}=0\) in \eqref{e:def_Aalpha}. Otherwise, we recall $\sfC(\cdot;(x,s),(y,t))$ from \Cref{def:CE}, and take \(J^{(1)}\) to be
the single-contour term in \eqref{e:single_term}, with
\begin{align}\label{e:defC}
        \sfC=\sfC(E_\fT;(x,s),(y,t))
\end{align}
where $E_\fT$ is chosen as follows.
\begin{enumerate}
\item If \(\nabla H^*=(0,0)\) on \(\fT\), then the portion of the arctic
boundary contained in \(\fT\) has slopes in \([1,\infty]\) and is mapped by
\eqref{e:emb} to
$
        [a_{i-1},b_i]\subset\cC(\bR).
$
In this case, choose any
$
        \max\{x-s,y-t\}<E_\fT<\min\{x,y\}
        $.
        
\item If \(\nabla H^*=(1,0)\) on \(\fT\), then the portion of the arctic
boundary contained in \(\fT\) has slopes in \([-\infty,0]\) and is mapped by
\eqref{e:emb} to
$
        [b_i,\infty_i]\subset\cC(\bR).
$
In this case, choose any
$
        E_\fT>\max\{x,y\}.
$

\item If \(\nabla H^*=(1,-1)\) on \(\fT\), then the portion of the arctic
boundary contained in \(\fT\) has slopes in \([0,1]\) and is mapped by
\eqref{e:emb} to
$
        [\infty_i,a_i]\subset\cC(\bR).
$
In this case, choose any
$
        E_\fT<\min\{x-s,y-t\}.
$
\end{enumerate}

Equivalently, we may write
\begin{align}\begin{split}\label{e:J1form}
&J^{(1)}((x,s),(y,t))
=
\bm1\!\left(
\fC(\fT;\fN_{(y,t)})\cap\fN_\alpha\neq\emptyset
\right)J_\fT((x,s),(y,t)),\\
&J_\fT((x,s),(y,t)):=
\frac{n}{2\pi \ri}
\int_{\sfC(E_\fT;(x,s), (y,t))}
P_{ns}(nz,nx)\,Q_{nt}(nz,ny)\,\rd z 
\end{split}\end{align}

In particular, if \(\fN_\alpha\) is a liquid neighborhood or
\[
\fN_\alpha\cap\fT=\emptyset,
\]
then the charts associated with \((x,s)\) are disjoint from those associated
with \((y,t)\). In either case, the preceding construction gives
\(J^{(1)}=0\).

The following lemma is an immediate consequence of
\Cref{c:single_contour,r:order_condition}.

\begin{lemma}\label{l:JT_nonvanish}
Let \(J_{\fT}((x,s),(y,t))\) be defined as in \eqref{e:J1form}.
\begin{enumerate}
\item
If \(\nabla H^*=(1,0)\) on \(\fT\), then
\(J_{\fT}((x,s),(y,t))\neq0\) only if
\[
t>s,
\qquad
x\geq y.
\]

\item
If \(\nabla H^*=(1,-1)\) on \(\fT\), then
\(J_{\fT}((x,s),(y,t))\neq0\) only if
\[
t>s,
\qquad
y-t\geq x-s.
\]

\item
If \(\nabla H^*=(0,0)\) on \(\fT\), then
\(J_{\fT}((x,s),(y,t))\neq0\) only if
\[
y-t\geq x-s,
\qquad
x\geq y.
\]
\end{enumerate}
\end{lemma}

\subsubsection{$\fN_{(y,t)}$ intersects $\fT_A$ and $\fT_B$}
Assume that $\fN_{(y,t)}$ lies on the shared boundary of two curvilinear triangles $\fT_A$ and $\fT_B$.

There are several possibilities for $(x,s)\in \fN_\al$.

If
$
\fN_\alpha\cap (\fT_A\cup\fT_B)=\emptyset,
$
including the case that $\fN_\alpha$ is a liquid neighborhood, 
then the charts associated with $(x,s)$ are disjoint from those associated with $(y,t)$. In this case, we take
$J^{(1)}=0$
in \eqref{e:def_Aalpha}.

If $\fN_\alpha$ intersects $\fT_A$ but not $\fT_B$, then we define $J^{(1)}$ as in \eqref{e:J1form} with $\fT=\fT_A$.

If $\fN_\alpha$ intersects $\fT_B$ but not $\fT_A$,  then we define $J^{(1)}$ as in \eqref{e:J1form} with $\fT=\fT_B$.

Otherwise, if $\fN_\alpha$ also lies on the shared boundary of $\fT_A$ and $\fT_B$, i.e.
\[
\fN_\alpha\cap (\fT_A\cap \fT_B)\neq\emptyset,
\]
then we may define $J^{(1)}$ as in \eqref{e:J1form} using either $\fT=\fT_A$ or $\fT=\fT_B$, provided that the local paths in \eqref{e:all_term} are chosen consistently relative to the same curvilinear triangle (recall \Cref{r:choose_set}).

\begin{figure}
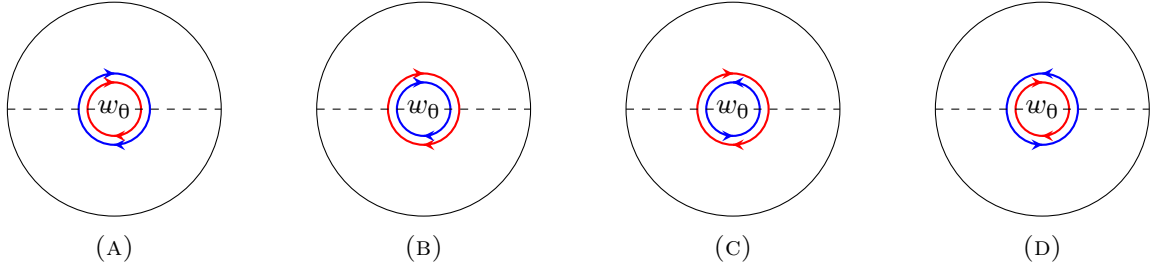

			\begin{subfigure}[t]{0.24\textwidth}
			\centering

		\caption{}

	\end{subfigure}
	\begin{subfigure}[t]{0.24\textwidth}
			\centering
			%
		\caption{}

	\end{subfigure}
			\begin{subfigure}[t]{0.24\textwidth}
			\centering
			%
		\caption{}

	\end{subfigure}
	\begin{subfigure}[t]{0.24\textwidth}
			\centering
			%
		\caption{}

	\end{subfigure}
		\caption{
\label{f:contour_nesting}
Difference of contours in \eqref{e:TAcontour} and \eqref{e:TBcontour}.}
	\end{figure}

\begin{lemma}\label{c:change_triangle}
Let \(\fT_A\) and \(\fT_B\) be adjacent curvilinear triangles. Suppose that we
make one of the following two consistent choices of contours.

First, choose the single contour \(\sfC\) in \eqref{e:single_term} according
to \eqref{e:defC} with \(\fT=\fT_A\), and choose
\begin{align}\label{e:TAcontour}
        \sfC^{\rm d}=\sfC^{\rm d}(w_0),
        \qquad
        \sfC^{\rm a}=\sfC^{\rm a}(w_0)
\end{align}
in \eqref{e:all_term} as constructed in
\Cref{s:vertical_tangent,s:unit_slope_tangent,s:horizontal_tangent,s:vertical_frozen_neighborhood,s:unit_slope_frozen_neighborhood,s:horizontal_frozen_neighborhood}
relative to \(\fT_A\).

Alternatively, choose the single contour \(\sfC\) in \eqref{e:single_term}
according to \eqref{e:defC} with \(\fT=\fT_B\), and choose
\begin{align}\label{e:TBcontour}
        \sfC^{\rm d}=\sfC^{\rm d}(w_0),
        \qquad
        \sfC^{\rm a}=\sfC^{\rm a}(w_0)
\end{align}
relative to \(\fT_B\).

Then the sum of the corresponding single-contour integral and double-contour
integral is the same for the two choices.

Consequently, if both \(\fN_\alpha\) and \(\fN_{(y,t)}\) intersect the shared
boundary of \(\fT_A\) and \(\fT_B\), that is,
\[
        \fN_\alpha\cap(\fT_A\cap\fT_B)\neq\emptyset,
        \qquad
        \fN_{(y,t)}\cap(\fT_A\cap\fT_B)\neq\emptyset,
\]
then constructing the contours in \eqref{e:all_term} and \eqref{e:J1form} using either \(\fT=\fT_A\) or \(\fT=\fT_B\)
produces the same approximate kernel
$
        A_\alpha((x,s),(y,t)).
$
\end{lemma}

\begin{proof}
We prove \Cref{c:change_triangle} in the case where the shared boundary of
\(\fT_A\) and \(\fT_B\) is a vertical tangent line. The cases where the
shared boundary is a unit-slope or horizontal tangent line are treated in the
same way, so we omit them.

In the vertical case, we have \(\nabla H^*=(0,0)\) on \(\fT_A\). Hence we may
choose
\[
        \max\{x-s,y-t\}<E_A<\min\{x,y\},
\]
and define
\begin{align}
J_A
:=
\frac{n}{2\pi \ri}
\int_{\sfC(E_A;(x,s),(y,t))}
P_{ns}(nz,nx)\,Q_{nt}(nz,ny)\,\rd z .
\end{align}
On \(\fT_B\), we have \(\nabla H^*=(1,0)\). Hence we may choose
\[
        E_B>\max\{x,y\},
\]
and define
\begin{align}
J_B
:=
\frac{n}{2\pi \ri}
\int_{\sfC(E_B;(x,s),(y,t))}
P_{ns}(nz,nx)\,Q_{nt}(nz,ny)\,\rd z .
\end{align}

The difference between the two choices of the single contour,
$
        \sfC(E_A;(x,s),(y,t))$
and $
        \sfC(E_B;(x,s),(y,t)),
$
is therefore
\begin{equation}\label{e:residual}
J_B-J_A
=
\frac{n}{2\pi \ri}
\oint_{\omega}
P_{ns}(nz,nx)\,Q_{nt}(nz,ny)\,\rd z ,
\end{equation}
where \(\omega\) is a counterclockwise contour enclosing the interval
$
        [\min\{x,y\},\max\{x,y\}].
$

We now compare the corresponding double-contour terms. In the present case,
the shared vertical tangent line corresponds to
$
        w_0=z_0=b_i
$
for some \(1\leq i\leq d\). Hence
\[
        I_+=I_-=I_i.
\]
Moreover, the residue at \(w=z\) of the double-contour integrand in
\eqref{e:all_term} is
\[
        \frac{n}{2\pi \ri}
        P_{ns}(nz,nx)\,Q_{nt}(nz,ny),
\]
which is precisely the integrand of the single-contour term
\eqref{e:residual}.

By the contour constructions in
\Cref{s:vertical_tangent,s:unit_slope_tangent,s:horizontal_tangent,s:vertical_frozen_neighborhood,s:unit_slope_frozen_neighborhood,s:horizontal_frozen_neighborhood},
the choices in \eqref{e:TAcontour} and \eqref{e:TBcontour} differ only by the
nesting of the local circular pieces. All possible nestings are shown in
\Cref{f:contour_nesting}, up to reversing the orientations of both the blue
and red paths.

Passing from the contour choice associated with \(\fT_A\) to the one
associated with \(\fT_B\) amounts to deforming the local descent
piece (blue circle) through the local ascent piece (red circle). During this deformation, the
only pole crossed is the pole at \(w=z\). By the residue computation above,
the resulting change in the double-contour term is exactly the negative of
\eqref{e:residual}, with the sign determined by the orientations of the two
local circular pieces. Thus the change in the double-contour term cancels the
change
$
        J_B-J_A
$
in the single-contour term. Therefore the sum of the single-contour and
double-contour contributions is the same for the two choices. This proves the
first statement of \Cref{c:change_triangle}.

For the second statement, suppose that both \(\fN_\alpha\) and
\(\fN_{(y,t)}\) intersect the shared boundary of \(\fT_A\) and \(\fT_B\).
Then, by \Cref{l:boundary_piece_in_C},
\[
        \fN_\alpha\cap\fC(\fT_A;\fN_{(y,t)})\neq\emptyset,
        \qquad
        \fN_\alpha\cap\fC(\fT_B;\fN_{(y,t)})\neq\emptyset.
\]
Hence, if we construct the contour in \eqref{e:J1form} using
\(\fT=\fT_A\), we obtain the single-contour term \(J_A\), while using
\(\fT=\fT_B\) gives \(J_B\). The first part of the lemma shows that the
corresponding sums of single-contour and double-contour terms are equal.
Therefore the resulting approximate kernel
$
        A_\alpha((x,s),(y,t))
$
is independent of whether the construction is made with \(\fT_A\) or
\(\fT_B\). This proves the second statement.
\end{proof}

\subsection{The region $\fC(\fT;y,t)$ and interlacing condition}\label{s:deffC2}
Fix a point $(y,t)\in \fT$  for some curvilinear triangle \(\fT\). We define
the pointwise analogue of \(\fC(\fT;\fN_{(y,t)})\), denoted by
\[
\fC(\fT;y,t)\subset\fT.
\]
The construction is the same as that of
\(\fC(\fT;\fN_{(y,t)})\) in \Cref{s:deffC}, except that all side conditions
are imposed directly on the point \((y,t)\), rather than on the neighborhood
\(\fN_{(y,t)}\).

Recall from \Cref{p:associate_critical_points} that \((y,t)\) is associated
with two or three real critical points, counted with multiplicity. Among
these, we consider the distinct critical-point locations
\(z_c\in\cC(\bR)\) for which the tangent line \(L(z_c)\) is tangent to the
portion of \(\fA\) contained in \(\fT\).

For each such \(z_c\), recall the descent and ascent cuts
\(\ell_-(z_c;\fT)\) and \(\ell_+(z_c;\fT)\) from
\eqref{e:l0}. As in the construction of
\(\fC(\fT;\fN_{(y,t)})\), we retain only those potential ascent critical
points \(z_c\) satisfying
\begin{align}\label{e:selectwc}
(y,t)\in\ell_+(z_c;\fT).
\end{align}
Equivalently, \((y,t)\) lies in the ascent part of
\(L(z_c)\cap\fT\). 

We define \(\fC(\fT;y,t)\) by applying the construction of
\(\fC(\fT;\fN_{(y,t)})\) with the replacements
\[
\fN_{(y,t)}
\quad\longrightarrow\quad
(y,t),
\qquad
z_0
\quad\longrightarrow\quad
z_c,
\qquad
\ell_-(z_0;\fT)
\quad\longrightarrow\quad
\ell_-(z_c;\fT).
\]
Thus, the cuts consist of the descent cuts \(\ell_-(z_c;\fT)\) satisfying
\eqref{e:selectwc}, together with the boundary cuts specified below. Every
side condition is tested at the point \((y,t)\).

\[
\begin{array}{c|c|c|c}
\hline
\nabla H^*
&
\text{boundary tangent lines}
&
\text{side of descent cut}
&
\text{additional boundary cuts}
\\
\hline
(1,0)
&
L_\infty,\ L_0
&
\text{below/right}
&
\begin{gathered}
L_\infty\cap\fT
\text{ if }(y,t)\text{ lies to the right of }L_\infty,\\
L_0\cap\fT
\text{ if }(y,t)\text{ lies on or below }L_0
\end{gathered}
\\[2.5ex]
\hline
(1,-1)
&
L_1,\ L_0
&
\text{below/left}
&
\begin{gathered}
L_1\cap\fT
\text{ if }(y,t)\text{ lies to the left of }L_1,\\
L_0\cap\fT
\text{ if }(y,t)\text{ lies on or below }L_0
\end{gathered}
\\[2.5ex]
\hline
(0,0)
&
L_\infty,\ L_1
&
\text{above/right}
&
\begin{gathered}
L_\infty\cap\fT
\text{ if }(y,t)\text{ lies to the right of }L_\infty,\\
L_1\cap\fT
\text{ if }(y,t)\text{ lies to the left of }L_1
\end{gathered}
\\
\hline
\end{array}
\]

The lines
\(L_\infty,L_0,L_1\) are the boundary tangent lines determined by the
boundary pieces of \(\fT\), as in \Cref{s:deffC}.

The descent cuts \(\ell_-(z_c;\fT)\), together with the additional boundary
cuts specified in the table, divide \(\fT\) into pieces, which we take to be
closed. If \((y,t)\notin\fA\cap\fT\), then \((y,t)\) lies on none of these
cuts and hence belongs to a unique piece. We define
\(\fC(\fT;y,t)\) to be this piece.

If \((y,t)\in\fA\cap\fT\), then \((y,t)\) may belong to the boundaries of
several pieces. In this case, we define \(\fC(\fT;y,t)\) using the same
limiting convention as in the definition of
\(\fC(\fT;\fN_{(y,t)})\), with \((y,t)\) approached through the interior
of \(\fT\). More explicitly, for each admissible center \(z_c\), let
\((y',t')\) be the point at which \(L(z_c)\) is tangent to the arctic
boundary.

\begin{itemize}
\item If there is a single admissible descent cut, centered at \(z_c\), and
\[
\ell_-(z_c;\fT)=\{(y',t')\},
\]
then this cut does not further divide \(\fT\), and
\(\fC(\fT;y,t)\) is the unique piece containing \((y,t)\).

\item If there is a single admissible descent cut, centered at \(z_c\), and
\[
\ell_+(z_c;\fT)=\{(y',t')\},
\]
then the limiting piece degenerates to
\[
\fC(\fT;y,t)
=
\ell_-(z_c;\fT).
\]

\item In all remaining configurations, we define
\(\fC(\fT;y,t)\) to be the unique closed piece containing \((y,t)\) and a
nontrivial portion of the corresponding ascent cut
\(\ell_+(z_c;\fT)\) for every admissible center \(z_c\). If more than one
admissible descent cut is present, the corresponding selection conditions
are imposed simultaneously.
\end{itemize}

Now regard \((x,s)\) as varying continuously, and suppose that it crosses a
nontrivial cut \(\ell_-(z_c;\fT)\), where \(z_c\) satisfies
\eqref{e:selectwc}. At the crossing, one of the descent critical points \(w_c\)
associated with \((x,s)\) coincides with \(z_c\). The critical point
\(w_c=z_c\) carries a steepest-descent path, while \(z_c\), viewed as
a critical point associated with \((y,t)\), carries a
steepest-ascent path. Thus, the corresponding descent and ascent paths meet
at the common critical point. As \((x,s)\) crosses the cut, the moving
descent critical point \(w_c\) passes from one side of \(z_c\) to the other,
changing the relative order of the relevant critical points. This motivates
the following definition.

\begin{definition}\label{d:interlacing}
Take two lattice points \((x,s), (y,t)\in\fP\cap\bZ^2/n\). Define
\begin{align}\label{e:def_interlacing}
\cI((x,s),(y,t))
:=
\begin{cases}
0, & (y,t)\in\fL,\\[1mm]
\bm1\!\left(
(x,s)\in
\fC\bigl(\operatorname{Cell}^{\rw}(y,t);y,t\bigr)
\right),
& (y,t)\in\fP\setminus\fL,
\end{cases}
\end{align}
where \(\operatorname{Cell}^{\rw}(y,t)\) is the cell assignment introduced
in \Cref{s:assign_cell}.

We say that the critical points associated with \((x,s)\)
\emph{interlace} those associated with \((y,t)\) if
\[
\cI((x,s),(y,t))=1.
\]
\end{definition}

Roughly speaking, \(\cI((x,s),(y,t))=1\) means that the descent critical
points associated with \((x,s)\) interlace with the ascent critical points
associated with \((y,t)\).

\begin{lemma}\label{l:interlace_critical_point}
Let \((y,t)\in\fP\cap\bZ^2/n\) represent a white triangle, and let
\((x,s)\in\fP\cap\bZ^2/n\) represent a blue triangle. Suppose that
\((y,t)\in\fP\setminus\fL\), set
\[
\fT:=\operatorname{Cell}^{\rw}(y,t),
\]
and assume that
\begin{align}\label{e:IJ}
\cI((x,s),(y,t))J_{\fT}((x,s),(y,t))\neq0.
\end{align}
Then \((x,s)\in\fT\). Moreover, for every
\[
z_c^{\rm d}\in\operatorname{Crit}^{\rm d}(y,t^-),
\]
there exists
\[
w_c\in\operatorname{Crit}^{\rm d}(x,s)
\]
such that
\begin{align}\label{e:Srelation}
\Re S(w_c;y,t)
\geq
\Re S(z_c^{\rm d};y,t).
\end{align}
If, in addition, \(z_c^{\rm d}\) is bounded away from
\(\operatorname{Crit}^{\rm d}(x,s)\), then there exists \(\fc'>0\) such
that
\begin{align}\label{e:Srelation2}
\Re S(w_c;y,t)
\geq
\Re S(z_c^{\rm d};y,t)+\fc'.
\end{align}
\end{lemma}

\begin{proof}[Proof of \Cref{l:interlace_critical_point}]
We prove the statement under the assumption that
\[
\nabla H^*=(1,0)
\]
on \(\fT\). The other two cases follow from the same argument, using the
symmetry in \Cref{f:symmetry}.

In this case, all critical points corresponding to tangencies to the portion
of the arctic boundary contained in \(\fT\) lie on
\([b_i,\infty_i]\subset\cC(\bR)\), and the slopes of the corresponding
tangent lines lie in \([-\infty,0]\). 
Moreover, by \eqref{e:vert_tangent_diff}, if \(y>b_i\), then the branch cut
\([b_i,y]\) is a source for the steepest-descent flow of
\(\Re S(\,\cdot\,;y,t)\). Let the horizontal extended side of \(\fT\) be
supported on the line \(s=s_0\). By \eqref{e:hor_tangent_diff}, if
\(t>s_0\), then \(\infty_i\) is a sink for the same steepest-descent flow.

Since \(\cI((x,s),(y,t))\) is an indicator, \eqref{e:IJ} implies that
\begin{align}\label{e:nonzero_cond}
\cI((x,s),(y,t))=1
\qquad\text{and}\qquad
J_{\fT}((x,s),(y,t))\neq0.
\end{align}
By \Cref{d:interlacing}, the first relation in
\eqref{e:nonzero_cond} is equivalent to
$
(x,s)\in\fC(\fT;y,t).
$
Since
$
\fC(\fT;y,t)\subseteq\fT,
$
it follows that \((x,s)\in\fT\). Moreover, by
\Cref{l:JT_nonvanish}, the second relation in
\eqref{e:nonzero_cond} implies
\begin{align}\label{e:xy_order}
x\geq y,
\qquad
s<t.
\end{align}

We first prove \eqref{e:Srelation} when \((y,t)\) lie
neither on the arctic boundary nor on an extended side, and all relevant
critical points are nondegenerate. In particular,
\[
\operatorname{Crit}^{\rm d}(y,t^-)
=
\operatorname{Crit}^{\rm d}(y,t).
\]
Fix
$
z_c^{\rm d}\in\operatorname{Crit}^{\rm d}(y,t).
$
Then \(z_c^{\rm d}\) lies on the real arc
\begin{align}\label{e:interval}
[\max\{b_i,y\},\infty_i]\subset\cC(\bR).
\end{align}
The ascent critical points associated with \((y,t)\), together with the two
boundary points \(\max\{b_i,y\}\) and \(\infty_i\), divide the real arc
\eqref{e:interval} into consecutive closed subarcs.
Let
\([\alpha,\beta]\) be the subinterval containing \(z_c^{\rm d}\). Since the
real critical points of \(S(\,\cdot\,;y,t)\) alternate between ascent and
descent critical points, \(z_c^{\rm d}\) is the unique descent critical
point in \([\alpha,\beta]\). Moreover,
\(\Re S(\,\cdot\,;y,t)\) is strictly decreasing toward
\(z_c^{\rm d}\) from one side and strictly increasing away from it on the
other. Consequently,
\begin{align}\label{e:Sminimum}
\Re S(w;y,t)
\geq
\Re S(z_c^{\rm d};y,t),
\qquad
w\in[\alpha,\beta],
\end{align}
with equality only when \(w=z_c^{\rm d}\). It therefore suffices to show
that \([\alpha,\beta]\) contains a descent critical point associated with
\((x,s)\).

Choose a continuous path from \((y,t)\) to \((x,s)\) with the following
properties:
\begin{enumerate}
\item
The path remains in \(\fC(\fT;y,t)\) and, except at its initial point,
lies in
\begin{align}\label{e:xxregion}
\{(x',s')\in\fT:x'\geq y,\ s'<t\}.
\end{align}

\item
Whenever the path crosses a vertical extended side, it does so from left to
right, and whenever it crosses a horizontal extended side, it does so from
above to below. If the terminal point \((x,s)\) lies on a horizontal
extended side, the path approaches \((x,s)\) from above.
\end{enumerate}

Such a path exists by the geometric construction of
\(\fC(\fT;y,t)\) in \Cref{s:deffC}. Indeed, when
\(\nabla H^*=(1,0)\), the relevant descent cuts lie below, or
equivalently to the right of, their tangency points, while the additional
boundary cuts ensure that the extended sides are crossed only in the
directions stated above. The path may be chosen, except at its endpoints, to avoid the boundary cuts of
\(\fC(\fT;y,t)\).

Starting from \(z_c^{\rm d}\), track the corresponding descent critical
point as the spatial point moves along this path. Away from the extended
sides, the critical point varies continuously by the implicit function
theorem and remains a descent critical point. At an extended-side crossing,
we use the corresponding one-sided continuation described in
\Cref{r:change_critical} and its symmetric versions. In the crossing
directions chosen above, the tracked descent critical point does not
disappear.

The tracked critical point cannot cross an ascent critical point associated
with \((y,t)\). Indeed, if it coincided with such an ascent critical point,
then the moving spatial point would lie on the corresponding descent cut
used in the construction of \(\fC(\fT;y,t)\). This contradicts the choice
of the path. Thus, the tracked critical point cannot leave
\([\alpha,\beta]\) through an internal endpoint.

It remains to rule out escape through an endpoint of the arc
\eqref{e:interval}. If \(y>b_i\), then every critical point associated with
a point \((x',s')\) in the region \eqref{e:xxregion} satisfies
\[
w_c=x'-s'\chi(w_c)>x'\geq y,
\]
because \(\chi(w_c)<0\). Hence the tracked
critical point cannot escape through the left endpoint \(y\).

The remaining boundary points \(b_i\) and \(\infty_i\) correspond to the
vertical and horizontal extended sides, respectively. A critical point can
reach either of these boundary points only when the moving spatial point
crosses the corresponding extended side or when the terminal point
\((x,s)\) lies on that side. Again, by \Cref{r:change_critical}, neither
crossing a vertical extended side from left to right nor approaching a
terminal point on a vertical extended side from the left can cause the
tracked descent critical point to disappear through \(b_i\). Similarly,
neither crossing a horizontal extended side from above to below nor
approaching a terminal point on a horizontal extended side from above can
cause it to disappear through \(\infty_i\).

Thus, the tracked critical point persists along the entire path and
terminates at a descent critical point
$
w_c\in\operatorname{Crit}^{\rm d}(x,s).
$
Moreover,
$
w_c\in[\alpha,\beta].
$
Applying \eqref{e:Sminimum} gives
\[
\Re S(w_c;y,t)
\geq
\Re S(z_c^{\rm d};y,t),
\]
which proves \eqref{e:Srelation} in the nondegenerate case.

The general case follows by approximation. Perturb \((y,t)\) downward
through generic points of \(\operatorname{int}\fT\), as in the definition
of \(\operatorname{Crit}^{\rm d}(y,t^-)\), so that the perturbed point lies
neither on the arctic boundary nor on an extended side. The perturbation may
be chosen sufficiently small that the relations in \eqref{e:xy_order} are
preserved and \((x,s)\) remains in the corresponding component of the
interlacing region.

Apply the preceding argument to each perturbed configuration. By the
limiting convention in \eqref{e:defCdown}, the corresponding initial
descent critical points converge to \(z_c^{\rm d}\). Passing to the limit and using the continuity of the critical values gives
\[
\Re S(w_c;y,t)
\geq
\Re S(z_c^{\rm d};y,t).
\]

Finally, suppose that \(z_c^{\rm d}\) is bounded away from
\(\operatorname{Crit}^{\rm d}(x,s)\). Then the critical point \(w_c\)
constructed above is bounded away from \(z_c^{\rm d}\). By the strict
monotonicity of \(\Re S(\,\cdot\,;y,t)\) on the two sides of
\(z_c^{\rm d}\) in \([\alpha,\beta]\), the inequality in
\eqref{e:Sminimum} is strict away from \(z_c^{\rm d}\). Compactness of the
relevant closed subinterval therefore gives a constant \(\fc'>0\) such
that
\[
\Re S(w_c;y,t)
\geq
\Re S(z_c^{\rm d};y,t)+\fc'.
\]
This proves \eqref{e:Srelation2} and completes the proof.
\end{proof}

\begin{figure} 
  \begin{subfigure}{0.23\textwidth}
    \centering

  \end{subfigure}
  \caption{The region $\fC(\fT;y,t)$ when $(y,t)$ is in four quadrants.}
  \label{f:Cuts}
\end{figure}

\section{Revisit Critical Points}
\label{s:revisit_critical}

In this section, we collect some refined
properties of the critical points and critical values of the tiling action function.

\subsection{Critical points associated with neighborhoods}

\begin{definition}\label{def:neighborhood_ascent_descent_critical}
Let \((x,s)\in\fP\), and let \(\fN\) be a neighborhood associated with
\((x,s)\), together with its collection of charts, as constructed in
\Cref{p:construct_neighborhood1}. Deform the local descent contours
\(\sfC^{\rm d}(w_0)\) carried by these charts as in
\Cref{s:critical_bulk,s:critical_point}. Each nonempty resulting contour is a
steepest-descent path \(\sfD^{\rm d}(w_c)\). We denote by
\[
\operatorname{Crit}^{\rm d}(x,s;\fN)
\]
the set of descent critical points \(w_c\) obtained in this way.

Similarly, deform the associated local ascent contours
\(\sfC^{\rm a}(z_0)\). Each nonempty resulting contour is a
steepest-ascent path \(\sfD^{\rm a}(z_c)\). We denote by
\[
\operatorname{Crit}^{\rm a}(x,s;\fN)
\]
the set of ascent critical points \(z_c\) obtained in this way.
\end{definition}

Some circular local contours contract during the deformation and therefore
do not produce genuine critical points. The following lemma identifies the
cases in which this occurs.

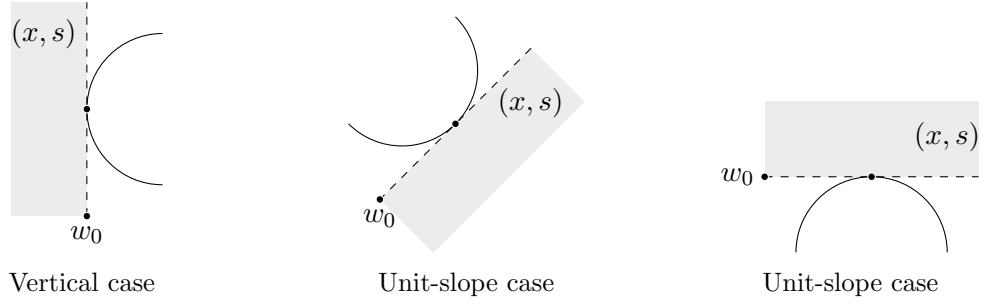
\begin{figure}			
		\begin{subfigure}[t]{0.3\textwidth}

		\centering
			\begin{tikzpicture}
			 \fill[gray!15]
  (-1,{-sqrt(2)})--(0,{-sqrt(2)})--(0,{sqrt(2)})--(-1,{sqrt(2)})
  -- cycle;

			\draw (1,1) arc (90:180:1);
			\draw[] (0,0) arc (180:270:1);
			\draw[dashed] (0,{-sqrt(2)})--(0,{sqrt(2)});
			%\draw[dashed] (1,{-sqrt(2)})--(-1,{-sqrt(2)});
			
			\draw[white, fill=black]  (0,0) circle (0.05);

			\draw[white, fill=black]  (0,{-sqrt(2)}) circle (0.05);
			\draw[](0,{-sqrt(2)}) node[below]{$w_0$};
			
			   			\draw[](0,1) node[left]{$(x,s)$};

			\end{tikzpicture}
				
			\caption*{Vertical case}
			
			\end{subfigure}
			\begin{subfigure}[t]{0.3\textwidth}

		\centering
			\begin{tikzpicture}[rotate=135]
			 \fill[gray!15]
  (-1,{-sqrt(2)})--(0,{-sqrt(2)})--(0,{sqrt(2)})--(-1,{sqrt(2)})
  -- cycle;

			\draw (1,1) arc (90:180:1);
			\draw[] (0,0) arc (180:270:1);
			\draw[dashed] (0,{-sqrt(2)})--(0,{sqrt(2)});
			%\draw[dashed] (1,{-sqrt(2)})--(-1,{-sqrt(2)});
			
			\draw[white, fill=black]  (0,0) circle (0.05);

			\draw[white, fill=black]  (0,{sqrt(2)}) circle (0.05);
			\draw[](0,{sqrt(2)}) node[below]{$w_0$};
			
			   			\draw[](-0.1,-0.5) node[right]{$(x,s)$};

			\end{tikzpicture}
				
			\caption*{Unit-slope case}
			
			\end{subfigure}
	\begin{subfigure}[t]{0.3\textwidth}

		\centering
			\begin{tikzpicture}[rotate=-90]
			 \fill[gray!15]
  (-1,{-sqrt(2)})--(0,{-sqrt(2)})--(0,{sqrt(2)})--(-1,{sqrt(2)})
  -- cycle;

			\draw (1,1) arc (90:180:1);
			\draw[] (0,0) arc (180:270:1);
			\draw[dashed] (0,{-sqrt(2)})--(0,{sqrt(2)});
			%\draw[dashed] (1,{-sqrt(2)})--(-1,{-sqrt(2)});
			
			\draw[white, fill=black]  (0,0) circle (0.05);

			\draw[white, fill=black]  (0,{-sqrt(2)}) circle (0.05);
			\draw[](0,{-sqrt(2)}) node[left]{$w_0$};
			
			   			\draw[](-0.2,1) node[above]{$(x,s)$};

			\end{tikzpicture}
				
			\caption*{Unit-slope case}
			
			\end{subfigure}

\caption{
Cases the integrand \eqref{e:w_integral} is holomorphic in $w$.}
	\label{f:vanish_contour}
	\end{figure}

\begin{lemma}\label{l:holomorphic}
Fix $(x,s)\in \bZ^2/n\cap \fP$. Consider the \(w\)-integral in the double-contour integral
\eqref{e:all_term},
\begin{align}\label{e:w_integral}
\int_{\sfC^{\rm d}(w_0)}
P_{ns}(nw,nx)\,I_+(w)\,(\cdots)\,\rd w,
\end{align}
where \(\sfC^{\rm d}(w_0)\) is the local descent contour carried by a
chart \(\fU\) centered at \(w_0\). The contour
\(\sfC^{\rm d}(w_0)\) is circular in each of the following cases:
\begin{enumerate}
\item
\(\fU\) is a vertical, unit-slope, or horizontal tangent frozen chart
carrying a local descent contour;

\item
\(\fU\) is a vertical-tangent chart and the arctic boundary is locally
tangent from the right; \(\fU\) is a unit-slope-tangent chart and the
arctic boundary is locally tangent from the right; or \(\fU\) is a
horizontal-tangent chart and the arctic boundary is locally tangent from
below.
\end{enumerate}

Moreover, the integrand in \eqref{e:w_integral} is holomorphic in \(w\) in
a sufficiently small neighborhood of \(w_0\) under the following respective
conditions:
\begin{enumerate}
\item
in the vertical case, \((x,s)\) lies to the left of the tangent line;

\item
in the unit-slope case, \((x,s)\) lies to the right of the tangent
line;

\item
in the horizontal case, \((x,s)\) lies on or above the tangent line.
\end{enumerate}
Consequently, in these cases the circular contour can be contracted and
makes no contribution to the integral. See \Cref{f:vanish_contour}.

The analogous statement holds for the \(z\)-integral
\[
\int_{\sfC^{\rm a}(z_0)}
Q_{nt}(nz,ny)\,I_-(z)^{-1}\,(\cdots)\,\rd z
\]
in \eqref{e:twoint}, where \(\sfC^{\rm a}(z_0)\) is the local ascent
contour. In this statement, descent is replaced by ascent, and the
corresponding tangency-side and point-location conditions are reversed.
\end{lemma}

\begin{proof}
The claim follows from the local analysis in
\Cref{s:vertical_tangent,s:unit_slope_tangent,s:horizontal_tangent,s:vertical_frozen_neighborhood,s:unit_slope_frozen_neighborhood,s:horizontal_frozen_neighborhood}.
\end{proof}

\subsection{Extra critical points}
\label{r:spillover}

In \Cref{p:construct_neighborhood1}, each point \((x,s)\in\fP\) is
assigned an open neighborhood \(\fN\), together with a finite collection of
pairwise disjoint \(\fc\)-charts. The union of these charts contains all
critical points associated with \((x,s)\) in
\Cref{p:associate_critical_points}. After deforming the local descent
contours carried by these charts, we obtain the neighborhood-dependent set
\[
\operatorname{Crit}^{\rm d}(x,s;\fN).
\]
By construction,
\begin{align}\label{e:Crit_inclusion}
\operatorname{Crit}^{\rm d}(x,s)
\subseteq
\operatorname{Crit}^{\rm d}(x,s;\fN).
\end{align}
The inclusion in \eqref{e:Crit_inclusion} may be strict: the chosen charts
may contain additional descent critical points that are not included in the
intrinsic set \(\operatorname{Crit}^{\rm d}(x,s)\).

This discrepancy arises because the charts are chosen so that the local
geometric configuration remains stable throughout the entire neighborhood
\(\fN\), whereas \Cref{p:associate_critical_points} selects critical points
according to the liquid region or the curvilinear triangle assigned to the
particular point \((x,s)\). The additional descent critical points in
\[
\operatorname{Crit}^{\rm d}(x,s;\fN)
\setminus
\operatorname{Crit}^{\rm d}(x,s)
\]
can arise in the following two cases.

\begin{enumerate}
\item
\textbf{Liquid spillover.}
Suppose that \((x,s)\in\fL\) and that \(\fN\) is an arctic or tangent
neighborhood. In this case, \Cref{p:associate_critical_points} associates
with \((x,s)\) only the pair of complex-conjugate critical points in the
liquid region.

The chart collection in \Cref{p:construct_neighborhood1}, however, must
remain valid throughout \(\fN\). In addition to the arctic or tangent
charts containing the two liquid critical points near the distinguished
arctic or tangency point, the collection may therefore retain an additional
regular, tangent, or cusp frozen chart. This chart may contain a real
descent critical point corresponding to another tangent line to the arctic
boundary. Such a critical point belongs to
\[
\operatorname{Crit}^{\rm d}(x,s;\fN)
\setminus
\operatorname{Crit}^{\rm d}(x,s).
\]

\item
\textbf{Boundary spillover.}
Suppose that \((x,s)\in\fP\setminus\fL\) lies in a tangent or interface
frozen neighborhood that intersects the boundary shared by two adjacent
curvilinear triangles \(\fT_A\) and \(\fT_B\). Assume, for instance, that
\[
(x,s)\in \fT_A,
\qquad
(x,s)\notin\fT_B,
\]
so that \((x,s)\) does not lie on the shared boundary. Then
\Cref{p:associate_critical_points} associates with \((x,s)\) only those
critical points whose points of tangency lie on the portion of the arctic
boundary contained in \(\fT_A\).

On the other hand, the chart collection in
\Cref{p:construct_neighborhood1} may retain at most one additional regular,
tangent, or cusp frozen chart containing a descent critical point whose
point of tangency lies on the portion of the arctic boundary contained in
\(\fT_B\). Since \((x,s)\not\in \fT_B\), this critical point
is not contained in \(\operatorname{Crit}^{\rm d}(x,s)\), but it may belong
to \(\operatorname{Crit}^{\rm d}(x,s;\fN)\).
\end{enumerate}

We will show that the terms indexed by these additional critical points make
negligible contributions to the neighborhood-dependent double-contour integrals \eqref{e:all_term}. The following lemma establishes the required gap in
the tiling action between the additional descent critical points and the
critical points selected in \Cref{p:associate_critical_points}.

\begin{figure}
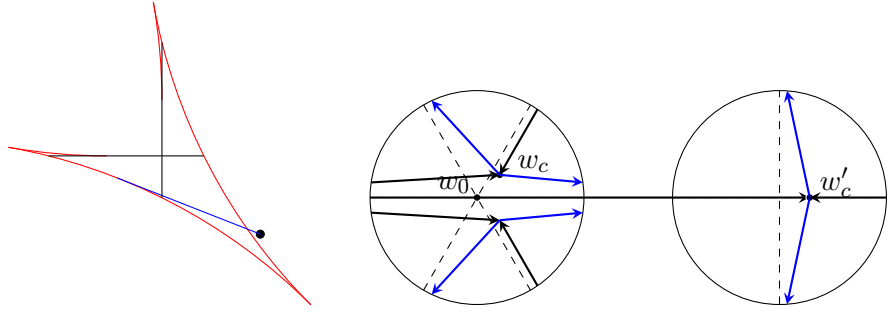

\begin{subfigure}{0.32\textwidth}
    \centering

				
			\end{subfigure}
			\caption{In the left panel, arctic neighborhood with an extra frozen chart; in the right panel, the corresponding local descent paths.
}\label{f:extra_chart1}
\end{figure}

\begin{figure}
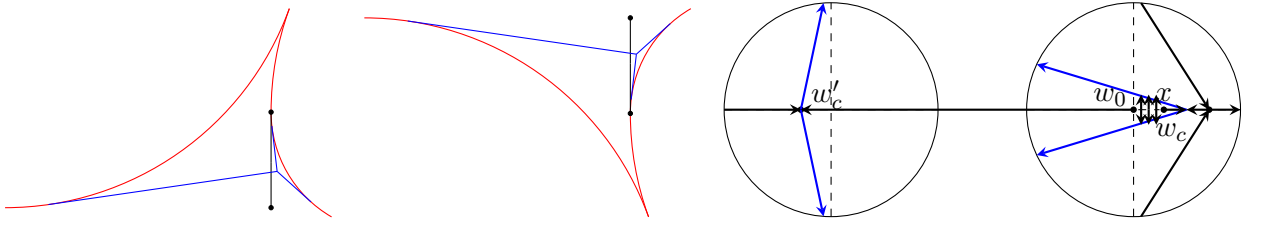

	 \begin{subfigure}{0.28\textwidth}
	%

			\end{subfigure}
		
				\caption{
				In the left and middle panel, $(x,s)$ in a tangent neighborhood is associated with an extra frozen chart; in the right panel, the corresponding local descent paths.
}\label{f:extra_chart2}
\end{figure}

\begin{lemma}\label{l:extra_frozen_action_gap}
Consider one of the following configurations.

\begin{enumerate}
\item
\textbf{Arctic neighborhood.}
Suppose that \((x,s)\in\fN\), where \(\fN\) is an arctic neighborhood.
Then \((x,s)\) is associated with an arctic chart \(\fU\) centered at
\(w_0\), which contains one or two descent critical points. Let \(w_c\) be
any one of them. Suppose, in addition, that the chart collection associated
with \(\fN\) contains a frozen chart with a descent critical point \(w'_c\);
see \Cref{f:extra_chart1}.

\item
\textbf{Tangent neighborhood.}
Suppose that \((x,s)\in\fN\), where \(\fN\) is a tangent neighborhood
satisfying
\[
\fN\cap\fT_A\cap\fT_B\neq\emptyset.
\]
Then \((x,s)\) is associated with a tangent chart \(\fU\) centered at
\(w_0\), corresponding to the extended side separating \(\fT_A\) and
\(\fT_B\). Suppose that \(\fU\) contains a descent critical point \(w_c\)
and that the associated chart collection also contains a frozen chart with
a descent critical point \(w'_c\); see \Cref{f:extra_chart2}.

\item
\textbf{Interface frozen neighborhood.}
Suppose that \((x,s)\in\fN\), where \(\fN\) is an interface frozen
neighborhood satisfying
\[
\fN\cap\fT_A\cap\fT_B\neq\emptyset.
\]
Then \((x,s)\) is associated with a tangent frozen chart \(\fU\) centered
at \(w_0\), corresponding to the extended side separating \(\fT_A\) and
\(\fT_B\). Suppose that \(\fU\) contains a descent critical point \(w_c\)
and that the associated chart collection also contains a frozen chart with
a descent critical point \(w'_c\).
\end{enumerate}

There exists \(\fc'>0\) such that, in each of the above cases,
\begin{align}\label{e:extra_frozen_action_gap}
\Re S(w_c;x,s)
\geq
\Re S(w'_c;x,s)+\fc'.
\end{align}

The analogous statement holds for additional ascent critical points, with
the inequality reversed.
\end{lemma}

\begin{proof}[Proof of \Cref{l:extra_frozen_action_gap}]
Suppose first that \((x,s)\) belongs to an arctic neighborhood \(\fN\).
Then \((x,s)\) is associated with an arctic chart centered at \(w_0\) and
an additional frozen chart containing \(w'_c\). The point \(w'_c\) is a
local minimum of \(\Re S(\,\cdot\,;x,s)\) along the real axis. Without loss
of generality, assume that
$
w'_c>w_0.
$
The other case is analogous. Along the real steepest-descent path from
\(w_0+\fc\) to \(w'_c\), the function
\(\Re S(\,\cdot\,;x,s)\) is decreasing; see the right panel of
\Cref{f:extra_chart1}. Hence
\begin{align}\label{e:extra_arctic_monotonicity}
\Re S(w_0+\fc;x,s)
\geq
\Re S(w'_c;x,s).
\end{align}

On the other hand, \(w_c\) is contained in the arctic chart centered at
\(w_0\), as in \Cref{c:arctic_critical}. By
\eqref{e:arctic_cubic},
\begin{align}\label{e:extra_arctic_local_gap}
\Re S(w_c;x,s)
\geq
\Re S(w_0+\fc;x,s)+\fc'
\end{align}
for some \(\fc'>0\). Combining
\eqref{e:extra_arctic_monotonicity} and
\eqref{e:extra_arctic_local_gap} gives
\[
\Re S(w_c;x,s)
\geq
\Re S(w'_c;x,s)+\fc'.
\]

We next consider a tangent neighborhood. Without loss of generality, assume
that
\[
\nabla H^*=(0,0)\quad\text{on }\fT_A,
\qquad
\nabla H^*=(1,0)\quad\text{on }\fT_B.
\]
The other configurations follow by symmetry. The point \((x,s)\) is
associated with a tangent chart centered at \(w_0\) and an additional
frozen chart containing \(w'_c\). The point \(w'_c\) is again a local
minimum of \(\Re S(\,\cdot\,;x,s)\) along the real axis. Without loss of
generality, assume that
\[
w'_c<w_0.
\]
The other case is analogous. Along the real steepest-descent path from
\(w_0-\fc\) to \(w'_c\), the function
\(\Re S(\,\cdot\,;x,s)\) is decreasing; see
\Cref{f:extra_chart2}. Therefore,
\begin{align}\label{e:extra_tangent_monotonicity}
\Re S(w_0-\fc;x,s)
\geq
\Re S(w'_c;x,s).
\end{align}

The local configuration in the tangent chart is described in
\Cref{f:tangent1}. By assumption, the chart contains a descent critical
point \(w_c\), and it is adapted to the tangent chart centered at \(w_0\),
as in \Cref{c:tangent_critical1}. By \eqref{e:tangent_S},
\begin{align}\label{e:extra_tangent_local_gap}
\Re S(w_c;x,s)
\geq
\Re S(w_0-\fc;x,s)+\fc'.
\end{align}
Combining \eqref{e:extra_tangent_monotonicity} and
\eqref{e:extra_tangent_local_gap} yields
\[
\Re S(w_c;x,s)
\geq
\Re S(w'_c;x,s)+\fc'.
\]

Finally, suppose that \((x,s)\) belongs to an interface frozen
neighborhood. The argument is identical to that for a tangent neighborhood,
with the tangent chart replaced by the tangent frozen chart. So we omit. 
\end{proof}

\subsection{Deformation of critical points}

\begin{lemma}\label{l:motion_descent_critical}
Fix a neighborhood \(\fN_\al\) and one of its associated charts \(\fU\).
Consider the descent critical points associated with
\((x,s)\in\fN_\al\) that are contained in \(\fU\).

As long as \((x,s)\) does not cross the arctic boundary or an extended
side, each nondegenerate descent critical point can be labeled by a
continuous, in fact real-analytic, function
$
w_c=w_c(x,s).
$

Suppose that \((x,s)\) crosses the arctic boundary at a regular point that
is neither a tangent location nor a cusp. Then two critical-point branches
coalesce at the arctic boundary and split after crossing it.

Now suppose that \((x,s)\) crosses the vertical extended side
$
x=b_i
$
transversely from right to left through a point
$
(b_i,s_0)
$
that is neither the tangent location nor the cusp-turning point. Then
exactly one of the following two cases occurs:
\begin{enumerate}
\item
If \(\fU\) is a tangent or tangent frozen chart carrying a local descent
path, then one descent critical point contained in \(\fU\) disappears as
\((x,s)\) crosses the side from right to left.

\item
If \(\fU\) is a cusp-turning chart, then the vertical extended side
separates two adjacent curvilinear triangles. The descent critical point
selected by the corresponding curvilinear triangle changes from the
right-hand branch to the left-hand branch. Writing
\[
w_c^+
:=
\lim_{\varepsilon\downarrow0}
w_c(b_i+\varepsilon,s_0),
\qquad
w_c^-
:=
\lim_{\varepsilon\downarrow0}
w_c(b_i-\varepsilon,s_0),
\]
one has $w_c^-=w_c(b_i,s_0)$ and 
\begin{align}\label{e:critical_jump_action}
\Re S(w_c^+;b_i,s_0)
\geq
\Re S(w_c^-;b_i,s_0).
\end{align}
\end{enumerate}

The corresponding statements for unit-slope and horizontal extended sides
follow by symmetry. For a unit-slope side, one replaces the coordinate
\(x\) by \(x-s\) and crosses from left to right. For a horizontal side, one
crosses from below to above.
\end{lemma}

\begin{proof}
Suppose first that \(w_c\) is a nondegenerate critical point and that
\((x,s)\) lies away from the arctic boundary and all extended sides. Then
\[
w_c-x+s\chi(w_c)=0,
\qquad
1+s\chi'(w_c)\neq0.
\]
The implicit-function theorem therefore gives a unique real-analytic branch
$
w_c=w_c(x,s)
$
in a neighborhood of \((x,s)\). This proves the first statement.

Suppose next that \((x,s)\) crosses the arctic boundary at a regular point
\((x_0,s_0)\), and let \(w_0\) be the corresponding critical point. It
follows from the discussion in \Cref{s:critical_arctic} that, on one side
of the arctic boundary, the two critical points are distinct and real,
whereas on the other side they form a complex-conjugate pair. At the arctic
boundary, they coalesce at \(w_0\).

We now consider a transverse crossing of the vertical extended side
\(x=b_i\) from right to left. By
\eqref{e:geometric_descent_ascent}, a tangent line whose point of tangency
lies to the left of \((x,s)\) corresponds to a descent critical point,
whereas one whose point of tangency lies to the right corresponds to an
ascent critical point.

Suppose first that \(\fU\) is a tangent or tangent frozen chart carrying a
local descent path. It follows from the discussion in
\Cref{s:vertical_tangent,s:vertical_frozen_neighborhood} that, as
\((x,s)\) approaches \(x=b_i\) from the right, the corresponding tangent
line converges to the vertical supporting line, and the corresponding
descent critical point converges to \(b_i\). At the regular point
\((b_i,s_0)\) of the extended side, however, \(b_i\) is spurious because
of the logarithmic cancellation:
\[
S'(b_i;b_i,s_0)\neq0.
\]
After crossing to the left, the corresponding critical point becomes an
ascent critical point. Hence the descent critical-point branch terminates
at the side and disappears.

Suppose next that \(\fU\) is a cusp-turning chart. In this case, the
vertical extended side separates two adjacent curvilinear triangles. There
are three tangent lines from \((b_i,s_0)\): one with tangency point to the
left, the vertical supporting line, and one with tangency point to the
right. They correspond, respectively, to a descent critical point
\(\xi_c^{\rm d}\), the spurious critical point \(b_i\), and an ascent
critical point \(\xi_c^{\rm a}\). Moreover,
\begin{align}\label{e:order_relation}
\Re S(\xi_c^{\rm d};b_i,s_0)
\leq
\Re S(b_i;b_i,s_0)
\leq
\Re S(\xi_c^{\rm a};b_i,s_0).
\end{align}

It follows from the discussion in \Cref{s:vertical_tangent} that, as
\((x,s)\) approaches \(x=b_i\) from the right, the selected descent
critical point is the one whose point of tangency lies on the portion of
the arctic boundary contained in the curvilinear triangle on the right.
Thus,
\[
w_c^+=b_i.
\]
After crossing to the left, the selected descent critical point becomes
\[
w_c^-=\xi_c^{\rm d}.
\]
The inequality \eqref{e:critical_jump_action} now follows from
\eqref{e:order_relation}.

The unit-slope and horizontal cases follow by the corresponding symmetries,
using \(x-s\) as the transverse coordinate in the unit-slope case and the
local coordinate at infinity in the horizontal case.
\end{proof}

\begin{lemma}\label{l:critical_value_continuity}
Fix a neighborhood \(\fN_\al\) and one of its associated charts \(\fU\).
There exists \(\delta_0>0\) such that the following holds. Let
\[
0<\delta\leq\delta_0,
\qquad
(x,s),(x',s')\in\fN_\al,
\qquad
\|(x,s)-(x',s')\|_2\leq\delta.
\]
Suppose that \((x,s)\) and \((x',s')\) can be joined by a path contained in
\(\fN_\al\) that does not intersect any extended side, and that
\(w_c(x,s)\) and \(w_c(x',s')\) lie on the same continuous local branch of
associated critical points in \(\fU\). Then
\begin{align}\label{e:critical_value_continuity}
\left|
S(w_c(x,s);x,s)
-
S(w_c(x',s');x',s')
\right|
\lesssim
\delta\ln\frac{2}{\delta}.
\end{align}
If the  critical points are uniformly bounded away
from the tangent locations, the error can be improved to $\delta$.
\end{lemma}

\begin{proof}
By \eqref{e:critical_value_derivatives_gff}, along the critical-point branch
and for \((u,v)\) on the path connecting \((x',s')\) to \((x,s)\), one has
\begin{align}\label{e:critical_value_differential}
\rd S(w_c(u,v);u,v)
=
-\ln f(w_c(u,v))\,\rd u
-\ln\bigl(1-\chi(w_c(u,v))\bigr)\,\rd v.
\end{align}

We will use the elementary bound
\begin{align}\label{e:ulog_modulus}
\int_a^b
\left(
1+\left|\ln|u-u_0|\right|
\right)
|\rd u|
\lesssim
|a-b|\ln\frac{2}{|a-b|}.
\end{align}

Suppose first that the critical-point branch remains away from the tangent
locations. Then \(\chi(w_c)\) is bounded away from \(0\), \(1\), and
\(\infty\). Hence both logarithms in
\eqref{e:critical_value_differential} are uniformly bounded, and
\[
\left|
S(w_c(x,s);x,s)
-
S(w_c(x',s');x',s')
\right|
\lesssim
\|(x,s)-(x',s')\|_2
\lesssim
\delta.
\]

It remains to consider the three types of tangent charts.

\smallskip
\noindent
\emph{Vertical tangent chart.}
Suppose that the chart is centered at \(b_i\), corresponding to a vertical
tangent point. By
\eqref{e:wc_vertical_tangent}, \eqref{e:wc_vertical_cusp}, and
\eqref{e:wcest}, along the critical point branch,
\[
|\chi(w_c(u,v))|
\asymp
|w_c(u,v)-b_i|
\gtrsim
|u-b_i|.
\]
Since \(f=\chi/(1-\chi)\), it follows that
\begin{align}\label{e:vertical_log_bound}
|\ln f(w_c(u,v))|
&\lesssim
1+\left|\ln|u-b_i|\right|,
\qquad
\left|\ln\bigl(1-\chi(w_c(u,v))\bigr)\right|
&\lesssim 1.
\end{align}
Integrating \eqref{e:critical_value_differential} along the path from $(x,s)$ to $(x',s')$, and using
\eqref{e:ulog_modulus}, we obtain
\begin{align}
&\left|
S(w_c(x,s);x,s)
-
S(w_c(x',s');x',s')
\right|
\lesssim
\int_{x'}^x
\left(
1+\left|\ln|u-b_i|\right|
\right)
|\rd u|
+
|s-s'|
\lesssim
\delta\ln\frac{2}{\delta}.
\end{align}

\smallskip
\noindent
\emph{Unit-slope tangent chart.}
Suppose that the chart is centered at \(a_i\), corresponding to a
unit-slope tangent point. By
\eqref{e:wc_unit_slope_tangent}, \eqref{e:wc_unit_slope_cusp}, and
\eqref{e:wc_unit_slope_frozen},
\[
|1-\chi(w_c(u,v))|
\asymp
|w_c(u,v)-a_i|
\gtrsim
|u-v-a_i|.
\]
Since \(f=\chi/(1-\chi)\), it follows that
\begin{align}\label{e:unit_log_bound}
|\ln f(w_c(u,v))|
&\lesssim
1+\left|\ln|u-v-a_i|\right|,
\qquad
|\ln\chi(w_c(u,v))|
&\lesssim 1.
\end{align}
Moreover, \eqref{e:critical_value_differential} can be rewritten as
\begin{align}
\rd S(w_c(u,v);u,v)
=
-\ln f(w_c(u,v))\,\rd(u-v)
-\ln\chi(w_c(u,v))\,\rd v.
\end{align}
Thus, integrating along the path from $(x,s)$ to $(x',s')$, and applying
\eqref{e:ulog_modulus}, gives
\begin{align*}
&\left|
S(w_c(x,s);x,s)
-
S(w_c(x',s');x',s')
\right|
\lesssim
\int_{x'-s'}^{x-s}
\left(
1+\left|\ln|u-v-a_i|\right|
\right)
|\rd (u-v)|
+
|s-s'|
\lesssim
\delta\ln\frac{2}{\delta}.
\end{align*}

\smallskip
\noindent
\emph{Horizontal tangent chart.}
Suppose that the chart is centered at \(\infty_i\), corresponding to a
horizontal tangent location \((x_0,s_0)\). By
\eqref{e:wc_horizontal_tangent}, \eqref{e:wc_horizontal_cusp}, and
\eqref{e:twc_bound},
\[
|\chi(w_c(u,v))|
\asymp
|w_c(u,v)|
\lesssim
\frac{1}{|v-s_0|}.
\]
It follows that
\begin{align}\label{e:horizontal_log_bound}
\left|
\ln\bigl(1-\chi(w_c(u,v))\bigr)
\right|
\lesssim
1+\left|\ln|v-s_0|\right|.
\end{align}
Since \(f=\chi/(1-\chi)\), the quantity \(\ln f(w_c(u,v))\) remains
uniformly bounded. Therefore, integrating along the path from $(x,s)$ to $(x',s')$,
we obtain
\begin{align}
&\left|
S(w_c(x,s);x,s)
-
S(w_c(x',s');x',s')
\right|
\lesssim
\int_{s'}^s
\left(
1+\left|\ln|v-s_0|\right|
\right)
|\rd v|
+
|x-x'|
\lesssim
\delta\ln\frac{2}{\delta}.
\end{align}
\end{proof}

\section{Estimates of the Integrand}\label{s:Integrand_est}

In this section, we collect several asymptotic estimates for the integrands of the single-contour integral \eqref{e:single_term} and the double-contour integral \eqref{e:all_term}.

\subsection{Single-contour integral}

In this section we collect some basic estimates on the asymptotics of the integrand of the single-contour integral \eqref{e:single_term}.

\begin{lemma}\label{l:PQ_bound}
Given  \(z\in \bC\), define
\begin{align}
\delta:=\operatorname{dist}\left(z,[\min\{x,y\}, \max\{x,y\}]\cup [\min\{x-s,y-t\}, \max\{x-s,y-t\}]\right).
\end{align}
If $\delta\gg n^{-1}$, then
\begin{align}\label{e:PQess}
P_{ns}(nz,nx) Q_{nt}(nz,ny)
=
\frac{\sqrt{st}\,
e^{n (S(z;x,s)-S(z;y,t))+\OO(1/(\delta n))}}
{n\sqrt{x-z}\sqrt{z-(x-s)}\sqrt{y-z}\sqrt{z-(y-t)}}.
\end{align}

If \(z\) is in a small neighborhood of
$
[\min\{x,y\}, \max\{x,y\}],
$
then $\delta=\dist(z, [\min\{x,y\}, \max\{x,y\}])$, and
the following estimates hold. If \(x<y\), then
\begin{align}\label{e:PQubb}
|P_{ns}(nz,nx) Q_{nt}(nz,ny)|
\leq
\frac{Cn e^{n \Re[S(z;x,s)-S(z;y,t)]}}
{\sqrt{(n|x-z|+1)(n|z-(x-s)|+1)(n|y-z|+1)(n|z-(y-t)|+1)}}.
\end{align}
If \(x\ge y\) and \(\delta>0\), then
\begin{align}\label{e:PQubb2}
|P_{ns}(nz,nx) Q_{nt}(nz,ny)|
\leq
\frac{C\min\{1,n\delta\}^{-1} n e^{n \Re[S(z;x,s)-S(z;y,t)]}}
{\sqrt{(n|x-z|+1)(n|z-(x-s)|+1)(n|y-z|+1)(n|z-(y-t)|+1)}}.
\end{align}

If \(z\) is in a small neighborhood of
$
[\min\{x-s,y-t\}, \max\{x-s,y-t\}],
$
then $\delta=\dist(z, [\min\{x-s,y-t\}, \max\{x-s,y-t\}])$, and the following estimates hold. If \(y-t<x-s\), then \eqref{e:PQubb} holds.
If \(y-t\ge x-s\) and \(\delta>0\), then \eqref{e:PQubb2} holds.
\end{lemma}

We need the following product-integral comparison lemma. Its proof is postponed to the end of this section. 
\begin{lemma}\label{l:product_integral}
For any positive integers \(n,\ell\geq 1\), and any complex number
\(z\in \bC\setminus [1/n,\ell/n]\), let
\[
\delta:=\operatorname{dist}\left(z,[{1}/{n},{\ell}/{n}]\right)>0.
\]
If \(\delta\gg n^{-1}\), we have
\begin{align}\label{e:asymplog}
\prod_{j=1}^{\ell}\left(z-\frac jn\right)
\exp\left(
-n\int_{1/(2n)}^{\ell/n+1/(2n)}\log(z-x)\,\rd x
\right)
=
1+\OO\left(\frac{1}{n\delta}\right).
\end{align}
More generally, there exists a universal constant \(C>0\) such that 
\begin{align}\label{e:sum_approx_integral}
C^{-1}\min\{1,n\delta\}
\le
\left|
\prod_{j=1}^{\ell}\left(z-\frac jn\right)
\exp\left(
-n\int_{1/(2n)}^{\ell/n+1/(2n)}\log(z-x)\,\rd x
\right)
\right|
\le C.
\end{align}
Here \(\log\) denotes the principal logarithm, with argument in
\((-\pi,\pi]\). If \(z\) lies in the interval of integration, the integral
is understood as an improper integral.
\end{lemma}

\begin{proof}[Proof of \Cref{l:PQ_bound}]
By definition,
\begin{align}\label{e:PQ_gamma_decomposition}
P_{ns}(nz,nx)Q_{nt}(nz,ny)
&=
\frac{\Gamma(ns+1)}{\Gamma(nt)}
\frac{\Gamma(n(y-z))}{\Gamma(n(x-z)+1)}
\frac{\Gamma(n(z-(y-t)))}{\Gamma(n(z-(x-s))+1)}.
\end{align}

We will only prove the case in which $z$ is closer to $[\min\{x,y\}, \max\{x,y\}]$, than to $[\min\{x-s,y-t\}, \max\{x-s,y-t\}]$. In this case $\delta=\dist(z, [\min\{x,y\}, \max\{x,y\}])$.

Suppose first that \(x<y\). By the upper bound in
\Cref{l:product_integral},
\begin{align}
\left|
\frac{\Gamma(n(y-z))}{\Gamma(n(x-z)+1)}
\right|
&\leq
C n^{n(y-x)-1}
\exp\left(
n\int_{x+\frac1{2n}}^{y-\frac1{2n}}
\log|u-z|\,\rd u
\right)\\
&=
C n^{n(y-x)-1}
\exp\left(
n\Re\int_x^y\log(u-z)\,\rd u
\right)\\
&\quad\times
\exp\left(
-n\int_x^{x+\frac1{2n}}\log|u-z|\,\rd u
\right)
\exp\left(
-n\int_{y-\frac1{2n}}^y\log|u-z|\,\rd u
\right).
\end{align}
Using the elementary estimates
\begin{align}\label{e:half_cell_log_bound}
\max\left\{
\exp\left(
-n\int_{r-\frac1{2n}}^r\log|u-z|\,\rd u
\right),
\exp\left(
-n\int_r^{r+\frac1{2n}}\log|u-z|\,\rd u
\right)
\right\}
\leq
\frac{C\sqrt n}{\sqrt{n|r-z|+1}}.
\end{align}
we obtain
\begin{align}\label{e:first_gamma_direct}
\left|
\frac{\Gamma(n(y-z))}{\Gamma(n(x-z)+1)}
\right|
\leq
\frac{
C n^{n(y-x)}
\exp\left(
n\Re\int_x^y\log(u-z)\,\rd u
\right)}
{\sqrt{(n|x-z|+1)(n|y-z|+1)}}.
\end{align}

Suppose next that \(x\geq y\). The lower bound in \Cref{l:product_integral} gives
\begin{align}
\left|
\frac{\Gamma(n(y-z))}{\Gamma(n(x-z)+1)}
\right|
&\leq
\frac{C n^{-n(x-y)-1}}{\min\{1,n\delta\}}
\exp\left(
-n\int_{y-\frac1{2n}}^{x+\frac1{2n}}
\log|u-z|\,\rd u
\right)\\
&=
\frac{C n^{-n(x-y)-1}}{\min\{1,n\delta\}}
\exp\left(
n\Re\int_x^y\log(u-z)\,\rd u
\right)\\
&\quad\times
\exp\left(
-n\int_{y-\frac1{2n}}^y\log|u-z|\,\rd u
\right)
\exp\left(
-n\int_x^{x+\frac1{2n}}\log|u-z|\,\rd u
\right).
\end{align}
Thus,
\begin{align}\label{e:first_gamma_inverse}
\left|
\frac{\Gamma(n(y-z))}{\Gamma(n(x-z)+1)}
\right|
\leq
\frac{
C n^{n(y-x)}
\exp\left(
n\Re\int_x^y\log(u-z)\,\rd u
\right)}
{\min\{1,n\delta\}
\sqrt{(n|x-z|+1)(n|y-z|+1)}}.
\end{align}
If $\delta=\dist(z, [\min\{x,y\}, \max\{x,y\}])\gg n^{-1}$, then \Cref{l:product_integral} gives
\begin{align}\label{e:far_gamma_asymptotic}
\frac{\Gamma(n(y-z))}{\Gamma(n(x-z)+1)}
=
\frac{
 n^{n(y-x)-1}
\exp\left(
n\int_x^y\log(u-z)\,\rd u
+\OO(1/(\delta n))\right)}
{\sqrt{x-z}\sqrt{y-z}}.
\end{align}

Stirling's formula gives
\begin{align}\label{e:second_gamma_asymptotic}
\frac{\Gamma(n(z-(y-t)))}
{\Gamma(n(z-(x-s))+1)}
=
\frac{
n^{n((x-s)-(y-t))-1}
\exp\left(
n\int_{y-t}^{x-s}\log(z-u)\,\rd u
+\OO\left(\frac{1}{n}\right)
\right)}
{\sqrt{(z-(x-s))(z-(y-t))}}.
\end{align}
and 
\begin{align}\label{e:gamma_st_asymptotic}
\frac{\Gamma(ns+1)}{\Gamma(nt)}
=
\sqrt{st}\,
n^{n(s-t)+1}
\exp\left(
n(s\log s-t\log t-s+t)
+\OO\left(\frac1n\right)
\right).
\end{align}

By the definition of the tiling action \(S\) in \eqref{e:def_action_C},
\begin{align*}
S(z;x,s)-S(z;y,t)=s\log s-t\log t-s+t
+\int_x^y\log(u-z)\,\rd u
+\int_{y-t}^{x-s}\log(z-u)\,\rd u.
\end{align*}
Combining \eqref{e:PQ_gamma_decomposition}, \eqref{e:far_gamma_asymptotic}, \eqref{e:second_gamma_asymptotic} and \eqref{e:gamma_st_asymptotic}
therefore yields
\begin{align*}
P_{ns}(nz,nx)Q_{nt}(nz,ny)
&=
\frac{\sqrt{st}\,
\exp\left(
n(S(z;x,s)-S(z;y,t))
+\OO\left(\frac{1}{n\delta}\right)
\right)}
{n\sqrt{x-z}\sqrt{z-(x-s)}
\sqrt{y-z}\sqrt{z-(y-t)}},
\end{align*}
which is  \eqref{e:PQess}.

Combining \eqref{e:PQ_gamma_decomposition},
\eqref{e:first_gamma_direct}, \eqref{e:second_gamma_asymptotic} and \eqref{e:gamma_st_asymptotic} gives
\begin{align*}
|P_{ns}(nz,nx)Q_{nt}(nz,ny)|
&\leq
\frac{
Cn e^{n\Re[S(z;x,s)-S(z;y,t)]}}
{\sqrt{
(n|x-z|+1)(n|z-(x-s)|+1)
(n|y-z|+1)(n|z-(y-t)|+1)
}},
\end{align*}
which is \eqref{e:PQubb}.

Similarly, combining \eqref{e:PQ_gamma_decomposition},
\eqref{e:first_gamma_inverse}, \eqref{e:second_gamma_asymptotic} and \eqref{e:gamma_st_asymptotic} gives
\begin{align*}
|P_{ns}(nz,nx)Q_{nt}(nz,ny)|
&\leq
\frac{
C\min\{1,n\delta\}^{-1}
n e^{n\Re[S(z;x,s)-S(z;y,t)]}}
{\sqrt{
(n|x-z|+1)(n|z-(x-s)|+1)
(n|y-z|+1)(n|z-(y-t)|+1)
}},
\end{align*}
which is \eqref{e:PQubb2}.

\end{proof}

\begin{proof}[Proof of \Cref{l:product_integral}]
The midpoint expansion gives, uniformly in \(j\),
\begin{align}
n\int_{(j-\frac12)/n}^{(j+\frac12)/n} \log(z-x)\,\rd x
&=
\log\left(z-\frac jn\right)
+
\OO\left(\sup_{x\in[(j-\frac12)/n,(j+\frac12)/n]}
\frac{1}{|z-x|^2n^2}\right)\\
&=
\log\left(z-\frac jn\right)
+\OO\left(\frac{1}{|z-j/n|^2n^2}\right),
\end{align}
provided that \(\delta\gg n^{-1}\). Indeed, since
$|z- j/n|\geq\delta\gg n^{-1}$,
every \(x\in[(j-\frac12)/n,(j+\frac12)/n]\) satisfies
\[
|z-x|
\geq
\left|z-\frac jn\right|-\frac{1}{2n}
\asymp
\left|z-\frac jn\right|.
\]
Summing the midpoint expansions over \(1\leq j\leq\ell\) gives
\[
\sum_{j=1}^{\ell}\log\left(z-\frac jn\right)
-
n\int_{1/(2n)}^{\ell/n+1/(2n)}\log(z-x)\,\rd x
=
\OO\left(\frac{1}{n\delta}\right).
\]
The claim \eqref{e:asymplog} follows by exponentiating.

To prove \eqref{e:sum_approx_integral}, we first rescale. Put \(\zeta=nz\), then 
\[
\left|
\prod_{j=1}^{\ell}\left(z-\frac jn\right)
\exp\left(
-n\int_{1/2n}^{\ell/n+1/2n}\log(z-x)\,\rd x
\right)
\right|=
\frac{\left|\prod_{j=1}^{\ell}(\zeta-j)\right|}
{\exp\left(\int_{1/2}^{\ell+1/2}\log|\zeta-y|\,\rd y\right)}=:R_\ell(\zeta).
\]

For \(u\in\mathbb C\setminus\{0\}\), define
\[
q(u):=
\log|u|-\int_{-1/2}^{1/2}\log|u-y|\,\rd y.
\]
Then, by decomposing the integral into unit cells,
\[
\log R_\ell(\zeta)
=
\sum_{j=1}^{\ell} q(\zeta-j).
\]

By Taylor expansion, there exists some \(C>0\) such that for any
\(|u|\geq 1/2\),
\begin{align}\label{e:qbound}
\frac{C}{|u|^2}\geq q(u)\geq -\frac{C}{|u|^2}.
\end{align}

For every \(u\in\mathbb C\), we have that  
\[
-\log 2-1=\int_{-1/2}^{1/2}\log|y|\,\rd y
\leq \int_{-1/2}^{1/2}\log|u-y|\,\rd y
\le
\log\left(|u|+\frac12\right).
\]
Hence, for every \(u\in\mathbb C\setminus\{0\}\),
\begin{align}\label{e:lowup}
e^{q(u)}
\ge
\frac{|u|}{|u|+1/2}\geq \frac23\min\{1,|u|\}.
\end{align}
Moreover, if \(|u|\leq 1/2\), then
\begin{align}\label{e:ubound}
e^{q(u)}\leq 2e|u|\leq e.
\end{align}

Choose \(j_0\in\{1,\dots,\ell\}\) such that
\[
|\zeta-j_0|
=
\min_{1\le j\le \ell}|\zeta-j|
\geq n\delta. 
\]
Then, since \(j_0\) is closest to \(\zeta\) and \(|j-j_0|\geq 1\), the remaining indices \(j\ne j_0\) satisfy
\[
|\zeta-j|\geq 1/2.
\]
Moreover, from \eqref{e:lowup} we have
\[
e^{q(\zeta-j_0)}
\ge
\frac23\min\{1,n\delta\}.
\]

For the upper bound, using the above estimates \eqref{e:qbound} and \eqref{e:ubound}, we have
\begin{align}
\log R_\ell(\zeta)=\sum_{j=1}^{\ell}
q(\zeta-j)&\leq e+\sum_{\substack{1\le j\le \ell\\ |\zeta-j|\ge1/2}}
q(\zeta-j) \leq e+
\sum_{\substack{1\le j\le \ell\\ |\zeta-j|\ge1/2}}
\frac{C}{|\zeta-j|^2} \\
&\leq e+2C\left(\frac{1}{(1/2)^2}+\frac{1}{(3/2)^2}+\cdots\right)
\leq e+16C.
\end{align}

For the lower bound we have 
\begin{align}
R_\ell(\zeta)
&=
e^{q(\zeta-j_0)}
\exp\left(\sum_{\substack{1\le j\le \ell\\ j\ne j_0}}q(\zeta-j)\right) \\
&\geq
\frac23\min\{1,n\delta\}
\exp\left(
-\sum_{\substack{1\le j\le \ell\\ j\ne j_0}}
\frac{C}{|\zeta-j|^2}
\right) \geq
\frac23\min\{1,n\delta\}e^{-16C}.
\end{align}

Combining the above estimates,  we obtain
\[
e^{e+16C}\geq R_\ell(\zeta)
\ge
\left(\frac23\min\{1,n\delta\}\right)
e^{-16C}.
\]
Changing the value of the universal constant \(C\) proves the lemma.
\end{proof}

\subsection{Double-contour integral}
In this section we collect some basic estimates on the asymptotics of the integrand of the double-contour integral \eqref{e:all_term}.

\begin{lemma}\label{l:PIi_bound}
Given any complex number \(w\in \bC\) and denote
\begin{align}\label{e:PIest}
\delta:=&\operatorname{dist}\left(w,[\min\{b_i+1/(2n),x\}, \max\{b_i+1/2n,x\}]\right.\\
&\qquad\qquad \left.\cup [\min\{a_i-1/(2n),x-s\}, \max\{a_i-1/(2n),x-s\}]\right).
\end{align}
If $\delta\gg n^{-1}$, then
\begin{align}
P_{ns}(nw,nx) I_i(w)= \frac{\sqrt{s}}{\sqrt{2\pi n}\sqrt{x-w}\sqrt{w-(x-s)}} e^{n S(w;x,s)+\OO(1/(\delta n))}.
\end{align}

If \(w\) is in a small neighborhood of
$
[\min\{b_i+1/(2n),x\}, \max\{b_i+1/2n,x\}],
$
then $\delta=\dist(w, [\min\{b_i+1/(2n),x\}, \max\{b_i+1/2n,x\}])$, and
the following estimates hold.
If \(x\leq b_i-1/(2n)\), then
\begin{align}\label{e:PIupbb1}
|P_{ns}(nw,nx) I_i(w)|\leq \frac{C\sqrt{n}}{\sqrt{(n|x-w|+1)(n|w-(x-s)|+1)}} e^{n \Re[S(w;x,s)]}.
\end{align}
If \(x\geq b_i+1/(2n)\) and \(\delta>0\), then
\begin{align}\label{e:PIupbb2}
|P_{ns}(nw,nx) I_i(w)|\leq  \frac{C\sqrt{n}}{\min\{1,n\delta\}\sqrt{(n|x-w|+1)(n|w-(x-s)|+1)}}  e^{n \Re[S(w;x,s)]}.
\end{align}

If \(w\) is in a small neighborhood of
$
[\min\{a_i-1/(2n),x-s\}, \max\{a_i-1/(2n),x-s\}],
$
then $\delta=\dist(w, [\min\{a_i-1/(2n),x-s\}, \max\{a_i-1/(2n),x-s\}])$, and
the following estimates hold.
If \(x-s\geq a_i+1/(2n)\), then \eqref{e:PIupbb1} holds.
If \(x-s\leq a_i-1/(2n)\) and \(\delta>0\), then \eqref{e:PIupbb2} holds.
\end{lemma}

\begin{lemma}\label{l:QIi_bound}
Given any complex number \(z\in \bC\) and denote
\begin{align}
\delta:=&\operatorname{dist}\left(z,[\min\{b_i-1/(2n),y\}, \max\{b_i-1/2n,y\}]\right.\\
&\qquad \qquad \left.\cup [\min\{a_i+1/(2n),y-t\}, \max\{a_i+1/(2n),y-t\}]\right).
\end{align}
If $\delta\gg n^{-1}$, then
\begin{align}
Q_{nt}(nz,ny) I^{-1}_i(z)= \frac{\sqrt{2\pi t}}{\sqrt{ n}\sqrt{y-z}\sqrt{z-(y-t)}} e^{-n S(z;y,t)+\OO(1/(\delta n))}.
\end{align}

If \(z\) is in a small neighborhood of
$
[\min\{b_i-1/(2n),y\}, \max\{b_i-1/2n,y\}],
$
then $\delta=\dist(z, [\min\{b_i-1/(2n),y\}, \max\{b_i-1/2n,y\}])$, and
the following estimates hold.
If \(y\geq b_i+1/(2n)\), then
\begin{align}\label{e:QIupbb1}
|Q_{nt}(nz,ny) I^{-1}_i(z)|\leq \frac{C\sqrt{n}}{\sqrt{(n|y-z|+1)(n|z-(y-t)|+1)}} e^{-n \Re[S(z;y,t)]}.
\end{align}
If \(y\leq b_i-1/(2n)\) and \(\delta>0\), then
\begin{align}\label{e:QIupbb2}
|Q_{nt}(nz,ny) I^{-1}_i(z)|\leq  \frac{C\sqrt{n}}{\min\{1,n\delta\}\sqrt{(n|y-z|+1)(n|z-(y-t)|+1)}}e^{-n \Re[S(z;y,t)]}.
\end{align}

If \(z\) is in a small neighborhood of
$
[\min\{a_i+1/(2n),y-t\}, \max\{a_i+1/(2n),y-t\}],
$
then $\delta=\dist(w, [\min\{a_i+1/(2n),y-t\}, \max\{a_i+1/(2n),y-t\}])$, and
the following estimates hold.
If \(y-t\leq a_i-1/(2n)\), then \eqref{e:QIupbb1} holds.
If \(y-t\geq a_i+1/(2n)\) and \(\dist(z, [y, b_i-1/(2n)])>0\), then \eqref{e:QIupbb2} holds.
\end{lemma}

\begin{lemma}\label{l:arctic_S}
Let $\fU$ be a liquid, arctic or cusp chart, and let
$w,z\in \fU$. Then, uniformly for $w,z\in \fU$, we have
\begin{align}\label{e:PI+_bound}
P_{ns}(nw,nx) I_+(w)
=
\frac{\sqrt{s}}{\sqrt{2\pi n}\sqrt{x-w}\sqrt{w-(x-s)}}
\exp\left(n S(w;x,s)+\OO(1/n)\right),
\end{align}
and
\begin{align}\label{e:QI-_bound}
Q_{ns}(nz,nx) I_-^{-1}(z)
=
\frac{\sqrt{2\pi t}}{\sqrt{n}\sqrt{x-z}\sqrt{z-(x-s)}}
\exp\left(-n S(z;x,s)+\OO(1/n)\right).
\end{align}
\end{lemma}

\begin{lemma}\label{l:vertical_tangent_S}
%Adopt the notation of either
%\Cref{c:tangent_critical1,l:vertical_tangent_steepest}
%or
%\Cref{c:cusp_turning_critical1,c:vertical_cusp_steepest}.
Let \(\fU\) be a vertical tangent, cusp-turning or frozen chart centered at \(w_0=b_i\).
Then \eqref{e:PI+_bound} and \eqref{e:QI-_bound} hold uniformly for
\(w,z\in\fU\) such that \(w\) and \(z\) are bounded away from
\[
[\min\{b_i,x\},\max\{b_i,x\}]
\qquad\text{and}\qquad
[\min\{b_i,y\},\max\{b_i,y\}],
\]
respectively.

Moreover, if $(x,s)$ is adapted to $\fU$, then for every \(w\in\sfD^{\rm d}(w_c)\),
\begin{align}\label{e:PIupbb1_copy}
|P_{ns}(nw,nx)I_i(w)|
\leq
\frac{C\sqrt{n}}
{\sqrt{n|x-w|+1}\sqrt{n|w-(x-s)|+1}}
e^{n\Re S(w;x,s)}.
\end{align}
If, in addition, \(x\geq b_i+1/(2n)\), then
\begin{align}\label{e:Dist1}
\dist\bigl(\sfD^{\rm d}(w_c),[b_i,x]\bigr)
\gtrsim \frac{1}{n}.
\end{align}

For every \(z\in\sfD^{\rm a}(z_c)\),
\begin{align}\label{e:QIupbb1_copy}
|Q_{nt}(nz,ny)I_i^{-1}(z)|
\leq
\frac{C\sqrt{n}}
{\sqrt{n|y-z|+1}\sqrt{n|z-(y-t)|+1}}
e^{-n\Re S(z;y,t)}.
\end{align}
If, in addition, \(y\leq b_i-1/(2n)\), then
\begin{align}\label{e:Dist2}
\dist\bigl(\sfD^{\rm a}(z_c),[y,b_i]\bigr)
\gtrsim \frac{1}{n}.
\end{align}
\end{lemma}

\begin{lemma}\label{l:unit_slope_tangent_S}
Let \(\fU\) be a unit-slope tangent, cusp-turning or frozen chart centered at $w_0=a_i$. Then \eqref{e:PI+_bound} and \eqref{e:QI-_bound} hold uniformly for
$w,z\in\fU$ bounded away from
$
[\min\{b_i,x-s\},\max\{b_i,x-s\}]
$ and $
[\min\{b_i,y-t\},\max\{b_i,y-t\}],
$
respectively.

Moreover, \eqref{e:PIupbb1_copy} holds for every
\(w\in\sfD^{\rm d}(w_c)\). If, in addition,
\(x-s\leq a_i-1/(2n)\), then 
\begin{align}
\label{e:Dist12} 
\dist(\sfD^{\rm d}(w_c), [x-s, a_i])\gtrsim \frac{1}{n}.
\end{align}
Likewise, \eqref{e:QIupbb1_copy} holds for every
\(z\in\sfD^{\rm a}(z_c)\). If, in addition,
\(y-t\geq a_i+1/(2n)\), then 
\begin{align}
\label{e:Dist22} 
\dist(\sfD^{\rm a}(z_c), [ a_i, y-t])\gtrsim \frac{1}{n}.
\end{align}
\end{lemma}

\begin{lemma}\label{l:horizontal_tangent_S}
Let
$\wt\fU$ be a horizontal tangent, cusp-turning or frozen chart centered at $0$ and associated with
the horizontal tangency point $(x_0,s_0)$. For
$\wt w,\wt z\in\wt\fU\setminus\{0\}$, set
\[
w=x_0-\frac{1}{\wt w},
\qquad
z=x_0-\frac{1}{\wt z}.
\]
Then \eqref{e:PI+_bound} and \eqref{e:QI-_bound} hold uniformly for
$\wt w,\wt z\in\wt\fU\setminus\{0\}$, with $w$ and $z$ defined as above.
\end{lemma}

\begin{proof}[Proof of \Cref{l:PIi_bound}]
Combining \eqref{e:defI2} and \eqref{e:Pterm}, we have
\begin{align}\begin{split}\label{e:Pexp}
&P_{ns}(nw,nx) I_i(w)
=\Gamma(ns+1)\frac{\Gamma(n(b_i-w)+1/2)}{\Gamma(n(x-w)+1)}\frac{\Gamma(n(w-a_i)+1/2)}{\Gamma(n(w+s-x)+1)}e^{-n\int_0^w \ln f(u)\,\rd u}\,\\
&\times
\exp\Bigl(-n\bigl((b_i-w)\ln(b_i-w)+(w-a_i)\ln(w-a_i)\bigr)-n(b_i-a_i)\ln(n/e)-\ln(2\pi)\Bigr).
\end{split}\end{align}

We will only prove the case in which $w$ is closer to $[\min\{b_i+1/(2n),x\}, \max\{b_i+1/(2n),x\}]$, than to $[\min\{a_i-1/(2n),x-s\}, \max\{a_i-1/(2n),x-s\}]$. In this case $\delta=\dist(w, [\min\{b_i+1/(2n),x\}, \max\{b_i+1/(2n),x\}])$.

We now estimate the two gamma ratios. First consider
\[
\frac{\Gamma(n(b_i-w)+1/2)}{\Gamma(n(x-w)+1)}.
\]
If \(x\le b_i-1/(2n)\), then this quotient is a finite product. Applying \Cref{l:product_integral} and gives
\begin{align}
\left|
\frac{\Gamma(n(b_i-w)+1/2)}{\Gamma(n(x-w)+1)}
\right|
%\le
%\frac{C}{\sqrt{n|x-w|}}\exp\left(
%n\int_{x}^{b_i}\log|n(u-w)|\,\rd u
%\right).
&\leq C n^{n(b_i-x)-1/2}\exp\left(
n\Re\int_{x+1/2n}^{b_i}\log(u-w)\,\rd u
\right)\\
&=
C n^{n(b_i-x)-1/2} \exp\left(
n\Re \int_{x}^{b_i}\log (u-w)\,\rd u
\right)
e^{-n  \int^{x+1/2n}_{x}\log |u-w|\,\rd u }
\end{align}
Using the identity
\[
\int_x^{b_i}\log(u-w)\,\rd u
=
(b_i-w)\log(b_i-w)-(b_i-w)
-(x-w)\log(x-w)+(x-w),
\]
and the upper bound \eqref{e:half_cell_log_bound}
\begin{align}
e^{-n  \int_{x}^{x+1/2n}\log |u-w|\,\rd u }\leq \frac{C}{\sqrt{|x-w|+1/n}}=\frac{C\sqrt{n}}{\sqrt{n|x-w|+1}}
\end{align}
we obtain
\begin{align}\label{e:Gammag1}
\left|
\frac{\Gamma(n(b_i-w)+1/2)}{\Gamma(n(x-w)+1)}
e^{-n(b_i-w)\log(b_i-w)}
\right|
\le
\frac{C\,n^{n(b_i-x)}e^{-n(b_i-x)}}{\sqrt{n|x-w|+1}}
e^{-n\Re[(x-w)\log(x-w)]}.
\end{align}
On the other hand, if \(x\ge b_i+1/(2n)\), the same argument applied to the
inverse product gives
\begin{align}\label{e:Gammag2}
\left|
\frac{\Gamma(n(b_i-w)+1/2)}{\Gamma(n(x-w)+1)}
e^{-n(b_i-w)\log(b_i-w)}
\right|
\le
\frac{Cn^{n(b_i-x)}e^{-n(b_i-x)}}{\min\{1,n\dist(w,[b_i+1/(2n),x])\}}
\frac{e^{-n\Re[(x-w)\log(x-w)]}}{\sqrt{n|x-w|+1}}.
\end{align}

If $\delta=\dist(w, [\min\{b_i+1/(2n),x\}, \max\{b_i+1/(2n),x\}])\gg n^{-1}$, by \Cref{l:product_integral} the following estimates hold
\begin{align}\label{e:Gammag25}
\frac{\Gamma(n(b_i-w)+1/2)}{\Gamma(n(x-w)+1)}
e^{-n(b_i-w)\log(b_i-w)}
=
\frac{n^{n(b_i-x)}e^{-n(b_i-x)}}{\sqrt{n(x-w)}}
e^{-n(x-w)\log(x-w)+\OO(1/(\delta n))}.
\end{align}

For $w$ in a small neighborhood of $[\min\{b_i+1/(2n),x\}, \max\{b_i+1/2n,x\}]$, then it is bounded away from \(x-s\) and \(a_i\), Stirling's formula gives
\begin{align}\label{e:Gammag3}
\frac{\Gamma(n(w-a_i)+1/2)}{\Gamma(n(w+s-x)+1)}
e^{-n(w-a_i)\log(w-a_i)}
=
\frac{n^{n(x-s-a_i)}e^{-n(x-s-a_i)}}{\sqrt{n(w+s-x)}}
e^{-n (w+s-x)\log(w+s-x) +\OO(1/n)},
\end{align}
and 
\begin{align}\label{e:Gammag4}
\Gamma(ns+1)
=
\sqrt{2\pi ns}\,
\exp\left(n(s\ln n+s\log s-s)+\OO(1/n)\right).
\end{align}

Combining the last four estimates \eqref{e:Gammag1}, \eqref{e:Gammag2}, \eqref{e:Gammag25}, \eqref{e:Gammag3},\eqref{e:Gammag4} with the  formula \eqref{e:Pexp} for
\(P_{ns}(nw,nx)I_i(w)\), we obtain that 
If $\delta\gg n^{-1}$, then
\begin{align}
P_{ns}(nw,nx) I_i(w)= \frac{\sqrt{s}}{\sqrt{2\pi n}\sqrt{x-w}\sqrt{w-(x-s)}} e^{n S(w;x,s)+\OO(1/(\delta n))}.
\end{align}
This gives \eqref{e:PIest}.
If \(x\le b_i-1/(2n)\), 
\[
|P_{ns}(nw,nx) I_i(w)|
\leq
\frac{C\sqrt{s}}{\sqrt{2\pi (n|x-w|+1)|w-(x-s)|}}
e^{n \Re[S(w;x,s)]}.
\]
This gives \eqref{e:PIupbb1}.
If \(x\ge b_i+1/(2n)\), then
\[
|P_{ns}(nw,nx) I_i(w)|
\leq
\frac{C}{\min\{1,n\delta \}}
\frac{\sqrt{s}}{\sqrt{2\pi (n|x-w|+1)|w-(x-s)|}}
e^{n \Re[S(w;x,s)]}.
\]
This gives \eqref{e:PIupbb2}.

\end{proof}

\begin{proof}[Proof of \Cref{l:QIi_bound}]
The estimates for \(Q_{nt}(nz,ny)I_i^{-1}(z)\) follow from the same argument as that in \Cref{l:PIi_bound}, so we omit them. 
\end{proof}

\begin{proof}[Proof of \Cref{l:arctic_S}]
We prove only \eqref{e:PI+_bound}; the proof of \eqref{e:QI-_bound} is identical, using the corresponding estimates for \(Q\) and \(I_-\).

Recall that
\[
w_0=x_0-s_0\chi(w_0).
\]
By \eqref{e:arctic_bounds}, the function
\(\chi(w)\) is bounded away from both \(0\) and \(1\) for all \(w\in \fU\). Thus exactly one of the following three cases holds.

\begin{enumerate}
\item Suppose that \(\chi(w_0)\in(0,1)\). Then \(w_0\in(a_i,b_i)\) for some
\(1\leq i\leq d\). In this case \(\fU\) is bounded away from
\[
(-\infty,x-s]\cup[x,\infty),
\]
and \(I_+=I\). Therefore \eqref{e:PI+_bound} follows directly from
\Cref{l:PIQI_bound}.

\item Suppose that \(\chi(w_0)\in(1,\infty)\). Then \(w_0\in(-\infty,a_i)\) for some
\(1\leq i\leq d\). In this case \(\fU\) is bounded away from the interval
\[
[\min\{a_i,x-s\},\max\{a_i,x-s\}],
\]
and \(I_+=I_i\). Hence \eqref{e:PI+_bound} follows from
\Cref{l:PIi_bound}.

\item Suppose that \(\chi(w_0)\in(-\infty,0)\). Then \(w_0\in(b_i,\infty)\) for some
\(1\leq i\leq d\). In this case \(\fU\) is bounded away from the interval
\[
[\min\{b_i,x\},\max\{b_i,x\}],
\]
and \(I_+=I_i\). Again, \eqref{e:PI+_bound} follows from
\Cref{l:PIi_bound}.
\end{enumerate}
\end{proof}

Before proving \Cref{l:vertical_tangent_S}, we first record the following lemma.
\begin{lemma}\label{l:repulsion}
Fix $b<x$, and let $g$ be analytic in a neighborhood of $[b,x]$ and satisfy
$
g(\overline{z})=\overline{g(z)}.
$
There exists a sufficiently small $\delta>0$ such that, for
$0<r\leq \delta(x-b)$, the following holds.

Let
$
R=\{z\in\bC:\dist(z,[b,x])\leq r\}.
$
Then, on $\partial R$, the vector
\begin{align}\label{e:vector}
-\overline{\left(\ln\frac{z-x}{z-b}+g(z)\right)}
\end{align}
points outward from $R$.
\end{lemma}

\begin{proof}[Proof of \Cref{l:repulsion}]
We use the branch of the logarithm on $\bC\setminus[b,x]$ that is real
on $(x,\infty)$. The boundary of $R$ consists of
\[
\left\{x+re^{\ri\theta}:-\frac{\pi}{2}\leq\theta\leq\frac{\pi}{2}\right\},
\quad
\left\{b-re^{\ri\theta}:-\frac{\pi}{2}\leq\theta\leq\frac{\pi}{2}\right\},
\]
and
\[
\{E\pm \ri r:b\leq E\leq x\}.
\]

On the right semicircle, write $z=x+re^{\ri\theta}$. Its outward unit
normal is $e^{\ri\theta}$, and
\begin{align*}
\Re\left[
-\overline{\left(\ln\frac{z-x}{z-b}+g(z)\right)}
e^{-\ri\theta}
\right]& =
\left(\ln\frac{|z-b|}{r}-\Re g(z)\right)\cos\theta
+\left(\theta-\arg(z-b)+\Im g(z)\right)\sin\theta \\
&=
\left(\ln\frac{x-b}{r}+\OO(1)\right)\cos\theta
+\theta\sin\theta
+\OO\left(\left(\frac{r}{x-b}+r\right)\sin^2\theta\right)>0
\end{align*}
for sufficiently small $\delta$, since
$\cos\theta\geq0$ and $\theta\sin\theta\geq\sin^2\theta$.
The same argument applies to the left semicircle.

On the upper boundary, write $z=E+\ri r$, where $b\leq E\leq x$.
Since $g$ is real on the real axis, $\Im g(E+\ri r)=\OO(r)$, and hence
\begin{align*}
\Im\left[
-\overline{\left(\ln\frac{z-x}{z-b}+g(z)\right)}
\right]
&=
\arg(z-x)-\arg(z-b)+\Im g(z) \\
&=
\arctan\frac{x-E}{r}
+\arctan\frac{E-b}{r}
+\OO(r) \geq
\arctan\frac{x-b}{r}-\OO(r)>0.
\end{align*}
Thus, the vector points upward on the upper boundary. The claim on the
lower boundary follows from the symmetry under complex conjugation.
\end{proof}

\begin{proof}[Proof of \Cref{l:vertical_tangent_S}]
In this setting, \eqref{e:PI+_bound} and \eqref{e:QI-_bound} follow from
\Cref{l:PIi_bound,l:QIi_bound}.

It remains to prove \eqref{e:PIupbb1_copy}; the proof of
\eqref{e:QIupbb1_copy} is identical, so we omit. Since
$x\in\bZ/n$ and $b_i\in\bZ'/n$, either
\[
x\leq b_i-\frac{1}{2n}
\qquad\text{or}\qquad
x\geq b_i+\frac{1}{2n}.
\]
In the first case, \eqref{e:PIupbb1_copy} follows from
\eqref{e:PIupbb1}. We therefore assume that
$x\geq b_i+1/(2n)$.

By \eqref{e:vert_tangent_diff}, in a neighborhood of $b_i$,
\[
S'(w;x,s)
=
\ln\frac{w-x}{w-b_i}
+\text{a conjugation-symmetric analytic function}.
\]
Moreover, \eqref{e:wc_vertical_tangent} and
\eqref{e:wc_vertical_cusp} imply that
\[
\dist(w_c,[b_i,x])\gtrsim x-b_i.
\]
Applying \Cref{l:repulsion} with a sufficiently small tubular radius
proportional to $x-b_i$, we see that the steepest-descent path issuing
from $w_c$ cannot enter this tubular neighborhood. Hence, for every
$w\in\sfD^{\rm d}(w_c)$,
\[
\dist(w,[b_i,x])
\gtrsim x-b_i
\geq \frac{1}{2n}.
\]
The estimate \eqref{e:PIupbb1_copy} now follows from
\eqref{e:PIupbb2}.
\end{proof}

\begin{proof}[Proof of \Cref{l:unit_slope_tangent_S}]
The proof is identical to those for \Cref{l:vertical_tangent_S}, so we omit.
\end{proof}

\begin{proof}[Proof of \Cref{l:horizontal_tangent_S}]
The proof is identical to those for \Cref{l:arctic_S}, so we omit.
\end{proof}

\chapter{Gaussian Free Field}

\section{Global Kernel Approximation}
\label{s:global_kernel_approximate}

In this section, we construct a single global approximation to the inverse Kasteleyn matrix by gluing together the local kernel ansatz. We then prove that the resulting error is uniformly small throughout the polygonal domain.

\subsection{Kernel condition}
In this section, we check that $A_\al((x,s),(y,t))$ as defined in \eqref{e:def_Aalpha} satisfies the same recursion and boundary condition as the inverse Kasteleyn matrix. 

\begin{proposition}\label{l:Aaleq}
Fix $(y,t)\in \bZ^2/n\cap \fP$ bounded away from ramification points. We recall the kernel ansatz \(A_\al=J^{(1)}+J^{(2)}\) from \eqref{e:def_Aalpha}. Then for any \((x,s)\in \fN_\al\), the following holds
\begin{align}\label{e:Aeq}
&A_\al((x,s),(y,t))
 - A_\al((x,s-1/n),(y,t))
 - A_\al((x-1/n,s-1/n),(y,t))= \delta_{(x,s),(y,t)}.
\end{align}
If \((x,s)\) corresponds to a blue boundary triangle of \(\fP\), recall
from \eqref{def:boundary_triangle}, then
\begin{align}\label{e:Aboundar}
A_\al((x,s),(y,t))=0 .
\end{align}
\end{proposition}

\begin{proof}[Proof of \Cref{l:Aaleq}]
The identity \eqref{e:Aeq} follows from the same argument as the proof of
\Cref{l:Aeq}. 

It remains to prove the boundary condition \eqref{e:Aboundar}. We prove the
case where \((x,s)\) corresponds to a blue triangle adjacent to a vertical side
\([\zeta,\zeta']\) of \(\fP\) on the left; see Panel (A) of
\Cref{f:boundary_case}. The cases where \((x,s)\) corresponds to a blue
triangle adjacent to a unit-slope side on the right, or to a horizontal side
on top, are analogous.

Recall from \Cref{s:tangent_location} that the vertical side
\([\zeta,\zeta']\) determines a critical point \(b_i\in\cC(\bR)\). The
unit-slope side passing through \(\zeta'\) determines the critical point
\(a_{i-1}\). The horizontal side passing through \(\zeta\) determines the
point \(\infty_i\), which, after the change of variables, gives a chart
centered at \(0\).

The corner \(\zeta\) belongs to a curvilinear triangle \(\fT_B\), on which
\(\nabla H^*=(1,0)\), while the corner \(\zeta'\) belongs to a curvilinear
triangle \(\fT_A\), on which \(\nabla H^*=(0,0)\). Since \((x,s)\) is a
boundary blue triangle adjacent to \([\zeta,\zeta']\), the neighborhood
\(\fN_\al\) intersects \(\fT_A\), \(\fT_B\), or both.

If \(\fN_\al\cap \fT_B\neq\emptyset\), there are several cases:
\begin{enumerate}
\item In the cases in the first row of \Cref{f:curvilinear_triangle},
\(\fN_\al\) is associated with a tangent or tangent frozen chart centered at
\(b_i\).

\item In the cases in the second row of \Cref{f:curvilinear_triangle}, if
\(\fN_\al\) intersects the horizontal side, then \(\fN_\al\) is associated
with a tangent frozen chart centered at \(0\).

\item The cases in the third row of \Cref{f:curvilinear_triangle} are
impossible.

\item In the cases in the fourth row of \Cref{f:curvilinear_triangle},
\(\fN_\al\) is associated with a tangent or tangent frozen chart centered at
\(b_i\). If \(\fN_\al\) also intersects the horizontal side, then
\(\fN_\al\) is also associated with a tangent frozen chart centered at \(0\).
\end{enumerate}
Similarly, if \(\fN_\al\cap \fT_A\neq\emptyset\), then \(\fN_\al\) is
associated with tangent or tangent frozen charts centered at \(b_i\), or with
tangent frozen charts centered at \(a_{i-1}\).

In all of these charts, the descent contour for the first variable is a small
circular contour. Moreover, by \Cref{l:holomorphic} in the double-contour integral from
\eqref{e:all_term},
\[
\int_{\sfC^{\rm d}} P_{ns}(nw,nx)\, I_i(w)\,(\cdots)\,\rd w,
\]
the \(w\)-integrand is holomorphic in the corresponding local coordinate near
the center of this circular descent contour. More precisely, for charts centered at \(b_i\),
this follows from \(x<b_i\). For charts centered at \(0\), this follows from
\(s\geq t_i^{(0)}\). For charts centered at \(a_{i-1}\), this follows from
\(x-s>a_{i-1}\). Therefore, whenever the circular blue descent contour can be
contracted without crossing a red ascent contour, the double-contour
contribution vanishes.

Recall from \Cref{s:single_contour} that \(J^{(1)}\) is the
single-contour integral whose form depends on the location of
\(\fN_{(y,t)}\), and recall \(J_{\fT}\) from \eqref{e:J1form}. Denote
\begin{align}
J_A((x,s),(y,t))
&:=J_{\fT_A}((x,s),(y,t)),
\qquad
J_B((x,s),(y,t))
:=J_{\fT_B}((x,s),(y,t)).
\end{align}
Then, by \Cref{l:JT_nonvanish}, \(J_A\neq0\) only if
\begin{align}\label{e:JA_nonvanish}
y-t\geq x-s,
\qquad
x\geq y,
\end{align}
and \(J_B\neq0\) only if
\begin{align}\label{e:JB_nonvanish}
t>s,
\qquad
x\geq y.
\end{align}

We now consider the possible positions of \(\fN_{(y,t)}\).
\smallskip

\noindent
\textbf{Case 1.}
Assume that
\[
\fN_{(y,t)}\cap(\fT_A\cup\fT_B)=\emptyset .
\]
Then \(J^{(1)}=0\). Moreover, by \Cref{l:chart_criterion}, \(\fN_\al\) and \(\fN_{(y,t)}\) do not share a
concentric chart. Hence the circular blue descent contour can be contracted to
the empty contour without crossing any red ascent contour. Therefore the
double-contour contribution vanishes, and so
$
A_\al((x,s),(y,t))=0 .
$

\smallskip

\noindent
\textbf{Case 2.}
Assume that
\[
\fN_{(y,t)}\cap(\fT_A\cup\fT_B)\neq\emptyset
\qquad\text{and}\qquad
y>b_i .
\]
Then $y>b_i>x$ and by \eqref{e:JA_nonvanish} and \eqref{e:JB_nonvanish}, both possible single-contour integrals vanish:
\[
J_A=J_B=0 .
\]
For the double-contour contribution, we can choose the contour configuration
relative to whichever of \(\fT_A\) or \(\fT_B\) is relevant, so that the
circular blue descent contour can be contracted to the empty contour without
crossing a red ascent contour. Hence the double-contour contribution also
vanishes. Therefore
$
A_\al((x,s),(y,t))=0 .
$

\smallskip

\noindent
\textbf{Case 3.}
Assume that
\[
\fN_{(y,t)}\cap(\fT_A\cup\fT_B)\neq\emptyset
\qquad\text{and}\qquad
y<b_i .
\]
There are two geometric subcases.

\smallskip

\noindent
\emph{Subcase 3a.}
Suppose that \((y,t)\) lies above the unit-slope side passing through
\(\zeta'\). Then \(\zeta'\) is a concave corner of \(\fP\). Moreover,
\(y-t<x-s\), and hence, by \eqref{e:JA_nonvanish},
\[
J_A=0.
\]

First assume that \(\fN_{(y,t)}\cap\fT_A\neq\emptyset\), for instance, in
the configuration shown in Panel (B) of \(\Cref{f:boundary_case}\). 

If
\(\fN_\al\) and \(\fN_{(y,t)}\) do not share a concentric chart, then
\(\fN_{(y,t)}\) does not meet the shared boundary
\(\fT_A\cap\fT_B\); otherwise, they would share the chart centered at
\(b_i\). We use the contour configuration relative to \(\fT_A\). Hence
\[
J^{(1)}=J_A=0.
\]
Moreover, the circular blue descent contour can be contracted to the empty
contour, so the double-contour contribution vanishes.

Suppose now that \(\fN_\al\) and \(\fN_{(y,t)}\) share a concentric chart.
Then the shared chart is centered at \(b_i\) or \(a_{i-1}\). Suppose first
that it is centered at \(b_i\). If \(\fN_{(y,t)}\) intersects \(\fT_A\)
but not \(\fT_B\), then, since the associated chart carries ascent paths,
\(\fN_{(y,t)}\) lies below the tangency location, and we are in the
configuration of \Cref{f:vertical_tangent2}. In this case, however, the
vertical side lies above the tangency location, contradicting the
assumption that \((y,t)\) lies above the unit-slope side passing through
\(\zeta'\).

Therefore, \(\fN_{(y,t)}\) meets the shared boundary
\(\fT_A\cap\fT_B\), and the vertical side lies below the tangency location.
If \(\fN_\al\cap\fT_A\neq\emptyset\), we use the contour configuration
relative to \(\fT_A\), so
\[
J^{(1)}=J_A=0.
\]
The relative configuration of the blue descent contour and the red ascent
contour is as in Panel (A) of \(\Cref{f:vertical_tangent1}\), Panel (A) of
\(\Cref{f:vertical_tangent4}\), or Panel (A) of
\(\Cref{f:c_vertical_cusp2}\). In each of these configurations, the
circular blue descent contour can be contracted to the empty contour, so
the double-contour contribution also vanishes.

If \(\fN_\al\) intersects only \(\fT_B\), we use the contour configuration
relative to \(\fT_B\), so
\[
J^{(1)}=J_B.
\]
We switch the circular blue descent contour with the circular red ascent
contour. After this switch, the contour configuration is relative to
\(\fT_A\) rather than \(\fT_B\). By \Cref{c:change_triangle}, the
single-contour integral changes from \(J^{(1)}=J_B\) to
\(J^{(1)}=J_A\). Since \(J_A=0\) in the present subcase, the
single-contour contribution vanishes. After the switch, the circular blue
descent contour can be contracted to the empty contour, so the
double-contour contribution also vanishes.

The same argument applies when the shared chart is centered at
\(a_{i-1}\).

It remains in this subcase to consider the possibility that
\[
\fN_{(y,t)}\cap\fT_B\neq\emptyset
\qquad\text{but}\qquad
\fN_{(y,t)}\cap\fT_A=\emptyset.
\]
This can occur only in the configuration shown in Panel (C) of
\(\Cref{f:boundary_case}\).

If \(\fN_\al\) and \(\fN_{(y,t)}\) do not share a concentric chart, then
\[
\fN_\al\cap\fC(\fT_B;\fN_{(y,t)})=\emptyset,
\]
and \(J^{(1)}=0\). Moreover, the circular blue descent contour can be
contracted to the empty contour, so the double-contour contribution
vanishes.

If instead \(\fN_\al\) and \(\fN_{(y,t)}\) share the concentric chart
centered at \(b_i\), then the relative configuration is as in Panel (B) of
\(\Cref{f:vertical_tangent1}\). We use the contour configuration relative
to \(\fT_B\), so
\[
J^{(1)}=J_B.
\]
As above, we switch the circular blue descent contour with the circular red
ascent contour.  After the switch, the
single-contour integral changes to
\(J^{(1)}=J_A=0\), and  the double-contour contribution
also vanishes.

Thus,
$
A_\al((x,s),(y,t))=0
$
throughout Subcase 3a.

\smallskip

\noindent
\emph{Subcase 3b.}
Suppose that \((y,t)\) lies on or below the horizontal side passing through
\(\zeta\). Then \(\zeta\) is a concave corner of \(\fP\). Moreover,
\(s\geq t\), and hence, by \eqref{e:JB_nonvanish},
\[
J_B=0.
\]

First assume that \(\fN_{(y,t)}\cap\fT_B\neq\emptyset\). If
\(\fN_\al\) and \(\fN_{(y,t)}\) do not share a concentric chart, then
\(\fN_{(y,t)}\) does not meet the shared boundary
\(\fT_A\cap\fT_B\); otherwise, they would share the chart centered at
\(b_i\). We use the contour configuration relative to \(\fT_B\). Hence
\[
J^{(1)}=J_B=0.
\]
Moreover, the circular blue descent contour can be contracted to the empty
contour, so the double-contour contribution vanishes.

Suppose now that \(\fN_\al\) and \(\fN_{(y,t)}\) share a concentric chart.
Then the shared chart is centered at \(b_i\) or \(\infty_i\). Suppose first
that the shared chart is centered at \(b_i\). If \(\fN_{(y,t)}\) intersects
\(\fT_B\) but not \(\fT_A\), then, since the associated chart carries
ascent paths, \(\fN_{(y,t)}\) lies above the tangency location, and we are
in the configuration of \Cref{f:vertical_tangent1}. In this case, however,
the vertical side lies below the tangency location, contradicting the
assumption that \((y,t)\) lies on or below the horizontal side passing
through \(\zeta\).

Therefore, \(\fN_{(y,t)}\) meets the shared boundary
\(\fT_A\cap\fT_B\), and the vertical side lies above the tangency location.
If \(\fN_\al\cap\fT_B\neq\emptyset\), we use the contour configuration
relative to \(\fT_B\), so
\[
J^{(1)}=J_B=0.
\]
The relative configuration of the blue descent contour and the red ascent
contour is as in Panel (B) of \(\Cref{f:vertical_tangent2}\), Panel (B) of
\(\Cref{f:vertical_tangent3}\), or Panel (B) of
\(\Cref{f:c_vertical_cusp1}\). In each of these configurations, the
circular blue descent contour can be contracted to the empty contour, so
the double-contour contribution also vanishes.

If \(\fN_\al\) intersects \(\fT_A\) but not \(\fT_B\), we use the contour
configuration relative to \(\fT_A\), so
\[
J^{(1)}=J_A.
\]
We switch the circular blue descent contour with the circular red ascent
contour. After this switch, the contour configuration is relative to
\(\fT_B\) rather than \(\fT_A\). By \Cref{c:change_triangle}, the
single-contour integral changes from \(J^{(1)}=J_A\) to
\(J^{(1)}=J_B\). Since \(J_B=0\) in the present subcase, the
single-contour contribution vanishes. After the switch, the circular blue
descent contour can be contracted to the empty contour, so the
double-contour contribution also vanishes.

The same argument applies when the shared chart is centered at
\(\infty_i\).

It remains in this subcase to consider the possibility that
\[
\fN_{(y,t)}\cap\fT_A\neq\emptyset
\qquad\text{but}\qquad
\fN_{(y,t)}\cap\fT_B=\emptyset.
\]
If \(\fN_\al\) and \(\fN_{(y,t)}\) do not share a concentric chart, then
\[
\fN_\al\cap\fC(\fT_A;\fN_{(y,t)})=\emptyset,
\]
and \(J^{(1)}=0\). Moreover, the circular blue descent contour can be
contracted to the empty contour, so the double-contour contribution
vanishes.

If instead \(\fN_\al\) and \(\fN_{(y,t)}\) share the concentric chart
centered at \(b_i\), then the relative configuration is as in Panel (A) of
\(\Cref{f:vertical_tangent2}\). We use the contour configuration relative
to \(\fT_A\), so
\[
J^{(1)}=J_A.
\]
As above, we switch the circular blue descent contour with the circular red
ascent contour. After the switch,  the
single-contour integral changes to
\(J^{(1)}=J_B=0\), and the double-contour contribution
also vanishes.

Thus,
$
A_\al((x,s),(y,t))=0
$
throughout Subcase 3b.

\smallskip

The three cases above exhaust all possible relative positions of
\(\fN_{(y,t)}\). This proves \eqref{e:Aboundar}.
\end{proof}

\begin{figure}
    \definecolor{bluetri}{rgb}{0.68, 0.85, 0.9}

    \begin{subfigure}{0.28\textwidth}
    \centering

      \caption{}
  \end{subfigure}
\begin{subfigure}{0.28\textwidth}
    \centering
      %
       \caption{}
  \end{subfigure}
		 \begin{subfigure}{0.4\textwidth}
    \centering
      %
  
       \caption{}    
  \end{subfigure}
  \caption{}
  \label{f:boundary_case}
\end{figure}

%
%
%\subsection{Liquid point}
%
%For any $(x_0,s_0)\in \fL$ which is bounded away from the arctic boundary and critical point.  Let  $w_0=x_0-s_0\chi(x_0, s_0)$, 

\subsection{Standard form of approximate kernel}

Let \((y,t)\in\fP\cap\bZ^2/n\) represent a white triangle, and let
\((x,s)\in\fP\cap\bZ^2/n\) represent a blue triangle. Suppose that
\[
(y,t)\in\fN_{(y,t)}
\qquad\text{and}\qquad
(x,s)\in\fN_\alpha.
\]
Recall the approximate kernel \(A_\alpha((x,s),(y,t))\) defined in
\eqref{e:def_Aalpha}.

When deforming the local descent and ascent contours to their corresponding
steepest-descent and steepest-ascent paths, one must account for the
residues at \(w=z\). The following propositions give the resulting standard
form. We postpone its proof to \Cref{s:liquid_standard_form} and \Cref{s:non-liquid_standard_form}.

The contour deformation naturally produces a possibly larger collection of
critical points that depends on the chosen neighborhoods. Recall the sets
\[
\operatorname{Crit}^{\rm d}(x,s;\fN_\alpha)
\qquad\text{and}\qquad
\operatorname{Crit}^{\rm a}(y,t;\fN_{(y,t)})
\]
from \Cref{def:neighborhood_ascent_descent_critical}.

\begin{proposition}\label{p:standard_form}
Let \((y,t)\in\fP\cap\bZ^2/n\) represent a white triangle, and let
\((x,s)\in\fP\cap\bZ^2/n\) represent a blue triangle. Suppose that
\[
(y,t)\in\fN_{(y,t)}
\qquad\text{and}\qquad
(x,s)\in\fN_\alpha.
\]
Throughout the formulas below, \(w_c\) ranges over
\(\operatorname{Crit}^{\rm d}(x,s;\fN_\alpha)\), and \(z_c\) ranges over
\(\operatorname{Crit}^{\rm a}(y,t;\fN_{(y,t)})\), as defined in
\Cref{def:neighborhood_ascent_descent_critical}. There exists \(\fc'>0\)
such that the following statements hold.

Suppose first that neither \(\fN_\alpha\) nor \(\fN_{(y,t)}\) is a ramification
neighborhood. Then
\begin{align}\label{e:Aal}
A_\alpha((x,s),(y,t))
&=
\wh J^{(0)}+\wh J^{(1)}+\wh J^{(2)}
+\sum_{w_c,z_c}
\OO\left(
e^{-\fc'n}
e^{n\Re\left[S(w_c;x,s)-S(z_c;y,t)\right]}
\right).
\end{align}
If \((y,t)\in\fP\setminus\fL\), set
$
\fT:=\operatorname{Cell}^{\rw}(y,t)
$
and define
\begin{align}
\wh J^{(0)}
:=
\cI((x,s),(y,t))J_{\fT}((x,s),(y,t)),
\end{align}
where \(J_{\fT}((x,s),(y,t))\) and
\(\cI((x,s),(y,t))\) are defined in \eqref{e:J1form} and
\Cref{d:interlacing}, respectively. If \((y,t)\in\fL\), we set
$
\wh J^{(0)}:=0.
$
Furthermore,
\begin{align}
\label{e:J1exp}
\wh J^{(1)}
&:=
-\sum_{\xi}
\sgn(\xi)\,
\frac{n}{2\pi\ri}
\int_{\sfD(\xi; (x,s),(y,t))}
P_{ns}(nz,nx)\,Q_{nt}(nz,ny)
\,\rd z,
\\
\label{e:J2exp}
\wh J^{(2)}
&:=
\sum_{w_c,z_c}
\frac{n}{(2\pi\ri)^2}
\int_{\sfD^{\rm a}(z_c)}
\int_{\sfD^{\rm d}(w_c)}
P_{ns}(nw,nx)\,Q_{nt}(nz,ny)
\frac{I_+(w)}{I_-(z)}
\frac{\sqrt{\phi'(w)}\sqrt{\phi'(z)}}
{\phi(w)-\phi(z)}
\,\rd w\,\rd z.
\end{align}
The sum over \(\xi\) ranges over all nonreal intersection points of the
deformed descent paths \(\sfD^{\rm d}(w_c)\) and ascent paths
\(\sfD^{\rm a}(z_c)\). The sign \(\sgn(\xi)\) is determined by the local
orientation of the crossing, and the steepest-descent path $\sfD(\xi; (x,s),(y,t))$ is defined in \Cref{s:pof_single_integral}. For each pair \((w_c,z_c)\), the factors
\(I_+\) and \(I_-\) are determined by the corresponding descent and ascent
charts.

Suppose next that at least one of \(\fN_\alpha\) and
\(\fN_{(y,t)}\) is a ramification neighborhood. Then \eqref{e:Aal} remains
valid with
\[
\wh J^{(0)}=0
\]
and with the following modifications.

If \(\fN_\alpha\) is a ramification neighborhood, in \eqref{e:J1exp} and \eqref{e:J2exp}, we replace
\[
P_{ns}(nw,nx), \sfD^{\rm d}(w_c), I_+(w), \phi(w), \phi'(w)\quad
\text{by}\quad 
P_{n(s+\ft)}(nw,nx),\sfD^{\rm d}(w_{c,\ft}), I_\ft(w),
\phi_\ft(w), \phi_\ft'(w),
\]
respectively. Similarly, if $\fN_{(y,t)}$ is a ramification
neighborhood, in \eqref{e:J1exp} and \eqref{e:J2exp}, we replace 
\[
Q_{nt}(nz,ny), \sfD^{\rm a}(z_c), I_-(z), \phi(z), \phi'(z)
\quad
\text{by}\quad 
Q_{n(t+\ft)}(nz,ny), \sfD^{\rm a}(z_{c,\ft}), I_\ft(z),
\phi_\ft(z), \phi_\ft'(z),
\]
respectively. Here $w_{c,\ft}$ and $z_{c,\ft}$ denote the corresponding critical points in the ramification charts.

\end{proposition}

\begin{proposition}\label{p:standard_form2}
Adopt notations and assumptions in \Cref{p:standard_form}. 
Assume that \(\fN_\alpha\) is a liquid or ramification neighborhood, and $(x,s)$ is bounded away from $(y,t)$. If $\fN_{(y,t)}$ is not a ramification neighborhood,
 then
\begin{align}\label{e:Aal_ramification}
A_\alpha((x,s),(y,t))
&=
\sum_{w_c,z_c}
\frac{(1+\OO(1/n))}{(2\pi\ri)^2}
e^{nS(w_c;x,s)}
\int_{\sfD^{\rm a}(z_c)}
\frac{Q_{nt}(nz,ny)}{I_-(z)}
\frac{\sqrt{\phi'(z)}\sqrt{\del_x\phi(x,s)}}
{\phi(w_c)-\phi(z)}
\,\rd z
\nonumber\\
&\quad+
\sum_{w_c,z_c}
\OO\left(
e^{-\fc'n}
e^{n\Re\left[S(w_c;x,s)-S(z_c;y,t)\right]}
\right).
\end{align}
If $\fN_{(y,t)}$ is a ramification neighborhood,
 then
\begin{align}\label{e:Aal_ramification2}
A_\alpha((x,s),(y,t))
&=
\sum_{w_c,z_{c,\ft}}
\frac{(1+\OO(1/n))}{(2\pi\ri)^2}
e^{nS(w_c;x,s)}
\int_{\sfD^{\rm a}(z_{c,\ft})}
\frac{Q_{n(t+\ft)}(nz,ny)}{I_\ft(z)}
\frac{\sqrt{\phi_\ft'(z)}\sqrt{\del_x\phi(x,s)}}
{\phi(w_c)-\phi_\ft(z)}
\,\rd z
\nonumber\\
&\quad+
\sum_{w_c,z_c}
\OO\left(
e^{-\fc'n}
e^{n\Re\left[S(w_c;x,s)-S(z_c;y,t)\right]}
\right).
\end{align}

\end{proposition}

The standard forms in \Cref{p:standard_form} depend on the neighborhoods
\(\fN_\al\) and \(\fN_{(y,t)}\). Recall from
\Cref{def:ascent_descent_critical} the sets of descent and ascent critical
points
\[
\operatorname{Crit}^{\rm d}(x,s),
\qquad
\operatorname{Crit}^{\rm a}(y,t).
\]
By construction,
\begin{align}\label{e:critical_set_inclusion}
\operatorname{Crit}^{\rm d}(x,s)
&\subseteq
\operatorname{Crit}^{\rm d}(x,s;\fN_\al),
\qquad
\operatorname{Crit}^{\rm a}(y,t)
\subseteq
\operatorname{Crit}^{\rm a}(y,t;\fN_{(y,t)}).
\end{align}
The following proposition shows that the terms indexed by the additional
critical points in
\[
\operatorname{Crit}^{\rm d}(x,s;\fN_\al)
\setminus
\operatorname{Crit}^{\rm d}(x,s),
\qquad
\operatorname{Crit}^{\rm a}(y,t;\fN_{(y,t)})
\setminus
\operatorname{Crit}^{\rm a}(y,t)
\]
are negligible. Moreover, the steepest-descent and steepest-ascent paths may
be truncated to sufficiently small neighborhoods of their corresponding
critical points. In both cases, the discarded contributions can be absorbed
into the exponentially small error terms. We postpone its proof in \Cref{s:largerset_proof}.

\begin{proposition}\label{p:largerset}
Under the assumptions of \Cref{p:standard_form}, the formulas
\eqref{e:Aal} and \eqref{e:Aal_ramification} remain valid after restricting
\(w_c\) to
\[
\operatorname{Crit}^{\rm d}(x,s)
\]
and \(z_c\) to
\[
\operatorname{Crit}^{\rm a}(y,t).
\]
If \((x,s)\) belongs to an arctic neighborhood, then \(w_c\) may be further
restricted to the critical points contained in the arctic chart, and the
critical point contained in the additional frozen chart may be omitted.
The analogous statement holds for \(z_c\) if \((y,t)\) belongs to an arctic
neighborhood.

Moreover, each steepest-descent or steepest-ascent path appearing in
\eqref{e:Aal} and \eqref{e:Aal_ramification}  may be truncated to a sufficiently
small neighborhood of its corresponding critical point. The contributions
from the omitted critical points and the discarded portions of the paths
are absorbed into the exponentially small error terms in
\eqref{e:Aal} and \eqref{e:Aal_ramification}.
\end{proposition}

\subsection{Compatibility}
It may happen that a point
\((x,s)\) belongs to two distinct neighborhoods:
\[
(x,s)\in\fN_\alpha\cap\fN_\beta,
\qquad
\alpha\neq\beta.
\]
We show that the corresponding local formulas are compatible on the
overlap. More precisely, we estimate
\[
A_\alpha((x,s),(y,t))-A_\beta((x,s),(y,t)).
\]

\begin{figure}
\centering
\includegraphics[scale=0.3]{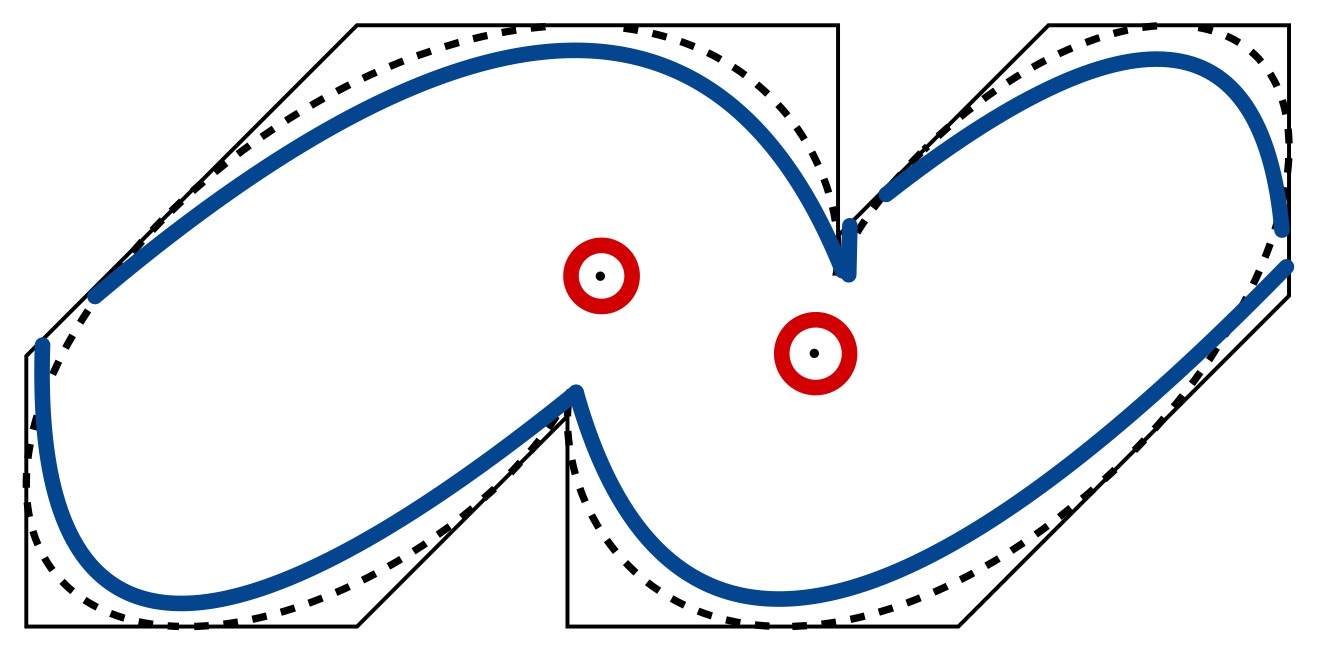}
\caption{The exceptional overlap regions. The region \(\fO_1\) is shown
in red, and the region \(\fO_2\) is shown in blue.}
\label{f:12gon}
\end{figure}

There are two types of overlaps that require separate treatment.

\noindent\textbf{Case 1: Ramification--liquid overlap \(\fO_1\).}

We say that \((x,s)\in\fO_1\) if
\[
(x,s)\in\fN_\alpha\cap\fN_\beta,
\]
where, after interchanging \(\alpha\) and \(\beta\) if necessary,
\(\fN_\alpha\) is a ramification neighborhood and \(\fN_\beta\) is a
regular liquid neighborhood.

By construction, \(\fO_1\) is a finite union of small annular regions
surrounding the ramification points. Each ramification point lies in the
hole of the corresponding annulus; see the red regions in
\Cref{f:12gon}. In terms of critical points, the descent critical points
associated with such an overlap lie in the corresponding annular overlap
between the ramification and regular liquid charts on \(\cC\).

\noindent\textbf{Case 2: \(I_i\)--\(I\) overlap \(\fO_2\).}

We say that \((x,s)\in\fO_2\) if
\[
(x,s)\in\fN_\alpha\cap\fN_\beta
\]
for a pair of overlapping neighborhoods
\((\fN_\alpha,\fN_\beta)\) satisfying the following conditions, after
interchanging \(\alpha\) and \(\beta\) if necessary.

The neighborhood \(\fN_\alpha\) is associated with a chart \(\fU\)
centered at a point
\[
w_0\in[b_i,a_i],
\]
and \(\fU\) contains a descent critical point
\[
w_c\notin[b_i,a_i]
\]
associated with \((x,s)\). Thus, in the corresponding double-contour
formula for \(\fN_\alpha\), one has
\[
I_+=I_i.
\]
The overlapping neighborhood \(\fN_\beta\) is associated with another
chart \(\fU'\), which also contains \(w_c\), but is centered at a point
\[
w_0'\notin[b_i,a_i].
\]
In the corresponding double-contour formula for \(\fN_\beta\), one has
\[
I_+=I.
\]

Thus, \(\fO_2\) consists of points \((x,s)\) associated with a descent
critical point \(w_c\) lying outside, but sufficiently close to
\([b_i,a_i]\),  that it is contained in both charts \(\fU\) and \(\fU'\).
Such an overlap can occur only in the following three cases.

\begin{figure}
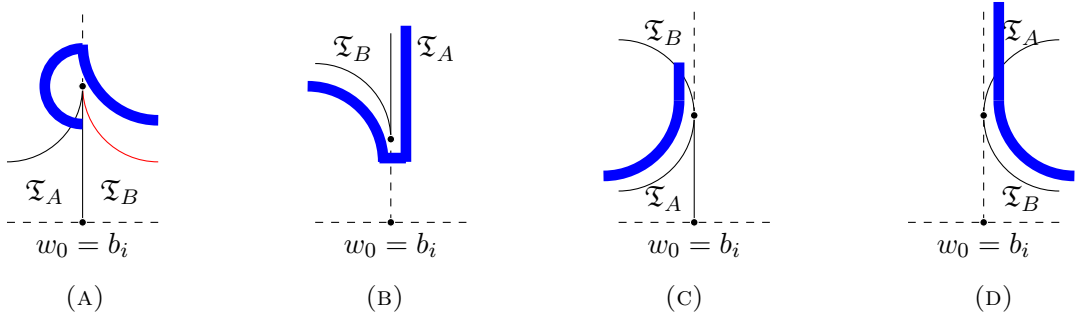
			
\begin{subfigure}[t]{0.24\textwidth}
			\centering

			\caption{}
			\end{subfigure}%
			\begin{subfigure}[t]{0.24\textwidth}
			\centering
			%
			\caption{}
			\end{subfigure}%				
	\begin{subfigure}[t]{0.24\textwidth}
			\centering
			%
			\caption{}
		\end{subfigure}	
		\begin{subfigure}[t]{0.24\textwidth}

		\centering
			%
			\caption{}
	
			\end{subfigure}

	\caption{
\label{f:overlap_region}
Examples of $I_i-I$ overlap region. The blue bands
indicate portions of \(\fO_2\).}
	\end{figure}

\begin{enumerate}
\item
\emph{Overlap along the arctic boundary.}
The neighborhood \(\fN_\alpha\) is an arctic, tangent, cusp, or
cusp-turning neighborhood and intersects a curvilinear triangle
\(\fT_B\) such that
\begin{align}\label{e:DH}
\nabla H^*
\in
\{(1,0),(1,-1)\}
\qquad\text{on }\fT_B.
\end{align}
In this case, the overlapping neighborhood \(\fN_\beta\) is a liquid
neighborhood.

Geometrically, this part of \(\fO_2\) consists of narrow bands inside the
liquid region \(\fL\), running along the portions of the arctic boundary
contained in curvilinear triangles satisfying \eqref{e:DH}; see
\Cref{f:overlap_region}.

\item
\emph{Overlap along a vertical extended side.}
The chart \(\fU\) is centered at \(w_0=b_i\), and \(\fN_\alpha\) lies
on a vertical extended side separating two adjacent curvilinear triangles
\(\fT_A\) and \(\fT_B\), where
\[
\nabla H^*=(0,0)
\quad\text{on }\fT_A,
\qquad
\nabla H^*=(1,0)
\quad\text{on }\fT_B.
\]
The overlapping neighborhood satisfies
\[
\fN_\beta\subset\fT_A\cup\fL,
\qquad
\fN_\beta\cap\fT_B=\emptyset.
\]

The chart \(\fU'\), centered at some \(w_0'<b_i\), is either an arctic
chart, in which case \(\fN_\beta\) is an arctic neighborhood, or a
frozen chart,  in which case \(\fN_\beta\) is a frozen neighborhood. In the frozen case, let \(L(w_0')\) be the tangent line
corresponding to \(w_0'\). This line is tangent to the portion of the
arctic boundary contained in \(\fT_A\) at a point
\[
(x_0',s_0')\in\fA
\]
and has slope in \((1,\infty)\). Since \(\fU'\) contains the descent
critical point \(w_c\), by the geometric criterion \eqref{e:geometric_descent_ascent}, the neighborhood \(\fN_\beta\) lies above the
tangency point \((x_0',s_0')\) along \(L(w_0')\).

The relevant tangent or cusp-turning point divides the vertical extended
side into an upper and a lower part. The corresponding portion of
\(\fO_2\) is a narrow band contained in \(\fT_A\cup\fL\), running along
the upper part of the vertical extended side; see Panels~(B) and~(D) of
\Cref{f:overlap_region}.

\item
\emph{Overlap along a unit-slope extended side.}
The chart \(\fU\) is centered at \(w_0=a_i\), and \(\fN_\alpha\) lies
on a unit-slope extended side separating two adjacent curvilinear triangles
\(\fT_A\) and \(\fT_B\), where
\[
\nabla H^*=(0,0)
\quad\text{on }\fT_A,
\qquad
\nabla H^*=(1,-1)
\quad\text{on }\fT_B.
\]
The overlapping neighborhood satisfies
\[
\fN_\beta\subset\fT_A\cup\fL,
\qquad
\fN_\beta\cap\fT_B=\emptyset.
\]

The chart \(\fU'\), centered at some \(w_0'>a_i\), is either an arctic
chart, in which case \(\fN_\beta\) is an arctic neighborhood, or a
frozen chart, in which case \(\fN_\beta\) is a frozen neighborhood. In the frozen case, let \(L(w_0')\) be the tangent line
corresponding to \(w_0'\). This line is tangent to the portion of the
arctic boundary contained in \(\fT_A\) at a point
\[
(x_0',s_0')\in\fA
\]
and has slope in \((1,\infty)\). As above, \(\fN_\beta\) lies above the
tangency point \((x_0',s_0')\) along \(L(w_0')\).

The relevant tangent or cusp-turning point divides the unit-slope extended
side into an upper and a lower part. The corresponding portion of
\(\fO_2\) is a narrow band contained in \(\fT_A\cup\fL\), running along
the upper part of the unit-slope extended side.
\end{enumerate}

Finally, we set
\begin{align}\label{e:def_overlap_region}
\fO:=\fO_1\cup\fO_2.
\end{align}
The regions in \(\fO_1\) are shown in red in \Cref{f:12gon}, and those
in \(\fO_2\) are shown in blue. The discussion above leads to the following lemma.

\begin{lemma}\label{p:overlap_critical_region}
There exists a fixed region
\[
\fW\subset\cC,
\]
depending only on the chosen neighborhood construction, with the following
properties:
\begin{enumerate}
\item
\(\fW\) is a finite union of small annular regions surrounding the
ramification points of \(\cC\) and thin stretched annular regions
surrounding the real arcs
\[
[b_i,a_i]\subset\cC(\bR),
\qquad
1\leq i\leq d.
\]

\item
Every \((x,s)\in\fO\) belongs to a liquid, frozen, or arctic neighborhood,
and the descent critical points responsible for the overlap belong to
\(\fW\).  More precisely:
\begin{enumerate}
\item
if \((x,s)\in\fO_1\), then every descent critical point contained in both
the relevant ramification chart and the relevant regular liquid chart lies
in an annular component of \(\fW\) surrounding the corresponding
ramification point;

\item
if \((x,s)\in\fO_2\), then the common descent critical point
\[
w_c\in\fU\cap\fU'
\]
lies in a stretched annular component of \(\fW\) surrounding
\([b_i,a_i]\).
\end{enumerate}
\end{enumerate}
\end{lemma}

Let \((x,s)\in\fP\), and let \(w_c\) be an associated formal critical
point. We introduce a localization length \(\Delta(w_c)\), depending on
whether \((x,s)\) lies in a liquid, frozen, or arctic neighborhood.
\begin{definition}\label{def:localization_length}
Let \((x,s)\in\fP\), and let \(w_c\) be an associated formal critical
point.
\begin{enumerate}
\item
If \((x,s)\) lies in a liquid or frozen neighborhood, then
\[
\Delta(w_c):=n^{-1/2}.
\]

\item
If \((x,s)\) lies in an arctic neighborhood and \(w_c\) is contained in
the corresponding arctic chart, then
\begin{align}\label{e:arctic_D}
\Delta(w_c)
:=
\min\left\{
\frac{1}{n^{1/2}\dist((x,s),\fA)^{1/4}},
\frac{1}{n^{1/3}}
\right\}.
\end{align}
If \(w_c\) is contained in a frozen chart, we still set
\begin{align}\label{e:arctic_D2}
\Delta(w_c):=n^{-1/2}\lesssim \min\left\{
\frac{1}{n^{1/2}\dist((x,s),\fA)^{1/4}},
\frac{1}{n^{1/3}}
\right\}.
\end{align}

\item
In all other cases, we simply set
\[
\Delta(w_c):=1.
\]
\end{enumerate}
\end{definition}
Let \((x,s),(y,t)\in\bZ^2/n\) represent white triangles contained in
\(\fP\). We introduce the following quantity, which will be used to bound the difference of approximate kernels in \Cref{p:Bbound}.
\begin{align}\label{e:intro_B}
B((x,s),(y,t))
:=
\sum_{w_c,z_c}
\left(
\frac{\bm1((x,s)\in\fO, w_c\in \fW)}{n}+e^{-\fc n}
\right)
\Delta(w_c)\Delta(z_c)
e^{n\Re[S(w_c;x,s)-S(z_c;y,t)]},
\end{align}
where \(w_c\) ranges through
\(\operatorname{Crit}^{\rm d}(x,s^-)\), and \(z_c\) ranges through
\(\operatorname{Crit}^{\rm a}(y,t)\).

The following proposition states that if
\(\fN_\al\cap\fN_\beta\neq\emptyset\), then the local formulas
\(A_\al\) and \(A_\beta\) are compatible on the overlap, with their
difference bounded by $B$ as in \eqref{e:intro_B}. We defer its proof to
\Cref{s:compatible}.
\begin{proposition}\label{p:Bbound}
Let \((x,s),(y,t)\in\bZ^2/n\) represent white triangles contained in
\(\fP\), and suppose that
\[
(x,s)\in\fN_\al\cap\fN_\beta.
\]
We have
\begin{align}\begin{split}\label{e:Bbound}
&\left|A_\al\left(((x,s),(y,t))\right)-A_\beta\left(((x,s),(y,t))\right)\right|\lesssim (\ln n)^{10}B((x,s),(y,t)),\\
&\left|A_\al\left(\left(x-\frac{1}{n},s-\frac{1}{n}\right);(y,t)\right)-A_\beta\left(\left(x-\frac{1}{n},s-\frac{1}{n}\right);(y,t)\right)\right|\lesssim (\ln n)^{10} B((x,s),(y,t)),\\
&\left|A_\al\left(\left(x,s-\frac{1}{n}\right);(y,t)\right)-A_\beta\left(\left(x,s-\frac{1}{n}\right);(y,t)\right)\right|\lesssim (\ln n)^{10} B((x,s),(y,t)).
\end{split}\end{align}

\end{proposition}

The following proposition gives a convolution-type upper bound for
\(B(\,\cdot\,;\,\cdot\,)\).

\begin{proposition}\label{p:B_est}
Let \((x,s),(y,t)\in\bZ^2/n\) represent white triangles contained in
\(\fP\). Then, summing over all \((u,v)\in\bZ^2/n\) representing white
triangles contained in \(\fP\), we have
\begin{align}
\sum_{(u,v)}
B((x,s);(u,v))B((u,v);(y,t))
\lesssim
B((x,s),(y,t)).
\end{align}
\end{proposition}

\begin{proof}[Proof of \Cref{p:B_est}]
Let
\[
w_c\in\operatorname{Crit}^{\rm d}(x,s^-),
\qquad
\xi_c^{\rm d}\in\operatorname{Crit}^{\rm d}(u,v^-),
\]
and
\[
\xi_c^{\rm a}\in\operatorname{Crit}^{\rm a}(u,v),
\qquad
z_c\in\operatorname{Crit}^{\rm a}(y,t).
\]
Then, by \Cref{p:critical_compare},
\begin{align}
\Re S(\xi_c^{\rm a};u,v)
\geq
\Re S(\xi_c^{\rm d};u,v).
\end{align}

We expand the product and estimate it as follows:
\begin{align}
&\sum_{(u,v)}
B((x,s);(u,v))B((u,v);(y,t))
\notag\\
&=
\sum_{(u,v)}
\sum_{w_c,\xi_c^{\rm a}}
\left(
\frac{\bm1((x,s)\in\fO, w_c\in \fW)}{n}+e^{-\fc n}
\right)
\Delta(w_c)\Delta(\xi_c^{\rm a})
e^{n\Re[S(w_c;x,s)-S(\xi_c^{\rm a};u,v)]}
\notag\\
&\qquad\times
\sum_{\xi_c^{\rm d},z_c}
\left(
\frac{ \bm1((u,v)\in\fO, \xi^{\rm d}_c\in \fW)}{n}+e^{-\fc n}
\right)
\Delta(\xi_c^{\rm d})\Delta(z_c)
e^{n\Re[S(\xi_c^{\rm d};u,v)-S(z_c;y,t)]}
\notag\\
&\lesssim
\sum_{w_c,z_c}
\left(
\frac{\bm1((x,s)\in\fO, w_c\in \fW)}{n}+e^{-\fc n}
\right)
\Delta(w_c)\Delta(z_c)
e^{n\Re[S(w_c;x,s)-S(z_c;y,t)]}
\notag\\
&\qquad\times
\sum_{(u,v)}
\sum_{\xi_c^{\rm a},\xi_c^{\rm d}}
\left(
\frac{ \bm1((u,v)\in\fO)}{n}+e^{-\fc n}
\right)
\Delta(\xi_c^{\rm a})\Delta(\xi_c^{\rm d}).
\end{align}
The first factor on the right-hand side is
\(B((x,s),(y,t))\). It remains to bound the last factor.

For \((u,v)\in\fO\) in a liquid or frozen chart, we have
\[
\Delta(\xi_c^{\rm a}),
\Delta(\xi_c^{\rm d})
\lesssim n^{-1/2},
\]
and hence
\begin{align}
\sum_{\substack{(u,v)\in\fO\\
\mathrm{liquid/frozen}}}
\sum_{\xi_c^{\rm a},\xi_c^{\rm d}}
\frac{1}{n}
\Delta(\xi_c^{\rm a})\Delta(\xi_c^{\rm d})
\lesssim 1.
\end{align}
For \((u,v)\in\fO\) in an arctic chart, we use
\eqref{e:arctic_D} and \eqref{e:arctic_D2}. Thus,
\begin{align}
\sum_{\substack{(u,v)\in\fO\\
\mathrm{arctic}}}
\sum_{\xi_c^{\rm a},\xi_c^{\rm d}}
\frac{1}{n}
\Delta(\xi_c^{\rm a})\Delta(\xi_c^{\rm d})
&\lesssim
\sum_{\substack{(u,v)\in\fO\\
\mathrm{arctic}}}
\frac{1}{n}
\min\left\{
\frac{1}{n\dist((u,v),\fA)^{1/2}},
\frac{1}{n^{2/3}}
\right\}
\\
&\lesssim
\frac{1}{n^{2/3}}
+
\sum_{1\leq i\leq n}
\min\left\{
\frac{1}{n(i/n)^{1/2}},
\frac{1}{n^{2/3}}
\right\}
\lesssim 1.
\end{align}
Finally, since all localization lengths are bounded by \(1\) and the
numbers of ascent and descent critical points are uniformly bounded,
\begin{align}
\sum_{(u,v)}
\sum_{\xi_c^{\rm a},\xi_c^{\rm d}}
e^{-\fc n}
\Delta(\xi_c^{\rm a})\Delta(\xi_c^{\rm d})
\lesssim
n^2e^{-\fc n}
\lesssim 1.
\end{align}
Combining the preceding estimates gives
\begin{align}\label{e:sumoveruv}
\sum_{(u,v)}
\sum_{\xi_c^{\rm a},\xi_c^{\rm d}}
\left(
\frac{ \bm1((x,s)\in\fO)}{n}+e^{-\fc n}
\right)
\Delta(\xi_c^{\rm a})\Delta(\xi_c^{\rm d})
\lesssim 1,
\end{align}
and the claim follows.
\end{proof}

\subsection{Approximate global kernel}

Take a partition of unity \(\{\rho_\al\}\) adapted to the neighborhoods
\(\{\fN_\al\}\). Let \((y,t)\in\bZ^2/n\) represent a white triangle
contained in \(\fP\), and let \((x,s)\in\bZ^2/n\) represent a blue triangle. We define the following global kernel:
\begin{align}\label{e:defA}
A((x,s),(y,t))
=
\sum_{\al}\rho_\al(x,s)A_\al((x,s),(y,t)).
\end{align}
If \((x,s)\) represents a blue boundary triangle of \(\fP\), as defined in
\Cref{def:boundary_triangle}, then, by \eqref{e:Aboundar}, we have
\begin{align}
A((x,s),(y,t))=0.
\end{align}

The following theorem states that $A((x,s),(y,t))$ is a good approximation of the inverse Kasteleyn matrix
\begin{theorem}\label{l:global_approximation}
Let \((y,t)\in\bZ^2/n\) represent a white triangle, and let
\((x,s)\in\bZ^2/n\) represent a blue triangle, both contained in \(\fP\).
Then
\begin{align}\label{e:global_approx}
\left|K^{-1}((x,s),(y,t))-A((x,s),(y,t))\right|
\lesssim
\frac{(\ln n)^{10}}{n}(A|B|)((x,s),(y,t)).
\end{align}
\end{theorem}

\Cref{l:global_approximation} follows from the following proposition, which
states that \(A((x,s),(y,t))\) satisfies the same recursion as the 
inverse Kasteleyn matrix up to a small error.

\begin{proposition}\label{l:global_kernel}
Let \((x,s),(y,t)\in\bZ^2/n\) represent white triangles contained in
\(\fP\). Then
\begin{align}
\begin{split}\label{e:OK}
(KA)((x,s),(y,t))
&=
\sum_\al
\Biggl[
\rho_\al(x,s)A_\al((x,s),(y,t))
-\rho_\al\left(x,s-\frac{1}{n}\right)
A_\al\left(\left(x,s-\frac{1}{n}\right);(y,t)\right)\\
&\qquad
-\rho_\al\left(x-\frac{1}{n},s-\frac{1}{n}\right)
A_\al\left(\left(x-\frac{1}{n},s-\frac{1}{n}\right);(y,t)\right)
\Biggr]\\
&=
\delta_{(x,s),(y,t)}-\cE((x,s),(y,t)),
\end{split}
\end{align}
where
\begin{align}
\left|\cE((x,s),(y,t))\right|
\lesssim
\frac{(\ln n)^{10}}{n}B((x,s),(y,t)).
\end{align}
\end{proposition}

\Cref{l:global_approximation} follows from the following proposition, which states that $A((x,s),(y,t))$ satisfies the same boundary condition as the inverse Kasteleyn matrix, and recursion with small error.

\begin{proof}[Proof of \Cref{l:global_approximation}]
Multiplying \eqref{e:OK} on the left by \(K^{-1}\), we obtain
\begin{equation}\label{e:resolvent}
A=K^{-1}-K^{-1}\cE=K^{-1}(I-\cE).
\end{equation}

Thanks to \eqref{p:B_est},  there exists a constant \(C>0\) such that, for every \(m\ge 1\),
\begin{equation}\label{e:Em_bound}
|(\cE^m)((x,s),(y,t))|\le \left(\frac{C(\ln n)^{10}}{n}\right)^m\,B((x,s),(y,t)).
\end{equation}
In particular, for \(n\) sufficiently large, the Neumann series
\[
(I-\cE)^{-1}=I+\cE+\cE^2+\cdots
\]
converges entrywise absolutely. Hence \eqref{e:resolvent} implies
\[
K^{-1}=A(I-\cE)^{-1},
\]
and therefore
\[
K^{-1}-A=A\sum_{m\ge 1}\cE^m.
\]

Using \eqref{e:Em_bound}, we obtain
\begin{align*}
&\phantom{{}={}}|(K^{-1}-A)((x,s),(y,t))|
=\Bigl|\Bigl(A\sum_{m\ge 1}\cE^m\Bigr)((x,s),(y,t))\Bigr|\le \sum_{m\ge 1}|(A\cE^m)((x,s),(y,t))| \\
&\lesssim \sum_{m\ge 1}\left(\frac{C(\ln n)^{10}}{n}\right)^m (|A|B)((x,s),(y,t))
\lesssim \frac{(\ln n)^{10}}{n} (|A|B)((x,s),(y,t)),
\end{align*}
as claimed.
\end{proof}

\begin{proof}[Proof of \Cref{l:global_kernel}]

For the simplicity of notation, in this proof, we simply write $A_\al(\cdot ;y,t)=A_\al(\cdot)$, and $B(\cdot ;y,t)=B(\cdot)$.
Since $\{\fN_\al\}_{\al}$ covers the polygon $\fP$, there exists an index $\beta$ such that $(x,s), (x-1/n,s-1/n), (x,s-1/n)\in \fN_\beta$, provided $n$ is sufficiently large. 

Since the partition of unity \(\{\rho_\alpha\}\) is adapted to the neighborhoods
\(\{\fN_\alpha\}\), if \(\rho_\alpha(x,s)\neq 0\), then \((x,s)\in \fN_\alpha\). Hence, by \Cref{l:Aaleq},
\begin{align}\label{e:kernel_eq}
\delta_{(x,s),(y,t)}=A_\al(x,s)-A_\al\left(x,s-\frac{1}{n}\right)-A_\al\left(x-\frac{1}{n},s-\frac{1}{n}\right)
\end{align}
We can then take difference of \eqref{e:OK} and \eqref{e:kernel_eq} multiplied by $\rho_\al$
\begin{align}\begin{split}\label{e:OK2}
(K A)(x,s)&=\delta_{(x,s),(y,t)}+\sum_\al \left(\rho_\al(x,s)-\rho_\al\left(x,s-\frac{1}{n}\right)\right)A_\al\left(x,s-\frac{1}{n}\right)\\
&+\left(\rho_\al(x,s)-\rho_\al\left(x-\frac{1}{n},s-\frac{1}{n}\right)\right)A_\al\left(x-\frac{1}{n},s-\frac{1}{n}\right)
\end{split}\end{align}

Since $\{\rho_\al\}$ is a partition of unity, 
\begin{align}\label{e:partition}
\sum_\al \left(\rho_\al(x,s)-\rho_\al\left(x,s-\frac{1}{n}\right)\right)=\sum_\al \left(\rho_\al(x,s)-\rho_\al\left(x-\frac{1}{n},s-\frac{1}{n}\right)\right)=0.
\end{align}
then we can multiply \eqref{e:partition} by $A_\beta(x,s-1/n)$ and $A_\beta(x-1/n, s-1/n)$, and take difference with \eqref{e:OK2}, 
\begin{align}\begin{split}\label{e:OK3}
(K A)(x,s)&=\sum_\al \left(\rho_\al(x,s)-\rho_\al\left(x,s-\frac{1}{n}\right)\right)\left(A_\al\left(x,s-\frac{1}{n}\right)-A_\beta\left(x,s-\frac{1}{n}\right)\right)\\
&+\left(\rho_\al(x,s)-\rho_\al\left(x-\frac{1}{n},s-\frac{1}{n}\right)\right)\left(A_\al\left(x-\frac{1}{n},s-\frac{1}{n}\right)-A_\beta\left(x-\frac{1}{n},s-\frac{1}{n}\right)\right).
\end{split}\end{align}

We now estimate the two sums on the right-hand side of \eqref{e:OK3}. Fix \(\alpha\). If none of the points
$(x,s), (x-1/n, s-1/n), (x,s-1/n)$
belongs to \(\operatorname{supp}\rho_\alpha\), then
\begin{align}\label{e:twocase2}
\rho_\alpha(x,s)-\rho_\alpha\left(x,s-\frac{1}{n}\right)
=
\rho_\alpha(x,s)-\rho_\alpha\left(x-\frac{1}{n},s-\frac{1}{n}\right)
=0.
\end{align}
Otherwise, at least one of these three points belongs to \(\operatorname{supp}\rho_\alpha\). Since the partition of unity is adapted to the cover, we may assume that, for each \(\alpha\), a fixed small neighborhood of \(\operatorname{supp}\rho_\alpha\) is contained in \(\fN_\alpha\). Therefore, for \(n\) sufficiently large, all three points above belong to \(\fN_\alpha\). Then $(x,s), (x-1/n, s-1/n), (x,s-1/n)\in \fN_\al \cap \fN_\beta$, and we have
\begin{align}\begin{split}\label{e:twocase1}
&\left|\rho_\al(x,s)-\rho_\al\left(x,s-\frac{1}{n}\right)\right|,\quad \left|\rho_\al(x,s)-\rho_\al\left(x-\frac{1}{n},s-\frac{1}{n}\right)\right|\lesssim \frac{1}{n},\\
&\left|A_\al\left(x,s-\frac{1}{n}\right)-A_\beta\left(x,s-\frac{1}{n}\right)\right|,\quad 
\left|A_\al\left(x-\frac{1}{n},s-\frac{1}{n}\right)-A_\beta\left(x-\frac{1}{n},s-\frac{1}{n}\right)\right|\lesssim (\ln n)^{10} B(x,s),
\end{split}\end{align}
where the first line follows from Lipschtiz properties of $\rho_\al$, and the second row follows from \Cref{p:Bbound}.

Combining \eqref{e:twocase1} and \eqref{e:twocase2} with \eqref{e:OK3}, and using that the cover is finite, we conclude that
\begin{align}
\left|(K A)(x,s)-\delta_{(x,s),(y,t)}\right|
\lesssim \frac{(\ln n)^{10}}{n}B(x,s).
\end{align}
This proves the claim.
\end{proof}

\subsection{Estimates for the inverse Kasteleyn matrix}\label{subsec:edge-facet}

Throughout this section, \((x,s)\in\bZ^2/n\) represents a blue triangle,
whereas \((y,t)\in\bZ^2/n\) represents a white triangle. We state estimates
for the inverse Kasteleyn matrix
\[
K^{-1}((x,s),(y,t))
\]
in three regimes: when both arguments lie in the liquid region, when they
are close to the arctic boundary \(\fA\), and when they lie in the frozen
region.

We first introduce an exceptional arctic set \(\fA_{\rm ex}\) with the
following property: if \((x,s)\) is bounded away from \(\fA_{\rm ex}\) and
either lies in a small neighborhood of \(\fL\) or is contained in a
curvilinear triangle \(\fT\) on which \(\nabla H^*=(0,0)\), then it does
not belong to the overlap region \(\fO\), defined in
\eqref{e:def_overlap_region}.

\begin{definition}[Exceptional arctic set]\label{def:exA}
The \emph{exceptional arctic set}
\(\fA_{\mathrm{ex}}\subseteq\fA\) consists of:
\begin{enumerate}
\item the portions of \(\fA\) contained in curvilinear triangles on which
\[
\nabla H^*\in\{(1,0),(1,-1)\};
\]
\item the cusp points and the endpoints of extended sides lying on portions
of \(\fA\) contained in curvilinear triangles on which
\[
\nabla H^*=(0,0).
\]
\end{enumerate}
\end{definition}

\begin{proposition}[Liquid-region estimate]\label{p:liquid}
Fix \(\fb>0\), and let \(\delta>0\) be sufficiently small. Uniformly for
\((x,s),(y,t)\in\fL\cap\bZ^2/n\) satisfying
\[
\dist((x,s),\fA)\geq n^{-\delta},
\qquad
\dist((y,t),\fA)\geq n^{-\delta},
\]
and 
\[
\|(x,s)-(y,t)\|_2\geq\fb,
\qquad 
\dist((x,s),\fA_{\mathrm{ex}})\geq\fb,
\qquad
\dist((y,t),\fA_{\mathrm{ex}})\geq\fb,
\]
let \(w_{c},\overline{w_{c}}\) and \(z_c,\overline{z_c}\) be the
complex-conjugate critical points associated with \((x,s)\) and \((y,t)\),
respectively. Then
\begin{align}\begin{split}\label{e:liquid_kernel}
K^{-1}((x,s),(y,t))
=
-\frac{1}{2\pi n\ri}\Bigg(&
\bigl(1+\OO(n^{-1+2\delta})\bigr)
\frac{
\sqrt{\del_x\phi(x,s)}\sqrt{\del_y\phi(y,t)}\,
e^{n(S(w_{c};x,s)-S(z_c;y,t))}}
{\phi(w_{c})-\phi(z_c)}
\\
&-
\bigl(1+\OO(n^{-1+2\delta})\bigr)
\frac{
\overline{\sqrt{\del_x\phi(x,s)}}\sqrt{\del_y\phi(y,t)}\,
e^{n(S(\overline{w_{c}};x,s)-S(z_c;y,t))}}
{\phi(\overline{w_{c}})-\phi(z_c)}
\\
&+
\bigl(1+\OO(n^{-1+2\delta})\bigr)
\frac{
\sqrt{\del_x\phi(x,s)}\overline{\sqrt{\del_y\phi(y,t)}}\,
e^{n(S(w_{c};x,s)-S(\overline{z_c};y,t))}}
{\phi(w_{c})-\phi(\overline{z_c})}
\\
&-
\bigl(1+\OO(n^{-1+2\delta})\bigr)
\frac{
\overline{\sqrt{\del_x\phi(x,s)}}
\overline{\sqrt{\del_y\phi(y,t)}}\,
e^{n(S(\overline{w_{c}};x,s)-S(\overline{z_c};y,t))}}
{\phi(\overline{w_{c}})-\phi(\overline{z_c})}
\Bigg).
\end{split}\end{align}
\end{proposition}

\begin{proposition}[Arctic-boundary estimate]\label{p:arctic}
Fix \(\fb>0\). Uniformly for \((x,s),(y,t)\in \bZ^2/n\) satisfying
\[
\dist((x,s),\fL)\leq n^{-\delta},
\qquad
\dist((y,t),\fL)\leq n^{-\delta},
\]
and
\[
\|(x,s)-(y,t)\|_2\geq\fb,
\qquad 
\dist((x,s),\fA_{\mathrm{ex}})\geq\fb,
\qquad
\dist((y,t),\fA_{\mathrm{ex}})\geq\fb,
\]
we have
\[
\bigl|K^{-1}((x,s),(y,t))\bigr|
\lesssim  \frac{(\ln n)^{10}\sum_{w_c,z_c} e^{n\Re\left(S(w_c;x,s)-S(z_c;y,t)\right)}}{n\sqrt{ (\dist((x,s), \fA)^{1/2}+n^{-1/3})(\dist((y,t), \fA)^{1/2}+n^{-1/3})}}
.
\]
Here the sum ranges over the descent critical points \(w_c\)  and  ascent critical points \(z_c\) contained in the corresponding liquid or arctic charts.
\end{proposition}

\begin{proposition}[Frozen-region estimate]\label{p:frozen}
Fix \(\fb>0\),  a curvilinear triangle $\fT$ on which $\nabla H^*=(0,0)$, and let \(\delta>0\) be sufficiently small.  Then uniformly
for \((x,s)\in \fT\cap\bZ^2/n\) satisfying
\[
\dist((x,s),\fA)\geq n^{-\delta},\quad
\dist((x,s),\fA_{\mathrm{ex}})\geq\fb,
\]
we have
\begin{align}
&K^{-1}((x,s),(x,s))=1+\OO(
e^{-n^{1-2\delta}}),\\
&K^{-1}((x,s-1/n),(x,s)), K^{-1}((x-1/n,s-1/n),(x,s))=\OO(
e^{-n^{1-2\delta}})
\end{align}
\end{proposition}

\section{Convergence to the Gaussian Free Field}
\label{s:gff_moment_convergence}

Fix pairwise distinct points
\[
 (x_1,s_1),\ldots,(x_k,s_k)\in\fL,
 \qquad
 u_i:=\phi(x_i,s_i)\in\bC_+.
\]
Let
\[
 \cG_{\bC_+}(u,v)
 :=-\frac{1}{2\pi}
 \log\left|\frac{u-v}{u-\overline v}\right|
\]
be the Dirichlet Green function on the upper half-plane.  In this section, using \Cref{p:liquid}, \Cref{p:arctic} and \Cref{p:frozen} as input, we will prove
\begin{equation}\label{e:height_moment_raw_limit}
 \lim_{n\to\infty}
 \bE\!\left[\prod_{i=1}^k H_n^\circ(x_i,s_i)\right]
 =
 \begin{cases}
 \displaystyle
 \sum_{\pi\in\mathcal P_k}
 \prod_{\{i,j\}\in\pi}
 \frac1\pi\,\cG_{\bC_+}(u_i,u_j),
 &k\text{ even},\\[1.1em]
 0,&k\text{ odd},
 \end{cases}
\end{equation}
where \(\mathcal P_k\) is the set of pairings of
\(\{1,\ldots,k\}\).  Multiplication by \(\pi^{k/2}\) gives
\eqref{eq:moment-conv}. 
In the rest of this section, we sketch the proof of \Cref{thm:moment-GFF}, which follows from the standard kernel computations, see \cite{petrov2015asymptotics,kenyon2008height}.

\subsection{Height paths and the exact cumulant formula}

Choose a curvilinear triangle \(\fT\) on which
\(\nabla H^*=(0,0)\).  Inside \(\fA\cap\fT\), choose a compact regular
subarc which is disjoint from \(\fA_{\mathrm{ex}}\) from \eqref{def:exA}.  For each \(i\), choose a polygonal path \(\gamma_i\), made of a
bounded number of segments in the three lattice directions, which starts on
\(\partial\fP\cap\fT\), crosses this regular arctic subarc once and
transversely, and ends at \((x_i,s_i)\).  The paths can be chosen pairwise
disjoint, positively separated, and a positive distance from
\(\fA_{\mathrm{ex}}\).  Then for points on these paths, \Cref{p:liquid}, \Cref{p:arctic} and \Cref{p:frozen} hold.

We denote the set of lozenges by $\mathscr M$. If the dual edge
\(e=\rw\rb\), with \(\rw\) white and \(\rb\) blue, crosses \(\gamma_{i}\), let
\(\epsilon_i(e)\in\{-1,1\}\) be its signed contribution to the height
increment.  The height at the initial boundary point is deterministic, so
\begin{equation}\label{e:height_path_sum_gff}
 H_n^\circ(x_i,s_i)
 =\sum_{e=\rw\rb\in\gamma_{i}}
 \epsilon_i(e)
 \left(\bm1_{\{e\in\mathscr M\}}-\bP(e\in\mathscr M)\right).
\end{equation}
For the three positively oriented lattice steps, the height convention (recall from \eqref{e:height1} and \eqref{e:height2}) and
the Kasteleyn weights give
\begin{equation}\label{e:height_kasteleyn_signs}
\begin{array}{c|c|c|c|c}
\text{path step}&\rb-\rw&\epsilon_i(\rw\rb)&K(\rw,\rb)&\epsilon_i(\rw\rb)K(\rw,\rb)\\
\hline
(1/n,0)&(0,0)&-1&1&-1\\
(0,-1/n)&(-1/n,-1/n)&1&-1&-1\\
(1/n,1/n)&(0,-1/n)&1&-1&-1.
\end{array}
\end{equation}
 Reversing a path step
reverses both \(\epsilon_i(\rw\rb)\).

For distinct edges \(e_r=\rw_r\rb_r\), the dimer correlation formula is
\[
 \bP(e_1,\ldots,e_m\in\mathscr M)
 =\left(\prod_{r=1}^mK(\rw_r,\rb_r)\right)
 \det\!\left(K^{-1}(\rb_r,\rw_q)\right)_{r,q=1}^m.
\]
Equivalently,
\begin{align}\label{e:cumulant}
\bE\prod_{r=1}^m
\left(1+v_r\bm1_{\{e_r\in\mathscr M\}}\right)
=
\det\!\left(
\delta_{rq}
+
v_rK(\rw_r,\rb_r)K^{-1}(\rb_r,\rw_q)
\right)_{r,q=1}^m.
\end{align}
Taking logarithms and using the trace--log expansion
\[
\log\det(I+A)
=
\operatorname{Tr}\log(I+A)
=
\sum_{\ell\geq1}
\frac{(-1)^{\ell+1}}{\ell}\operatorname{Tr}(A^\ell),
\]
we obtain the joint cumulants of the edge indicators. Since
\(e_1,\ldots,e_m\) are distinct, the coefficient of
\(v_1\cdots v_m\) in the logarithm of the left-hand side of
\eqref{e:cumulant} is
\begin{align}\label{e:cumulant2}
&\kappa\left(
\bm1_{\{e_1\in\mathscr M\}},\ldots,
\bm1_{\{e_m\in\mathscr M\}}
\right)=
(-1)^{m-1}
\sum_{\sigma\in\mathfrak C_m}
\prod_{r=1}^m
K(\rw_r,\rb_r)K^{-1}(\rb_r,\rw_{\sigma(r)}),
\end{align}
where \(\mathfrak C_m\) denotes the set of permutations of
\(\{1,\ldots,m\}\) consisting of a single cycle of length \(m\).
Indeed, only the term with \(\ell=m\) in the trace--log expansion can
contribute the monomial \(v_1\cdots v_m\). Its sign is
\((-1)^{m+1}=(-1)^{m-1}\), and the factor \(1/m\) is canceled because each
cycle is counted \(m\) times, once for each cyclic rotation. For \(m\geq2\),
the same formula holds after centering the edge indicators.

Let \(\kappa_n(i_1,\ldots,i_m)\) be the joint cumulant of the centered
heights at \((x_{i_1},s_{i_1}), (x_{i_2},s_{i_2}),\cdots, (x_{i_m},s_{i_m})\).  Multilinearity and
\eqref{e:height_path_sum_gff} therefore give the exact identity
\begin{align}\label{e:cumulant_cycle_formula}
 \kappa_n(i_1,\ldots,i_m)
 &=(-1)^{m-1}\sum_{\sigma\in\mathfrak C_m}
 \sum_{e_1=\rw_1\rb_1\in\gamma_{i_1}}\cdots
 \sum_{e_m=\rw_m\rb_m\in\gamma_{i_m}}
 \notag\\[-0.2em]
 &\hspace{2.1cm}\times
 \prod_{r=1}^m
 \epsilon_{i_r}(e_r)K(\rw_r,\rb_r)
 K^{-1}(\rb_r,\rw_{\sigma(r)}).
\end{align}
Also \(\kappa_n(i)=0\) by centering.

\subsection{Removing the frozen part and the arctic collar}

Take \(\eta=n^{-\delta}\), truncate \(\gamma_i\) at its entering the shrunk liquid region
\(\{\zeta\in \fL: \dist(\zeta,\fA)\geq \eta\}\), retaining the terminal part ending at
\((x_i,s_i)\); denote this terminal part by \(\gamma_i^\eta\).  The paths
\(\gamma_i^\eta\) lie inside \(\fL\) and remain pairwise
separated.  Let \(\kappa_n^\eta(i_1,\ldots,i_m)\) denote
\eqref{e:cumulant_cycle_formula} with every \(\gamma_{i_r}\) replaced by
 \(\gamma_{i_r}^\eta\).  We shall prove
\begin{equation}\label{e:truncate_error}
\limsup_{n\to\infty}
 \left|\kappa_n(i_1,\ldots,i_m)
       -\kappa_n^\eta(i_1,\ldots,i_m)\right|=0.
\end{equation}

We now prove \eqref{e:truncate_error}.    On the
frozen part of a path at distance at least \(\eta=n^{-\delta}\) from \(\fA\), the frozen
kernel estimate \Cref{p:frozen} says that the occupation probability of each crossed edge
is within \(\OO(e^{-n^{1-2\delta}})\) of either \(0\) or \(1\).  Therefore its
centered indicator has \(L^1\)-norm \(\OO(e^{-n^{1-2\delta}})\).  The
moment--cumulant formula implies that a fixed-order cumulant containing one
such centered indicator is also \(\OO_m(e^{-n^{1-2\delta}})\).  Since there are only
\(\OO(n^m)\) edge tuples, all tuples containing an edge in this frozen
part contribute \(\oo(1)\).

It remains to control a two-sided collar of the arctic curve, which consists of points at distance $\eta=n^{-\delta}$ from $\fA$. We first show for edges $e_1=\rw_1\rb_1, \cdots, e_r=\rw_r\rb_r$ in this region
\begin{align}
n\sum_{r=1}^m \Re[S(w_c(\rb_r); \rb_r)-S(z_c(\rw_{\sigma(r)}); \rw_{\sigma(r)})]
=n\sum_{r=1}^m \Re[S(w_c(\rb_r); \rb_r)-S(z_c(\rw_{r}); \rw_{r})]\leq Cm
\end{align}
Here we used that the real part of every descent
critical value is no larger than the real part of every ascent critical
value from \Cref{p:critical_compare}:
\begin{align}\label{e:descent_small_ascent}
\Re[S(w_c(w_r); w_r)-S(z_c(w_{r}); w_{r})]\leq 0.
\end{align}
Moreover, from \Cref{l:critical_value_continuity}, all critical values are uniformly Lipschitz near the chosen regular arctic subarc:
\begin{align}\label{e:lipschitz}
n|\Re[S(w_c(\rb_r); \rb_r)-S(z_c(\rw_{r}); \rw_{r})]|\lesssim1.
\end{align}

 Multiplying the arctic kernel
bounds \Cref{p:arctic} around the cycle therefore gives
\begin{equation}\label{e:arctic_cycle_bound}
 \prod_{r=1}^m
 \left|K(\rw_r,\rb_r)K^{-1}(\rb_r,\rw_{\sigma(r)})\right|
 \lesssim \frac{(\ln n)^{10m}}{n^m}
 \prod_{r=1}^m
 \frac1{\dist(e_r,\fA)^{1/2}+n^{-1/3}}.
\end{equation}
Indeed, after cyclic relabeling each selected edge supplies the two
quarter-power factors belonging to its blue and white endpoints.

Because every path crosses \(\fA\) transversely, its edges in an
\(\eta\)-collar can be indexed so that their distances are comparable to
\(j/n\), up to a bounded shift.  Consequently,
\begin{align}\label{e:arctic_weight_sum}
 \frac1n
 \sum_{\substack{e\in\gamma_{i}\\
                  \dist(e,\fA)\le \eta}}
 \frac1{\dist(e,\fA)^{1/2}+n^{-1/3}}\le C\left[
 \frac1n\sum_{0\le j\le n^{1/3}}n^{1/3}
 +\frac1n\sum_{n^{1/3}<j\le C\eta n}
 \left(\frac nj\right)^{1/2}
 \right]
 \le C\eta^{1/2}.
\end{align}
The same normalized weighted sum over the whole fixed arctic neighborhood
is uniformly bounded.  Summing \eqref{e:arctic_cycle_bound} over edge tuples
and using a union bound for which selected edge lies in the collar shows that
all tuples containing a collar edge contribute
\(\OO(\eta^{1/2})\).
 Taking 
\(n\to\infty\) proves
\eqref{e:truncate_error}.

\subsection{The liquid cycle sum}

Differentiating \(x=z+s\chi(z)\) and composing with \(\phi\) gives
\begin{equation}\label{e:phi_directional_derivatives_gff}
 \partial_s\phi=-\chi\,\partial_x\phi,
 \qquad
 (\partial_x+\partial_s)\phi=(1-\chi)\partial_x\phi.
\end{equation}
We recall from \eqref{e:critical_value_derivatives_gff},
\begin{equation}\label{e:critical_value_derivatives_gff2}
 \partial_xS(w_c;x,s)=-\ln f,
 \qquad
 \partial_sS(w_c;x,s)=-\ln(1-\chi),
 \qquad
 (\partial_x+\partial_s)S(w_c;x,s)=-\ln\chi.
\end{equation}

The product of the two square-root factors in the liquid kernel \eqref{e:liquid_kernel} at the same upper  half plane critical point is
\(\partial_x\phi\); on the lower half plane it is
\(-\partial_x\overline\phi\).  Indeed, \eqref{e:critical_value_derivatives_gff2} 
 gives, for adjacent \(\rb,\rw\), and critical points both on the upper half plane 
\begin{equation}\label{e:three_local_saddle_factors}
\begin{aligned}
\sqrt{\partial_x\phi(\rb)}\sqrt{\partial_y\phi(\rw)}
 e^{n(S(w_c;\rb)-S(z_c;\rw))}
=
\begin{cases}
\partial_x\phi(\rw)+\OO(n^{-1}),
 &\rb-\rw=(0,0),\\
\chi(w)\partial_x\phi(\rw)+\OO(n^{-1}),
 &\rb-\rw=(-1/n,-1/n),\\
(1-\chi(w))\partial_x\phi(\rw)+\OO(n^{-1}),
 &\rb-\rw=(0,-1/n).
\end{cases}
\end{aligned}
\end{equation}
In the second row
\[
 n\bigl(S(w_c;\rb)-S(z_c;\rw)\bigr)
 =-(\partial_x+\partial_s)S(w_c;\rw)+\OO(n^{-1})
 =\ln\chi(\rw)+\OO(n^{-1}),
\]
and the third row is identical with \(\chi\) replaced by \(1-\chi\).
The first row has coincident coordinate labels.  The terms with both critical points in the lower half plane is
the complex conjugate of \eqref{e:three_local_saddle_factors}.

Let \(e=\rw\rb\) cross a straight portion of a truncated path, and let
\(p,p'\) be the consecutive vertices along the paths.  By
\eqref{e:height_kasteleyn_signs},
\eqref{e:phi_directional_derivatives_gff}, and Taylor expansion,
\eqref{e:three_local_saddle_factors} is exactly the statement
\begin{align}\label{e:local_differential_gff}
 &-\frac1{2\pi n\ri}\epsilon_i(e)K(\rw,\rb)
 \sqrt{\partial_x\phi(\rb)}\sqrt{\partial_y\phi(\rw)}
 e^{n(S(w_c;\rb)-S(z_c;\rw))}
 =\frac{\phi(p')-\phi(p)}{2\pi\ri}+\OO(n^{-2}).
\end{align}
For the lower branch, the right side is
\(-\bigl(\overline{\phi(p')}-\overline{\phi(p)}\bigr)/(2\pi\ri)
+\OO(n^{-2})\).

To bound the terms, some critical points from upper half plane, some from lower half plane, we shall also use the following elementary oscillation consequence of
\eqref{e:critical_value_derivatives_gff2}.  If the outgoing and incoming
saddles at a selected edge are conjugate rather than matched, cyclic
regrouping leaves the phase
\[
 e^{\pm2n\ri\Im S(w_c(x,s);x,s)}.
\]
Along the three possible straight lattice directions, by \eqref{e:critical_value_derivatives_gff2} the ratios of
successive phases are, up to a multiplicative \(1+\OO(n^{-1})\),
\[
 e^{-2\ri\arg f},
 \qquad
 e^{2\ri\arg(1-\chi)},
 \qquad
 e^{-2\ri\arg\chi},
\]
respectively.  In the liquid region, \(f,\chi\in\bC_-\) and
\(1-\chi\in\bC_+\). Moreover, the truncated paths remain at distance at least \(n^{-\delta}\)
from the regular arctic boundary and are bounded away from the exceptional
arctic set. The square-root behavior near the arctic
boundary therefore gives
\[
\min\left\{
\left|1-e^{-2\ri\arg f}\right|,
\left|1-e^{2\ri\arg(1-\chi)}\right|,
\left|1-e^{-2\ri\arg\chi}\right|
\right\}
\geq
c n^{-\delta/2}.
\]
The oscillatory sum over one path variable is
\(\oo(n)\), rather than \(\OO(n)\).

Define the oriented chain
\begin{equation}\label{e:def_Lambda_eta}
 \Lambda_i^\eta
 :=\phi(\gamma_i^\eta)-\overline{\phi(\gamma_i^\eta)}.
\end{equation}
The first part is oriented from the truncation point to \(u_i\); the
conjugate part has the opposite orientation.  As \(\eta=n^{-\delta}\) approaches $0$, the
two parts join on the real axis and form a contour from
\(\overline{u_i}\) to \(u_i\).  The paths were chosen so that these chains
are pairwise disjoint.

Fix \(m\ge2\), distinct indices \(i_1,\ldots,i_m\), and
\(\sigma\in\mathfrak C_m\).  We claim that
\begin{align}\label{e:cycle_riemann_sum}
&\lim_{n\to\infty}
 \sum_{e_1=\rw_1\rb_1\in\gamma_{i_1}^\eta}\cdots
 \sum_{e_m=\rw_m\rb_m\in\gamma_{i_m}^\eta}
 \prod_{r=1}^m
 \epsilon_{i_r}(e_r)K(\rw_r,\rb_r)
 K^{-1}(\rb_r,\rw_{\sigma(r)})
 \notag\\
&\qquad=
 \frac1{(2\pi\ri)^m}
 \int_{\Lambda_{i_1}^\eta}\cdots
 \int_{\Lambda_{i_m}^\eta}
 \prod_{r=1}^m
 \frac{\rd\zeta_r}{\zeta_r-\zeta_{\sigma(r)}}.
\end{align}

To prove this, expand the four saddle terms in every liquid kernel.  In the
cyclic product, group at each selected edge the source factor from the
outgoing kernel with the target factor from the incoming kernel.  If the two
factors use the same branch at every edge, then at each edge one may choose
independently the upper or lower branch.  By
\eqref{e:local_differential_gff}, these \(2^m\) choices are precisely the
Riemann-sum expansion of the right side of
\eqref{e:cycle_riemann_sum}: the upper choice supplies
\(\rd\phi/(2\pi\ri)\), and the lower choice supplies
\(-\rd\overline\phi/(2\pi\ri)\).

Replacing the blue or white triangle coordinate in a Cauchy denominator by
the adjacent path vertex changes it by \(\OO(n^{-1})\).  Since points
on distinct paths are uniformly separated, the resulting error is
\(\OO(n^{-m-1})\) for each edge tuple and hence
\(\OO(n^{-1})\) after summation.  Replacing one local factor in
\eqref{e:local_differential_gff} by its \(\OO(n^{-2})\) error gives the
same bound.  Tuples containing one of the finitely many path corners also
contribute \(\OO(n^{-1})\).

If the two factors are conjugate at some edge, sum first over that edge.
The preceding oscillation estimate makes this sum \(\oo(n)\); the other
\(m-1\) path sums are \(\OO(n^{m-1})\), while the \(m\) kernels contribute
\(n^{-m}\).  Thus every incompatible saddle assignment is \(\oo(1)\).
This proves \eqref{e:cycle_riemann_sum}.

Combining \eqref{e:cumulant_cycle_formula} and
\eqref{e:cycle_riemann_sum}, for \(\eta=n^{-\delta}\) we obtain
\begin{equation}\label{e:cumulant_integral_limit_eta}
\kappa_n^\eta(i_1,\ldots,i_m)
 -\frac{(-1)^{m-1}}{(2\pi\ri)^m}
 \sum_{\sigma\in\mathfrak C_m}
 \int_{\Lambda_{i_1}^\eta}\cdots
 \int_{\Lambda_{i_m}^\eta}
 \prod_{r=1}^m
 \frac{\rd\zeta_r}{\zeta_r-\zeta_{\sigma(r)}}=\oo(1).
\end{equation}

\subsection{Evaluation of the cumulants}

For \(m=2\), the only cycle is the transposition, so
\begin{equation}\label{e:two_cumulant_eta}
\kappa_n^\eta(i,j)
-\frac1{4\pi^2}
 \int_{\Lambda_i^\eta}\int_{\Lambda_j^\eta}
 \frac{\rd\zeta\,\rd\omega}{(\zeta-\omega)^2}=\oo(1).
\end{equation}
As \(\eta\downarrow0\), the two contours have endpoints
\((\overline{u_i},u_i)\) and \((\overline{u_j},u_j)\).  Since they are
disjoint, one may integrate successively using a continuous logarithm along
the contours.  Endpoint evaluation gives
\begin{align*}
 \lim_{\eta\downarrow0}
 \int_{\Lambda_i^\eta}\int_{\Lambda_j^\eta}
 \frac{\rd\zeta\,\rd\omega}{(\zeta-\omega)^2}
 &=\log\frac{(u_i-u_j)(\overline{u_i}-\overline{u_j})}
 {(u_i-\overline{u_j})(\overline{u_i}-u_j)}=2\log\left|\frac{u_i-u_j}{u_i-\overline{u_j}}\right|.
\end{align*}
The cross ratio is positive, so there is no logarithm ambiguity.  Together
with \eqref{e:truncate_error}, this proves
\begin{equation}\label{e:second_cumulant_limit}
 \lim_{n\to\infty}\kappa_n(i,j)
 =-\frac1{2\pi^2}
 \log\left|\frac{u_i-u_j}{u_i-\overline{u_j}}\right|
 =\frac1\pi\cG_{\bC_+}(u_i,u_j).
\end{equation}

For \(m\ge3\), the elementary cyclic Cauchy identity is
\begin{equation}\label{e:cauchy_cycle_identity}
 \sum_{\sigma\in\mathfrak C_m}
 \prod_{r=1}^m\frac1{\zeta_r-\zeta_{\sigma(r)}}=0
 \qquad
 (\zeta_1,\ldots,\zeta_m\text{ pairwise distinct}).
\end{equation}

For fixed \(\eta\), the chains in
\eqref{e:cumulant_integral_limit_eta} are pairwise disjoint, so
\eqref{e:cauchy_cycle_identity} applies pointwise under the integrals.  Hence
\[
\kappa_n^\eta(i_1,\ldots,i_m)=\oo(1),
 \qquad m\ge3.
\]
Using \eqref{e:truncate_error}, we conclude that
\begin{equation}\label{e:higher_cumulants_vanish}
 \lim_{n\to\infty}\kappa_n(i_1,\ldots,i_m)=0,
 \qquad m\ge3.
\end{equation}

Finally, the moment--cumulant formula gives
\[
 \bE\!\left[\prod_{i=1}^kH_n^\circ(x_i,s_i)\right]
 =\sum_{\Pi}\prod_{B\in\Pi}\kappa_n(B),
\]
where \(\Pi\) ranges over the set partitions of \(\{1,\ldots,k\}\).
Singleton blocks vanish by centering, blocks of size at least three vanish by
\eqref{e:higher_cumulants_vanish}, and a block \(\{i,j\}\) of size two
converges to \(\pi^{-1}\cG_{\bC_+}(u_i,u_j)\) by
\eqref{e:second_cumulant_limit}.  Only pairings survive, which proves
\eqref{e:height_moment_raw_limit} and hence \eqref{eq:moment-conv}.

\chapter{Compatibility of Approximate Kernels}

\section{Estimates of the Contour Integrals}\label{s:cintegral}
In this section we derive estimates for the single-contour integral \eqref{e:single_term} and the double-contour integral \eqref{e:all_term}.

\subsection{Nodal length bound}

Let $f$ be holomorphic in the disk $B_R:=\{z\in\mathbb C:|z|<R\}$. Define the \emph{vector-valued (Almgren) frequency}
\begin{equation}\label{eq:freq-def}
\mathcal N_f(r):=\frac{r\,D(r)}{H(r)}
=\frac{2r\displaystyle\int_{B_r}|f'(z)|^2}{\displaystyle\int_{\partial B_r}|f-f(0)|^2},
\qquad 0<r<R.
\end{equation}
It is classical that $\mathcal N_f(r)$ is nondecreasing in $r$, and that in terms of the Taylor series
$f(z)=\sum_{k\ge0} a_k z^k$ one has the exact identity
\begin{equation}\label{eq:freq-avg}
\mathcal N_f(r)=
\frac{\sum_{k\ge1} k\,|a_k|^2 r^{2k}}{\sum_{k\ge1} |a_k|^2 r^{2k}},
\end{equation}
i.e. $\mathcal N_f(r)$ is the weighted average of degrees present at scale $r$.

\medskip

For any $c\in \bR$, write $\Gamma_c:=\{z\in B_R:\Im f(z)=c\}$ for the nodal set of $\Im f-c$. We denote by $\mathcal H^1$ the one-dimensional Hausdorff measure (arc length). The following statement that the nodal length can be controlled by frequency, which is a special case of \cite[Theorem 3.1']{lin1991nodal}. 
\begin{lemma}\label{lem:nodal-length-frequency}
There exists a universal constant $C>0$ such that for every $0<\rho\le R/2$,
\begin{equation}\label{eq:nodal-length}
\mathcal H^1\!\big(\Gamma_c\cap B_\rho\big)\;\leq\; C\,\rho\,\mathcal N_f(R).
\end{equation}
\end{lemma}

\begin{proof}
Take $u(z)=\Im[f(2\rho z)]-c$, which is a harmonic function on $B_1$. It follows from \cite[Theorem 3.1']{lin1991nodal} (by taking $n=2$ and $g$ is the Euclidean metric), 
\begin{align}
\mathcal H^1\!\big(z\in B_1: u(z)=0\big)\;\leq\; C \frac{\int_{B_1} |\nabla u|^2}{\int_{\del B_1}|u|^2}
\end{align}
We notice that
\begin{align}
\int_{B_1} |\nabla u|^2
=\int_{B_1}|\nabla \Im[f(2\rho z)]|^2
=\int_{B_{2\rho}}|\nabla \Im[f(z)]|^2= \int_{B_{2\rho}}|f'(z)|^2,
\end{align}
and 
\begin{align}
\int_{\del B_1}|u|^2
=\int_{\del B_1}| \Im[f(2\rho z)-f(0)]|^2
=\frac{1}{2\rho}\int_{\del B_{2\rho}}|\Im[f( z)-c]|^2=\frac{1}{4\rho} \int_{\del B_{2\rho}}|f|^2
\end{align}
The statement \eqref{eq:nodal-length} follows 
\begin{align}
\mathcal H^1\!\big(\Gamma_c\cap B_\rho\big)=\rho \mathcal H^1\!\big(z\in B_1: u(z)=0\big)\;\leq\; C\rho \frac{\int_{B_1} |\nabla u|^2}{\int_{\del B_1}|u|^2}
=C\rho \cN_f(2\rho)\leq C\rho\cN_f(R).
\end{align}

\end{proof}

\begin{lemma}[Uniform local length bound for two logarithmic singularities]\label{l:length_bound}
Let $g$ be holomorphic and nonconstant in a neighborhood of $0$. Then there exist
$r_g>0$ and $C_g<\infty$, depending only on $g$, such that for every
$0<\rho<r_g$, every $\varepsilon\in\mathbb C\setminus\{0\}$, and every $c\in\mathbb R$,
\[
\mathcal H^1\!(
\Gamma_{\varepsilon, c}
)
\le C_g \rho ,\quad \Gamma_{\varepsilon, c}:=\left\{
z\in B_\rho\setminus\{0,-\varepsilon\}:
\Im\left[\log\frac{z}{z+\varepsilon}+g(z)\right]=c
\right\}
\]
Here the logarithm may be taken in any local branch; changing the branch only
changes $c$ by a multiple of $2\pi$.
\end{lemma}

\begin{proof}
First suppose that $|\varepsilon|\ge 4\rho$. Then, in $B_\rho$,
\[
\log\frac{z}{z+\varepsilon}+g(z)=\log z+G_\varepsilon(z),
\qquad
G_\varepsilon(z):=g(z)-\log(z+\varepsilon),
\]
with $G_\varepsilon$ holomorphic. For $|z|<\rho$,
\[
|zG_\varepsilon'(z)|
\le |z g'(z)|+\frac{|z|}{|z+\varepsilon|}
\le \oo_{r_g}(1)+\frac13
\le \frac12 .
\]
Writing $z=se^{\ri\theta}$, the level equation is
\[
\theta+\Im G_\varepsilon(se^{\ri\theta})=c.
\]
Since
\[
\partial_\theta\big(\theta+\Im G_\varepsilon(se^{\ri\theta})\big)
=
1+\Im\!\left(\ri zG_\varepsilon'(z)\right)\big|_{z=se^{\ri\theta}}
\ge \frac12,
\]
the level set consists of a uniformly bounded number of radial graphs
$\theta=\theta(s)$. Differentiating implicitly gives
\[
|s\theta'(s)|\le C,
\]
and therefore its length in $B_\rho$ is at most
\[
\int_0^\rho \sqrt{1+s^2|\theta'(s)|^2}\,ds
\le C\rho .
\]

It remains to consider the case $|\varepsilon|<4\rho$. Put
$a:=|\varepsilon|$. We cover $0$ and $-\varepsilon$ by two balls
$B_{a/4}(0)$ and $B_{a/4}(-\varepsilon)$. By the first case, and by the same
argument after translating near $-\varepsilon$,
\[
\mathcal H^{1}(\Gamma_{\varepsilon,c}\cap B_{a/4}(0)),
\quad
\mathcal H^{1}(\Gamma_{\varepsilon,c}\cap B_{a/4}(-\varepsilon))
\le C a .
\]

Outside the two balls, we decompose
\[
B_\rho\setminus \big(B_{a/4}(0)\cup B_{a/4}(-\varepsilon)\big)
\subset
\bigcup_j A_j,
\]
where
\[
A_j:=\left\{z:\dist(z,\{0,-\varepsilon\})\in [s_j,2s_j]\right\},
\qquad
s_j:=2^{j-2}a,
\]
and the union is over those $j\ge0$ with $s_j<\rho$.

We can then cover $A_j$ with a uniformly finite number of balls
$B_{s_j/4}(z_0)$. Moreover, the enlarged balls $B_{s_j/2}(z_0)$ have distance
at least $s_j/2$ from $\{0,-\varepsilon\}$.

One can then check that, after translating to $z_0$ and choosing a local branch,
the frequency of
\[
\log\frac{z}{z+\varepsilon}+g(z)
\]
is uniformly bounded on the enlarged ball $B_{s_j/2}(z_0)$. Applying
\eqref{eq:nodal-length} gives
\[
\mathcal H^{1}(\Gamma_{\varepsilon,c}\cap B_{\kappa s_j}(z_0))
\le C' s_j .
\]
Since only uniformly many balls are needed to cover $A_j$, it follows that
\[
\mathcal H^{1}(\Gamma_{\varepsilon,c}\cap A_j)\le C s_j .
\]

Combining this with the core estimate gives
\[
\mathcal H^1(\Gamma_{\varepsilon,c}\cap B_\rho)
\le C|\varepsilon|+C\sum_{j:s_j<\rho}s_j
\le 8C\rho,
\]
because $|\varepsilon|<4\rho$. The proof is complete.
\end{proof}

\subsection{Single-contour integrals}

\begin{proposition}\label{p:single_denominator_bound}
Fix $\xi\in \bC$, and let
$
\sfD(\xi)=\sfD(\xi;(x,s),(y,t))
$
be a steepest-descent path for
$
\Re\bigl[S(z;x,s)-S(z;y,t)\bigr],
$
as recalled in \Cref{s:pof_single_integral}. Then
\begin{align}\label{e:single_denominator_bound}
\int_{\sfD(\xi)}
\frac{n\,|\rd z|}
{
\sqrt{
(1+n|x-z|)
(1+n|z-(x-s)|)
(1+n|y-z|)
(1+n|z-(y-t)|)
}
}
\lesssim
\frac{\ln n}{n}.
\end{align}
\end{proposition}

\begin{proof}
Let
\[
\mathcal V:=\{x,x-s,y,y-t\}.
\]
Choose \(R>0\) sufficiently large so that \(B_R(0)\) contains
\(\mathcal V\).

Since
\[
|x-(x-s)|=s,
\qquad
|y-(y-t)|=t,
\]
and \(s,t>0\) are fixed, for every \(z\), at least one of
\(|z-x|\) and \(|z-(x-s)|\) is bounded below by a positive constant, and
the same holds for \(|z-y|\) and \(|z-(y-t)|\). Consequently, on
\(B_R(0)\),
\begin{align}\label{e:single_endpoint_comparison}
\frac{n^2}{
\sqrt{
(1+n|x-z|)
(1+n|z-(x-s)|)
(1+n|y-z|)
(1+n|z-(y-t)|)
}
}
\lesssim
\sum_{a\in\mathcal V}
\frac{n}{1+n|z-a|}.
\end{align}

Next we show that for each \(a\in\mathcal V\), a dyadic decomposition around \(a\),  gives
\begin{align}\label{e:single_endpoint_bound}
\int_{\sfD(\xi)\cap B_R(0)}
\frac{n\,|\rd z|}{1+n|z-a|}
\lesssim
\ln n.
\end{align}
Indeed, the portion of the contour on which
\[
\frac{2^{j-1}}{n}<|z-a|\leq\frac{2^j}{n}
\]
has length \(\OO(2^j/n)\), while the integrand is \(\OO(n/2^j)\). Each
nonempty dyadic region therefore contributes \(\OO(1)\), and there are
\(\OO(\ln n)\) such regions. It follows from
\eqref{e:single_endpoint_comparison} that the contribution from
\(\sfD(\xi)\cap B_R(0)\) is \(\OO(\ln n)\).

It remains to control the unbounded portion of the contour, if present.
Outside \(B_R(0)\), one has
\[
|z-a|\asymp |z|,
\qquad
a\in\mathcal V.
\]
Therefore,
\[
\frac{n^2}{
\sqrt{
(1+n|x-z|)
(1+n|z-(x-s)|)
(1+n|y-z|)
(1+n|z-(y-t)|)
}
}
\lesssim
\frac{1}{|z|^2}.
\]
Another dyadic decomposition, together with
\(\Cref{l:length_bound}\), gives
\[
\int_{\sfD(\xi)\setminus B_R(0)}
\frac{|\rd z|}{|z|^2}
\lesssim 1.
\]
Combining the bounded and unbounded portions proves
\eqref{e:single_denominator_bound}.
\end{proof}

\begin{lemma}\label{l:single_contour_integral}
Let \(\xi\in\bC\) satisfy
\[
\dist(\xi,[y,x])\gtrsim\frac{1}{n}
\qquad\text{if }x\geq y,
\]
and
\[
\dist(\xi,[x-s,y-t])\gtrsim\frac{1}{n}
\qquad\text{if }y-t\geq x-s.
\]
Let
$
\sfD(\xi)=\sfD(\xi;(x,s),(y,t))
$
be a steepest-descent path for
$
\Re\bigl[S(z;x,s)-S(z;y,t)\bigr],
$
as recalled in \Cref{s:pof_single_integral}. If \(\xi\in\bR\), then the
corresponding steepest-descent path is understood as a boundary value
obtained by perturbing \(\xi\) infinitesimally into the upper or lower half-plane.
Then
\begin{align}\label{e:single_contour_integral}
\left|
\frac{n}{2\pi\ri}
\int_{\sfD(\xi)}
P_{ns}(nz,nx)Q_{nt}(nz,ny)\,\rd z
\right|
\lesssim
(\ln n)
e^{n\Re[S(\xi;x,s)-S(\xi;y,t)]}.
\end{align}
\end{lemma}

\begin{lemma}\label{l:intersection_descent}
Let
\[
\xi\in
\sfD^{\rm d}(w_c)\cap\sfD^{\rm a}(z_c)\setminus\bR
\]
be a nonreal intersection point of the steepest-descent path
\(\sfD^{\rm d}(w_c)\) and the steepest-ascent path
\(\sfD^{\rm a}(z_c)\). Let
$
\sfD(\xi)=\sfD(\xi;(x,s),(y,t))
$
be a steepest-descent path for
$
\Re\bigl[S(z;x,s)-S(z;y,t)\bigr],
$
as recalled in \Cref{s:pof_single_integral}.  If \(x\geq y\), then
\begin{align}\label{e:distxi}
\dist(\xi,[y,x])
\gtrsim
\frac{1}{n}.
\end{align}
If \(y-t\geq x-s\), then
\begin{align}\label{e:distxi_shifted}
\dist(\xi,[x-s,y-t])
\gtrsim
\frac{1}{n}.
\end{align}
This verifies the assumptions in \Cref{l:single_contour_integral}, and hence
\begin{align}\label{e:single_contour_integral2}
\left|
\frac{n}{2\pi\ri}
\int_{\sfD(\xi)}
P_{ns}(nz,nx)Q_{nt}(nz,ny)\,\rd z
\right|
\lesssim
(\ln n)
e^{n\Re[S(\xi;x,s)-S(\xi;y,t)]}\leq (\ln n)
e^{n\Re[S(w_c;x,s)-S(z_c;y,t)]}.
\end{align}
\end{lemma}

\begin{proof}[Proof of \Cref{l:single_contour_integral}]
Since \(\sfD(\xi)\) is a steepest-descent path for
$
\Re\bigl[S(z;x,s)-S(z;y,t)\bigr],
$
one has
\begin{align}\label{e:single_phase_bound}
\Re\bigl[S(z;x,s)-S(z;y,t)\bigr]
\leq
\Re\bigl[S(\xi;x,s)-S(\xi;y,t)\bigr]
\end{align}
for every \(z\in\sfD(\xi)\).

The \(n^{-1}\)-scale pole-avoiding perturbations used below are made only
at points where the phase difference is analytic. Hence, they change the
real part of the phase by \(\OO(n^{-1})\), and therefore change the
exponential factor in \eqref{e:single_contour_integral} by at most a
constant factor.

We first show that, after these \(n^{-1}\)-scale pole-avoiding
perturbations,
\begin{align}\label{e:single_PQ_global}
|P_{ns}(nz,nx)Q_{nt}(nz,ny)|
\lesssim
\frac{
n e^{n\Re[S(z;x,s)-S(z;y,t)]}
}{
\sqrt{
(1+n|x-z|)
(1+n|z-(x-s)|)
(1+n|y-z|)
(1+n|z-(y-t)|)
}
}.
\end{align}
Away from the two intervals
\[
[\min\{x,y\},\max\{x,y\}],
\qquad
[\min\{x-s,y-t\},\max\{x-s,y-t\}],
\]
the estimate \eqref{e:single_PQ_global} follows directly from
\eqref{e:PQess}. It remains to consider neighborhoods of these intervals.
We prove \eqref{e:single_PQ_global} for \(z\) close to
\([\min\{x,y\},\max\{x,y\}]\). The other interval
\([\min\{x-s,y-t\},\max\{x-s,y-t\}]\) is treated in the same way, so we
omit it.

Suppose first that \(x<y\). Then \eqref{e:PQubb} gives
\eqref{e:single_PQ_global} in a neighborhood of \([x,y]\).

Suppose next that \(x>y\). In this case,
\begin{align}
\partial_z\bigl[S(z;x,s)-S(z;y,t)\bigr]
=
\ln\frac{z-x}{z-y}
+
\ln\frac{z-(y-t)}{z-(x-s)}.
\end{align}
The second logarithm is analytic in a neighborhood of \([y,x]\) and
satisfies the symmetry under complex conjugation required in
\Cref{l:repulsion}. Hence, for a sufficiently small fixed \(c>0\), the
steepest-descent vector field points outward on the boundary of
\[
\left\{
z\in\bC:
\dist(z,[y,x])\leq\frac{c}{n}
\right\}.
\]
Since the starting point lies outside this region, uniqueness of the
gradient flow implies that
\begin{align}\label{e:upper_source_separation}
\dist(z,[y,x])
\geq
\frac{c}{n},
\qquad
z\in\sfD(\xi).
\end{align}
Thus, in the case \(x>y\), the additional factor in
\eqref{e:PQubb2} satisfies
\[
\min\left\{
1,
n\dist(z,[y,x])
\right\}^{-1}
\lesssim 1
\]
along \(\sfD(\xi)\). Therefore, \eqref{e:PQubb2} again gives
\eqref{e:single_PQ_global}.

It remains to treat the degenerate case \(x=y\). The integrand has a
simple pole at \(z=x\), whereas the phase difference is analytic there.
Indeed,
\[
S(z;x,s)-S(z;x,t)
\]
extends analytically across \(z=x\) and vanishes at \(z=x\). We recall from
\Cref{s:pof_single_integral} that, if \(\xi\in\bR\), the corresponding
steepest-descent path is understood after an infinitesimal perturbation of
\(\xi\) into the upper half-plane and does not pass through \(z=x\).

If the path comes within distance \(c/n\) of \(x\), we deform the
corresponding portion to an arc of the circle \(|z-x|=c/n\) in the same
open half-plane. Since the phase difference is analytic at \(z=x\), on the
detour arc,
\[
\Re\bigl[S(z;x,s)-S(z;x,t)\bigr]
\leq
\Re\bigl[S(\xi;x,s)-S(\xi;x,t)\bigr]
+
\OO\left(\frac{1}{n}\right).
\]
Thus,
\begin{align}\label{e:detour_phase_bound}
e^{n\Re[S(z;x,s)-S(z;x,t)]}
\lesssim
e^{n\Re[S(\xi;x,s)-S(\xi;x,t)]}
\end{align}
on the detour arc. We continue to denote the resulting contour by
\(\sfD(\xi)\). Since this contour stays at distance
\(\gtrsim n^{-1}\) from \(x\), \eqref{e:single_PQ_global} holds along the
resulting contour.

Using \eqref{e:single_PQ_global}, \eqref{e:single_phase_bound}, and
\eqref{e:detour_phase_bound}, we obtain
\begin{align}
&\phantom{{}={}}\left|
\frac{n}{2\pi\ri}
\int_{\sfD(\xi)}
P_{ns}(nz,nx)Q_{nt}(nz,ny)\,\rd z
\right|\\
&\lesssim
\int_{\sfD(\xi)}
\frac{n^2
e^{n\Re[S(\xi;x,s)-S(\xi;y,t)]}|\rd z|}
{
\sqrt{
(1+n|x-z|)
(1+n|z-(x-s)|)
(1+n|y-z|)
(1+n|z-(y-t)|)
}
}.
\end{align}
By \Cref{p:single_denominator_bound}, whose proof applies unchanged to the
locally deformed contour,
\begin{align}
&\phantom{{}={}}n^2
\int_{\sfD(\xi)}
\frac{|\rd z|}
{
\sqrt{
(1+n|x-z|)
(1+n|z-(x-s)|)
(1+n|y-z|)
(1+n|z-(y-t)|)
}
}\\
&=
n
\int_{\sfD(\xi)}
\frac{n\,|\rd z|}
{
\sqrt{
(1+n|x-z|)
(1+n|z-(x-s)|)
(1+n|y-z|)
(1+n|z-(y-t)|)
}
}
\lesssim
\ln n.
\end{align}
This proves \eqref{e:single_contour_integral}.
\end{proof}

\begin{proof}[Proof of \Cref{l:intersection_descent}]
We prove \eqref{e:distxi}; the proof of
\eqref{e:distxi_shifted} is analogous.

When \(\dist(\xi,[y,x])\gtrsim1\), there is nothing to prove. Otherwise,
the concentric charts containing \(\sfD^{\rm d}(w_c)\) and
\(\sfD^{\rm a}(z_c)\) are of vertical tangent, cusp-turning, or tangent
frozen type. Let \(b_i\) denote the vertical tangency location associated
with the relevant charts. Recall from
\Cref{l:vertical_tangent_S} that, if
$
x\geq b_i+{1}/({2n}),
$
then
\begin{align}\label{e:Dist1}
\dist\bigl(
\sfD^{\rm d}(w_c),[b_i,x]
\bigr)
\gtrsim
\frac{1}{n}.
\end{align}
Similarly, if
$
y\leq b_i-{1}/({2n}),
$
then
\begin{align}\label{e:Dist2}
\dist\bigl(
\sfD^{\rm a}(z_c),[y,b_i]
\bigr)
\gtrsim
\frac{1}{n}.
\end{align}

Since \(nx,ny\in\bZ\) and \(nb_i\in\bZ+1/2\), neither \(x\) nor \(y\)
can equal \(b_i\). We distinguish three cases.

If \(b_i<y\leq x\), then
$
[y,x]\subset[b_i,x].
$
Since \(\xi\in\sfD^{\rm d}(w_c)\), \eqref{e:Dist1} gives
\[
\dist(\xi,[y,x])
\geq
\dist(\xi,[b_i,x])
\gtrsim
\frac{1}{n}.
\]

If \(y\leq x<b_i\), then
$
[y,x]\subset[y,b_i].
$
Since \(\xi\in\sfD^{\rm a}(z_c)\), \eqref{e:Dist2} gives
\[
\dist(\xi,[y,x])
\geq
\dist(\xi,[y,b_i])
\gtrsim
\frac{1}{n}.
\]

Finally, if \(y<b_i<x\), then
$
[y,x]=[y,b_i]\cup[b_i,x].
$
Since
$
\xi\in
\sfD^{\rm d}(w_c)\cap\sfD^{\rm a}(z_c),
$
the estimates \eqref{e:Dist1} and \eqref{e:Dist2} imply
\[
\dist(\xi,[y,x])
=
\min\left\{
\dist(\xi,[y,b_i]),
\dist(\xi,[b_i,x])
\right\}
\gtrsim
\frac{1}{n}.
\]
This proves \eqref{e:distxi}.

The proof of \eqref{e:distxi_shifted} is identical, using the corresponding
estimates in
\Cref{l:unit_slope_tangent_S} and decomposing
\([x-s,y-t]\) according to the position of the unit-slope tangent location
\(a_i\).
\end{proof}

\subsection{Double-contour integrals}

We first record one-contour estimates that will be used repeatedly. We recall the localization length $\Delta(\cdot)$ from \Cref{def:localization_length}

\begin{lemma}\label{l:contour_integral}
Let \(\sfD^{\rm d}(w_c)\) be a steepest-descent contour contained in a
chart \(\fU\).

If \(\fU\) is centered at a finite point, then
\begin{align}\label{e:one_contour_finite}
\int_{\sfD^{\rm d}(w_c)}
\frac{\sqrt n\,|\rd w|}
{\sqrt{(1+n|x-w|)(1+n|w-(x-s)|)}}
\lesssim
\frac{1}{\sqrt n}.
\end{align}
If, in addition, \(\sfS^{\rm d}(w_c)\) is the truncated version of
\(\sfD^{\rm d}(w_c)\), then
\begin{align}\label{e:one_contour_finite2}
\int_{\sfS^{\rm d}(w_c)}
\frac{\sqrt n\,|\rd w|}
{\sqrt{(1+n|x-w|)(1+n|w-(x-s)|)}}
\lesssim
\frac{(\ln n)^5\Delta(w_c)}{\sqrt n}.
\end{align}
If \(\fU\) is a liquid, arctic, regular frozen, or cusp frozen chart, then
the stronger bound
\begin{align}\label{e:one_contour_finite_regular}
\int_{\sfS^{\rm d}(w_c)}
\frac{\sqrt n\,|\rd w|}
{\sqrt{(1+n|x-w|)(1+n|w-(x-s)|)}}
\lesssim
\frac{\Delta(w_c)}{\sqrt n}
\end{align}
holds.

Suppose that \(\fU\) is centered at \(\infty_i\), for some \(i\). After
passing to the local coordinate \(\wt w\) at \(\infty_i\), the contour
\(\sfD^{\rm d}(w_c)\) becomes
\(\wt\sfD^{\rm d}(\wt w_c)\), and
\begin{align}\label{e:one_contour_infinity}
\int_{\wt\sfD^{\rm d}(\wt w_c)}
\frac{|\rd\wt w|}{\sqrt n}
\lesssim
\frac{1}{\sqrt n}.
\end{align}
If \(\wt\sfS^{\rm d}(\wt w_c)\) is the truncated version of
\(\wt\sfD^{\rm d}(\wt w_c)\), then
\begin{align}\label{e:one_contour_infinity2}
\int_{\wt\sfS^{\rm d}(\wt w_c)}
\frac{|\rd\wt w|}{\sqrt n}
\lesssim
\frac{(\ln n)^5\Delta(w_c)}{\sqrt n}.
\end{align}

The analogous estimates hold for a steepest-ascent contour
\(\sfD^{\rm a}(z_c)\) and its truncated version
\(\sfS^{\rm a}(z_c)\), with \(x,w,s\) replaced by \(y,z,t\).
\end{lemma}

\begin{proof}
Suppose first that \(\fU\) is centered at a finite point. Since
$
|x-(x-s)|=s
$
and \(s>0\) is fixed, at least one of \(|x-w|\) and
\(|w-(x-s)|\) is bounded below by \(s/2\). On the portion of the contour
where
$
|x-w|\leq |w-(x-s)|,
$
one has
$
|w-(x-s)|\geq {s}/{2},
$
and hence
\begin{align}
\frac{\sqrt n}
{\sqrt{(1+n|x-w|)(1+n|w-(x-s)|)}}
\lesssim
\frac{1}{\sqrt{1+n|x-w|}}.
\end{align}
On the remaining portion of the contour,
\begin{align}
\frac{\sqrt n}
{\sqrt{(1+n|x-w|)(1+n|w-(x-s)|)}}
\lesssim
\frac{1}{\sqrt{1+n|w-(x-s)|}}.
\end{align}
It is therefore enough to show that, for
\[
a\in\{x,x-s\},
\]
\begin{align}\label{e:one_endpoint_integral}
\int_{\sfD^{\rm d}(w_c)}
\frac{|\rd w|}{\sqrt{1+n|w-a|}}
\lesssim
\frac{1}{\sqrt n}.
\end{align}

Decompose the contour into the sets on which
\[
|w-a|\leq\frac1n
,\qquad 
\frac{2^{j-1}}{n}
<
|w-a|
\leq
\frac{2^j}{n},
\qquad
j\geq1.
\]
By \Cref{l:length_bound}, the portion of the contour in the \(j\)-th
region has length \(\OO(2^j/n)\), while
\[
\frac{1}{\sqrt{1+n|w-a|}}
\lesssim
2^{-j/2}
\]
there. Since \(\fU\) is a fixed bounded chart, the largest nonempty
dyadic index satisfies \(2^j\lesssim n\). Therefore,
\begin{align}
\int_{\sfD^{\rm d}(w_c)}
\frac{|\rd w|}{\sqrt{1+n|w-a|}}
&\lesssim
\frac1n
+
\frac1n\sum_j2^{j/2}
\lesssim
\frac1{\sqrt n}.
\end{align}
This proves \eqref{e:one_contour_finite}.

Suppose next that \(\fU\) is a liquid, arctic, regular-frozen, or
cusp-frozen chart. By the construction of the truncated contour,
\[
\operatorname{length}\bigl(\sfS^{\rm d}(w_c)\bigr)
\lesssim
\Delta(w_c).
\]
Moreover, \(\sfS^{\rm d}(w_c)\) remains uniformly bounded away from
\(x\) and \(x-s\). Hence
\[
\frac{\sqrt n}
{\sqrt{(1+n|x-w|)(1+n|w-(x-s)|)}}
\lesssim
\frac1{\sqrt n},
\qquad
w\in\sfS^{\rm d}(w_c),
\]
and therefore
\[
\int_{\sfS^{\rm d}(w_c)}
\frac{\sqrt n\,|\rd w|}
{\sqrt{(1+n|x-w|)(1+n|w-(x-s)|)}}
\lesssim
\frac{\Delta(w_c)}{\sqrt n}.
\]
This proves \eqref{e:one_contour_finite_regular}.

It remains to consider a tangent frozen chart centered at a finite point.
We treat a vertical tangent frozen chart centered at \(b_i\); the
unit-slope case is identical after replacing \(x\) by \(x-s\) and \(b_i\)
by \(a_i\). The contour \(\sfS^{\rm d}(w_c)\) is uniformly bounded away
from \(x-s\), and hence
\begin{align}\label{e:vt_integral}
\int_{\sfS^{\rm d}(w_c)}
\frac{\sqrt n\,|\rd w|}
{\sqrt{(1+n|x-w|)(1+n|w-(x-s)|)}}
\lesssim
\int_{\sfS^{\rm d}(w_c)}
\frac{|\rd w|}{\sqrt{1+n|x-w|}}.
\end{align}
We recall the construction of \(\sfS^{\rm d}(w_c)\) from
\Cref{l:local_descent_deformation}.

If
\[
|x-b_i|
\leq
\frac{(\ln n)^3}{n},
\]
then
\[
\operatorname{length}\bigl(\sfS^{\rm d}(w_c)\bigr)
\lesssim
\frac{(\ln n)^5}{n}.
\]
Thus,
\begin{align}
\eqref{e:vt_integral}
\lesssim
\frac{(\ln n)^5}{n}
\lesssim
\frac{(\ln n)^5\Delta(w_c)}{\sqrt n},
\end{align}
where we used $\Delta(w_c)=n^{-1/2}$.

If
\[
|x-b_i|
>
\frac{(\ln n)^3}{n},
\]
then
\[
\dist\bigl(\sfS^{\rm d}(w_c),x\bigr)
\asymp
|x-b_i|
\]
and
\[
\operatorname{length}\bigl(\sfS^{\rm d}(w_c)\bigr)
\lesssim
(\ln n)^3
\left(\frac{|x-b_i|}{n}\right)^{1/2}.
\]
Consequently,
\begin{align}
\eqref{e:vt_integral}
&\lesssim
\frac{
\operatorname{length}\bigl(\sfS^{\rm d}(w_c)\bigr)
}{
\sqrt{n|x-b_i|}
}
\lesssim
\frac{(\ln n)^3}{n}
\lesssim
\frac{(\ln n)^5\Delta(w_c)}{\sqrt n}.
\end{align}
This proves \eqref{e:one_contour_finite2}.

Suppose finally that \(\fU\) is centered at \(\infty_i\). By
\Cref{l:length_bound}, the transformed contour
\(\wt\sfD^{\rm d}(\wt w_c)\) has uniformly bounded length. Therefore,
\[
\int_{\wt\sfD^{\rm d}(\wt w_c)}
\frac{|\rd\wt w|}{\sqrt n}
\lesssim
\frac1{\sqrt n},
\]
which proves \eqref{e:one_contour_infinity}.

The construction of the truncated contour gives
\[
\operatorname{length}
\bigl(\wt\sfS^{\rm d}(\wt w_c)\bigr)
\lesssim
(\ln n)^5\Delta(w_c).
\]
Consequently,
\[
\int_{\wt\sfS^{\rm d}(\wt w_c)}
\frac{|\rd\wt w|}{\sqrt n}
\lesssim
\frac{(\ln n)^5\Delta(w_c)}{\sqrt n},
\]
which proves \eqref{e:one_contour_infinity2}. The ascent estimates follow
in the same way.
\end{proof}

We next record the behavior of the Cauchy factor in local coordinates near
the points at infinity.

\begin{lemma}\label{l:phi_local_coordinate_infinity}
Let \(w\) and \(z\) lie in a sufficiently small chart centered at
\(\infty_i\), and let \(\wt w\) and \(\wt z\) be the corresponding local
coordinates, vanishing at \(\infty_i\). Then
\begin{align}\label{e:phi_local_coordinate_infinity}
\frac{
\sqrt{|\phi'(w)|\,|\phi'(z)|}\,
|\rd w|\,|\rd z|
}{
|w|\,|z|\,|\phi(w)-\phi(z)|
}
\asymp
\frac{
|\rd\wt w|\,|\rd\wt z|
}{
|\wt w-\wt z|
}.
\end{align}

Suppose that \(w\) lies in a sufficiently small chart centered at
\(\infty_i\), while \(z\) lies in a finite chart disjoint from this chart.
Then
\begin{align}\label{e:phi_local_coordinate_infinity_disjoint}
\frac{
\sqrt{|\phi'(w)|\,|\phi'(z)|}\,
|\rd w|\,|\rd z|
}{
|w|\,|\phi(w)-\phi(z)|
}
\asymp
|\rd\wt w|\,|\rd z|.
\end{align}

Finally, suppose that \(w\) and \(z\) lie in two disjoint charts centered
at \(\infty_i\) and \(\infty_j\), respectively. Then
\begin{align}\label{e:phi_local_coordinate_two_infinities}
\frac{
\sqrt{|\phi'(w)|\,|\phi'(z)|}\,
|\rd w|\,|\rd z|
}{
|w|\,|z|\,|\phi(w)-\phi(z)|
}
\asymp
|\rd\wt w|\,|\rd\wt z|.
\end{align}
\end{lemma}

\begin{proof}
Since \(\wt w\) is a local coordinate vanishing at \(\infty_i\), one has
\begin{align}\label{e:infinity_coordinate_relation}
|w|
\asymp
\frac{1}{|\wt w|},
\qquad
|\rd w|
\asymp
\frac{|\rd\wt w|}{|\wt w|^2}.
\end{align}
The same estimates hold for \(z\) and \(\wt z\).

We first prove \eqref{e:phi_local_coordinate_infinity}. Suppose that
\(i=1\), so that \(p_1=\infty\). By \eqref{e:inftoinf},
\[
|\phi(w)|
\asymp
|w|,
\qquad
|\phi'(w)|
\asymp1,
\]
and similarly for \(z\). Moreover,
\[
\left|
\frac1{\phi(w)}
-
\frac1{\phi(z)}
\right|
\asymp
|\wt w-\wt z|,
\qquad
\left|\frac1{\phi(w)}\right|
\asymp
|\wt w|.
\]
It follows that
\begin{align}\label{e:phi_difference_infinity1}
|\phi(w)-\phi(z)|
\asymp
\frac{
|\wt w-\wt z|
}{
|\wt w|\,|\wt z|
}.
\end{align}
Combining this with \eqref{e:infinity_coordinate_relation} gives
\begin{align}
\frac{
\sqrt{|\phi'(w)|\,|\phi'(z)|}\,
|\rd w|\,|\rd z|
}{
|w|\,|z|\,|\phi(w)-\phi(z)|
}
\asymp
\frac{|\rd\wt w|\,|\rd\wt z|}
{|\wt w-\wt z|}.
\end{align}

Suppose next that \(2\leq i\leq d\), so that \(p_i\in\bR\). By
\eqref{e:inftofinite},
\[
|\phi(w)-p_i|
\asymp
|\wt w|,
\qquad
|\phi'(w)|
\asymp
|\wt w|^2,
\]
and similarly for \(z\). Hence
\[
|\phi(w)-\phi(z)|
\asymp
|\wt w-\wt z|,
\]
and \eqref{e:phi_local_coordinate_infinity} follows again from
\eqref{e:infinity_coordinate_relation}.

We next prove \eqref{e:phi_local_coordinate_infinity_disjoint}. Since \(z\)
lies in a fixed finite chart disjoint from the chart centered at
\(\infty_i\),
\[
|\phi'(z)|\asymp1.
\]
If \(i=1\), then
\[
|\phi(w)-\phi(z)|
\asymp
|\phi(w)|
\asymp
|w|
\asymp
|\wt w|^{-1},
\qquad
|\phi'(w)|\asymp1.
\]
Therefore,
\[
\frac{
\sqrt{|\phi'(w)|}\,|\rd w|
}{
|w|\,|\phi(w)-\phi(z)|
}
\asymp
|\rd\wt w|.
\]
If \(2\leq i\leq d\), then the disjointness of the charts gives
\[
|\phi(w)-\phi(z)|
\asymp1,
\]
while
\[
\sqrt{|\phi'(w)|}
\asymp
|\wt w|.
\]
Using \eqref{e:infinity_coordinate_relation}, we again obtain
\[
\frac{
\sqrt{|\phi'(w)|}\,|\rd w|
}{
|w|\,|\phi(w)-\phi(z)|
}
\asymp
|\rd\wt w|.
\]
Multiplying by
\[
\sqrt{|\phi'(z)|}\,|\rd z|
\asymp
|\rd z|
\]
proves \eqref{e:phi_local_coordinate_infinity_disjoint}.

It remains to prove \eqref{e:phi_local_coordinate_two_infinities}. If
\(p_i,p_j\in\bR\), then \(p_i\neq p_j\), and hence
\[
|\phi(w)-\phi(z)|
\asymp1.
\]
Using
\[
\sqrt{|\phi'(w)|}
\asymp
|\wt w|,
\qquad
\sqrt{|\phi'(z)|}
\asymp
|\wt z|,
\]
together with \eqref{e:infinity_coordinate_relation}, gives
\eqref{e:phi_local_coordinate_two_infinities}. If one of \(p_i,p_j\) is
\(\infty\), then
\[
|\phi(w)-\phi(z)|
\asymp
|\wt w|^{-1}
\]
for the corresponding local coordinate \(\wt w\), and the same
cancellation gives \eqref{e:phi_local_coordinate_two_infinities}.
\end{proof}

We next estimate the singular interaction between the two contours.

\begin{lemma}\label{l:weighted_contour_integral}
Suppose that the steepest-descent contour
\(\sfD^{\rm d}(w_c)\) and the steepest-ascent contour
\(\sfD^{\rm a}(z_c)\) are contained either in disjoint charts or in a
common concentric chart. In the latter case, assume that the two contours
do not share a nontrivial arc. Then
\begin{align}\label{e:weighted_contour_integral}
&\int_{\sfD^{\rm a}(z_c)}
\int_{\sfD^{\rm d}(w_c)}
\frac{
n\sqrt{|\phi'(w)|\,|\phi'(z)|}\,
|\rd w|\,|\rd z|
}{
\sqrt{(1+n|x-w|)(1+n|w-(x-s)|)}
\sqrt{(1+n|y-z|)(1+n|z-(y-t)|)}
\,|\phi(w)-\phi(z)|
}
\notag\\
&\qquad\lesssim
\frac{(\ln n)^2}{n}.
\end{align}
If the two contours are contained in disjoint charts, then the stronger
bound \(\OO(n^{-1})\) holds.
\end{lemma}

\begin{proof}
Suppose first that the two contours are contained in disjoint charts.

If both charts are centered at finite points, then
\[
\frac{\sqrt{|\phi'(w)|\,|\phi'(z)|}}
{|\phi(w)-\phi(z)|}
\lesssim1.
\]
Applying \eqref{e:one_contour_finite} to the two contours gives
\[
\int_{\sfD^{\rm a}(z_c)}
\int_{\sfD^{\rm d}(w_c)}
\frac{
n\sqrt{|\phi'(w)|\,|\phi'(z)|}\,
|\rd w|\,|\rd z|
}{
\sqrt{(1+n|x-w|)(1+n|w-(x-s)|)}
\sqrt{(1+n|y-z|)(1+n|z-(y-t)|)}
\,|\phi(w)-\phi(z)|
}
\lesssim
\frac1n.
\]

Suppose next that the descent contour lies in a chart centered at
\(\infty_i\), while the ascent contour lies in a finite disjoint chart.
In the chart at infinity,
\[
\sqrt{(1+n|x-w|)(1+n|w-(x-s)|)}
\asymp
n|w|.
\]
Using \eqref{e:phi_local_coordinate_infinity_disjoint}, the double integral
is bounded by
\[
\left(
\int_{\wt\sfD^{\rm d}(\wt w_c)}
|\rd\wt w|
\right)
\left(
\int_{\sfD^{\rm a}(z_c)}
\frac{|\rd z|}
{\sqrt{(1+n|y-z|)(1+n|z-(y-t)|)}}
\right)
\lesssim
\frac1n,
\]
where we used \eqref{e:one_contour_finite} in the second factor. The case in
which the ascent contour is centered at infinity is identical.

If both contours lie in two disjoint charts centered at points at infinity,
then
\[
\sqrt{(1+n|x-w|)(1+n|w-(x-s)|)}
\asymp
n|w|
\]
and
\[
\sqrt{(1+n|y-z|)(1+n|z-(y-t)|)}
\asymp
n|z|.
\]
By \eqref{e:phi_local_coordinate_two_infinities}, the double integral is
bounded by
\[
\frac1n
\left(
\int_{\wt\sfD^{\rm d}(\wt w_c)}
|\rd\wt w|
\right)
\left(
\int_{\wt\sfD^{\rm a}(\wt z_c)}
|\rd\wt z|
\right)
\lesssim
\frac1n.
\]
This proves the stronger estimate in the disjoint case.

Suppose now that the two contours are contained in a common concentric
chart centered at a finite point. Since \(s,t>0\) are fixed, we split
\(\sfD^{\rm d}(w_c)\) according to whether \(w\) is closer to \(x\) or
to \(x-s\), and split \(\sfD^{\rm a}(z_c)\) according to whether \(z\)
is closer to \(y\) or to \(y-t\). Thus, for each
\[
a\in\{x,x-s\},
\qquad
b\in\{y,y-t\},
\]
it suffices to estimate
\begin{align}\label{e:generic_weighted_contour}
\int_{\sfD^{\rm a}(z_c)}
\int_{\sfD^{\rm d}(w_c)}
\frac{|\rd w|\,|\rd z|}
{\sqrt{1+n|w-a|}\sqrt{1+n|z-b|}
\,|\phi(w)-\phi(z)|}.
\end{align}
Indeed, on the portion of the descent contour where \(w\) is closer to
\(a\), the other endpoint factor is bounded below by a constant times
\(\sqrt n\), and the same holds for the ascent contour.

Since \(\phi\) is a uniformizing coordinate, it is bi-Lipschitz on the
fixed finite chart. Moreover, \Cref{l:length_bound} gives the same linear
length bound for the images of the two contours under \(\phi\). Set
\[
J_{\rm d}(z)
:=
\int_{\sfD^{\rm d}(w_c)}
\frac{|\rd w|}{|\phi(w)-\phi(z)|},
\qquad
J_{\rm a}(w)
:=
\int_{\sfD^{\rm a}(z_c)}
\frac{|\rd z|}{|\phi(w)-\phi(z)|}.
\]
Using
\[
\frac{1}{\sqrt{AB}}
\leq
\frac{1}{2A}+\frac{1}{2B},
\qquad
A,B>0,
\]
we obtain
\begin{align}\label{e:split_square_root_weights}
&\int_{\sfD^{\rm a}(z_c)}
\int_{\sfD^{\rm d}(w_c)}
\frac{|\rd w|\,|\rd z|}
{\sqrt{1+n|w-a|}\sqrt{1+n|z-b|}
\,|\phi(w)-\phi(z)|}
\notag\\
&\quad\leq
\frac12
\int_{\sfD^{\rm a}(z_c)}
\frac{J_{\rm d}(z)}{1+n|z-b|}\,|\rd z|
+
\frac12
\int_{\sfD^{\rm d}(w_c)}
\frac{J_{\rm a}(w)}{1+n|w-a|}\,|\rd w|.
\end{align}
We estimate the first term; the second is identical.

For almost every \(z\in\sfD^{\rm a}(z_c)\), let
\[
\delta
:=
\dist\left(
\phi(z),\phi(\sfD^{\rm d}(w_c))
\right)>0.
\]
A dyadic decomposition of \(\phi(\sfD^{\rm d}(w_c))\) into the regions
\[
2^k\delta
\leq
|\phi(w)-\phi(z)|
<
2^{k+1}\delta
\]
and the linear length bound give
\begin{align}\label{e:log_potential_pointwise}
J_{\rm d}(z)
\lesssim
1+
\ln\frac{C}{
\dist\left(
\phi(z),\phi(\sfD^{\rm d}(w_c))
\right)}.
\end{align}

Because the two contours are real analytic and do not share a nontrivial
arc, they have only finitely many intersections, each of finite order.
Since the family of contours is fixed, the number and orders of these
intersections are uniformly bounded. Near each intersection, the logarithm
on the right-hand side of \eqref{e:log_potential_pointwise} is bounded by
a constant times \(1+|\ln r|\), where \(r\) is the arclength distance to
the intersection. Away from the intersections it is uniformly bounded.
Consequently, for every measurable set
\(E\subset\sfD^{\rm a}(z_c)\) of positive arclength,
\begin{align}\label{e:log_potential_set}
\int_E J_{\rm d}(z)\,|\rd z|
\lesssim
|E|
\left(
1+\ln\frac{C}{|E|}
\right).
\end{align}

Decompose \(\sfD^{\rm a}(z_c)\) into
\[
E_0
:=
\left\{
z\in\sfD^{\rm a}(z_c):|z-b|\leq n^{-1}
\right\}
\]
and, for \(j\geq1\),
\[
E_j
:=
\left\{
z\in\sfD^{\rm a}(z_c):
\frac{2^{j-1}}{n}<|z-b|\leq\frac{2^j}{n}
\right\}.
\]
Only \(\OO(\ln n)\) of these sets are nonempty, and
\[
|E_j|\lesssim\frac{2^j}{n},
\qquad
\frac{1}{1+n|z-b|}
\lesssim2^{-j},
\qquad
z\in E_j.
\]
Hence, by \eqref{e:log_potential_set},
\begin{align}
\int_{\sfD^{\rm a}(z_c)}
\frac{J_{\rm d}(z)}{1+n|z-b|}\,|\rd z|
&\lesssim
\sum_{E_j\neq\emptyset}
2^{-j}|E_j|
\left(
1+\ln\frac{C}{|E_j|}
\right)
\notag\\
&\lesssim
\frac{1}{n}
\sum_{E_j\neq\emptyset}
\left(
1+\ln\frac{Cn}{2^j}
\right)
\lesssim
\frac{(\ln n)^2}{n}.
\end{align}
Together with \eqref{e:split_square_root_weights}, this proves
\eqref{e:weighted_contour_integral} in a common finite chart.

Suppose finally that the common concentric chart is centered at
\(\infty_i\). Then
\[
\sqrt{(1+n|x-w|)(1+n|w-(x-s)|)}
\asymp n|w|,
\]
and similarly,
\[
\sqrt{(1+n|y-z|)(1+n|z-(y-t)|)}
\asymp n|z|.
\]
Therefore, the left-hand side of
\eqref{e:weighted_contour_integral} is bounded by
\begin{align}
\frac1n
\int_{\sfD^{\rm a}(z_c)}
\int_{\sfD^{\rm d}(w_c)}
\frac{
\sqrt{|\phi'(w)|\,|\phi'(z)|}\,
|\rd w|\,|\rd z|
}{
|w|\,|z|\,|\phi(w)-\phi(z)|
}.
\end{align}
By \eqref{e:phi_local_coordinate_infinity}, this is comparable to
\[
\frac1n
\int_{\wt\sfD^{\rm a}(\wt z_c)}
\int_{\wt\sfD^{\rm d}(\wt w_c)}
\frac{|\rd\wt w|\,|\rd\wt z|}
{|\wt w-\wt z|}.
\]
The transformed contours are real analytic and do not share a nontrivial
arc. The same finite-order intersection argument as above shows that the
last double integral is uniformly bounded. Hence the contribution is
\(\OO(n^{-1})\), which completes the proof.
\end{proof}

The following lemma gives an upper bound for the double-contour integral.

\begin{lemma}\label{l:double_bound}
Let \(\sfD^{\rm d}(w_c)\) be a steepest-descent contour and let
\(\sfD^{\rm a}(z_c)\) be a steepest-ascent contour. Suppose that the two
contours are contained either in disjoint charts or in a common concentric
chart and, in the latter case, do not share a nontrivial arc. Then
\begin{align}\begin{split}\label{e:PQint}
&\phantom{{}={}}\left|
\frac{n}{(2\pi\ri)^2}
\int_{\sfD^{\rm a}(z_c)}
\int_{\sfD^{\rm d}(w_c)}
P_{ns}(nw,nx)\,Q_{nt}(nz,ny)\,
\frac{I_+(w)}{I_-(z)}\,
\frac{\sqrt{\phi'(w)}\sqrt{\phi'(z)}}
{\phi(w)-\phi(z)}\,
\rd w\,\rd z
\right|
\\
&\lesssim
\int_{\sfD^{\rm a}(z_c)}
\int_{\sfD^{\rm d}(w_c)}
e^{n\Re[S(w;x,s)-S(z;y,t)]}
\frac{\sqrt{|\phi'(w)|\,|\phi'(z)|}}
{|\phi(w)-\phi(z)|}
\\
&\times
\frac{n^2\,|\rd w|\,|\rd z|}
{
\sqrt{(1+n|x-w|)(1+n|w-(x-s)|)}
\sqrt{(1+n|y-z|)(1+n|z-(y-t)|)}
}.
\end{split}\end{align}
Moreover,
\begin{align}\label{e:PQint2}
&\left|
\frac{n}{(2\pi\ri)^2}
\int_{\sfD^{\rm a}(z_c)}
\int_{\sfD^{\rm d}(w_c)}
P_{ns}(nw,nx)\,Q_{nt}(nz,ny)\,
\frac{I_+(w)}{I_-(z)}\,
\frac{\sqrt{\phi'(w)}\sqrt{\phi'(z)}}
{\phi(w)-\phi(z)}\,
\rd w\,\rd z
\right|
\notag\\
&\qquad\lesssim
(\ln n)^2
e^{n\Re[S(w_c;x,s)-S(z_c;y,t)]}.
\end{align}

If the two contours are contained in disjoint charts, or if they are
contained in the same chart and \(w_c\) and \(z_c\) are uniformly
separated in the local coordinate of that chart, then
\begin{align}\label{e:PQint3}
&\left|
\frac{n}{(2\pi\ri)^2}
\int_{\sfD^{\rm a}(z_c)}
\int_{\sfD^{\rm d}(w_c)}
P_{ns}(nw,nx)\,Q_{nt}(nz,ny)\,
\frac{I_+(w)}{I_-(z)}\,
\frac{\sqrt{\phi'(w)}\sqrt{\phi'(z)}}
{\phi(w)-\phi(z)}\,
\rd w\,\rd z
\right|
\notag\\
&\qquad\lesssim
(\ln n)^{10}\Delta(w_c)\Delta(z_c)
e^{n\Re[S(w_c;x,s)-S(z_c;y,t)]}.
\end{align}
\end{lemma}

\begin{proof}[Proof of \Cref{l:double_bound}]
By
\Cref{l:PIQI_bound,l:vertical_tangent_S,l:unit_slope_tangent_S,l:horizontal_tangent_S}, for
\[
w\in\sfD^{\rm d}(w_c),
\qquad
z\in\sfD^{\rm a}(z_c),
\]
we have
\begin{align}\label{e:PQII}
\left|P_{ns}(nw,nx)I_+(w)\right|
&\lesssim
\frac{\sqrt n}
{\sqrt{(1+n|x-w|)(1+n|w-(x-s)|)}}
e^{n\Re S(w;x,s)},
\notag\\
\left|Q_{nt}(nz,ny)I_-^{-1}(z)\right|
&\lesssim
\frac{\sqrt n}
{\sqrt{(1+n|y-z|)(1+n|z-(y-t)|)}}
e^{-n\Re S(z;y,t)}.
\end{align}
Substituting \eqref{e:PQII} into the double-contour integral gives
\eqref{e:PQint}.

Since \(\sfD^{\rm d}(w_c)\) is a steepest-descent contour for
\(\Re S(\,\cdot\,;x,s)\), while \(\sfD^{\rm a}(z_c)\) is a
steepest-ascent contour for \(\Re S(\,\cdot\,;y,t)\), one has
\[
\Re S(w;x,s)
\leq
\Re S(w_c;x,s),
\qquad
\Re S(z;y,t)
\geq
\Re S(z_c;y,t).
\]
Therefore,
\begin{align}\label{e:PQ_action_bound}
\Re[S(w;x,s)-S(z;y,t)]
\leq
\Re[S(w_c;x,s)-S(z_c;y,t)].
\end{align}
The double integral on the right-hand side of \eqref{e:PQint}, after the
exponential factor is removed using \eqref{e:PQ_action_bound}, is \(n\)
times the left-hand side of \eqref{e:weighted_contour_integral}. Therefore,
\eqref{e:PQint2} follows from \Cref{l:weighted_contour_integral}.

We now prove \eqref{e:PQint3}. Whenever a truncated contour is defined,
replace \(\sfD^{\rm d}(w_c)\) by \(\sfS^{\rm d}(w_c)\), and similarly
replace \(\sfD^{\rm a}(z_c)\) by \(\sfS^{\rm a}(z_c)\). If no truncation
is made, we retain the full contour; in that case, the corresponding
localization length equals \(1\).

By
\Cref{c:bulk_steepest,c:arctic_steepest,l:local_descent_deformation},
for every
\[
w\in
\sfD^{\rm d}(w_c)\setminus\sfS^{\rm d}(w_c),
\]
one has
\[
e^{n\Re S(w;x,s)}
\leq
e^{n\Re S(w_c;x,s)}
e^{-\fc'(\ln n)^2}.
\]
The analogous estimate holds on the discarded portion of the ascent
contour:
\[
e^{-n\Re S(z;y,t)}
\leq
e^{-n\Re S(z_c;y,t)}
e^{-\fc'(\ln n)^2}.
\]
Hence, by the proof of \eqref{e:PQint2}, the contribution of all discarded
portions is bounded by
\[
(\ln n)^2
e^{-\fc'(\ln n)^2}
e^{n\Re[S(w_c;x,s)-S(z_c;y,t)]},
\]
and  this contribution is absorbed into the
right-hand side of \eqref{e:PQint3}.

It remains to estimate the retained portions. Suppose first that both
charts are centered at finite points. If the charts are disjoint, or if
\(w_c\) and \(z_c\) are uniformly separated in a common chart, then, after
choosing the truncations sufficiently small,
\[
\frac{\sqrt{|\phi'(w)|\,|\phi'(z)|}}
{|\phi(w)-\phi(z)|}
\lesssim1
\]
on the retained contours. Therefore,
\begin{align}\begin{split}
&\phantom{{}={}}\int_{\sfS^{\rm a}(z_c)}
\int_{\sfS^{\rm d}(w_c)}
\frac{
n^2\,|\rd w|\,|\rd z|
}{
\sqrt{(1+n|x-w|)(1+n|w-(x-s)|)}\sqrt{(1+n|y-z|)(1+n|z-(y-t)|)}}
\\
&=
n
\left(
\int_{\sfS^{\rm d}(w_c)}
\frac{\sqrt n\,|\rd w|}
{\sqrt{(1+n|x-w|)(1+n|w-(x-s)|)}}
\right)
\left(
\int_{\sfS^{\rm a}(z_c)}
\frac{\sqrt n\,|\rd z|}
{\sqrt{(1+n|y-z|)(1+n|z-(y-t)|)}}
\right)
\\
&\quad\lesssim
(\ln n)^{10}\Delta(w_c)\Delta(z_c),
\end{split}\end{align}
where we used \Cref{l:contour_integral}. If one of the contours is not
truncated, the corresponding estimate follows from
\eqref{e:one_contour_finite} and the fact that its localization length is
\(1\).

Suppose next that both contours lie in a common chart centered at
\(\infty_i\), and that \(\wt w_c\) and \(\wt z_c\) are uniformly
separated. Using
\eqref{e:phi_local_coordinate_infinity} and the endpoint estimates at
infinity, the retained double integral is comparable to
\[
\int_{\wt\sfS^{\rm a}(\wt z_c)}
\int_{\wt\sfS^{\rm d}(\wt w_c)}
\frac{|\rd\wt w|\,|\rd\wt z|}
{|\wt w-\wt z|}.
\]
The denominator is uniformly bounded away from zero on sufficiently small
truncated contours. By \eqref{e:one_contour_infinity2}, the lengths of the
two transformed contours are bounded by
\[
(\ln n)^5\Delta(w_c)
\qquad\text{and}\qquad
(\ln n)^5\Delta(z_c),
\]
respectively. Hence the last double integral is bounded by
\[
(\ln n)^{10}\Delta(w_c)\Delta(z_c).
\]

Suppose that the descent contour lies in a chart centered at
\(\infty_i\), while the ascent contour lies in a disjoint finite chart.
Using \eqref{e:phi_local_coordinate_infinity_disjoint}, the retained
double integral is bounded by
\begin{align}
&n
\left(
\int_{\wt\sfS^{\rm d}(\wt w_c)}
|\rd\wt w|
\right)
\left(
\int_{\sfS^{\rm a}(z_c)}
\frac{|\rd z|}
{\sqrt{(1+n|y-z|)(1+n|z-(y-t)|)}}
\right)
\notag\\
&\qquad\lesssim
(\ln n)^{10}\Delta(w_c)\Delta(z_c).
\end{align}
The case in which the ascent contour is centered at infinity is identical.

Finally, if the two contours lie in two disjoint charts centered at points
at infinity, then \eqref{e:phi_local_coordinate_two_infinities} reduces
the retained double integral to
\[
\int_{\wt\sfS^{\rm a}(\wt z_c)}
\int_{\wt\sfS^{\rm d}(\wt w_c)}
|\rd\wt w|\,|\rd\wt z|,
\]
which is again bounded by
\[
(\ln n)^{10}\Delta(w_c)\Delta(z_c).
\]
Combining the retained and discarded estimates with
\eqref{e:PQ_action_bound} proves \eqref{e:PQint3}.
\end{proof}

\section{Standard Form in the Liquid Case}\label{s:liquid_standard_form}
Let \((y,t)\in\fP\cap\bZ^2/n\) represent a white triangle, and let
\((x,s)\in\fP\cap\bZ^2/n\) represent a blue triangle.

In this section, we prove \Cref{p:standard_form,p:standard_form2} when
\(\fN_{(y,t)}\) is a liquid or ramification neighborhood.

\begin{proposition}
Assume that \(\fN_{(y,t)}\) is a liquid or ramification neighborhood.
Then \((y,t)\) is associated with a pair of complex-conjugate liquid or
ramification charts, respectively. Suppose that \((x,s)\in\fN_\al\).
There are two cases:
\begin{enumerate}
\item
The charts associated with \((y,t)\in\fN_{(y,t)}\) are disjoint from all
charts associated with \(\fN_\alpha\).

\item
The two liquid or ramification charts associated with
\((y,t)\in\fN_{(y,t)}\) have the same centers as the corresponding liquid
or ramification charts associated with \(\fN_\alpha\).
\end{enumerate}
In either case, the conclusions of
\Cref{p:standard_form,p:standard_form2} hold.
\end{proposition}

\subsection{Deforming double-contour integral}

As in
\Cref{s:critical_point}, we associate with \((x,s)\) a collection of charts
carrying local descent contours \(\sfC^{\rm d}(w_0)\). We deform each such
contour to  steepest-descent paths \(\sfD^{\rm d}(w_c)\), where \(w_c\) are
 descent critical points associated with \((x,s)\).
Similarly, the charts associated with \((y,t)\) carry local ascent contours
\(\sfC^{\rm a}(z_0)\), which we deform to steepest-ascent paths
\(\sfD^{\rm a}(z_c)\), where \(z_c\) are ascent critical points associated
with \((y,t)\).

During these deformations, one must account for the residue at \(w=z\).
The following proposition identifies the combination of the double integral
and the corresponding residue contribution that remains invariant under the
contour deformations.

\begin{proposition}[Invariance]\label{p:invariance_contour}
Let $U\subset \bC$ be a simply connected domain, and let
$\phi$ be biholomorphic on $U$. Fix a branch of \(\sqrt{\phi'}\) on \(U\). Let \(f\) be a holomorphic
function on \(U\times U\), and let \(\sfC_w\) and \(\sfC_z\) be finite
unions of oriented, piecewise smooth simple paths in \(U\) with fixed
endpoints. Assume that \(\sfC_w\) and \(\sfC_z\) have only finitely many
intersection points, all of which are transverse and away from the
endpoints.
For each $\xi\in \sfC_w\cap\sfC_z$, let
$\sgn(\xi)\in\{\pm1\}$ denote the local intersection number of $\sfC_w$ and
$\sfC_z$ at $\xi$ (recall from \Cref{d:positive_negative}).

Fix $\xi_0\in U\cup\{\infty\}$ and consider the quantity
\begin{align}\label{e:invariant}
\frac{n}{(2\pi \ri)^2}\int_{\sfC_w}\rd w\int_{\sfC_z}\rd z\,f(w,z)\frac{\sqrt{\phi'(w)}\sqrt{\phi'(z)}}{\phi(w)-\phi(z)}
+\sum_{\xi\in \sfC_w\cap\sfC_z}\sgn(\xi)\cdot \frac{n}{2\pi \ri}\int_{\xi_0}^{\xi} f(\zeta,\zeta)\,\rd\zeta.
\end{align}
Then \eqref{e:invariant} is invariant under any deformation of
$(\sfC_w,\sfC_z)$ that keeps the endpoints fixed, keeps intersections away
from the endpoints, and preserves transversality except at finitely many
ordinary tangencies.
\end{proposition}

\begin{proof}[Proof of \Cref{p:invariance_contour}]
It suffices to check the claim for local deformations. Away from the diagonal
$w=z$, the kernel
\[
f(w,z)\frac{\sqrt{\phi'(w)}\sqrt{\phi'(z)}}{\phi(w)-\phi(z)}\,\rd w\,\rd z
\]
is holomorphic. Hence the double integral can change only when the local
deformation crosses the diagonal, that is, only through the motion, creation,
or annihilation of intersection points of $\sfC_w$ and $\sfC_z$.

Near the diagonal, since $\phi$ is biholomorphic,
\[
\phi(w)-\phi(z)=\phi'(z)(w-z)+\OO((w-z)^2),
\qquad
\sqrt{\phi'(w)}\sqrt{\phi'(z)}=\phi'(z)+\OO(w-z).
\]
Therefore
\[
\Res_{w=z}
\left[
f(w,z)\frac{\sqrt{\phi'(w)}\sqrt{\phi'(z)}}{\phi(w)-\phi(z)}
\right]
=f(z,z).
\]
Thus the only possible change in the double integral is the residue
contribution along the portion of the diagonal swept out by the local move.

We now consider the three possible local moves. Throughout, $\Delta$ denotes
the new value minus the old value.

\begin{enumerate}
\item \emph{A single intersection point moves.}
Suppose an intersection point moves continuously from $\xi$ to $\xi'$, with
the same local intersection number:
$
\sgn(\xi)=\sgn(\xi')$.
Then the residue computation gives
\[
\Delta\!\left(
\frac{n}{(2\pi \ri)^2}
\int_{\sfC_w}\rd w\int_{\sfC_z}\rd z\,
f(w,z)\frac{\sqrt{\phi'(w)}\sqrt{\phi'(z)}}{\phi(w)-\phi(z)}
\right)
=
-\sgn(\xi)\cdot \frac{n}{2\pi \ri}
\int_{\xi}^{\xi'} f(\zeta,\zeta)\,\rd\zeta .
\]
On the other hand, the correction term changes by
\[
\sgn(\xi)\cdot \frac{n}{2\pi \ri}
\int_{\xi_0}^{\xi'} f(\zeta,\zeta)\,\rd\zeta
-
\sgn(\xi)\cdot \frac{n}{2\pi \ri}
\int_{\xi_0}^{\xi} f(\zeta,\zeta)\,\rd\zeta
=
\sgn(\xi)\cdot \frac{n}{2\pi \ri}
\int_{\xi}^{\xi'} f(\zeta,\zeta)\,\rd\zeta .
\]
The two changes cancel.

\item \emph{A pair of intersections is created.}
Suppose a deformation creates two new intersections $\xi,\xi'$. Then their
local intersection numbers are opposite:
$
\sgn(\xi')=-\sgn(\xi)$.
The residue contribution to the double integral is
\[
\Delta\!\left(
\frac{n}{(2\pi \ri)^2}
\int_{\sfC_w}\rd w\int_{\sfC_z}\rd z\,
f(w,z)\frac{\sqrt{\phi'(w)}\sqrt{\phi'(z)}}{\phi(w)-\phi(z)}
\right)
=
\sgn(\xi)\cdot \frac{n}{2\pi \ri}
\int_{\xi}^{\xi'} f(\zeta,\zeta)\,\rd\zeta .
\]
The correction term changes from $0$ to
\[
\sgn(\xi)\cdot \frac{n}{2\pi \ri}
\int_{\xi_0}^{\xi} f(\zeta,\zeta)\,\rd\zeta
-\sgn(\xi)\cdot \frac{n}{2\pi \ri}
\int_{\xi_0}^{\xi'} f(\zeta,\zeta)\,\rd\zeta
=
-\sgn(\xi)\cdot \frac{n}{2\pi \ri}
\int_{\xi}^{\xi'} f(\zeta,\zeta)\,\rd\zeta .
\]
Again the two changes cancel.

\item \emph{A pair of intersections is annihilated.}
This is the reverse of the preceding move, so the same cancellation applies.
\end{enumerate}

Since any endpoint-preserving generic deformation can be decomposed into these
local moves, and no other change is possible, the quantity
\eqref{e:invariant} is invariant.
\end{proof}

\subsection{Concentric case}
We start with the case that both $\fN_\al$ and $\fN_{(y,t)}$ are liquid neighborhood, and share concentric liquid charts. For the pair of complex-conjugate liquid charts centered at $w_0\in \bC_+$ and
$\overline{w_0}\in \bC_-$, the kernel ansatz \eqref{e:def_Aalpha} is given by
\begin{align}
A_\al((x,s),(y,t))&=T_1+T_2,
\end{align}
where
\begin{align}\begin{split}\label{e:liquid-liquid}
T_1&:=\frac{n}{2\pi\ri}\int_{\sfC(w_0; (x,s),(y,t))}
P_{ns}(nz,nx)\,Q_{nt}(nz,ny)\,\rd z\\
&+\frac{n}{(2\pi \ri)^2}
\left(\int_{\sfC^{\rm a}(w_0)}\!\!\int_{\sfC^{\rm d}(w_0)}
+\int_{\sfC^{\rm a}(\overline{w_0})}\!\!\int_{\sfC^{\rm d}(\overline{w_0})}
\right)
P_{ns}(nw,nx)\,Q_{nt}(nz,ny)\,
\frac{I(w)}{I(z)}\,
\frac{\sqrt{\phi'(w)}\sqrt{\phi'(z)}}{\phi(w)-\phi(z)}\,
\rd w\,\rd z,\\
T_2&:=\frac{n}{(2\pi \ri)^2}
\left(\int_{\sfC^{\rm a}(w_0)}\!\!\int_{\sfC^{\rm d}(\overline{w_0})}
+\int_{\sfC^{\rm a}(\overline{w_0})}\!\!\int_{\sfC^{\rm d}(w_0)}\right)
P_{ns}(nw,nx)\,Q_{nt}(nz,ny)\,
\frac{I(w)}{I(z)}\,
\frac{\sqrt{\phi'(w)}\sqrt{\phi'(z)}}{\phi(w)-\phi(z)}\,
\rd w\,\rd z.
\end{split}\end{align}

By \Cref{c:bulk}, \((x,s)\) is associated with a pair of
complex-conjugate critical points \(w_c,\overline{w_c}\) and is adapted to
the corresponding liquid charts, while \((y,t)\) is associated with a pair
of complex-conjugate critical points \(z_c,\overline{z_c}\) and is adapted
to the corresponding liquid charts.

We start with $T_2$. Notice that the liquid charts centered at $w_0$ and
$\overline{w_0}$ are disjoint, so we can deform
$\sfC^{\rm d}(w_0)$ and $\sfC^{\rm d}(\overline{w_0})$ to the steepest
descent paths $\sfD^{\rm d}(w_c)$ and $\sfD^{\rm d}(\overline{w_c})$,
respectively, and deform $\sfC^{\rm a}(w_0)$ and
$\sfC^{\rm a}(\overline{w_0})$ to the steepest ascent paths
$\sfD^{\rm a}(z_c)$ and $\sfD^{\rm a}(\overline{z_c})$, respectively, as
introduced in \Cref{c:bulk_steepest}. Then
\begin{align}\begin{split}\label{e:liquid-liquid1}
T_2
&=\frac{n}{(2\pi \ri)^2}
\left(\int_{\sfD^{\rm a}(z_c)}\!\!\int_{\sfD^{\rm d}(\overline{w_c})}
+\int_{\sfD^{\rm a}(\overline{z_c})}\!\!\int_{\sfD^{\rm d}(w_c)}\right)
P_{ns}(nw,nx)\,Q_{nt}(nz,ny)\,
\frac{I(w)}{I(z)}\,
\frac{\sqrt{\phi'(w)}\sqrt{\phi'(z)}}{\phi(w)-\phi(z)}\,
\rd w\,\rd z\\
&\quad
+\OO\left(e^{-\fc' n}e^{n\Re[S(w_c;x,s)-S(z_c;y,t)]}\right),
\end{split}
\end{align}
where the error comes from the integrals over the sub-arcs, using
\eqref{e:bulk_subarc}.

Next we study $T_1$. We recall from \Cref{d:defCw} that there are two cases
for $\sfC(w_0; (x,s),(y,t))$:
\begin{enumerate}
\item If \(s\ge t\), let $\xi_0=\infty$. Then \(\sfC(w_0; (x,s),(y,t))\) consists of a
path from \(w_0\) to \(\xi_0\) in the upper half-plane, and one path from
\(\xi_0\) to \(\overline{w_0}\) in the lower half-plane.

\item If \(s<t\), let $\xi_0$ be a point in
$(\max\{x-s,y-t\},\,\min\{x,y\})$. Then \(\sfC(w_0; (x,s),(y,t))\) consists of
a path from \(w_0\) to \(\xi_0\) in the upper half-plane, and one path from
\(\xi_0\) to \(\overline{w_0}\) in the lower half-plane.
\end{enumerate}
Moreover the two paths
$\sfC^{\rm a}(w_0)$ and $\sfC^{\rm d}(w_0)$ intersect negatively, and the two paths
$\sfC^{\rm a}(\overline{w_0})$ and $\sfC^{\rm d}(\overline{w_0})$ intersect positively.

Therefore $T_1$ is in the form of \eqref{e:invariant}. More precisely, the
double-contour integral over $\sfC^{\rm a}(w_0)$ and $\sfC^{\rm d}(w_0)$, where the
two paths $\sfC^{\rm a}(w_0)$ and $\sfC^{\rm d}(w_0)$ intersect negatively,
together with the single-contour integral from $w_0$ to $\xi_0$, is in the form of
\eqref{e:invariant}; the double-contour integral over
$\sfC^{\rm a}(\overline{w_0})$ and $\sfC^{\rm d}(\overline{w_0})$, where the
two paths $\sfC^{\rm a}(\overline{w_0})$ and
$\sfC^{\rm d}(\overline{w_0})$ intersect positively, together with the
single-contour integral from $\xi_0$ to $\overline{w_0}$, is also in the form of
\eqref{e:invariant}.

Hence, by \Cref{p:invariance_contour}, similarly to the disjoint case, we
can deform $\sfC^{\rm d}(w_0)$ and $\sfC^{\rm d}(\overline{w_0})$ to the
steepest descent paths $\sfD^{\rm d}(w_c)$ and
$\sfD^{\rm d}(\overline{w_c})$, respectively, and deform
$\sfC^{\rm a}(w_0)$ and $\sfC^{\rm a}(\overline{w_0})$ to the steepest
ascent paths $\sfD^{\rm a}(z_c)$ and $\sfD^{\rm a}(\overline{z_c})$,
respectively, as introduced in \Cref{c:bulk_steepest}, and compensate by the
residual:
\begin{align}\begin{split}\label{e:liquid-liquid2}
T_1
&=\frac{n}{(2\pi \ri)^2}
\left(\int_{\sfD^{\rm a}(z_c)}\!\!\int_{\sfD^{\rm d}(w_c)}
+\int_{\sfD^{\rm a}(\overline{z_c})}\!\!\int_{\sfD^{\rm d}(\overline{w_c})}\right)
P_{ns}(nw,nx)\,Q_{nt}(nz,ny)\,
\frac{I(w)}{I(z)}\,
\frac{\sqrt{\phi'(w)}\sqrt{\phi'(z)}}{\phi(w)-\phi(z)}\,
\rd w\,\rd z\\
&\quad
-\sum_{\xi}\sgn(\xi)\cdot \frac{n}{2\pi \ri}
\int_{\sfC(\xi;(x,s),(y,t))}
P_{ns}(nz,nx)\,Q_{nt}(nz,ny)\,\rd z
+\OO\left(e^{-\fc' n}e^{n\Re[S(w_c;x,s)-S(z_c;y,t)]}\right),
\end{split}
\end{align}
where the sum over $\xi$ is over all intersections of
$\sfD^{\rm d}(w_c)\cap \sfD^{\rm a}(z_c)$, and we used the fact that the
integral paths from $\xi$ to $\xi_0$ and from $\xi_0$ to $\overline{\xi}$
together give $\sfC(\xi;(x,s),(y,t))$. Finally, using
\Cref{l:deform_single_contour}, we can further deform the contours
$\sfC(\xi;(x,s),(y,t))$ to the steepest descent contours
$\sfD(\xi;(x,s),(y,t))$ over all intersections of
$\sfD^{\rm d}(w_c)\cap \sfD^{\rm a}(z_c)$ and $\sfD^{\rm d}(\overline{w_c})\cap \sfD^{\rm a}(\overline{z_c})$. The claim of \Cref{p:standard_form} in this case follows by
combining \eqref{e:liquid-liquid1} and \eqref{e:liquid-liquid2}.

The case in which both \(\fN_\al\) and \(\fN_{(y,t)}\) are ramification
neighborhoods and share concentric ramification charts can be proved in the
same way after replacing
\begin{align}\label{e:replaceP}
P_{ns}(nw,nx), \sfD^{\rm d}(w_c), I(w), \phi(w), \phi'(w)\quad
\text{by}\quad 
P_{n(s+\ft)}(nw,nx),\sfD^{\rm d}(w_{c,\ft}), I_\ft(w),
\phi_\ft(w), \phi_\ft'(w),
\end{align}
respectively, and replacing
\begin{align}\label{e:replaceQ}
Q_{nt}(nz,ny), \sfD^{\rm a}(z_c), I(z), \phi(z), \phi'(z)
\quad
\text{by}\quad 
Q_{n(t+\ft)}(nz,ny), \sfD^{\rm a}(z_{c,\ft}), I_\ft(z),
\phi_\ft(z), \phi_\ft'(z),
\end{align}
respectively. Here \(w_{c,\ft}\) and \(z_{c,\ft}\) denote the corresponding
critical points in the ramification charts. The error terms in
\eqref{e:liquid-liquid1} and \eqref{e:liquid-liquid2} take the form
\begin{align}\label{e:error_liquid}
\OO\left(
e^{-\fc'n}
e^{n\Re\left[
S_\ft(w_{c,\ft};x,s+\ft)
-
S_\ft(z_{c,\ft};y,t+\ft)
\right]}
\right).
\end{align}
By \Cref{c:change_time}, we have
\begin{align}\label{e:replaceS}
S_\ft(w_{c,\ft};x,s+\ft)
&=
S(w_c;x,s),
\qquad
S_\ft(z_{c,\ft};y,t+\ft)
=
S(z_c;y,t).
\end{align}
Thus, \eqref{e:error_liquid} can be replaced by
\begin{align}
\OO\left(
e^{-\fc'n}
e^{n\Re\left[S(w_c;x,s)-S(z_c;y,t)\right]}
\right).
\end{align}
The claim of \Cref{p:standard_form} in this case follows.

\subsection{Disjoint case}
We first consider the case in which \(\fN_{(y,t)}\) is a liquid
neighborhood and \(\fN_\al\) is not a ramification neighborhood.

If the charts associated with \(\fN_{(y,t)}\) are disjoint from all charts
associated with \(\fN_\alpha\), then the single-contour integral satisfies
$
J^{(1)}\equiv0.
$

Each relevant pair of charts associated with \((x,s)\) and \((y,t)\),
centered at \(w_0\) and \(z_0\), respectively, gives a double-contour
integral term as in \eqref{e:all_term}:
\begin{equation}\label{e:all_term2}
\frac{n}{(2\pi\ri)^2}
\int_{\sfC^{\rm a}(z_0)}\!\!\int_{\sfC^{\rm d}(w_0)}
P_{ns}(nw,nx)\,Q_{nt}(nz,ny)\,
\frac{I_+(w)}{I(z)}\,
\frac{\sqrt{\phi'(w)}\sqrt{\phi'(z)}}{\phi(w)-\phi(z)}\,
\rd w\,\rd z.
\end{equation}
We can deform the descent contour \(\sfC^{\rm d}(w_0)\) to a possibly
disconnected union of steepest-descent paths
\(\sfD^{\rm d}(w_c)\), and the ascent contour
\(\sfC^{\rm a}(z_0)\) to the steepest-ascent path
\(\sfD^{\rm a}(z_c)\), as introduced in \Cref{s:critical_point}. Here
the \(w_c\)'s and \(z_c\) are the critical points associated with
\((x,s)\) and \((y,t)\), respectively; the liquid chart contains only one
critical point \(z_c\). We obtain
\begin{align}
\begin{split}\label{e:Aal_term_dis_liquid}
\eqref{e:all_term2}
&=
\sum_{w_c}\frac{n}{(2\pi\ri)^2}
\int_{\sfD^{\rm a}(z_c)}\!\!\int_{\sfD^{\rm d}(w_c)}
P_{ns}(nw,nx)\,Q_{nt}(nz,ny)\,
\frac{I_+(w)}{I(z)}\,
\frac{\sqrt{\phi'(w)}\sqrt{\phi'(z)}}{\phi(w)-\phi(z)}\,
\rd w\,\rd z\\
&\quad+
\sum_{w_c}
\OO\left(
e^{-\fc'n}
e^{n\Re[S(w_c;x,s)-S(z_c;y,t)]}
\right).
\end{split}
\end{align}
The claim of \Cref{p:standard_form} in this case follows by summing over
all relevant pairs of charts.

The remaining cases, in which \(\fN_\al\) is a ramification neighborhood,
\(\fN_{(y,t)}\) is a ramification neighborhood, or both are ramification
neighborhoods, follow from exactly the same argument after making the
replacements in \eqref{e:replaceP} or \eqref{e:replaceQ}, or both, as
appropriate, and using the relations \eqref{e:replaceS} for the error
terms.

\section{Standard Form in the Non-Liquid Case}\label{s:non-liquid_standard_form}

There are two possibilities. Either \(\fN_{(y,t)}\) intersects a unique
curvilinear triangle \(\fT\), in which case
$
(y,t)\in\fL$
or $
\operatorname{Cell}^{\rw}(y,t)=\fT$;
or \(\fN_{(y,t)}\) intersects two adjacent curvilinear triangles
\(\fT_A\) and \(\fT_B\). In the latter case, after relabeling if necessary,
we may assume that
$
(y,t)\in\fL$
or $
\operatorname{Cell}^{\rw}(y,t)=\fT_B$,
and we set \(\fT:=\fT_B\). Thus, in either case, we fix a curvilinear
triangle \(\fT\) such that
\[
\fN_{(y,t)}\cap\fT\neq\emptyset,
\qquad
(y,t)\in\fL
\quad\text{or}\quad
\operatorname{Cell}^{\rw}(y,t)=\fT.
\]

If \(\fN_{(y,t)}\) is not a liquid or ramification neighborhood, then, by
\Cref{c:concentric}, the chart collections associated with
\(\fN_{(y,t)}\) and \(\fN_\al\) contain at most one concentric pair.

\begin{proposition}\label{l:compare_Aalpha}
Under the assumptions above, exactly one of the following two possibilities
occurs:
\begin{enumerate}
\item
No chart associated with \(\fN_\al\) carrying a local descent path is
concentric with a chart associated with \(\fN_{(y,t)}\) carrying a local
ascent path.

\item
There exists a descent cut \(\ell_-(w_0;\fT)\) appearing in the
construction of \(\fC(\fT;\fN_{(y,t)})\) such that
\[
\fN_\al
\cap
\ell_-(w_0;\fT)\}
\neq\emptyset.
\]
In this case, a chart associated with \(\fN_\al\) carrying a local descent
path and a chart associated with \(\fN_{(y,t)}\) carrying a local ascent
path are concentric, with common center \(w_0\). Every other pair consisting
of a chart associated with \(\fN_\al\) and a chart associated with
\(\fN_{(y,t)}\) is disjoint.
\end{enumerate}
In either case, the conclusion of \Cref{p:standard_form} holds.
\end{proposition}

If \(\fN_\al\) is a ramification neighborhood, then all charts associated
with \(\fN_{(y,t)}\) are disjoint from the charts associated with
\(\fN_\al\). The same argument as in the liquid case therefore gives
\Cref{p:standard_form}. In the following, we assume that \(\fN_\al\)
is not a ramification neighborhood.

\subsection{Estimates on the single-contour integral}

We recall \(J_\fT=J_{\fT}((x,s),(y,t))\) and
\(\cI((x,s),(y,t))\) from \eqref{e:J1form} and
\Cref{d:interlacing}. The following statements give estimates on the single-contour integral $J^{(1)}$ from \eqref{e:def_Aalpha}
\begin{lemma}[Disjoint case]\label{l:JTsmall}
Let \(\fT\) be a curvilinear triangle such that
\[
\fN_{(y,t)}\cap\fT\neq\emptyset, \quad \operatorname{Cell}^{\rm w}(y,t)\neq \fT.
\]
Suppose that no chart associated with \(\fN_\al\) carrying a local descent path is
concentric with a chart associated with \(\fN_{(y,t)}\) carrying a local
ascent path. Then there exists a constant \(\fc'>0\) such that
\begin{align}\label{e:JTsmall}
\bm1\!\left(
\fN_\al\cap\fC(\fT;\fN_{(y,t)})\neq\emptyset
\right)
|J_\fT|
=
\sum_{w_c,z_c}
\OO\left(
e^{-\fc'n}
e^{n\Re[S(w_c;x,s)-S(z_c;y,t)]}
\right).
\end{align}
Here \(w_c\) ranges over 
\(\operatorname{Crit}(x,s;\fN_\al)\), and \(z_c\) ranges over \(\operatorname{Crit}(y,t;\fN_{(y,t)})\). 
\end{lemma}

\begin{lemma}[Disjoint case]\label{l:JTterm}
Let \(\fT\) be a curvilinear triangle, and assume that
\[
\operatorname{Cell}^{\rw}(y,t)=\fT.
\]
Suppose that no chart associated with \(\fN_\al\) carrying a local descent path is
concentric with a chart associated with \(\fN_{(y,t)}\) carrying a local
ascent path. Then there exists a
constant \(\fc'>0\) such that
\begin{align}\label{e:JTterm}
&\bm1\!\left(
\fN_\al\cap\fC(\fT;\fN_{(y,t)})\neq\emptyset
\right)
J_\fT
=
\cI((x,s),(y,t))J_\fT
+
\sum_{w_c,z_c}
\OO\left(
e^{-\fc'n}
e^{n\Re[S(w_c;x,s)-S(z_c;y,t)]}
\right).
\end{align}
Here \(w_c\) ranges over 
\(\operatorname{Crit}(x,s;\fN_\al)\), and \(z_c\) ranges over \(\operatorname{Crit}(y,t;\fN_{(y,t)})\). 
\end{lemma}

\begin{lemma}[Concentric case]\label{p:concentric_setting}
Suppose that a chart associated with \(\fN_\al\), carrying a local descent
path, and a chart associated with \(\fN_{(y,t)}\), carrying a local ascent
path, are concentric, with common center \(w_0\), and that \(w_0\)
corresponds to a tangency or cusp-turning location. The tangent line
\(L(w_0)\) separates two curvilinear triangles \(\fT_A\) and \(\fT_B\).

By the convention in \Cref{s:single_contour}, we choose the local descent
and ascent contours relative to one of the two triangles
\(\fT\in\{\fT_A,\fT_B\}\), according to the local configuration. Then
\begin{align}
J^{(1)}=J_\fT.
\end{align}
\end{lemma}

\begin{lemma}\label{p:J1vanish}
Adopt the assumptions of \Cref{p:concentric_setting}, and assume further
that
\[
(x,s)\in\fL
\qquad\text{or}\qquad
(y,t)\in\fL.
\]
Fix \(\fT\in\{\fT_A,\fT_B\}\), and assume that, for the local configuration
of the descent and ascent contours relative to \(\fT\), one of the local
contours
\[
\sfC^{\rm d}(w_0),
\qquad
\sfC^{\rm a}(w_0)
\]
can be deformed to the empty contour without crossing the other. Then there
exists a constant \(\fc'>0\) such that
\begin{align}\label{e:J1vanish}
|J_\fT|
=
\sum_{w_c,z_c}
\OO\left(
e^{-\fc'n}
e^{n\Re[S(w_c;x,s)-S(z_c;y,t)]}
\right).
\end{align}
Here \(w_c\) ranges over
\(\operatorname{Crit}^{\rm d}(x,s;\fN_\al)\), and \(z_c\) ranges over
\(\operatorname{Crit}^{\rm a}(y,t;\fN_{(y,t)})\).
\end{lemma}

\begin{lemma}[Concentric case]\label{p:concentric_setting}

Suppose that a chart associated with \(\fN_\al\) carrying a local descent
path and a chart associated with \(\fN_{(y,t)}\) carrying a local ascent
path are concentric, with common center
\(w_0\), and that \(w_0\) corresponds to a tangent or cusp-turning
location. The tangent line $L(w_0)$ separates two curvilinear triangles $\fT_A$ and $\fT_B$. 

Then it is impossible for one of \(\fN_\al\) and
\(\fN_{(y,t)}\) to be contained in \(\fT_A\cup\fL\) and the other to be
contained in \(\fT_B\cup\fL\).
By the convention in \eqref{s:single_contour}, we are using local contours for $\fT\in \{\fT_A,\fT_B\}$ for possible local configurations of the descent and ascent contours
Then 
\begin{align}
J^{(1)}=J_\fT.
\end{align}
\end{lemma}

\begin{lemma}\label{p:J1vanish}
Adopt assumptions in \eqref{p:concentric_setting}.
and sssume further that
\[
(x,s)\in\fL\quad  \text{or} \quad (y,t)\in \fL.
\]
Fix $\fT\in \{\fT_A, \fT_B\}$, and assume that for the local configurations of the descent and ascent contours for $\fT$,  one of the local contours
\[
\sfC^{\rm d}(w_0),
\qquad
\sfC^{\rm a}(w_0)
\]
can be deformed to \(\emptyset\) without crossing the other. Then there
exists a constant \(\fc'>0\) such that
\begin{align}\label{e:J1vanish}
|J_\fT|
=
\sum_{w_c,z_c}
\OO\left(
e^{-\fc'n}
e^{n\Re[S(w_c;x,s)-S(z_c;y,t)]}
\right).
\end{align}
Here \(w_c\) ranges over 
\(\operatorname{Crit}(x,s;\fN_\al)\), and \(z_c\) ranges over \(\operatorname{Crit}(y,t;\fN_{(y,t)})\). 
\end{lemma}

Without loss of generality, assume that
\[
\nabla H^*=(1,0)
\quad\text{on}\quad \fT.
\]
By \Cref{t:frozen_structure}, the two straight boundary pieces of \(\fT\)
determine a vertical tangent line \(L_\infty\) and a horizontal tangent line
\(L_0\). These two lines meet at the vertex \(\zeta\) of \(\fP\).
Moreover, all tangent lines to the portion of the arctic boundary contained
in \(\fT\) have slopes in \([-\infty,0]\), and this portion of the arctic
boundary is parametrized by
\[
[b_i,\infty]\subset\cC(\bR).
\]

\begin{lemma}[Ordering condition]
If $(y,t)\in \fT$, recall from \eqref{l:JT_nonvanish}, we have that if \begin{align}\label{e:orderc}
x< y \quad \text{or}\quad s\geq t 
\end{align} then $J_\fT=0$. This is what we will refer to as ordering condition.
\end{lemma}

\begin{lemma}[Phase separation]\label{l:phase_separation}
Suppose that the descent chart associated with \(\fN_\al\), centered at
\(w_0\), is disjoint from the ascent chart associated with
\(\fN_{(y,t)}\), centered at \(z_0\). Let \(w_c\) be a descent critical point
in the former chart, and let \(z_c\) be an ascent critical point in the latter
chart.

Assume that there exist real points \(\xi\) in the chart centered at \(w_0\)
and \(\xi'\) in the chart centered at \(z_0\), and a constant \(\fc'>0\),
such that one of the following two alternatives holds.

\smallskip
\noindent
\emph{First alternative.}
\[
        S(\xi;x,s)\leq \Re S(w_c;x,s),
        \qquad
        S(\xi';y,t)\geq \Re S(z_c;y,t)+\fc',
\]
 the gradient flow of \(S(\cdot;y,t)\) runs from \(\xi\) to \(\xi'\), and $\xi$ satisfies the assumptions in \Cref{l:single_contour_integral}

\smallskip
\noindent
\emph{Second alternative.}
\[
        S(\xi;x,s)\leq \Re S(w_c;x,s)-\fc',
        \qquad
        S(\xi';y,t)\geq \Re S(z_c;y,t),
\]
 the gradient flow of \(S(\cdot;x,s)\) runs from \(\xi\) to \(\xi'\), and $\xi'$ satisfies the assumptions in \Cref{l:single_contour_integral}

Then
\[
     |J_\fT|\leq e^{n (S(w_c;x,s)-S(z_c;y,t))-\fc'n}
\]
\end{lemma}

\begin{proof}
We prove the claim only under the first alternative. The proof under the
second alternative is identical, so we omit it.

By \Cref{l:single_contour_integral}, we have
\[
|J_\fT|
\leq
(\ln n)
e^{n\Re[S(\xi;x,s)-S(\xi;y,t)]}.
\]
The first alternative implies
\[
\begin{aligned}
\Re[S(\xi;x,s)-S(\xi;y,t)]
&\leq
\Re S(w_c;x,s)-\Re S(\xi';y,t) \\
&\leq
\Re S(w_c;x,s)-\Re S(z_c;y,t)-\fc'.
\end{aligned}
\]
After decreasing \(\fc'>0\) if necessary, the factor \(\ln n\) can be
absorbed into the exponentially small factor, and therefore
\[
|J_\fT|
\leq
e^{-\fc'n}
e^{n\Re[S(w_c;x,s)-S(z_c;y,t)]}.
\]
\end{proof}

\begin{proof}[Proof of \Cref{l:JTsmall}]
Under the assumption that
$
\operatorname{Cell}^{\rw}(y,t)\neq\fT,
$
there are two cases:
\begin{enumerate}
\item
The point \((y,t)\) lies in \(\fL\), and \(\fN_{(y,t)}\) is an arctic,
cusp, tangent, or cusp-turning neighborhood.

\item
The neighborhood \(\fN_{(y,t)}\) is a tangent, cusp-turning, or interface
frozen neighborhood that intersects two curvilinear triangles
\(\fT_A\) and \(\fT_B\), with
\[
\operatorname{Cell}^{\rw}(y,t)=\fT_A,
\qquad
\fT=\fT_B.
\]
\end{enumerate}

We discuss all cases one by one.

In this case, \(\fN_{(y,t)}\) is an arctic, cusp, tangent, cusp-turning or interface frozen
neighborhood. It is associated with an arctic, cusp, tangent, cusp-turning, or interface frozen
chart centered at \(z_0\). Let \(\ell_-(z_0;\fT)\) be the corresponding descent
cut. This cut is tangent to the arctic boundary at a point
\((y_0',t_0')\), where
\[
        z_0=y_0'-t_0'\chi(z_0).
\]
In the cusp and cusp-turning cases, the cut \(\ell_-(z_0;\fT)\) may degenerate to
a point.

Recall from \Cref{l:number_of_cuts} and its proof that, in addition to
\(\ell_-(z_0;\fT)\), there may be one further admissible descent cut, which we
denote by \(\ell_-(z_0';\fT)\). Such an additional cut can occur only if
\begin{itemize}
\item \(\fN_{(y,t)}\) intersects the second quadrant;
\item if \(\fN_{(y,t)}\cap L_\infty\neq\emptyset\), then
      \(\fN_{(y,t)}\) lies above the tangency location on \(L_\infty\);
\item if \(\fN_{(y,t)}\cap L_0\neq\emptyset\), then
      \(\fN_{(y,t)}\) lies to the left of the tangency location on \(L_0\).
\end{itemize}
If \(\fN_{(y,t)}\) contains the tangency location on \(L_\infty\) or on
\(L_0\), then the existence of the additional cut depends on the local
configuration. In all remaining cases, no additional admissible descent cut
occurs.

\medskip
\noindent
\textbf{Case 1: Arctic and cusp neighborhoods.}

In this case, \(\fN_{(y,t)}\) is associated with an arctic or cusp chart centered at \(z_0\). Let \(\ell_-(z_0;\fT)\) be the corresponding descent
cut. This cut is tangent to the arctic boundary at a point
\((y_0',t_0')\), where
\[
        z_0=y_0'-t_0'\chi(z_0).
\]
In the cusp case, the cut \(\ell_-(z_0;\fT)\) may degenerate to
a point.

Recall from \Cref{l:number_of_cuts} and its proof that, in addition to
\(\ell_-(z_0;\fT)\), there may be one further admissible descent cut, which we
denote by \(\ell_-(z_0';\fT)\).

For each tangent line, we view the line in its downward-pointing direction
in order to label its two sides as left and right. At a non-cusp point, the
arctic boundary lies locally on one of these two sides of the tangent line;
see Panels~(A) and~(B) of \Cref{f:CT}.

\smallskip
\noindent
\emph{Subcase 1: Cusp neighborhoods and the left-side arctic
configuration.}

Suppose that \(\fN_{(y,t)}\) is a cusp neighborhood, or that it is an
arctic neighborhood for which the arctic boundary lies locally to the left
of \(\ell_-(z_0;\fT)\), when viewed in its downward-pointing direction.
In the arctic case,
$
S'''(z_0;y_0',t_0')<0;
$
see Panel (A) of \Cref{f:arctic_boundary}.

If \(\fN_{(y,t)}\) intersects the second quadrant, then the boundary of
\(\fC(\fT;\fN_{(y,t)})\) contains another descent cut
\(\ell_-(z_0';\fT)\); see Panel~(C) of \Cref{f:CT}. In the cusp case,
the two descent cuts coincide, so \(z_0=z_0'\).

Otherwise, \(\fN_{(y,t)}\) lies either to the left or to the right of
\(L_\infty\). The former possibility is incompatible with the present
configuration, since inspection of \Cref{f:curvilinear_triangle} shows
that the arctic boundary would then lie locally to the right of
\(\ell_-(z_0;\fT)\). Hence \(\fN_{(y,t)}\) lies to the right of
\(L_\infty\), and \(L_\infty\) is the other boundary cut of
\(\fC(\fT;\fN_{(y,t)})\); see Panel~(D) of \Cref{f:CT}.

We first treat the configuration in Panel~(C) of \Cref{f:CT}. As one
moves along the arctic boundary from the tangency point
\((y_0',t_0')\) in the ascending direction, the corresponding
downward-pointing tangent ray varies continuously and sweeps out the entire
region \(\fC(\fT;\fN_{(y,t)})\). Therefore, there exists \(w_0\) such
that
\[
\fN_\al\cap\ell_-(w_0;\fT)\neq\emptyset.
\]
Thus, \(\fN_\al\) is associated with a local descent chart centered at
\(w_0\).

By the chart-disjointness assumption,
\begin{align}\label{e:cclow0}
z_0'<w_0<z_0.
\end{align}
In particular, this rules out the cusp configuration \(z_0'=z_0\).
Moreover, the fixed separation of the two charts gives
$
w_0\leq z_0-\fc.
$

Let \(z_c\) be an ascent critical point associated with
\(\fN_{(y,t)}\) and contained in the chart centered at \(z_0\). Since
\((y,t)\in\fL\), there is a complex-conjugate pair of such critical
points, and we choose either one. By the cubic estimate
\eqref{e:arctic_cubic} and the inequality
\(S'''(z_0;y_0',t_0')<0\),
\begin{align}\label{e:z_gap_left}
\Re S(z_0-\fc;y,t)
\geq
\Re S(z_c;y,t)+\fc'
\end{align}
for some \(\fc'>0\).

The chart associated with \(\fN_\al\), centered at \(w_0\), is an
arctic, cusp, or frozen chart. In the cusp case, the cusp points upward and
the corresponding fourth derivative is negative; see Panel (C) of \Cref{f:cusp}. We
distinguish two possibilities.

\begin{enumerate}
\item
Suppose first that the chart contains a real descent critical point \(w_c\).
Then the relation
$
w_c-x=-s\chi(w_c)\gtrsim1,
$
together with \eqref{e:z_gap_left}, verifies the assumptions of
\Cref{l:phase_separation} with
\[
\xi=w_c,
\qquad
\xi'=z_0-\fc.
\]
This gives
\begin{align}\label{e:J1_dd1}
|J_\fT|
\leq
e^{-\fc'n}
e^{n\Re[S(w_c;x,s)-S(z_c;y,t)]}.
\end{align}

\item
Otherwise, \((x,s)\in\fL\), and the chart centered at \(w_0\) is an
arctic or cusp chart. There is then a complex-conjugate pair of critical
points; denote either one by \(w_c\). By \eqref{e:arctic_cubic} or
\eqref{e:cusp_quartic}, there exists
\[
\xi\in\{w_0-\fc,w_0+\fc\}
\]
such that
\begin{align}\label{e:w_gap_left}
\Re S(w_c;x,s)
\geq
\Re S(\xi;x,s)+\fc'
\end{align}
for some \(\fc'>0\). Moreover,
$
\xi-x\geq w_0-\fc-x\gtrsim1.
$
By \eqref{e:z_gap_left} and \eqref{e:w_gap_left}, the assumptions of
\Cref{l:phase_separation} hold with this choice of \(\xi\) and with
\(\xi'=z_0-\fc\). Applying \Cref{l:phase_separation} gives
\eqref{e:J1_dd1}.
\end{enumerate}

We next treat the configuration in Panel~(D) of \Cref{f:CT}. The same
sweeping argument shows that \(\fN_\al\) is associated with a local
descent chart centered at some \(w_0\), where
\[
b_i\leq w_0<z_0.
\]
The equality \(w_0=b_i\) occurs precisely when
\[
\ell_-(w_0;\fT)=L_\infty.
\]

In this configuration, \(\fN_{(y,t)}\) is either an arctic neighborhood
with
\[
S'''(z_0;y_0',t_0')<0,
\]
or a cusp neighborhood whose cusp points downward and whose corresponding
fourth derivative is positive; see Panel (F) of \Cref{f:cusp}. The associated chart
contains an ascent critical point \(z_c\). In both cases,
\eqref{e:z_gap_left} holds, by \eqref{e:arctic_cubic} in the arctic
case and by \eqref{e:cusp_quartic} in the cusp case.

If \(w_0>b_i\), the same argument as in Panel~(C) gives
\eqref{e:J1_dd1}. It remains to treat the endpoint case \(w_0=b_i\).
Then \(\ell_-(w_0;\fT)=L_\infty\). If \(x>b_i\), the chart
centered at \(b_i\) contains a real descent critical point \(w_c>b_i\), and by \eqref{e:wc_vertical_tangent}, \eqref{e:wc_vertical_cusp} and \eqref{e:wcest}, we have $w_c-x\gtrsim |x-b_i|\geq 1/(2n)$. 
The same phase-separation argument using \Cref{l:phase_separation} gives \eqref{e:J1_dd1}.

If instead \(x\leq b_i\), while
\((y,t)\in\fN_{(y,t)}\) lies to the right of \(L_\infty\), then
\(y>x\). Hence the ordering condition \eqref{e:orderc} gives
\[
J_\fT=0,
\]
contrary to our standing assumption. This proves \eqref{e:JTsmall} in
Subcase~1.

\smallskip
\noindent
\emph{Subcase 2: The arctic boundary lies locally on the right of
\(\ell_-(z_0;\fT)\).}

We now assume that the arctic boundary lies locally on the right of
\(\ell_-(z_0;\fT)\). In the arctic case, this means \(S'''>0\); see
Panel (D) of \Cref{f:arctic_boundary}.

If \(\fN_{(y,t)}\) intersects the second quadrant, then the boundary of
\(\fC(\fT;\fN_{(y,t)})\) contains another descent cut, denoted
\(\ell_-(z_0';\fT)\). Otherwise, \(\fN_{(y,t)}\) lies either above or below
\(L_0\). The former possibility is incompatible with the present assumption,
since inspection of \Cref{f:curvilinear_triangle} shows that the arctic
boundary would then lie locally on the left. Hence \(\fN_{(y,t)}\) lies below
\(L_0\), and the horizontal cut \(L_0\) is the other boundary cut of
\(\fC(\fT;\fN_{(y,t)})\).

As before, moving along the arctic boundary from \((y_0',t_0')\) in the
ascending direction, the downward-pointing tangent ray sweeps the whole
region \(\fC(\fT;\fN_{(y,t)})\). Therefore, if
\[
        \fN_\al\cap\fC(\fT;\fN_{(y,t)})\neq\emptyset,
\]
then there exists a parameter \(w_0\) such that
\[
        \fN_\al\cap\ell_-(w_0)\neq\emptyset.
\]
Thus \(\fN_\al\) is associated with a local descent chart centered at \(w_0\).

By chart disjointness, if the second boundary cut is \(\ell_-(z_0')\), then
\[
        z_0<w_0<z_0'.
\]
If the second boundary cut is \(L_0\), then
\begin{align}\label{e:cclow}
        z_0<w_0\leq\infty.
\end{align}
The equality \(w_0=\infty\) occurs exactly when
\(\ell_-(w_0)=L_0\).

If \(w_0<\infty\), the same phase-separation calculation using \Cref{l:phase_separation},
with \(z_0-\fc\) replaced by \(z_0+\fc\), gives \eqref{e:J1_dd1}. The
same conclusion holds when \(w_0=\infty\) and \((x,s)\) lies below \(L_0\).
If instead \(w_0=\infty\) and \((x,s)\) lies above or on \(L_0\), while
\((y,t)\in\fN_{(y,t)}\) lies below \(L_0\), then the ordering condition
\eqref{e:orderc} gives
\[
        J_\fT=0.
\]
This completes Case 1.

\begin{figure} 
\begin{subfigure}{0.15\textwidth}
\centering

\caption{}
\end{subfigure}
\begin{subfigure}{0.16\textwidth}
\centering
%
\caption{}
\end{subfigure}
\begin{subfigure}{0.35\textwidth}
    \centering
    %
    \caption{}
\end{subfigure}
\begin{subfigure}{0.3\textwidth}
    \centering
    %
    \caption{}
\end{subfigure}
  \caption{Configurations}
  \label{f:CT}
\end{figure}

\begin{figure}[ht]
\begin{subfigure}{0.23\textwidth}
\centering
%
\caption{}
\end{subfigure}
\caption{$\fN_{(y,t)}$ is a tangent or cusp-turning neighborhood.}\label{f:cuts}
\end{figure}

\medskip
\noindent
\textbf{Case 2: Tangent and cusp-turning neighborhoods.}
By our assumption, $(y,t)\in \fL\cup \fT_A$.
We first consider the case in which \(\fN_{(y,t)}\) is a vertical tangent or
cusp-turning neighborhood. 
It is associated with a descent cut \(\ell_-(z_0;\fT)\), which lies on the vertical tangent line
\(L_\infty\), and $z_0=b_i$. There
are four possible configurations, shown in \Cref{f:cuts}.

In Panel (A) of \Cref{f:cuts}, we have 
\[
        \fC(\fT;\fN_{(y,t)})=\ell_-(z_0;\fT).
\]
Thus, on the event
\[
        \fN_\al\cap\fC(\fT;\fN_{(y,t)})\neq\emptyset,
\]
we have
\[
        \fN_\al\cap\ell_-(z_0;\fT)\neq\emptyset.
\]
This contradicts the standing assumption that the descent chart associated
with \(\fN_\al\) is disjoint from the ascent chart associated with
\(\fN_{(y,t)}\). 

In Panel (B) of \Cref{f:cuts}, if
\[
        \fN_\al\cap L_\infty\neq\emptyset,
\]
then chart disjointness forces \(\fN_\al\) to be above of \(\fN_{(y,t)}\). Hence the ordering condition
\eqref{e:orderc} gives
\[
        J_\fT=0.
\]
If instead \(\fN_\al\) does not intersect \(L_\infty\), then it lies either above $\fN_{(y,t)}$, or $x\leq b_i$. 
Since $(y,t)\in \fL\cap \fT_A$, we have $y\geq b_i$.  Thus by the same ordering condition \eqref{e:orderc}, in both cases we have
\[
        J_\fT=0.
\]

In Panel (C) of \Cref{f:cuts}, \(\fN_{(y,t)}\) is associated with an
additional descent cut \(\ell_-(z_0';\fT)\), where
\[
b_i=z_0<z_0'.
\]
As in \eqref{e:cclow0}, there exists \(w_0\) satisfying
\[
z_0<w_0<z_0'
\]
such that
\[
\fN_\al\cap\ell_-(w_0;\fT)\neq\emptyset,
\]
and the associated chart contains a descent critical point \(w_c\).
Since \((y,t)\in\fL\cap\fT_A\), the chart centered at \(z_0\) contains an
ascent critical point \(z_c\). The claim \eqref{e:J1_dd1} then follows from
the phase-separation estimate in \Cref{l:phase_separation}.

In Panel (D) of \Cref{f:cuts}, \(\ell_-(z_0;\fT)\) degenerates to a
point. If
\[
\fN_\al\cap L_\infty\neq\emptyset,
\]
then the disjointness of the charts forces \(\fN_\al\) to lie above
\(\fN_{(y,t)}\). Hence the ordering condition \eqref{e:orderc} gives
\[
J_\fT=0.
\]
Suppose instead that \(\fN_\al\) does not intersect \(L_\infty\). Then
either \(\fN_\al\) lies above \(\fN_{(y,t)}\), or \(x\leq b_i\). In the
first case, \eqref{e:orderc} again gives \(J_\fT=0\). If
\((y,t)\in\fT_A\) and \(y\geq b_i\), then, in the second case,
\(x\leq b_i\leq y\), so the same ordering condition implies that
\(J_\fT=0\).

It remains to consider the case in which either \((y,t)\in\fL\), or
\((y,t)\in\fT_A\) and \(y<b_i\). In both cases, the chart centered at
\(z_0\) contains an ascent critical point. The claim
\eqref{e:J1_dd1} then follows by the same phase-separation argument as in
Panel (C).

This completes Case 2.

\medskip
\noindent
\textbf{Case 3: Interface frozen neighborhoods.}
This is similar to the tangent neighborhood case, so we omit.

\end{proof}

\begin{figure}

\begin{subfigure}{0.23\textwidth}
\centering

\end{subfigure}
  \caption{Configurations}
  \label{f:CT2}
\end{figure}

\begin{proof}[Proof of \Cref{l:JTterm}]
In this case, we use the pointwise curvilinear region
\(\fC(\fT;y,t)\). Since \(\fN_\al\) and \(\fN_{(y,t)}\) do not share a
concentric chart, the definitions of the regions
\(\fC(\fT;\fN_{(y,t)})\) and \(\fC(\fT;y,t)\) in
\Cref{s:deffC,s:deffC2}, together with the ordering condition
\eqref{e:orderc}, imply that
\[
\bm1\bigl(\fN_\al\cap\fC(\fT;\fN_{(y,t)})\neq\emptyset\bigr)J_\fT
=
\bm1\bigl(\fN_\al\cap\fC(\fT;y,t)\neq\emptyset\bigr)J_\fT.
\]
Therefore, it is enough to prove
\begin{align}\label{e:J1_small}
\bm1\bigl(
\fN_\al\cap\fC(\fT;y,t)\neq\emptyset,\,
(x,s)\notin\fC(\fT;y,t)
\bigr)|J_\fT|
=
\sum_{w_c,z_c}
\OO\left(
e^{-\fc'n}
e^{n\Re[S(w_c;x,s)-S(z_c;y,t)]}
\right).
\end{align}
Indeed, the complementary condition
\((x,s)\in\fC(\fT;y,t)\) gives precisely the interlacing contribution.

There are three possible ways in which the indicator in \eqref{e:J1_small}
can be nonzero:
\begin{enumerate}
\item \(\fN_\al \cap \fC(\fT;y,t)\neq \emptyset\) and \((x,s)\in\fL\);

\item \(\fN_\al \cap \fC(\fT;y,t)\neq \emptyset\) and
$
        (x,s)\in\fT\setminus\fC(\fT;y,t).
$
This can happen only if \(\fN_{(y,t)}\) is on the right of \(L_\infty\) and
\((x,s)\) is on the left of \(L_\infty\), or if \(\fN_{(y,t)}\) is below
\(L_0\) and \((x,s)\) is above \(L_0\).

\item \((x,s)\) belongs to an adjacent curvilinear triangle. More precisely,
write
$
        \fT_B:=\fT,
$
and let \(\fT_A\) be the adjacent curvilinear triangle containing \((x,s)\).
Thus
$
        (x,s)\in\fT_A,
 $ and $
        (x,s)\notin\fT_B.
$
In this case \(\fN_\al\) is a tangent, cusp-turning, or interface frozen
neighborhood, and
\[
        \fN_\al\cap(\fT_A\cap\fT_B)\neq\emptyset.
\]
Moreover, it is associated with a tangent, cusp-turning, or interface frozen
chart centered at \(w_0\), which carries local descent paths.
\end{enumerate}

In the first case, the same detailed phase estimate as in \Cref{l:JTsmall} gives the
bound \eqref{e:J1_small}. In the second case, we have either \(x<y\) or
\(s> t\), and hence the ordering condition \eqref{e:orderc} gives
$
        J_\fT=0.
$
It remains to treat the adjacent-triangle case.

\medskip
\noindent
\textbf{Case 1: The adjacent triangles share a vertical boundary.}

Suppose first that \(\fT_A\) and \(\fT_B\) share a vertical boundary piece,
corresponding to \(w_0=b_i\).

Since
\[
        \fN_\al\cap(\fT_A\cap\fT_B)\neq\emptyset,
\]
by the classification in \Cref{f:curvilinear_triangle}, either \(\fN_\al\)
lies above \(L_0\), or it lies below \(L_0\).

If \(\fN_\al\) lies above \(L_0\), in the configurations of the first and
fourth rows of \Cref{f:curvilinear_triangle}, there are two possibilities as
in the first two panels of \Cref{f:CT2}. Then \((x,s)\in\fT_A\) lies to the
left of \(L_\infty\), so \(x<w_0\). In these cases,
$
        \fN_\al \cap \fC(\fT;y,t)\neq \emptyset
$
and chart disjointness imply that \(\fN_{(y,t)}\) lies in the first quadrant,
so \(y\geq w_0\). Thus
\[
        x<w_0\leq y,
\]
and the ordering condition \eqref{e:orderc} gives
$
        J_\fT=0.
$

If \(\fN_\al\) lies below \(L_0\), in the configurations of the second and
third rows of \Cref{f:curvilinear_triangle}, then it either contains the
tangent or cusp-turning location along \(L_\infty\), or lies above that
location; see the third panel of \Cref{f:CT2}.

We first rule out the possibility that \(\fN_\al\) is a cusp-turning
neighborhood. If \(\fN_\al\) is a cusp-turning neighborhood, then the
condition
$
        \fN_\al\cap\fC(\fT;y,t)\neq\emptyset
$
would force \(\fN_{(y,t)}\) to be a cusp-turning neighborhood as well. This
would be the concentric-chart case, contrary to the standing disjointness
assumption. Hence \(\fN_\al\) is either a tangent neighborhood or an
interface-frozen neighborhood.

We now consider the position of \(\fN_{(y,t)}\) relative to \(L_\infty\).

\begin{enumerate}
\item If \(\fN_{(y,t)}\) lies to the right of \(L_\infty\), then either
$
        \fN_\al\cap\fC(\fT;y,t)=\emptyset,
$
in which case the indicator in \eqref{e:J1_small} vanishes, or
\((x,s)\in\fT_A\) lies above \((y,t)\). In the latter case, we have
$
        s>t,
$
and hence the ordering condition \eqref{e:orderc} gives
$
        J_\fT=0.
$

\item If \(\fN_{(y,t)}\) lies on the vertical line \(L_\infty\), then there are
two possibilities. Either it lies on or above the tangent or cusp-turning
location along \(L_\infty\), which is the concentric-chart case and is
excluded by the standing disjointness assumption; or it lies below the tangent
or cusp-turning location. In the latter case, \((x,s)\in\fT_A\) lies above
\((y,t)\in\fT\), so
$
        s>t,
$
and the ordering condition \eqref{e:orderc} gives
$
        J_\fT=0.
$

\item Finally, suppose that \(\fN_{(y,t)}\) lies to the left of
\(L_\infty\). Then it is associated with a cut \(\ell_-(z_0;\fT)\), the arctic
boundary lies locally on the right of this cut, and the corresponding
chart contains an ascent critical point \(z_c\); see the third panel of
\Cref{f:CT2}.

Since \((x,s)\in\fT_A\) lies to the right of \(L_\infty\), we are in the
situations shown in Panels (B) and (D) of \Cref{f:tangent1}, or in Panel (B)
of \Cref{f:side}. Thus the chart centered at \(w_0\) contains a descent
critical point \(w_c<w_0\). By \eqref{e:tangent_S} and
\eqref{e:tangent_frozen_S},
\begin{align}\label{e:tangent_gap_reorganized}
        \Re S(w_c;x,s)
        \geq
        \Re S(w_0+\fc;x,s)+\fc'.
\end{align}
Moreover, \((x,s)\) is associated with another frozen chart centered at some
\(w_0'\), with
\[
        w_0<z_0<w_0',
\]
and the gradient flow of \(S(\cdot;x,s)\) runs from \(w_0+\fc\) to \(w_0'\).
Then the assumptions of \Cref{l:phase_separation} hold with
$
        \xi=w_0+\fc,
$ and $
        \xi'=z_c.
$
Therefore
\begin{align}\label{e:J1bb2}
        |J_\fT|
        \leq
        e^{-\fc'n}
        e^{n\Re[S(w_c;x,s)-S(z_c;y,t)]}.
\end{align}
This is the desired estimate.
\end{enumerate}

\medskip
\noindent
\textbf{Case 2: The adjacent triangles share a horizontal boundary.}

It remains to suppose that \(\fT_A\) and \(\fT_B\) share a horizontal boundary
piece, corresponding to $w_0=\infty$. By the classification in \Cref{f:curvilinear_triangle}, either
\(\fN_\al\) lies to the left of \(L_\infty\), or it lies to the right of
\(L_\infty\).

If \(\fN_\al\) lies to the left of \(L_\infty\), in the configurations of the
second and fourth rows of \Cref{f:curvilinear_triangle}, then
\[
        \fN_\al \cap \fC(\fT;y,t)\neq \emptyset
\]
and chart disjointness imply that \(\fN_{(y,t)}\) lies in the third quadrant,
and \((y,t)\) is on or below \(L_0\). Then
$
        s\geq t,
$
and the ordering condition \eqref{e:orderc} gives
$
        J_\fT=0.
$

If \(\fN_\al\) lies to the right of \(L_\infty\), in the configurations of the
first and third rows of \Cref{f:curvilinear_triangle}, then it either contains
the tangent or cusp-turning location along \(L_0\), or lies to the left of
that location. Again, the disjointness assumption rules out the possibility
that \(\fN_\al\) is a cusp-turning neighborhood. Hence \(\fN_\al\) is either a
tangent neighborhood or an interface-frozen neighborhood.

We now consider the position of \(\fN_{(y,t)}\) relative to \(L_0\), as in the
vertical-boundary case.

\begin{enumerate}
\item If \(\fN_{(y,t)}\) lies on or below \(L_0\), then the disjointness
assumption implies that
$
        x<y,
$
and the ordering condition \eqref{e:orderc} gives
$
        J_\fT=0.
$

\item If \(\fN_{(y,t)}\) lies above \(L_0\), then the phase-separation estimate
from \Cref{l:phase_separation} gives \eqref{e:J1bb2}.
\end{enumerate}

The claim \eqref{e:JTterm} follows from combining Case 1 and Case 2.
\end{proof}

\begin{proof}[Proof of \Cref{p:concentric_setting}]
We assume that the tangent line \(L(w_0)\) is vertical and denote the
tangency location by \((y_0,t_0)\); see Panels (C) and (D) of
\Cref{f:l}. The cases in which the tangent line is horizontal or of unit-slope are analogous. On \(\fT_A\), we have
$
\nabla H^*=(0,0),
$
and on \(\fT_B\), we have
$
\nabla H^*=(1,0).
$

Each of \(\fN_{(y,t)}\) and \(\fN_\al\) is of one of the following
types: a tangent or cusp-turning neighborhood, an interface frozen
neighborhood, a  neighborhood disjoint from \(\fT_B\), or a 
neighborhood disjoint from \(\fT_A\). If \(\fN_{(y,t)}\) is 
disjoint from \(\fT_B\), then
\[
\fN_{(y,t)}\subset\fT_A\cup\fL,
\]
and it lies below the tangency location \((y_0,t_0)\). If
\(\fN_{(y,t)}\) is disjoint from \(\fT_A\), then
\[
\fN_{(y,t)}\subset\fT_B\cup\fL,
\]
and it lies above \((y_0,t_0)\).

For \(\fN_\al\), the same two containments hold, but the vertical
positions are reversed.

In summary, each of the neighborhoods \(\fN_\al\) and
\(\fN_{(y,t)}\) is of one of the following types:
\begin{align}\begin{split}\label{e:configuration}
&\text{a tangent or cusp-turning neighborhood, an interface frozen
neighborhood,} \\
&\text{a neighborhood in \(\fT_B\cup\fL\), or a 
neighborhood in \(\fT_A\cup\fL\).}
\end{split}\end{align}

Next we show that it is impossible for one of \(\fN_\al\) and
\(\fN_{(y,t)}\) to be contained in \(\fT_A\cup\fL\) and the other to be
contained in \(\fT_B\cup\fL\). Indeed, in this case, the two neighborhoods
would lie on the same side of the tangency location \((y_0,t_0)\). On the
other hand, one neighborhood intersects the segment
\((\fT_A\cap L(w_0))\setminus (\fT_A\cap \fT_B)\), while the other intersects the segment
\((\fT_B\cap L(w_0))\setminus (\fT_A\cap \fT_B)\), so they must lie on opposite sides of
\((y_0,t_0)\), a contradiction.

Therefore, if either \(\fN_\al\) or \(\fN_{(y,t)}\) is contained in
\(\fT_A\cup\fL\), we take \(\fT=\fT_A\). If either is contained in
\(\fT_B\cup\fL\), we take \(\fT=\fT_B\). In the remaining cases, we may
take \(\fT\) to be either \(\fT_A\) or \(\fT_B\).

By \Cref{l:boundary_piece_in_C}, we always have
\begin{align}
\fN_\al\cap\fC(\fT;\fN_{(y,t)})\neq\emptyset.
\end{align}
Consequently,
\begin{align}
J^{(1)}
&=
\bm1\bigl(
\fN_\al\cap\fC(\fT;\fN_{(y,t)})\neq\emptyset
\bigr)J_\fT
=
J_\fT.
\end{align}
\end{proof}

\begin{proof}[Proof of \Cref{p:J1vanish}]
We recall the possible types of the neighborhoods \(\fN_\al\) and
\(\fN_{(y,t)}\) from \eqref{e:configuration}. In all cases, we use the
contours in Panel (A) or Panel (B) of
\[
\Cref{f:vertical_tangent1},\quad
\Cref{f:vertical_tangent2},\quad
\Cref{f:vertical_tangent3},\quad
\Cref{f:vertical_tangent4},\quad
\Cref{f:c_vertical_cusp1},\quad
\Cref{f:c_vertical_cusp2}
\]
for the possible local configurations of the descent and ascent contours.

By symmetry between \((x,s)\) and \((y,t)\), we only treat the case
\[
(y,t)\in\fL.
\]
The case \((x,s)\in\fL\) can be treated in the same way, so we omit it.
Throughout the rest of the proof, each case will be handled by one of the
following two mechanisms.

First, let
\[
J_A:=J_{\fT_A},
\qquad
J_B:=J_{\fT_B};
\]
recall \eqref{e:J1form}. By \Cref{l:JT_nonvanish}, \(J_A\neq0\) only if
\begin{align}\label{e:JA_nonvanish}
x\geq y,
\qquad
y-t\geq x-s,
\end{align}
and \(J_B\neq0\) only if
\begin{align}\label{e:JB_nonvanish}
x\geq y,
\qquad
t>s.
\end{align}

Second, if the descent chart associated with \(\fN_\al\) contains a
descent critical point \(w_c\), the ascent chart associated with
\(\fN_{(y,t)}\) contains an ascent critical point \(z_c\), and the two
charts satisfy the phase-separation condition in
\Cref{l:phase_separation}, then
\begin{align}\label{e:J1phase}
|J_\fT|
\leq
e^{-\fc'n}
e^{n\Re[S(w_c;x,s)-S(z_c;y,t)]}.
\end{align}

\medskip
\noindent
\textbf{Case 1. \(\fN_{(y,t)}\) is a tangent or cusp-turning
neighborhood.}
See the first panel of \Cref{f:concentric_chart}.

In this case, the chart associated with \((y,t)\), centered at \(w_0\),
contains a complex-conjugate pair of ascent critical points
\(z_c,\overline{z_c}\). The local ascent contour
\(\sfC^{\rm a}(w_0)\) is deformed to
$
\sfD^{\rm a}(z_c)\cup\sfD^{\rm a}(\overline{z_c}).
$
By assumption, the local descent contour \(\sfC^{\rm d}(w_0)\) associated
with \((x,s)\) can be deformed to the empty contour without crossing the
local ascent contour.

This can happen only if
\[
x\leq w_0
\]
and \(\fN_\al\) is an interface frozen neighborhood, a neighborhood
contained in \(\fT_A\cup\fL\), or a  neighborhood contained in
\(\fT_B\cup\fL\). We now treat all possible configurations one by one.

\begin{enumerate}
\item\label{i:cc1}
Consider the configuration in Panel (A) of
\Cref{f:vertical_tangent1}. In this case, we have
\(\fT=\fT_A\). Moreover, \(y<w_0\), and the two ascent paths
\(\sfD^{\rm a}(z_c)\) and \(\sfD^{\rm a}(\overline{z_c})\) meet at a
point
\[
\zeta\in[y,w_0].
\]

Since the local descent contour \(\sfC^{\rm d}(w_0)\) associated with
\((x,s)\) can be deformed to the empty contour, there are three
possibilities, which we discuss one by one.

First, suppose that \(\fN_\al\) is an interface frozen neighborhood, with
\(x\leq w_0\) and \((x,s)\in\fT_A\). In this case, \(\fN_\al\) is
associated with another frozen chart centered at \(w_0'<w_0\), which
contains a descent critical point \(w_c\) corresponding to a tangent line
from \((x,s)\) to the portion of the arctic boundary contained in
\(\fT_A\); see the first panel of \Cref{f:concentric_chart}. The charts
centered at \(w_0'\) and \(w_0\) are disjoint. Hence we may apply the
phase-separation estimate \Cref{l:phase_separation}, taking
\[
\xi=w_c,
\qquad
\xi'=\zeta.
\]
Indeed, the gradient flow of \(\Re S(\,\cdot\,;y,t)\) runs from \(w_c\)
toward \(w_0\). Therefore, \eqref{e:J1phase} holds.

Second, suppose that \(\fN_\al\) is a neighborhood contained in
\(\fT_B\cup\fL\) and lies below \((y_0,t_0)\). In this case, \(t>s\), so
the two conditions in \eqref{e:JA_nonvanish} cannot hold simultaneously.
Therefore,
\[
J_A=0.
\]

Third, suppose that \(\fN_\al\) is a frozen neighborhood contained in
\(\fT_A\cup\fL\) and lies above the tangency location \((y_0,t_0)\).
This is impossible in the configuration of
\Cref{f:vertical_tangent1}.

\item\label{i:cc2}
Consider the configuration in Panel (B) of
\Cref{f:vertical_tangent2}. In this case, we have
\(\fT=\fT_B\). Moreover, \(y<w_0\), and the two ascent paths
\(\sfD^{\rm a}(z_c)\) and \(\sfD^{\rm a}(\overline{z_c})\) meet at a
point
\[
\zeta\in[y,w_0].
\]
Again, we discuss the three possibilities for \(\fN_\al\) one by one.

First, suppose that \(\fN_\al\) is an interface frozen neighborhood, with
\(x\leq w_0\) and \((x,s)\in\fT_B\). In this case, \(\fN_\al\) is
associated with another frozen chart centered at \(w_0'>w_0\), which
contains a descent critical point \(w_c\) corresponding to a tangent line
from \((x,s)\) to the portion of the arctic boundary contained in
\(\fT_B\). The charts centered at \(w_0'\) and \(w_0\) are disjoint, so
the same phase-separation estimate as in \eqref{e:J1phase} applies.

Second, suppose that \(\fN_\al\) is a frozen neighborhood contained in
\(\fT_A\cup\fL\) and lies above \((y_0,t_0)\). In this case, \(s>t\), so
\eqref{e:JB_nonvanish} gives
\[
J_B=0.
\]

Third, suppose that \(\fN_\al\) is a frozen neighborhood contained in
\(\fT_B\cup\fL\) and lies below the tangency location \((y_0,t_0)\).
This is impossible in the configuration of
\Cref{f:vertical_tangent2}.

\item\label{i:cc3}
Consider the configuration in Panel (B) of
\Cref{f:vertical_tangent3}. In this case, \(\fT=\fT_B\). Since
\((y,t)\in\fL\), we have \(y>w_0\), and hence
\[
x\leq w_0<y.
\]
Therefore, by \eqref{e:JB_nonvanish},
\[
J_B=0.
\]

\item\label{i:cc4}
Consider the configuration in Panel (A) of
\Cref{f:vertical_tangent4}. In this case, \(\fT=\fT_A\), and the geometry
again gives
\[
x<y.
\]
Therefore, by \eqref{e:JA_nonvanish},
\[
J_A=0.
\]

\item\label{i:cc5}
Consider the configuration in Panel (B) of
\Cref{f:c_vertical_cusp1}. This case is handled by the same argument as in
\Cref{i:cc2}.

\item
Consider the configuration in Panel (A) of
\Cref{f:c_vertical_cusp2}. This case is handled by the same argument as in
\Cref{i:cc1}.
\end{enumerate}

\medskip
\noindent
\textbf{Case 2. \(\fN_{(y,t)}\) is a neighborhood contained in
\(\fT_B\cup\fL\).}
See the second panel of \Cref{f:concentric_chart}.

In this case, \(\fN_{(y,t)}\) lies above the tangency location
\((y_0,t_0)\). Since \((y,t)\in\fL\), the neighborhood
\(\fN_{(y,t)}\) contains part of the arctic boundary. Hence it is
associated with another chart centered at \(z_0'>w_0\), which contains a
complex-conjugate pair of ascent critical points; see the second panel of
\Cref{f:concentric_chart}.

This can happen only in the configurations appearing in
\[
\Cref{f:vertical_tangent1},\quad
\Cref{f:vertical_tangent2},\quad
\Cref{f:vertical_tangent3},\quad
\Cref{f:c_vertical_cusp1}.
\]
Equivalently, the possible triangles \(\fT_B\) are those shown in the
second and fourth panels of the first and fourth rows of
\Cref{f:curvilinear_triangle}, together with the configurations in the
third row.

We first claim that \(\fN_\al\) must intersect \(\fT_B\). Otherwise,
\[
\fN_\al\subset\fT_A\cup\fL,
\]
which contradicts \Cref{p:concentric_setting}. Thus, \(\fN_\al\) is 
a tangent or cusp-turning neighborhood, an interface frozen neighborhood,
or a neighborhood contained in \(\fT_B\cup\fL\). In all cases,
\(\fN_\al\) lies below \(\fN_{(y,t)}\). Therefore, \(t>s\), and the two
conditions in \eqref{e:JA_nonvanish} cannot hold simultaneously. Hence
\[
J_A=0.
\]

We now treat all remaining configurations one by one.
\begin{enumerate}
\item
Consider the configuration in Panel (B) of
\Cref{f:vertical_tangent1}. By assumption,
\(\sfC^{\rm a}(w_0)\) can be deformed to the empty contour, so
\(y\geq w_0\). If
$
x\leq w_0\leq y,
$
then by \eqref{e:JB_nonvanish},
\[
J_B=0.
\]
Otherwise, \(x>w_0\). In the configuration of
\Cref{f:vertical_tangent1}, we then have \((x,s)\in\fT_B\), and the chart
centered at \(w_0\) associated with \((x,s)\) contains a real descent
critical point \(w_c>w_0\). This chart and the chart centered at \(z_0'\)
associated with \((y,t)\) satisfy the phase-separation condition in
\Cref{l:phase_separation}. Hence \eqref{e:J1phase} holds with
\(\fT=\fT_B\):
\begin{align}\label{e:JBBB}
|J_B|
\leq
e^{-\fc'n}
e^{n\Re[S(w_c;x,s)-S(z_c;y,t)]}.
\end{align}

\item
Consider the configurations in Panel (B) of
\(\Cref{f:vertical_tangent2}\), \(\Cref{f:vertical_tangent3}\), or
\(\Cref{f:c_vertical_cusp1}\). By assumption,
\(\sfC^{\rm d}(w_0)\) can be deformed to the empty contour, so
\(x\leq w_0\). Moreover, in the configuration of
\(\Cref{f:vertical_tangent2}\), \(\fN_\al\) is not a tangent
neighborhood, while in the configuration of
\(\Cref{f:c_vertical_cusp1}\), it is not a cusp-turning neighborhood. In
all cases, the neighborhood \(\fN_\al\) is associated with another frozen
chart centered at some
\[
w_0'<z_0',
\]
which contains a descent critical point \(w_c\); see the second panel of
\Cref{f:concentric_chart}.

If
$
x\leq w_0\leq y,
$
then by \eqref{e:JB_nonvanish},
\[
J_B=0.
\]
Otherwise, \(y<w_0\). The chart centered at \(z_0'\), associated with
\((y,t)\), contains an ascent critical point \(z_c\). The charts centered
at \(w_0'\) and \(z_0'\) satisfy the phase-separation condition. Therefore,
\eqref{e:J1phase} holds with \(\fT=\fT_B\), and \eqref{e:JBBB} follows.
\end{enumerate}
This completes Case 2.

\medskip
\noindent
\textbf{Case 3. \(\fN_{(y,t)}\) is a neighborhood contained in
\(\fT_A\cup\fL\).}

This case is symmetric to Case 2, so we omit it.
\end{proof}

	\begin{figure}
	 \begin{subfigure}{0.4\textwidth}


  \end{subfigure}
		
				\caption{
				$\fN_\al$ and $\fN_{(y,t)}$ are associated with concentric charts.
}\label{f:concentric_chart}
\end{figure}

\subsection{Disjoint case}

Suppose that every chart associated with \(\fN_{(y,t)}\) is disjoint from
every chart associated with \(\fN_\al\). We consider the case in which
\(\fN_\al\) is not a ramification neighborhood. The case in which
\(\fN_\al\) is a ramification neighborhood can be proved in the same way,
so we omit it.

We may deform each local descent contour \(\sfC^{\rm d}(w_0)\) to the
steepest-descent paths \(\sfD^{\rm d}(w_c)\), and each local ascent contour
\(\sfC^{\rm a}(z_0)\) to the steepest-ascent paths
\(\sfD^{\rm a}(z_c)\) introduced in \Cref{s:critical_point}. This gives
\begin{align}\begin{split}
&\phantom{{}={}}A_\al((x,s),(y,t))
=
J^{(1)}
+\sum_{w_c,z_c}
\OO\left(
e^{-\fc'n}
e^{n\Re[S(w_c;x,s)-S(z_c;y,t)]}
\right),
\\
&+
\sum_{w_c,z_c}
\frac{n}{(2\pi \ri)^2}
\int_{\sfD^{\rm a}(z_c)}
\!\!\int_{\sfD^{\rm d}(w_c)}
P_{ns}(nw,nx)\,Q_{nt}(nz,ny)\,
\frac{I_+(w)}{I_-(z)}\,
\frac{\sqrt{\phi'(w)}\sqrt{\phi'(z)}}
{\phi(w)-\phi(z)}\,
\rd w\,\rd z,
\end{split}\end{align}
where \(w_c\) ranges over
\(\operatorname{Crit}^{\rm d}(x,s;\fN_\al)\), and \(z_c\) ranges over
\(\operatorname{Crit}^{\rm a}(y,t;\fN_{(y,t)})\).

By \Cref{l:JTsmall,l:JTterm}, according as \((y,t)\in\fL\) or
\(\operatorname{Cell}^{\rw}(y,t)=\fT\), we have
\begin{align}
J^{(1)}
&=
\wh J^{(0)}
+
\sum_{w_c,z_c}
\OO\left(
e^{-\fc'n}
e^{n\Re[S(w_c;x,s)-S(z_c;y,t)]}
\right),
\end{align}
where \(\wh J^{(0)}\) is as defined in \Cref{p:standard_form}. Substituting
this estimate into the preceding expansion proves \Cref{p:standard_form}
in the disjoint case.

\subsection{Concentric case}

Suppose that the neighborhoods \(\fN_\al\) and \(\fN_{(y,t)}\) share
concentric charts centered at \(w_0\in\cC(\bR)\). The shared center \(w_0\)
determines the tangent line
\begin{align}
L(w_0)
:=
\{(y',t'):\ w_0=y'-t'\chi(w_0)\},
\end{align}
which is tangent to the arctic boundary at some point
\((y_0',t_0')\in\fA\).

The double-contour integrals corresponding to pairs of disjoint charts can
be treated in exactly the same way as in the disjoint case. Thus,
\Cref{p:standard_form} follows from the following proposition.

\begin{proposition}\label{p:concentric}
The sum of the single-contour integral term and the double-contour integral
associated with the shared concentric charts,
\begin{align}
\begin{split}\label{e:sumoftwo}
J^{(1)}
+
\frac{n}{(2\pi\ri)^2}
\int_{\sfC^{\rm a}(w_0)}
\!\!\int_{\sfC^{\rm d}(w_0)}
P_{ns}(nw,nx)\,Q_{nt}(nz,ny)\,
\frac{I_+(w)}{I_-(z)}\,
\frac{\sqrt{\phi'(w)}\sqrt{\phi'(z)}}
{\phi(w)-\phi(z)}\,
\rd w\,\rd z,
\end{split}
\end{align}
has the form \eqref{e:invariant}. We may deform
\(\sfC^{\rm d}(w_0)\) to the steepest-descent paths
\(\sfD^{\rm d}(w_c)\), and \(\sfC^{\rm a}(w_0)\) to the
steepest-ascent paths \(\sfD^{\rm a}(z_c)\) introduced in
\Cref{s:critical_point}. This gives
\begin{align}
\begin{split}\label{e:sumoftwo2}
&\phantom{{}={}}\eqref{e:sumoftwo}
=
\wh J^{(0)}+\wh J^{(1)}
+
\sum_{w_c,z_c}
\OO\left(
e^{-\fc'n}
e^{n\Re[S(w_c;x,s)-S(z_c;y,t)]}
\right)\\
&+
\sum_{w_c,z_c}
\frac{n}{(2\pi\ri)^2}
\int_{\sfD^{\rm a}(z_c)}
\!\!\int_{\sfD^{\rm d}(w_c)}
P_{ns}(nw,nx)\,Q_{nt}(nz,ny)\,
\frac{I_+(w)}{I_-(z)}\,
\frac{\sqrt{\phi'(w)}\sqrt{\phi'(z)}}
{\phi(w)-\phi(z)}\,
\rd w\,\rd z
\end{split}
\end{align}
where \(\wh J^{(0)}\) and \(\wh J^{(1)}\) are as defined in
\Cref{p:standard_form}, \(w_c\) ranges over the descent critical points in
\(\operatorname{Crit}^{\rm d}(x,s;\fN_\al)\) contained in the shared chart,
and \(z_c\) ranges over the ascent critical points in
\(\operatorname{Crit}^{\rm a}(y,t;\fN_{(y,t)})\) contained in the shared
chart.
\end{proposition}

In the rest of this section we prove \Cref{p:concentric}. There are several cases for $w_0$.

\noindent\textbf{Case 1: \(w_0\) corresponds to a regular arctic point $(x_0,s_0)$; see panel (A) of \Cref{f:l}. }

The chart associated with \((y,t)\) is either an arctic chart, or a frozen chart with
$S''(w_0;\cdot)<0$,
and in either case carries a local ascent path. 
Similarly, the chart associated with \((x,s)\) is either an arctic chart, or a frozen chart with
$
S''(w_0;\cdot)>0$,
and carries a local descent path.

Without loss of generality, we assume that at the tangency location $(x_0,s_0)\in \fT$ for some curvilinear triangle $\fT$ and satisfies
\[
S'''(w_0;x_0,s_0)<0.
\]
We recall from \Cref{f:frozen_path1} all possible local configurations of the relevant paths. 
After a slight perturbation, the local descent paths (blue) lie to the left of the local ascent paths (red).

Assume first that \((x,s),(y,t)\notin \fL\), and let \(w_c\in \bR\) be the associated descent critical point for $(x,s)$,  and \(z_c\) the associated ascent critical point for $(y,t)$. 
If \(w_c\leq z_c\), namely in the interlacing case $\cI((x,s),(y,t))=1$, we may deform the local descent path \(\sfD^{\rm d}(w_0)\) to a contour
passing through \(w_c\), and the local ascent path
\(\sfD^{\rm a}(w_0)\) to a contour passing through \(z_c\). Since
\(w_c\leq z_c\), these contours may be chosen to be disjoint when
\(w_c<z_c\) and to intersect only at their common critical point when
\(w_c=z_c\). Using \Cref{p:invariance_contour}, we may then further deform them, separately in the upper and lower half-planes, to the corresponding steepest descent $\sfD^{\rm d}(w_c)$ and ascent contours $\sfD^{\rm a}(z_c)$, introduced in \Cref{s:path_analysis}.
This yields \eqref{e:sumoftwo2} with
\begin{align}\label{e:wJ1_case1}
\wt J^{(0)}+\wh J^{(1)}
:=\cI((x,s),(y,t))J_\fT-\sum_{\xi}\sgn(\xi)\frac{n}{2\pi \ri}
\int_{\sfD(\xi)}
P_{ns}(nz,nx)\,Q_{nt}(nz,ny)\,\rd z.
\end{align}

\begin{figure}
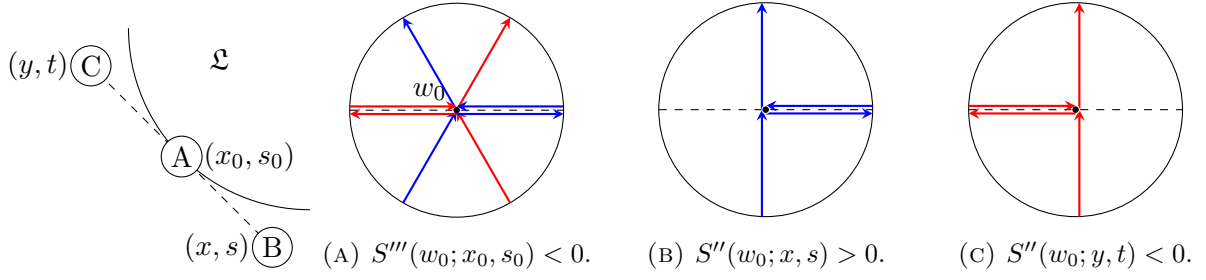

	\begin{subfigure}{0.24\textwidth}
	
		\centering

		\end{center}
		\caption{$S'''(w_0;x_0,s_0)<0$.}
		\end{subfigure}	
		\begin{subfigure}{0.24\textwidth}
			\begin{center}		
			%
			\end{center}
			\caption{$S''(w_0;x,s)>0$.}
			\end{subfigure}
			\begin{subfigure}{0.24\textwidth}
			\begin{center}		
			%
				
			\end{center}
			\caption{$S''(w_0;y,t)<0$.}
			\end{subfigure}

				\caption{\label{f:deform_frozen_path1}  
Deformed local path associated with a regular arctic point.
}
\end{figure}

We now treat all remaining cases: either at least one of \((x,s)\) and \((y,t)\) belongs to \(\fL\), or else \((x,s),(y,t)\notin \fL\) and \(w_c> z_c\). 
In these situations, we augment the local descent path $\sfD^{\rm d}(w_0)$ by adding two short segments: one from \(w_0+\fc+0\ri\) to \(w_0+0\ri\) in the upper half-plane, and one from \(w_0-0\ri\) to \(w_0+\fc-0\ri\) in the lower half-plane, as in  \Cref{f:deform_frozen_path1}. 
Thus the amended local descent contour consists of one component in the upper half-plane and one in the lower half-plane. 
We amend the local ascent contour in the same way.

With this convention, the upper-half-plane descent and ascent contours intersect at \(w_0+0\ri\) with negative orientation, while the corresponding lower-half-plane contours intersect at \(w_0-0\ri\) with positive orientation. 
Viewed in this way, the sum of the single-integral term and the double-integral term in \eqref{e:sumoftwo} is exactly of the form covered by \Cref{p:invariance_contour}. 
We may therefore deform \(\sfC^{\rm d}(w_0)\) and \(\sfC^{\rm a}(w_0)\) to the steepest descent and ascent contours $\sfD^{\rm d}(w_c)$ and $\sfD^{\rm a}(z_c)$.

Consequently, if \((x,s)\in\fL\) or \((y,t)\in\fL\), or if \((x,s),(y,t)\notin\fL\) and \(w_c> z_c\), then \eqref{e:sumoftwo2} holds with
\begin{align}\label{e:tJ1new}
\wh J^{(0)}=0,\quad \wh J^{(1)}
:=
-\sum_{\xi}\sgn(\xi)\frac{n}{2\pi \ri}
\int_{\sfD(\xi)}
P_{ns}(nz,nx)\,Q_{nt}(nz,ny)\,\rd z,
\end{align}
This matches with \eqref{e:sumoftwo2} by noticing $\cI((x,s),(y,t))=0$.

\noindent\textbf{Case 2: \(w_0\) corresponds to  a cusp location \((x_0,s_0)\); see panel (B) of \Cref{f:l}. }

The chart associated with \((y,t)\) is either a cusp chart, or a cusp frozen chart satisfying
$
S''(w_0;\cdot)<0$,
and in either case carries local ascent paths. 
Similarly, the chart associated with \((x,s)\) is either a cusp chart, or a cusp frozen chart satisfying
$
S''(w_0;\cdot)>0$, and carries a local descent path.

Without loss of generality, we assume that at the cusp location $(x_0,s_0)\in \fT$ for some curvilinear triangle $\fT$ and satisfies
\[
S''''(w_0;x_0,s_0)<0.
\]
We recall from \Cref{f:frozen_path2} all possible local configurations of the relevant paths. 
After a slight perturbation, the local descent path (blue) lies between the two local ascent paths (red).

Assume first that \((x,s),(y,t)\notin \fL\), and let \(w_c\) and \(z_{c,1}\leq z_{c,2}\) denote the associated frozen critical points. 
If
\[
z_{c,1}\leq w_c\leq z_{c,2},
\]
that is, in the interlacing case $\cI((x,s),(y,t))=1$, then we may deform the local descent path
\(\sfC^{\rm d}(w_0)\) to a contour passing through \(w_c\), and the local
ascent path \(\sfC^{\rm a}(w_0)\) to a pair of contours passing through
\(z_{c,1}\) and \(z_{c,2}\), respectively. Since
\(z_{c,1}\leq w_c\leq z_{c,2}\), these contours may be chosen so that the
descent contour intersects the ascent contours only at a common critical
point, which can occur only when \(w_c=z_{c,1}\) or \(w_c=z_{c,2}\).
Using \Cref{p:invariance_contour}, we may then further deform them, separately in the upper and lower half-planes, to the corresponding steepest descent contour \(\sfD^{\rm d}(w_c)\) and ascent contours \(\sfD^{\rm a}(z_{c,1})\cup \sfD^{\rm a}(z_{c,2})\). 
This yields \eqref{e:sumoftwo2}, with \(\wh J^{(0)}, \wh J^{(1)}\) given by \eqref{e:wJ1_case1}.

\begin{figure}
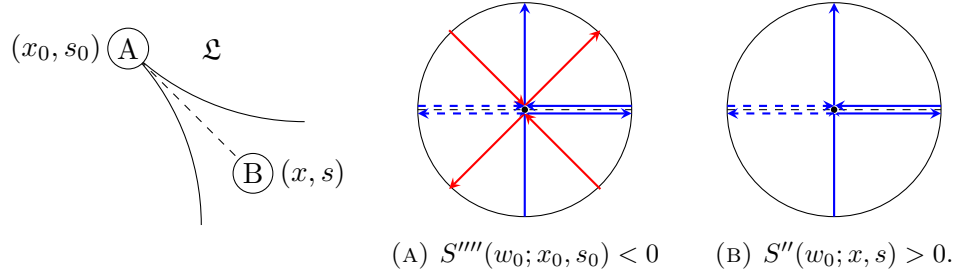


  \begin{subfigure}{0.3\textwidth}
    \begin{center}

			\end{center}\caption{$S''''(w_0;x_0, s_0)<0$}
			\end{subfigure}
			\begin{subfigure}{0.24\textwidth}
			\begin{center}		
			%
			\end{center}\caption{$S''(w_0;x,s)>0$.}
			\end{subfigure}  
  
	 \caption{\label{f:deform_frozen_path2}  
Deformed local path associated with a cusp point.
}
	 \end{figure}

We now treat all remaining cases: either at least one of \((x,s)\) and \((y,t)\) belongs to \(\fL\), or else \((x,s),(y,t)\notin \fL\) and \(w_c> z_{c,2}\) or \(w_c< z_{c,1}\). 
In these situations, we regard the local ascent contour as consisting of one component in the upper half-plane and one component in the lower half-plane. 
We then augment the local descent path by adding two short segments: either
\[
[w_0+\fc+0\ri,\,w_0+0\ri]
\quad\text{and}\quad
[w_0-0\ri,\,w_0+\fc-0\ri],
\]
or
\[
[w_0-\fc+0\ri,\,w_0+0\ri]
\quad\text{and}\quad
[w_0-0\ri,\,w_0-\fc-0\ri],
\]
as in \Cref{f:deform_frozen_path2}. 
Thus the amended local descent contour consists of one component in the upper half-plane and one component in the lower half-plane.

With this convention, the upper-half-plane descent and ascent contours intersect at \(w_0+0\ri\) with negative orientation, while the corresponding lower-half-plane contours intersect at \(w_0-0\ri\) with positive orientation. 
Viewed in this way, the sum of the single-integral term and the double-integral term in \eqref{e:sumoftwo} is exactly of the form covered by \Cref{p:invariance_contour}. 
We may therefore deform \(\sfC^{\rm d}(w_0)\) and \(\sfC^{\rm a}(w_0)\) to the steepest descent and ascent contours introduced in \Cref{s:path_analysis}. 
Hence \eqref{e:sumoftwo2} holds with \(\wh J^{(0)}, \wh J^{(1)}\) given by \eqref{e:tJ1new}.

\begin{figure}
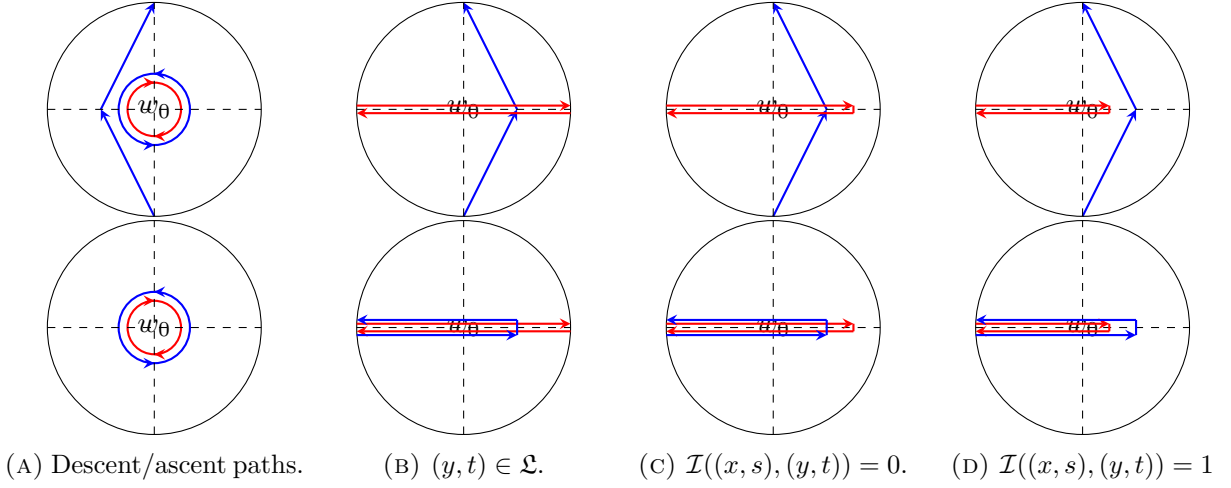
						
		\begin{subfigure}[t]{0.24\textwidth}
			\centering
			%
		\caption{Descent/ascent paths.}

	\end{subfigure}
	\begin{subfigure}[t]{0.24\textwidth}
			\centering
			%
		\caption{$(y,t)\in \fL$.}

	\end{subfigure}
\begin{subfigure}[t]{0.24\textwidth}
			\centering
			%
		\caption{$\cI((x,s),(y,t))=0$.}

	\end{subfigure}
	\begin{subfigure}[t]{0.24\textwidth}
			\centering
			%
		\caption{$\cI((x,s),(y,t))=1$}

	\end{subfigure}

%	\begin{subfigure}[t]{0.24\textwidth}
%			\centering
%			%
%		\caption{$z_c>w_c$}
%
%	\end{subfigure}

	\caption{
\label{f:concentric_tangent1}
Local contour configurations for concentric charts centered at $w_0$ from Panel (B) of \Cref{f:vertical_tangent1}. First row: $\fN_\al$ is a tangent neighborhood; Second row: $\fN_\al$ is a frozen neighborhood.}
	\end{figure}

\begin{figure}
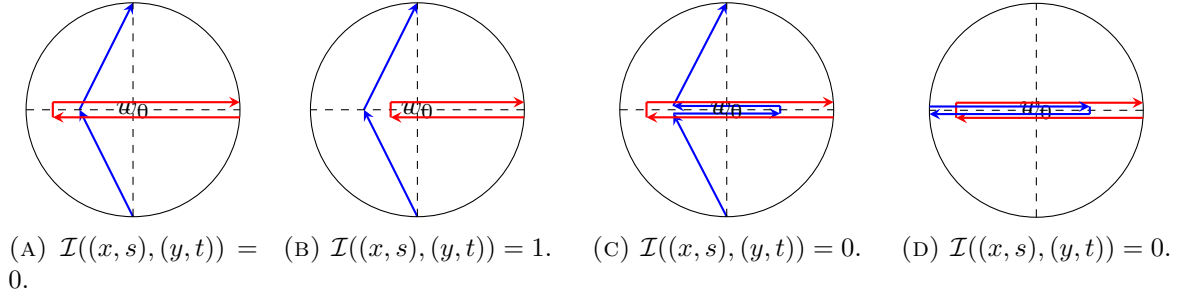
						
		\begin{subfigure}[t]{0.2\textwidth}
			\centering
			%
			\caption{$\cI((x,s),(y,t))=0$.}

	\end{subfigure}
	\begin{subfigure}[t]{0.24\textwidth}
			\centering
			%
		\caption{$\cI((x,s),(y,t))=1$.}

	\end{subfigure}
	\begin{subfigure}[t]{0.24\textwidth}
			\centering
			%
		\caption{$\cI((x,s),(y,t))=0$.}
	\end{subfigure}
	
%	\begin{subfigure}[t]{0.24\textwidth}
%			\centering
%			%
%		\caption{$z_c>w_c$}
%
%	\end{subfigure}

	\caption{
\label{f:concentric_tangent2}
Local contour configurations for concentric charts centered at $w_0$ from Panel (A) of \Cref{f:vertical_tangent1}.}
	\end{figure}

\noindent\textbf{Case 3: \(w_0\) corresponds to  a tangent location $(y_0,t_0)$.} 

Without loss of generality, assume that the tangent line is vertical and
separates the curvilinear triangles \(\fT_A\) and \(\fT_B=\fT\), as in
Panel (C) of \Cref{f:l}. In this case, the local configuration at the
tangency location \((y_0,t_0)\) is as in
\Cref{f:vertical_tangent1}.

We recall the possible types of the neighborhoods \(\fN_\al\) and
\(\fN_{(y,t)}\) from \eqref{e:configuration}. Each of
\(\fN_{(y,t)}\) and \(\fN_\al\) is of one of the following four types: a
tangent neighborhood, an interface frozen neighborhood, a neighborhood
contained in \(\fT_A\cup\fL\), or a neighborhood contained in
\(\fT_B\cup\fL\).

If \(\fN_{(y,t)}\) is contained in \(\fT_A\cup\fL\), then it lies below
the tangency location \((y_0,t_0)\). If \(\fN_{(y,t)}\) is contained in
\(\fT_B\cup\fL\), then it lies above \((y_0,t_0)\). If \(\fN_\al\) is
contained in \(\fT_A\cup\fL\), then it lies above \((y_0,t_0)\), which is
impossible in the configuration of \Cref{f:vertical_tangent1}. If
\(\fN_\al\) is contained in \(\fT_B\cup\fL\), then it lies below
\((y_0,t_0)\).

In summary, the neighborhood \(\fN_{(y,t)}\) is of one of the following
types:
\begin{align}\label{e:configurationY}
&\text{a tangent neighborhood, an interface frozen neighborhood,} \notag\\
&\text{a neighborhood contained in \(\fT_B\cup\fL\), or a neighborhood
contained in \(\fT_A\cup\fL\).}
\end{align}
The neighborhood \(\fN_\al\) is of one of the following types:
\begin{align}\label{e:configurationX}
&\text{a tangent neighborhood, an interface frozen neighborhood, or a
neighborhood contained in \(\fT_B\cup\fL\).}
\end{align}

The main difference from the regular arctic and cusp cases is that, in the
present tangent case, the local ascent and descent paths
$
        \sfD^{\rm a}(w_0),
$ and $
        \sfD^{\rm d}(w_0)$
may deform to the empty set. We first recall when this happens.

By \Cref{c:PIproperty} and \Cref{l:holomorphic}, if
$
        y\geq w_0,
$
then the integrand is holomorphic in \(z\) in a small neighborhood of \(w_0\),
and the (red) circular ascent path  deforms to \(\emptyset\). Similarly, if
$
        x\leq w_0,
$
then the integrand is holomorphic in \(w\) in a small neighborhood of \(w_0\),
and the (blue) local descent path deforms to \(\emptyset\).

Hence, by \Cref{l:vertical_tangent_steepest} and
\Cref{l:local_descent_deformation}, if 
\begin{align}
y\geq w_0,
\end{align}
 then the tangent or
frozen chart centered at \(w_0\), associated with \((y,t)\), contains no
ascent critical point, and
$
        \sfC^{\rm a}(w_0)
$
deforms to \(\emptyset\). Similarly
 if \(\fN_\al\) is a frozen neighborhood and
\[
x\leq w_0,
\]
 then the frozen chart centered at \(w_0\), associated with
\((x,s)\), contains no descent critical point, and
$
        \sfC^{\rm d}(w_0)
$
deforms to \(\emptyset\).

In all other cases, \(\sfC^{\rm d}(w_0)\) deforms to a union of steepest
descent paths \(\sfD^{\rm d}(w_c)\), passing through descent critical points
\(w_c\), and \(\sfC^{\rm a}(w_0)\) deforms to a union of steepest ascent paths
\(\sfD^{\rm a}(z_c)\), passing through ascent critical points \(z_c\).

\medskip
\noindent
\emph{Case {\rm i}: both local ascent and descent paths deform to nonempty
steepest ascent/descent contours.}

We first classify the possible descent configurations associated with
\((x,s)\).

\begin{enumerate}
\item\label{i:xscase1}
Suppose that \(\fN_\al\) is a tangent neighborhood and
$
        (x,s)\in\fL.
$
Then \(x\leq w_0\), there is a pair of complex conjugate critical points $w_{c}$ and $\overline{w_c}$. The local descent path
\(\sfC^{\rm d}(w_0)\) deforms to
$
        \sfD^{\rm d}(w_{c})\cup \sfD^{\rm d}(\overline{w_{c}}),
$
and the two components meet at a point
$
        \xi\in [x,y_0];
$
see Panel (F) of \Cref{f:tangent2}.

\item\label{i:xscase2}
Suppose that \(\fN_\al\) is a tangent neighborhood,
$
        x\leq w_0,
$ and $
        (x,s)\in\fT_A.
$
Then the local descent path \(\sfC^{\rm d}(w_0)\) deforms to
$\sfD^{\rm d}(w_c)$,
with $
        w_c<w_0.
$

\item\label{i:xscase3}
In the remaining cases, either \(\fN_\al\) is a tangent neighborhood and
\((x,s)\in\fT_B\), or \(\fN_\al\) is not a tangent neighborhood and
\(x>w_0\).
In the latter case, by \eqref{e:configurationX}, either $\fN_\al$ is an interface frozen neighborhood, or a  neighborhood in  \( \fT_B\cup \fL\). 
The local descent path
\(\sfC^{\rm d}(w_0)\) deforms to
$
        \sfD^{\rm d}(w_c)$,
with $
        w_c>w_0.
$
\end{enumerate}

We now discuss the possible positions of \((y,t)\).

\begin{enumerate}
\item
Suppose that \(\fN_{(y,t)}\) is a tangent neighborhood and
$
        (y,t)\in\fL.
$
Then  there is a pair of complex conjugate critical points $z_{c}$ and $\overline{z_c}$, and  the local ascent path \(\sfC^{\rm a}(w_0)\) deforms to
$
        \sfD^{\rm a}(z_{c})\cup \sfD^{\rm a}(\overline{z_{c}}),
$ and we are in the non-interlacing case, so
$
        \cI((x,s),(y,t))=0.
$

We use the contours in Panel (B) of \Cref{f:vertical_tangent1} for the possible
local configurations of the descent and ascent contours. Then by \Cref{l:boundary_piece_in_C}, $\fN_\al \cap \fC(\fT_B;\fN_{(y,t)})\neq \emptyset$ and 
\begin{align}\label{e:J1}
J^{(1)}=\bm1(\fN_\al \cap \fC(\fT_B;\fN_{(y,t)})\neq \emptyset)J_B=J_B.
\end{align}

We modify the local ascent and descent paths as in Panel (B) of
\Cref{f:concentric_tangent1}, so that the blue descent path passes through
the descent critical point $w_c$. In this form, the sum of the single-contour integral term
and the double-contour integral term in \eqref{e:sumoftwo} is exactly of the type
covered by \Cref{p:invariance_contour}.

\item
Suppose either that \(\fN_{(y,t)}\) is a tangent neighborhood with
$
(y,t)\in\fT_B,
$
or that \(\fN_{(y,t)}\) is not a tangent neighborhood and
\[
y<w_0,
\qquad
(y,t)\in\fT_B\cup\fL.
\]
In these cases, \((y,t)\) is associated with an ascent critical point
$
        z_c>w_0.
$
We use the contours in Panel (B) of \Cref{f:vertical_tangent1} for the possible
local configurations of the descent and ascent contours, and \eqref{e:J1} holds.

In the cases \Cref{i:xscase1}, \Cref{i:xscase2}, or \Cref{i:xscase3} with
$
        w_c<z_c,
$
we are in the non-interlacing case, so
$
        \cI((x,s),(y,t))=0.
$
We modify the local ascent and descent paths as in Panel (C) of
\Cref{f:concentric_tangent1}, so that the blue descent path passes through
\(\xi\) in \Cref{i:xscase1}, or through \(w_c\) in \Cref{i:xscase2} and
\Cref{i:xscase3}, while the red ascent path passes through \(z_c>w_0\). In
this form, the sum of the single-contour integral term and the double-contour integral term in
\eqref{e:sumoftwo} is exactly of the type covered by
\Cref{p:invariance_contour}.

In \Cref{i:xscase3} with
$
        w_c\geq z_c,
$
we are in the interlacing case, so
$
        \cI((x,s),(y,t))=1.
$
We may deform the local descent path to a contour passing through \(w_c\), and
the local ascent path to a contour passing through \(z_c\); see Panel (D) of
\Cref{f:concentric_tangent1}. Since
\(z_c\leq w_c\), these contours can be chosen to intersect only at their common critical point when
\(w_c=z_c\).

\item
Suppose that either \(\fN_{(y,t)}\) is a tangent neighborhood with
$
        (y,t)\in\fT_A,
$
or that \(\fN_{(y,t)}\) is not a tangent neighborhood and 
\begin{align}
        y<w_0, \qquad 
        (y,t)\in\fT_A\cup\fL.
\end{align}
In these cases, \((y,t)\) is associated with an ascent critical point
$
        z_c<y_0.
$

If $\fN_{\al}$ is a neighborhood in $\fT_B\cup \fL$ and $x>w_0$ as in \Cref{i:xscase3}, we use the contours in Panel (B) of \Cref{f:vertical_tangent1} for the possible local configurations of the descent and ascent contours. Since $y\leq w_0<x$, the ordering condition \eqref{e:JB_nonvanish} gives,
\begin{align}
J^{(1)}=J_B=0.
\end{align}
Moreover, this is the noninterlacing case, and
\(\cI((x,s);(y,t))=0\). We may deform the local descent path to a contour
passing through \(w_c\), and the local ascent path to a contour passing
through \(z_c\); see Panel (D) of \Cref{f:concentric_tangent1}.

In the remaining cases, we use the contours in Panel (A) of \Cref{f:vertical_tangent1} for the possible
local configurations of the descent and ascent contours. Then by \Cref{l:boundary_piece_in_C}, $\fN_\al \cap \fC(\fT_A;\fN_{(y,t)})\neq \emptyset$ and 
\begin{align}\label{e:J1JA}
J^{(1)}=\bm1(\fN_\al \cap \fC(\fT_A;\fN_{(y,t)})\neq \emptyset)J_A=J_A.
\end{align}

We now split according to the three possibilities
\Cref{i:xscase1}--\Cref{i:xscase3}.

\smallskip
\noindent
\emph{Subcase \Cref{i:xscase1}.} We are in the non-interlacing case, so
$
        \cI((x,s),(y,t))=0.
$
Since \(x\leq w_0\), the integrand is holomorphic in \(w\) in a small
neighborhood of \(y_0\), and the blue circular descent contour deforms to
\(\emptyset\). There are two possibilities.

\begin{enumerate}
\item
If \(x\geq y\), then
$
        z_c<y\leq x.
$
We modify the local ascent and descent paths as in Panel (A) of
\Cref{f:concentric_tangent2}, so that the blue descent path passes through
\(\zeta\in[x,w_0]\), and the red ascent path passes through \(z_c\) with $z_c<\zeta$. In this
form, the sum of the single-integral term and the double-integral term in
\eqref{e:sumoftwo} is exactly of the type covered by
\Cref{p:invariance_contour}.

\item
If \(x<y\), then the ordering condition \Cref{e:JA_nonvanish} gives
$
        J^{(1)}=J_A=0.
$
If \(\zeta\geq z_c\), we modify the local ascent and descent paths as in Panel
(A) of \Cref{f:concentric_tangent2}. If \(\zeta\leq z_c\), we modify them as in
Panel (B). In either case, the blue descent path passes through \(\zeta\), and
the red ascent path passes through \(z_c\).
\end{enumerate}

\smallskip
\noindent
\emph{Subcase \Cref{i:xscase2}.}
Again \(x\leq w_0\), so the integrand is holomorphic in \(w\) in a small
neighborhood of \(y_0\), and the circular descent contour deforms to
\(\emptyset\). There are two possibilities.

\begin{enumerate}
\item
If
$
        w_c>z_c,
$
then we are in the non-interlacing case,
$
        \cI((x,s),(y,t))=0.
$
We modify the local ascent and descent paths as in Panel (A) of
\Cref{f:concentric_tangent2}. The blue descent path passes through \(w_c\),
and the red ascent path passes through \(z_c\).

\item
If
$
        w_c\leq z_c\leq y_0,
$
then we are in the interlacing case,
$
        \cI((x,s),(y,t))=1.
$
We may deform the local descent path to a contour passing through \(w_c\), and
the local ascent path to a contour passing through \(z_c\); see Panel (B) of
\Cref{f:concentric_tangent2}. Since \(w_c\leq z_c\), these contours can be
chosen to remain disjoint.
\end{enumerate}

\smallskip
\noindent
\emph{Subcase \Cref{i:xscase3}.}
In this case,
$
        \cI((x,s),(y,t))=0.
$
We modify the local ascent and descent paths as in Panel (C) of
\Cref{f:concentric_tangent2} when \(\fN_\al\) is a tangent neighborhood, and
as in Panel (D) when \(\fN_\al\) is an interface frozen neighborhood. The blue descent
path passes through \(w_c>w_0\), while the red ascent path passes through
\(z_c<w_0\).
\end{enumerate}
%\end{enumerate}

In all cases discussed above, \Cref{p:invariance_contour} allows us to deform
the local contours, separately in the upper and lower half-planes, to the
corresponding steepest descent contour \(\sfD^{\rm d}(w_c)\) and steepest
ascent contour \(\sfD^{\rm a}(z_c)\). This gives \eqref{e:sumoftwo2}, with
\(\wh J^{(0)}\) and \(\wh J^{(1)}\) given by \eqref{e:wJ1_case1}.

\noindent\emph{Case {\rm ii}: one of \(\sfD^{\rm d}(w_0)\) or
\(\sfD^{\rm a}(w_0)\) deforms to \(\emptyset\).}

We need to show that
\begin{align}\label{e:J1vanish}
\eqref{e:sumoftwo}
=
\cI((x,s),(y,t))J_\fT
+
\OO\left(
\sum_{w_c,z_c}
e^{n\Re[S(w_c;x,s)-S(z_c;y,t)]-\fc'n}
\right).
\end{align}

We recall the possibilities of $\fN_{(y,t)}$ and $\fN_\al$ from \eqref{e:configurationY} and \eqref{e:configurationX}. 

We start with the case
$
        y\geq w_0.
$
Then the tangent or frozen chart centered at \(w_0\), associated with
\((y,t)\), contains no ascent critical point, and
$
        \sfD^{\rm a}(w_0)
$
deforms to \(\emptyset\).

\begin{enumerate}
\item
Suppose first that \(y\geq w_0\) and
$
(y,t)\in\fT_B\cup\fL.
$
Then \((y,t)\) is also associated with a frozen chart centered at \(z_0'>w_0\),
which contains an ascent critical point \(z_c\); see Panel (A) of
\Cref{f:tangent_t}. Recall from \eqref{e:configurationX} that
\(\fN_\al\) is either a tangent neighborhood, an interface frozen
neighborhood, or a neighborhood contained in
\(\fT_B\cup\fL\). In all cases, \Cref{l:boundary_piece_in_C} gives
$
\fN_\al\cap\fC(\fT_B;\fN_{(y,t)})\neq\emptyset.
$
We use the contours in Panel (B) of \Cref{f:vertical_tangent1} for the
local descent and ascent contours. Then 
\begin{align}
J^{(1)}=  \bm1\!\left(
        \fN_\al\cap\fC(\fT_B;\fN_{(y,t)})\neq\emptyset
        \right)J_B=J_B
\end{align}
 Since the local ascent paths can be
deformed to the empty contour, the double-contour integral vanishes. We
now distinguish three possibilities.

\begin{enumerate}
\item
If \(x\leq w_0\), then
$
        x\leq w_0\leq y.
$
Hence, by \eqref{e:JB_nonvanish},
$
        J_B=0.
$
Therefore
\[
        0
        =
        J^{(1)}=
        J_B
        =
        \cI((x,s),(y,t))J_B.
\]

\item
If \(x>w_0\) and \((x,s)\in\fT_B\), then, since the region
\(\fC(\fT_B;y,t)\) is bounded by the vertical tangent line \(L_\infty\) and
the cut \(\ell_-(z_c;\fT_B)\), the interlacing condition holds:
\[
(x,s)\in\fC(\fT_B;\fN_{(y,t)}).
\]
Thus,
\(\cI((x,s);(y,t))=1\), and
\[
J^{(1)}
=
J_B
=
\cI((x,s);(y,t))J_B.
\]

\item
If \(x>w_0\) and \((x,s)\in\fL\), then \Cref{p:J1vanish} gives
\[
        J^{(1)}
        =
        J_B
        =\OO\left( \sum_{w_c, z_c}
        e^{n\Re[S(w_c;x,s)-S(z_c;y,t)]-\fc'n}\right).
\]
\end{enumerate}

\item \label{i:case2}
Suppose next that \(\fN_{(y,t)}\) is a neighborhood in
$
       \fT_A\cup\fL.
$
Then \(\fN_{(y,t)}\) lies below the tangency location \((y_0,t_0)\); see
Panel (B) of \Cref{f:tangent_t}. We use the contours in Panel (A) of
\Cref{f:vertical_tangent1} for the local descent and ascent contours.

By \Cref{p:concentric_setting}, we have
$
\fN_\al\cap\fT_A\neq\emptyset.
$
Together with \eqref{e:configurationX}, this implies that \(\fN_\al\) is
either a tangent neighborhood or an interface frozen neighborhood. Then, by
\Cref{l:boundary_piece_in_C}, we have
$
\fN_\al\cap\fC(\fT_A;\fN_{(y,t)})\neq\emptyset,
$
and hence
\begin{align}
J^{(1)}
=
\bm1\bigl(
\fN_\al\cap\fC(\fT_A;\fN_{(y,t)})\neq\emptyset
\bigr)J_A
=
J_A.
\end{align}

In this case \(\fN_\al\) lies above \(\fN_{(y,t)}\). In particular,
$
        s\geq t,
$
and hence, by \eqref{e:JB_nonvanish},
$
        J_B=0.
$
By \Cref{c:change_triangle},
\begin{align}\begin{split}\label{e:sumzero}
&\phantom{{}={}}
J_A
+
(\text{double-contour integral relative to \(\fT_A\)}) \\
&=
J_B
+
(\text{double-contour integral relative to \(\fT_B\)})
=
0.
\end{split}\end{align}
Here the last equality follows because \(J_B=0\), and, if we use the contours
in Panel (B) of \Cref{f:vertical_tangent1}, the local ascent paths deform to
\(\emptyset\), so the local double-contour integral vanishes.

Moreover, in this configuration the vertical tangent line \(L_\infty\) is a
boundary cut of \(\fC(\fT_A;y,t)\). Hence the interlacing condition fails:
$
        (x,s)\notin\fC(\fT_A;y,t),
$, and we have $
        \cI((x,s),(y,t))=0.
$
This gives \eqref{e:J1vanish}.
\end{enumerate}

It remains to discuss the case in which \(\fN_\al\) is not a tangent neighborhood
and
$
        x\leq w_0.
$
Then the frozen chart centered at \(w_0\), associated with \((x,s)\), contains
no descent critical point, and
$
        \sfD^{\rm d}(w_0)
$
deforms to \(\emptyset\). We may further assume that
$
        y<w_0,
$
since the complementary case \(y\geq w_0\) was already treated.

\begin{enumerate}
\item
Suppose first that
$
        x\leq w_0,
$ and $
        \fN_\al$ is a neighborhood contained in $\fT_B\cup\fL.
$
Then \(\fN_\al\) lies below the tangency location; see Panel (C) of
\Cref{f:tangent_t}. 

By \Cref{p:concentric_setting},
$
        \fN_{(y,t)}\cap\fT_B\neq\emptyset.
$
We use the contours in Panel (B) of
\Cref{f:vertical_tangent1} for the local descent and ascent contours.
Then by \Cref{l:boundary_piece_in_C}, we have $ \fN_\al\cap\fC(\fT_B;\fN_{(y,t)})\neq\emptyset$, and 
\begin{align}
       \quad J^{(1)}=\bm1(\fN_\al\cap\fC(\fT_B;\fN_{(y,t)}))J_B=J_B.
\end{align}

In this case, \((y,t)\in\fN_{(y,t)}\) lies above \(\fN_\al\). In particular,
$
        t>s,
$
and hence, by \eqref{e:JA_nonvanish},
$
        J_A=0.
$
By the same argument as in \eqref{e:sumzero}, we obtain
\[
        J_B
        +
        (\text{double-contour integral relative to \(\fT_B\)})
        =
        0.
\]
Moreover, either \((y,t)\notin\fT_B\), or \((y,t)\in\fT_B\) and $y<w_0$.  In the
latter case,
$
        (x,s)\notin\fC(\fT_B;y,t).
$
Thus, in both cases, the interlacing condition fails:
$
        \cI((x,s),(y,t))=0.
$
This gives \eqref{e:J1vanish}.

\item
Suppose next that \((x,s)\in\fN_\al\) is an interface frozen neighborhood.
Since \(x\leq w_0\), we must have
$
        (x,s)\in\fT_A.
$
In this case, \(\fN_\al\) is associated with another frozen chart centered at
\(w_0'\). This chart contains a descent critical point \(w_c\), corresponding
to a tangent line from \((x,s)\) to the portion of the arctic boundary
contained in \(\fT_A\).

If \(\fN_{(y,t)}\) is a neighborhood contained in \(\fT_B\cup\fL\),
then this case has already been treated above in \Cref{i:case2}. We therefore assume that
\(\fN_{(y,t)}\) is either a tangent or interface frozen neighborhood, or a
 neighborhood contained in \(\fT_A\cup\fL\).

In all remaining cases, by \Cref{l:boundary_piece_in_C}, 
we have $ \fN_\al\cap\fC(\fT_A;\fN_{(y,t)})\neq\emptyset$, and 
\begin{align}
       \quad J^{(1)}=\bm1(\fN_\al\cap\fC(\fT_A;\fN_{(y,t)}))J_A=J_A.
\end{align}
we use the contours in Panel (A) of
\Cref{f:vertical_tangent1} for the local descent and ascent contours. 
Then the
local descent paths deform to \(\emptyset\), and the double-contour integral
vanishes. 

There are three subcases.

\begin{enumerate}
\item
If \((y,t)\in\fL\), then
$
        \cI((x,s),(y,t))=0.
$
By \Cref{p:J1vanish},
\[
        J^{(1)}
        =
        J_A
        =\OO\left( \sum_{w_c, z_c}
        e^{n\Re[S(w_c;x,s)-S(z_c;y,t)]-\fc'n}\right).
\]
\item
If \(\fN_{(y,t)}\) is a tangent neighborhood and
$
        (y,t)\in\fT_B,
$
then
$
        \cI((x,s),(y,t))=0.
$
Moreover, in this case \(t>s\), and therefore, by \eqref{e:JA_nonvanish},
$
        J_A=0.
$

\item
If
$
        (y,t)\in\fT_A,
$ and $
        y<w_0,
$
then \((y,t)\) is associated with an ascent critical point \(z_c\) in the
chart centered at \(w_0\). The charts centered at \(w_0'\) and \(w_0\) are
disjoint, so
$
        w_c<z_c.
$
Thus we are in the interlacing regime:
$
        (x,s)\in\fC(\fT_A;y,t).
$
Therefore $  \cI((x,s),(y,t))=1$ and 
\[
        J^{(1)}
        =
       J_A
        =
        \cI((x,s),(y,t))J_A.
\]
\end{enumerate}
\end{enumerate}

\begin{figure}
	 \begin{subfigure}{0.32\textwidth}

    \caption{}
  \end{subfigure}
  	 \begin{subfigure}{0.32\textwidth}
      %
\caption{}
  \end{subfigure}
    	 \begin{subfigure}{0.32\textwidth}
      %
\caption{}
  \end{subfigure}

  	\caption{
	Left panel: local geometry;	Middle panel: $\fN_{(y,t)}\in \fT_A\cup \fL$; Right panel $\fN_{\al}\in\fT_B\cup \fT$.
}\label{f:tangent_t}
\end{figure}

\noindent\textbf{Case 4: \(w_0\) corresponds to a cusp-turning location
\((y_0,t_0)\).}

Without loss of generality, assume that the tangent line is vertical and
separates \(\fT_A\) and \(\fT_B=\fT\); see Panel (D) of \Cref{f:l}.

Then \(\fN_{(y,t)}\) is either a cusp-turning neighborhood or a tangent frozen
neighborhood; in the latter case, the associated tangent frozen chart carries
local ascent paths. Similarly, \(\fN_\al\) is either a cusp-turning
neighborhood or a tangent frozen neighborhood; in the latter case, the
associated tangent frozen chart carries local descent paths.

Without loss of generality, assume that the local configuration at the
cusp-turning location \((y_0,t_0)\) is as in \Cref{f:c_vertical_cusp2}. We
discuss the case in which both \(\fN_\al\) and \(\fN_{(y,t)}\) are
cusp-turning neighborhoods. The remaining cases, namely cusp-turning--frozen,
frozen--cusp-turning, and frozen--frozen, are treated in the same way as in
the tangent case above, so we omit them. If both \(\fN_\al\) and \(\fN_{(y,t)}\) are
cusp-turning neighborhoods, by \Cref{l:boundary_piece_in_C}, we have
\begin{align}
\fN_\al \cap \fC(\fT_A; \fN_{(y,t)})\neq \emptyset, \quad \fN_\al \cap \fC(\fT_B; \fN_{(y,t)})\neq \emptyset.
\end{align}

Recall from \Cref{c:vertical_cusp_steepest} that if
\((y,t)\in\fP\setminus\fL\), then \(\sfC^{\rm a}(w_0)\) can be deformed to a
steepest ascent path \(\sfD^{\rm a}(z_c)\) passing through a real ascent
critical point \(z_c\). Similarly, if \((x,s)\in\fP\setminus\fL\), then
\(\sfC^{\rm d}(w_0)\) can be deformed to a steepest descent path
\(\sfD^{\rm d}(w_c)\) passing through a real descent critical point \(w_c\).

We first consider the cases in which at least one of \((x,s)\) or \((y,t)\)
belongs to \(\fL\). In these cases we have
\begin{align}
\cI((x,s),(y,t))=0.
\end{align}

\begin{figure}
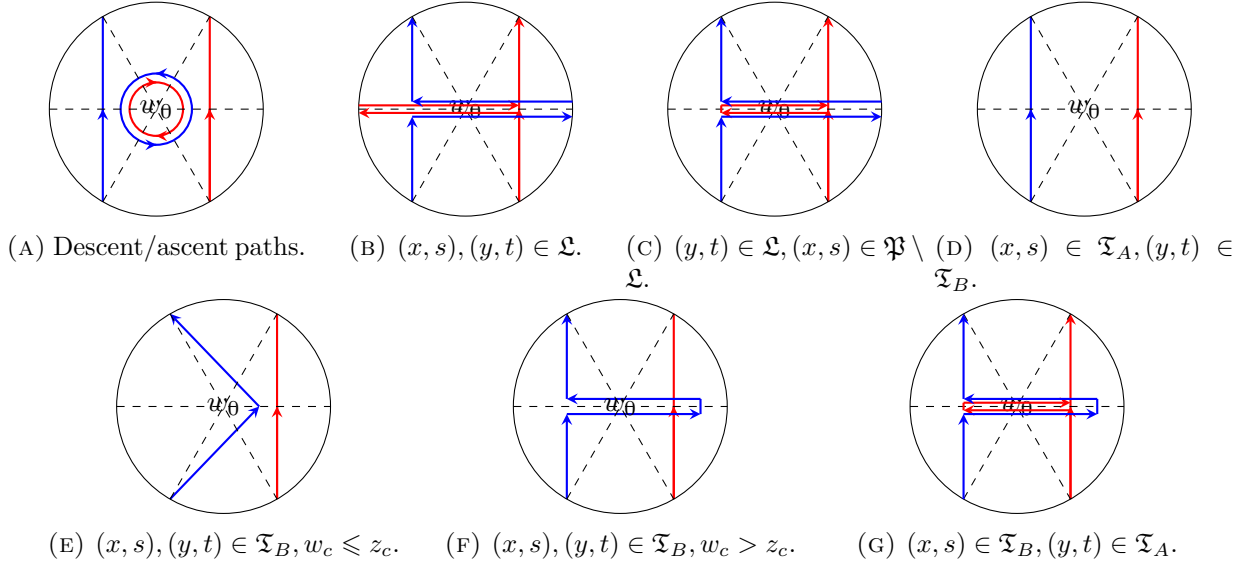
			
	\begin{subfigure}[t]{0.24\textwidth}
			\centering

		\caption{Descent/ascent paths.}

	\end{subfigure}
		\begin{subfigure}[t]{0.24\textwidth}
			\centering
			%
		\caption{ $(x,s),(y,t)\in \fL$.}
	
	\end{subfigure}
	\begin{subfigure}[t]{0.24\textwidth}
			\centering
			%
		\caption{$(y,t)\in \fL, (x,s)\in \fP\setminus \fL$.}
	\end{subfigure}
\begin{subfigure}[t]{0.24\textwidth}
			\centering
			%
		\caption{$(x,s)\in \fT_A, (y,t)\in \fT_B$.}

	\end{subfigure}

\begin{subfigure}[t]{0.31\textwidth}
			\centering
			%
		\caption{$(x,s),(y,t)\in \fT_B, w_c\leq z_c$.}

	\end{subfigure}
\begin{subfigure}[t]{0.31\textwidth}
			\centering
			%
		\caption{$(x,s),(y,t)\in \fT_B,w_c>z_c$.}

	\end{subfigure}
\begin{subfigure}[t]{0.31\textwidth}
			\centering
			%
		\caption{$(x,s)\in \fT_B,(y,t)\in \fT_A$.}
	
	\end{subfigure}

		\caption{
\label{f:concentric_cusp_turning}
Two concentric charts centered at $w_0$ which is a cusp-turning point.}
	\end{figure}

\begin{enumerate}
\item Suppose that \((x,s)\in\fL\) and \((y,t)\in\fL\). We use the contours in
Panel (B) of \Cref{f:c_vertical_cusp2} for the possible local configurations
of the descent and ascent contours. We then modify the local ascent and
descent paths as in Panel (B) of \Cref{f:concentric_cusp_turning}.

After this modification, the local descent and ascent contours may be regarded
as having one component in the upper half-plane and one component in the lower
half-plane. In this form, the sum of the single-contour integral term and the
double-contour integral term in \eqref{e:sumoftwo} is exactly of the type covered by
\Cref{p:invariance_contour}. 

\item Suppose that \((x,s)\in\fL\), but
\((y,t)\in\fP\setminus\fL\). We again use the contours in Panel (B) of
\Cref{f:c_vertical_cusp2} for the possible local configurations of the
descent and ascent contours. We then modify the local ascent and descent paths
as in Panel (C) of \Cref{f:concentric_cusp_turning}, where the red local
ascent path passes through the real ascent critical point \(z_c\).

As before, the sum of the single-integral term and the double-integral term in
\eqref{e:sumoftwo} is exactly of the type covered by
\Cref{p:invariance_contour}. 

\item Suppose that \((x,s)\in\fP\setminus\fL\) and \((y,t)\in\fL\). 
We can deform the contours in exactly the same way as the previous case. We omit the details.
\end{enumerate}
In all cases, we can further deform
\(\sfC^{\rm d}(w_0)\) and \(\sfC^{\rm a}(w_0)\) to the steepest descent and
ascent contours introduced in \Cref{c:vertical_cusp_steepest}. Hence by \Cref{p:invariance_contour}, \eqref{e:sumoftwo2} holds with
\(\wh J^{(0)}\) and \(\wh J^{(1)}\) given by \eqref{e:tJ1new}.

It remains to consider the case in which
$
        (x,s),(y,t)\in\fP\setminus\fL.
$
We split into four subcases according to the positions of \(x\) and \(y\)
relative to \(w_0\).

\begin{enumerate}
\item Suppose that
\[
        x\leq w_0,
        \qquad
        y\geq w_0.
\]
Then \((x,s)\in\fT_A\) and \((y,t)\in\fT_B\). We use the contours in Panel
(B) of \Cref{f:c_vertical_cusp2} for the possible local configurations of the
descent and ascent contours. In this case \(x\leq y\), so the ordering condition \eqref{e:JB_nonvanish} gives
$
        J^{(1)}=J_B=0.
$

In a small neighborhood of \(w_0\), the integrand is holomorphic in both \(z\)
and \(w\). Hence both the circular descent contour and the circular ascent
contour can be deformed to \(\emptyset\), leaving the configuration in Panel
(D) of \Cref{f:concentric_cusp_turning}. We can then further deform the
remaining descent and ascent contours to the steepest descent and ascent
contours introduced in \Cref{c:vertical_cusp_steepest}. These contours pass
through \(w_c\) and \(z_c\), respectively, with
\[
        w_c\leq w_0\leq z_c.
\]

\item Suppose that
\[
        x>w_0,
        \qquad
        y\geq w_0.
\]
Then \((x,s),(y,t)\in\fT_B\). We use the contours in Panel (B) of
\Cref{f:c_vertical_cusp2} for the possible local configurations of the
descent and ascent contours, and $J^{(1)}=J_B$.

In a small neighborhood of \(w_0\), the integrand is holomorphic in \(z\).
Hence the circular ascent contour can be deformed to \(\emptyset\). There are
two possibilities.

\begin{enumerate}
\item If
$
        z_c\geq w_c\geq w_0,
$
then we are in the interlacing case $\cI((x,s),(y,t))=1$, and
\[
        J^{(1)}=J_B=\cI((x,s),(y,t))J_B.
\]
In this case we deform the remaining local descent and ascent contours as in
Panel (E) of \Cref{f:concentric_cusp_turning}; they pass through \(w_c\) and \(z_c\),
respectively.

\item If
$
        w_c>z_c\geq w_0,
$
then we are in the non-interlacing case, and
$
        \cI((x,s),(y,t))=0.
$
In this case we deform the remaining local descent and ascent contours as in
Panel (F) of \Cref{f:concentric_cusp_turning}; they pass through \(w_c\) and \(z_c\),
respectively. In this form, the sum of the single-integral term and the
double-integral term in \eqref{e:sumoftwo} is exactly of the type covered by
\Cref{p:invariance_contour}.
\end{enumerate}

\item Suppose that
\[
        x\leq w_0,
        \qquad
        y<w_0.
\]
Then \((x,s),(y,t)\in\fT_A\). We use the contours in Panel (A) of
\Cref{f:c_vertical_cusp2} for the possible local configurations of the
descent and ascent contours. We can deform the contours in the same way as in the previous
case.
\item Suppose that
\[
        x>w_0,
        \qquad
        y<w_0.
\]
Then \((x,s)\in\fT_B\) and \((y,t)\in\fT_A\). We use the contours in Panel
(B) of \Cref{f:c_vertical_cusp2} for the possible local configurations of the
descent and ascent contours.

In this case we are in the non-interlacing regime $      w_c>w_0>z_c$, and 
$
        \cI((x,s),(y,t))=0.
$
We deform the local descent and ascent contours as in Panel (G) of
\Cref{f:c_vertical_cusp2}; they pass through \(w_c\) and \(z_c\),
respectively. In this form, the sum of the single-integral term and the
double-integral term in \eqref{e:sumoftwo} is exactly of the type covered by
\Cref{p:invariance_contour}.
\end{enumerate}

In all cases, we may further deform \(\sfC^{\rm d}(w_0)\) and
\(\sfC^{\rm a}(w_0)\) to the steepest descent and ascent contours introduced
in \Cref{c:vertical_cusp_steepest}. Hence by \Cref{p:invariance_contour}, \eqref{e:sumoftwo2} holds with
\(\wh J^{(0)}\) and \(\wh J^{(1)}\) given by \eqref{e:wJ1_case1}.

\section{Properties of the Standard Form and Compatibility}
\label{s:compatibility_proof}

In this section, we establish properties of the standard forms of the kernel ansatz introduced in \Cref{p:standard_form} and prove that the resulting approximation kernels are compatible in the sense of \Cref{p:Bbound}.

\subsection{Proof of \Cref{p:standard_form2}}

We prove \Cref{p:standard_form2} in the case where
\(\fN_{(y,t)}\) is not a ramification neighborhood. The case where
\(\fN_{(y,t)}\) is a ramification neighborhood can be proved in the same
way, so we omit it.

There are two possibilities: \(\fN_\al\) is either a liquid neighborhood
or a ramification neighborhood.

Suppose first that \(\fN_\al\) is a ramification neighborhood. In this
case, \((x,s)\) is associated with a complex-conjugate pair of ramification
charts, which are disjoint from the charts associated with \((y,t)\). By
\Cref{p:standard_form}, up to the error in \eqref{e:errorK}, the kernel
ansatz is given by
\begin{align}\label{e:ramK}
\sum_{w_c,z_c}
\frac{n}{(2\pi\ri)^2}
\int_{\sfD^{\rm a}(z_c)}
\!\!\int_{\sfD^{\rm d}(w_{c,\ft})}
P_{n(s+\ft)}(nw,nx)\,Q_{nt}(nz,ny)\,
\frac{I_\ft(w)}{I_-(z)}\,
\frac{\sqrt{\phi_\ft'(w)}\sqrt{\phi'(z)}}
{\phi_\ft(w)-\phi(z)}\,
\rd w\,\rd z.
\end{align}

Suppose next that \(\fN_\al\) is a liquid neighborhood. By
\Cref{p:standard_form}, up to the error in \eqref{e:errorK}, the kernel
ansatz is given by
\begin{align}
\begin{split}\label{e:liqK}
&-\sum_{\xi}
\sgn(\xi)\,
\frac{n}{2\pi\ri}
\int_{\sfD(\xi;(x,s),(y,t))}
P_{ns}(nz,nx)\,Q_{nt}(nz,ny)\,\rd z
\\
&\quad+
\sum_{w_c,z_c}
\frac{n}{(2\pi\ri)^2}
\int_{\sfD^{\rm a}(z_c)}
\!\!\int_{\sfD^{\rm d}(w_c)}
P_{ns}(nw,nx)\,Q_{nt}(nz,ny)\,
\frac{I(w)}{I_-(z)}\,
\frac{\sqrt{\phi'(w)}\sqrt{\phi'(z)}}
{\phi(w)-\phi(z)}\,
\rd w\,\rd z.
\end{split}
\end{align}
In both \eqref{e:ramK} and \eqref{e:liqK}, the sums range over
$
w_c\in\operatorname{Crit}^{\rm d}(x,s),
$, $
z_c\in\operatorname{Crit}^{\rm a}(y,t),
$
and \(w_{c,\ft}\in\cC_\ft\) is obtained from \(w_c\) by the time shift in
\Cref{c:change_time}. Both formulas hold up to an exponentially small error
of the form
\begin{align}\label{e:errorK}
\sum_{w_c,z_c}
\OO\left(
e^{-\fc'n}
e^{n\Re[S(w_c;x,s)-S(z_c;y,t)]}
\right).
\end{align}
In \eqref{e:liqK}, the sum over \(\xi\) ranges over all nonreal
intersection points of the deformed descent paths
\(\sfD^{\rm d}(w_c)\) and ascent paths \(\sfD^{\rm a}(z_c)\).

By the separation assumption in \Cref{p:standard_form2},  $(x,s)$ is bounded away from $(y,t)$, so the relevant critical points \(w_c\) and \(z_c\) are
uniformly separated. We may therefore slightly truncate the
steepest-descent and steepest-ascent paths in \eqref{e:liqK} so that they
remain uniformly separated and, in particular, are disjoint. By \Cref{p:largerset}, the discarded
portions of the paths are absorbed into the exponentially small error
\eqref{e:errorK}. We continue to denote the truncated contours by the same
symbols. The sum over \(\xi\) is then empty, and \eqref{e:liqK} reduces to
\begin{align}
\begin{split}\label{e:liqK2}
\sum_{w_c,z_c}
\frac{n}{(2\pi\ri)^2}
\int_{\sfD^{\rm a}(z_c)}
\!\!\int_{\sfD^{\rm d}(w_c)}
P_{ns}(nw,nx)\,Q_{nt}(nz,ny)\,
\frac{I(w)}{I_-(z)}\,
\frac{\sqrt{\phi'(w)}\sqrt{\phi'(z)}}
{\phi(w)-\phi(z)}\,
\rd w\,\rd z.
\end{split}
\end{align}

The following lemma evaluates the \(w\)-integrals in
\eqref{e:ramK} and \eqref{e:liqK2}. Applying \Cref{l:onebulk} term by term to \eqref{e:ramK} and
\eqref{e:liqK2}, and combining the resulting expressions with the error
estimate \eqref{e:errorK}, proves \Cref{p:standard_form2}.

\begin{lemma}\label{l:onebulk}
For every relevant pair of critical points \(w_c,z_c\), we have
\begin{align}\begin{split}
\label{e:asymp1}
&\frac{n}{(2\pi\ri)^2}
\int_{\sfD^{\rm a}(z_c)}
\!\!\int_{\sfD^{\rm d}(w_c)}
P_{ns}(nw,nx)\,Q_{nt}(nz,ny)\,
\frac{I(w)}{I_-(z)}\,
\frac{\sqrt{\phi'(w)}\sqrt{\phi'(z)}}
{\phi(w)-\phi(z)}\,
\rd w\,\rd z
\\
&\qquad=
-\frac{1+\OO(1/n)}{(2\pi)^2}
e^{nS(w_c;x,s)}
\int_{\sfD^{\rm a}(z_c)}
\frac{Q_{nt}(nz,ny)}{I_-(z)}\,
\frac{\sqrt{\phi'(z)}\sqrt{\del_x\phi(x,s)}}
{\phi(w_c)-\phi(z)}\,
\rd z,
\end{split}\end{align}
and
\begin{align}\begin{split}
\label{e:asymp2}
&\frac{n}{(2\pi\ri)^2}
\int_{\sfD^{\rm a}(z_c)}
\!\!\int_{\sfD^{\rm d}(w_{c,\ft})}
P_{n(s+\ft)}(nw,nx)\,Q_{nt}(nz,ny)\,
\frac{I_\ft(w)}{I_-(z)}\,
\frac{\sqrt{\phi_\ft'(w)}\sqrt{\phi'(z)}}
{\phi_\ft(w)-\phi(z)}\,
\rd w\,\rd z
\\
&\qquad=
-\frac{1+\OO(1/n)}{(2\pi)^2}
e^{nS(w_c;x,s)}
\int_{\sfD^{\rm a}(z_c)}
\frac{Q_{nt}(nz,ny)}{I_-(z)}\,
\frac{\sqrt{\phi'(z)}\sqrt{\del_x\phi(x,s)}}
{\phi(w_c)-\phi(z)}\,
\rd z.
\end{split}\end{align}
\end{lemma}

\begin{proof}[Proof of \Cref{l:onebulk}]
We first prove \eqref{e:asymp1}. To simplify the notation, write
$
S(w):=S(w;x,s).
$
By \Cref{l:PIQI_bound},
\begin{align}\label{e:PQasymp-used}
P_{ns}(nw,nx)I(w)
=
\frac{\sqrt{s}}
{\sqrt{2\pi n}\sqrt{x-w}\sqrt{w-(x-s)}}
e^{nS(w)}
\left(1+\OO\left(\frac1n\right)\right).
\end{align}
Thus, the \(w\)-dependent part of the integrand can be written as
\[
\frac{1}{\sqrt{2\pi n}}e^{nS(w)}
H(w,z)
\left(1+\OO\left(\frac1n\right)\right),
\]
where
\[
H(w,z)
:=
\frac{\sqrt{s}}
{\sqrt{x-w}\sqrt{w-(x-s)}}
\frac{Q_{nt}(nz,ny)}{I_-(z)}
\frac{\sqrt{\phi'(w)}\sqrt{\phi'(z)}}
{\phi(w)-\phi(z)}.
\]
The error is uniform for \(w\) in a fixed neighborhood of \(w_c\) and
\(z\in\sfD^{\rm a}(z_c)\). Indeed, by the preceding truncation, the
\(w\)- and \(z\)-contours remain uniformly separated.

Since \(w_c\) is a simple critical point of \(S\), we have
\[
S'(w_c)=0,
\qquad
S''(w_c)\neq0.
\]
Thus, near \(w_c\),
\begin{align}\label{e:S_expand_prop}
S(w)
=
S(w_c)
+\frac12S''(w_c)(w-w_c)^2
+\frac16S'''(w_c)(w-w_c)^3
+\OO(|w-w_c|^4).
\end{align}
Moreover, since the prefactor is analytic near \(w_c\),
\begin{align}\label{e:H_expand_prop}
H(w,z)
=
H(w_c,z)
+\del_1H(w_c,z)(w-w_c)
+\OO(|w-w_c|^2),
\end{align}
uniformly for \(z\in\sfD^{\rm a}(z_c)\).

Rescaling
\[
w=w_c+\frac{u}{\sqrt n}
\]
and applying the standard steepest-descent expansion gives
\begin{align}
&\phantom{{}+{}}\int_{\sfD^{\rm d}(w_c)}
\frac{1}{\sqrt{2\pi n}}e^{nS(w)}
H(w,z)\left(1+\OO\left(\frac1n\right)\right)\rd w  \\
& =
\frac{e^{nS(w_c)}}{\sqrt{2\pi n}}
\cdot
\frac{1}{\sqrt n}
\int_{\sqrt n(\sfD^{\rm d}(w_c)-w_c)}
e^{\frac12 S''(w_c)u^2}
H(w_c,z)\left(1+\frac{\OO_{\rm odd}(u)}{\sqrt n}+\OO\left(\frac{1}{n}\right)\right)\rd u  \\
& =
\frac{e^{nS(w_c)}}{n}
\frac{H(w_c,z)}{\sqrt{-S''(w_c)}}
\left(1+\OO\left(\frac{1}{n}\right)\right),
\end{align}
where $\OO_{\rm odd}(u)$ is a polynomial in $u$ containing only odd powers and it integrates to zero over the two local branches of the
steepest-descent contour.
Here the branch of \(\sqrt{-S''(w_c)}\) is determined by the orientation of
the steepest-descent contour.

Multiplying by the outside factor \(n/(2\pi\ri)^2\), we obtain
\begin{align}\begin{split}\label{e:liquid_PQint}
&\phantom{{}={}}\frac{n}{(2\pi\ri)^2}
\int_{\sfD^{\rm a}(z_c)}
\!\!\int_{\sfD^{\rm d}(w_c)}
P_{ns}(nw,nx)\,Q_{nt}(nz,ny)\,
\frac{I(w)}{I_-(z)}
\frac{\sqrt{\phi'(w)}\sqrt{\phi'(z)}}
{\phi(w)-\phi(z)}
\,\rd w\,\rd z
\\
&=
-\frac{1+\OO(1/n)}{(2\pi)^2}
e^{nS(w_c)}
\int_{\sfD^{\rm a}(z_c)}
\frac{H(w_c,z)}{\sqrt{-S''(w_c)}}\,\rd z.
\end{split}\end{align}
Substituting the definition of \(H\) gives
\begin{align}
\frac{H(w_c,z)}{\sqrt{-S''(w_c)}}
&=
\frac{Q_{nt}(nz,ny)}{I_-(z)}
\frac{\sqrt{\phi'(z)}\sqrt{\phi'(w_c)}}
{\phi(w_c)-\phi(z)}
\frac{\sqrt{s}}
{\sqrt{x-w_c}\sqrt{w_c-(x-s)}
\sqrt{-S''(w_c)}}.
\end{align}

Finally, by \eqref{e:Sdirection},
\begin{align}
-\frac{1}{S''(w_c)}
=
s\chi(w_c)\bigl(1-\chi(w_c)\bigr)
\frac{\del_x\phi(x,s)}{\phi'(w_c)}.
\end{align}
Since
$
w_c=x-s\chi(w_c),
$
we have
\[
x-w_c=s\chi(w_c),
\qquad
w_c-(x-s)=s\bigl(1-\chi(w_c)\bigr),
\]
and hence
\[
\del_x\phi(x,s)
=
\frac{s\phi'(w_c)}
{(x-w_c)(w_c-(x-s))[-S''(w_c)]}.
\]
Therefore,
\begin{align}
\frac{H(w_c,z)}{\sqrt{-S''(w_c)}}
=
\frac{Q_{nt}(nz,ny)}{I_-(z)}
\frac{\sqrt{\phi'(z)}\sqrt{\del_x\phi(x,s)}}
{\phi(w_c)-\phi(z)}.
\end{align}
Substituting this identity into \eqref{e:liquid_PQint} proves
\eqref{e:asymp1}.

We next prove \eqref{e:asymp2}. Applying the same argument on
\(\cC_\ft\) gives
\begin{align}
&\phantom{{}+{}}\frac{n}{(2\pi\ri)^2}
\int_{\sfD^{\rm a}(z_c)}
\!\!\int_{\sfD^{\rm d}(w_{c,\ft})}
P_{n(s+\ft)}(nw,nx)\,Q_{nt}(nz,ny)\,
\frac{I_\ft(w)}{I_-(z)}\,
\frac{\sqrt{\phi_\ft'(w)}\sqrt{\phi'(z)}}
{\phi_\ft(w)-\phi(z)}
\,\rd w\,\rd z
\\
&=
-\frac{1+\OO(1/n)}{(2\pi)^2}
e^{nS_\ft(w_{c,\ft};x,s+\ft)}
\int_{\sfD^{\rm a}(z_c)}
\frac{Q_{nt}(nz,ny)}{I_-(z)}\,
\frac{\sqrt{\phi'(z)}
\sqrt{\del_x\phi_\ft(x,s+\ft)}}
{\phi_\ft(w_{c,\ft})-\phi(z)}
\,\rd z.
\end{align}
By \Cref{c:change_time},
\begin{align}
S_\ft(w_{c,\ft};x,s+\ft)
=
S(w_c;x,s),
\quad
\phi_\ft(w_{c,\ft})
=
\phi(w_c),\quad 
\del_x\phi_\ft(x,s+\ft)
=
\del_x\phi(x,s).
\end{align}
Substituting these identities proves \eqref{e:asymp2}.
\end{proof}

\subsection{Proof of \Cref{p:largerset}}
\label{s:largerset_proof}

\begin{proof}[Proof of \Cref{p:largerset}]
We first prove that the formulas
\eqref{e:Aal} and \eqref{e:Aal_ramification} remain valid after restricting
\(w_c\) to
\[
\operatorname{Crit}^{\rm d}(x,s)
\]
and \(z_c\) to
\[
\operatorname{Crit}^{\rm a}(y,t).
\]
Moreover, if \((x,s)\) belongs to an arctic neighborhood, then \(w_c\) may be further
restricted to the critical points contained in the arctic chart, and the
critical point contained in the additional frozen chart may be omitted.
The analogous statement holds for \(z_c\) if \((y,t)\) belongs to an arctic
neighborhood.

Recall from \Cref{r:spillover} that a descent critical point \(w_c'\) to be
omitted is either an additional critical point satisfying
\[
w_c'
\in
\operatorname{Crit}^{\rm d}(x,s;\fN_\al)
\setminus
\operatorname{Crit}^{\rm d}(x,s),
\]
or the descent critical point contained in the additional frozen chart when
\(\fN_\al\) is an arctic neighborhood. Such a point can arise either in the
liquid-spillover case or in the boundary-spillover case in \Cref{l:extra_frozen_action_gap}. In both cases, by
\Cref{l:extra_frozen_action_gap}, there exists
\[
w_c\in\operatorname{Crit}^{\rm d}(x,s),
\]
which may be chosen to lie in the arctic chart in the arctic-neighborhood
case, such that
\begin{align}\label{e:Swcsmall}
\Re S(w_c;x,s)
\geq
\Re S(w_c';x,s)+\fc'
\end{align}
for some \(\fc'>0\).

Let
\[
\xi
\in
\sfD^{\rm d}(w_c')\cap\sfD^{\rm a}(z_c).
\]
Since \(\xi\) lies on the steepest-descent path from \(w_c'\) and the
steepest-ascent path from \(z_c\),
\[
\Re S(\xi;x,s)
\leq
\Re S(w_c';x,s),
\qquad
\Re S(\xi;y,t)
\geq
\Re S(z_c;y,t).
\]
Therefore,
\[
\Re\bigl[S(\xi;x,s)-S(\xi;y,t)\bigr]
\leq
\Re\bigl[S(w_c';x,s)-S(z_c;y,t)\bigr].
\]
By \Cref{l:intersection_descent,l:single_contour_integral}, the corresponding
single-contour integral is bounded by
\[
(\ln n)
e^{n\Re[S(w_c';x,s)-S(z_c;y,t)]}.
\]
By \Cref{l:double_bound}, the corresponding double-contour integral is
bounded by
\[
(\ln n)^2
e^{n\Re[S(w_c';x,s)-S(z_c;y,t)]}.
\]
Using \eqref{e:Swcsmall}, both bounds are at most
\[
(\ln n)^2
e^{-\fc'n}
e^{n\Re[S(w_c;x,s)-S(z_c;y,t)]}.
\]
After decreasing \(\fc'>0\), the factor \((\ln n)^2\) can be absorbed into
\(e^{-\fc'n}\). Thus, all terms indexed by the omitted descent critical
points can be absorbed into the exponentially small error terms.

Similarly, the single- and double-contour integrals indexed by each
additional ascent critical point can be absorbed into the exponentially
small error terms. This proves the first assertion.

We next show that each steepest-descent or steepest-ascent path appearing
in \eqref{e:Aal} and \eqref{e:Aal_ramification} may be truncated to a
sufficiently small neighborhood of its corresponding critical point.

Choose these neighborhoods sufficiently small. Since the real part of the
action decreases strictly along a steepest-descent path away from its
critical point and increases strictly along a steepest-ascent path away from
its critical point, there exists \(\fc'>0\) such that, whenever \(w\) lies
in the discarded portion of \(\sfD^{\rm d}(w_c)\), or \(z\) lies in the
discarded portion of \(\sfD^{\rm a}(z_c)\),
\begin{align}\label{e:truncated_action_gap}
\Re\bigl[S(w;x,s)-S(z;y,t)\bigr]
\leq
\Re\bigl[S(w_c;x,s)-S(z_c;y,t)\bigr]-\fc'.
\end{align}
Hence, by the same argument as in \Cref{l:double_bound}, the contribution
of the discarded portion of a double-contour integral is bounded by
\begin{align}
(\ln n)^2
e^{n\Re[S(w_c;x,s)-S(z_c;y,t)]-\fc'n}.
\end{align}

For a single-contour integral, suppose that an intersection point
\[
\xi
\in
\sfD^{\rm d}(w_c)\cap\sfD^{\rm a}(z_c)
\]
is removed after the two paths are truncated. Then \(\xi\) is bounded away
from at least one of \(w_c\) and \(z_c\). Consequently, either
\[
\Re S(\xi;x,s)
\leq
\Re S(w_c;x,s)-\fc'
\]
or
\[
\Re S(\xi;y,t)
\geq
\Re S(z_c;y,t)+\fc'.
\]
In either case,
\begin{align}
\Re\bigl[S(\xi;x,s)-S(\xi;y,t)\bigr]
\leq
\Re\bigl[S(w_c;x,s)-S(z_c;y,t)\bigr]-\fc'.
\end{align}
By \Cref{l:single_contour_integral}, the corresponding single-contour
integral is therefore bounded by
\begin{align}
(\ln n)
e^{n\Re[S(w_c;x,s)-S(z_c;y,t)]-\fc'n}.
\end{align}

After decreasing \(\fc'>0\) if necessary, the factors \(\ln n\) and
\((\ln n)^2\) can again be absorbed into the exponential decay. Thus, the
contributions of all discarded portions are absorbed into the exponentially
small error terms in \eqref{e:Aal_ramification} and \eqref{e:Aal}.
\end{proof}

\subsection{Proof of compatibility}\label{s:compatible}
The compatibility result \Cref{p:white_triangle_critical_value} is a consequence of the following statement.
\begin{proposition}\label{p:compatible_error}
Let \((y,t)\in\fP\cap\bZ^2/n\) represent a white triangle, and let
\((x,s)\in\fP\cap\bZ^2/n\) represent a blue triangle. Suppose that
\[
(y,t)\in\fN_{(y,t)}
\qquad\text{and}\qquad
(x,s)\in\fN_\al\cap\fN_\beta
\]
for some \(\alpha\neq\beta\). We recall $\fO$ and $\fW$ from \Cref{p:overlap_critical_region}. Then the following holds
\begin{align}\begin{split}\label{e:compatible_error}
&\phantom{{}={}}\left|A_\alpha((x,s),(y,t))-A_\beta((x,s),(y,t))\right|\\
&\lesssim
\sum_{w_c,z_c}
\left(\bm1\bigl((x,s)\in\fO, w_c\in \fW\bigr) \frac{(\ln n)^{10}\Delta(w_c)\Delta(z_c)}{n}
+
e^{-\fc'n}\right)
e^{n\Re\left[S(w_c;x,s)-S(z_c;y,t)\right]}
\end{split}\end{align}
where  \(z_c\) ranges over
\(\operatorname{Crit}^{\rm a}(y,t)\) and  \(w_c\) ranges over
\(\operatorname{Crit}^{\rm d}(x,s)\), as defined in
\Cref{def:ascent_descent_critical}.
\end{proposition}

\begin{proof}[Proof of \Cref{p:compatible_error}]
Suppose first that either both \(\fN_\al\) and \(\fN_\beta\) are
ramification neighborhoods, or neither is a ramification neighborhood and,
for every \(w_c\in\operatorname{Crit}^{\rm d}(x,s)\), the corresponding
\(I_+\)-factor is the same for \(\fN_\alpha\) and \(\fN_\beta\). Then
\Cref{p:standard_form} and \Cref{p:largerset} imply
\begin{align}\label{e:compatible_exp_error}
A_\alpha((x,s),(y,t))-A_\beta((x,s),(y,t))
=
\sum_{w_c,z_c}
\OO\left(
e^{-\fc'n}
e^{n\Re\left[S(w_c;x,s)-S(z_c;y,t)\right]}
\right).
\end{align}
Thus, \eqref{e:compatible_error} follows in this case.

It remains to consider the cases in which one of
\(\fN_\al,\fN_\beta\) is a ramification neighborhood, as in the
\emph{ramification--liquid overlap} case, or the \(I_+\)-assignment differs
between the two neighborhoods, as in the \emph{\(I_i-I\) overlap} case.

\medskip
\noindent\textbf{Ramification--liquid overlap: \((x,s)\in\fO_1\).}
Without loss of generality, assume that \(\fN_\al\) is a ramification
neighborhood and that \(\fN_\beta\) is a liquid neighborhood. Since
\((x,s)\in\fN_\beta\), the relevant critical points associated with
\((x,s)\) are uniformly separated from the ramification point. Moreover,
by construction, \(\fN_{(y,t)}\) is bounded away from all ramification
points, so the collection of charts associated with \((y,t)\) contains no
ramification chart. Consequently, the charts associated with
\(\fN_\alpha\) and those associated with \(\fN_{(y,t)}\) are disjoint.
In particular, every
$
w_c\in\operatorname{Crit}^{\rm d}(x,s)
$
is uniformly separated from every
$
z_c\in\operatorname{Crit}^{\rm a}(y,t).
$

By \eqref{e:Aal_ramification},
\begin{align}
A_\al((x,s),(y,t))
&=
\sum_{w_c,z_c}
\frac{(1+\OO(1/n))}{(2\pi\ri)^2}
e^{nS(w_c;x,s)}
\int_{\sfD^{\rm a}(z_c)}
\frac{Q_{nt}(nz,ny)}{I_-(z)}
\frac{\sqrt{\phi'(z)\,\del_x\phi(x,s)}}
{\phi(w_c)-\phi(z)}
\,\rd z
\nonumber\\
&\quad+
\sum_{w_c,z_c}
\OO\left(
e^{-\fc'n}
e^{n\Re\left[S(w_c;x,s)-S(z_c;y,t)\right]}
\right),
\end{align}
and 
\(A_\beta((x,s),(y,t))\) has the same leading expression. Subtracting the two
representations therefore gives
\begin{align}
A_\al((x,s),(y,t))-A_\beta((x,s),(y,t))
&=
\sum_{w_c,z_c}
\frac{\OO(1)}{n}
e^{nS(w_c;x,s)}
\int_{\sfD^{\rm a}(z_c)}
\frac{Q_{nt}(nz,ny)}{I_-(z)}
\frac{\sqrt{\phi'(z)\,\del_x\phi(x,s)}}
{\phi(w_c)-\phi(z)}
\,\rd z
\nonumber\\
&\quad+
\sum_{w_c,z_c}
\OO\left(
e^{-\fc'n}
e^{n\Re\left[S(w_c;x,s)-S(z_c;y,t)\right]}
\right).
\end{align}
Applying \eqref{l:double_bound}, we obtain
\begin{align}\label{e:O1_compatible}
\left|A_\al((x,s),(y,t))-A_\beta((x,s),(y,t))\right|
&\lesssim
\sum_{w_c,z_c}
\frac{(\ln n)^{10}\Delta(w_c)\Delta(z_c)}{n}
e^{n\Re\left[S(w_c;x,s)-S(z_c;y,t)\right]}
\nonumber\\
&\quad+
e^{-\fc'n}
\sum_{w_c,z_c}
e^{n\Re\left[S(w_c;x,s)-S(z_c;y,t)\right]}.
\end{align}

\medskip
\noindent\textbf{\(I_i-I\) overlap: \((x,s)\in\fO_2\).}
In this case, \((x,s)\) is associated with a descent critical point \(w_c\)
lying outside, but sufficiently close to, \([b_i,a_i]\), so that \(w_c\)
is contained in both a chart \(\fU\) centered at a point
\(w_0\in[b_i,a_i]\) and a chart \(\fU'\) centered at a point
\(w_0'\notin[b_i,a_i]\). In particular $w_c\in \fW$. Without loss of generality, assume that the
double-contour formula associated with \(\fN_\al\) uses \(I_+=I_i\),
whereas the formula associated with \(\fN_\beta\) uses \(I_+=I\).

The correction terms \(\wh J^{(0)}\) and \(\wh J^{(1)}\) are the same in the two representations
and therefore cancels. Moreover, the double-contour integrals associated
with all descent critical points outside \(\fU\cap\fU'\) also cancel.
Consequently, \Cref{p:standard_form} gives
\begin{align}\label{e:O2_difference}
&A_\al((x,s),(y,t))-A_\beta((x,s),(y,t))=
\sum_{w_c,z_c}
\OO\left(
e^{-\fc'n}
e^{n\Re\left[S(w_c;x,s)-S(z_c;y,t)\right]}
\right)\\
+
&\sum_{w_c, z_c}
\frac{\bm1(w_c\in \fW) n}{(2\pi\ri)^2}
\int_{\sfD^{\rm a}(z_c)}
\int_{\sfD^{\rm d}(w_c)}
P_{ns}(nw,nx)\,Q_{nt}(nz,ny)
\frac{I_i(w)-I(w)}{I_-(z)}
\frac{\sqrt{\phi'(w)}\sqrt{\phi'(z)}}
{\phi(w)-\phi(z)}
\,\rd w\,\rd z
\end{align}

We first observe that every
\(z_c\in\operatorname{Crit}^{\rm a}(y,t)\) is uniformly separated from
\(w_c\). Indeed, let \(\fV\) be a chart associated with \((y,t)\) that
contains \(z_c\). If \(z_c\) were not uniformly separated from \(w_c\),
then \(\fV\) would intersect both \(\fU\) and \(\fU'\). By the chart
construction in \Cref{p:construct_neighborhood2}, \(\fV\) would then have
to be concentric with both \(\fU\) and \(\fU'\). This would imply
\(w_0=w_0'\), contradicting
\[
w_0\in[b_i,a_i],
\qquad
w_0'\notin[b_i,a_i].
\]
Thus, \(w_c\) is uniformly separated from every ascent critical point
associated with \((y,t)\).

By \Cref{p:standard_form2}, after truncating the local steepest-descent and steepest-ascent paths, if
necessary, we may therefore assume that
\(\sfD^{\rm d}(w_c)\) and \(\sfD^{\rm a}(z_c)\) are uniformly separated
and that \(\sfD^{\rm d}(w_c)\) is disjoint from \([b_i,a_i]\). By
\Cref{l:PIQI_bound} and \Cref{l:PIi_bound}, uniformly for
\(w\in\sfD^{\rm d}(w_c)\),
\begin{align}\label{e:IiI_difference}
P_{ns}(nw,nx)I_i(w)
&=
P_{ns}(nw,nx)I(w)\left(1+\OO(1/n)\right).
\end{align}
Substituting \eqref{e:IiI_difference} into
\eqref{e:O2_difference} and applying \eqref{l:double_bound}, we obtain
\begin{align}\label{e:O2_compatible}
\left|A_\al((x,s),(y,t))-A_\beta((x,s),(y,t))\right|
&\lesssim
\sum_{w_c,z_c}\bm1(w_c\in \fW)
\frac{(\ln n)^{10}\Delta(w_c)\Delta(z_c)}{n}
e^{n\Re\left[S(w_c;x,s)-S(z_c;y,t)\right]}
\nonumber\\
&\quad+
e^{-\fc'n}
\sum_{w_c,z_c}
e^{n\Re\left[S(w_c;x,s)-S(z_c;y,t)\right]}.
\end{align}

Combining \eqref{e:compatible_exp_error},
\eqref{e:O1_compatible}, and \eqref{e:O2_compatible} proves
\eqref{e:compatible_error}.
\end{proof}

\begin{proposition}\label{p:white_triangle_critical_value}
Let \((x,s)\in\bZ^2/n\) represent a white triangle \(\rw\) contained in
\(\fP\), and let \((x',s')\) lie in the interior of \(\rw\). For any
\[
\zeta
\in
\left\{
(x,s),
\left(x-\frac1n,s-\frac1n\right),
\left(x,s-\frac1n\right)
\right\}
\cap\fP
\]
and any
\[
w_c\in\operatorname{Crit}^{\rm d}(u,v),
\]
there exists
\[
\xi_c^{\rm d}\in\operatorname{Crit}^{\rm d}(x',s')
\]
such that
\begin{align}\label{e:white_triangle_critical_value}
\Re S(w_c;\zeta)
\leq
\Re S(\xi_c^{\rm d};x',s')
+
\OO\left(\frac{\ln n}{n}\right).
\end{align}
If the corresponding critical-point branches remain uniformly bounded away
from the tangent locations, then the error in
\eqref{e:white_triangle_critical_value} can be improved to
\begin{align}\label{e:white_triangle_critical_value_regular}
\Re S(w_c;\zeta)
\leq
\Re S(\xi_c^{\rm d};x',s')
+
\OO\left(\frac1n\right).
\end{align}
\end{proposition}

\begin{proof}[Proof of \Cref{p:white_triangle_critical_value}]
The diameter of the white triangle \(\rw\) is of order \(n^{-1}\).
For each of the three possible choices of \(\zeta\), we can choose a path
from \((x',s')\) to \(\zeta\) of length \(\OO(n^{-1})\) that crosses an
extended side at most once. Denote a point along each path by
\((u,v)\). More precisely:
\begin{enumerate}
\item
The path from \((x',s')\) to \((x,s)\) can meet a horizontal extended
side only at the endpoint \((x,s)\).

\item
The path from \((x',s')\) to
$(x-1/n, s-1/n)$ can cross only the vertical extended side
\[
u=b_i=x-\frac{1}{2n},
\]
and, if it does, it crosses this side once from right to left.

\item
The path from \((x',s')\) to
$(x,s-1/n)$ can cross only the unit-slope extended side
\[
u-v=a_i=x-s+\frac{1}{2n},
\]
and, if it does, it crosses this side once from left to right.
\end{enumerate}

Fix \(w_c\in\operatorname{Crit}^{\rm d}(\zeta)\), and follow its
critical-point branch backward along the chosen path from \(\zeta\) to
\((x',s')\). Away from the extended sides, the branch varies continuously.
If the path crosses the arctic boundary, two critical-point branches may
coalesce and split, but the critical points and their critical values remain
continuous.

Suppose that the path crosses an extended side. In the tangent or
tangent frozen case, one descent critical point disappears in the direction
from \((x',s')\) to \(\zeta\). Therefore, no new descent critical point is
created in this direction, and every descent critical point at \(\zeta\)
can be followed backward to a descent critical point at \((x',s')\).

In the cusp-turning case, the selected descent critical point may jump from
one branch to another. By \Cref{l:motion_descent_critical}, and by its
unit-slope and horizontal analogues, the real part of the critical value
does not increase across the jump in the direction from \((x',s')\) to
\(\zeta\). Thus, if \(w_c^{\rm before}\) and \(w_c^{\rm after}\) denote
the descent critical points immediately before and after the crossing, then
\begin{align}\label{e:white_triangle_jump}
\Re S(w_c^{\rm after})
\leq
\Re S(w_c^{\rm before}),
\end{align}
where both critical values are evaluated at the crossing point. At a
horizontal side met only at the endpoint, this inequality is understood in
the corresponding one-sided sense.

It follows that \(w_c\) can be traced backward, allowing such a jump if
necessary, to a descent critical point $\xi_c^{\rm d}$ associated with $(x',s')$.

We now compare the critical values. Split the path at the extended-side
crossing, if one occurs. Each resulting subpath has length \(\OO(n^{-1})\)
and does not cross any extended side. Hence
\eqref{e:critical_value_continuity}, applied with
\(\delta=\OO(n^{-1})\), gives
\begin{align}\label{e:point_change}
\left|
\Re S(w_c(p);p)
-
\Re S(w_c(q);q)
\right|
\lesssim
\frac{\ln n}{n}
\end{align}
between the endpoints \(p,q\) of each subpath. Combining these estimates
with the favorable jump inequality \eqref{e:white_triangle_jump} yields
\begin{align}
\Re S(w_c;u,v)
\leq
\Re S(\xi_c^{\rm d};x',s')
+
\OO\left(\frac{\ln n}{n}\right),
\end{align}
which proves \eqref{e:white_triangle_critical_value}.

If the corresponding critical-point branches remain uniformly bounded away
from the tangent locations, by \Cref{l:critical_value_continuity}, we have the improved estimate \(\OO(n^{-1})\). This proves
\eqref{e:white_triangle_critical_value_regular}.

\end{proof}

\begin{proof}[Proof of \Cref{p:Bbound}]
Let \(\rw\) denote the white triangle represented by
\((x,s)\in\bZ^2/n\). Choose $(x',s')\in \rw$ and let \((x',s')\to(x,s)\) from below, as in the definition of
\(\operatorname{Crit}^{\rm d}(x,s^-)\) in \eqref{e:defCdown}. Fix
\[
\zeta
\in
\left\{
(x,s),
\left(x-\frac1n,s-\frac1n\right),
\left(x,s-\frac1n\right)
\right\}.
\]
We apply \Cref{p:compatible_error} at \(\zeta\) and estimate the two terms
on the right-hand side of \eqref{e:compatible_error} separately.

We first consider the term containing
\(\bm1((x,s)\in\fO)/n\). For each
\[
w_c(\zeta)
\in
\operatorname{Crit}^{\rm d}(\zeta;\fN_\al)
\cap \fW,
\]
the improved estimate
\eqref{e:white_triangle_critical_value_regular} gives a critical point
\[
\xi_c^{\rm d}(x',s')
\in
\operatorname{Crit}^{\rm d}(x',s';\fN_\al)
\]
such that
\begin{align}\label{e:first_term_action}
\Re S(w_c(\zeta);\zeta)
\leq
\Re S(\xi_c^{\rm d}(x',s');x',s')
+
\OO\left(\frac1n\right).
\end{align}
Moreover,
$
\Delta(w_c(\zeta))
\asymp
\Delta(\xi_c^{\rm d}(x',s')).
$
Therefore,
\begin{align}\label{e:first_termbb}
\frac{\Delta(w_c(\zeta))\Delta(z_c)}{n}
e^{n\Re\left[
S(w_c(\zeta);\zeta)-S(z_c;y,t)
\right]}
\lesssim
\frac{\Delta(\xi_c^{\rm d}(x',s'))\Delta(z_c)}{n}
e^{n\Re\left[
S(\xi_c^{\rm d}(x',s');x',s')-S(z_c;y,t)
\right]}.
\end{align}
Letting \((x',s')\to(x,s)\), these critical points converge to points in
\(\operatorname{Crit}^{\rm d}(x,s^-)\). Hence the first
term on the right-hand side of \eqref{e:compatible_error} is bounded by
$
(\ln n)^{10} B((x,s),(y,t)).
$

We next consider the exponentially small term in
\eqref{e:compatible_error}. By
\eqref{e:white_triangle_critical_value}, for every
$
w_c\in\operatorname{Crit}^{\rm d}(\zeta;\fN_\al)
$
there exists
$
\xi_c^{\rm d}
\in
\operatorname{Crit}^{\rm d}(x',s';\fN_\al)
$
such that
\[
\Re S(w_c;\zeta)
\leq
\Re S(\xi_c^{\rm d};x',s')
+
\OO\left(\frac{\ln n}{n}\right).
\]
Consequently,
\begin{align}\label{e:second_termbb}
\sum_{w_c,z_c}
e^{n\Re\left[
S(w_c;\zeta)-S(z_c;y,t)
\right]}
\lesssim
e^{\OO(\ln n)}
\sum_{\xi_c^{\rm d},z_c}
e^{n\Re\left[
S(\xi_c^{\rm d};x',s')-S(z_c;y,t)
\right]}.
\end{align}
Multiplying \eqref{e:second_termbb} by \(e^{-\fc'n}\), the 
factor \(e^{\OO(\ln n)}\) can be absorbed into the exponential after
decreasing \(\fc'>0\).

Letting \((x',s')\to(x,s)\) and applying
\Cref{p:critical_compare}, we may replace the descent critical points
on the right-hand side by those in
\(\operatorname{Crit}^{\rm d}(x,s^-)\). Since the localization lengths are
bounded below by a negative power of \(n\), the missing factors
\(\Delta(\xi_c^{\rm d})\Delta(z_c)\) can also be absorbed into the
exponentially small factor after decreasing \(\fc'>0\) once more. Thus, the
second term on the right-hand side of \eqref{e:compatible_error} is 
bounded by
$
B((x,s),(y,t)).
$

Since \(\zeta\) was arbitrary among the three points above, this proves
\eqref{e:Bbound}.
\end{proof}

\section{Estimates for the Inverse Kasteleyn Matrix}\label{s:final_kernel}
In this section, we prove
\Cref{p:liquid,p:arctic,p:frozen}. The proofs reduce to kernel estimates in
the following two cases:
\begin{enumerate}
\item
The points \((x,s)\) and \((y,t)\) lie in liquid, ramification, or 
arctic neighborhoods and are bounded away from each other.

\item
The points \((x,s)\) and \((y,t)\) both lie in frozen neighborhoods.
\end{enumerate}

Recall from \Cref{l:global_approximation} that
\begin{align}\label{e:Kexp}
K^{-1}((x,s);(y,t))
=
A((x,s);(y,t))
+
\OO\left(
\frac{(\ln n)^{10}}{n}
(|A|B)((x,s);(y,t))
\right).
\end{align}
For any \(\alpha\) such that \((x,s)\in\fN_\alpha\), the construction of
\(A((x,s);(y,t))\) in \eqref{e:defA}, together with \eqref{e:Bbound},
gives
\begin{align}\label{e:Aexp}
A((x,s);(y,t))
=
A_\alpha((x,s);(y,t))
+
\OO\left(
(\ln n)^{10}B((x,s);(y,t))
\right).
\end{align}
Moreover, by \Cref{p:standard_form},
\begin{align}\label{e:Aexp2}
A_\alpha((x,s);(y,t))
&=
\wh J^{(0)}((x,s);(y,t))
+
\wh J^{(1)}((x,s);(y,t))
+
\wh J^{(2)}((x,s);(y,t))
\notag\\
&\quad+
\OO\left(
(\ln n)^{10}B((x,s);(y,t))
\right).
\end{align}

\subsection{Liquid and arctic cases}
We recall the exceptional arctic set \(\fA_{\rm ex}\) from
\eqref{def:exA}. In this section, we fix a sufficiently small \(\fb>0\)
and points
$
(x,s),(y,t)\in\fP\cap\bZ^2/n
$
lying in a sufficiently small neighborhood of the liquid region \(\fL\)
and satisfying
\begin{align}\label{e:close}
\dist((x,s),\fL)\leq n^{-\delta},
\qquad
\dist((y,t),\fL)\leq n^{-\delta},
\end{align}
and
\begin{align}\label{e:far}
\|(x,s)-(y,t)\|_2&\geq\fb,
&
\dist((x,s),\fA_{\mathrm{ex}})&\geq\fb,
&
\dist((y,t),\fA_{\mathrm{ex}})&\geq\fb.
\end{align}
These are also among the assumptions of \Cref{p:liquid,p:arctic}. By the
defining property of \(\fA_{\rm ex}\), Then \((x,s)\) and \((y,t)\) are
contained in liquid, ramification, or arctic neighborhoods.

Fix \(\alpha\) such that \((x,s)\in\fN_\alpha\). By
\Cref{p:standard_form,p:largerset}, in the standard form of the approximate
kernel \(A_\alpha((x,s);(y,t))\), we may retain only the critical points
contained in the corresponding liquid, ramification, or arctic charts; the
contributions from any additional frozen charts are absorbed into the
exponentially small error terms. Moreover, by \eqref{e:far} and
\Cref{p:largerset}, after slightly truncating the steepest-descent and
steepest-ascent paths, we may assume that every relevant pair of paths is
uniformly separated. The discarded portions are again absorbed into the
exponentially small error terms.

For instance, if neither \(\fN_\alpha\) nor \(\fN_{(y,t)}\) is a
ramification neighborhood, then each double-contour integral in the
standard form of \(A_\alpha((x,s);(y,t))\) is of the form
\begin{align}\label{e:asymp}
\frac{n}{(2\pi\ri)^2}
\int_{\sfD^{\rm a}(z_c)}\!\!\int_{\sfD^{\rm d}(w_c)}
P_{ns}(nw,nx)\,Q_{nt}(nz,ny)\,
\frac{I(w)}{I(z)}\,
\frac{\sqrt{\phi'(w)}\sqrt{\phi'(z)}}
{\phi(w)-\phi(z)}\,
\rd w\,\rd z,
\end{align}
where \(w_c\) and \(z_c\) are the relevant descent and ascent critical
points, respectively.

\begin{lemma}\label{l:liquidregion}
Assume \eqref{e:close} and \eqref{e:far}. Then
\begin{align}\label{e:Aal_liquid_arctic_bound}
\left|A_\al((x,s);(y,t))\right|
\lesssim
\sum_{w_c,z_c}
(\ln n)^{10}\Delta(w_c)\Delta(z_c)
e^{n\Re[S(w_c;x,s)-S(z_c;y,t)]},
\end{align}
and
\begin{align}\label{e:K_liquid_arctic_bound}
|K^{-1}((x,s);(y,t))
-
A_\al((x,s);(y,t))|\lesssim
\sum_{w_c,z_c}
\frac{(\ln n)^{20}\Delta(w_c)\Delta(z_c)}{n}
e^{n\Re[S(w_c;x,s)-S(z_c;y,t)]}.
\end{align}
Here \(w_c\) ranges over the critical points in
\(\operatorname{Crit}^{\rm d}(x,s)\) contained in the corresponding liquid,
ramification, or arctic charts, and \(z_c\) ranges over the critical points
in \(\operatorname{Crit}^{\rm a}(y,t)\) contained in the corresponding
liquid, ramification, or arctic charts.
\end{lemma}

The following lemma gives a refined estimate for the double-contour
integral \eqref{e:asymp}and its analogues in the ramification cases.

\begin{lemma}\label{l:kernel_liquid}
Assume \eqref{e:close} and \eqref{e:far}. We assume in addition that $(x,s), (y,t)\in \fL$ and 
\[
\dist((x,s),\fA)\geq n^{-\delta},
\qquad
\dist((y,t),\fA)\geq n^{-\delta}.
\]
Suppose first that neither \(\fN_\alpha\) nor \(\fN_{(y,t)}\) is a
ramification neighborhood. Then, uniformly for every relevant pair
\((w_c,z_c)\),
\begin{align}\label{e:asymp2}
\eqref{e:asymp}=
-\frac{\left(1+\OO(n^{-1+2\delta})\right)e^{n(S(w_c;x,s)-S(z_c;y,t))}}{2\pi n\ri}
\begin{cases}
\dfrac{
\sqrt{\partial_x\phi(x,s)}
\sqrt{\partial_x\phi(y,t)}
}{\phi(x,s)-\phi(y,t)},
& w_c,z_c\in\bC_+,\\[3mm]
\dfrac{
\sqrt{\partial_x\phi(x,s)}
\overline{\sqrt{\partial_x\phi(y,t)}}
}{\phi(x,s)-\overline{\phi(y,t)}},
& w_c\in\bC_+,\ z_c\in\bC_-,\\[3mm]
-\dfrac{
\overline{\sqrt{\partial_x\phi(x,s)}}
\sqrt{\partial_x\phi(y,t)}
}{\overline{\phi(x,s)}-\phi(y,t)},
& w_c\in\bC_-,\ z_c\in\bC_+,\\[3mm]
-\dfrac{
\overline{\sqrt{\partial_x\phi(x,s)}}
\overline{\sqrt{\partial_x\phi(y,t)}}
}{\overline{\phi(x,s)}-\overline{\phi(y,t)}},
& w_c,z_c\in\bC_-.
\end{cases}
\end{align}

If \(\fN_\alpha\) is a ramification neighborhood, then
\eqref{e:asymp2} remains valid after replacing, on its left-hand side,
\[
P_{ns}(nw,nx), \sfD^{\rm d}(w_c), I(w), \phi(w), \phi'(w)\quad
\text{by}\quad 
P_{n(s+\ft)}(nw,nx),\sfD^{\rm d}(w_{c,\ft}), I_\ft(w),
\phi_\ft(w), \phi_\ft'(w),
\]
respectively. Similarly, if \(\fN_{(y,t)}\) is a ramification
neighborhood, then \eqref{e:asymp2} remains valid after replacing, on its
left-hand side,
\[
Q_{nt}(nz,ny), \sfD^{\rm a}(z_c), I_-(z), \phi(z), \phi'(z)
\quad
\text{by}\quad 
Q_{n(t+\ft)}(nz,ny), \sfD^{\rm a}(z_{c,\ft}), I_\ft(z),
\phi_\ft(z), \phi_\ft'(z),
\]
respectively. If both neighborhoods are ramification neighborhoods, both
sets of replacements are made. Here \(w_{c,\ft}\) and \(z_{c,\ft}\)
denote the corresponding critical points in the ramification charts. In
all cases, the right-hand side of \eqref{e:asymp2} remains unchanged.
\end{lemma}

\begin{proof}[Proof of \Cref{l:liquidregion}]
We prove the case in which neither \(\fN_\al\) nor \(\fN_{(y,t)}\) is a
ramification neighborhood. The cases involving ramification neighborhoods
follow from the same argument, after making the corresponding replacements
in \Cref{p:standard_form} and using \Cref{c:change_time}, so we omit them.

\medskip
\noindent\emph{Step 1: Preliminary bounds.}
We first record bounds for the three terms in the standard-form expansion
\eqref{e:Aexp2}. Fix
$
(u,v)\in\fP\cap\bZ^2/n.
$
By \Cref{l:single_contour_integral},
\begin{align}\label{e:frozen_J0}
\left|
\wh J^{(0)}((x,s);(u,v))
\right|
\lesssim
(\ln n)\,
\cI((x,s);(u,v))\,
e^{n\Re[S(\xi;x,s)-S(\xi;u,v)]}.
\end{align}
Whenever the interlacing condition holds, fix
$
\xi_c^{\rm d}\in
\operatorname{Crit}^{\rm d}(u,v^-).
$
By \Cref{l:interlace_critical_point}, there exists
$
w_c\in\operatorname{Crit}^{\rm d}(x,s)
$
such that
\[
\Re S(w_c;u,v)
\geq
\Re S(\xi_c^{\rm d};u,v).
\]
Applying \eqref{e:frozen_J0} with \(\xi=w_c\), we obtain
\begin{align}\label{e:J0bb}
\left|
\wh J^{(0)}((x,s);(u,v))
\right|
\lesssim
(\ln n)
e^{n\Re[S(w_c;x,s)-S(w_c;u,v)]}
\lesssim
(\ln n)
e^{n\Re[S(w_c;x,s)-S(\xi_c^{\rm d};u,v)]}.
\end{align}

By \Cref{l:intersection_descent,l:single_contour_integral},
\begin{align}\label{e:J1bb}
\left|
\wh J^{(1)}((x,s);(u,v))
\right|
\lesssim
\sum_{w_c,\xi_c^{\rm a}}
(\ln n)
e^{n\Re[S(w_c;x,s)-S(\xi_c^{\rm a};u,v)]},
\end{align}
and, by \Cref{l:double_bound},
\begin{align}\label{e:J2bb}
\left|
\wh J^{(2)}((x,s);(u,v))
\right|
\lesssim
\sum_{w_c,\xi_c^{\rm a}}
(\ln n)^2
e^{n\Re[S(w_c;x,s)-S(\xi_c^{\rm a};u,v)]}.
\end{align}

\medskip
\noindent\emph{Step 2: Proof of \eqref{e:Aal_liquid_arctic_bound}.}
We now take \((u,v)=(y,t)\). Since \((x,s)\) and \((y,t)\) are bounded
away from each other, \Cref{p:largerset} allows us to truncate the relevant
steepest-descent and steepest-ascent paths so that every relevant pair is
uniformly separated. The discarded portions are absorbed into the
exponentially small error in \eqref{e:Aexp2}. Consequently,
\begin{align}\label{e:ytJ1}
\wh J^{(1)}((x,s);(y,t))=0,
\end{align}
and the separated-contour estimate in \Cref{l:double_bound} gives
\begin{align}\label{e:J2_separated}
\left|
\wh J^{(2)}((x,s);(y,t))
\right|
\lesssim
\sum_{w_c,z_c}
(\ln n)^{10}
\Delta(w_c)\Delta(z_c)
e^{n\Re[S(w_c;x,s)-S(z_c;y,t)]}.
\end{align}

It remains to estimate \(\wh J^{(0)}((x,s);(y,t))\). If this term
vanishes, there is nothing to prove. Otherwise, both $(x,s)$ and $(y,t)$ are in the frozen region, are contained in arctic charts, the interlacing condition
holds. 
Let
\[
\xi_c^{\rm d}\in\operatorname{Crit}^{\rm d}(y,t)
\]
be the corresponding descent critical point, and let
\[
w_c\in\operatorname{Crit}^{\rm d}(x,s)
\]
be the critical point furnished by
\Cref{l:interlace_critical_point}. By \eqref{e:far}, these two critical
points are uniformly separated. Hence the strict part of
\Cref{l:interlace_critical_point} gives
\begin{align}\label{e:Syt0}
\Re S(w_c;y,t)
\geq
\Re S(\xi_c^{\rm d};y,t)+\fc'
\end{align}
for some \(\fc'>0\).

The corresponding arctic chart associated with \((y,t)\) contains an
ascent critical point \(z_c\) and the descent critical point
\(\xi_c^{\rm d}\). By \eqref{e:close} and
\eqref{e:Sdiff_at_critical},
\begin{align}\label{e:Syt1}
0
\leq
\Re\left[
S(z_c;y,t)-S(\xi_c^{\rm d};y,t)
\right]
\asymp
\dist((y,t),\fA)^{3/2}
\lesssim
n^{-3\delta/2}.
\end{align}
Thus, after decreasing \(\fc'>0\), \[
\Re S(w_c;y,t)
\geq
\Re S(z_c;y,t)+\frac{\fc'}{2}.
\]
Applying \eqref{e:frozen_J0} with \(\xi=w_c\) and absorbing the factor
\(\ln n\) into the exponential decay, we obtain
\begin{align}\label{e:J0_strict_compare}
\left|
\wh J^{(0)}((x,s);(y,t))
\right|
\lesssim
e^{n\Re[S(w_c;x,s)-S(z_c;y,t)]-\fc'n/3}.
\end{align}

Combining \eqref{e:ytJ1}, \eqref{e:J2_separated}, and
\eqref{e:J0_strict_compare} with the standard-form expansion
\eqref{e:Aexp2} proves \eqref{e:Aal_liquid_arctic_bound}. Here and below,
the exponentially small terms are absorbed into the stated bounds because
\(\Delta(w_c)\) and \(\Delta(z_c)\) are bounded below by negative powers of
\(n\).

\medskip
\noindent\emph{Step 3: Estimate of the convolution \((|A|B)\).}
By \eqref{e:Aexp} and \eqref{e:Aexp2}, we have
\begin{align}\label{e:AB_decomposition}
(|A|B)
&\lesssim
|\wh J^{(0)}|B
+
\bigl(|\wh J^{(1)}|+|\wh J^{(2)}|\bigr)B
+
(\ln n)^{10}B^2.
\end{align}
We estimate these three terms separately.

\smallskip
\noindent\emph{The \(\wh J^{(0)}\)-term.}
Suppose that \((u,v)\in\fO\) and that the descent critical point
$
\xi_c^{\rm d}\in\operatorname{Crit}^{\rm d}(u,v^-)
$
appearing in the overlap part of \(B\) belongs to \(\fW\). Since the
relevant critical points \(w_c\) associated with \((x,s)\) are uniformly
bounded away from \(\fW\), the strict part of
\Cref{l:interlace_critical_point} gives
\[
\Re S(w_c;u,v)
\geq
\Re S(\xi_c^{\rm d};u,v)+\fc'
\]
for some \(\fc'>0\). Applying \eqref{e:frozen_J0} with \(\xi=w_c\), we
obtain
\begin{align}\label{e:J0_strict_compare2}
\left|
\wh J^{(0)}((x,s);(u,v))
\right|
\lesssim
(\ln n)
e^{n\Re[S(w_c;x,s)-S(\xi_c^{\rm d};u,v)]-\fc'n}.
\end{align}

We now split \(B((u,v);(y,t))\) into its exponentially small part and its
overlap part. Using \eqref{e:J0bb} for the former and
\eqref{e:J0_strict_compare2} for the latter, we obtain
\begin{align}\begin{split}\label{e:J0B}
&\bigl(|\wh J^{(0)}|B\bigr)((x,s);(y,t))\lesssim
\sum_{(u,v)}
\sum_{w_c,\xi_c^{\rm d},z_c}
(\ln n)
e^{n\Re[S(w_c;x,s)-S(\xi_c^{\rm d};u,v)]}
e^{-\fc n}
\Delta(\xi_c^{\rm d})\Delta(z_c)
e^{n\Re[S(\xi_c^{\rm d};u,v)-S(z_c;y,t)]}
\\
&+
\sum_{(u,v)\in\fO}
\sum_{w_c,\xi_c^{\rm d},z_c}
(\ln n)
e^{n\Re[S(w_c;x,s)-S(\xi_c^{\rm d};u,v)]-\fc'n}
\frac{
\bm1(\xi_c^{\rm d}\in\fW)
\Delta(\xi_c^{\rm d})\Delta(z_c)
}{n}
e^{n\Re[S(\xi_c^{\rm d};u,v)-S(z_c;y,t)]}
\\
&\lesssim
\sum_{w_c,z_c}
(\ln n)
\left(
n^2e^{-\fc n}
+
ne^{-\fc'n}
\right)
e^{n\Re[S(w_c;x,s)-S(z_c;y,t)]}
\lesssim
\sum_{w_c,z_c}
e^{-\fc''n}
e^{n\Re[S(w_c;x,s)-S(z_c;y,t)]}
\end{split}\end{align}
for some \(\fc''>0\).

\smallskip
\noindent\emph{The \(\wh J^{(1)}\)- and \(\wh J^{(2)}\)-terms.}
Suppose that \((u,v)\in\fO\). If the neighborhoods associated with
\((x,s)\) and \((u,v)\) do not share a chart, then
\[
\wh J^{(1)}((x,s);(u,v))=0,
\]
and the separated-contour estimate in \Cref{l:double_bound} gives
\begin{align}\label{e:J2_separated2}
\left|
\wh J^{(2)}((x,s);(u,v))
\right|
\lesssim
\sum_{w_c,\xi_c^{\rm a}}
(\ln n)^{10}
\Delta(w_c)\Delta(\xi_c^{\rm a})
e^{n\Re[S(w_c;x,s)-S(\xi_c^{\rm a};u,v)]}.
\end{align}

If the two neighborhoods share concentric charts, then, under
\eqref{e:close} and \eqref{e:far}, the shared chart is the arctic chart.
In this case, the ascent critical point \(\xi_c^{\rm a}\) is uniformly
bounded away from \(\fW\), whereas the descent critical point
\(\xi_c^{\rm d}\) appearing in the overlap part of \(B\) lies in \(\fW\).
Thus, the strict part of \Cref{p:critical_compare} gives
\begin{align}\label{e:strict_critical_compare}
\Re S(\xi_c^{\rm d};u,v)
\leq
\Re S(\xi_c^{\rm a};u,v)-\fc'
\end{align}
for some \(\fc'>0\).

Using \eqref{e:J1bb} and \eqref{e:J2bb} for the exponentially small part of
\(B\), \eqref{e:J2_separated2} when the charts are disjoint, and
\eqref{e:strict_critical_compare} when the arctic chart is shared, we obtain
\begin{align}\begin{split}
&\phantom{{}={}}\bigl(
(|\wh J^{(1)}|+|\wh J^{(2)}|)B
\bigr)((x,s);(y,t))
\\
&\lesssim
\sum_{(u,v)}
\sum_{w_c,\xi_c^{\rm a},\xi_c^{\rm d},z_c}
(\ln n)^2
e^{n\Re[S(w_c;x,s)-S(\xi_c^{\rm a};u,v)]}
e^{-\fc n}
\Delta(\xi_c^{\rm d})\Delta(z_c)
e^{n\Re[S(\xi_c^{\rm d};u,v)-S(z_c;y,t)]}
\\
&\quad+
\sum_{\substack{(u,v)\in\fO\\
\text{disjoint chart}}}
\sum_{w_c,\xi_c^{\rm a},\xi_c^{\rm d},z_c}
(\ln n)^{10}
\Delta(w_c)\Delta(\xi_c^{\rm a})
e^{n\Re[S(w_c;x,s)-S(\xi_c^{\rm a};u,v)]}
\frac{\Delta(\xi_c^{\rm d})\Delta(z_c)}{n}
e^{n\Re[S(\xi_c^{\rm d};u,v)-S(z_c;y,t)]}
\\
&\quad+
\sum_{\substack{(u,v)\in\fO\\
\text{concentric charts}}}
\sum_{w_c,\xi_c^{\rm a},\xi_c^{\rm d},z_c}
(\ln n)^2
e^{n\Re[S(w_c;x,s)-S(\xi_c^{\rm a};u,v)]}
\frac{\Delta(\xi_c^{\rm d})\Delta(z_c)}{n}
e^{n\Re[S(\xi_c^{\rm d};u,v)-S(z_c;y,t)]}.
\end{split}\end{align}

By \Cref{p:critical_compare},
\[
\Re S(\xi_c^{\rm d};u,v)
\leq
\Re S(\xi_c^{\rm a};u,v).
\]
Thus, the first term is exponentially small. For the second term, we use
\begin{align}\label{e:sumoveruv_copy}
\sum_{(u,v)\in\fO}
\sum_{\xi_c^{\rm a},\xi_c^{\rm d}}
\frac{
\Delta(\xi_c^{\rm a})\Delta(\xi_c^{\rm d})
}{n}
\lesssim1,
\end{align}
which follows from \eqref{e:sumoveruv}. For the third term,
\eqref{e:strict_critical_compare} gives an additional factor
\(e^{-\fc'n}\). Consequently,
\begin{align}\label{e:J1J2B}
&\bigl(
(|\wh J^{(1)}|+|\wh J^{(2)}|)B
\bigr)((x,s);(y,t))\lesssim
\sum_{w_c,z_c}
\left(
(\ln n)^{10}\Delta(w_c)\Delta(z_c)
+
e^{-\fc''n}
\right)
e^{n\Re[S(w_c;x,s)-S(z_c;y,t)]}
\end{align}
for some \(\fc''>0\).

\smallskip
\noindent\emph{The \(B^2\)-term.}
By \Cref{p:B_est},
\[
B^2((x,s);(y,t))
\lesssim
B((x,s);(y,t)).
\]
Moreover, by \eqref{e:far} and the defining property of
\(\fA_{\rm ex}\), the point \((x,s)\) is bounded away from \(\fO\).
Therefore, only the exponentially small part of \(B((x,s);(y,t))\)
contributes, and
\begin{align}
B^2((x,s);(y,t))
\lesssim
B((x,s);(y,t))\lesssim
\sum_{w_c,z_c}
e^{-\fc n}
\Delta(w_c)\Delta(z_c)
e^{n\Re[S(w_c;x,s)-S(z_c;y,t)]}.
\end{align}

Combining \eqref{e:AB_decomposition}, \eqref{e:J0B},
\eqref{e:J1J2B}, and the preceding \(B^2\)-estimate gives
\begin{align}\label{e:AB_liquid_arctic_bound}
(|A|B)((x,s);(y,t))
\lesssim
\sum_{w_c,z_c}
\left(
(\ln n)^{10}\Delta(w_c)\Delta(z_c)
+
e^{-\fc'n}
\right)
e^{n\Re[S(w_c;x,s)-S(z_c;y,t)]}.
\end{align}

\medskip
\noindent\emph{Step 4: Proof of \eqref{e:K_liquid_arctic_bound}.}
By \eqref{e:Kexp}, \eqref{e:Aexp}, and \eqref{e:Aexp2},
\begin{align}
&K^{-1}((x,s);(y,t))-A_\al((x,s);(y,t))
\notag\\
&\qquad=
\OO\left(
(\ln n)^{10}B((x,s);(y,t))
+
\frac{(\ln n)^{10}}{n}
(|A|B)((x,s);(y,t))
\right).
\end{align}
The first term on the right-hand side is exponentially small because
\((x,s)\) is bounded away from \(\fO\). Substituting
\eqref{e:AB_liquid_arctic_bound} into the second term and absorbing all
exponentially small contributions, using again that
\(\Delta(w_c)\) and \(\Delta(z_c)\) are bounded below by negative powers
of \(n\), gives \eqref{e:K_liquid_arctic_bound}. This completes the proof.
\end{proof}

\begin{proof}[Proof of \Cref{l:kernel_liquid}]
We prove the case in which neither \(\fN_\al\) nor \(\fN_{(y,t)}\) is a
ramification neighborhood. The cases involving ramification neighborhoods
can be proved in the same way, so we omit them.

To ease the notation, we will simply write
\(S(w)=S(w;x,s)\) and \(S(z)=S(z;y,t)\).
Since \(w_c\) and \(z_c\) are bounded away from each other, we can truncate
the descent and ascent paths \(\sfD^{\rm d}(w_c)\) and
\(\sfD^{\rm a}(z_c)\) to \(\sfS^{\rm d}(w_c)\) and
\(\sfS^{\rm a}(z_c)\), as introduced in
\Cref{c:bulk_steepest,c:arctic_steepest}, so that they are bounded away
from each other. By the same argument as in the proof of
\Cref{l:double_bound}, the truncation error is bounded by
\[
(\ln n)^2
e^{-\fc'(\ln n)^2}
e^{n\Re[S(w_c;x,s)-S(z_c;y,t)]}.
\]

For liquid and arctic charts, \(w,z\) are bounded away from
\(x,x-s,y,y-t\). Then, by \Cref{l:PIQI_bound}, we have
\begin{align}\begin{split}
P_{ns}(nw,nx)I(w)
&=
\frac{\sqrt{s}}
{\sqrt{2\pi n}\sqrt{x-w}\sqrt{w-(x-s)}}
e^{nS(w)+\OO(1/n)},\\
Q_{nt}(nz,ny)I^{-1}(z)
&=
\frac{\sqrt{2\pi t}}
{\sqrt{n}\sqrt{y-z}\sqrt{z-(y-t)}}
e^{-nS(z)+\OO(1/n)}.
\end{split}\end{align}
Thus, the product of these two factors can be written as
\[
\frac{1}{n}e^{n(S(w)-S(z))+\OO(1/n)}\,G(w,z),
\]
where
\[
G(w,z):=
\frac{\sqrt{st}}
{\sqrt{x-w}\sqrt{w-(x-s)}\sqrt{y-z}\sqrt{z-(y-t)}}
\frac{\sqrt{\phi'(w)}\sqrt{\phi'(z)}}{\phi(w)-\phi(z)}.
\]

Since \(w_c\) is a critical point, we have \(S'(w_c)=0\). Hence, in a
neighborhood of \(w_c\),
\begin{align}\begin{split}\label{e:S_expand_w}
S(w)
&=
S(w_c)
+\frac12 S''(w_c)(w-w_c)^2
+\frac16 S'''(w_c)(w-w_c)^3
+\OO(|w-w_c|^4).
\end{split}\end{align}
An analogous expansion holds in a small neighborhood of \(z_c\).
Likewise, since all prefactors are analytic near \(w_c\) and \(z_c\),
\begin{equation}\label{e:G_expand_prop}
G(w,z)
=
G(w_c,z_c)
+G_1(w-w_c,z-z_c)
+G_2(w-w_c,z-z_c)
+\OO\!\bigl(|w-w_c|^3+|z-z_c|^3\bigr),
\end{equation}
where \(G_1\) is homogeneous of degree \(1\) and \(G_2\) is homogeneous
of degree \(2\) in \(w-w_c\) and \(z-z_c\).

If \(\sfD^{\rm d}(w_c)\) is contained in a liquid chart, then by \Cref{c:bulk} and \Cref{c:bulk_steepest}
\[
|S''(w_c;x,s)|\asymp 1,
\qquad
\operatorname{length}(\sfS^{\rm d}(w_c))
\lesssim \frac{\ln n}{\sqrt n}.
\]
If it is contained in an arctic chart, then by \Cref{c:arctic_critical} and \Cref{c:arctic_steepest}
\[
|S''(w_c;x,s)|
\asymp
\dist((x,s),\fA)^{1/2}
\gtrsim
n^{-\delta/2},
\]
and
\[
\operatorname{length}(\sfS^{\rm d}(w_c))
\asymp
\frac{\ln n}
{n^{1/2}\dist((x,s),\fA)^{1/4}}
\lesssim
\frac{\ln n}{n^{1/2-\delta/4}}.
\]
The analogous estimates hold for \(S''(z_c;y,t)\) and
\(\sfS^{\rm a}(z_c)\). In particular,
\[
n|S''(w_c;x,s)|^3
\gtrsim
n^{1-3\delta/2},
\qquad
n|S''(z_c;y,t)|^3
\gtrsim
n^{1-3\delta/2}.
\]

We rescale
\[
w=w_c+\frac{u}{\sqrt n},
\qquad
z=z_c+\frac{v}{\sqrt n}.
\]
Then, using \(\rd w=\rd u/\sqrt n\), \(\rd z=\rd v/\sqrt n\),
\eqref{e:S_expand_w}, and \eqref{e:G_expand_prop}, we obtain
\begin{align}\begin{split}\label{e:scaled_integral_prop}
\eqref{e:asymp}
&=
\frac{G(w_c,z_c)e^{\,n(S(w_c)-S(z_c))}}
{(2\pi\ri)^2n}
\int_{\sqrt n(\sfS^{\rm a}(z_c)-z_c)}
\rd v
\int_{\sqrt n(\sfS^{\rm d}(w_c)-w_c)}
\rd u\,
e^{\frac12S''(w_c)u^2-\frac12S''(z_c)v^2}\\
&\qquad\times
\left(
1+\frac1{\sqrt n}\mathcal O_{\mathrm{odd}}(u,v)
+\frac1n\OO\bigl((1+|u|+|v|)^6\bigr)
\right),
\end{split}\end{align}
where every monomial in
\(\mathcal O_{\mathrm{odd}}(u,v)\) is odd in at least one of \(u\) and
\(v\). After extending the rescaled contours to the corresponding full
steepest quadratic contours, with an exponentially small error, the Gaussian
weight
\[
e^{\frac12S''(w_c)u^2-\frac12S''(z_c)v^2}
\]
is even in \(u\) and \(v\). Therefore, every term in the coefficient of
\(n^{-1/2}\) integrates to zero. The Gaussian moments show that the first
nonvanishing correction is bounded by
\[
\OO\left(
\frac1n
+\frac{1}{n|S''(w_c;x,s)|^3}
+\frac{1}{n|S''(z_c;y,t)|^3}
\right)
=
\OO(n^{-1+3\delta/2})
=
\OO(n^{-1+2\delta}).
\]
Consequently,
\begin{align}\begin{split}\label{e:main_before_gauss_prop}
 \eqref{e:asymp}
&=
\frac{e^{\,n(S(w_c)-S(z_c))}}{(2\pi\ri)^2n}\,
G(w_c,z_c)\\
&\times
\int_{\sqrt n(\sfS^{\rm a}(z_c)-z_c)}
\rd v
\int_{\sqrt n(\sfS^{\rm d}(w_c)-w_c)}
\rd u\,
e^{\frac12S''(w_c)u^2-\frac12S''(z_c)v^2}
\left(1+\OO(n^{-1+2\delta})\right).
\end{split}\end{align}

We can then evaluate the Gaussian integrals along the steepest contours and
conclude that
\begin{align}\begin{split}\label{e:after_gauss_prop}
 \eqref{e:asymp}
=
-\frac{e^{\,n(S(w_c)-S(z_c))}}{2\pi n}\,
\frac{G(w_c,z_c)}
{\sqrt{-S''(w_c)}\sqrt{S''(z_c)}}
\left(1+\OO(n^{-1+2\delta})\right).
\end{split}\end{align}

We now simplify the prefactor. Using the critical-point relations
\[
w_c=x-s\chi(w_c),
\qquad
z_c=y-t\chi(z_c),
\]
we obtain
\begin{align}\label{e:G_critical}
G(w_c,z_c)
&=
\frac{\sqrt{\phi'(w_c)}\sqrt{\phi'(z_c)}}
{\phi(w_c)-\phi(z_c)}
\frac{1}{\sqrt{st}}
\frac{1}
{\sqrt{\chi(w_c)}\sqrt{1-\chi(w_c)}
 \sqrt{\chi(z_c)}\sqrt{1-\chi(z_c)}}.
\end{align}

By the square-root conventions in
\eqref{e:Sdirection}, \eqref{eq:integrandlocal}, and \eqref{e:vdlow}, we
have
\begin{align}\label{e:descent_saddle_factor}
\sqrt{\frac{1}{-S''(w_c)}}
&=
\frac{\sqrt{s}\sqrt{\chi(w_c)}\sqrt{1-\chi(w_c)}}
{\sqrt{\phi'(w_c)}}
\begin{cases}
\sqrt{\partial_x\phi(x,s)},
& w_c\in\bC_+,\\[1mm]
-\overline{\sqrt{\partial_x\phi(x,s)}},
& w_c\in\bC_-.
\end{cases}
\end{align}
Similarly, by \eqref{e:va},
\begin{align}\label{e:ascent_saddle_factor}
\sqrt{\frac{1}{S''(z_c)}}
&=
\frac{\sqrt{t}\sqrt{\chi(z_c)}\sqrt{1-\chi(z_c)}}
{\sqrt{\phi'(z_c)}}
\begin{cases}
\dfrac{\sqrt{\partial_x\phi(y,t)}}{\ri},
& z_c\in\bC_+,\\[3mm]
\dfrac{\overline{\sqrt{\partial_x\phi(y,t)}}}{\ri},
& z_c\in\bC_-.
\end{cases}
\end{align}

Combining \eqref{e:G_critical}--\eqref{e:ascent_saddle_factor}, all factors
involving \(s,t,\chi\), and \(\phi'\) cancel, and we obtain
\begin{align}\label{e:G_saddle_simplified}
&\frac{G(w_c,z_c)}
{\sqrt{-S''(w_c)}\sqrt{S''(z_c)}}
=
\frac{1}{\phi(w_c)-\phi(z_c)}
\begin{cases}
\dfrac{
\sqrt{\partial_x\phi(x,s)}
\sqrt{\partial_x\phi(y,t)}
}{\ri},
& w_c,z_c\in\bC_+,\\[3mm]
\dfrac{
\sqrt{\partial_x\phi(x,s)}
\overline{\sqrt{\partial_x\phi(y,t)}}
}{\ri},
& w_c\in\bC_+,\ z_c\in\bC_-,\\[3mm]
-\dfrac{
\overline{\sqrt{\partial_x\phi(x,s)}}
\sqrt{\partial_x\phi(y,t)}
}{\ri},
& w_c\in\bC_-,\ z_c\in\bC_+,\\[3mm]
-\dfrac{
\overline{\sqrt{\partial_x\phi(x,s)}}
\overline{\sqrt{\partial_x\phi(y,t)}}
}{\ri},
& w_c,z_c\in\bC_-.
\end{cases}
\end{align}
Substituting this into \eqref{e:after_gauss_prop} and using
\[
\phi(w_c)=\phi(x,s),
\qquad
\phi(z_c)=\phi(y,t),
\]
when \(w_c,z_c\in\bC_+\), we obtain \eqref{e:asymp2}.
\end{proof}

\subsection{Frozen case}

\begin{lemma}\label{l:frozenregion}
For any lattice points
$
(x,s),(y,t)\in\fP\cap\bZ^2/n
$ we have
\begin{align}\label{e:Aal_frozen_bound}
A_\al((x,s),(y,t))=\cI((x,s),(y,t))J_\fT((x,s),(y,t))+\OO\left(\sum_{w_c,z_c}
(\ln n)^{2}
e^{n\Re[S(w_c;x,s)-S(z_c;y,t)]}\right)
\end{align}
and 
\begin{align}\label{e:K_frozen_bound}
K^{-1}((x,s),(y,t))
=
A_\al((x,s),(y,t))
+
\OO\left(
\sum_{w_c,z_c}
(\ln n)^{12}
e^{n\Re[S(w_c;x,s)-S(z_c;y,t)]}
\right),
\end{align}
where \(w_c\) ranges over
\(\operatorname{Crit}^{\rm d}(x,s)\), and \(z_c\) ranges over
\(\operatorname{Crit}^{\rm a}(y,t)\).
\end{lemma}

\begin{proof}
The estimate \eqref{e:Aal_frozen_bound} follows by applying the
bounds \eqref{e:J1bb} and \eqref{e:J2bb} with \((u,v)=(y,t)\).

Since the localization length \(\Delta(\cdot)\leq1\) (recall from \Cref{def:localization_length}), we have the crude bound from the definition \eqref{e:intro_B}
\begin{align}\label{e:frozen_B_crude}
B((u,v);(y,t))
\lesssim
\sum_{\xi_c^{\rm d},z_c}
\frac{1}{n}
e^{n\Re[S(\xi_c^{\rm d};u,v)-S(z_c;y,t)]}.
\end{align}
where the sum is over
$
\xi_c^{\rm d}
\in
\operatorname{Crit}^{\rm d}(u,v^-),
$
and 
$
z_c\in
\operatorname{Crit}^{\rm a}(y,t).
$

Using \eqref{e:J0bb} and \eqref{e:frozen_B_crude}, we
obtain
\begin{align}
\bigl(|\wh J^{(0)}|B\bigr)((x,s),(y,t))
&\lesssim
\sum_{(u,v)}
\sum_{w_c,\xi_c^{\rm d},z_c}
(\ln n)
e^{n\Re[S(w_c;x,s)-S(\xi_c^{\rm d};u,v)]}
\times\frac{1}{n}
e^{n\Re[S(\xi_c^{\rm d};u,v)-S(z_c;y,t)]}
\notag\\
&\lesssim
\sum_{w_c,z_c}
n\ln n\,
e^{n\Re[S(w_c;x,s)-S(z_c;y,t)]}.
\end{align}
where we used that the number of lattice points \((u,v)\in\fP\cap\bZ^2/n\) is \(\OO(n^2)\),
and the numbers of associated critical points are uniformly bounded.

Similarly, by \Cref{p:critical_compare},
\[
\Re S(\xi_c^{\rm d};u,v)
\leq
\Re S(\xi_c^{\rm a};u,v).
\]
Therefore, using \eqref{e:J1bb}, \eqref{e:J2bb}  and \eqref{e:frozen_B_crude}, we obtain
\begin{align}
\bigl(|\wh J^{(1)}|B\bigr)((x,s),(y,t))
&\lesssim
\sum_{w_c,z_c}
n\ln n\,
e^{n\Re[S(w_c;x,s)-S(z_c;y,t)]},
\end{align}
and
\begin{align}
\bigl(|\wh J^{(2)}|B\bigr)((x,s),(y,t))
&\lesssim
\sum_{w_c,z_c}
n(\ln n)^2
e^{n\Re[S(w_c;x,s)-S(z_c;y,t)]}.
\end{align}

By \Cref{p:B_est} and \eqref{e:frozen_B_crude}, we have
\begin{align}
B^2((x,s),(y,t))
\lesssim
B((x,s),(y,t))
\lesssim
\sum_{w_c,z_c}
\frac1n
e^{n\Re[S(w_c;x,s)-S(z_c;y,t)]}.
\end{align}
It follows that
\begin{align}
(|A|B)((x,s),(y,t))
\lesssim
\sum_{w_c,z_c}
n(\ln n)^2
e^{n\Re[S(w_c;x,s)-S(z_c;y,t)]}.
\end{align}

Finally, using \Cref{l:global_approximation}
\[
K^{-1}
=
A
+
\OO\left(\frac{(\ln n)^{10}}{n}(|A|B)\right)
\]
and \eqref{e:Aexp}, 
we obtain
\begin{align}\begin{split}
K^{-1}((x,s),(y,t))
&=
A_\al((x,s),(y,t))
+
\OO\left((\ln n)^{10}B((x,s),(y,t))
+
\frac{(\ln n)^{10}}{n}(|A|B)((x,s),(y,t))\right)
\\
&=
A_\al((x,s),(y,t))
+
\OO\left(
\sum_{w_c,z_c}
(\ln n)^{12}
e^{n\Re[S(w_c;x,s)-S(z_c;y,t)]}
\right).
\end{split}\end{align}
This proves the claim \eqref{e:K_frozen_bound}.
\end{proof}

\subsection{Proofs of \Cref{p:liquid,p:arctic,p:frozen}}

\begin{proof}[Proof of \Cref{p:liquid}]
The claim follows from \Cref{l:kernel_liquid,l:liquidregion}, since the
error term in \Cref{l:liquidregion} is
\(\OO(n^{-1+2\delta})\) times the leading term and can therefore be
absorbed into the stated error.
\end{proof}

\begin{proof}[Proof of \Cref{p:arctic}]
The claim follows directly from
\Cref{l:kernel_liquid,l:liquidregion}.
\end{proof}

\begin{proof}[Proof of \Cref{p:frozen}]
Suppose first that \((x,s)\) lies in a frozen neighborhood. Then the
relevant descent critical point \(w_c\) and ascent critical point \(z_c\)
are uniformly bounded away from each other. Hence the strict part of
\Cref{p:critical_compare} gives
\begin{align}\label{e:agap1}
\Re S(z_c;x,s)
\geq
\Re S(w_c;x,s)+\fc'
\end{align}
for some \(\fc'>0\).

Suppose instead that \((x,s)\) lies in a regular arctic neighborhood and
\[
\dist((x,s),\fA)\geq n^{-\delta}.
\]
Then \Cref{c:arctic_steepest} implies that
\begin{align}\label{e:agap2}
\Re S(z_c;x,s)-\Re S(w_c;x,s)
&\gtrsim
\dist((x,s),\fA)^{3/2}
\geq
n^{-3\delta/2}.
\end{align}

By the continuity estimate in
\Cref{p:white_triangle_critical_value}, replacing
\(\Re S(w_c;x,s)\) by either
\[
\Re S\left(
w_c\left(x-\frac1n,s-\frac1n\right);
x-\frac1n,s-\frac1n
\right)
\]
or
\[
\Re S\left(
w_c\left(x,s-\frac1n\right);
x,s-\frac1n
\right)
\]
changes the critical value by at most
$\OO(\ln n/n)$.
For \(\delta>0\) sufficiently small, this error is
\(\oo(n^{-3\delta/2})\). Thus, the preceding
action-gap estimates \eqref{e:agap1} and \eqref{e:agap2} remain valid for all three adjacent blue triangles.

It follows from \Cref{l:frozenregion} that for $(x',s')\in \{(x,s), (x-1/n, s-1/n), (x,s-1/n)\}$
\begin{align}
K^{-1}((x',s'),(x,s))
&=
\cI(x',s';x',s')J_{\fT}(x',s';x,s)
+
\OO\left(e^{-n^{1-2\delta}}\right).
\end{align}

On \(\fT\), one has
$
\nabla H^*\equiv(0,0).
$
Thus, by \eqref{e:gH00},
\begin{align}
J_{\fT}(x,s;x,s)
=1,\quad 
J_{\fT}\left(x,s-\frac1n;x,s\right)
=
J_{\fT}\left(x-\frac1n,s-\frac1n;x,s\right)
=0
\end{align}
We therefore conclude that
\begin{align}
&K^{-1}((x,s),(x,s))=1+\OO(
e^{-n^{1-2\delta}}),\\
&K^{-1}((x,s-1/n),(x,s)), K^{-1}((x-1/n,s-1/n),(x,s))=\OO(
e^{-n^{1-2\delta}})
\end{align}
\end{proof}

\chapter{Appendices}

\section{Structure of Frozen Region}\label{s:frozen_structure}
\subsection{Limit shapes and arctic boundary} 
	\label{s:limits}
 We introduce the complex coordinates and complex derivatives
 \begin{align}
z=x+\ri s,\quad  \del_{\overline z} =\frac{\del_x +\ri \del_s}{2}, \quad \del_{z} =\frac{\del_x -\ri \del_s}{2}
 \end{align}
and recall the function $g: \fL \mapsto \bD$ from \eqref{e:defg}. %Moreover, for a sequence of points $u_k\in \fL$ converging to $(x,t)\in \del \fL$ with $0<t<1$, we have 
 %$f(u_k)\rightarrow \bR$ and $|g(u_k)|\rightarrow 1$. Similarly for a sequence of points $u_k\in \fL$ converging to $(a_i,0)\in \del \fL$ or $(b_i,1)\in \del \fL$, we have
% $\rho(u_k)\rightarrow \infty$, $f(u_k)\rightarrow \infty$ and we also have $|g(u_k)|\rightarrow 1$. Thus $g$ is a proper map from $\fL$ to $\bD$. 
We notice
 \begin{align}
  \del_{\overline z} f= \frac{\del_x f +\ri \del_s f}{2}=\frac{\del_x f -\ri (f+1)^{-1}f\del_x f}{2},\quad 
  \del_{z} f =\frac{\del_x f -\ri \del_s f}{2}=\frac{\del_x f +\ri (f+1)^{-1}f \del_x f }{2}
 \end{align}
 It follows that
  \begin{align}\label{e:Beltrami2}
   \del_{\overline z} f(z)=\frac{1-\ri (f(x,s)+1)^{-1}f(x,s)}{1+\ri (f(x,s)+1)^{-1}f(x,s)}\del_{z} f(z)=  g(z) \del_{z} f(z). 
 \end{align}
By \eqref{e:Beltrami2} and chain rule, we conclude that $g(z)$ satisfies the following universal Beltrami equation
 \begin{align}\label{e:Beltrami}
   \del_{\overline z} g(z)= g(z) \del_{z} g(z). 
 \end{align}

 The regularity properties of the universal Beltrami equation \eqref{e:Beltrami} has been intensively studied in \cite[Section 6]{astala2026dimer}, we record them below.

\begin{theorem}[{\cite{astala2026dimer}}] \label{t:g_behavior}
 Under the notation of \Cref{p},  the following holds for the function $g : \fL \to \mathbb{D}$ defined in \eqref{e:defg}. 
 \begin{enumerate}
\item The arctic curve $\fA=\del \fL$ is the real locus of an algebraic curve. It has finitely many singularities, and these are either first order (inward) cusps or tacnodes (double points). 

\item The map $g : \fL \to \mathbb{D}$ extends continuously to the closure $\overline{\fL}$.  The tangent vectors $\tau(\zeta)$ of the arctic curve and the boundary values $g(\zeta) \in \partial \mathbb{D}$ are related via the identity
  \[
    g(\zeta) = -\tau(\zeta)^2, \qquad \zeta \in \partial \fL \setminus \{\text{cusps}\}.
  \]
 At a cusp point $\zeta_0$,  as $\zeta\rightarrow \zeta_0$ along $\del \fL$ the tangent vectors $\tau(\zeta)$ of the arctic curve have a well defined limit, which is the direction of the cusp.

  \item 
 $\partial \fL$ is locally strictly convex and smooth, except at $d-2$ cusps, where $d$ is the degree of $g:\del \fL\mapsto \del \mathbb D$. More precisely for every $\zeta \in \partial \fL$ outside the cusps and tacnodes, $B_\varepsilon(\zeta) \cap \fL$ is strictly convex for $\varepsilon > 0$ small enough. At the tacnodes $\zeta \in \partial \fL$ the set $B_\varepsilon(\zeta) \cap \fL$ has two components, both convex.

 \item For every $\zeta \in \partial \fL$ outside the tacnodes, the map $(x,t)\mapsto g(x,t)$ is injective for $B_\varepsilon(\zeta) \cap \fL$ for $\varepsilon>0$ small enough. At the tacnodes $\zeta \in \partial \fL$ the set $B_\varepsilon(\zeta) \cap \fL$ has two components, and $g(x,t)$ is injective on each of them.

\end{enumerate}
\end{theorem}
\begin{proof}
By the first statement of \Cref{pla}, the arctic boundary $\del \fL$ is frozen. Then it follows from \cite[Theorem 2.6]{astala2026dimer} that $g : \fL \to \mathbb{D}$ is proper.
The first statement follows from \cite[Theorem 6.3]{astala2026dimer}; The second statement follows from \cite[Theorem 6.6, Theorem 5.16]{astala2026dimer}; The third statement is given in \cite[Theorem 6.6]{astala2026dimer}. The last statement follows from the decomposition of $g$ as in \cite[Theorem 6.1]{astala2026dimer}: $g=B\circ \widetilde g$, where $\widetilde g$ is a homeomorphism from $\fL$ to $\bD$, and $B$ is a finite Blaschke product that is locally biholomorphic in a neighborhood of each boundary point of $\bD$.

\end{proof}

It follows from the last statement of \Cref{t:g_behavior} that $g$ extends continuously to $\overline{\fL}$. Recall the relation between $f$ and $g$ from \eqref{e:defg}, we also have that $f$ extends to a continuous map from $\overline \fL$ to $\bC_-\cup \bR\cup \{\infty\}$.

 We introduce the map
\begin{align}\label{e:defchi}
\chi(x,s)=\frac{f(x,s)}{f(x,s)+1}
\end{align}
  	
\begin{proof}[Proof of \Cref{pa1}]
It follows from the last statement of \Cref{t:g_behavior} that $g$ extends continuously to $\overline{\fL}$. Recall the relation between $f$ and $g$ from \eqref{e:defg}, we also have that $f$ extends to a continuous map from $\overline \fL$ to $\bC_-\cup \bR\cup \{\infty\}$.  
 
 It follows from the second statement of \Cref{t:g_behavior} that the tangent vectors $\tau(\zeta)$ of the arctic boundary satisfies \begin{align}
\tau^2(\zeta)=- g(\zeta)=\frac{\ri (f(\zeta)+1)^{-1}f(\zeta)-1}{1+\ri (f(\zeta)+1)^{-1} f(\zeta)}=\frac{ (f(\zeta)+1)^{-1}f(\zeta)+\ri }{ (f(\zeta)+1)^{-1}f(\zeta)-\ri}
 \end{align}
 Thus the slope of $\tau(\zeta)$ is given by $(f(\zeta)+1)/f(\zeta)$. This gives the first statement in \Cref{pa1}.

Since \( \mathfrak{L} \) is simply connected, the arctic boundary \( \mathfrak{A} \) is connected and homeomorphic to a circle.  
Next, we investigate $f(x,s)$ as \( (x, s) \) moves along the arctic curve \( \mathfrak{A} \) in counter-clockwise order. %, $f(x,t)$ decreases until it diverges to $-\infty$, then wraps around to $\infty$ and continuous decreasing.  
The third statement in \Cref{t:g_behavior} asserts that \( \partial \fL \) is locally strictly convex. As \( (x, s) \) traverses the arctic curve \( \mathfrak{A} \) in counterclockwise order, the slope of the tangent vector, given by \((f(x, s) + 1)/{f(x, s)} \), increases monotonically until it diverges to 
$+\infty$, then wraps around to $-\infty$ and continues increasing. It follows that the complex slope map \( f: \fA \to \mathbb{R} \cup \{\infty\} \) winds around \( \mathbb{R} \cup \{\infty\} \) exactly \( d \) times in the negative orientation.

For \Cref{i:no_tacnode}, the first statement that $\fA$ has $d-2$ inward-pointing cusps follows from the third statement in \Cref{t:g_behavior}. We postpone the proof of the third statement that $\fA$ has no tacnodes  to \Cref{s:frozen}.

%Next we show that $\fA$ does not have tacnodes (double points). Otherwise the tacnode $\zeta$ is where  two segments $\Gamma, \Gamma'$ of arctic boundary are tangent to each other. If $\Gamma$ and $\Gamma'$ are different segments, then they separate the liquid region into disconnected components, which contradicts to that $\fL$ is connected. If $\Gamma=\Gamma'$ is the same segment, then $\Gamma$ contains loop from $\zeta$ to itself. And some point on the loop has slope $0$. This is impossible because only points with slope $0$ are given by these boundary point \eqref{e:bbpoint}.

%  By the third statement of \Cref{t:g_behavior}, the set $B(\zeta,\varepsilon)\cap \fL$ has two components. Their boundaries give by two pieces $\Gamma_1, \Gamma_2$ of the arctic curve, which are tangent to each other at the tacnode $\zeta$. Since when we move along the arctic curve, $f(x,t)$ decreases monotonically

\end{proof}

\subsection{Properties of the Riemann surface $Q$}
\label{s:Rsurface}

	\begin{proof}[Proof of \Cref{p:surface}]
	We recall from \eqref{e:defchi} and \eqref{e:emb} that for $(x,s)\in \fL$,
	\begin{align}
	\Im[\chi(x,s)]=\frac{\Im[f(x,s)]}{|f(x,s)+1|^2}<0, \quad \Im[z(x,s)]=\Im[x-s\chi(x,s)]=-s\Im[\chi(x,s)].
	\end{align}
	Hence, the map \eqref{e:emb} restricted to $\fL$ is an injection to its image, since we can recover $(x,s)$ from its image as
\begin{align}\label{e:txconstruct}
s=-\frac{{\rm Im}[z(x,s)]}{{\rm Im}\left[\chi(x,s)\right]},\quad x=z(x,s)+s \chi(x,s).
\end{align}
If $\del_x \chi(x_0,s_0)\neq 1/s_0$, then $\del_xz(x_0,s_0) =1-s_0\del_x \chi(x_0,s_0)\neq 0$, and $\del_s z(x,s)=-\chi(x,s)-s\del_s \chi(x,s)=-\chi(x,s)(1-s\del_x \chi(x,s))$
which follows from the complex Burgers equation \eqref{ftx} that
\begin{align}\label{ftx2}
\del_s \chi(x,s)=\frac{\del_s f(x,s)}{(f(x,s)+1)^2}=-\frac{f(x,s)\del_s f(x,s)}{(f(x,s)+1)^3}=-\chi(x,s)\del_x \chi(x,s).\end{align}
 Thus 
\begin{align}
\label{e:dtzdxz}
\frac{\del_s z(x_0,s_0)}{\del_x z(x_0,s_0)}=-\chi(x_0,s_0)\not\in \bR.
\end{align}
And the second coordinate of the map \eqref{e:emb}, $(x,s)\mapsto z(x,s)$ is a local homeomorphism and there exists a function $f_0$ such that $f(x,s)=f_0(z(x,s))$ in a small neighborhood of $(x_0,s_0)$ in $\fL$. 
By taking derivatives with respect to $x,s$ on both sides of $f(x,s)=f_0(z(x,s))$, 
\begin{align}\label{e:dxfdtf1}
\del_x f(x,s)&=\del_z f_0 (z(x,s)) \del_x z(x,s)+\del_{\overline z} f_0 (z(x,s)) \del_x \overline{z}(x,s)\\
\begin{split}\label{e:dxfdtf2}
\del_s f(x,s)&=-\chi(x,s) \del_x f(x,s)=\del_z f_0 (z(x,s)) \del_s z(x,s)+\del_{\overline z} f_0 (z(x,s)) \del_s \overline{z}(x,s)\\
&=-\del_z f_0 (z(x,s))\chi(x,s) \del_x z(x,s)-\del_{\overline z} f_0 (z(x,s)) \overline{\chi(x,s)}\del_x \overline{z}(x,s),
\end{split}
\end{align}
where in the last line we used \eqref{ftx2} and \eqref{e:dtzdxz}. By taking the sum of \eqref{e:dxfdtf1} multiplying $\chi(x,t)$ and \eqref{e:dxfdtf2}, we conclude that  $f_0$ satisfies the Cauchy--Riemann equations, i.e. $\del_{\overline z}f_0=0$. Thus $f_0$ is analytic. 

Similarly if $\del_x \chi(x_0,s_0)\neq 0$ (in particularly it is possible $\del_x \chi(x_0,s_0)=1/s_0$), then $\del_x f(x_0,s_0)\neq 0$, and we have 
\begin{align}
\label{e:dtzdxf3}
\frac{\del_s f(x_0,s_0)}{\del_x f(x_0,s_0)}=-f(x_0,s_0)\not\in \bR.
\end{align}
And the first coordinate of the map \eqref{e:emb}, $(x,s)\mapsto f(x,s)$ is a local homeomorphism and there exists a function $Q_0$ such that $z(x,s)=Q_0(f(x,s))$ in a small neighborhood of $(x_0,s_0)$ in $\fL$. The same as in \eqref{e:dxfdtf1} and \eqref{e:dxfdtf2}, one can check that $Q_0$ satisfies the Cauchy--Riemann equations, and is analytic. 
It follows that the image of the map \eqref{e:emb} restricted to  $\fL$ gives a Riemann surface: $\{(f,z)=(f(x,s),z(x,s))\in \mathbb C^2: (x,s)\in \fL\}$ with local charts given by $Q_0$ or $f_0$.

% implies $\del_t f(x,t)/\del_x f(x,t)=-f(x,t)\not\in \bR$. Therefore, the first coordinate of the map \eqref{e:emb}, $(x,t)\mapsto f(x,t)$ is a local homeomorphism and there exists a function $Q_0$ such that $z(x,t)=Q_0(f(x,t))$ in a small neighborhood of $\fU\subseteq \fL$ of $(x_0,t_0)$. The complex Burgers equation \eqref{ftfx} implies that $Q_0$ satisfies the Cauchy--Riemann equations. It follows that $Q_0$ is analytic, and the image of the map \eqref{e:emb} restricted to  $\fL$ gives a Riemann surface: $\{(f,z)=(f(x,t),z(x,t))\in \mathbb C^2: (x,t)\in \fL\}$ with local charts given by $Q_0$.

%Since $\fL$ is simply connected, \cite[Theorem 5.1]{DMCS} gives a decomposition of the map $(x,t)\mapsto f(x,t)$.  It implies 

By the third statement of \Cref{pa1}, $\fA$ does not contain tacnodes. For any $(x_0,s_0)\in \fA$, there exists a small neighborhood $\fU=B_\delta(x_0,s_0)$ of it, such that the map  $(x,s)\mapsto \chi(x,s)$ is injective for $(x,s)\in \fU\cap \overline\fL$ by the fourth statement in \Cref{t:g_behavior}. Thus the same argument as above, there exists a continuous map $Q_0$ such that $z(x,s)=Q_0(\chi(x,s))$ for  $(x,s)\in \fU\cap \overline{\fL}$, and $Q_0$ is analytic on $\chi(\fU\cap \fL)$. For \( (x,s) \in \fU \cap \partial \fL \), both \( z(x,s) \) and \( \chi(x,s) \) are real-valued (possibly infinite). Hence, by the Schwarz reflection principle, the function \( Q_0 \) can be extended to a real analytic function in a neighborhood of \( \chi(x_0, s_0) \).

In a small neighborhood of $\chi(x_0,s_0)$, we rewrite \eqref{e:txconstruct} as
\begin{align}\label{e:tz}
s=-\frac{\Im[z(x,s)]}{\Im[\chi(x,s)]}=-\frac{\Im[Q_0(\chi(x,s))]}{\Im[\chi(x,s)]}.
\end{align}

We  notice that if $(x_0,s_0)\in \fA$ is a tangent location with $f(x_0, s_0)=-1$ (the slope of the tangent is $0$), then $\chi(x_0, s_0)=\infty$ and $z(x_0,s_0)=\infty$. 

In the following we first discuss the case that $f(x_0, s_0)\neq -1$. In this case $\chi(x_0, s_0)\in \bR\setminus \{0\}$, by sending $(x,s)\in \fL$ to $(x_0,s_0)$, we conclude $Q_0'(\chi(x_0,s_0))=-s_0<0$. Thus in a small neighborhood of $\chi(x_0,s_0)$, $z=Q_0(\chi)$ is invertible, we denote its inverse as $\chi=Q_0^{-1}$: $\chi(x,s)=\chi(z(x,s))$. Using $\chi(z)$, we can locally parameterize the arctic boundary in the following way: for $z\in \bR$, 
\begin{align}\label{e:chid1}
-\frac{1}{s}=\chi'(z),\quad x=z+s\chi(z)=z-\frac{\chi(z)}{\chi'(z)}.
\end{align}
We can view $x,s$ as functions of $z$, by taking more derivatives
\begin{align}\label{e:chid2}
\del_z s =\frac{\chi''(z)}{(\chi'(z))^2},\quad \del_z x=\frac{\chi(z)\chi''(z)}{(\chi'(z))^2}, 
\end{align}
If $\chi''(z_0)=0$, both derivatives above vanishes at $z_0=z(x_0,s_0)$, and $(x_0,s_0)$ corresponds to a cusp location.  If \( \chi''(z_0) \neq 0 \), we can take the ratio \eqref{e:chid2} and differentiate further to obtain
\begin{align}
\partial_s x = \chi(z), \quad \partial_s^2 x = -\frac{1}{s^3 \chi''(z)}.
\end{align}
If \( \chi''(z_0) > 0 \), then \( \partial_s^2 x < 0 \), so locally the arctic boundary is concave in the \( x \)-direction.  
If \( \chi''(z_0) < 0 \), then \( \partial_s^2 x > 0 \), and locally the arctic boundary is convex in the \( x \)-direction.

Next we study the cusp case, i.e. $\chi''(z_0)=0$. We first show that $\chi'''(z_0)\neq 0$ by contradiction. Assume that $0=\chi'''(z_0)=\cdots=\chi^{(k-1)}(z_0)$ and $\chi^{(k)}(z_0)\neq 0$ for some $k\geq 4$. We notice $k\neq \infty$, otherwise $\chi'(z)\equiv -1/s_0$ is a constant, which is impossible. 

We can reconstruct the liquid region around $(x_0, s_0)$, by inverting the map $(x,s)\mapsto z(x,s)$ in a neighborhood of $z_0=z(x_0, s_0)$ in the following way:
Let $z=z_0+a+\ri b$ for $b> 0$, we can rewrite \eqref{e:tz} as
\begin{align}\label{e:1/texp}
-\frac{1}{s}=\frac{\Im[\chi(z)]}{\Im[z]}=\frac{1}{b}\Im\left[\chi(z_0)+\chi'(z_0)(z-z_0)+\chi^{(k)}(z_0)\frac{(z-z_0)^k}{k!}+\cdots \right].
\end{align}
We can reorganize it as
\begin{align}
\frac{b}{\chi^{(k)}(z_0)}\left(\frac{1}{s_0}-\frac{1}{s}\right)=\Im\left[\frac{(z-z_0)^k}{k!}+\cdots\right]:=\frac{\Im[\Phi^k(z-z_0)]}{k!},
\end{align}
where $\Phi(z)$ is a real analytic function in a neighborhood of $0$, and $\Phi'(0)=1$.
Locally around $0$, $\Phi(w)$ behaves like $w$, there are $k-1$ branches of solutions:  $z-z_0=\Phi^{-1}(re^{i\ell \pi/k})$ for $1\leq \ell \leq k-1$ and some $r>0$ such that $s=s_0$ and $\Im[z]>0$.

For these values of $z$ such that $s=s_0$, we can Taylor expand the second relation in \eqref{e:tz} as
\begin{align}
x&=z+s_0\chi(z)=x_0+ s_0 \left(\chi^{(k)}(z_0)\frac{(z-z_0)^k}{k!}+\cdots\right)
=x_0+ s_0\chi^{(k)}(z_0)\Re \left[\frac{(z-z_0)^k}{k!}+\cdots\right]\\
&=x_0+s_0\chi^{(k)}(z_0)\frac{\Re[\Phi^k(z-z_0)]}{k!}
=x_0+s_0 \chi^{(k)}(z_0) \frac{(-1)^\ell r^k}{k!} ,\quad r>0,\quad 1\leq \ell\leq k-1
\end{align}
If $k\geq 4$, the inverse map from $z$ to $(x,s)$, covers $(0,\varepsilon)\times \{s_0\}$ or $(-\varepsilon,0)\times \{s_0\}$ at least twice. This leads to a contradiction. We conclude that $\chi'''(z_0)\neq 0$, and for $z\in [z_0-\varepsilon, z_0+\varepsilon]\subset \bR$, \eqref{e:1/texp} gives
\begin{align}
-\frac{1}{s}&=\frac{\Im[\chi(z)]}{\Im[z]}=\frac{1}{b}\Im\left[\chi(z_0)+\chi'(z_0)(z-z_0)+\chi'''(z_0)\frac{(z-z_0)^3}{3!}+\cdots \right]\\
&=-\frac{1}{s_0}+\frac{\chi'''(z_0)a^2}{2}+\OO(a^3)
\end{align}

If $\chi'''(z_0)>0$, $s$ as a function of $z\in [z_0-\varepsilon, z_0+\varepsilon]$ first decreases to $s_0$ then increases, we have a cusp pointing downward.

If $\chi'''(z_0)<0$, $s$ as a function of $z\in [z_0-\varepsilon, z_0+\varepsilon]$ first increases to $s_0$ then decreases, we have a cusp pointing upward.

%\begin{align}
%\Im\left[f'''(z_0)\frac{(z-z_0)^2}{2}+f''''(z_0)\frac{(z-z_0)^3}{3!}+\cdots \right]>0
%\end{align}

Next we discuss the cases that $(x_0, s_0)\in \fA$ is a tangent location with $f(x_0, s_0)=-1$. In this case, $\chi(x_0,s_0)=\infty$. By the same argument as for the $f(x_0, s_0)\neq -1$ case, we can extend $Q_0$ to a real analytic function in a neighborhood of $\infty$. Because $Q_0(\chi(x,s))=z(x,s)=x-s\chi(x,s)$, by sending $(x,s)\rightarrow (x_0,s_0)$, we conclude that 
\begin{align}\label{e:Q0_infty}
Q_0(w)=-s_0w +\OO(1), \text{ as } w\rightarrow\infty.
\end{align}
Thus in a neighborhood of $\infty$, $z=Q_0(\chi)$ is invertible, we denote its inverse as $\chi=Q_0^{-1}$: $\chi(x,s)=\chi(z(x,s))$. From \eqref{e:Q0_infty}, $\chi(z)$ has the following Laurent form in a neighborhood of $\infty$ 
\begin{align}\label{e:expansion}
\chi(z)=\frac{x_0-z}{s_0}+\sum_{i\geq 1}\frac{c_i}{(x_0-z)^i},
\end{align}
The same as in \eqref{e:chid1}, using $\chi(z)$, we can locally parameterize the arctic boundary in the following way: for $z\in \bR$ in a neighborhood of $\infty$,
\begin{align}\label{e:chid1_copy}
-\frac{1}{s}=\chi'(z),\quad x=z+s\chi(z)=z-\frac{\chi(z)}{\chi'(z)}.
\end{align}

In \eqref{e:expansion}, if $c_1\neq 0$, then as $z\in \bR$ and $|z|\rightarrow \infty$ the first relation of \eqref{e:chid1} gives
\begin{align}
-\frac{1}{s}=-\frac{1}{s_0}+\frac{c_1}{(x_0-z)^2} +\OO\left(\frac{1}{|x_0-z|^3}\right).
\end{align}
Thus if $c_1>0$, then $s_0<s$, and around $(x_0, s_0)$ the arctic boundary is convex in the $s$ direction. 
Similarly, if $c_1<0$, then $s_0>s$, and around $(x_0, s_0)$ the arctic boundary is concave in the $s$ direction.

Next we assume that $0=c_1=c_2=\cdots=c_{k-1}$ and $c_k\neq 0$ for some $k\geq 1$. We notice $k\neq \infty$, otherwise $\chi'(z)\equiv -1/s_0$ is a constant, which is impossible.

We can reconstruct the liquid region around $(x_0, s_0)$, by inverting the map $(x,s)\mapsto z(x,s)$ in a neighborhood of $\infty$ in the following way:
Let $z\in \bC_+$, we can rewrite \eqref{e:tz} as
\begin{align}\label{e:1/texp}
-\frac{1}{s}=\frac{\Im[\chi(z)]}{\Im[z]}=\frac{1}{\Im[z]}\Im\left[\frac{-z+x_0}{s_0}+\sum_{i\geq k}\frac{c_i}{(x_0-z)^i}\right]
=-\frac{1}{s_0}+\frac{1}{\Im[z]}\Im\left[\sum_{i\geq k}\frac{c_i}{(x_0-z)^i}\right]
.
\end{align}
We can reorganize it as
\begin{align}\label{e:phik}
\frac{\Im[z]}{c_k}\left(\frac{1}{s_0}-\frac{1}{s}\right)=\Im\left[\frac{1}{(x_0-z)^k}+\sum_{i\geq k+1}\frac{c_i/c_k}{(x_0-z)^i}\cdots\right]:=\Im[\Phi^k(1/(x_0-z))],
\end{align}
where $\Phi(z)$ is a real analytic function in a neighborhood of $0$, and $\Phi'(0)=1$.
Locally around $0$, $\Phi(w)$ behaves like $w$, and maps upper half plane to upper half plane. 

Next we use \eqref{e:phik} to solve for $z\in \bC_+$ corresponding to $\{(x, s_0)\in \fL\}$. If we take $s=s_0$, 
in \eqref{e:phik}, then $\Phi^k(1/(x_0-z))\in \bR$ and $\Phi(1/(x_0-z))\in \bC_+$. It has $k-1$ branches of solutions:  $1/(x_0-z)=\Phi^{-1}(re^{i\ell \pi/k})$ with $1\leq \ell \leq k-1$ and  $r>0$.

For these values of $z$ corresponding to $\{(x, s_0)\in \fL\}$, we can Taylor expand the second relation in \eqref{e:tz} as
\begin{align}
x&=z+s_0\chi(z)=x_0+ s_0 \left(\sum_{i\geq k}\frac{c_i}{(x_0-z)^i}\right)
=x_0+ s_0c_k \Phi^k(1/(x_0-z))=
x_0+s_0 c_k (-1)^\ell r^k,
\end{align}
for any $r>0$ and $1\leq \ell\leq k-1$.  If $k\geq 3$, the inverse map from $z$ to $(x,s)$, covers both $(x_0,x_0+\varepsilon)\times \{s_0\}$ and $(x_0-\varepsilon,x_0)\times \{s_0\}$. This leads to a contradiction that the tangent line at $(x_0, s_0)$ has slope $0$. So $k=2$ and $c_2\neq 0$. 
For $z\in  \bR$ with $|z|\rightarrow \infty$, \eqref{e:chid1_copy} gives
\begin{align}
x=z-\frac{\chi(z)}{\chi'(z)}
=z-\frac{(x_0-z)/s_0+c_2/(x_0-z)^2+\OO(1/|(x_0-z)|^3)}{-1/s_0+2c_2/(x_0-z)^3+\OO(1/|x_0-z|^4)}=x_0+\frac{3c_2 s_0}{(x_0-z)^2}+\OO\left(\frac{1}{|x_0-z|^3}\right)
\end{align}
We conclude that if $c_2>0$, $x\geq x_0$ and we have a cusp pointing leftward.
If $c_2<0$, $x\leq  x_0$ and we have a cusp pointing rightward.

Since the image of the arctic boundary $\fA$ under the map \eqref{e:emb} is real, we can glue the image of the map \eqref{e:emb} restricted to the liquid region $\fL$, and its complex conjugate along the image of the arctic boundary to get an immersed Riemann surface in $\mathbb{CP}^2$:
$\{(f,z)=(f(x,s),z(x,s))\in \mathbb{CP}^2: (x,s)\in \overline\fL\}\cup\{(f,z)=(\overline{f(x,s)},\overline{z(x,s)})\in \mathbb{CP}^2: (x,s)\in \overline\fL\}$. Since $Q_0$ is real analytic, the local charts around any point $(f(x_0,s_0),z(x_0,s_0))$  for $(x_0,s_0)\in \fA$ on the arctic boundary are given by $z=Q_0(f/(f+1))$.
$f(x,s)$ with its complex conjugate  and $z(x,s)$ with its complex conjugate are two meromorphic function on this Riemann surface. By \cite[Theorem 5.8.1]{jost2013compact},  there exists a rational function $Q$ such that $Q(f(x,s), z(x,s))=0$ for any $(x,s)\in \overline{\fL}$, which gives \eqref{e:qfh}.

%{\color{red}The fourth statement follows from the facts that $(x, t) \in \mathfrak{A}$ if and only if $f_t (x) \in \mathbb{R}$, by \eqref{fh}, and that any root of $Q_0$ is real if and only if it is a double root, as $Q_0$ is real anaytic (see also the discussion at the end of \cite[Section 1.6]{LSCE})}.

\end{proof}

\begin{proof}[Proof of \Cref{l:derivechi}]
Since \(x=z+s\chi(z)\) with \(x\in\bR\), taking imaginary parts gives
\begin{equation}\label{e:Im_identity}
0=\Im z(x,s)+s\,\Im\chi(z(x,s))
\qquad\Longrightarrow\qquad
\frac1s=-\frac{\Im\chi(z(x,s))}{\Im z(x,s)}.
\end{equation}
Differentiating \(x=z+s\chi(z)\) with respect to \(x\) yields
\begin{equation}\label{e:dxz}
1=\partial_x z(x,s)\,\bigl(1+s\chi'(z(x,s))\bigr)
\qquad\Longrightarrow\qquad
\partial_x z(x,s)=\frac{1}{s}\,\frac{1}{1/s+\chi'(z(x,s))}.
\end{equation}
We now expand the denominator in \eqref{e:dxz}.

If $(x_0, s_0)$ is away from horizontal tangencies, using \eqref{e:Im_identity} we rewrite
\begin{equation}\label{e:denom_rewrite}
\frac1s+\chi'(z)=\chi'(z)-\frac{\Im\chi(z)}{\Im z}
=\chi'(z)-\frac{\Im[\chi(z)-\chi(z_0)]}{b},
\end{equation}
since \(\chi(z_0)\in\bR\) and \(\Im z=b\).
Taylor-expanding \(\chi\) at \(z_0\) gives
\[
\chi'(z)=\sum_{k\ge 1}\frac{\chi^{(k)}(z_0)}{(k-1)!}(a+\ri b)^{k-1},
\qquad
\chi(z)-\chi(z_0)=\sum_{k\ge 1}\frac{\chi^{(k)}(z_0)}{k!}(a+\ri b)^{k}.
\]
Substituting into \eqref{e:denom_rewrite} yields the exact series identity
\begin{equation}\label{e:series_identity}
\frac1s+\chi'(z)
=\sum_{k\ge 1}\frac{\chi^{(k)}(z_0)}{k!}\Bigl(k(a+\ri b)^{k-1}-\frac{\Im[(a+\ri b)^{k}]}{b}\Bigr).
\end{equation}
If \((x_0,s_0)\) is not a cusp location, then \(\chi''(z_0)\neq 0\) (by \Cref{p:surface}), and the \(k=2\) term in \eqref{e:series_identity} equals
\[
\frac{\chi''(z_0)}{2}\Bigl(2(a+\ri b)-\frac{2ab}{b}\Bigr)
=\ri\,\chi''(z_0)\,b.
\]
The remaining terms satisfy $\OO(b(|a|+b))$.
This proves \eqref{e:curve_reg}.

If \((x_0,s_0)\) is a cusp location, then \(\chi''(z_0)=0\) and \(\chi'''(z_0)\neq 0\) (again by \Cref{p:surface}).
Hence the leading contribution comes from \(k=3\) in \eqref{e:series_identity}:
\[
\frac{\chi'''(z_0)}{6}\Bigl(3(a+\ri b)^2-\frac{\Im[(a+\ri b)^3]}{b}\Bigr)=\chi'''(z_0) b\Bigl(\ri a-\frac{b}{3}\Bigr),
\]
and the remainder is
\(\OO\!\bigl(b(|a|+b)^2\bigr)\), proving \eqref{e:curve_cusp}.

If $(x_0, s_0)$ is a horizontal tangency location,
differentiating \(w=(x_0-z)^{-1}\) gives \(\partial_x w=w^2\,\partial_x z\), and combining with \eqref{e:dxz} yields
\begin{equation}\label{e:dxw}
\partial_x w(x,s)=\frac{1}{s}\,\frac{w(x,s)^2}{1/s+w(x,s)^2\widetilde\chi'(w(x,s))}.
\end{equation}
Using \eqref{e:Im_identity} and \(\Im z=\Im(x_0-1/w)=\Im(-1/w)={b}/(a^2+b^2)\), we obtain
\begin{equation}\label{e:rewrite_w}
\frac1s=-\frac{\Im\chi(z)}{\Im z}
=-\frac{\Im\widetilde\chi(w)}{b/(a^2+b^2)}
=-(a^2+b^2)\frac{\Im\widetilde\chi(w)}{b}.
\end{equation}
Therefore
\begin{equation}\label{e:denom_w}
\frac1s+w^2\widetilde\chi'(w)
=w^2\widetilde\chi'(w)-(a^2+b^2)\frac{\Im\widetilde\chi(w)}{b}.
\end{equation}

Write the Laurent expansion of $\widetilde\chi$ at $w=0$ as
\begin{equation}\label{e:chi_Laurent}
\widetilde\chi(w)=\frac{1}{t_0 w}+c_1 w+c_2 w^2+\OO(w^3),
\end{equation}
so that $c_1$ and $c_2$ are the coefficients of the $w$ and $w^2$ terms in the regular part.
Differentiating \eqref{e:chi_Laurent} gives
\[
\widetilde\chi'(w)=-\frac{1}{s_0 w^2}+c_1+2c_2 w+\OO(w^2),
\qquad
w^2\widetilde\chi'(w)=-\frac{1}{s_0}+c_1 w^2+2c_2 w^3+\OO(w^4).
\]
Moreover,
\[
\Im\widetilde\chi(w)=\Im\!\Bigl(\frac{1}{s_0 w}\Bigr)+\Im\!\bigl(c_1 w+c_2 w^2+\OO(w^3)\bigr)
=-\frac{b}{s_0(a^2+b^2)}+\Im\!\bigl(c_1 w+c_2 w^2+\OO(w^3)\bigr),
\]
hence
\[
(a^2+b^2)\frac{\Im\widetilde\chi(w)}{b}
=-\frac{1}{s_0}+(a^2+b^2)\frac{\Im\!\bigl(c_1 w+c_2 w^2+\OO(w^3)\bigr)}{b}.
\]
Substituting into \eqref{e:denom_w}, the constants $-1/s_0$ cancel and we obtain
\begin{equation}\label{e:denom_cancel}
\frac1s+w^2\widetilde\chi'(w)
= c_1 w^2+2c_2 w^3-(a^2+b^2)\frac{\Im\!\bigl(c_1 w+c_2 w^2+\OO(w^3)\bigr)}{b}
+\OO(w^4).
\end{equation}

\medskip
\noindent\textbf{Regular case: $c_1\neq 0$.}
Using $\Im(w)=b$ and $\Im(w^2)=2ab$, we have
\[
\Im\!\bigl(c_1 w+c_2 w^2+\OO(w^3)\bigr)=c_1 b+\OO\bigl(b(|a|+b)\bigr).
\]
Substituting into \eqref{e:denom_cancel} yields
\[
\frac1s+w^2\widetilde\chi'(w)
=c_1\bigl(w^2-(a^2+b^2)\bigr)+\OO\!\bigl(b(|a|+b)^2\bigr).
\]
Since
\[
w^2-(a^2+b^2)=(a+\ri b)^2-(a^2+b^2)=2\ri b\,(a+\ri b),
\]
we obtain
\[
\frac1s+w^2\widetilde\chi'(w)
=2c_1\,\ri b\,(a+\ri b)+\OO\!\bigl(b(|a|+b)^2\bigr)
=\del_w(w\wt \chi(w))|_{w=0}\,\ri b\,(a+\ri b)+\OO\!\bigl(b(|a|+b)^2\bigr)
\]
which is \eqref{e:curve_ht_reg}.

\medskip
\noindent\textbf{Cusp case: $c_1=0$ and $c_2\neq 0$.}
Now
\[
\Im\!\bigl(c_2 w^2+\OO(w^3)\bigr)=2c_2 ab+\OO\bigl(b(|a|+b)^2\bigr).
\]
Substituting into \eqref{e:denom_cancel} gives
\[
\frac1s+w^2\widetilde\chi'(w)
=2c_2\Bigl(w^3-a(a^2+b^2)\Bigr)+\OO\!\bigl(b(|a|+b)^3\bigr).
\]
Finally,
\[
w^3-a(a^2+b^2)=(a+\ri b)^3-a(a^2+b^2)
=-4ab^2+\ri b(3a^2-b^2)
=\ri b\,(a+\ri b)(3a+\ri b),
\]
and hence, using $6c_2=\del_w^2 (w\widetilde\chi(w))|_{w=0}$ ,
\[
\frac1s+w^2\widetilde\chi'(w)
=\del_w^2 (w\widetilde\chi(w))|_{w=0} \ri b\,(a+\ri b)(3a+\ri b)/3+\OO\!\bigl(b(|a|+b)^3\bigr),
\]
which is \eqref{e:curve_ht_cusp}.
\end{proof}

	\subsection{Frozen region}\label{s:frozen}
The following results due to \cite{astala2026dimer,de2010minimizers} indicates continuity properties for the gradient $\nabla H^*$ of the maximizer $H^*$ of $\mathcal{E}$ on $\mathfrak{P}$, as well as convexity properties for its arctic boundary. In what follows, for any direction $\omega \in \mathbb{R}^2 \setminus \big\{ (0, 0) \big\}$, the graph $G \subset \mathbb{R}^2$ of a (possibly discontinuous) function (whose domain is possibly disconnected or empty) is said to be convex (or concave) in the $\omega$ direction, if the following holds. Let $\rho_{\omega} : \mathbb{R}^2 \rightarrow \mathbb{R}^2$ denote the rotation such that $\rho_{\omega} (\omega) \in \mathbb{R}_{> 0} \cdot (0, 1)$ points vertically upwards. Then, each connected component of $\rho_{\omega} (G)$ is convex (or concave, respectively). 

\begin{proposition}[{\cite{astala2026dimer, de2010minimizers}}]\label{p:H*reg}
	
	The following statements hold. 
	
\begin{enumerate}
\item %On $\mathfrak{P} \setminus \mathfrak{L} (\mathfrak{P})$, $\nabla H^* (x, t)$ is piecewise constant, taking values in $\big\{ (0, 0), (1, 0), (1, -1) \big\}$ {\color{red} not sure if I need this statement}. 
On the frozen region $\fP\setminus \fL(\fP)$, $H^*$ is countably piecewise affine with gradient taking values in $\{ (0, 0), (1, 0), (1, -1)\}$.

\item If $\nabla H^*$ is continuous at $\zeta\in \fA$, then $\nabla H^*(\zeta)\in \{ (0, 0), (1, 0), (1, -1)\}$, and
\begin{enumerate}
\item $\nabla H^*(\zeta) = (0,0)$ if $f(\zeta) \in (0,\infty)$ and the slope of the arctic curve at $\zeta$ lies in $(1,\infty)$.
\item $\nabla H^*(\zeta) = (1,0)$ if $f(\zeta) \in (-1,0)$ and the slope of the arctic curve at $\zeta$ lies in $(-\infty,0)$.
\item $\nabla H^*(\zeta) = (1,-1)$ if $f(\zeta) \in (-\infty,-1)$ and the slope of the arctic curve at $\zeta$ lies in $(0,1)$.

\end{enumerate}

\item \label{i:disconnect_point} If $\nabla H^*$ is discontinuous at $\zeta_0\in\mathfrak{P}$, then there exists $\zeta_1\in \del \fP$ with the segment $[\zeta_0,\zeta_1]\subset \overline{\fP}\setminus \fL$, such that one of the following holds
\begin{enumerate}
\item $\zeta_0-\zeta_1=t\cdot (1,0)$ for some $t\in \bR$, and  $H^*(\zeta)=H^*(\zeta_1)+(1,0)\cdot (\zeta-\zeta_1)$ for any $\zeta\in [\zeta_0, \zeta_1]$.
\item $\zeta_0-\zeta_1=t\cdot (0,1)$ for some $t\in \bR$, and $H^*(\zeta)=H^*(\zeta_1)$  for any $\zeta\in [\zeta_0, \zeta_1]$ .
\item $\zeta_0-\zeta_1=t \cdot (1,1)$ for some $t\in \bR$, and $H^*(\zeta)=H^*(\zeta_1)$  for any $\zeta\in [\zeta_0, \zeta_1]$.
\end{enumerate}

%\item The gradient $\nabla H^*$ exists and is continuous on the set $\big\{ z \in \mathfrak P: m(z)<H^*(z)<M(z) \big\}$.
\item Fix a real number $c \in \mathbb{R}$; a vertex $p_0 \in \big\{ (0, 0), (1, 0), (1, -1) \big\} \in \overline{\mathcal{T}}$; and a direction $\omega \in \mathbb{R}^2 \setminus \big\{ (0, 0) \big\}$ such that
\begin{align}\label{e:omega}
\omega\cdot(p-p_0)>0, \quad \text{for all } p \in \overline{\mathcal T}\setminus \{p_0\}.
\end{align}

\noindent Let $S$ denote the interior of the set $\{ z \in \mathfrak P: H^*(z)=c+p_0\cdot z\}$; assume that $S\neq \emptyset$. Then $\partial S\cap \mathfrak P$ consists of the union of a convex graph (by above) and a concave graph (by below) in the $\omega$ direction.

\end{enumerate}
\end{proposition}

\begin{proof}

The first result appears in \cite[Theorem 2.5]{astala2026dimer}. For the second result, we prove the case $\nabla H^*(\zeta)=(0,0)$; the other cases can be treated analogously and are therefore omitted. In this case if we take a sequence of $\{\zeta_i\}_{i\geq 1}\in \fL$ approaching $\zeta$, then $\nabla H^*(\zeta_i)\rightarrow \nabla H^*(\zeta)=(0,0)$. By the relation \eqref{fh}, we have $\arg f(\zeta_i), \arg(f(\zeta_i)+1)\rightarrow 0$. It follows that $f(\zeta_i)\rightarrow f(\zeta)\in (0,\infty)$. Then first statement in \Cref{pa1} then implies that the slope of the arctic curve at $\zeta$ lies in $(1,\infty)$.
The third result follows from \cite[Theorem 1.3]{de2010minimizers} by noting that for $\zeta \in [\zeta_0, \zeta_1]$ we have $\nabla H(\zeta) \in \del \cT$, which implies $[\zeta_0, \zeta_1] \notin \fL$. Finally the fourth statement appears in  \cite[Theorem 4.2]{de2010minimizers}.
\end{proof}

The third statement in \eqref {pa1} that
the arctic boundary $\fA$ does not contain a tacnode(double points)
follows from \Cref{p:H*reg}.

\begin{proof}[Proof of \Cref{i:no_tacnode} in \Cref{pa1}]
By our assumption that the liquid region $\fA$ is connected, and the arctic boundary is also connected. 
Suppose, for contradiction, that a tacnode exists, where two pieces of the arctic boundary are tangent to each other. Then the arctic boundary contains a jordan curve from $\zeta^*$ to itself, and no other tacnodes (see \Cref{f:tacnode}), and $H^*$ is frozen inside the curve. We denote by $\fF$ the  region enclosed by the curve,  and its boundary $\del\fF\subset \fA$.

As $\zeta$ traverses the boundary curve $\del \fF$ from the tacnode back to itself, the slope of the tangent line to $\del \fF$ winds around $\bR\cup\{\infty\}$ at least once:
\begin{align}\label{e:all_tangent}
\{\text{slope of the tangent line to $\del \fF$ at $\zeta: \zeta\in \del \fF$}\}=\bR\cup\{\infty\}. 
\end{align}

We consider two cases for the behavior of $\nabla H^*$ on $\overline \fF\setminus\{\zeta^*\}$; either $\nabla H^*$ is continuous or discontinuous. If $\nabla H^*$ is continuous on $\overline \fF\setminus\{\zeta^*\}$, then $\nabla H^*$ is constant  on $\overline \fF\setminus\{\zeta^*\}$ and taking values in $\{(0,0), (1,0),(1,-1)\}$ by the first statement in \Cref{p:H*reg}. Without loss of generality, we assume $\nabla H^*\equiv (0,0)$ on $\overline \fF$. Then by the second statement in \Cref{p:H*reg} for any $\zeta\in \del \fF$, the slope of its tangent vector lies in $(1,\infty)$. It contradicts with \eqref{e:all_tangent}.

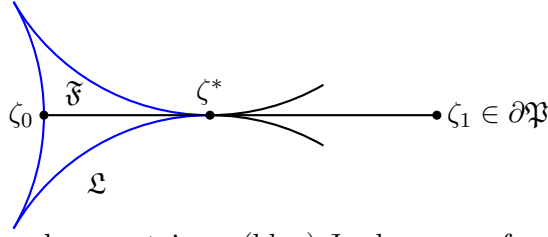
\begin{figure}
\centering
\vspace{-4em}
\begin{tikzpicture}[scale=3]
 \draw[thick,blue] (0,-1) ++(90:1)
       arc[start angle=90, end angle=150, radius=1];
      \draw[thick] (0,-1) ++(60:1)
       arc[start angle=60, end angle=90, radius=1];
 
  \draw[thick, blue] (0,1) ++(-150: 1)
        arc[start angle=-150, end angle=-90, radius=1];
    \draw[thick] (0,1) ++(-90: 1)
        arc[start angle=-90, end angle=-60, radius=1];
  
    \draw[thick, blue] (-1.732,0) ++(-30: 1)
        arc[start angle=-30, end angle=30, radius=1];

    \draw[thick] (-0.732,0) -- (1,0);

  % Labels
  \filldraw (0,0) circle (0.5pt) node[above] {$\zeta^*$};
    \filldraw  (-0.732,0)  circle (0.5pt) node[left] {$\zeta_0$};
    \filldraw  (1,0)  circle (0.5pt) node[right] {$\zeta_1\in \del \fP$};
    
     \node at (-0.6,0.1) {$\fF$};
    \node at (-0.5,-0.3) {$\fL$};
\end{tikzpicture}
\vspace{-5em}
\caption{The arctic boundary contains a (blue) Jordan curve from the tacnode $\zeta^*$ to itself.}
\label{f:tacnode}
\end{figure}

If $\nabla H^*$ is discontinuous at some point $\zeta_0 \in \overline{\fF} \setminus \{\zeta^*\}$, 
we first show that $\zeta_0$ lies on the tangent line to $\fA$ through $\zeta^*$. 
By the third statement of \Cref{p:H*reg}, there exists a segment connecting $\zeta_0$ to some 
$\zeta_1 \in \partial \fP$ that remains outside the liquid region. 
This segment $[\zeta_0, \zeta_1]$ must pass through the tacnode $\zeta^*$ and is tangent to $\fA$ at $\zeta^*$ (see \Cref{f:tacnode}). 
Again, by the third statement of \Cref{p:H*reg}, we may, without loss of generality, assume that 
 $\zeta_0-\zeta_1=t\cdot (1,0)$ for some $t\in \bR$, and for any $\zeta\in [\zeta_0, \zeta_1]$, \begin{align}\label{e:H_change}
 H^*(\zeta)=H^*(\zeta_1)+(1,0)\cdot (\zeta-\zeta_1).
 \end{align}
 
 We conclude that $[\zeta_0, \zeta_1]$ lies on the tangent line to $\fA$ through $\zeta^*$, which has slope~$0$. 
Let $\zeta_0$ denote the discontinuity point farthest from $\zeta^*$ along this tangent line. 
Then the segment $[\zeta_0, \zeta^*]$ contains all possible discontinuity points. 

The segment $[\zeta_0, \zeta^*]$ divides $\overline{\fF}$ into several disconnected pieces. 
On each piece, $\nabla H^*$ is constant (by continuity) and takes values in 
$\{(0,0), (1,0), (1,-1)\}$ by the first statement of \Cref{p:H*reg}. 
We now show that $\nabla H^*$ cannot be constantly $(0,0)$ on any piece. 
Indeed, otherwise $H^*$ would be constant, contradicting \eqref{e:H_change}. 
Hence, on $\overline{\fF}\setminus [\zeta_0, \zeta^*]$, the gradient $\nabla H^*$ can only take values in $\{(1,0), (1,-1)\}$. 

By the second statement of \Cref{p:H*reg}, it follows that for any 
$\zeta \in \partial \fF \setminus [\zeta_0, \zeta^*]$, the slope of the tangent vector lies in $(-\infty,0)\cup (0,1)$. 
This again contradicts \eqref{e:all_tangent}, since $\partial \fF \setminus [\zeta_0, \zeta^*]$ is obtained from $\partial \fF$ by removing only finitely many points.
\end{proof}

In the rest of this section, we prove \Cref{t:frozen_structure}.

\begin{figure}[ht]
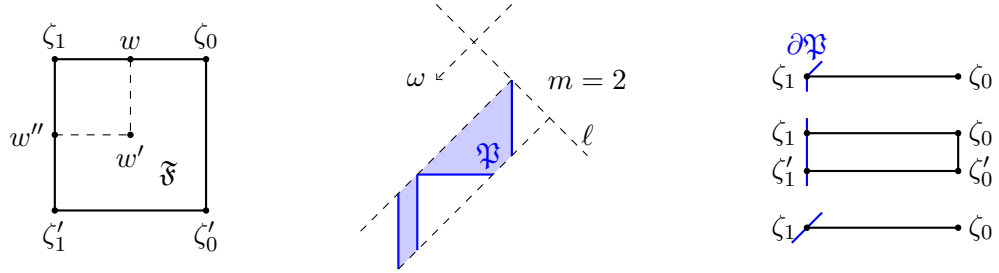

\centering
\begin{minipage}{0.30\textwidth}
\centering

\end{minipage}
\caption{The left panel illustrates \Cref{l:rec_height};
the middle panel is the region $
\{\, x + c \omega \in \fP : c \in \bR, \, x \in I \,\}
$;
the right panel depicts the three cases of $\zeta_1 \in \del\fP$.}
\label{f:zeta}
\end{figure}

\begin{lemma}\label{l:rec_height}
For a rectangle $\fF\in \overline\fP$ with vertices $\zeta_1, \zeta_0, \zeta_1', \zeta_0'$ (see  \Cref{f:zeta} left panel), if %the following holds
%\begin{enumerate}
%\item For any $\zeta\in (\zeta_1',\zeta_1]\cup [\zeta_1, \zeta_0)$, there exists a small $\varepsilon>0$, such that $B(\zeta, \varepsilon)\cap \fF\subset \overline \fP$
%\item 
for any $\zeta\in[\zeta_1',\zeta_1]\cup [\zeta_1, \zeta_0]$ or for any $\zeta\in[\zeta_1',\zeta'_0]\cup [\zeta_0', \zeta_0]$, the height function $H$ satisfies
\begin{align}\label{e:height_condition}
H(\zeta)=H(\zeta_1)+(1,0)\cdot (\zeta-\zeta_1).
\end{align}
Then \eqref{e:height_condition} holds for any $\zeta\in \fF$.
%\end{enumerate}
\end{lemma}
\begin{proof}
We prove the statement when \eqref{e:height_condition} holds for any $\zeta\in[\zeta_1',\zeta'_0]\cup [\zeta_0', \zeta_0]$. The other case can be proven in the same way, so we omit.

For any $w'\in \fF$, we denote $w\in [\zeta_0, \zeta_1]$ such that $(1,0)\cdot(w-w')=0$, and $w''\in [\zeta_1', \zeta_1]$ such $(0,1)\cdot(w''-w')=0$ then
\begin{align}
H(w)\leq H(w')\leq H(w'')+(1,0)\cdot(w'-w'')=H(\zeta_1)+(1,0)\cdot(w-\zeta_1)=H(w),
\end{align}
We conclude that $H(w')=H(\zeta_1)+(1,0)\cdot(w-\zeta_1)=H(\zeta_1)+(1,0)\cdot(w'-\zeta_1)$ for any $w'\in \fR$. 
\end{proof}

\begin{lemma}
On the frozen region $\fP\setminus \fL(\fP)$, $H^*$ is continuous and finitely piecewise affine with gradient taking values in $\big\{ (0, 0), (1, 0), (1, -1) \big\}$. Each affine piece is an open, simply connected region whose boundary consists of finitely many curves, each of which is, up to a rotation of the plane, the graph of a convex function.
\end{lemma}

\begin{proof}
We introduce the following set  $\cE$ of curves:
\begin{enumerate}
\item Lines starting from a corner vertex of $\fP$ and in one of the directions $(1,0), (0,1), (1,1)$.
There are finitely many such Lines.
\item We recall that each affine piece of the obstacle functions $M, m$ is in the form $\{c+p_0\cdot z\}$,  where $p_0 \in \{ (0, 0), (1, 0), (1, -1)\} $. There are finitely many such pieces. By the fourth statement of \Cref{p:H*reg}, let $S$ be the interior of the set $\{z\in\fP: H^*(z)=c+p_0 z\}$, and fix any direction $\omega$ such that \eqref{e:omega} holds. If $S\neq \emptyset$, then $\del S\cap \fP$ consists of the union of a convex graph (by above) and a concave graph (by below) in the $\omega$ direction. Next we show that $\del S\cap \fP$ consists of finite union of such curves (each of which is, up to a rotation of the plane, the graph of a convex function), and we collect them into $\cE$.

%Fix a line $\ell$ orthogonal to $\omega$, for any $z\in \ell$, we denote $N(z)$ the number of intervals of $\{z+c\omega\in \fP: c\in \bR\}$.  Since $\fP$ is a polygon, $N:\ell\mapsto \bZ_{\geq0}$ is bounded and finitely piecewise constant. Given an interval $I\subset \ell$ such that for any $z\in I$, $N(z)=n$ for some $n\geq 1$. Then $\{z+c\omega\in \fP: c\in \bR, z\in I\}$ consists of $n$ pieces. On each piece, the height function is increasing in the $\omega$ direction. So its intersection with $S$ consists of at most one convex curve as top boundary, and one concave curve as the bottom boundary. In total there are at most $2n$ curves. We conclude that $\del S\cap \fP$ consists of finite union of such convex/concave curves. 

Fix a line $\ell$ orthogonal to $\omega$ (see \Cref{f:zeta}, middle panel). 
The projections of the vertices of the polygon $\fP$ onto $\ell$ (in the $\omega$-direction) divide $\ell$ into finitely many intervals. 
Fix one such interval $I \subset \ell$. Then 
\[
\{\, x + c \omega \in \fP : c \in \bR, \, x \in I \,\}
\]
(if nonempty) decomposes into $m$ connected components, where $m \geq 1$. 
We denote one of these open components as $\fF$. 

In the $\omega$-direction, $\fF$ is bounded above by one edge of $\fP$ and below by another edge of $\fP$. The height function is linear on each edge of $\fP$. Hence for the top and bottom edges of $\fF$, either the entire edge lies in the set $\{z\in\fP: H^*(z)=c+p_0 z\}$, or it intersects $\{z\in\fP: H^*(z)=c+p_0 z\}$ in at most one point. Moreover, using \eqref{e:omega}, on $\fF$ the function $H^*(z) - p_0 \cdot z$ is monotone in the $\omega$-direction. Hence $S \cap \fF$ is bounded by a top boundary and a bottom boundary.

Consequently, the top boundary of $S \cap \fF$ is either a line segment lying on the top edge of $\fF$, or else it consists of at most two convex curves, separated by a point on the top edge of $\fF$.  
Similarly, the bottom boundary of $S \cap \fF$ consists of at most two concave curves.  
Altogether, there are at most $4m$ such curves.  
Summing over all subintervals of $\ell$, we conclude that $\partial S \cap \fP$ is a finite union of curves, each of which is, up to a rotation of the plane, the graph of a convex function.

\end{enumerate}

From the construction above, any two curves in $\cE$ intersect only at finitely many points. 
Consequently, the set $\fP \setminus (\overline{\fL(\fP)} \cup \cE )$ splits into finitely many connected open components. 
We now show that $\nabla H^*$ is continuous on each such component. 
Since $\nabla H^*$ is locally constant, it must take one of the values in $\{(0,0), (1,0), (1,-1)\}$ on each component. 
Hence $H^*$ is affine on every component, and therefore $H^*$ is finitely piecewise affine. Moreover, each affine piece is an open region whose boundary consists of finitely many curves from $\cE$.

Assume $\nabla H^*$ is  discontinuous at $\zeta_0\in \fP\setminus \overline{\fL(\fP)}$, by the third statement of \Cref{p:H*reg},  there exists $\zeta_1\in \del \fP$ with the segment $[\zeta_0,\zeta_1]\subset \overline{\fP}$ and $H^*$ is linear on $[\zeta_0, \zeta_1]$. Without loss of generality, we assume that $\zeta_0-\zeta_1=t\cdot(1,0)$ for some $t>0$, and  
\begin{align}\label{e:Hcondition}
H^*(\zeta)=H^*(\zeta_1)+(1,0)\cdot (\zeta-\zeta_1) \text{ for any } \zeta\in [\zeta_0, \zeta_1].
\end{align} 
There are three possibilities for $\zeta_1$, see \Cref{f:zeta} right panel:
\begin{enumerate}
\item $\zeta_1$ is a corner vertex of the polygon $\fP$.
\item $\zeta_1$ belongs to a vertical edge.
\item $\zeta_1$ belongs to an edge with slope $1$.
\end{enumerate}

In the first case, $\zeta_0$ belongs to a line starting from a corner vertex of $\fP$ and in one of the directions $(1,0), (0,1), (1,1)$, and $\zeta_0\in \cE$. The last two cases are analogous; hence we restrict attention to the second case. In this case there exists a thin rectangle $\fF\in \overline \fP$ with vertices $\zeta_1, \zeta_0, \zeta_1', \zeta_0'$ (see \Cref{f:zeta}), such that $[\zeta_1', \zeta_1]$ lies in a vertical edge of $\del \fP$. Then for $\zeta\in [\zeta_1', \zeta_1]$, $H^*(\zeta)=H^*(\zeta_1)=H(\zeta_1)+(1,0)\cdot (\zeta-\zeta_1)$. This together with \eqref{e:Hcondition} verifies the assumption of \Cref{l:rec_height} for $H=H^*$, we conclude $H^*(\zeta)=H^*(\zeta_1)+(1,0)\cdot (\zeta-\zeta_1)$ for any $\zeta\in \fF$.

For any $w'\in \fF$, we introduce $w''\in [\zeta'_1,\zeta_1]$ such that $(0,1)\cdot (w'-w'')=0$. The upper obstacle function $M$ satisfies
\begin{align}
H^*(w')\leq M(w')\leq M(w'')+(1,0)\cdot(w'-w'')=H^*(\zeta_1)+(1,0)\cdot(w'-\zeta_1)= H^*(w').
\end{align}
Thus $M(w')=H^*(w')$ for all $w'\in \fF$.

We denote the interior of the set $\{w\in \fP: H^*(w)=H^*(\zeta_1)+(1,0)\cdot(w-\zeta_1)\}$ by $S$, then $S$ contains the rectangle $\fF$, and $\zeta_0$ is on the boundary of $S$. In particular we also have $\zeta_0\in \cE$.

Finally, we prove that each affine piece is simply connected.  
Let $\fF$ be an affine piece with gradient $p_0$. Then 
$\fF$ is the interior of $\{z: H^*(z)=c+p_0\cdot z\}$. 
Take a direction $\omega$ satisfying \eqref{e:omega}.  
In this direction the function $H^*(z) - p_0 \cdot z$ is monotone.  

Suppose for contradiction that $\fF$ is not simply connected.  
Then $\fP \setminus \fF$ contains a bounded connected component, which we denote by $\frak O$.  
For any point $\zeta \in \frak O$, there exists a line segment $[\zeta_1,\zeta_2]\subset \overline{\fF}$ in the $\omega$-direction such that $\zeta \in [\zeta_1,\zeta_2]$.  
Since $H^*(z) - p_0 \cdot z$ is monotone in the $\omega$-direction and 
$
H^*(\zeta_1) - p_0 \cdot \zeta_1 
= H^*(\zeta_2) - p_0 \cdot \zeta_2 
= c$,
it follows that $H^*(\zeta) = c + p_0 \cdot \zeta$.  
Thus $\frak O \subset \{ z : H^*(z) = c + p_0 \cdot z \}$, contradicting the assumption that $\frak O$ is a hole of $\fF$.  

\end{proof}

\begin{lemma}\label{c:share_boundary}

If two affine pieces of $H^*$ with gradients $p_0\neq p_1 \in \{(0,0),(1,0),  (1,-1)\}$, 
share a portion of their boundary, then it consists of line segments $[\zeta_1,\zeta_0]$ satisfying
 $[\zeta_1,\zeta_0]$ is orthogonal to $p_0-p_1$.  Moreover, for any  $\zeta\in [\zeta_1,\zeta_0]$,     
    \begin{align}\label{e:Hexp}
        H^*(\zeta) \;=\; H^*(\zeta_1) + p_0 \cdot (\zeta - \zeta_1).
    \end{align}
\end{lemma}

\begin{proof}
We parametrize a connected piece of the shared boundary by $\gamma:[0,1]\mapsto \bR^2$. For any two points $\gamma(t), \gamma(t')$ with $t<t'$ on the shared boundary, since they are on the boundary of an affine piece of $H^*$ with gradient $p_0$, by continuity of $H^*$
\begin{align}\label{e:Hp0}
H^*(t')-H^*(t)=\int_{t}^{t'} p_0 \cdot \gamma'(u)\rd u.
\end{align}
By the same reasoning, we also have
\begin{align}\label{e:Hp1}
H^*(t')-H^*(t)=\int_{t}^{t'} p_1 \cdot \gamma'(u)\rd u
\end{align}
By comparing \eqref{e:Hp0} and \eqref{e:Hp1}, and sending $t'\rightarrow t$, we conclude that $\gamma'(t)$ is orthogonal to $p_0-p_1$. We conclude the shared boundary piece must be a line segment in the direction $\omega$, where $\omega$ is orthogonal to $p_0-p_1$. Moreover, \eqref{e:Hexp} holds.

\end{proof}

\begin{lemma}\label{c:boundary1}
Let $\fF$ be an affine piece of $H^*$ with gradient 
$
p_0 \in \{(0,0),\, (1,0),\, (1,-1)\},
$
and let 
$
\{p_1,p_2\} = \{(0,0),\, (1,0),\, (1,-1)\} \setminus \{p_0\}$.
Then the boundary of $\fF$ consists of finitely many curves of three types:
\begin{enumerate}
    \item portions of the arctic boundary whose slope belongs to $(1,\infty)$ if $p_0=(0,0)$; to $(-\infty,0)$ if $p_0=(1,0)$; to $(0,1)$ if $p_0=(1,-1)$.
    \item line segments orthogonal to $p_1-p_0$.
    \item line segments orthogonal to $p_2-p_0$.
\end{enumerate}
We refer to points where boundary pieces of different types meet as \emph{corner points} of $\fF$.  
\end{lemma}
\begin{proof}
Without loss of generality, we assume $p_0=(1,0)$, $p_1=(0,0)$ and $p_2=(1,-1)$. 
The boundary of $\fF$ consists of 
\begin{enumerate}
\item portions of the arctic boundary;
\item shared boundary with another affine piece;
\item boundary edge of the polygon $\fP$;
\end{enumerate}
For the first case, the slopes are described in the second statement of \Cref{p:H*reg}.
For the second case, by \Cref{c:share_boundary}, the shared boundary is a line segment oriented in the direction orthogonal to either 
$p_0 - p_1 = (1,0)$ or $p_0 - p_2 = (0,1)$.
For the third case, note that the boundary edges of $\fP$ lie in one of the directions 
$
(1,0), (0,1),  (1,1)$.
We claim that the boundary edges of $\fF$ can only be in the directions $(1,0)$ or $(0,1)$.  
Suppose, for the sake of contradiction, that the boundary of $\fF$ contains a segment 
$
[\zeta, \zeta']$ in the direction $(1,1)$.
By the definition of the boundary height function as in \Cref{p}, we then have
\begin{align}\label{e:H_equal}
H(\zeta') = H(\zeta).
\end{align}
On the other hand, since $[\zeta, \zeta']$ lies on the boundary of $\fF$, where $\nabla H^* = (1,0)$, it follows that
\begin{align}
H(\zeta') - H(\zeta) = (1,0)\cdot (\zeta'-\zeta) \neq 0.
\end{align}
This contradicts \eqref{e:H_equal}. Hence no boundary edge of $\fF$ can lie in the direction $(1,1)$.

\end{proof}

\begin{figure}
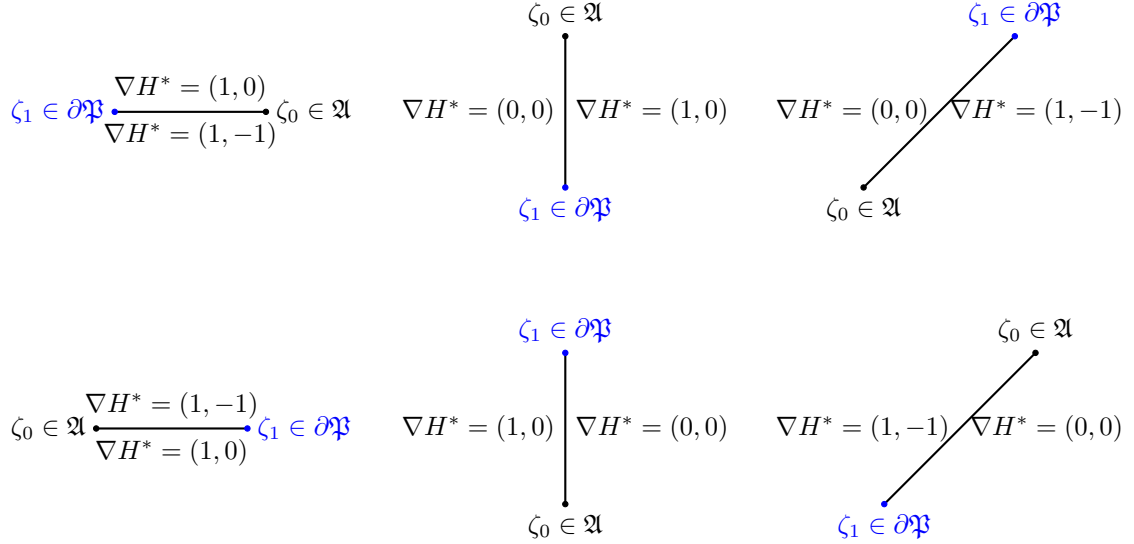

\centering
\begin{minipage}{0.30\textwidth}
\centering

\end{minipage}

\caption{For the cases in the first row, $H^* = M$ on $[\zeta_1, \zeta_0]$;  
for those in the second row, $H^* = m$ on $[\zeta_1, \zeta_0]$.}
\label{f:share_edge}
\end{figure}

\begin{lemma}\label{c:edge}
Let $\fF_0$ and $\fF_1$ be two affine pieces of $H^*$ with distinct gradients 
$p_0 \neq p_1 \in \{(0,0),(1,0),(1,-1)\}$. 
Then their shared boundary is a finite union of line segments. Moreover, each segment $[\zeta_1,\zeta_0]$ satisfies:
\begin{enumerate}
\item There are six cases in total. In the three cases shown in the first row of \Cref{f:share_edge}, we have $H^* = M$ on $[\zeta_1,\zeta_0]$, while in the three cases shown in the second row, we have $H^* = m$ on $[\zeta_1,\zeta_0]$.
\item There exists an edge $[\zeta_2,\zeta_1]$ of the polygon $\fP$ such that $[\zeta_1,\zeta_0]$ lies on the extension of $[\zeta_2,\zeta_1]$.
\end{enumerate}
\end{lemma}

\begin{proof}

%We now prove the third statement in \eqref{c:share_boundary}. By possibly partitioning the interval
%$[\zeta_1,\zeta_0]$ into subintervals, we may—without loss of
%generality—work on a subinterval (which we continue to denote by
%$[\zeta_1,\zeta_0]$) such that
%\begin{align}\label{e:z1z0}
%(\zeta_1,\zeta_0)\cap \partial \fP=\varnothing.
%\end{align}
%Equivalently, every $\zeta \in (\zeta_1,\zeta_0)$ lies strictly inside $\fP$ and away from the arctic boundary.  
%We also assume that the segment $[\zeta_1,\zeta_0]$ cannot be extended to a larger interval still satisfying %\eqref{e:z1z0}.

As illustrated in \Cref{f:share_edge}, six distinct cases arise, depending on the gradients $p_0, p_1$ and the relative positioning of the two affine pieces.  
We consider the case $p_0=(1,0)$ and $p_1=(1,-1)$, in which the shared boundary $[\zeta_1,\zeta_0]$ is horizontal. Moreover, the affine pieces $\fF_0$ and $\fF_1$, with $\nabla H^* = p_0 = (1,0)$ and  $\nabla H^* = p_1 = (1,-1)$ respectively, lie above and below the interval $[\zeta_1,\zeta_0]$.
This corresponds to the first panel in the first row of \Cref{f:share_edge}.  The remaining cases can be proved in the same way and are therefore omitted. In the rest of the proof, we will use these two properties for $\fF_0$.

\begin{enumerate}[label=P\arabic*]
\item \label{f1}By \Cref{c:boundary1}, the boundary $\fF_0$ consists of portions of the arctic boundary with slope in $(-\infty,0)$, together with horizontal and vertical segments.  
Thus, except for corner points, the boundary of $\fF_0$ has slope in $[-\infty,0]$.
\item \label{f2} By the fourth statement of \Cref{p:H*reg},  $\del\fF_0\cap \fP$ consists of a convex graph (by above) in the $\omega$ direction, which satisfies \eqref{e:omega} with $p_0=(1,0)$
\end{enumerate}

%We first show that $\zeta_1$ is a vertex of $\fP$. 

If $\zeta_1$ is not a corner point of $\fF_0$, we can continue moving to the left until reaching a corner point of $\fF_0$, which we denote by $w_0$. By \Cref{f1}, the boundary of $\fF_0$ has slope in $[-\infty,0]$. Hence, there exists $\zeta\in (\zeta_1,\zeta_0)$, and  a sufficiently small rectangle $w_0', w_0'', \zeta, w_0$ inside $\fF_0$. By the same argument, there exists a sufficiently small parallelogram $w_1, \zeta,w_1'',  w_1'$ inside $\fF_1$. See \Cref{f:shared_boundary}.

\begin{figure}[ht]
\centering
\begin{minipage}{0.33\textwidth}
\centering

\end{minipage}

\caption{$\fF_0$ contains a rectangle $w_0', w_0'', \zeta, w_0$,  and $\fF_1$ contains a parallelogram $w_1, \zeta,w_1'',  w_1'$.}
\label{f:shared_boundary}
\end{figure}

By our assumption that the segment $[\zeta_1,\zeta_0]$ cannot be extended to a larger interval, there are three possibilities corresponding to three panels in \Cref{f:shared_boundary}
\begin{enumerate}
\item $w_0=w_1=\zeta_1$
\item $w_0$ is on the left of $w_1=\zeta_1$
\item $w_1$ is on the left of $w_0=\zeta_1$.

\end{enumerate}
We discuss the above three cases one by one. In the first case,  by the third statement of \Cref{t:g_behavior}, the arctic boundary $\fA$ is locally strictly convex, except for inward cusps. 
Since the angle $\angle w_0' \zeta_1 w_1'=120^\circ<180^\circ$, we conclude that  $\zeta_1\not\in\fA$. Because $\zeta_1$ is a corner point for both $\fF_0$ and $\fF_1$, $[w_0', \zeta_1]$ belongs to the boundary of $\fF_0$ and  $[w_1', \zeta_1]$ belongs to the boundary of $\fF_1$.

By \Cref{f2}, $\del\fF_0\cap \fP$ consists of a convex graph (by above) in the $\omega$ direction. This convex graph contains $[\zeta_1, \zeta]$, so it does not contain $[w_0', \zeta_1]$. Hence, $[w_0', \zeta_1]\in \del \fP$. By the same argument $[w_1', \zeta_1]\in \del \fP$. However, this contradicts with our \Cref{p} that $\fP$ is formed by segments
with slopes $0, 1, \infty$ cyclically repeated as we follow the boundary in the
counterclockwise direction. 

% and $\zeta_1$ is in the interior of the frozen region. 
In the second case that $w_0$ is on the left of $w_1=\zeta_1$. By the same argument as in the first case, we have that $[w_1', \zeta_1]\in \del \fP$. By \Cref{p},  $\fP$ is formed by segments
with slopes $0, 1, \infty$ cyclically repeated as we follow the boundary in the
counterclockwise direction. Thus after $[w_1', \zeta_1]$, the next edge of $\fP$ is horizontal, we denote it as $[\zeta_2,\zeta_1]$. Then $[\zeta_1, \zeta_0]$ lies on the extension of $[\zeta_2, \zeta_1]$; see the first row of \Cref{f:P_boundary}.

%Then  $\zeta_2\in [w_0,\zeta_1]$, because $\del\fF_0\cap \fP$ consists of a convex graph (by above) in the $\omega$ satisfying \eqref{e:omega} with $p_0=(1,0)$. 

In the third case that $w_1$ is on the left of $w_0=\zeta_1$, by the same argument as in the second case, we have that $[w_0', \zeta_1]\in \del \fP$, and $\fP$ has an horizontal edge $[\zeta_2,\zeta_1]$ with  $\zeta_2\in [w_1,\zeta_1]$; see the second row of \Cref{f:P_boundary}.

In both cases, $\zeta_1\in \del \fP$ and for any  $\zeta\in [\zeta_1, \zeta_0]$, $H^*$ satisfies
\begin{align}
M(\zeta)\leq h(\zeta_1)+(1,0)\cdot(\zeta-\zeta_1)=H^*(\zeta)\leq M(\zeta).
\end{align}
Hence $H^*$ agrees with the obstacle function $M$ on $[\zeta_1, \zeta_0]$.

\begin{figure}[ht]
\centering
\begin{minipage}{0.4\textwidth}
\centering

\end{minipage}

\caption{$\fF_0$ contains a thin rectangle $w_0', w_0'', \zeta, w_0$,  and $\fF_1$ contains a thin parallelogram $w_1, \zeta,w_1'',  w_1'$.}
\label{f:P_boundary}
\end{figure}

\end{proof}

\begin{lemma}\label{c:boundary2}
For any $\zeta\in \fA$, the following holds
\begin{enumerate}
\item If $\zeta$ is not a tangent location, then $\zeta\in \fP$ and $\nabla H^*$ is continuous at $\zeta$.
\item If $\zeta\in \fP$ is a tangent location,  there exist two affine pieces $\fF_0$ and $\fF_1$, whose shared boundary $\ell$ is tangent to $\fA$ at $\zeta$, and for sufficiently small $\delta>0$, $B_\delta(\zeta)\cap( \fP\setminus \overline{\fL})\subset  \fF_0\cup  \fF_1 \cup \ell$.
\item If $\zeta\in \del \fP$ is a tangent location, then an edge of $\fP$ is tangent to $\fA$ at $\zeta$. Moreover there exist two affine pieces $\fF_0$ and $\fF_1$, and for sufficiently small $\delta>0$, $B_\delta(\zeta)\cap (\fP\setminus \overline{\fL})\subset \fF_0\cup \fF_1$.
\end{enumerate}
\end{lemma}

\begin{figure}[ht]
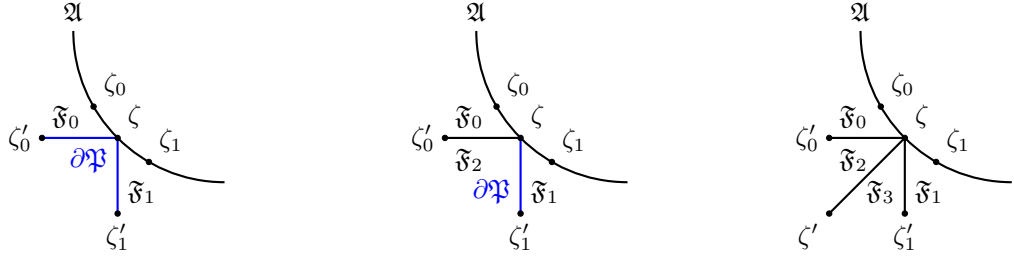

\centering
\begin{minipage}{0.33\textwidth}
\centering

\end{minipage}

\caption{The slope of the tangent line at $\zeta\in \fA$ lies in $(0,\infty)$. }
\label{f:arc_point}
\end{figure}

\begin{proof}

%We can take two points $\zeta_0, \zeta_1\in \fA$ such that $\zeta$ lies on the piece of arctic boundary from $\zeta_0$ to $\zeta_1$ as we follow $\fA$ in the counterclockwise direction. Moreover, $\nabla H^*$ is continuous on the arctic boundary pieces $\fA(\zeta_0, \zeta)$(from $\zeta_0$ to $\zeta$) and $\fA(\zeta, \zeta_1)$ (from $\zeta$ to $\zeta_1$); see \Cref{f:arc_point}.
%
%
%Since $\zeta\in \fA$ is not a tangent location, without loss of generality we assume that the slope of the tangent line at each point of the arctic boundary piece $\fA(\zeta_0, \zeta_1)$ (from $\zeta_0$ to $\zeta_1$) is in $(0,\infty)$. 

We choose two points $\zeta_0, \zeta_1 \in \fA$ such that $\zeta$ lies on the portion of the arctic boundary from $\zeta_0$ to $\zeta_1$ when traversing $\fA$ in the counterclockwise direction, which we denote as $\fA(\zeta_0,\zeta_1)$. Moreover, the gradient $\nabla H^*$ is continuous along the boundary pieces $\fA(\zeta_0,\zeta)$ (from $\zeta_0$ to $\zeta$) and $\fA(\zeta,\zeta_1)$ (from $\zeta$ to $\zeta_1$); see \Cref{f:arc_point}.  

Since $\zeta \in \fA$ is not a tangency point, we may assume without loss of generality that the slope of the tangent line at each point of $\fA(\zeta_0,\zeta_1)$ lies in $(0,\infty)$.

There exist affine pieces $\fF_0, \fF_1$ such that $\fA(\zeta_0, \zeta)\subset \del \fF_0$ and $\fA(\zeta, \zeta_1)\subset \del\fF_1$ respectively. Moreover, by the second statement of \Cref{p:H*reg}, $\nabla H^*=(1,0)$ on $\fF_0$ and $\fF_1$. If $\fF_0=\fF_1$ are the same affine piece, whose boundary contains the whole arc $\fA(\zeta_0, \zeta_1)$(from $\zeta_0$ to $\zeta_1$), then $\zeta\in \fP$ and $\nabla H^*$ is continuous at $\zeta$.  Otherwise $\fF_0\neq \fF_1$, and $\zeta$ is a corner point of both $\fF_0$ and $\fF_1$. By \Cref{f1} and the third statement in \Cref{pa1} that there is no tacnode, there exists short horizontal segments $[\zeta_0',\zeta]\subset \del \fF_0$, and a vertical segment $[\zeta_1',\zeta]\subset \del \fF_1$. We notice that $\zeta$ cannot be a cusp point. We have several cases:
\begin{enumerate}
\item 
If $[\zeta_0',\zeta], [\zeta_1', \zeta]\subset \del\fP$ (see the first panel in \Cref{f:arc_point}), this contradicts with \Cref{p} that $\fP$ is formed by segments
with slopes $0, 1, \infty$ cyclically repeated as we follow the boundary in the
counterclockwise direction. 

\item If $[\zeta_0',\zeta]\subset \fP$ and $[\zeta_1', \zeta]\subset  \del\fP$ (see the second panel in \Cref{f:arc_point}), 
 then $[\zeta_0', \zeta]$ belongs to a share boundary between $\fF_0$ and another  affine piece $\fF_2$ with $\nabla H^*=(1,-1)$ (here we used the second statement of \Cref{c:share_boundary}). In this case, $H^*=M$ on $[\zeta_0', \zeta]$. Since $\zeta\in \del \fP$, for any $\zeta'\in [\zeta_0', \zeta]$
 \begin{align}
m(\zeta')\geq h(\zeta)+(1,0)\cdot(\zeta'-\zeta)=H^*(\zeta')\geq m(\zeta').
\end{align}
Thus on $[\zeta_0', \zeta]$, $H^*=m=M$. This contradicts with \Cref{a:asump} that the trivial set is empty.

\item If $[\zeta_1',\zeta]\subset \fP$ and $[\zeta_0', \zeta]\subset  \del\fP$, we get a contradiction by the same argument as in the second case, so we omit.

\item
If $[\zeta_0',\zeta], [\zeta_1', \zeta]\subset \fP$ (see the third panel in \Cref{f:arc_point}), there exist affine pieces $\fF_2$ with $\nabla H^*=(1,-1)$, and $\fF_3$ with $\nabla H^*=(0,0)$ such that $[\zeta_0',\zeta]\subset \del \fF_0\cap \del \fF_2$ and $[\zeta_1',\zeta]\subset \del \fF_1\cap \del \fF_3$. Moreover, by \Cref{f1}, there exists a short segment $[\zeta',\zeta]\subset  \del \fF_2\cap \del \fF_3$ with slope $1$. By the second statement of \Cref{c:share_boundary}, $H^*=M$ on $[\zeta_0',\zeta]$, and $H^*=m$ on $[\zeta_2',\zeta]$. Hence $H^*(\zeta)=m(\zeta)=M(\zeta)$, this contradicts with \Cref{a:asump} that the trivial set is empty. 
\end{enumerate}

If $\zeta \in \fA$ is a tangency point, we may assume without loss of generality that the tangent at $\zeta$ has slope $\infty$.  
Then, by the second statement of \Cref{pa1}, the tangent along $\fA(\zeta_0,\zeta)$ has slope in $(0,1)$, while the tangent along $\fA(\zeta,\zeta_1)$ has slope in $(-\infty,0)$. There exist affine pieces $\fF_0, \fF_1$ such that $\fA(\zeta_0, \zeta)\subset \del \fF_0$ and $\fA(\zeta, \zeta_1)\subset \del\fF_1$ respectively. Moreover, by the second statement of \Cref{p:H*reg}, $\nabla H^*=(0,0)$ on $\fF_0$,  and $\nabla H^*=(1,0)$ on $\fF_1$.
See the first and the second panels of \Cref{f:arc_boundary_point}.

We first discuss the case that $\zeta\in \fP$. By \Cref{c:boundary1}, the potential boundary segments of affine pieces meeting at $\zeta$ consist of the coordinate-axis rays from $\zeta$ that lie outside the liquid region, together with the portion of the arctic boundary through $\zeta$.  
We enumerate these axis-aligned rays as
$
[\zeta,\zeta_0'],  [\zeta,\zeta_1'],  [\zeta,\zeta_2']$, and $ [\zeta,\zeta_3']$,
see the first panel of \Cref{f:arc_boundary_point}.  
Together with the arctic boundary, they exhaust all possible boundary pieces of affine regions $\fF$ such that $\zeta \in \partial \fF$.

If $\zeta$ is also a cusp point, the only possible ray is $[\zeta,\zeta_1']$, as shown in the second panel of \Cref{f:arc_boundary_point}.  
In this case, the shared boundary $[\zeta,\zeta_1'] \subset \ell$ of $\fF_0$ and $\fF_1$ is tangent to $\fA$ at $\zeta$.  
Moreover, for sufficiently small $\delta>0$, we have
$B_\delta(\zeta)\cap( \fP\setminus \overline{\fL})\subset  \fF_0\cup  \fF_1 \cup \ell$ .

If $\zeta$ is not a cusp location, there are several possibilities
\begin{enumerate}
\item $[\zeta, \zeta_0']\subset \del \fF_0\cap \del \fF_1$, then the shared boundary $[\zeta, \zeta_0']\subset \ell$ of $\fF_0, \fF_1$ is tangent to $\fA$ at $\zeta$, and for sufficiently small $\delta>0$, $B_\delta(\zeta)\cap( \fP\setminus \overline{\fL})\subset  \fF_0\cup  \fF_1 \cup \ell$. See the third panel of \Cref{f:arc_boundary_point}.

\item $[\zeta, \zeta_1']\subset \del \fF_0\cap \del \fF_1$,  then the shared boundary $[\zeta, \zeta_1']\subset \ell$ of $\fF_0, \fF_1$ is tangent to $\fA$ at $\zeta$, and for sufficiently small $\delta>0$, $B_\delta(\zeta)\cap( \fP\setminus \overline{\fL})\subset  \fF_0\cup  \fF_1 \cup \ell$. See the fourth panel of \Cref{f:arc_boundary_point}.

\item  $[\zeta, \zeta_0']\subset \del \fF_0$ and $[\zeta, \zeta_1']\subset \del \fF_1$. In this case, by the second statement of \Cref{c:share_boundary}, $[\zeta, \zeta_0']$ is shared by $\fF_0$ and another affine piece with gradient $(1,0)$. Hence $H^*=m$ on $[\zeta, \zeta_0']$. By the same argument $H^*=M$ on $[\zeta, \zeta_1']$. It follows that $H^*(\zeta)=m(\zeta)=M(\zeta)$, which contradicts with \Cref{a:asump} that the trivial set is empty.

\item $(\zeta, \zeta_0')\subset \fF_0$ and $[\zeta, \zeta_1']\subset \del \fF_1$. In this case, either $[\zeta, \zeta_1']\subset \del \fF_0$, which has been discussed in the second case, or $[\zeta, \zeta_2']\subset \del \fF_0$. Then by \Cref{c:share_boundary}, $[\zeta, \zeta_2']$ is shared by $\fF_0$ and another affine piece $\fF_2$ with gradient $(1,-1)$. Then we must have that $[\zeta, \zeta_1']\subset \del \fF_1\cap \del \fF_2$, which is impossible by \Cref{c:share_boundary}.

\item $[\zeta, \zeta_0']\subset \del \fF_0$ and $(\zeta, \zeta_1')\subset  \fF_1$. This is impossible, by the same argument as in the fourth statement. 

\item $(\zeta,\zeta_0') \subset \fF_0$ and $(\zeta,\zeta_1') \subset \fF_1$.  
By \Cref{c:boundary1}, this implies that $\fF_0$ contains the angular region $\angle \zeta_0'\zeta\zeta_2'$ and $\fF_1$ contains the angular region $\angle \zeta_3'\zeta\zeta_1'$.  
Since these two regions overlap, such a configuration is impossible.

\end{enumerate}

Next we discuss the case that $\zeta\in \del \fP$. %If $\zeta$ is on the interior of an edge, the edge must be vertical and tangent to $\fA$ at $\zeta$. 
%Otherwise $\zeta$ is a vertex of the polygon $\fP$. 
In this case, $
[\zeta,\zeta_0'],  [\zeta,\zeta_1'],  [\zeta,\zeta_2']$, and $ [\zeta,\zeta_3']$, also exhaust all possible boundary edges of $\fP$ adjacent to $\zeta$. There are several possibilities
\begin{enumerate}
\item  $[\zeta, \zeta_0'], [\zeta, \zeta_1']\subset \del \fP$. In this case $\zeta\in (\zeta_0', \zeta_1')\subset \del \fP$, and the edge is tangent to $\fA$ at $\zeta$. Moreover, for sufficiently small $\delta>0$, $B_\delta(\zeta)\cap (\fP\setminus \overline{\fL})\subset \fF_0\cup \fF_1$. See the first panel of \Cref{f:arc_edge_point}.

\item $[\zeta, \zeta_0']\subset \del \fP$ and $(\zeta, \zeta_1')\subset \fP$. In this case, $\zeta$ is a vertex of $\fP$, and  $[\zeta, \zeta_3']\subset \del \fP$, because by \Cref{p} that $\fP$ is formed by segments with slopes $0, 1, \infty$ cyclically repeated as we follow the boundary in the counterclockwise direction. In this case the edge of $\fP$ containing $[\zeta, \zeta_0']$ is tangent to $\fA$ at $\zeta$, and $\fF_0$ contains the whole angular region $\angle \zeta_0 \zeta \zeta_0'$.

For $\fF_1$, either it contains the whole angular region $\angle \zeta_3' \zeta \zeta_1$, or $[\zeta, \zeta_1']\subset\del \fF_1$, is shared by $\fF_1$  and another affine piece $\fF_2$ with gradient $(0,0)$, by \Cref{c:share_boundary}.  In this second case, $H^*=M$ on $[\zeta, \zeta_1']$. Since $\zeta\in \del \fP$, for any $\zeta'\in [\zeta, \zeta_1']$
 \begin{align}
m(\zeta')\geq m(\zeta)=H^*(\zeta)=H^*(\zeta')\geq m(\zeta').
\end{align}
Thus on $[\zeta, \zeta_1']$, $H^*=m=M$. This contradicts with \Cref{a:asump} that the trivial set is empty. 

We conclude that $\fF_1$ contains the whole angular region $\angle \zeta_3' \zeta \zeta_1$, and for sufficiently small $\delta>0$, $B_\delta(\zeta)\cap (\fP\setminus \overline{\fL})\subset \fF_0\cup \fF_1$. See the second panel of \Cref{f:arc_edge_point}.

\item $(\zeta, \zeta_0')\subset \fP$ and $[\zeta, \zeta_1']\subset\del \fP$. By the same argument in the second case,  $\zeta$ is a vertex of $\fP$, $[\zeta, \zeta_2'], [\zeta, \zeta_1']\subset\del \fP$, and the edge of $\fP$ containing $[\zeta, \zeta_1']$ is tangent to $\fA$ at $\zeta$.  Moreover,  $\fF_0$ contains the whole angular region $\angle \zeta_0 \zeta \zeta_2'$, $\fF_1$ contains the whole angular region $\angle \zeta_1' \zeta \zeta_1'$. It follows that for sufficiently small $\delta>0$, $B_\delta(\zeta)\cap (\fP\setminus \overline{\fL})\subset \fF_0\cup \fF_1$. See the third panel of \Cref{f:arc_edge_point}.

\item $(\zeta, \zeta_0'), (\zeta, \zeta_1')\subset \fP$. In this case, by the same argument as in the second case, we must have that $\fF_1$ contains the angular region $\angle \zeta_3' \zeta \zeta_1$, and $\fF_0$ contains the angular region $\angle \zeta_0 \zeta \zeta_2'$, which is impossible.

 %$[\zeta, \zeta_2'], [\zeta,\zeta_3']\subset \fP$. By the same argument as in the second case, $\fF_1$ contains the angular region $\angle \zeta_2' \zeta \zeta_1$. In particular $[\zeta, \zeta_2']\subset \del \fF_1$, which contradicts with \Cref{e:boundary}. 

\end{enumerate}

\end{proof}

\begin{figure}[ht]
\centering
\begin{minipage}{0.23\textwidth}
\centering
%
\end{minipage}

\caption{The neighborhood of a tangent location $\zeta\in \fA\cap \fP$. }
\label{f:arc_boundary_point}
\end{figure}

\begin{figure}[ht]
\centering

\begin{minipage}{0.3\textwidth}
\centering
%

\end{minipage}

\caption{The neighborhood of a tangent location $\zeta\in \fA\cap \del\fP$. }
\label{f:arc_edge_point}
\end{figure}

\begin{lemma}\label{c:boundary3}.
Let $\fF$ be an affine piece of $H^*$. 
Then each boundary piece from a corner point $\zeta_1$ to $\zeta_0$ is either:
\begin{enumerate}
\item a portion of the arctic boundary $\fA$, with $\zeta_1$ and $\zeta_0$ being two consecutive tangency points along $\fA$; or
\item a line segment $[\zeta_1,\zeta_0]$ orthogonal to $p_1-p_0$ or $p_2-p_0$, tangent to $\fA$ at $\zeta_0$, and satisfying one of the following:
  \begin{enumerate}
  \item $[\zeta_1,\zeta_0] \subset [\zeta_1,\zeta_2]$ for some edge $[\zeta_1,\zeta_2]$ of $\fP$;
  \item $[\zeta_1,\zeta_2] \subset [\zeta_1,\zeta_0)$ for some edge $[\zeta_1,\zeta_2]$ of $\fP$;
  \item $[\zeta_1,\zeta_0]$ lies on the extension of an edge $[\zeta_2,\zeta_1]$ of $\fP$.
  \end{enumerate}
\end{enumerate}
\end{lemma}

\begin{figure}[ht]
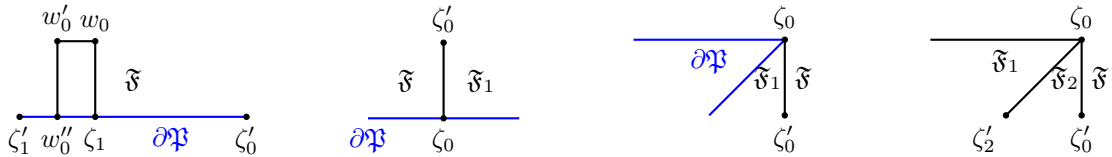

\centering
\begin{minipage}{0.23\textwidth}
\centering


\end{minipage}

\caption{The neighborhoods of $\zeta_1$ and $\zeta_0$.}
\label{f:segment_case}
\end{figure}

\begin{proof}
The first statement follows from the first statement of \Cref{c:boundary1}, and the first statement of \Cref{c:boundary2}.

Next we prove the second statement. Without loss of generality, we assume that on $\fF$, $\nabla H^*=(1,0)$,  $[\zeta_1, \zeta_0]$ is a horizontal segment and $\fF$ lies above it. The other cases can be proven in the same way, so we omit. 
% In the rest of the proof, we will use these two properties for $\fF$.
%\begin{enumerate}[label=P\arabic*]
%\item \label{f11}By \Cref{c:boundary1}, the boundary $\fF$ consists of portions of the arctic boundary with slope in $(-\infty,0)$, together with horizontal and vertical segments.  
%Thus, except for corner points, the boundary of $\fF$ has slope in $[-\infty,0]$.
%\item \label{f21} By the fourth statement of \Cref{p:H*reg},  $\del\fF\cap \fP$ consists of a convex graph (by above) in the $\omega$ direction, which satisfies \eqref{e:omega} with $p_0=(1,0)$
%\end{enumerate}

First we notice that $[\zeta_1,\zeta_0]\cap \del \fP$ is connected. Otherwise, say there exists $(\zeta_1',\zeta_0')\subset [\zeta_1,\zeta_0]\cap\fP$ and $\zeta_1',\zeta_0'\in \del \fP$. Then by \Cref{c:share_boundary}, $[\zeta_1',\zeta_0']$ is on the shared boundary of $\fF$ and another affine piece with gradient $(1,-1)$. Moreover, by the first statement of \Cref{c:edge} $H^*=M$ on $[\zeta_1', \zeta_0']$. Moreover, since $\zeta_0'\in \del \fP$, for any $\zeta\in [\zeta_1', \zeta_0']$
    \begin{align}
        m(\zeta) \geq h(\zeta'_0) + (1,0)\cdot (\zeta - \zeta'_0)=H^*(\zeta)\geq m(\zeta).
    \end{align}
Thus on $[\zeta_1, \zeta_0]$, $H^*=m=M$. This contradicts with \Cref{a:asump} that the trivial set is empty. 

There are several cases for $[\zeta_1,\zeta_0]\cap \del \fP$. 
\begin{enumerate}
\item There exists an horizontal edge $[\zeta_1', \zeta_2]$ of $\fP$, such that $\zeta_1\in [\zeta_1', \zeta_2)$. In this case we first show that $\zeta_1=\zeta_1'$. Otherwise $\zeta_1 \in (\zeta_1', \zeta_2)$. By \Cref{c:boundary1}, except for corner points, the boundary of $\fF$ has slope in $[-\infty,0]$. There exists a short vertical segment $[w_0, \zeta_1]\subset \overline \fF$. We can then take a small rectangle $\zeta_1, w_0, w'_0,  w_0''$ inside $\overline \fP$ with $w_0''\in (\zeta_1', \zeta_1)$. See the first panel of \Cref{f:segment_case}.
Then $[w_0'', \zeta_1]\subset [\zeta_1', \zeta_2]$, and  for $\zeta\in [w_0, \zeta_1]\cap [w_0'', \zeta_1]$ the height function satisfies
\begin{align}\label{e:heightbb}
H^*(\zeta)=H^*(\zeta_1)+(1,0)\cdot(\zeta-\zeta_1).
\end{align}
This verifies the assumptions in \Cref{l:rec_height}, and it follows that \eqref{e:heightbb} holds for any $\zeta$ in the rectangle $\zeta_1, w_0, w_0', w_0''$. In particular, $\fF$ contains this rectangle. But this contradicts to our assumption that $\zeta_1$ is a corner point of $\fF$.

\item There exists $\zeta_0'\in (\zeta_1, \zeta_0]$ such that $(\zeta_1, \zeta_0')\in \fP$. In this case $(\zeta_1, \zeta_0')$ is on the shared boundary of $\fF$ with another affine piece with gradient $(1,-1)$. By \Cref{c:edge}, there exists an edge $[\zeta_2, \zeta_1]$ of the polygon $\fP$ such that $[\zeta_1, \zeta_0]$ lies on the extension of $[\zeta_2, \zeta_1]$. Because $\zeta_1\in \del \fP$ and $[\zeta_1,\zeta_0]\cap \del \fP$ is connected, we have $(\zeta_1, \zeta_0]\in \fP$.

\end{enumerate}

We remark that in both cases above, $\zeta_1$ is also a vertex of the polygon $\fP$.

Next we show that $[\zeta_1, \zeta_0]$ is tangent to $\fA$ at $\zeta_0$. Combining the discussions above, there are three cases
\begin{enumerate}

  \item $[\zeta_1,\zeta_0] \subset [\zeta_1,\zeta_2]$ for some edge $[\zeta_1,\zeta_2]$ of $\fP$. Next we prove by contradiction that $\zeta_0\in \fA$. Assume that $\zeta_0 \notin \fA$. Since $\zeta_0$ is a corner point of $\fF$, there exists a short vertical segment $[\zeta_0,\zeta_0'] \subset \partial \fF \cap \fP$. If the segment $[\zeta_0,\zeta_0']$ points upward (see the second panel of \Cref{f:segment_case}), then it is on the boundary of $\fF$ and another affine piece with gradient $(0,0)$ with $\fF$ on the left. By \Cref{c:edge}, $H^*=m$ on  $[\zeta_0,\zeta_0']$. Since $\zeta_0\in \del \fP$,  for any $\zeta\in [\zeta_0, \zeta'_0]$,
\begin{align}
M(\zeta)\leq M(\zeta_0)=H^*(\zeta_0)=H^*(\zeta).
\end{align}
Thus on $[\zeta_0, \zeta'_0]$, $H^*=m=M$. This contradicts with \Cref{a:asump} that the trivial set is empty. 
If the segment $[\zeta_0, \zeta_0']$ points downward (see the third panel of \Cref{f:segment_case}), then it is on the boundary of $\fF$ and another affine piece with gradient $(0,0)$ with $\fF$ on the right.   By \Cref{c:edge}, $H^*=M$ on  $[\zeta_0,\zeta_0']$. Since $\zeta_0\in \del \fP$,  for any $\zeta\in [\zeta_0, \zeta'_0]$,
\begin{align}
m(\zeta)\geq m(\zeta_0)=H^*(\zeta_0)=H^*(\zeta).
\end{align}
Thus on $[\zeta_0, \zeta'_0]$, $H^*=m=M$. This contradicts with \Cref{a:asump} that the trivial set is empty. 

We conclude that $\zeta_0\in \fA$. Then by the third statement of \Cref{c:boundary2}, $[\zeta_1, \zeta_0]$ is tangent to $\fA$ at $\zeta_0$. 
  
  \item $[\zeta_1,\zeta_2] \subset [\zeta_1,\zeta_0)$ for some edge $[\zeta_1,\zeta_2]$ of $\fP$. In this case $\zeta_0\in \fP$, and $[\zeta_1, \zeta_0]$ is on the shared boundary of $\fF$ and another affine piece $\fF_1$ with gradient $(1,-1)$. 
  
  Next we prove by contradiction that $\zeta_0\in \fA$. Assume that $\zeta_0 \notin \fA$. Since $\zeta_0$ is a corner point of $\fF$, there exists a short vertical segment $[\zeta_0,\zeta_0'] \subset \partial \fF \cap \fP$. By the same argument as above, it is impossible that the segment $[\zeta_0,\zeta_0']$ points upward.

If the segment $[\zeta_0, \zeta_0']$ points downward, then it is on the boundary of $\fF$ and another affine piece $\fF_2$ with gradient $(0,0)$ with $\fF$ on the right. Then $\fF_1$ and $\fF_2$ also share a piece of boundary $[\zeta_0,\zeta_2']$ with slope $1$. See the fourth panel of \Cref{f:segment_case}. By \Cref{c:edge}, $H^*=m$ on  $[\zeta_0,\zeta_2']$, and $H^*=M$ on $[\zeta_1, \zeta_0]$. Thus $H^*(\zeta_0)=m(\zeta_0)=M(\zeta_0)$. This contradicts with \Cref{a:asump} that the trivial set is empty. 

We conclude that $\zeta_0\in \fA$. Then by the second statement of \Cref{c:boundary2}, $[\zeta_1, \zeta_0]$ is tangent to $\fA$ at $\zeta_0$.

  \item $[\zeta_1,\zeta_0]$ lies on the extension of an edge $[\zeta_2,\zeta_1]$ of $\fP$. The same as in the second case, $[\zeta_1, \zeta_0]$ is tangent to $\fA$ at $\zeta_0$.

\end{enumerate}

%Next we prove the second statement. Without loss of generality, assume that on $\fF$ we have $\nabla H^* = (1,0)$, and $[\zeta_1,\zeta_0]$ is a horizontal line segment with $\fF$ lying above it. The other cases are analogous and will be omitted.
%
%There are several cases
%\begin{enumerate}
%\item $[\zeta_1,\zeta_0] \subset [\zeta'_1,\zeta_2]$ for some edge $[\zeta'_1,\zeta_2]$ of $\fP$
%\end{enumerate}
%
%
%
%
%
%Next we show the right endpoint $\zeta_0\in \fP$. Otherwise $\zeta_0\in\del \fP$. Then for any $\zeta\in [\zeta_1, \zeta_0]$,
%\begin{align}
%m(\zeta)\geq h(\zeta_0)+(1,0)\cdot(\zeta-\zeta_0)=H^*(\zeta)\geq m(\zeta).
%\end{align}
%Thus on $[\zeta_1, \zeta_0]$, $H^*=m=M$. This contradicts with \Cref{a:asump} that the trivial set is empty. 
%Next we prove by contradiction that $\zeta_0\in \fA$. Otherwise, by our construction, the segment $[\zeta_1, \zeta_0]$ cannot be extended to a larger interval still satisfying \eqref{e:z1z0}. Then $\zeta_0$ must be a corner point of $\fF_0$ or $\fF_1$. Without loss of generality, we assume that $\zeta_0$ is a corner point of $\fF_0$. By \Cref{c:boundary1}, the boundary $\fF_0$ consists of portions of the arctic boundary with slope in $(-\infty,0)$, together with horizontal and vertical segments. Thus if $\zeta_0\not\in \fA$, there exists a short vertical segments $[\zeta_0,\zeta_0']\in \del \fF_0\cap \fP$. 
%

\end{proof}

\begin{lemma}\label{l:curvilinear_triangle}
The boundary of each affine piece $\fF$ of $H^*$ consists of two line segments 
$[\zeta,\zeta_1]$ and $[\zeta,\zeta_2]$, tangent to $\fA$ at $\zeta_1$ and $\zeta_2$, respectively, 
together with the portion of the arctic boundary between $\zeta_1$ and $\zeta_2$, 
where $\zeta_1$ and $\zeta_2$ are consecutive tangency points along $\fA$. 
Moreover, $\zeta$ is a vertex of the polygon $\fP$, and the segments 
$[\zeta,\zeta_1]$ and $[\zeta,\zeta_2]$ satisfy the second statement of \Cref{c:boundary3}.

Let the cone $\angle \zeta_1 \zeta \zeta_2$ denote the region containing a small neighborhood 
$B_\delta(\zeta)\cap \fP$ of $\zeta$. This cone has angle greater than $180^\circ$ if the corner $\zeta$ of $\fP$ is not convex. The number of cusp points on the arctic boundary from $\zeta_1$ to $\zeta_2$ is given by
\begin{align}\begin{split}\label{e:cusp_num}
\bm{1}\bigl(\text{$\zeta$ is a concave corner of $\fP$}\bigr)
&+ \bm{1}\bigl(\text{$\fA$ is tangent to $\zeta_1$ on the exterior side of $\angle \zeta_1 \zeta \zeta_2$}\bigr)\\
&+ \bm{1}\bigl(\text{$\fA$ is tangent to $\zeta_2$ on the exterior side of $\angle \zeta_1 \zeta \zeta_2$}\bigr).
\end{split}\end{align}
\end{lemma}

\begin{proof}

By \Cref{c:boundary3}, as we follow the boundary of $\fF$ clockwise from corner point to corner point, we begin with a portion of the arctic boundary from $\zeta_1$ to $\zeta_2$, as in the first statement of \Cref{c:boundary3}, where both $\zeta_1,\zeta_2 \in \fA$ are tangency points.  
The next piece is a line segment $[\zeta_3,\zeta_2]$ satisfying the second statement of \Cref{c:boundary3}, where $\zeta_2$ is a vertex of $\fP$.  
Note that it is possible for $[\zeta_3,\zeta_2]$ to degenerate to a single point.
The following piece is a line segment $[\zeta_3,\zeta_4]$, again satisfying the second statement of \Cref{c:boundary3}, where $\zeta_4 \in \fA$ is a tangency point.  
If $\zeta_4 \neq \zeta_1$, then the next piece, from $\zeta_4$ to $\zeta_5$, is another portion of the arctic boundary, as in the first statement of \Cref{c:boundary3}, with $\zeta_5 \in \fA$ a tangency point, and the process continues.

If the boundary of $\fF$ contains two separate piece of arctic boundary, then $\fP\setminus \overline \fF$ consists of at least two disconnected component, and each of them contains a liquid region. This contradicts to \Cref{a:asump} that the liquid region is connected. Thus in the above procedure, $\zeta_4=\zeta_1$.

Next we prove \eqref{e:cusp_num}. Without loss of generality, we assume that on $\fF$, $\nabla H^*=(1,0)$.
Fix a direction $\omega$ satisfying \eqref{e:omega} with $p_0=(1,0)$. Then the segment $[\zeta, \zeta_1]$ is vertical, and $[\zeta, \zeta_2]$ is horizontal. As we follow the boundary of $\fF$ clockwise once, the tangent direction rotates by $360^\circ$. From $\zeta_1$ to $\zeta_2$, the tangent direction rotates by
\begin{align}
-90^\circ + 180^\circ \cdot \text{\# cusps on $\fA(\zeta_1, \zeta_2)$}.
\end{align}
If $\fA$ is tangent to $\zeta_1$ (or $\zeta_2$) on the interior side of $\angle \zeta_1 \zeta \zeta_2$, the tangent 
direction rotates by $180^\circ$ at $\zeta_1$ (or $\zeta_2$). 
If $\fA$ is tangent to $\zeta_1$ (or $\zeta_2$) on the exterior side of $\angle \zeta_1 \zeta \zeta_2$, the tangent direction does not change. Finally if $\zeta$ is a convex corner of $\fP$, the tangent direction rotates by $90^\circ$ at $\zeta$. Otherwise, if $\zeta$ is a concave corner of $\fP$, the tangent direction rotates by $-90^\circ$ at $\zeta$. We obtain the following relation
\begin{align}\begin{split}
360^\circ &=-90+ 180^\circ \cdot \text{\# cusps on $\fA(\zeta_1, \zeta_2)$}+(90-180^\circ \cdot \bm{1}(\text{$\zeta$ is a concave corner of $\fP$}))\\
&+(180^\circ- 180^\circ\cdot \bm{1}(\text{$\fA$ is tangent to $\zeta_1$ on the exterior side of $\angle \zeta_1 \zeta \zeta_2$}))\\
&+(180^\circ- 180^\circ\cdot \bm{1}(\text{$\fA$ is tangent to $\zeta_2$ on the exterior side of $\angle \zeta_1 \zeta \zeta_2$})).
\end{split}\end{align}
The claim \eqref{e:cusp_num} follows from rearranging the above expression.

\end{proof}

\begin{proof}[Proof of \Cref{t:frozen_structure}]
\Cref{t:frozen_structure} follows from \Cref{l:curvilinear_triangle}.
\end{proof}

\section{Steepest Descent Paths Analysis}
\label{s:path_analysis}
In this section, we prove the results stated in \Cref{s:critical_bulk,s:critical_point} concerning deformations of local descent and ascent paths.

\subsection{Steepest descent/ascent paths}
We recall from \eqref{e:def_action} and \eqref{e:critical}
\begin{align}\begin{split}\label{e:S_copy}
&S(w;x,s)=s\ln s-(x-w)\ln (x-w)-(s-x+w)\ln (s-x+w) -\int_0^w \ln f(u)\rd u,\\
&S'(w;x,s)= \ln\!\left(\frac{w-x}{x-s-w}\right)-\ln f(w).
\end{split}\end{align}
The critical points are characterized by
\begin{align}
x=w+s \chi(w).
\end{align}
The steepest--descent paths are the integral curves of the vector field
\(-\overline{S'(w;x,s)}\). If \(\gamma(t)\) is such a path, then by the chain
rule
\begin{align}\label{e:gradient_flow}
\frac{\rd}{\rd t}S(\gamma(t);x,s)
= S'(\gamma(t);x,s)\,\gamma'(t)
= S'(\gamma(t);x,s)\bigl(-\overline{S'(\gamma(t);x,s)}\bigr)
= -\lvert S'(\gamma(t);x,s)\rvert^{2}\leq 0.
\end{align}
Consequently,
\[
\frac{\rd}{\rd t}\Im S(\gamma(t);x,s)=0,
\quad
\frac{\rd}{\rd t}\Re S(\gamma(t);x,s)=-\lvert S'(\gamma(t);x,s)\rvert^{2}\le 0,
\]
i.e. \(\Im S\) is constant and \(\Re S\) is nonincreasing along steepest--descent paths.

\begin{proposition}[Exit directions near a critical point]\label{prop:exit}
Fix $m\geq 1$, $d\in \bC\setminus \{0\}$ and a base point \(w_0\). Assume that, uniformly on $|w-w_0|=r$, for some small $\delta\leq |d|r^m$,
\begin{align}\label{e:Sclose}
\left|S(w;x,s)-S(w_c;x,s) -d (w-w_0)^{m}\right|\leq \delta.
\end{align}

Let \(\gamma\) be any steepest--descent or steepest--ascent curve of
\(S(\,\cdot\,;x,s)\) through \(w_c\).
Then \(\gamma\) first intersects the circle \(|w-w_0|=r\) at points in the narrow arcs
\begin{align}\label{e:Aj}
A^{(2m)}_j(w_0, r)=\left\{w_0+r e^{\ri \theta}: \left|\theta-\frac{(j+1)\pi-\arg d}{m}\right|
\leq\frac{\pi\delta}{2m |d|r^{m}}\right\},
\quad j=0,1,\dots,2m-1.
\end{align}
Moreover, at those intersection points $w\in A^{(2m)}_j(w_0,r)$,
\begin{align}\label{e:final}
\left|\Re\!\big(S(w;x,s)-S(w_c;x,s)\big)
-(-1)^{j+1}\,|d|\,r^{\,m}\right|\leq 3\delta,
\end{align}
so the sign agrees with that predicted by
the model term \(d(w-w_0)^m\).
\end{proposition}

\begin{proof}[Proof of \Cref{prop:exit}]
%Fix \(w=w_0+r e^{\ri \theta}\). Subtract and add \(S(w;x_0,s_0)\) and use
%\eqref{eq:param-stability} and \eqref{eq:Taylor}:
%\begin{align}\label{e:Sdiff}
%S(w;x,s)-S(w_c;x,s)
%=a r^{m}e^{\ri m\theta}+\OO(\delta).
%\end{align}
Along any steepest--descent or steepest--ascent curve through \(w_c\),
\(\Im S(\cdot;x,s)\) is constant, hence the intersection points $w=w_0+r e^{\ri \theta}$ on \(|w-w_0|=r\) must satisfy
\begin{align}\begin{split}\label{e:ImSdiff}
\big|\Im\big(d\,r^{\,m}e^{\ri m\theta}\big)\big|
=\left|\Im[S(w;x,s)-S(w_c;x,s) -d (w-w_0)^{m}]\right|\leq  \delta,
\end{split}\end{align}
where the first equality follows from $\Im[S(w;x,s)]=\Im[S(w_c;x,s)]$, and the inequality follows from taking the imaginary part of \eqref{e:Sclose}.

Dividing by \(|d|\,r^{\,m}\), \eqref{e:ImSdiff} gives  
\[
\left|\sin(m\theta+\arg d)\right|
\le \delta/(|d|r^{m})\leq 1,
\]
where the last inequality follows from our assumption that $\delta\leq |d|r^m$.
This confines \(\theta\) to \(2m\) short arcs 
\begin{align}\label{e:diffsmall}
\left|m\theta+\arg d-(j+1)\pi \right|
 \leq\frac{\pi\delta}{2 |d|r^{m}},\quad j=0,1,2,\cdots, 2m-1.
\end{align}
The above arcs correspond to the sets $A_j^{(2m)}$ defined in \eqref{e:Aj}.
On those arcs, for $w=w_0+r e^{\ri \theta}\in A_j^{(2m)}(w_0, r)$, we have
\begin{align}\label{e:real_part0}
\left|\Re[d r^m e^{\ri m \theta}]-(-1)^{j+1} |d|r^m\right|
=\left|(-1)^{j+1} |d|r^m \big(\cos(m\theta+\arg d-(j+1)\pi)-1\big)\right|
\leq \pi \delta/2,
\end{align}
where the first equality follows from taking the real part, and the inequality follows from $|\cos(x)-1|\leq |x|$ and \eqref{e:diffsmall}. Hence
\begin{align}\begin{split}\label{e:real_part}
&\phantom{{}={}}\left|\Re\!\big(S(w;x,s)-S(w_c;x,s)\big)
-(-1)^{j+1}\,|d|\,r^{\,m}\right|\\
&\leq \left|\Re\!\big(S(w;x,s)-S(w_c;x,s)-d r^m e^{\ri m \theta}\big)\right|
+\left|\Re[d r^m e^{\ri m \theta}]-(-1)^{j+1} |d|r^m\right|\\
&\leq \delta+\pi\delta/2\leq 3\delta.
\end{split}\end{align}
where the first inequality follows from the triangle inequality, and the second inequality follows from \eqref{e:real_part0} and \eqref{e:Sclose}. This gives \eqref{e:final}, and the ``even'' arcs
\(j=0,2, 4,\dots,2m-2\) correspond to descent, with negative real increment, while the
``odd'' arcs correspond to ascent. 
\end{proof}

Next, we verify that the assumptions of \Cref{prop:exit} hold in all cases
discussed in  \Cref{s:critical_bulk} and \Cref{s:critical_point}. There are two cases.

\begin{enumerate}
  \item In all cases discussed in \Cref{s:critical_bulk} and \Cref{s:critical_point}, except when
  $(x_0,s_0)$ lies on a horizontal extended side and $\fU$ is a
  tangent/cusp-turning/frozen chart, the point $w_0$ is a critical point of
  finite order $m\geq 1$. That is,
  \begin{equation}\label{eq:Taylor}
    S(w;x_0,s_0)
    =
    S(w_0;x_0,s_0)
    + d\,(w-w_0)^m
    + \OO\!\bigl(|w-w_0|^{m+1}\bigr),
    \qquad |d|\asymp 1.
  \end{equation}

  \item Suppose the line through $(x_0,s_0)$ is tangent to the arctic curve
  at a horizontal tangency point $(x_0',s_0)$; see
  \Cref{s:horizontal_tangent} and \Cref{s:horizontal_frozen_neighborhood}. In the
  case of \Cref{s:horizontal_tangent}, we have
  $(x'_0,s_0)=(x_0,s_0)$. The critical point of
  $S(\,\cdot\,;x_0,s_0)$ lies at infinity. Introducing the change of
  variables
  \[
    \wt w=\frac{1}{x_0'-w},
  \]
  as in \eqref{e:change_coordinate} or \eqref{e:frozen_change_coordinate}, the point $\wt w=0$ becomes a
  critical point of finite order $m\geq 1$. Equivalently,
  \begin{equation}\label{eq:Taylor2}
    \wt S(\wt w;x_0,s_0)
    =
    \wt S(0;x_0,s_0)
    + d\,\wt w^m
    + \OO\!\bigl(\wt w^{m+1}\bigr),
    \qquad |d|\asymp 1.
  \end{equation}
\end{enumerate}

The following lemma gives estimates for
$S(w;x,s)-S(w_c;x,s)$ and
$\wt S(\wt w;x,s)-\wt S(\wt w_c;x,s)$ when $(x,s)$ is sufficiently close to
$(x_0,s_0)$.

\begin{lemma}\label{c:S_diff}
Assume \eqref{eq:Taylor} in the finite critical point case and
\eqref{eq:Taylor2} in the horizontal tangency case. Fix small
$\delta,r>0$, and suppose that
\[
\|(x,s)-(x_0,s_0)\|_2\leq \delta.
\]
Then the following estimates hold.

\begin{enumerate}
\item For any $|w-w_0|=r$ and $|w_c-w_0|\leq r$,
\begin{align}\label{e:Scritical1}
S(w;x,s)-S(w_c;x,s)
=
d(w-w_0)^m
+
\OO\!\bigl(
  |d|\,|w_c-w_0|^m
  + r^{m+1}
  + \delta\log(1/\delta)
\bigr),
\end{align}
provided
\[
|s_0|+|x_0|+|w_0|\leq 1/\delta .
\]

\item For any $|\wt w|=r$ and $0<|\wt w_c|\leq r$,
\begin{align}\label{e:Scritical2}
\wt S(\wt w;x,s)-\wt S(\wt w_c;x,s)
=
d \wt w^m
+
\OO\!\bigl(
  |d|\,|\wt w_c|^m
  + r^{m+1}
  + \delta
  + |s-s_0|\log(1/|\wt w_c|)
\bigr),
\end{align}
provided
\[
|x_0-x_0'|r\leq \frac14,
\qquad
|s_0|r\leq \frac14 .
\]
\end{enumerate}
\end{lemma}

\begin{proof}[Proof of \Cref{c:S_diff}]
We first record a useful estimate. Let $u,v$ be complex numbers lying in the
same half-plane, and assume that $|u-v|\leq \delta$. Then
\begin{equation}\label{e:u-v0}
  |u\log u-v\log v|
  \leq
  C\delta\bigl(\log(1/\delta)+\log(1+|u|)\bigr).
\end{equation}

Indeed, if $|u|\leq 2\delta$, then $|v|\leq 3\delta$, and therefore
\begin{equation}\label{e:u-v1}
  |u\log u-v\log v|
  \leq |u\log u|+|v\log v|
  \leq C\delta\log(1/\delta).
\end{equation}
On the other hand, if $|u|\geq 2\delta$, then the straight segment from $u$ to
$v$ remains in the same half-plane and stays a distance at least $|u|/2$ from
the origin. Thus, fixing a branch of $\log$ on this half-plane,
\begin{equation}\label{e:u-v2}
  |u\log u-v\log v|
  =
  \left|\int_u^v (\log \zeta+1)\,\rd\zeta\right|
  \leq
  C\delta(1+|\log u|).
\end{equation}
Combining the two cases gives \eqref{e:u-v0}.

We now prove \eqref{e:Scritical1}. For
$\|(x,s)-(x_0,s_0)\|_2\leq \delta$ and $|u-w_0|\leq r$, we have
\[
|s-s_0|\leq \delta,
\qquad
|(x-u)-(x_0-u)|\leq \delta,
\qquad
|(s-x+u)-(s_0-x_0+u)|\leq 2\delta .
\]
Recalling \eqref{e:S_copy}, we get
\begin{align}
\begin{split}\label{e:Sdiff}
S(u;x,s)-S(u;x_0,s_0)
&=
\bigl(s\log s-(x-u)\log(x-u)
       -(s-x+u)\log(s-x+u)\bigr) \\
&\quad
-\bigl(s_0\log s_0-(x_0-u)\log(x_0-u)
       -(s_0-x_0+u)\log(s_0-x_0+u)\bigr).
\end{split}
\end{align}
Applying \eqref{e:u-v0} to the three logarithmic terms gives
\begin{align}\label{e:Sdiff1}
|S(u;x,s)-S(u;x_0,s_0)|
&\leq
C\delta\bigl(
  \log(1/\delta)
  +\log(1+|s_0|+|x_0|+|w_0|)
\bigr) \leq 2C\delta\log(1/\delta),
\end{align}
provided we take $\delta$ small enough  such that 
$1+|s_0|+|x_0|+|w_0|\leq 1/\delta$.

Applying \eqref{e:Sdiff1} with $u=w$ and $u=w_c$, we obtain
\[
S(w;x,s)-S(w_c;x,s)
=
S(w;x_0,s_0)-S(w_c;x_0,s_0)
+
\OO\!\bigl(\delta\log(1/\delta)\bigr).
\]
By \eqref{eq:Taylor},
\begin{align}
S(w;x_0,s_0)-S(w_c;x_0,s_0)
&=
d(w-w_0)^m
-d(w_c-w_0)^m 
+\OO\!\bigl(r^{m+1}+|w_c-w_0|^{m+1}\bigr).
\end{align}
Since $|w_c-w_0|\leq r$, this gives \eqref{e:Scritical1}.

It remains to prove \eqref{e:Scritical2}. From \eqref{e:S_copy}, after the
change of variables $u=1/(x_0'-\wt w)$, we have
\begin{align}
\begin{split}\label{e:Sdiff2}
&\phantom{{}={}}
\wt S(u;x,s)-\wt S(u;x_0,s_0) \\
&=
s\log s-s_0\log s_0 
+(x_0-x_0'+1/u)\log(x_0-x_0'+1/u)
 -(x-x_0'+1/u)\log(x-x_0'+1/u) \\
&
+(s_0-x_0+x_0'-1/u)\log(s_0-x_0+x_0'-1/u) 
-(s-x+x_0'-1/u)\log(s-x+x_0'-1/u).
\end{split}
\end{align}
For $|u|\leq r$, the assumptions
$
|x_0-x_0'|r\leq 1/4$,
and $
|s_0|r\leq1/4
$
imply
$
|x_0-x_0'|\,|u|\leq 1/4$,
and 
$
|s_0-x_0+x_0'|\,|u|\leq \frac12$.
Therefore,
\begin{align}
\begin{split}\label{e:two_term}
&\phantom{{}={}}
(x_0-x_0'+1/u)\log(x_0-x_0'+1/u)
-(x-x_0'+1/u)\log(x-x_0'+1/u) \\
&=
(x_0-x)\log(1/u)+\OO(|x-x_0|),
\end{split}
\end{align}
and similarly
\begin{align}
\begin{split}\label{e:two_term2}
&\phantom{{}={}}
(s_0-x_0+x_0'-1/u)\log(s_0-x_0+x_0'-1/u) 
-(s-x+x_0'-1/u)\log(s-x+x_0'-1/u) \\
&=
(s_0-s+x-x_0)\log(1/u)
+\OO(|x-x_0|+|s-s_0|).
\end{split}
\end{align}
Substituting \eqref{e:two_term} and \eqref{e:two_term2} into
\eqref{e:Sdiff2}, we obtain
\begin{align}\label{e:Sdiff3}
\wt S(u;x,s)-\wt S(u;x_0,s_0)
=
s\log s-s_0\log s_0
+
(s_0-s)\log(1/u)
+
\OO(|x-x_0|+|s-s_0|).
\end{align}
Using again $|s_0|r\leq 1/4$, we may bound the first term by the same type of
logarithmic estimate. Hence, for $0<|u|\leq r$,
\begin{align}\label{e:Sdiff4}
|\wt S(u;x,s)-\wt S(u;x_0,s_0)|
\leq
C\bigl(\delta+|s-s_0|\log(1/|u|)\bigr).
\end{align}

Applying \eqref{e:Sdiff4} with $u=\wt w$ and $u= \wt w_c$, and using $|\wt w|=r$ and
$|\wt w_c|\leq r$, gives
\begin{align}
\wt S(\wt w;x,s)-\wt S(\wt w_c;x,s)
&=
\wt S(\wt w;x_0,s_0)-\wt S(\wt w_c;x_0,s_0) 
+\OO\!\bigl(\delta+|s-s_0|\log(1/|\wt w_c|)\bigr).
\end{align}
Finally, \eqref{eq:Taylor2} gives
\[
\wt S(\wt w;x_0,s_0)-\wt S(\wt w_c;x_0,s_0)
=
d \wt w^m
-d \wt w_c^m
+
\OO(r^{m+1}+|\wt w_c|^{m+1}).
\]
Since $|\wt w_c|\leq r$, this proves \eqref{e:Scritical2}.
\end{proof}

\subsection{Perturbation estimates}

In this section, we collect several results on perturbations of linear, quadratic, and cubic equations that will be used to analyze the critical-point equation. Their proofs are deferred to \Cref{s:perturb}.

\begin{lemma}\label{lem:bulk-perturbation}
Given $b\in \bC$ and a power series $\sum_{k\ge 2} c_k z^k$ that is analytic on the open disk $\{z\in \bC:|z|<\fc_0\}$, consider the equation
\begin{equation}\label{eq:bulk_F}
F(z;b):=-z-b+\sum_{k\ge 2} c_k z^k=0.
\end{equation}
There exists a constant $\fc>0$ depending on $\fc_0$ and the power series such that the following holds. In the following statements, all implicit constants depend only on $\fc_0$ and the power series.

If $|b|\leq \fc/2$,
then the equation \eqref{eq:bulk_F} has exactly one solution $z_0$ in the disk $|z|\le \fc$, and
\begin{align}\label{e:bulk_small_real_z0}
z_0=-b+\OO\left(|b|^2\right).
\end{align}
If $b\in \bR$ and $c_k\in\bR$ for all $k$, then $z_0\in \bR$.
\end{lemma}

\begin{lemma}\label{lem:edge-perturbation}
Given $a,b\in \bR$ and a power series $\sum_{k\ge 3} c_k z^k$ with real coefficients that is analytic on the open disk $\{z\in \bC:|z|<\fc_0\}$, consider the equation
\begin{equation}\label{eq:F}
F(z;a,b):=z^2-a z-b+\sum_{k\ge 3} c_k z^k=0.
\end{equation}
There exists a constant $\fc>0$ depending on $\fc_0$ and the power series such that the following holds. In the following statements, all implicit constants depend only on $\fc_0$ and the power series.

If $\rho_0:=\max\{|a|/4,\,|b|^{1/2}\}\le \fc$, the equation \eqref{eq:F} has exactly two solutions (counting multiplicity) in the disk $|z|\le 4\fc$, and both are bounded in absolute value by $4\rho_0$, and lower bounded by $\min\{|b|/a, \rho_0\}/8$.

Moreover, there exist $b'=b'(a)\asymp -a^2$ and $z'=z'(a)\asymp a$, such that $F(\,\cdot\,;a,b')$ has a double real root at $z'$.

\begin{enumerate}
\item \textbf{Quadratic--constant balance ($b\ge b'$).} $F(\,\cdot\,;a,b)=0$ has two real roots $z_0\le z'\le z_1$ such that
\begin{align}\label{e:edge_realroot}
z'-z_0,\; z_1-z' \asymp \sqrt{b-b'}.
\end{align}
%Moreover
%\begin{enumerate}
%\item If $a\geq 0$ and $|b|\geq |b'|/2$ then $|z_0|, |z_1|\asymp \rho_0$;
%\item If $a\geq 0$ and $|b|\leq |b'|/2$ then $z_0|\asymp-b/a$,
%and $z_1\asymp \rho_0$;
%\item If $a\leq 0$ and $|b|\geq |b'|/2$ then $|z_0|, |z_1|\asymp \rho_0$;
%\item If $a\leq 0$ and $|b|\leq |b'|/2$ then $z_0\asymp -\rho_0$ and $z_1\asymp -b/a$.
%\end{enumerate}

\item \textbf{Linear--quadratic balance ($b<b'$).}
There is a complex conjugate pair $z_{\pm}$ with $|z_\pm|\asymp \rho_0$ and
\begin{align}\label{e:edge_complexroot}
|z_\pm-z'|\asymp |\Im z_\pm|\asymp \sqrt{b'-b}.
\end{align}
\end{enumerate}
\end{lemma}

\begin{lemma}\label{lem:cubic-perturbation}
Given $a,b\in \bR$ and a power series $\sum_{k\ge 4} c_k z^k$ with real coefficients that is analytic on the open disk $\{z\in \bC:|z|<\fc_0\}$, consider the equation
\begin{equation}\label{eq:cusp_F}
F(z;a,b):=z^3-a z-b+\sum_{k\ge 4} c_k z^k=0.
\end{equation}

There exists a constant $\fc>0$ depending on $\fc_0$ and the power series such that the following holds. In the following statements, all implicit constants depend only on $\fc_0$ and the power series.

%If $|b|^{1/3}\leq \fc \leq \sqrt{|a|/4}$,
%then the equation \eqref{eq:cusp_F} has exactly one solution $z_0$ in the disk $|z|\le \fc$, this root is real, and
%\begin{align}\label{e:cusp_small_real_z0}
%z_0=-\frac{b}{a}+\OO\left(\frac{b^3}{|a|^4}\right).
%\end{align}

If $\rho_0:=\max\{\sqrt{|a|/4},\,|b|^{1/3}\}\le \fc$, then
the equation \eqref{eq:cusp_F} has exactly three solutions (counting multiplicity) in the disk $|z|\le 4\fc$, and all are bounded in absolute value by $4\rho_0$, and lower bounded by $\min\{|b|/a, \rho_0\}/4$.

Moreover, the qualitative location of the three roots is as follows.
\begin{enumerate}
\item If $a<0$ (the ``one real + one complex pair'' regime), there are two subregimes:
\begin{enumerate}
\item \textbf{Linear--cubic balance} ($|a|\gtrsim |b|^{2/3}$). There is one real root
\begin{align}\label{e:realroot1}
z_0\asymp -\frac{b}{a},
\end{align}
and a complex conjugate pair $z_{\pm}$ with
\begin{align}\label{e:complexroot1}
|z_\pm|\asymp \sqrt{|a|},\quad |\Im z_\pm|\asymp \sqrt{|a|}.
\end{align}
% Moreover, if $|a|/4\geq |b|^{2/3}$, the real root $z_0$ satisfies the refined estimate 
%\begin{align}\label{e:cusp_small_real_z0}
%z_0=-\frac{b}{a}+\OO\left(\frac{b^3}{|a|^4}\right).
%\end{align}

\item \textbf{Cubic--constant balance} ($|a|\lesssim |b|^{2/3}$).
There is one real root of size $z_0\asymp b^{1/3}$ and a complex conjugate pair $z_{\pm}$ with
\begin{align}\label{e:realroot2}
|z_\pm|\asymp |b|^{1/3},\quad |\Im z_\pm|\asymp |b|^{1/3}.
\end{align}
%All roots \(z\) satisfy
%\(|F'(z;a,b)|\asymp |b|^{2/3}\).
\end{enumerate}

\item If $a>0$, there exist $b'_\pm=b'_\pm(a)$ and $z'_\pm=z'_\pm(a)$ with
\[
b_-'\asymp -a^{3/2},\quad b_+'\asymp a^{3/2}, \quad z_+'\asymp \sqrt{a}, \quad z_-'\asymp -\sqrt{a},
\]
such that $F(\,\cdot\,;a,b'_{\pm})$ has a double real root at $z'_{\mp}$.

\begin{enumerate}
\item \textbf{Linear--cubic balance} ($b'_-\le b\le b'_+$). The equation $F(\,\cdot\,;a,b)=0$ has three real roots $z_0\le z_1\le z_2$ such that
\begin{align}\label{e:realroot3}
z_0\asymp -\sqrt{a}, \quad z_1\asymp -\frac{b}{a},\quad
z_2\asymp \sqrt{a}.
\end{align}
%If $0\le b\le b'_+$, then $|F'(z_2;a,b)|\asymp a$, and
%\begin{align}\label{e:close_critical}
%z_1-z'_-,\; z'_- -z_0\asymp (b'_+ - b)^{1/2}a^{-1/4},\quad |F'(z_0;a,b)|, |F'(z_1;a,b)|\asymp (b'_+ - b)^{1/2}a^{1/4}.
%\end{align}
%If $b'_-\le b\le 0$, then $|F'(z_0;a,b)|\asymp a$, and
%\begin{align}
%z_2-z'_+,\; z'_+-z_1\asymp (b - b'_-)^{1/2}a^{-1/4},\quad |F'(z_1;a,b)|, |F'(z_2;a,b)|\asymp (b - b_-')^{1/2}a^{1/4}.
%\end{align}
%Moreover, if $|a|/4\geq |b|^{2/3}$, the real root $z_1$ satisfies the refined estimate \eqref{e:cusp_small_real_z0}.

\item \textbf{Cubic--constant balance} ($b< b'_-$ or $b> b'_+$).
There is one real root $z_0\asymp b^{1/3}$ and a complex conjugate pair $z_{\pm}$ with $|z_\pm|\asymp |b|^{1/3}$.
% If $b> b'_+$, then
%\begin{align}\label{e:complexroot2}
%|z_\pm-z'_-|\asymp |\Im z_\pm|\asymp \min \left\{|b|^{1/3}, \sqrt{\frac{b-b'_+}{\sqrt{a}}}\right\}.
%\end{align}
%If $b< b'_-$, then
%\begin{align}
%|z_\pm-z'_+|\asymp |\Im z_\pm|\asymp \min \left\{|b|^{1/3}, \sqrt{\frac{b'_--b}{\sqrt{a}}}\right\}.
%\end{align}
%In both cases, the real root satisfies
%\[
%|F'(z_0;a,b)|\asymp |b|^{2/3},
%\]
%while the complex roots satisfy
%\[
%\color{red}
%|F'(z_\pm;a,b)|\asymp |b|^{1/3}|\Im z_\pm|.
%\]
\end{enumerate}
\end{enumerate}
\end{lemma}

\subsection{Liquid chart}
\label{s:liquid_chart_proof}

We recall from \Cref{s:critical_bulk} that, given a  liquid chart
centered at $(x_0,s_0)$ and setting $w_0 = x_0 - s_0 \chi(w_0) \in \bC_+$, the 
point $w_0$ is a nondegenerate saddle point of $S(w;x_0,s_0)$, with $m=2$ in 
\Cref{prop:exit}. For any $(x,s)$ in this regular liquid neighborhood, we study 
$S(w;x,s)$ in a small chart of $w_0$.

\begin{proof}[Proof of \Cref{c:bulk}]
By the first statement of \Cref{p:surface}, after decreasing $\fc$ if necessary,
the Riemann surface $\cC$ can be parametrized over $\fU$ as $(f(w),w)$, and
both $f$ and $\chi$ are holomorphic on $\fU$. Moreover,
\begin{align}\begin{split}\label{e:chip}
\Im \chi(w_0)
&= \Im \chi(x_0,s_0)
 = \Im\left[\frac{f(x_0,s_0)}{f(x_0,s_0)+1}\right]
 = \frac{\Im f(x_0,s_0)}{|f(x_0,s_0)+1|^2}
 \asymp -1, \\
\Im w_0
&= \Im\bigl[x_0-s_0\chi(w_0)\bigr]
 = -s_0\Im \chi(w_0)\asymp 1,\quad |1+s_0\chi'(w_0)|\asymp 1,
\end{split}\end{align}
where the last statement is from \eqref{e:derchi_liquid}.
Thus, after possibly decreasing $\fc$ again, the bounds in
\eqref{e:bulk_bounds} hold uniformly for all $w\in \fU$.

The critical points $w_c$ of $S(w;x,s)$ satisfy
\begin{align}\label{e:bulk_x_exp}
x= w+s\chi(w) = w_0+s\chi(w_0)+(w-w_0)
   + s\sum_{k\geq 1}\frac{a_k(w-w_0)^k}{k!},
\end{align}
where $a_k=\chi^{(k)}(w_0)$. Since $x_0=w_0+s_0\chi(w_0)$, this gives
\begin{align}
(x-x_0)-(s-s_0)\chi(w_0)
= (1+s\chi'(w_0))(w-w_0)
  + s\sum_{k\geq 2}\frac{a_k(w-w_0)^k}{k!}.
\end{align}
Hence, for
$|s-s_0|\leq \delta$ with $\delta$ sufficiently small, we also have
$|1+s\chi'(w_0)|\asymp 1$. Therefore the critical point equation can be written as
\begin{align}
\frac{(x-x_0)-(s-s_0)\chi(w_0)}{1+s\chi'(w_0)}
-(w-w_0)
-\sum_{k\geq 2}
\frac{s a_k(w-w_0)^k}{(1+s\chi'(w_0))k!}
=0.
\end{align}
Equivalently, setting
\begin{align}\label{e:defb}
z=w-w_0,
\qquad
b=\frac{-(x-x_0)+(s-s_0)\chi(w_0)}{1+s\chi'(w_0)},
\end{align}
we obtain an equation of the form
\begin{align}\label{e:bulk_F}
F(z;b):=-z-b+\sum_{k\geq 2}c_k z^k=0,
\quad
c_k=-\frac{s a_k}{(1+s\chi'(w_0))k!}.
\end{align}
Thus \Cref{lem:bulk-perturbation} implies that $F(z;b)$ has exactly one zero in the disk
$|z|\leq \fc$, and  this zero satisfies
\begin{align}\label{e:S_diff_critical}
|w_c-w_0|
\lesssim |b|
\lesssim \|(x,s)-(x_0,s_0)\|_2
\leq \delta,
\end{align}
where the first statement follows from
\eqref{e:bulk_small_real_z0}, and the second statement follows from \eqref{e:chip}, \eqref{e:defb} and $|1+s\chi'(w_0)|\asymp 1$. This gives \eqref{e:wcbb}. 

We now recall the expansion of $S(w;x_0,s_0)$ and $S''(w;x_0,s_0)$ from
\eqref{e:liquid_neighborhood}. Then \eqref{e:chip} implies
\begin{align}
|d|=\frac12| S''(w_0;x_0,s_0)|\asymp 1. 
\end{align}
Combining \eqref{e:S_diff_critical} with
\eqref{e:Scritical1}, we obtain, uniformly for $|w-w_0|=r$,
\begin{align}
S(w;x,s)-S(w_c;x,s)
=
d(w-w_0)^2
+\OO\bigl(\delta^2+r^3+\delta\log(1/\delta)\bigr)
.
\end{align}
Taking $r=\fc$ and then choosing $\fc$ and $\delta$ sufficiently small gives
\eqref{e:bulk_err}.
\end{proof}

\begin{proof}[Proof of \Cref{c:bulk_steepest}]
By \Cref{c:bulk}, uniformly for $|w-w_0|=\fc$,
\begin{align}\label{e:bulk_steepest_expansion}
|S(w;x,s)-S(w_c;x,s)
-
d(w-w_0)^2|
\leq \frac{|d|\fc^2}{100},
\quad
d=\frac12 S''(w_0;x_0,s_0),
\end{align}
with $|d|\asymp 1$. The hypotheses of
\Cref{prop:exit} hold with $(m, \delta, r)$ taken to be $(2, d\fc^2/100, \fc)$.

Therefore two steepest--descent and two steepest--ascent branches emanate from
$w_c$, forming the usual cross pattern of a nondegenerate saddle. Moreover,
these four trajectories intersect the circle $|w-w_0|=\fc$ alternately through
the arcs $A^{(4)}_0(w_0;\fc),\cdots,A^{(4)}_3(w_0;\fc)$, and uniformly for $w\in A^{(4)}_j(w_0;\fc)$,
\begin{align}\label{e:final_copy}
\left|\Re [S(w;x,s)-S(w_c;x,s)]
-
(-1)^{j+1} |d|\fc^2\right|\leq \frac{3|d|\fc^2}{100}
\end{align}

\noindent\textbf{Replacing \(\mathsf C^{\rm d}(w_0)\) by \(\mathsf D^{\rm d}(w_c)\).}
Let $\mathsf D^{\rm d}(w_c)$ denote the union of the two steepest--descent trajectories
starting at $w_c$, stopped at their first exit from the disk
$
\{w: |w-w_0|\leq \fc\}$.
By the preceding paragraph, these two trajectories exit transversely through
the boundary arcs $A^{(4)}_0(w_0;\fc)$ and $A^{(4)}_2(w_0;\fc)$. Since the integrand is
holomorphic in the deformation region, Cauchy's theorem allows us to deform
$\mathsf C^{\rm d}(w_0)$ to $\mathsf D^{\rm d}(w_c)$, together with subarcs of
$A^{(4)}_0(w_0;\fc)\cup A^{(4)}_2(w_0;\fc)$.

For $w\in A^{(4)}_0(w_0;\fc)\cup A^{(4)}_2(w_0;\fc)$, \eqref{e:final_copy} gives
\[
\Re[S(w;x,s)-S(w_c;x,s)]
\leq - |d|\fc^2
+\frac{3|d|\fc^2}{100}\leq -\fc',
\]
 for some $\fc'>0$. 

%\[
%\OO\!\left(\bigl|e^{nS(w_c;x,s)}\bigr|e^{-n\fc'}\right),
%\]
%and is exponentially small.
\noindent\textbf{Replacing \(\mathsf D^{\rm d}(w_c)\) by \(\mathsf S^{\rm d}(w_c)\)}
Let
$
r_n:={\ln n}/{\sqrt n}$.
For $n$ sufficiently large, $r_n\ll \fc$, and by \eqref{e:S_diff_critical} $w_c=w_0+\OO(\delta)$, the
disk $\{w:|w-w_c|\le r_n\}$ is contained in $\{w:|w-w_0|\le \fc\}$. We recall that 
$\mathsf S(w_c)$ are the portions of the two steepest--descent trajectories
starting at $w_c$ and stopped at their first exit from the disk
$
\{w: |w-w_c|\leq r_n\}$.

Locally around $w_c$, the function $S(w;x,s)$ is holomorphic.
We recall the expansion of $S(w;x,s)$ and $S''(w;x,s)$ from
\eqref{e:liquid_neighborhood}. 
\begin{align}\label{e:exp_Sbulk}
S(w;x,s)
=
S(w_c;x,s)
+d(w-w_c)^2
+\OO(|w-w_c|^3),
\end{align}
where 
\begin{align}
d:=S''(w_c;x,s)
=
-\frac{\chi'(w_c)+1/s}{\chi(w_c)(1-\chi(w_c))}
=-(1+\OO(\delta))\frac{\chi'(w_0)+1/s}{\chi(w_0)(1-\chi(w_0))}
\asymp 1,
\end{align}
where we used $w_c=w_0+\OO(\delta)$ from \eqref{e:S_diff_critical}, and \eqref{e:chip}. 

Applying \Cref{prop:exit} once more, now with $(m, \delta, r)$ taken to be $(2, Cr_n^3, r_n)$. 
We conclude that the two steepest--descent branches from $w_c$ exit the circle
$|w-w_c|=r_n$ through the arcs $A^{(4)}_0(w_c;r_n)$ and $A^{(4)}_2(w_c;r_n)$. For any
point on these arcs,
\begin{align}\label{eq:bulk-drop0}
n\,\Re[S(w;x,s)-S(w_c;x,s)]
=
-n\bigl(|d|r_n^2+\OO(r_n^3)\bigr) \leq -\fc'(\ln n)^2,
\end{align}
for some $\fc'>0$. Indeed, $nr_n^2=(\ln n)^2$, while
$nr_n^3=\oo((\ln n)^2)$.

Since $\Re S(w;x,s)$ decreases along steepest--descent trajectories as one
moves away from $w_c$, the bound \eqref{eq:bulk-drop0} holds throughout
$\mathsf D^{\rm d}(w_c)\setminus \mathsf S^{\rm d}(w_c)$. Therefore for all \(w\in \mathsf D^{\rm d}(w_c)\setminus
\mathsf S^{\rm d}(w_c)\),
\[
e^{n\Re[S(w;x,s)]}
\le
e^{n\Re[S(w_c;x,s)]}e^{-\fc'(\ln n)^2}.
\]

\end{proof}

\subsection{Arctic chart}\label{s:arctic_chart_proof}

We recall from \Cref{s:critical_arctic} that, given a regular arctic
neighborhood centered at $(x_0,s_0)\in\fA$ and setting
$
w_0=x_0-s_0\chi(w_0)\in \cC(\bR)$,
the point $w_0$ is a cubic saddle point of $S(w;x_0,s_0)$. Thus $m=3$ in
\Cref{prop:exit}. For any $(x,s)$ in this regular arctic neighborhood, we
study $S(w;x,s)$ in a small neighborhood of $w_0$.

%
%{\color{red}
%More generally, let $(x_0,s_0)\in \fA$ be a non-horizontal tangency point, that
%is, assume $f(x_0,s_0)\neq -1$. By the second statement of \Cref{p:surface},
%there exists an analytic function $\chi(w)$ defined in a small neighborhood of
%$w_0$, and $1/\chi(w_0)$ is the slope of the tangent line to the arctic curve
%at $(x_0,s_0)$. For a fixed value of $s$, set
%\[
%x_0'(s):=w_0+s\chi(w_0).
%\]
%Then $(x_0',s)$ lies on the tangent line through
%$(x_0,s_0)$; see the first panel of \Cref{f:cubic_saddle}.
%}

\begin{proof}[Proof of \Cref{c:arctic_critical}]
By the second statement of \Cref{p:surface}, after decreasing $\fc$ if necessary,
the Riemann surface $\cC$ can be parametrized over $\fU$ as $(f(w),w)$, and
both $f$ and $\chi$ are holomorphic on $\fU$. Moreover, since $(x_0,s_0)$ is
bounded away from cusp and tangency points, we have
\begin{align}\label{e:chip_arctic}
|\chi(w_0)|,\ |1-\chi(w_0)|,\ |\chi''(w_0)| \asymp 1 .
\end{align}
Also, since $w_0$ is a cubic critical point of $S(\,\cdot\,;x_0,s_0)$, we have
\[
x_0=w_0+s_0\chi(w_0),
\qquad
\chi'(w_0)=-\frac1{s_0}.
\]
Thus, after possibly decreasing $\fc$ again, the bounds in \eqref{e:arctic_bounds} hold uniformly for all $w\in \fU$.
And \eqref{e:arctic_cubic} follows from \eqref{e:derSthird}.

Recall from \eqref{e:slope} that $1/\chi(w_0)$ is the slope of the tangent line
to the arctic curve at $(x_0,s_0)$. For a fixed value of $s$, set
\begin{align}\label{e:x_0'}
x_0'=x_0'(s):=w_0+s\chi(w_0).
\end{align}
Then $(x_0',s)$ lies on the tangent line through $(x_0,s_0)$; see the first
panel of \Cref{f:cubic_saddle}.

Let $a_k=\chi^{(k)}(w_0)$ and set $z=w-w_0$. Critical points $w_c$ of
$S(w;x,s)$ satisfy
$
x=w+s\chi(w)$.
Using $\chi'(w_0)=-1/s_0$, we expand
\begin{align}\label{e:arc_x_exp}
w+s\chi(w)
&=
w_0+s\chi(w_0)
+\left(1-\frac{s}{s_0}\right)(w-w_0)
+\sum_{k\geq 2}\frac{s a_k(w-w_0)^k}{k!} \notag \\
&=
x_0'
-\frac{s-s_0}{s_0}(w-w_0)
+\sum_{k\geq 2}\frac{s a_k(w-w_0)^k}{k!}.
\end{align}
Hence the critical point equation is equivalent to
\begin{align}
0
&=w+s\chi(w)-x =
-(x-x_0')
-\frac{s-s_0}{s_0}(w-w_0)
+\frac{s a_2}{2}(w-w_0)^2
+\sum_{k\geq 3}\frac{s a_k(w-w_0)^k}{k!}.
\end{align}
Multiplying by $2/(s a_2)$ gives
\begin{align}
0
=
-\frac{2(x-x_0')}{s a_2}
-\frac{2(s-s_0)}{a_2ss_0}(w-w_0)
+(w-w_0)^2
+\sum_{k\geq 3}\frac{2a_k(w-w_0)^k}{a_2k!}.
\end{align}
Thus the equation has the form
\begin{align}\label{e:arctic_F}
F(z;a,b):=z^2-a z-b+\sum_{k\geq 3}c_k z^k=0,
\qquad
z=w-w_0,
\end{align}
where
\[
a=\frac{2(s-s_0)}{a_2ss_0},
\qquad
b=\frac{2(x-x_0')}{s a_2},
\qquad
c_k=\frac{2a_k}{a_2k!}.
\]

For $\|(x,s)-(x_0,s_0)\|_2\leq \delta$, we have
\begin{align}
|a|\lesssim |s-s_0|\leq \delta, \quad
|b|\lesssim |x-x_0'|
\leq |x-x_0|+|x_0-x_0'|
= |x-x_0|+|(s-s_0)\chi(w_0)|
\lesssim \delta,
\end{align}
where we used \eqref{e:x_0'} and $x_0=w_0+s_0\chi(w_0)$.

Therefore, for $\delta$ sufficiently small, \Cref{lem:edge-perturbation}
implies that $F(z;a,b)$ has exactly two zeros in the disk $|z|\leq \fc$,
counted with multiplicity. Translating back to $w$, the corresponding critical
points satisfy
\begin{align}\label{e:edge_S_diff_critical}
|w_c-w_0|
\lesssim |a|+\sqrt{|b|}
\lesssim \sqrt{\delta}.
\end{align}

We recall the expansion of $S(w;x_0,s_0)$ and $S'''(w;x_0,s_0)$ from \eqref{e:arctic_neighborhood} and \eqref{e:derSthird2}. 
Then \eqref{e:chip_arctic} implies
\begin{align}
|d|=\frac16| S'''(w_0;x_0,s_0)|\asymp 1. 
\end{align}
The estimate \eqref{e:edge_S_diff_critical}, together with \eqref{e:Scritical1}, implies that uniformly for $|w-w_0|=r$,
\begin{align}
S(w;x,s)-S(w_c;x,s)
=
d(w-w_0)^3
+\OO\bigl(\delta^{3/2}+r^{4}+\delta\ln(1/\delta)\bigr),
\end{align}
Taking $r=\fc$ and then choosing $\fc$ and $\delta$ sufficiently small gives
\eqref{e:arctic_err}.

\end{proof}

\begin{proof}[Proof of \Cref{c:arctic_steepest}]

By \Cref{c:arctic_critical}, uniformly for $|w-w_0|=\fc$,
\begin{align}
|S(w;x,s)-S(w_c;x,s)-
d(w-w_0)^3|
\leq |d|\fc^3/100,
\quad
d:=S'''(w_0;x_0,s_0)/6.
\end{align}
with $|d|\asymp 1$. The hypotheses of
\Cref{prop:exit} hold with $(m, \delta, r)$ taken to be $(3,  |d|\fc^3/100, \fc)$.

\noindent\textbf{Replacing \(\mathsf C^{\rm d}(w_0)\) by \(\mathsf D^{\rm d}(w_c)\).}
We begin with \(d=S'''(w_0;x_0,s_0)/6<0\). There are three configurations.

\begin{enumerate}
\item \textbf{Arctic Boundary}
If $(x,s)\in \fA$, there is a single real critical point \(w_c\) with multiplicity two; see the second panel of
\Cref{f:cubic_saddle}. It is a cubic critical point of $S(\cdot;x,s)$. Hence, from \(w_c\) emanate three steepest--descent and three
steepest--ascent branches, alternating in angle. By
Proposition~\ref{prop:exit}, these six trajectories intersect the circle
\(|w-w_0|=\fc\) in the arcs \(A^{(6)}_0(w_0,\fc),A^{(6)}_1(w_0,\fc),\cdots,A^{(6)}_5(w_0,\fc)\) in alternating order. For $d<0$, one descent branch runs along the real axis from \(w_c\) to $w_0+\fc$, and one ascent branch runs along
the real axis from the opposite intersection $w_0-\fc$ to \(w_c\).

We define $\mathsf D^{\rm d}(w_c)$  to be the portions of
non-real steepest--descent trajectories, starting at \(w_c\) up to their first exit from
the disk $|w-w_0|\leq \fc$. Within the disk $|w-w_0|\leq \fc$, the two non-real steepest--descent trajectories from
\(w_c\) exit transversely through the boundary arcs \(A^{(6)}_2(w_0;\fc)\) and
\(A^{(6)}_4(w_0;\fc)\). By Cauchy's theorem, we may deform \(\mathsf C^{\rm d}(w_0)\) to
\(\mathsf D^{\rm d}(w_c)\), together with subarcs of \(A^{(6)}_2(w_0;\fc)\cup A^{(6)}_4(w_0;\fc)\).
By \eqref{e:final}, for any $w$ on $A^{(6)}_2(w_0;\fc)\cup A^{(6)}_4(w_0;\fc)$, 
\begin{align}\label{e:difS1}
n\,\Re\!\big[S(w;x,s)-S(w_c;x,s)\big]
\leq
-\,n\Bigl(|d|\,\fc^{3}-3|d|\fc^3/100\Bigr)
\le -\,n\fc'.
\end{align}

\item \textbf{Frozen Region}
If $(x,s)$ is in the frozen region, there are two simple saddles \(w_{c,1}<w_{c,2}\); see the third panel of
\Cref{f:cubic_saddle}. Each has two descent and two ascent branches, forming the
usual cross pattern of a simple saddle. Along the real axis, 
the gradient flow goes from $w_0-\fc$ to $w_{c,1}$, and from $w_{c,2}$ to $w_0+\fc$. $w_{c,1}$ is a descent critical point, and $w_{c,2}$ is an ascent critical point. The remaining branches are off the real axis. Since gradient flows
do not intersect, the two non-real descent branches from \(w_{c,1}\) exit the
circle through arcs \(A^{(6)}_2(w_0,\fc)\) and \(A^{(6)}_4(w_0,\fc)\), while the two complex ascent
branches from \(w_{c,2}\) exit through arcs \(A^{(6)}_1(w_0,\fc)\) and \(A^{(6)}_5(w_0,\fc)\), as shown in
the third panel of \Cref{f:cubic_saddle}.

For $w_c\in\{w_{c,1}, w_{c,2}\}$, we define $\mathsf D^{\rm d}(w_c)$ to be the portions of non-real
steepest--descent trajectories,  starting at \(w_c\) up to their first exit from the disk $|w-w_0|\leq \fc$. In this case $\mathsf D^{\rm d}(w_{c, 2})=\emptyset$.
By the same argument as in the first case, we may deform \(\mathsf C^{\rm d}(w_0)\) to
\(\mathsf D^{\rm d}(w_{c,1})\), and \eqref{e:difS1} holds.

\item \textbf{Liquid Region}
If $(x,s)\in \fL$, there are two complex conjugate critical points \(w_{c,+}\) and \(w_{c,-}=\overline{w_{c,+}}\);
see the fourth panel of \Cref{f:cubic_saddle}. By real-analyticity, the flow is
symmetric with respect to the real axis. It suffices to describe the upper
half-plane: from \(w_{c,+}\) there emanate two descent and two ascent branches,
all contained in the upper half-plane; they cannot cross the real axis, which itself is
a flow line from left to right. By
Proposition~\ref{prop:exit}, these four branches meet \(|w-w_0|=\fc\) through
the four alternating arcs \(A^{(6)}_0(w_0,\fc),A^{(6)}_1(w_0,\fc),A^{(6)}_2(w_0,\fc),A^{(6)}_3(w_0,\fc)\). The configuration from
\(w_{c,-}\) is the complex conjugate of that from \(w_{c,+}\), as shown in
the fourth panel of \Cref{f:cubic_saddle}.

For $w_c\in\{w_{c,+}, w_{c,-}\}$,  we define $\mathsf D^{\rm d}(w_c)$, to be the portions of
steepest--descent trajectories starting at \(w_c\) up to their first exit from the disk $|w-w_0|\leq \fc$.
By the same argument as in the first case, we may deform \(\mathsf C^{\rm d}(w_0)\) to
the union of the $\mathsf D(w_c)$'s, and \eqref{e:difS1} holds.

\end{enumerate}

When \(d=S'''(w_0;x_0,s_0)/6>0\), the same three configurations occur, but all gradient directions
are reversed. The claims can be proven in the same way, so we omit.

\noindent\textbf{Replacing \(\mathsf D^{\rm d}(w_c)\) by \(\mathsf S^{\rm d}(w_c)\).}
We recall $F(\cdot;a,b)$ from \eqref{e:arctic_F}. From \Cref{lem:edge-perturbation}, there exists $b'$ such that $F(\cdot;a,b')$ has a double root at $z'$. Denote $w'=z'+w_0$ and $x'=w'+s\chi(w')$. Then $b'=2(x'-x_0')/(sa_2)$, and $w'$ is a double root of $w+s\chi(w)-x'=0$. In particular, $\chi'(w')=-1/s$. By the second statement of \Cref{p:surface}, $(x',s)\in \fA$ lies on the arctic boundary, and 
\begin{align}
|x-x'|\asymp \dist((x,s), \fA).
\end{align}

We introduce the shifted coordinates
\begin{align}\label{e:shift_coordinate}
  \sfa = s - s_0, \quad
  \sfb = x - x_0', \quad
  \sfb' = x' - x_0'.
\end{align}
Then $|\sfa|\asymp |a|$, $|\sfb|\asymp |b|$, and $|\sfb-\sfb'|\asymp |b-b'|$. Moreover, \Cref{lem:edge-perturbation} implies
\begin{align}\label{e:w_c-w'}
|w'-w_0|\asymp |\sfa|, \quad
|\sfb'|\asymp |\sfa|^2,\quad
|w_c-w'|\asymp \sqrt{|\sfb'-\sfb|}.
\end{align}
By our assumption $\|(x,s)-(x_0,s_0)\|_2\leq  \delta$, so $|\sfa|,|\sfb|,|\sfb'|\lesssim \delta$. Together with the above estimate, we also have $|w'-w_0|, |w_c-w_0|\lesssim \sqrt{\delta}$.
We recall the estimates of $\chi$ from \eqref{e:chi_property}. For $|w-w_0|\lesssim \sqrt{\delta}$, we have
\begin{align}\begin{split}\label{e:chi_est}
&\chi(w)(1-\chi(w))=(1+\OO(\sqrt{\delta}))\chi(w_0)(1-\chi(w_0)),
\quad
\chi'(w)=(1+\OO(\sqrt{\delta}))\chi'(w_0),\\
&\chi''(w)=(1+\OO(\sqrt{\delta}))\chi''(w_0),
\quad
|\chi(w)(1-\chi(w))|,\ |\chi''(w)|\asymp 1. 
\end{split}\end{align}

We discuss three cases separately
\begin{enumerate}
\item
\textbf{Close to Arctic case $|\sfb-\sfb'|\leq (\ln n)^2/n^{2/3}$.}
In this case, $(x,s)$ is close to the arctic boundary point $(x',s)\in \fA$. Locally around $w'$, $S(w;x',s)$ is holomorphic and we have the Taylor expansion
\begin{align}\label{e:arc_S_exp1}
S(w;x',s)=S(w';x',s)+d(w-w')^3+\OO(|w-w'|^4),
\end{align}
where $d=S'''(w';x',s)/6\in \bR$, and by \eqref{e:chi_est} we have
\begin{align}\label{e:arc_S_exp2}
S'''(w';x',s)
=
-\frac{\chi''(w')}{\chi(w')(1-\chi(w'))},
\qquad
|d|\asymp |S'''(w';x',s)|
\asymp
\frac{|\chi''(w_0)|}{|\chi(w_0)(1-\chi(w_0))|}
\asymp 1.
\end{align}

Let
$
r_n:=(\ln n)^2/n^{1/3}$.
We define the local steepest--descent set $\mathsf S^{\rm d}(w_c)$ at \(w_c\) to be the portions of
steepest--descent trajectories starting at \(w_c\) up to their first exit from
the disk $\{w: |w-w'|\leq r_n\}$.  
For any $w=w'+r_n e^{\ri\theta}$, we have
\begin{align}\begin{split}\label{e:S_arc_wdiff}
&\phantom{{}={}}S(w;x,s)-S(w_c;x,s)\\
&=S(w;x',s)-S(w_c;x',s)
+\int_{w_c}^w
\left[
\ln\left(\frac{u-x}{x-s-u}\right)
-\ln\left(\frac{u-x'}{x'-s-u}\right)
\right]\rd u\\
&=d(w-w')^3
+\OO\bigl(|w_c-w'|^3+|w-w'|^4+|w_c-w'|^4+|w-w_c||x-x'|\bigr)\\
&=d(w-w')^3
+\OO\bigl(|\sfb-\sfb'|^{3/2}+r_n^4+r_n|\sfb-\sfb'|\bigr)\\
&=d(w-w')^3+\OO\bigl((\ln n)^4/n\bigr),
\end{split}\end{align}
where the first equality follows from \eqref{e:S_copy}; the second follows from \eqref{e:arc_S_exp1}, and we bound the integrand by $\OO(|x-x'|)$; the third follows from \eqref{e:shift_coordinate} and \eqref{e:w_c-w'}; the fourth follows from the bound $|\sfb-\sfb'|\leq (\ln n)^2/n^{2/3}$ and the choice $r_n=(\ln n)^2/n^{1/3}$. 
The estimate \eqref{e:S_arc_wdiff} verifies the assumptions in \Cref{prop:exit} with $(m,\delta,r)$ taken to be $(3,(\ln n)^4/n,r_n)$. 
It follows that the two steepest--descent branches from \(w_c\) belonging to $\mathsf D^{\rm d}(w_c)$ exit the circle
\(\{w: |w-w'|=r_n\}\) at points on \(A^{(4)}_2(w';r_n)\) and \(A^{(4)}_4(w';r_n)\). At any such
point,
\begin{equation}\label{eq:sd-drop0}
n\,\Re\!\big(S(w;x,s)-S(w_c;x,s)\big)
=
-\,n|d|\,r_n^{3}+\OO((\ln n)^4)
\le -\,\fc'\,(\ln n)^{6}.
\end{equation}
Moreover, since $\Re[S(w;x,s)]$ decreases along the steepest--descent, \eqref{eq:sd-drop0}
holds for all \(w\in\mathsf D^{\rm d}(w_c)\setminus \mathsf S^{\rm d}(w_c)\).
Consequently, for  \(w\in \mathsf D^{\rm d}(w_c)\setminus \mathsf S^{\rm d}(w_c)\) we have
\begin{align}
e^{n\Re[S(w;x,s)]}
\leq
e^{n\Re[S(w_c;x,s)]}e^{-\fc'(\ln n)^6}.
\end{align}
\item 
\textbf{Liquid Case with $|\sfb-\sfb'|\geq (\ln n)^2/n^{2/3}$.}
In this case there are two complex conjugate critical points $w_c$ and $\overline{w}_c$. Let
\[
r_n:=\frac{\ln n}{n^{1/2}|\sfb-\sfb'|^{1/4}}
\leq \frac{(\ln n)^{1/2}}{n^{1/3}}
\ll \sqrt{|\sfb-\sfb'|}.
\]
We define the local steepest--descent set $\mathsf S^{\rm d}(w_c)$ at \(w_c\) to be the portions of
steepest--descent trajectories starting at \(w_c\) up to their first exit from
the disk $\{w: |w-w_c|\leq r_n\}$.
In the following, we show that we can replace $\mathsf D^{\rm d}(w_c)$ by $\mathsf S^{\rm d}(w_c)$. The statement for $\overline w_c$ can be proven in the same way, so we omit it.

Locally around $w_c$, $S(w;x,s)$ is holomorphic and we have the Taylor expansion
\begin{align}\begin{split}\label{e:edge_S_exp1}
S(w;x,s)
&=S(w_c;x,s)+d(w-w_c)^2+\OO(|w-w_c|^3),
\end{split}\end{align}
where
\begin{align}\begin{split}\label{e:edge_S_exp2}
S''(w_c;x,s)
=
-\frac{\chi'(w_c)+1/s}{\chi(w_c)(1-\chi(w_c))},
\qquad
d=S''(w_c;x,s)/2.
\end{split}\end{align}
Notice that $\chi'(w')=-1/s$. We can further estimate $d$ as
\begin{align}\begin{split}\label{e:dbound}
|d|
&\asymp |S''(w_c;x,s)|
=
\frac{|\chi'(w_c)+1/s|}{|\chi(w_c)(1-\chi(w_c))|}
\asymp |\chi'(w_c)+1/s|\\
&\asymp |\chi'(w_c)-\chi'(w')|
\asymp |\chi''(w')(w_c-w')|
\asymp |\chi''(w_0)(w_c-w')|
\asymp |w_c-w'|
\asymp \sqrt{|\sfb-\sfb'|},
\end{split}\end{align}
where we used \eqref{e:chi_est} and \eqref{e:w_c-w'}. Since $\sqrt{|\sfb-\sfb'|}\gg r_n$, \eqref{e:dbound} implies $|d|r_n^2\asymp \sqrt{|\sfb-\sfb'|}r_n^2 \gg r_n^3$. This verifies the assumptions of \Cref{prop:exit} with $(m,\delta,r)$ taken to be $(2,r_n^3,r_n)$. We conclude that
the two steepest--descent branches from \(w_c\) exit the circle
\(\{w: |w-w_c|=r_n\}\) at points on \(A^{(4)}_0(w_c;r_n)\) and \(A^{(4)}_2(w_c;r_n)\). At any such
point,
\begin{equation}\label{eq:sd-arc_drop1}
n\,\Re\!\big(S(w;x,s)-S(w_c;x,s)\big)
=
-\,n\Bigl(|d|\,r_n^{2}+\OO(r_n^{3})\Bigr)
\le -\,\fc'\,(\ln n)^{2}.
\end{equation}
Moreover, since $\Re[S(w;x,s)]$ decreases along the steepest--descent, \eqref{eq:sd-arc_drop1}
holds for all \(w\in\mathsf D^{\rm d}(w_c)\setminus \mathsf S^{\rm d}(w_c)\).
Consequently, for  \(w\in \mathsf D^{\rm d}(w_c)\setminus \mathsf S^{\rm d}(w_c)\) we have
\begin{align}
e^{n\Re[S(w;x,s)]}
\leq
e^{n\Re[S(w_c;x,s)]}e^{-\fc'(\ln n)^2}.
\end{align}

\item 
\textbf{Frozen Case $|\sfb-\sfb'|\geq (\ln n)^2/n^{2/3}$.}
In this case there are two real solutions $w_{c,1}<w'<w_{c,2}$ with
\begin{align}\label{e:wc1wc2}
|w_{c,1}-w'|,\ |w_{c,2}-w'|\asymp \sqrt{|\sfb-\sfb'|}.
\end{align}
Without loss of generality, we assume $S'''(w_0;x_0,s_0)<0$, recall the first row of \Cref{f:cubic_saddle}; the case $S'''(w_0;x_0,s_0)>0$ can be proven in the same way, so we omit it. In this case, we let $w_c=w_{c,1}$ is a descent critical point. Then the same expansion \eqref{e:edge_S_exp1} holds with $d=S''(w_c;x,s)/2$, and 
\begin{align}\label{e:Svalue}
S''(w_c;x,s)
&\asymp
-\frac{\chi''(w_0)(w_c-w')}{\chi(w_0)(1-\chi(w_0))}
\asymp
S'''(w_0;x_0,s_0)(w_c-w')
\asymp
-(w_c-w')
\asymp
\sqrt{|\sfb-\sfb'|},
\end{align}
where we used \eqref{e:chi_est} and \eqref{e:w_c-w'}. 
By the same argument as in the liquid case, we can replace $\mathsf D^{\rm d}(w_c)$ by $\mathsf S^{\rm d}(w_c)$.
\end{enumerate}

Finally, to prove \eqref{e:Sdiff_at_critical}, we first record the following elementary consequence of Taylor expansion. Suppose
that \(F'(c_1)=F'(c_2)=0\). Then
\begin{align}\label{e:Taylor_critical_difference}
F(c_2)-F(c_1)
=
\frac{F''(c_1)(c_2-c_1)^2}{6}
+\OO(|c_2-c_1|^4).
\end{align}
Indeed, by Taylor expansion,
\begin{align}
F(c_2)-F(c_1)
&=
\frac{F''(c_1)(c_2-c_1)^2}{2}
+\frac{F'''(c_1)(c_2-c_1)^3}{6}
+\OO(|c_2-c_1|^4),
\\
0=F'(c_2)-F'(c_1)
&=
F''(c_1)(c_2-c_1)
+\frac{F'''(c_1)(c_2-c_1)^2}{2}
+\OO(|c_2-c_1|^3).
\end{align}
Subtracting
\[
\frac{c_2-c_1}{3}\bigl(F'(c_2)-F'(c_1)\bigr)
\]
from \(F(c_2)-F(c_1)\) gives \eqref{e:Taylor_critical_difference}.

We now apply \eqref{e:Taylor_critical_difference} this with
\[
F(\cdot)=S(\cdot;x,s), \qquad c_1=w_c, \qquad c_2=w_c'.
\]
This gives
\begin{align}\label{e:Svalueat}
S(w_c';x,s)-S(w_c;x,s)
=
\frac{1}{6}S''(w_c;x,s)(w_c'-w_c)^2
+\OO(|w_c'-w_c|^4).
\end{align}

Substituting \eqref{e:wc1wc2} and \eqref{e:Svalue} into
\eqref{e:Svalueat}, we obtain
\begin{align}
\Re\!\bigl[S(w_c';x,s)-S(w_c;x,s)\bigr]
=
\frac{1}{6}S''(w_c;x,s)(w_c'-w_c)^2
+\OO(|w_c'-w_c|^4)
\asymp
|\sfb-\sfb'|^{3/2}.
\end{align}
This proves \eqref{e:Sdiff_at_critical}, by recalling from \eqref{e:shift_coordinate} that $\sfb-\sfb'\asymp x-x'(s)$.

\end{proof}

\subsubsection{Orientation of $\sfD^{\rm d}(w_c)$}
\label{s:orientation}
For $(x,s)\in \fL$ in an arctic chart $\fU$ centered at $w_0$, the
orientation of the corresponding steepest--descent contour
$\sfD^{\rm d}(w_c)$, namely whether it travels from
$A^{(6)}_0(w_0;\fc)$ to $A^{(6)}_2(w_0;\fc)$ or in the reverse direction,
can be read off from its local orientation at the critical point $w_c$. 

By \eqref{e:arc_direction}, as $(x,s)\in\fL$ approaches a regular arctic
point $(x_0,s_0)$, the direction of $S''(w_c;x,s)$ is close to
$\ri S'''(w_0;x_0,s_0)$, which is in the direction of $\pm\ri$. Without
loss of generality, assume
\[
    S'''(w_0;x_0,s_0)<0 .
\]
The case $S'''(w_0;x_0,s_0)>0$ follows by the same argument. Under this
assumption, $S''(w_c;x,s)$ is close to the negative imaginary direction.
Thus the two steepest--descent directions from $w_c$ are close to
\[
    e^{3\ri\pi/4},\qquad -e^{3\ri\pi/4},
\]
while the two steepest--ascent directions are close to
\[
    e^{\ri\pi/4},\qquad -e^{\ri\pi/4}.
\]

The arcs $A^{(6)}_j(w_0;\fc)$, $0\leq j\leq 5$, are small pieces of
$|w-w_0|=\fc$ centered in the directions $j\pi/3$. By the exit statement
above, see Panel (C) of
\Cref{f:cubic_saddle}, the four steepest--ascent/descent trajectories from $w_c$ exit
\[
    \{w: |w-w_0|\leq \fc\}
\]
transversely through the boundary arcs
\[
    A^{(6)}_0(w_0;\fc),\,
    A^{(6)}_1(w_0;\fc),\,
    A^{(6)}_2(w_0;\fc),\,
    A^{(6)}_3(w_0;\fc)
\]
in cyclic order, since distinct steepest trajectories cannot intersect
away from a critical point.

We now show that the steepest--descent trajectory with initial direction
$e^{3\ri\pi/4}$ exits through $A^{(6)}_2(w_0;\fc)$, while the steepest--
descent trajectory with initial direction $-e^{3\ri\pi/4}$ exits through
$A^{(6)}_0(w_0;\fc)$. This is the same local picture as in Panel (C) of
\Cref{f:cubic_saddle}.

In the arctic chart $\fU$, we have
\begin{align}
\begin{split}\label{e:Sprime_expanded_local_simplified}
S'(w;x,s)
&=S'(w;x_0,s_0)+\OO(\delta)  =\frac{S'''(w_0;x_0,s_0)}{2}(w-w_0)^2
  +\OO\!\bigl(\delta+|w-w_0|^3\bigr),
\\
S''(w;x,s)
&=S''(w;x_0,s_0)+\OO(\delta)  =S'''(w_0;x_0,s_0)(w-w_0)
  +\OO\!\bigl(\delta+|w-w_0|^2\bigr).
\end{split}
\end{align}
By Rouch\'e's theorem, for $\delta$ sufficiently small the function
$S''(\,\cdot\,;x,s)$ has exactly one zero in $\fU$; denote it by
$w_\ast$. Since the tiling action function $S$ satisfies Schwarz reflection,
\[
    S(\overline w;x,s)=\overline{S(w;x,s)},
\]
the zero set of $S''(\,\cdot\,;x,s)$ is invariant under complex
conjugation. Since the zero in $\fU$ is unique, it follows that
\[
    w_\ast\in \bR .
\]

In the arctic chart $\fU$, the set
\[
    \{\,w\in\fU:\ \Im S'(w;x,s)=0\,\}
\]
is a finite union of real-analytic arcs. It contains the real axis in
$\fU$. There is also a non-real arc $\gamma$ which passes through
$w_c$ and $\overline{w_c}$, and this arc intersects the real axis at the
turning point $w_\ast$. By
\eqref{e:Sprime_expanded_local_simplified}, since
$S'''(w_0;x_0,s_0)<0$, the arc $\gamma$ is nearly vertical. It exits
$\fU$ through the upper boundary close to $w_0+\fc\ri$ and through the
lower boundary close to $w_0-\fc\ri$.

Expanding $S'$ about the critical point $w_c$ gives
\begin{equation}\label{e:Sprime_at_wc_simplified}
    S'(w;x,s)
    =
    S''(w_c;x,s)(w-w_c)+\OO(|w-w_c|^2).
\end{equation}
As $(x,s)$ approaches $(x_0,s_0)$, $S''(w_c;x,s)$ is close to the negative
imaginary direction; equivalently,
\[
    \arg S''(w_c;x,s)=-\frac{\pi}{2}+\oo(1).
\]
Thus \eqref{e:Sprime_at_wc_simplified} implies that the level set
$\Im S'(w;x,s)=0$ through $w_c$ is nearly vertical.

Consequently, the arc $\gamma$ has two local branches at $w_c$: one leaves
$w_c$ upward and the other leaves $w_c$ downward. The upward branch exits
$\fU$ through the upper boundary close to $w_0+\fc\ri$, while the downward
branch first connects to the turning point $w_\ast$ on the real axis.
Along the upward branch, $S'(w;x,s)>0$ as a real number, so a steepest--
descent trajectory crosses this branch from right to left. Along the
downward branch, $S'(w;x,s)<0$, so a steepest--descent trajectory crosses
this branch from left to right.

Suppose, for contradiction, that the steepest--descent trajectory with
initial direction $e^{3\ri\pi/4}$ exits through
$A^{(6)}_0(w_0;\fc)$ instead of $A^{(6)}_2(w_0;\fc)$. Starting from
$w_c$, this trajectory must cross the arc $\gamma$ before reaching
$A^{(6)}_0(w_0;\fc)$. If its first crossing is on the upward branch of
$\gamma$, then it crosses that branch from left to right, contradicting
the crossing direction described above.

It remains to consider the case where the first crossing is on the
downward branch of $\gamma$. Then the portion of the descent trajectory
from $w_c$ to this first crossing, together with the corresponding portion
of the downward branch of $\gamma$, traps the steepest--ascent trajectory
with initial direction $-e^{\ri\pi/4}$. To exit through one of the
remaining boundary arcs, this ascent trajectory would have to cross either
the descent trajectory or the downward branch of $\gamma$ in the forbidden
direction. This is impossible: distinct steepest trajectories do not
intersect away from critical points, and ascent trajectories cross the
downward branch of $\gamma$ in the direction opposite to steepest descent.
This contradiction proves the claim.

Therefore the steepest--descent trajectory with initial direction
$e^{3\ri\pi/4}$ exits through $A^{(6)}_2(w_0;\fc)$, and the one with
initial direction $-e^{3\ri\pi/4}$ exits through $A^{(6)}_0(w_0;\fc)$.
Hence, with the local orientation whose tangent at $w_c$ is
$e^{3\ri\pi/4}$, the contour $\sfD^{\rm d}(w_c)$ runs from
$A^{(6)}_0(w_0;\fc)$ to $A^{(6)}_2(w_0;\fc)$; reversing the local
orientation reverses the orientation of the contour.

\subsection{Cusp chart}\label{s:cusp_chart_proof}

We recall from \Cref{s:critical_cusp} that, given a regular cusp chart
centered at $(x_0,s_0)$ and setting $w_0 = x_0 - s_0 \chi(w_0) \in \bR$, the 
point $w_0$ is a quartic saddle point of $S(w;x_0,s_0)$, with $m=4$ in 
\Cref{prop:exit}. For any $(x,s)$ in this regular cusp neighborhood, we study 
$S(w;x,s)$ in a small neighborhood of $w_0$.

\begin{proof}[Proof of \Cref{c:cusp_critical}]
By the second statement of \Cref{p:surface}, after decreasing $\fc$ if necessary,
the Riemann surface $\cC$ can be parametrized over $\fU$ as $(f(w),w)$, and
both $f$ and $\chi$ are holomorphic on $\fU$. Moreover, since $(x_0,s_0)$ is
bounded away from tangency points, we have
\begin{align}\label{e:chip_cusp}
|\chi(w_0)|,\ |1-\chi(w_0)|,\ |\chi'''(w_0)| \asymp 1, \quad .
\end{align}
Also, since $w_0$ is a cubic critical point of $S(\,\cdot\,;x_0,s_0)$, we have
\[
x_0=w_0+s_0\chi(w_0),
\quad
\chi'(w_0)=-\frac1{s_0},\quad \chi''(w_0)=0.
\]
Thus, after possibly decreasing $\fc$ again, the bounds in \eqref{e:cusp_bounds} hold uniformly for all $w\in \fU$.
And \eqref{e:cusp_quartic} follows from \eqref{e:derSfourth}.

Critical points $w_c$ of $S(w;x,s)$ satisfy the following equation 
\begin{align}\label{e:x_exp2}
w_0+s\chi(w_0)+(x-x_0')=x=w+s\chi(w)=w+ s\chi(w_0)-\frac{s(w-w_0)}{s_0}+ \sum_{k\geq 3}\frac{sa_k(w-w_0)^k}{k!}
\end{align}
We can rewrite it as
\begin{align}
-\frac{6(x-x_0')}{s a_3}-\frac{6(s-s_0)}{a_3 s_0}(w-w_0)+(w-w_0)^3+\sum_{k\geq 4}\frac{6a_k(w-w_0)^k}{a_3k!}=0,
\end{align}
which is in the following form
\begin{align}\label{e:cusp_F}
F(z;a,b):=z^3-a z-b+\sum_{k\ge 4} c_k z^k=0,\quad z=w-w_0, \quad a=\frac{6(s-s_0)}{a_3 s_0}, \quad b=\frac{6(x-x_0')}{s a_3}.
\end{align}

For $\|(x,s)-(x_0, s_0)\|_2\leq \delta$, by \Cref{lem:cubic-perturbation}, we have 
\begin{align}\label{e:cusp_critical}
|w_c-w_0|\lesssim |s-s_0|^{1/2}+|x-x'|^{1/3}\lesssim |s-s_0|^{1/2}+(|x-x_0|+|(s-s_0)\chi(w_0)|)^{1/3}\lesssim \delta^{1/3}.
\end{align}

We recall the expansion of $S(w;x_0,s_0)$ and $S''''(w;x_0,s_0)$ from \eqref{e:cusp_neighborhood} and \eqref{e:derSfourth2}. 
Then \eqref{e:chip_cusp} implies
\begin{align}
|d|=\frac1{24}| S''''(w_0;x_0,s_0)|\asymp 1. 
\end{align}
The estimate \eqref{e:cusp_critical}, together with \eqref{e:Scritical1}, implies that uniformly for $|w-w_0|=r$,
\begin{align}
S(w;x,s)-S(w_c;x,s)
=
d(w-w_0)^3
+\OO\bigl(\delta^{3/2}+r^{4}+\delta\ln(1/\delta)\bigr),
\end{align}
Taking $r=\fc$ and then choosing $\fc$ and $\delta$ sufficiently small gives
\eqref{e:arctic_err}.

\end{proof}

\begin{proof}[Proof of \Cref{l:cusp_steepest}]

By \Cref{c:cusp_critical}, uniformly for $|w-w_0|=\fc$,
\begin{align}
|S(w;x,s)-S(w_c;x,s)-
d(w-w_0)^4|
\leq |d|\fc^4/100,
\quad
d:=S''''(w_0;x_0,s_0)/24.
\end{align}
with $|d|\asymp 1$. The hypotheses of
\Cref{prop:exit} hold with $(m, \delta, r)$ taken to be $(4,  |d|\fc^4/100, \fc)$.

We begin with \(d=S''''(w_0;x_0, s_0)/24<0\).
There are four configurations:

\begin{enumerate}
\item \textbf{Cusp Location.}
If $(x,s)=(x_0,s_0)$, there is a single critical point \(w_c=w_0\) with multiplicity three; see panel (A) of
\Cref{f:quartic_saddle}. At $w_c$, $S(\cdot; x,s)$  vanishes up to the quartic term.
Hence, from \(w_0\) emanate four steepest–descent and four
steepest–ascent branches, alternating in angle. By Proposition~\ref{prop:exit}, these
eight trajectories intersect the circle \(|w-w_0|=\fc\) in the arcs
\(A^{(8)}_0(w_0, \fc),A^{(8)}_1(w_0, \fc),\cdots ,A^{(8)}_7(w_0, \fc)\) in alternating order. In particular, when $d=S''''(w_0; x_0, s_0)/24<0$, the two
descent branches run along the real axis from \(w_0\) to
\(w_0\pm\fc\).

We define $\mathsf D^{\rm d}(w_0)$ the local steepest–descent set at \(w_0\) to be the portions of
non-real steepest–descent trajectories, starting at \(w_0\) up to their first exit from
the disk $|w-w_0|\leq \fc$.  Within the disk $|w-w_0|\leq \fc$, the two non-real steepest–descent trajectories from
\(w_c\) exit transversely through the boundary arcs \(A^{(8)}_2(w_0;\fc)\) and
\(A^{(8)}_6(w_0;\fc)\). By Cauchy’s theorem, we may deform \(\mathsf C^{\rm d}(w_0)\) to
\(\mathsf D^{\rm d}(w_0)\) together with (sub)arcs of \(A^{(8)}_2(w_0;\fc)\cup A^{(8)}_6(w_0;\fc)\).
By \eqref{e:final}, for any $w$ on $A^{(8)}_2(w_0;\fc)\cup A^{(8)}_6(w_0;\fc)$, 
\begin{align}\label{e:cusp_difS1}
n\,\Re\!\big[S(w;x,s)-S(w_c;x,s)\big]
\leq 
-\,n\Bigl(|d|\,\fc^{4}- |d|\fc^4/100\Bigr)
\le -\,n\fc'.
\end{align}

\item \textbf{Arctic Boundary.}
If $(x,s)\in \fA$ but $(x,s)\neq (x_0,s_0)$, there are two distinct real critical points \(w_{c,1}<w_{c,2}=w_{c,3}\),
where (without loss of generality) \(w_{c,2}\) has multiplicity two; see panel (B) of \Cref{f:quartic_saddle}. From \(w_{c,1}\) emanate two
descent and two ascent branches; from \(w_{c,2}\) emanate three descent and
three ascent branches, alternating in angle. Along the real axis, the gradient
flow goes from \(w_{c,1}\) to \(w_0-\fc\) and to \(w_{c,2}\), and from
\(w_{c,2}\) to \(w_0+\fc\). The remaining branches are off the real axis:
the two non-real ascent branches from \(w_{c,1}\) meet the circle \(|w-w_0|=\fc\)  through
arcs \(A^{(8)}_3(w_0, \fc)\) and \(A^{(8)}_5(w_0,\fc)\); the four non-real descent/ascent branches from
\(w_{c,2}\) meet the circle alternately through \(A^{(8)}_1(w_0, \fc),A^{(8)}_2(w_0, \fc),A^{(8)}_6(w_0, \fc),A^{(8)}_7(w_0, \fc)\), as shown
in panel (B) of \Cref{f:quartic_saddle}. (This follows from that flows do not
intersect.)

For $w_c\in\{w_{c,1},w_{c,2}\}$, we define $\mathsf D^{\rm d}(w_c)$ to be the portions of non-real
steepest–descent trajectories starting at \(w_{c}\) up to their first exit from the disk $|w-w_0|\leq \fc$.
In this case, $w_{c,1}$ is an ascent critical point, and $\mathsf D^{\rm d}(w_{c,1})=\emptyset$.
By the same argument as in the first case, we may deform \(\mathsf C^{\rm d}(w_0)\) to
\(\mathsf D^{\rm d}(w_{c,2})\), and \eqref{e:cusp_difS1} holds.

\item \textbf{Frozen Region.}
If $(x,s)$ is in the frozen region, there are three simple saddles \(w_{c,1}<w_{c,2}<w_{c,3}\); see panel (C)
of \Cref{f:quartic_saddle}. Each has two descent and two ascent branches, forming
the usual alternating pattern. Along the real axis, the flow goes from \(w_{c,1}\) to
\(w_0-\fc\) and to \(w_{c,2}\), and from \(w_{c,3}\) to \(w_{c,2}\) and to
\(w_0+\fc\). Off the real axis, the two non-real ascent branches from
\(w_{c,1}\) meet the circle through \(A^{(8)}_3(w_0, \fc)\) and \(A^{(8)}_5(w_0,\fc)\); the two non-real
ascent branches from \(w_{c,3}\) meet through \(A^{(8)}_1(w_0, \fc)\) and \(A^{(8)}_7(w_0,\fc)\); and the two
non-real descent branches from \(w_{c,2}\) meet through \(A^{(8)}_2(w_0,\fc)\) and \(A^{(8)}_6(w_0,\fc)\),
as depicted in panel (C) of \Cref{f:quartic_saddle}. (Again, this follows from that flows do not
intersect.)

For $w_c\in\{w_{c,1}, w_{c,2}, w_{c,3}\}$, we define $\mathsf D^{\rm d}(w_c)$ to be the portions of non-real
steepest–descent trajectories starting at \(w_{c}\) up to their first exit from the disk $|w-w_0|\leq \fc$.
In this case, $w_{c,2}$ is the only descent critical point, and $\mathsf D^{\rm d}(w_{c,2})\neq \emptyset$.
By the same argument as in the first case, we may deform \(\mathsf C^{\rm d}(w_0)\) to
\(\mathsf D^{\rm d}(w_{c,2})\), and \eqref{e:cusp_difS1} holds.

\item \textbf{Liquid Region.}
If $(x,s)\in \fL$, there is one real saddle \(w_{c,1}\) and a conjugate pair
\(w_{c,\pm}\), with \(\overline{w_{c,+}}=w_{c,-}\); see panel (D) of
\Cref{f:quartic_saddle}. Along the real axis, the two real branches of the flow
run from \(w_{c,1}\) to \(w_0-\fc\) and to \(w_0+\fc\). By real-analyticity,
the configuration is symmetric with respect to the real axis, so it suffices to
describe the upper half-plane. From \(w_{c,+}\) emanate two descent and two
ascent branches, forming the usual alternating pattern; these four branches
remain in the upper half-plane, since flow lines cannot cross the real axis. By
Proposition~\ref{prop:exit}, two ascent branches and one descent branch meet
\(|w-w_0|=\fc\) through
$
A^{(8)}_1(w_0,\fc)$,
$A^{(8)}_2(w_0,\fc)$,
and $A^{(8)}_3(w_0,\fc)$
in alternating order. The remaining descent branch either exits through
\(A^{(8)}_0(w_0,\fc)\) or \(A^{(8)}_4(w_0,\fc)\), or connects to the real
critical point \(w_{c,1}\); the dashed arrow in panel (D) of
\Cref{f:quartic_saddle} indicates one such possibility. The flows from
\(w_{c,-}\) are the complex conjugates of those from \(w_{c,+}\).

For \(w_c\in\{w_{c,-},w_{c,+}\}\), define \(\mathsf D^{\rm d}(w_c)\) to be the
portions of the steepest--descent trajectories starting at \(w_c\), stopped at
their first exit from the disk \(|w-w_0|\leq \fc\), or when they connect to
\(w_{c,1}\). By the same argument as in the first case, we may deform
\(\mathsf C^{\rm d}(w_0)\) to the union of these
\(\mathsf D^{\rm d}(w_c)\), up to boundary arcs on which
\eqref{e:cusp_difS1} holds.

\end{enumerate}
%Moreover, in all cases $w\in A_j$, 
%\begin{align}\label{e:cubic_final_copy}
%\Re\!\big(S(w;x,s)-S(w_c;x,s)\big)
%=(-1)^j\,|d|\,r^{\,m}+\OO\!\big(\delta+r^{\,m+1}\big),
%\end{align}
When \(d=S''''(w_0;x_0, s_0)/24>0\), the same four configurations occur, but all gradient directions
are reversed (second row of \Cref{f:quartic_saddle}).

\end{proof}

\subsection{Regular frozen chart and cusp frozen chart}
\label{s:frozen_chart_proof}.
We recall from \Cref{s:frozen_neighborhood} that, for a regular frozen chart or
a cusp frozen chart centered at $(x_0,s_0)$, the point
$
w_0=x_0-s_0\chi(w_0)\in \bR$
is a nondegenerate saddle point of $S(w;x_0,s_0)$. Thus $m=2$ in
\Cref{prop:exit}. For any $(x,s)$ in this chart, we study
$S(w;x,s)$ in a small neighborhood of $w_0$.

\begin{proof}[Proof of \Cref{c:regular_frozen_critical} and \Cref{c:cusp_frozen_critical}]
As in the regular liquid region, the critical point equation can be written in
the form \eqref{e:bulk_F}. In the present case, $w_0\in \bR$. Moreover, the
coefficients in \eqref{e:bulk_F} are real, and $b\in \bR$. Hence by \Cref{lem:bulk-perturbation}, $F(z;b)$ has a unique real zero in the relevant
disk. Therefore $w_c\in \bR$, and
\begin{align}\label{e:frozen_wc_est}
|w_c-w_0|
\lesssim |b|
\lesssim \|(x,s)-(x_0,s_0)\|_2
\leq \delta .
\end{align}
The rest of the argument is the same as in \Cref{c:bulk} for the regular liquid
region, so we omit it.
\end{proof}

\begin{proof}[Proof of \Cref{c:frozen_steepest}]
We repeat the argument from the proof of \Cref{c:bulk_steepest} for the regular
liquid region; see \Cref{f:bulk_critical}. In the present case,
\begin{align}\label{e:S''_exp}
S''(w_c;x,s)
=
-\frac{\chi'(w_c)+1/s}{\chi(w_c)(1-\chi(w_c))}
\in \bR,
\qquad
S''(w_0;x_0,s_0)
=
-\frac{\chi'(w_0)+1/s_0}{\chi(w_0)(1-\chi(w_0))}
\in \bR\setminus\{0\}.
\end{align}
For $\delta$ sufficiently small, \eqref{e:frozen_wc_est} and
\eqref{e:S''_exp} imply that $S''(w_c;x,s)$ has the same sign as
$S''(w_0;x_0,s_0)$. Hence there are two cases.

\begin{enumerate}
  \item If $S''(w_0;x_0,s_0)>0$, then $S''(w_c;x,s)>0$, the vertical
  direction through $w_c$ is a descent direction, and $w_c$ is a descent critical point. Moreover,
  $\mathsf C^{\rm a}(w_0)=\emptyset$, and we can deform
$
    \mathsf C^{\rm d}(w_0)
$
  to $\mathsf D^{\rm d}(w_c)$, together with finitely many arcs of total length
$\OO(1)$, on which \eqref{e:bulk_subarc} holds.

  \item If $S''(w_0;x_0,s_0)<0$, then $S''(w_c;x,s)<0$, the vertical
  direction through $w_c$ is an ascent direction, and $w_c$ is an ascent critical point. Moreover,
  $\mathsf C^{\rm d}(w_0)=\emptyset$, and we can deform
$
    \mathsf C^{\rm a}(w_0)
   $
  to $\mathsf D^{\rm a}(w_c)$, together with finitely many arcs of total length
$\OO(1)$, on which \eqref{e:bulk_subarc} holds.

\end{enumerate}
The rest of the argument is the same as in \Cref{c:bulk_steepest} for the regular liquid
region, so we omit it.
\end{proof}

\subsection{Vertical tangent chart and cusp-turning chart}
\label{s:vertical_tangent_chart_proof}.

%In this section, we fix a tangent neighborhood 
%centered at $(x_0,s_0)\in \fA$, which is not a cusp location. For any $(x,s)$ in this tangent neighborhood, we study 
%$S(w;x,s)$ in a small neighborhood of $w_0$. There are two cases, either $(x_0, s_0)$ is a vertical/unit-slope tangent location, or it is a horizontal tangent location. 

\begin{proof}[Proof of \Cref{c:tangent_critical1}]

We recall the expansion of $S(w;x_0, s_0)$ from \eqref{e:vert_tangent}. The main difference from the regular liquid region is that $S(w;x,s)$ is not analytic around $w_0$. In fact,  from \eqref{e:S_copy}, we have
\begin{align}\begin{split}\label{e:tangent_sdiff}
S(w; x,s) &= S(w; x_0,s_0)+(s\ln s-(x-w)\ln (x-w)-(s-x+w)\ln (s-x+w))\\
&-(s_0\ln s_0-(x_0-w)\ln (x_0-w)-(s_0-x_0+w)\ln (s_0-x_0+w)).
\end{split}\end{align}
Notice that $w_0=x_0$, from the above expression, $S(w;x,s)$ has a cut $[\min\{x_0, x\}, \max\{x_0, x\}]$. For $\dist(w,[\min\{x_0, x\}, \max\{x_0, x\}]) \gtrsim |x-x_0|$, we have the following estimates for the derivatives of \eqref{e:tangent_sdiff}
%\begin{align}
%S'(w; x,s) &= S'(w; x_0,s_0)+\ln \frac{x-w}{w-(x-s)}-\ln \frac{x_0-w}{w-(x_0-s_0)} 
%\end{align}
\begin{align}\begin{split}\label{e:tangent_S_der}
\del_w^k S(w;x,s)
&=\del_w^k S(w;x_0,s_0)+\del_w^{k-1} \ln\left(\frac{x-w}{x_0-w}\frac{w-(x_0-s_0)}{w-(x-s)}\right)\\
&=\del_w^k \ln\left(\frac{x-w}{x_0-w}\right)+\OO(1)\lesssim \frac{|x-x_0|}{|w-x_0|^k}+1,
\end{split}\end{align}

In this case, the same as for the regular arctic region, the critical point $w_c$ also satisfy the equation \eqref{e:arctic_F}:
\begin{align}
F(z;a,b):=z^2-a z-b+\sum_{k\geq 3}c_k z^k=0,
\qquad
z=w-w_0,
\end{align}
where
\[
a=\frac{2(s-s_0)}{a_2ss_0},
\qquad
b=\frac{2(x-x_0)}{s a_2},
\qquad
c_k=\frac{2a_k}{a_2k!},\quad a_k=\chi^{(k)}(w_0)
\]

By \Cref{lem:edge-perturbation}, and using
\(\delta\gtrsim |s-s_0|+|x-x_0|\), we have
\begin{align}\label{e:tangent}
\begin{split}
|w_c-w_0|
&\lesssim \max\{|a|,|b|^{1/2}\}
 \lesssim |s-s_0|+\sqrt{|x-x_0|}
 \lesssim \sqrt{\delta}, \\
|w_c-w_0|
&\gtrsim
\min\left\{\frac{|b|}{|a|},\max\{|a|,|b|^{1/2}\}\right\}  
\gtrsim
\min\left\{
\frac{|x-x_0|}{|s-s_0|}, |x-x_0|^{1/2}
\right\}
\gtrsim
\frac{|x-x_0|}{\sqrt{\delta}}.
\end{split}
\end{align}
Consequently,
\[
\dist\bigl(w_c,[\min\{x_0,x\},\max\{x_0,x\}]\bigr)
\gtrsim |x-x_0|,
\]
and \eqref{e:tangent_S_der} holds for $w=w_c$.

We now recall the expansion of \(S(w;x_0,s_0)\) from
\eqref{e:vert_tangent}.
Then \eqref{e:c_1} implies
\begin{align}
d:=\frac{1}{2}S''(w_0;x_0,s_0),\quad 
|d|=\frac{1}{2}|c_1 s_0|\asymp 1. 
\end{align}

Combining \eqref{e:tangent} with
\eqref{e:Scritical1}, we obtain, for \(|w-w_0|=r\),
\begin{align}
S(w;x,s)-S(w_c;x,s)
=
d(w-w_0)^2
+
\OO\bigl(\delta+r^3+\delta\ln(1/\delta)\bigr).
\end{align}
Taking $r=\fc$ and then choosing $\fc$ and $\delta$ sufficiently small gives
\eqref{e:tangent_err}.

\end{proof}

\begin{proof}[Proof of \Cref{c:cusp_turning_critical1}]

%The tangent through $(x_0,s_0)$ is vertical so $w_0=x_0$ and $\chi(w_0)=\chi(x_0, s_0)=0$.  
%We recall the expansion of $S(w;x_0, s_0)$ from \eqref{e:vert_tangent}. The main difference from the regular liquid region is that $S(w;x,s)$ is not analytic around $w_0$. In fact,  from \eqref{e:S_copy}, we have
%\begin{align}\begin{split}\label{e:tangent_sdiff}
%S(w; x,s) &= S(w; x_0,s_0)+(s\ln s-(x-w)\ln (x-w)-(s-x+w)\ln (s-x+w))\\
%&-(s_0\ln s_0-(x_0-w)\ln (x_0-w)-(s_0-x_0+w)\ln (s_0-x_0+w)).
%\end{split}\end{align}
%Notice that $w_0=x_0$, from the above expression, locally around $w_0$, $S(w;x,s)$ has a cut $[\min\{x_0, x\}, \max\{x_0, x\}]$. For $\dist(w,[\min\{x_0, x\}, \max\{x_0, x\}]) \gtrsim |x-x_0|$, we have the following estimates for the derivatives of \eqref{e:tangent_sdiff} from \eqref{e:tangent_S_der}
%%\begin{align}
%%S'(w; x,s) &= S'(w; x_0,s_0)+\ln \frac{x-w}{w-(x-s)}-\ln \frac{x_0-w}{w-(x_0-s_0)} 
%%\end{align}
%\begin{align}\begin{split}
%\del_w^k S(w;x,s)=\del_w^k \ln\left(\frac{x-w}{x_0-w}\right)+\OO(1)\lesssim \frac{|x-x_0|}{|w-x_0|^k}+1,
%\end{split}\end{align}

In this case, the same as for the regular cusp region, the critical point $w_c$ also satisfy the equation \eqref{e:cusp_F}:
\begin{align}
F(z;a,b):=z^3-a z-b+\sum_{k\ge 4} c_k z^k=0,\quad z=w-w_0
\end{align}
where
\[
a=\frac{6(s-s_0)}{a_3 s_0}, \qquad b=\frac{6(x-x_0)}{s a_3},
\qquad
c_k=\frac{6a_k}{a_3k!},\quad a_k=\chi^{(k)}(w_0)
\]

By \Cref{lem:cubic-perturbation}, we have 
\begin{align}\begin{split}\label{e:vertical_cusp}
&|w_c-w_0|\lesssim |s-s_0|^{1/2}+|x-x_0|^{1/3}\lesssim \delta^{1/3},\\
&|w_c-w_0|\gtrsim \min\left\{\frac{|x-x_0|}{|s-s_0|}, |x-x_0|^{1/3}\right\}\gtrsim \delta^{-2/3}|x-x_0|.
\end{split}\end{align}
In particular $\dist(w_c,[\min\{x_0, x\}, \max\{x_0, x\}]) \gtrsim |x-x_0|$. 

We now recall the expansion of \(S(w;x_0,s_0)\) from
\eqref{e:vert_tangent}.
Then \eqref{e:c_2} implies
\begin{align}
d:=\frac{1}{6}S'''(w_0;x_0,s_0),\quad 
|d|=\frac{1}{3}|c_2 s_0|\asymp 1. 
\end{align}

Combining \eqref{e:vertical_cusp} with
\eqref{e:Scritical1}, we obtain, for \(|w-w_0|=r\),
\begin{align}\label{e:vert_cusp}
S(w; x,s) - S(w_c; x,s)=d(w-w_0)^3+\OO(\delta+r^{4}+\delta  \ln(1/\delta )),\quad d:=S'''(w_0;x_0, s_0)/6.
\end{align}
Taking $r=\fc$ and then choosing $\fc$ and $\delta$ sufficiently small gives
\eqref{e:cusp_turning_err}.

\end{proof}

\begin{proof}[Proof of \Cref{l:vertical_tangent_steepest}]

By \Cref{c:tangent_critical1}, uniformly for $|w-w_0|=\fc$,
\begin{align}
|S(w;x,s)-S(w_c;x,s)-
d(w-w_0)^2|
\leq |d|\fc^2/100,
\quad
d:=S''(w_0;x_0,s_0)/2.
\end{align}
with $|d|\asymp 1$. The hypotheses of
\Cref{prop:exit} hold with $(m, \delta, r)$ taken to be $(2, |d|\fc^2/100, \fc)$.

We recall from \eqref{e:tangent_S_der} when $k=1$, $S'(w;x,s)$ contains the following non-analytic term
\begin{align}
\ln \frac{(x-w)}{x_0-w}.
\end{align}
The tangent location $(x_0, s_0)$ is on an vertical extended side of $\fP$. By our convention  $nx_0\in \bZ'=\bZ+1/2$, and $nx\in \bZ$, so either $x>x_0$ or $x<x_0$. 
If $x>x_0$, the interval $[x_0, x]$ is a source term, the steepest–descent paths point outward; 
If $x<x_0$  the interval $[x, x_0]$ is a source term, the steepest–descent paths point inward.

%; and if $x=x_0$ it vanishes and $S(w;x,s)$ is analytic locally around $w_0$. 

There are several cases. We begin with \(d=S''(w_0;x_0,s_0)<0\). This case corresponds to \Cref{f:vertical_tangent3} or \Cref{f:vertical_tangent4}, and $\mathsf C^{\rm d}(w_0)$ is a contour surrounding $w_0$.

\begin{enumerate}

\item \textbf{Arctic Boundary}
If $(x,s)\in \fA$, $w_c\in \bR$ is the only critical point of $S(\cdot;x,s)$ with multiplicity two. 
In this case,  $x > x_0$, so that the interval $[x_0,x]$ is a source: the steepest–descent paths point outward; see panel~(A) and (B) of
\Cref{f:tangent1}. 

There are two cases $(x,s)\in \fT_A$ or  $(x,s)\in \fT_B$. We will only discuss the case that $(x,s)\in \fT_B$, then $s<s_0$, and $w_c$ is determined by the intersection of the tangent line through $(x,s)$ with the $s=0$ axis, so $w_c>x$; see panel~(A) of
\Cref{f:tangent1}. The case that $(x,s)\in \fT_B$, see panel~(B) of
\Cref{f:tangent1} can be proven in the same way, so we omit. 

 At $w_c$, the expansion of $S(\cdot; x,s)$ vanishes up to the cubic term. Hence, from $w_c$ there emanate three steepest–descent and three steepest–ascent branches, alternating in angle. One descent branch and one ascent branch lie along the real axis, from $w_c$ to $w_0+\fc$ and from $x$ to $w_c$, respectively.

By Proposition~\ref{prop:exit}, the remaining two descent branches intersect the circle
$|w-w_0|=\fc$ in the arcs $A^{(4)}_0(w_0,\fc)$ or $A^{(4)}_2(w_0,\fc)$, while the remaining two ascent branches either originate from the source interval $[x_0,x]$ or intersect the circle $|w-w_0|=\fc$ in the arcs $A^{(4)}_1(w_0,\fc)$ or $A^{(4)}_3(w_0,\fc)$. Since these six steepest–descent and steepest–ascent branches alternate in angle and do not intersect, the only consistent configuration is the one shown in panel~(A) of \Cref{f:tangent1}: two non-real descent branches intersect the circle $|w-w_0|=\fc$ in the arc $A^{(4)}_2(w_0,\fc)$, and the two non-real ascent branches intersect the circle $|w-w_0|=\fc$ in the arcs $A^{(4)}_1(w_0,\fc)$ and $A^{(4)}_3(w_0,\fc)$, respectively.

We define the local steepest–descent contour $\mathsf D^{\rm d}(w_c)$ at $w_c$ to be the union of the non-real steepest–descent trajectories starting at $w_c$, truncated at their first exit from the disk $|w-w_0|\le \fc$. Within this disk, both non-real steepest–descent trajectories from $w_c$ exit transversely through the boundary arc $A^{(4)}_2(w_0;\fc)$. By Cauchy’s theorem, we may deform the contour $\mathsf C^{\rm d}(w_0)$ to $\mathsf D^{\rm d}(w_c)$ together with a (sub)arc of $A^{(4)}_2(w_0;\fc)$. Thanks to \eqref{e:final_copy}, for any $w$ on $A^{(4)}_2(w_0;\fc)$ we have
\begin{align}\label{e:difS_tangent1}
n\,\Re\!\big(S(w;x,s)-S(w_c;x,s)\big)
= -\,n\Bigl(|d|\,\fc^{2}-|d|\fc^2/100\Bigr)
\le -\,n\fc'.
\end{align}

\item \textbf{Frozen region}
If $(x,s)$ is in the frozen region, there are two cases $x>x_0$ or $x<x_0$.
\begin{enumerate}
\item If $x>x_0$, the interval $[x_0, x]$ is a source: the steepest-descent paths point outward.

There are two cases $(x,s)\in \fT_A$ or  $(x,s)\in \fT_B$. We will only discuss the case that $(x,s)\in \fT_B$, then $s<s_0$, and there are two real critical points determined by the intersection of the tangent line through $(x,s)$ with the $s=0$ axis, so $x<w_{c,1}<w_{c,2}$; see panel~(C) of
\Cref{f:tangent1}. The case that $(x,s)\in \fT_B$, see panel~(D) of
\Cref{f:tangent1} can be proven in the same way, so we omit.

Each of the two critical points $w_{c,1}, w_{c,2}$ has two descent and two ascent branches, forming the
usual cross pattern of a simple saddle. Along the real axis, 
the gradient flow goes from  $x$ to $w_{c,1}$, and from $w_{c,2}$ to $w_{c,1}$ and $w+\fc$. The remaining branches are off the real axis. Since gradient flows
do not intersect, the two non-real descent branches from \(w_{c,1}\) exit the
circle through arcs \(A^{(4)}_2(w_0, \fc)\), while the two non-real ascent
branches from \(w_{c,2}\) exit through arcs \(A^{(4)}_1(w_0, \fc)\) and \(A^{(4)}_3(w_0, \fc)\), as shown in
panel (C) of \Cref{f:cubic_saddle}.

Let $w_c\in \{w_{c,1}, w_{c,2}\}$, we define $\mathsf D^{\rm d}(w_c)$ to be the portions of non-real
steepest–descent trajectories starting at \(w_{c}\) up to their first exit from the disk $|w-w_0|\leq \fc$. In this case, $w_{c,2}$ is an ascent critical point, so $\mathsf D^{\rm d}(w_{c,2})=\emptyset$. By the same argument as in the first case, we may deform \(\mathsf C^{\rm d}(w_0)\) to the union of these
\(\mathsf D^{\rm d}(w_{c})\), and \eqref{e:difS_tangent1} holds.

\item If $x<x_0$,  the interval $[x_0, x]$ is a sink: the steepest-descent paths point inward, and there are two real critical points $w_{c,1}<x<w_{c,2}$; see panel (E) of \Cref{f:tangent1}.
Each has two descent and two ascent branches, forming the
usual cross pattern of a simple saddle. Along the real axis, 
the descent branches goes from  $w_{c,1}$ to $x$ and $w_0-\fc$, and from $w_{c,2}$ to $w_0$ and $w_0+\fc$. The remaining branches are off the real axis. Since gradient flows
do not intersect, the two non-real ascent branches from \(w_{c,1}\) and $w_{c,2}$ exit the
circle through arcs \(A^{(4)}_1(w_0, \fc)\) and \(A^{(4)}_3(w_0, \fc)\), as shown in
panel (E) of \Cref{f:cubic_saddle}.

Let $w_c\in \{w_{c,1}, w_{c,2}\}$, we define $\mathsf D^{\rm d}(w_c)$ to be the portions of non-real
steepest–descent trajectories starting at \(w_{c}\) up to their first exit from the disk $|w-w_0|\leq \fc$. In this case, both critical points are ascent critical points, so $\mathsf D^{\rm d}(w_{c,1})=\mathsf D^{\rm d}(w_{c,2})=\emptyset$. 

The key observation is that when $x<x_0$, by \Cref{c:PIproperty} ($b_i$ there is our $x_0$), the integrand
\begin{align}
P_{ns}(nw,nx)I_i(w)
\end{align} 
from \eqref{e:twoint} is holomorphic in a small neighborhood of $x_0$, so the circular contour 
  \(\mathsf C^{\rm d}(w_0)\) deforms to the empty set, which is also the union of these
\(\mathsf D^{\rm d}(w_{c})\).

\end{enumerate}

\item \textbf{Liquid Region}
If $(x,s)\in \fL$, there are two simple saddle points, denoted $w_{c,+}$ and $w_{c,-} = \overline{w_{c,+}}$. Again $x > x_0$, so that the interval $[x_0,x]$ is a source: the steepest–descent paths point outward; see panel~(F) of \Cref{f:tangent1}. By real-analyticity, the flow is symmetric with respect to the real axis, so it suffices to describe the configuration in the upper half-plane. From $w_{c,+}$ there emanate two steepest–descent and two steepest–ascent branches, alternating in angle.

By Proposition~\ref{prop:exit}, the two descent branches intersect the circle $|w-w_0|=\fc$ in the arcs $A^{(4)}_0(w_0,\fc)$ or $A^{(4)}_2(w_0,\fc)$, while the two ascent branches either originate from the source interval $[x_0,x]$ or intersect the circle $|w-w_0|=\fc$ in the arc $A^{(4)}_1(w_0,\fc)$. Since these four steepest–descent and steepest–ascent branches alternate in angle and do not intersect, the only consistent configuration is the one shown in panel~(F) of \Cref{f:tangent1}: the two descent branches intersect the circle $|w-w_0|=\fc$ in the arcs $A^{(4)}_0(w_0,\fc)$ and $A^{(4)}_2(w_0,\fc)$, one ascent branch intersects the circle in the arc $A^{(4)}_1(w_0,\fc)$, and the other ascent branch originates from the source interval $[x_0,x]$.

For each critical point $w_c \in \{w_{c,+},w_{c,-}\}$, we define the local steepest–descent contour $\mathsf{D}(w_c)$ to be the union of the steepest–descent trajectories starting at $w_c$, truncated at their first exit from the disk $|w-w_0|\le \fc$.  By the same argument as in the first case, we may deform \(\mathsf C^{\rm d}(w_0)\) to the union of these
\(\mathsf D^{\rm d}(w_{c})\), and \eqref{e:difS_tangent1} holds.

\end{enumerate}

Assume that \(d=S''(w_0;x_0,s_0)>0\). We discuss the case shown in
\Cref{f:vertical_tangent1}; the other case, shown in
\Cref{f:vertical_tangent2}, can be proved in the same way, so we omit it.

The steepest--descent/ascent configurations are the same as before, but all
gradient directions are reversed; see \Cref{f:tangent2}. In this case,
\(\mathsf C^{\rm d}(w_0)\) consists of two pieces: a path from
\(w_0-\ri\fc\) to \(w_0+\ri\fc\) passing to the left of \(w_0\), and a circle
surrounding \(w_0\). 

\begin{enumerate}

\item \textbf{Arctic Boundary.}
In the configuration shown in panel~(A), where
\((x,s)\in \fA\cap\fT_B\), we first deform \(\mathsf C^{\rm d}(w_0)\) to a
path from \(w_0-\ri\fc\) to \(w_0+\ri\fc\) passing to the right of \(w_0\), and
then to the union of the two non-real steepest--descent paths issuing from
\(w_c\).

In the configuration shown in panel~(B), where
\((x,s)\in \fA\cap\fT_A\) and \(x<x_0\), \Cref{c:PIproperty} implies that the
integrand
$
P_{ns}(nw,nx)I_i(w)
$
appearing in \eqref{e:twoint} is holomorphic in a small neighborhood of
\(w_0\); here the \(b_i\) in \Cref{c:PIproperty} corresponds to our \(w_0\).
Thus the circular component of \(\mathsf C^{\rm d}(w_0)\) deforms to the empty
set, and the remaining part deforms to \(\mathsf D^{\rm d}(w_c)\).

\item \textbf{Frozen Region.}
In the configuration shown in panel~(C), where \((x,s)\in\fT_B\), we first
deform \(\mathsf C^{\rm d}(w_0)\) to a path from \(w_0-\ri\fc\) to
\(w_0+\ri\fc\) passing to the right of \(w_0\), and then to
\(\mathsf D^{\rm d}(w_{c,2})\). In this case, \(w_{c,1}\) is an ascent
critical point, and \(\mathsf D^{\rm d}(w_{c,1})=\emptyset\).

In the configuration shown in panel~(D), where \((x,s)\in\fT_A\) and
\(x<x_0\), the circular component of \(\mathsf C^{\rm d}(w_0)\) again deforms
to the empty set, and the remaining part deforms to
\(\mathsf D^{\rm d}(w_{c,1})\). In this case, \(w_{c,2}\) is an ascent
critical point, and \(\mathsf D^{\rm d}(w_{c,2})=\emptyset\).

In the configuration shown in panel~(E), where \((x,s)\in\fT_B\), we keep only
the critical point \(w_{c,2}\), which originates from a tangent line from
\((x,s)\) to the portion of the arctic boundary contained in \(\fT_B\), and
ignore the critical point \(w_{c,1}\), which corresponds to a tangent line from
\((x,s)\) to the portion of the arctic boundary contained in \(\fT_A\). This is
consistent with \Cref{p:associate_critical_points}, where the critical points
are assigned. Then we first
deform \(\mathsf C^{\rm d}(w_0)\) to a path from \(w_0-\ri\fc\) to
\(w_0+\ri\fc\) passing to the right of \(w_0\), and then to
\(\mathsf D^{\rm d}(w_{c,2})\).

\item \textbf{Liquid Region.}
In the configuration shown in panel~(F), where \((x,s)\in\fL\), we deform
\(\mathsf C^{\rm d}(w_0)\) to the two descent paths from \(w_{c,+}\), one
running to the interval \([x_0,x]\) and the other to the arc
\(A^{(4)}_0(w_0,\fc)\), together with the symmetric descent paths from
\(w_{c,-}\).

\end{enumerate}

\end{proof}

\begin{proof}[Proof of \Cref{c:vertical_cusp_steepest}]

By \Cref{c:cusp_turning_critical1}, uniformly for $|w-w_0|=\fc$,
\begin{align}\label{e:Sddf}
|S(w;x,s)-S(w_c;x,s)-
d(w-w_0)^3|
\leq |d|\fc^3/100,
\quad
d:=S'''(w_0;x_0,s_0)/6.
\end{align}
with $|d|\asymp 1$. The hypotheses of
\Cref{prop:exit} hold with $(m, \delta, r)$ taken to be $(3,  |d|\fc^3/100, \fc)$.

%
%From \eqref{e:S_copy} and the Taylor expansion \eqref{e:vert_tangent} (with $c_1=0$)
%\begin{align}\begin{split}\label{e:tangent_S_exp}
%&\phantom{{}={}}S(w;x,s)-S(w_c;x,s)\\
%&=S(w;x_0,s_0)-S(w_c;x_0,s_0)+\int_{w_c}^w \ln\left(\frac{u-x}{x-s-u}\right)-\ln\left(\frac{u-x_0}{x_0-s-u}\right)\rd u\\
%&=d(w-w_0)^3 +\int_{w_c}^w \ln \frac{u-x}{u-x_0}\rd u+ \OO(|w_c-w_0|^3+|w-w_0|^4+|w_c-w_0|^4+|w-w_c||x-x_0|),
%\end{split}\end{align}
%where the error term is analytic in a small neighborhood of $w_0$.

%For $|w-w_0|\gg |x-x_0|$, we recall the estimates for the derivatives from \eqref{e:tangent_S_der},
%\begin{align}\begin{split}\label{e:tangent_S_der_copy}
%\del_w^k S(w;x,s)\lesssim \frac{|x-x_0|}{|w-x_0|^k}+1.
%\end{split}\end{align}

We recall from \eqref{e:tangent_S_der} when $k=1$, $S'(w;x,s)$ contains the following non-analytic term
\begin{align}
\ln \frac{(x-w)}{x_0-w}.
\end{align}
If $x>x_0$, the interval $[x_0, x]$ is a source term, the steepest–descent paths point outward; 
If $x<x_0$  the interval $[x, x_0]$ is a source term, the steepest–descent paths point inward; and if $x=x_0$ it vanishes and $S(w;x,s)$ is analytic locally around $w_0$.

We begin with the case  \(d=S'''(w_0;x_0,s_0)/6<0\). 
This case corresponds to \Cref{f:c_vertical_cusp2}. The contour
$\mathsf C^{\rm d}(w_0)$ consists of two pieces: a path from
$w_0+e^{-2\pi \ri/3}\fc$ to $w_0+e^{2\pi \ri/3}\fc$ passing to the
left of $w_0$, and a circular piece surrounding $w_0$, as in
\Cref{f:vertical_cusp1}.

For any critical point $w_c$, we define $\mathsf D^{\rm d}(w_c)$ to be
the union of the non-real steepest-descent trajectories starting at
$w_c$, stopped at their first exit from the disk $|w-w_0|\leq \fc$ or
their first intersection with the real axis. If no such non-real
descent trajectory emanates from $w_c$, we set
$\mathsf D^{\rm d}(w_c)=\emptyset $.

We shall repeatedly use the following estimate. Let $w_c$ denote the
critical point relevant to the contour deformation under consideration.
For every $w$ lying on one of the auxiliary subarcs contained in
$
    A_0^{(6)}(w_0,\fc)$,
    $A_2^{(6)}(w_0,\fc)$,
  and $ A_4^{(6)}(w_0,\fc)$,
the cubic expansion estimate \eqref{e:final} gives,
uniformly in the cases considered below,
\begin{align}\label{e:difS_cusp1}
n\,\Re\!\bigl(S(w;x,s)-S(w_c;x,s)\bigr)
&\leq 
-n\Bigl(|d|\,\fc^{3}
-3|d|\fc^3/100\Bigr)\leq -n\fc'
\end{align}
on all such auxiliary arc. 

We now distinguish the three possible locations of $(x,s)$.

\begin{enumerate}
\item \textbf{Arctic curve.}
Suppose first that $(x,s)\in\fA$.

Assume first that $x>x_0$, so that $(x,s)\in\fT_B$. Then the interval
$[x_0,x]$ is a source: the steepest-descent paths point outward. There
are two distinct real critical points
$
    w_{c,1}<x_0<w_{c,2},
$
where $w_{c,2}$ has multiplicity two; see Panel (A) of
\Cref{f:vertical_cusp1}. From $w_{c,2}$ emanate three descent and three
ascent branches, while from $w_{c,1}$ emanate two descent and two ascent
branches, alternating in angle. Along the real axis, the flow goes from
$w_0-\fc$ and $w_0$ to $w_{c,1}$, from $x$ to $w_{c,2}$, and from
$w_{c,2}$ to $w_0+\fc$.

The remaining branches are non-real. The four non-real branches from
$w_{c,2}$, two descent and two ascent, meet the circle
$|w-w_0|=\fc$ alternately through
$
    A^{(6)}_1(w_0,\fc)$,
    $A^{(6)}_2(w_0,\fc)$,
    $A^{(6)}_4(w_0,\fc)$,
    and $A^{(6)}_5(w_0,\fc)$.
And the two non-real descent branches from $w_{c,1}$ exit
through $A^{(6)}_2(w_0,\fc)$ and $A^{(6)}_4(w_0,\fc)$. This follows
from the fact that gradient-flow trajectories do not intersect.

By \Cref{p:associate_critical_points}, we assign to $(x,s)$ only those
critical points whose tangent lines from $(x,s)$ are tangent to the
portion of the arctic boundary contained in $\fT_B$. Therefore
$w_{c,2}$ is assigned to $(x,s)$, whereas $w_{c,1}$ is not.

In this case, the contour $\mathsf C^{\rm d}(w_0)$ can first be
deformed to a path from $w_0+e^{-2\pi \ri/3}\fc$ to
$w_0+e^{2\pi \ri/3}\fc$ passing to the right of $w_0$, and then to
$\mathsf D^{\rm d}(w_{c,2})$, together with possible subarcs of
$A^{(6)}_2(w_0,\fc)$ and $A^{(6)}_4(w_0,\fc)$. These subarc
contributions are exponentially small by \eqref{e:difS_cusp1}.

Assume next that $x<x_0$, so that $(x,s)\in\fT_A$. Then the interval
$[x,x_0]$ is a sink. There are two distinct real critical points
$
    w_{c,1}<x_0<w_{c,2}$,
where $w_{c,1}$ has multiplicity two; see Panel (B) of
\Cref{f:vertical_cusp1}. From $w_{c,1}$ emanate three descent and three
ascent branches, while from $w_{c,2}$ emanate two descent and two ascent
branches, alternating in angle. Along the real axis, the flow goes from
$w_0-\fc$ to $w_{c,1}$, from $w_{c,1}$ to $x$, and from $x_0$ and
$w_0+\fc$ to $w_{c,2}$.

The remaining branches are non-real. The four non-real branches from
$w_{c,1}$, two descent and two ascent, meet the circle
$|w-w_0|=\fc$ alternately through
$
    A^{(6)}_1(w_0,\fc)$,
    $A^{(6)}_2(w_0,\fc)$,
    $A^{(6)}_4(w_0,\fc)$,
    and $A^{(6)}_5(w_0,\fc)$.
And the two non-real ascent branches from $w_{c,2}$ exit
through $A^{(6)}_1(w_0,\fc)$ and $A^{(6)}_5(w_0,\fc)$.

By \Cref{p:associate_critical_points}, only tangent lines to the
portion of the arctic boundary contained in $\fT_A$ are relevant.
Therefore $w_{c,1}$ is assigned to $(x,s)$, whereas $w_{c,2}$ is not.

Moreover, in this case the integrand of the first integral in
\eqref{e:twoint} is holomorphic in a small neighborhood of $w_0$.
Therefore the circular component of $\mathsf C^{\rm d}(w_0)$ deforms to
the empty contour. The remaining component of $\mathsf C^{\rm d}(w_0)$
can be deformed to $\mathsf D^{\rm d}(w_{c,1})$, together with possible
subarcs of $A^{(6)}_2(w_0,\fc)$ and $A^{(6)}_4(w_0,\fc)$. These subarc
contributions are exponentially small by \eqref{e:difS_cusp1}.

\item \textbf{Liquid region.}
Suppose next that $(x,s)\in\fL$.

Assume first that $x>x_0$. Then $[x_0,x]$ is a source. There is one
real saddle $w_{c,1}<w_0$ and a conjugate pair of non-real critical
points $w_{c,\pm}$, with
$
    \overline{w_{c,+}}=w_{c,-}$;
see Panel (C) of \Cref{f:vertical_cusp1}. Along the real axis, the
descent branches go from $w_0-\fc$ and $w_0$ to $w_{c,1}$, and from
$x$ to $w_0+\fc$.

By real-analyticity, the configuration is symmetric with respect to the
real axis, so it suffices to describe the upper half-plane. From
$w_{c,+}$ emanate two descent and two ascent branches, forming the usual
alternating pattern. These four branches remain in the upper half-plane,
since gradient-flow trajectories cannot cross the real axis. By
Proposition~\ref{prop:exit}, two descent branches and one ascent branch
meet $|w-w_0|=\fc$ through
$
    A^{(6)}_0(w_0,\fc)$,
   $A^{(6)}_1(w_0,\fc)$,
   and $A^{(6)}_2(w_0,\fc)$
in alternating order. The remaining ascent branch either exits through
$A^{(6)}_3(w_0,\fc)$ or connects to the real critical point $w_{c,1}$
or to the source $[x_0,x]$. The flows from $w_{c,-}$ are the complex
conjugates of those from $w_{c,+}$.

It follows that $\mathsf C^{\rm d}(w_0)$ can be deformed to
$
    \mathsf D^{\rm d}(w_{c,+})
    \cup
    \mathsf D^{\rm d}(w_{c,-}),
$
together with possible subarcs of $A^{(6)}_0(w_0,\fc)$ and
$A^{(6)}_2(w_0,\fc)$. These subarc contributions are exponentially
small by \eqref{e:difS_cusp1}.

Assume next that $x<x_0$. Then $[x,x_0]$ is a sink. There is one real
saddle $w_{c,1}>w_0$ and a conjugate pair of non-real critical points
$w_{c,\pm}$, with
$
    \overline{w_{c,+}}=w_{c,-};
$
see Panel (D) of \Cref{f:vertical_cusp1}. Along the real axis, the
descent branches go from $w_0-\fc$ to $x$, and from $w_{c,1}$ to $w_0$
and $w_0+\fc$.

Again, by symmetry, it suffices to describe the upper half-plane. From
$w_{c,+}$ emanate two descent and two ascent branches, alternating in
angle, and all four branches remain in the upper half-plane. By
Proposition~\ref{prop:exit}, two ascent branches and one descent branch
meet $|w-w_0|=\fc$ through
$
    A^{(6)}_1(w_0,\fc)$,
   $A^{(6)}_2(w_0,\fc)$,
   and $A^{(6)}_3(w_0,\fc)$
in alternating order. The remaining descent branch either exits through
$A^{(6)}_0(w_0,\fc)$ or connects to the real critical point $w_{c,1}$
or to the sink $[x,x_0]$. The flows from $w_{c,-}$ are the complex
conjugates of those from $w_{c,+}$.

Thus $\mathsf C^{\rm d}(w_0)$ can be deformed to
$
    \mathsf D^{\rm d}(w_{c,+})
    \cup
    \mathsf D^{\rm d}(w_{c,-}),
$
together with possible subarcs of $A^{(6)}_0(w_0,\fc)$ and
$A^{(6)}_2(w_0,\fc)$. These subarc contributions are exponentially
small by \eqref{e:difS_cusp1}.

\item \textbf{Frozen region.}
Finally suppose that $(x,s)$ lies in the frozen region.

Assume first that $x>x_0$. Then $[x_0,x]$ is a source. There are three
real critical points
$
    w_{c,1}<w_{c,2}<x_0<w_{c,3};
$
see Panel (E) of \Cref{f:vertical_cusp1}. Each critical point has two
descent and two ascent branches, forming the usual cross pattern of a
simple saddle. Along the real axis, the descent branches go from
$w_0-\fc$ and $w_0$ to $w_{c,1}$, from $x$ to $w_{c,2}$, and from
$w_{c,3}$ to $w_{c,2}$ and $w_0+\fc$.

The remaining branches are non-real. Since gradient-flow trajectories
do not intersect, the non-real descent branches from $w_{c,1}$ and
$w_{c,2}$ exit the circle through
$A^{(6)}_2(w_0,\fc)$ and $A^{(6)}_4(w_0,\fc)$, while the non-real
ascent branches from $w_{c,3}$ exit through
$A^{(6)}_1(w_0,\fc)$ and $A^{(6)}_5(w_0,\fc)$.

By \Cref{p:associate_critical_points}, the critical points assigned to
$(x,s)$ are $w_{c,2}$ and $w_{c,3}$, but not $w_{c,1}$. However,
$
    \mathsf D^{\rm d}(w_{c,3})=\emptyset$,
because the non-real branches from $w_{c,3}$ are ascent branches. Hence
only $w_{c,2}$ contributes to the descent contour. Therefore
$\mathsf C^{\rm d}(w_0)$ can be deformed to
$\mathsf D^{\rm d}(w_{c,2})$, together with possible subarcs of
$A^{(6)}_2(w_0,\fc)$ and $A^{(6)}_4(w_0,\fc)$. These subarc
contributions are exponentially small by \eqref{e:difS_cusp1}.

Assume finally that $x<x_0$. Then $[x,x_0]$ is a sink. There are three
real critical points, ordered as in Panel (F) of
\Cref{f:vertical_cusp1}; in particular,
$
    w_{c,1}<w_{c,2}<x<w_{c,3}.
$
Each critical point has two descent and two ascent branches, forming
the usual cross pattern of a simple saddle. Along the real axis, the
descent branches go from $w_0-\fc$ to $w_{c,1}$, from $w_{c,2}$ to
$w_{c,1}$ and $x$, and from $w_{c,3}$ to $w_0$ and $w_0+\fc$.

The remaining branches are non-real. Since gradient-flow trajectories
do not intersect, the two non-real descent branches from $w_{c,1}$ exit
the circle through $A^{(6)}_2(w_0,\fc)$ and $A^{(6)}_4(w_0,\fc)$,
whereas the non-real ascent branches from $w_{c,2}$ and $w_{c,3}$ exit
through $A^{(6)}_1(w_0,\fc)$ and $A^{(6)}_5(w_0,\fc)$.

By \Cref{p:associate_critical_points}, the critical points assigned to
$(x,s)$ are $w_{c,1}$ and $w_{c,2}$, but not $w_{c,3}$. However,
$
    \mathsf D^{\rm d}(w_{c,2})=\emptyset,
$
because the non-real branches from $w_{c,2}$ are ascent branches. Thus
only $w_{c,1}$ contributes to the descent contour.

Moreover, the integrand of the first integral in \eqref{e:twoint} is
holomorphic in a small neighborhood of $w_0$. Hence the circular
component of $\mathsf C^{\rm d}(w_0)$ deforms to the empty contour. The
remaining component of $\mathsf C^{\rm d}(w_0)$ can be deformed to
$\mathsf D^{\rm d}(w_{c,1})$, together with possible subarcs of
$A^{(6)}_2(w_0,\fc)$ and $A^{(6)}_4(w_0,\fc)$. These subarc
contributions are exponentially small by \eqref{e:difS_cusp1}.
\end{enumerate}

This completes the proof in the case $d<0$. 
It remains to consider the case
$
    d={S'''(w_0;x_0,s_0)}/{6}>0$.
Then the local cubic term has the opposite sign, so the corresponding
gradient-flow diagram is obtained from the one above by reversing all
arrows; equivalently, the steepest-descent and steepest-ascent
directions are interchanged. See \Cref{f:vertical_cusp2}. The
classification of the critical points, the assignment rule from
\Cref{p:associate_critical_points}, the contour deformations, and the
exponential estimates on the auxiliary arcs are therefore identical
after this reversal. Hence the same argument proves the desired contour
deformation in the case $d>0$ as well. This completes the proof.

\end{proof}

\subsection{Vertical tangent frozen chart}
\label{s:vertical_frozen_chart_proof}

\begin{proof}[Proof of \Cref{c:vertical_frozen_critical1}]

Critical points $w_c$ of $S(w;x,s)$ satisfy the following equation 

\begin{align}
x=w+s\chi(w)=w_0+ s\chi(w_0)+(w-w_0)+ \sum_{k\geq 1}\frac{sa_k(w-w_0)^k}{k!}.
\end{align}
We can rewrite it as
\begin{align}
\frac{(x-x_0)}{1+s\chi'(w_0)}- (w-w_0)-\sum_{k\geq 2}\frac{a_k(w-w_0)^k}{(1+s\chi'(w_0))k!}=0,
\end{align}
which is in the following form: let $z=w-w_0$ 
\begin{align}
F(z;a,b):=-z-b+\sum_{k\ge 2} c_k z^k=0,\quad  \quad b=-\frac{x-x_0}{1+s\chi'(w_0)}.
\end{align}

For $\|(x,s)-(x_0,s_0)\|_2\leq \delta$, we have $|1+s\chi'(w_0)|\asymp 1$. By \Cref{lem:bulk-perturbation}, we have 
\begin{align}
w_c-w_0=w_c-x_0=-b+\OO(b^2)=\frac{x-x_0}{1-s/s_0'}+\OO(|x-x_0|^2).
\end{align}
Then either $1-s/s_0'\asymp 1$
\begin{align}
|w_c-x_0|-|x-x_0|=\frac{(s/s_0')|x-x_0|}{1-s/s_0'}+\OO(|x-x_0|^2)\asymp |x-x_0|.
\end{align}
Or $s/s_0'-1\asymp 1$
\begin{align}
w_c-w_0=w_c-x_0=-\frac{x-x_0}{s/s_0'-1}+\OO(\delta^2)\asymp -(x-x_0).
\end{align}
In both cases we have 
\begin{align}\label{e:critical_distance}
\dist(w_c,[\min\{x_0, x\}, \max\{x_0, x\}]) \asymp |x-x_0|\lesssim \delta.
\end{align}

We recall the expansion of $S(w;x_0, s_0)$ from \eqref{e:S_extend}. Then \eqref{e:hor_c_2} implies
\begin{align}
d:= S'(w_0;x_0,s_0)=\ln(s_0'/s_0),\quad 
|d|=\asymp 1. 
\end{align}
The estimate \eqref{e:critical_distance}, together with \eqref{e:Scritical1}, we conclude that for $|w-w_0|=r$,
\begin{align}
S(w; x,s) - S(w_c; x,s)=d(w-w_0)+\OO(\delta+r^{2}+\delta  \ln(1/\delta )),\quad d:=S'(w_0;x_0, s_0).
\end{align}
Taking $r=\fc$ and then choosing $\fc$ and $\delta$ sufficiently small gives
\eqref{e:tangent_frozen_err}.

\end{proof}

\begin{proof}[Proof of \Cref{l:local_descent_deformation}]

By \Cref{c:vertical_frozen_critical1}, uniformly for $|w-w_0|=\fc$,
\begin{align}\label{e:final_copy2}
|S(w; x,s) - S(w_c; x,s)-d(w-w_0)|\leq |d|\fc/100,\quad d:=S'(w_0;x_0, s_0).
\end{align}
with $|d|\asymp 1$. The hypotheses of
\Cref{prop:exit} hold with $(m, \delta, r)$ taken to be $(1,  |d|\fc/100, \fc)$.

\noindent\textbf{Replacing \(\mathsf C^{\rm d}(w_0)\) by \(\mathsf D^{\rm d}(w_c)\).}

We first record the four possible local configurations.

\begin{enumerate}
\item[(A)] Suppose that $s>s_0$ and $x<x_0$, as in Panel (A) of
\Cref{f:side}. Then the interval $[x,x_0]$ is a sink: the
steepest-descent paths point inward. The critical point $w_c$ has two
descent branches and two ascent branches, forming the usual cross
pattern of a simple saddle. Along the real axis, the gradient flow goes
from $w_0-\fc$ to $x$, and from $w_c$ to $w_0$ and $w_0+\fc$. The
remaining branches are non-real. The two non-real ascent branches from
$w_c$ exit the circle through arcs $A^{(2)}_1(w_0,\fc)$. Thus $w_c$ is an
ascent critical point for the present descent-contour construction, and
we set
$
    \sfD^{\rm d}(w_c)=\emptyset $.

\item[(B)] Suppose that $s>s_0$ and $x>x_0$, as in Panel (B) of
\Cref{f:side}. Then the interval $[x_0,x]$ is a source: the
steepest-descent paths point outward. There is one real critical point
$w_c<x_0$. It has two descent branches and two ascent branches, again
forming the usual cross pattern of a simple saddle. Along the real axis,
the gradient flow goes from $w_0-\fc$ and $w_0$ to $w_c$, and from $x$
to $w_0+\fc$. The remaining branches are non-real. The two non-real
descent branches from $w_c$ exit the circle through arcs
$A^{(2)}_0(w_0,\fc)$. Thus $w_c$ is a descent critical point. We define
$\sfD^{\rm d}(w_c)$ to be the union of the portions of the
steepest-descent trajectories starting at $w_c$ up to their first exit
from the disk $|w-w_0|\leq \fc$.

\item[(C)] Suppose that $s<s_0$ and $x<x_0$, as in Panel (C) of
\Cref{f:side}. Then the interval $[x,x_0]$ is a sink. There is one
ascent critical point $w_c<x_0$. Along the real axis, the gradient flow
goes from $w_c$ to $w_0$ and $w_0-\fc$, and from $w_0+\fc$ to $w_0$.
The remaining branches are non-real. The two non-real ascent branches
from $w_c$ exit the circle through arcs $A^{(2)}_1(w_0,\fc)$. Hence we set
$
    \sfD^{\rm d}(w_c)=\emptyset$ .

\item[(D)] Suppose that $s<s_0$ and $x>x_0$, as in Panel (D) of
\Cref{f:side}. Then the interval $[x_0,x]$ is a source. There is one
descent critical point $w_c>x_0$. Along the real axis, the gradient
flow goes from $w_0$ to $w_0-\fc$, and from $x$ and $w_0+\fc$ to
$w_c$. The remaining branches are non-real. The two non-real descent
branches from $w_c$ exit the circle through arcs $A^{(2)}_0(w_0,\fc)$. We
define $\sfD^{\rm d}(w_c)$ to be the union of the portions of the
steepest-descent trajectories starting at $w_c$ up to their first exit
from the disk $|w-w_0|\leq \fc$.
\end{enumerate}

We now distinguish two geometric cases.

\begin{enumerate}
\item Suppose first that $(x_0,s_0)$ lies in the closure of a single
curvilinear triangle $\fT$, and not on the common boundary of two
adjacent curvilinear triangles. Without loss of generality, assume that
on $\fT$ we have
\[
    \nabla H^*=(1,0).
\]
The argument splits according to whether $s>s_0$ or $s<s_0$.

\begin{enumerate}
\item Suppose that $s>s_0$. In this case,
$\mathsf C^{\rm d}(w_0)=\emptyset$, since by
\eqref{e:frozen_S''3} we have $S''<0$.

In the configuration of Panel (A) of
\Cref{f:side}, we have already set
$\sfD^{\rm d}(w_c)=\emptyset$. In the configuration of Panel (B) of
\Cref{f:side}, the
tangent line corresponding to $w_c$ is not tangent to the portion of the
arctic boundary contained in $\fT$, whose slopes lie in
$[-\infty,0]$. By \Cref{p:associate_critical_points}, we assign to
$(x,s)$ only those critical points whose tangent lines from $(x,s)$ are
tangent to the portion of the arctic boundary contained in $\fT$.
Therefore $w_c$ is not assigned to $(x,s)$ in this case. This is consistent with
$
    \sfC^{\rm d}(w_0)=\emptyset$ .

\item Suppose that $s<s_0$. In this case,
$\mathsf C^{\rm d}(w_0)$ consists of a circle surrounding $w_0$, since
by \eqref{e:frozen_S''3} we have $S''>0$.

In the configuration of Panel (C) of
\Cref{f:side}, we have $x<x_0$. The tangent line
corresponding to $w_c$ is not tangent to the portion of the arctic
boundary contained in $\fT$, whose slopes lie in $[-\infty,0]$.
Therefore, by \Cref{p:associate_critical_points}, the critical point
$w_c$ is not assigned to $(x,s)$.

Moreover, when $x<x_0$, \Cref{c:PIproperty} applies with the point
$b_i$ there equal to our $x_0$. Hence the integrand
\[
    P_{ns}(nw,nx) I_i(w)
\]
appearing in \eqref{e:twoint} is holomorphic in a small neighborhood of
$x_0$. Therefore the circular contour $\sfC^{\rm d}(w_0)$ can be
deformed to the empty contour.

It remains to consider the configuration of Panel (D) of
\Cref{f:side}. By Cauchy's
theorem, the contour $\sfC^{\rm d}(w_0)$ can be deformed to
$\sfD^{\rm d}(w_c)$ together with a subarc of $A^{(2)}_0(w_0;\fc)$. For
$w\in A^{(2)}_0(w_0;\fc)$, \eqref{e:final} gives
\begin{align}\label{e:difS_tangent3}
n\,\Re[S(w;x,s)-S(w_c;x,s)]
&=
-n\Bigl(|d|\,\fc-3|d|\fc/100\Bigr)
\leq -n\fc'
\end{align}
for all $w\in A^{(2)}_0(w_0;\fc)$.
\end{enumerate}

\item Suppose next that $(x_0,s_0)$ lies on the common boundary of two
adjacent curvilinear triangles $\fT_A$ and $\fT_B$, with
\[
    \nabla H^*=(0,0) \quad \text{on } \fT_A,
    \qquad
    \nabla H^*=(1,0) \quad \text{on } \fT_B;
\]
see \Cref{f:adjacent_curvilinear_triangle}. In this case,
$\mathsf C^{\rm d}(w_0)$ always consists of a circle surrounding
$w_0$.

In the configurations of Panels (A) and (C) of
\Cref{f:side}, we have $x<x_0$ and
$\sfD^{\rm d}(w_c)=\emptyset$. By \Cref{c:PIproperty}, again with
$b_i$ there equal to our $x_0$, the integrand in the first integral in
\eqref{e:twoint} is holomorphic in a small neighborhood of $x_0$.
Therefore the circular contour $\sfC^{\rm d}(w_0)$ can be deformed to
the empty contour.

In the configuration of Panel (B) of
\Cref{f:side}, we have $(x,s)\in\fT_A$. Hence, by
\Cref{p:associate_critical_points}, the critical point $w_c$ is
assigned to $(x,s)$. By Cauchy's theorem, the contour
$\sfC^{\rm d}(w_0)$ can be deformed to $\sfD^{\rm d}(w_c)$, up to
possible subarcs of $A^{(2)}_0(w_0;\fc)$. The contribution of any such subarc
is exponentially small by the estimate \eqref{e:difS_tangent3}.

In the configuration of Panel (D) of
\Cref{f:side}, we have $(x,s)\in\fT_B$. Again, by
\Cref{p:associate_critical_points}, the critical point $w_c$ is
assigned to $(x,s)$. The same contour deformation and the same estimate
\eqref{e:difS_tangent3} show that $\sfC^{\rm d}(w_0)$ may be replaced,
up to an exponentially small contribution, by $\sfD^{\rm d}(w_c)$.
\end{enumerate}

\noindent\textbf{Replacing \(\mathsf D^{\rm d}(w_c)\) by \(\mathsf S^{\rm d}(w_c)\).}
We introduce the shifted coordinate
\begin{align}\label{e:shift_coordinate}
  \sfb = x - x_0, \quad
\end{align}

We discuss two cases separately
\begin{enumerate}

\item \textbf{Close to the tangent line} $|\sfb|\leq (\ln n)^3/n$.
Let
$r_n:=(\ln n)^5 /n$. 
We define the local steepest–descent set $\mathsf S(w_c)$ at \(w_c\) to be the portions of
steepest–descent trajectories starting at \(w_0\) up to their first exit from
the disk $\{w: |w-w_0|\leq r_n\}$. 

We recall the expression \eqref{e:S_extend}, integrating it
\begin{align}\begin{split}\label{e:S_extend2}
S(w;(x_0,s))-S(w_0;(x_0,s))
&=\ln \frac{s'_0}{s_0}(w-w_0)+\OO((w-w_0)^2).
\end{split}\end{align}

On the circle $\{|w-w_0|=r_n\}$, we have
\begin{align}\begin{split}\label{e:tangent_F_exp1}
S(w;x,s)-S(w_{c};x,s)
&=S(w;(x_0,s))-S(w_{c};(x_0,s))\\
&+\int_{w_c}^w \ln \frac{u-x}{u-(x-s)} -\ln \frac{u-x_0}{u-(x_0-s)}\rd u\\
&=\ln(s_0/s)(w-w_0)+\int_{w_c}^w \ln \frac{u-x}{u-x_0}\rd u\\
&+\OO(|\ln(s_0/s)||w_c-w_0|+|w_c-w_0|^2+r_n^2+r_n|x-x_0|)
\end{split}\end{align}
where the first statement follows from \eqref{e:S_copy};
and the second statement follows from the Taylor expansion \eqref{e:S_extend2}, and $|\ln ((u-(x_0-s_0))/(u-(x-s)))|\lesssim |x-x_0|$. For the integral in \eqref{e:tangent_F_exp1}, we can bound it as
\begin{align}\begin{split}\label{e:t_integral_bound}
\int_{w_c}^w \ln \frac{u-x}{u-x_0}\rd u
&=\left.(u-x)\ln (u-x)- (u-x_0)\ln(u-x_0)\right|_{w_c}^w\\
&=\left.(u-x)\ln \frac{u-x}{u-x_0}+(x-x_0)\ln (u-x_0)\right|_{w_c}^w\\
&=\OO( |x-x_0|\ln(1/|w_c-x_0|))\lesssim |\sfb|\ln(1/|\sfb|)\lesssim (\ln n)^4/n
\end{split}\end{align}
where the first two statements from direct computation; in the third statement we used \eqref{e:wcest}, $|w_c-x_0|\asymp|x-x_0|$, and $|w-x_0|=r_n\gg  |\sfb|=|x-x_0|$, and bounding $|\ln(u-x)/(u-x_0)|\lesssim |x-x_0|/|u-x_0|\asymp |x-x_0|/|u-x|$; the last two statements also from direct computation. By plugging \eqref{e:t_integral_bound} into \eqref{e:tangent_F_exp1}, we conclude 
\begin{align}\begin{split}\label{e:frozen_tangent}
S(w;x,s)-S(w_c;x,s)
&=\ln(s_0/s)(w-w_0)+\OO((\ln n)^4/n)+\OO(|w_c-w_0|+r_n^2+r_n|\sfb|)\\
&=\ln(s_0/s)(w-w_0)+  \OO((\ln n)^{4}/n)
\end{split}\end{align}

The estimate \eqref{e:frozen_tangent} verifies the assumptions in \Cref{prop:exit} with $(m, \delta, r)$ taken to be $(1, (\ln n)^4/n, r_n)$. 
It follows that the two steepest–descent branches from \(w_c\) exit the circle
\(\{w: |w-w_0|=r_n\}\) at points on \(A^{(2)}_1(w_0;r_n)\).  At any such
point,
\begin{equation}\label{eq:sd-drop0}
n\,\Re\!\big(S(w; x,s)-S(w_c; x,s)\big)
= -\,n|d|\,r_n+\OO((\ln n)^4)
\le -\,\fc'\,(\ln n)^{5}.
\end{equation}
Moreover, since $\Re[S(w; x,s)]$ decreases along the steepest–descent, \eqref{eq:sd-drop0}
holds for all \(w\in\mathsf D(w_c)\setminus \mathsf S(w_c)\).
Consequently, for  \(w\in \mathsf D(w_c)\setminus \mathsf S(w_c)\) we have
\[
e^{n\Re[S(w;x,s)]}
\leq
e^{n\Re[S(w_c;x,s)]}e^{-\fc'(\ln n)^5}.
\]

\item \textbf{Other Frozen} $|\sfb|\geq (\ln n)^3/n$.

There is another critical point $|w_c-w_0|\asymp |\sfb|$ and $|\chi'(w_c)+1/s|\asymp 1$. 
Locally around $w_c$, $S(w_c; (x,s))$ is holomorphic and we have the Taylor expansion 
\begin{align}\begin{split}\label{e:vert_side_S0}
S(w;x,s)-S(w_{c};x,s)
&=d (w-w_{c})^2 +\OO\left(\frac{|\sfb||w-w_c|^3}{|w_c-w_0|^3}\right)
\end{split}\end{align}
where the error term is from \eqref{e:tangent_S_der}, and $d=S''(w_c; (x,s))/2$. 
Thanks to \eqref{e:derSsecond} and \eqref{e:chi_v_tangent}
\begin{align}\begin{split}\label{e:vert_side_S1}
&|d|\asymp|S''(w_{c};x,s)|=\frac{|\chi'(w_{c})+1/s|}{|\chi(w_{c})(1-\chi(w_{c})|}\asymp \frac{1}{|w_c-w_0|},
\end{split}
\end{align}

By plugging \eqref{e:vert_side_S1} into \eqref{e:vert_side_S0}, we get
\begin{align}\begin{split}\label{e:vert_side_S2}
S(w;x,s)-S(w_{c};x,s)
&=\Theta\left(1/|\sfb|\right)(w-w_{c})^2+\OO\left(\frac{|w-w_c|^3}{|\sfb|^2}\right),
\end{split}\end{align}

Let
$r_n:=(\ln n)^2 (|\sfb|/n)^{1/2}$. 
We define the local steepest–descent set $\mathsf S(w_c)$ at \(w_c\) to be the portions of
steepest–descent trajectories starting at \(w_c\) up to their first exit from
the disk $\{w: |w-w_c|\leq r_n\}$.  
For any $w=w_c+r_n e^{\ri\theta}$, we have
\begin{align}\begin{split}
&\phantom{{}={}}S(w;x,s)-S(w_c;x,s)=\Theta(1/|\sfb|)(w-w_c)^2+\OO(|w-w_c|^3/|\sfb|^2)\\
&=\Theta(1/|\sfb|)(w-w_c)^2+\OO((\ln n)^3/(|\sfb|^{1/2}n^{3/2}))
=\Theta(1/|\sfb|)(w-w_c)^2+\OO((\ln n)^{3/2}/n)
\end{split}\end{align}
The above estimate verifies the assumptions in \Cref{prop:exit} with $(m, \delta, r)$ taken to be $(2, (\ln n)^{3/2}/n, r_n)$. We can therefore replace $\sfD^{\rm d}(w_c)$ by $\sfS^{\rm d}(w_c)$ as in the
close to tangent case.
\end{enumerate}

\end{proof}

\subsection{Horizontal tangent chart and cusp-turning chart}
\label{s:horizontal_tangent_chart_proof}
\begin{proof}[Proof of \Cref{c:tangent_critical3}]
Let $(x_0,s_0)\in \fA$ be a horizontal tangency point which is not a cusp
point. Then $f(x_0,s_0)=-1$ and $\chi(x_0,s_0)=\infty$. As in
\eqref{e:chi-expand-horizontal}, we change coordinates by setting
\begin{align}\label{e:change_variable}
\wt w=\frac{1}{x_0-w},
\qquad
\wt S(\wt w;x,s):=S(w;x,s),
\qquad
\wt\chi(\wt w):=\chi(w)
=
\frac{1}{s_0\wt w}+\sum_{i\geq 1}c_i\wt w^i .
\end{align}
Since $(x_0,s_0)$ is not a cusp point, $c_1\neq 0$. In this new coordinate,
 from \eqref{e:hor_tangent} we have 
 \begin{align} \label{e:hor_tangent_copy}
 \wt S'(\wt w;x_0,s_0)=s_0^2 c_1\wt w +(c_1s_0^3+c_2s_0^2)\wt w^2+\OO(\wt w^3). 
 \end{align}
and the critical point of $\wt S(\wt w;x_0,s_0)$ is at $\wt w=0$. Then \eqref{e:hor_tangent} and \eqref{e:hor_c_1} implies
\begin{align}
d:=\frac{1}{2}\wt S''(0;x_0,s_0),\quad 
|d|=\frac{1}{2}|c_1 s^2_0|\asymp 1. 
\end{align}

Fix $\|(x,s)-(x_0,s_0)\|_2\leq \delta$. We study $\wt S(w;x,s)$ in a small
neighborhood of $0$. From \eqref{e:critical_point2}, the critical points of
$\wt S(w;x,s)$ satisfy
\begin{align}\label{e:x_exp3}
0
=
s\wt w\wt\chi(\wt w)-(x-x_0)\wt w-1
=
\frac{s}{s_0}
+\sum_{i\geq 1}sc_i\wt w^{i+1}
-(x-x_0)\wt w
-1 .
\end{align}
Equivalently,
\begin{align}
-\frac{s_0-s}{c_1ss_0}
-\frac{x-x_0}{c_1s}\wt w
+\wt w^2
+\sum_{i\geq 2}\frac{c_i}{c_1}\wt w^{i+1}
=0.
\end{align}
Thus the critical point equation has the form
\begin{equation}\label{eq:F}
F(\wt w;a,b):=
\wt w^2-a\wt w-b+\sum_{k\geq 3}c_k\wt w^k=0,
\qquad
a=\frac{x-x_0}{c_1s},
\qquad
b=\frac{s_0-s}{c_1ss_0}.
\end{equation}
For $\delta$ sufficiently small, \Cref{lem:edge-perturbation} implies that
$F(\wt w;a,b)$ has exactly two zeros in a small disk around $0$, counted with
multiplicity. Moreover, for any such zero $\wt w_c$,
\begin{align}\label{e:h-tangent}
\frac{|s-s_0|}{\sqrt{\delta}}
\lesssim
\min\left\{
\frac{|s-s_0|}{|x-x_0|},
|s-s_0|^{1/2}
\right\}
\lesssim
|\wt w_c|
\lesssim
|x-x_0|+\sqrt{|s-s_0|}
\lesssim
\sqrt{\delta}.
\end{align}
Here, if $x=x_0$, the first term inside the minimum is interpreted as $+\infty$.
In particular, $|\wt w_c|\gtrsim |s-s_0|$.

We now apply \eqref{e:Scritical2}. Since $x_0'=x_0$ in the present case, the
condition $|x_0-x_0'|r\leq 1/4$ is automatic; we choose $r$ small enough so
that $|s_0|r\leq 1/4$. For $|\wt w|=r$, \eqref{e:Scritical2} gives
\begin{align}
\wt S(\wt w;x,s)-\wt S(\wt w_c;x,s)
=
d \wt w^2
+
\OO\!\left(
|\wt w_c|^2+r^3+\delta+|s-s_0|\log(1/|\wt w_c|)
\right),
\quad
d:=\frac12 \wt S''(0;x_0,s_0).
\end{align}
Using \eqref{e:h-tangent}, we have $|w_c|^2\lesssim \delta$ and
$
|s-s_0|\log(1/|\wt w_c|)
\lesssim
\delta\log(1/\delta).
$
Hence, for $|\wt w|=r$,
\begin{align}
\wt S(\wt w;x,s)-\wt S(\wt w_c;x,s)
=
d \wt w^2
+
\OO\!\bigl(\delta\log(1/\delta)+r^3\bigr).
\end{align}
Taking $r=\fc$ and then choosing $\fc$ and $\delta$ sufficiently small gives
\eqref{e:tangent_err}.

%The estimate \eqref{e:h-tangent}, together with \Cref{c:S_diff}, we conclude that the statements in \Cref{prop:exit} hold with $(m, \delta, r)$ taken to be $(2, \delta\ln (1/\delta)+\fc^3, \fc)$.  Moreover, in this case $d\in\mathbb{R}\setminus\{0\}$ (recall from \eqref{eq:Taylor}). 

\end{proof}

\begin{proof}[Proof of \Cref{c:cusp_turning_critical3}]

For any horizontal cusp location $(x_0, s_0)\in\fA$, the same as in \eqref{e:hor_tangent_copy}, and notice that $c_1=0$,  we have 
\begin{align}\begin{split}\label{e:hor_cusp}
\wt S'({\wt w}; x_0, s_0)=c_2 s_0^2 {\wt w}^2 +\OO({\wt w}^3).
\end{split}\end{align}

Fix $|x-x_0|+|s-s_0|\leq \delta$, we study the $\wt S({\wt w};x,s)$ in a small neighborhood of $0$. From \eqref{e:critical_point2}, the critical points of
$\wt S(w;x,s)$ satisfy
\begin{align}
0
=
s\wt w\wt\chi(\wt w)-(x-x_0)\wt w-1
=
\frac{s}{s_0}
+\sum_{i\geq 1}sc_i\wt w^{i+1}
-(x-x_0)\wt w
-1 .
\end{align}
Using that  $c_1=0$ and $c_2\neq 0$, we can rewrite the above equation as
\begin{align}
-\frac{s_0-s}{c_2ss_0}
-\frac{x-x_0}{c_2s}\wt w
+\wt w^3
+\sum_{i\geq 3}\frac{c_i}{c_2}\wt w^{i+1}
=0,
\end{align}
which is in the following form
\begin{equation}\label{eq:F}
F(\wt w;a,b):={\wt w}^3-a {\wt w}-b+\sum_{k\ge3} c_k {\wt w}^k=0,\quad a=\frac{x-x_0}{c_2s}, \quad b=\frac{(s_0-s)}{c_2ss_0}
\end{equation}

By \Cref{lem:cubic-perturbation}, we have 
\begin{align}\label{e:h-tangent}
\frac{|s-s_0|}{\delta^{2/3}}\lesssim \min\left\{\frac{|s-s_0|}{|x-x_0|}, |s-s_0|^{1/3}\right\}\lesssim |{\wt w}_c|\lesssim |x-x_0|^{1/2}+|s-s_0|^{1/3}\lesssim  \delta^{1/3}.
\end{align}
In particular $|{\wt w}_c| \gtrsim |s-s_0|$.

We recall the expansion of $\wt S({\wt w};x_0, s_0)$ from \eqref{e:hor_cusp}. Then \eqref{e:hor_c_2} implies
\begin{align}
d:=\frac{1}{6}\wt S'''(0;x_0,s_0),\quad 
|d|=\frac{1}{3}|c_2 s^2_0|\asymp 1. 
\end{align}

The estimate \eqref{e:h-tangent}, together with \eqref{e:Scritical2}, we conclude that for $|{\wt w}|=r$,
\begin{align}
\wt S({\wt w}; x,s) - \wt S({\wt w}_c; x,s)=d{\wt w} ^3+\OO(\delta+r^{4}+\delta  \ln(1/\delta )),\quad d:=\wt S'''({\wt w}_0;x_0, s_0)/6.
\end{align}
Taking $r=\fc$ and then choosing $\fc$ and $\delta$ sufficiently small gives
\eqref{e:hor_cusp_turning_err}.
\end{proof}

%The estimate \eqref{e:h-tangent}, together with \Cref{c:S_diff}, we conclude that the statements in \Cref{prop:exit} hold with $(m, \delta, r)$ taken to be $(2, \delta\ln (1/\delta)+\fc^3, \fc)$.  Moreover, in this case $d\in\mathbb{R}\setminus\{0\}$ (recall from \eqref{eq:Taylor}). 

\begin{figure}
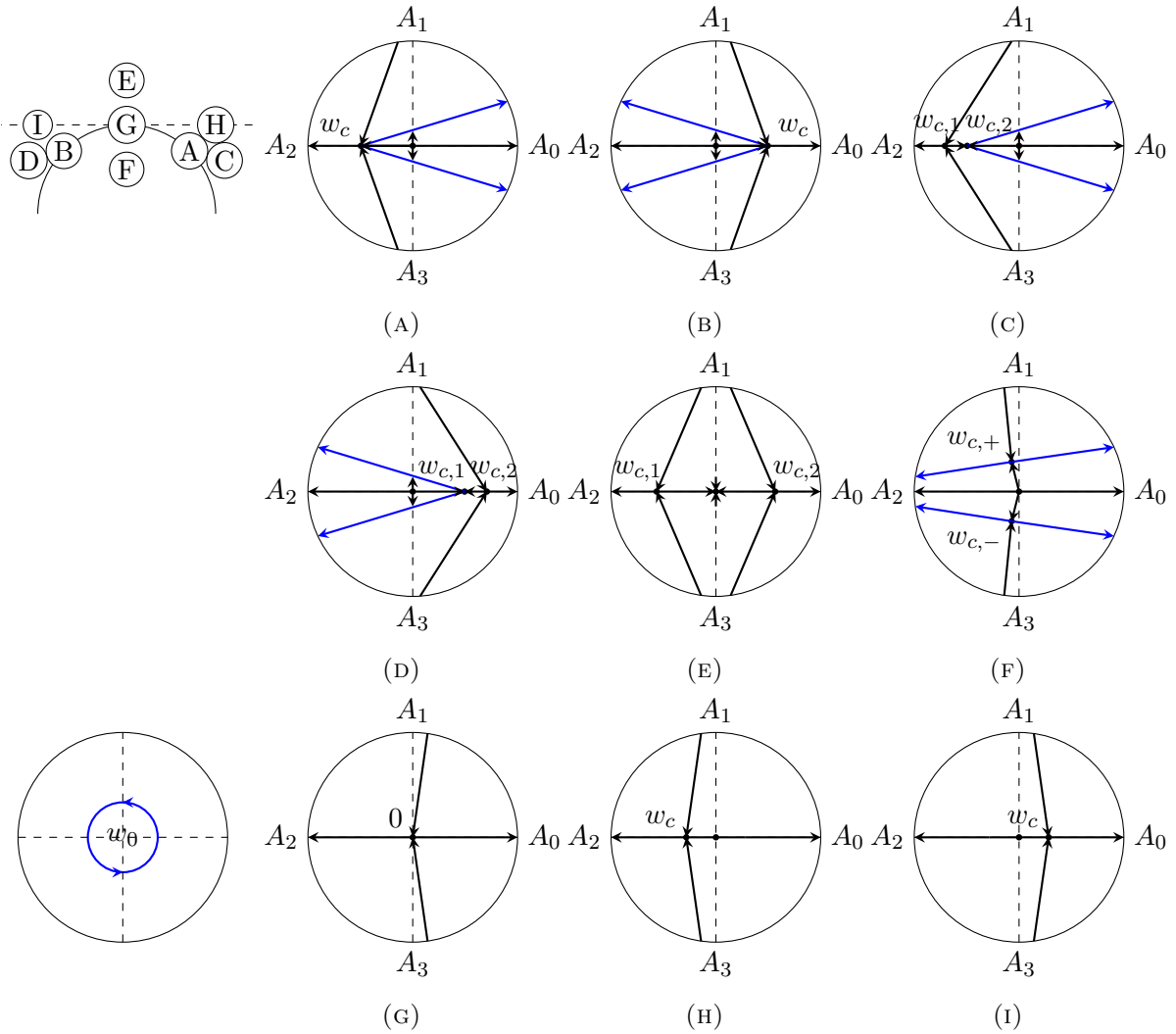

			\begin{subfigure}[t]{0.2\textwidth}
	
			\centering
		

				
			\caption{}
			\end{subfigure}
			\begin{subfigure}[t]{0.24\textwidth}
			
			%
				
			\caption{}
			\end{subfigure}
			\begin{subfigure}[t]{0.24\textwidth}
			
			%
				
			\caption{}
			\end{subfigure}
			
			\begin{subfigure}[t]{0.2\textwidth}
			\phantom{-}
			\end{subfigure}
			\begin{subfigure}[t]{0.24\textwidth}
			
			%
				
			\caption{}
			\end{subfigure}
			\begin{subfigure}[t]{0.24\textwidth}
			
			%
			\caption{}
			\end{subfigure}

			\begin{subfigure}[t]{0.2\textwidth}
			\centering
			%
				
		\caption{}
			\end{subfigure}

				\caption{
		Gradient-flow structure near a critical point $w_0$ corresponding to horizontal tangent region, with \(d=S''(w_0;x_0,s_0)/2<0\).}
	\label{f:h_tangent}
\end{figure}

\begin{proof}[Proof of the first statement in \Cref{c:horizontal_tangent_steepest}]

By \Cref{c:tangent_critical3}, uniformly for $|\wt w|=\fc$,
\begin{align}
|\wt S(\wt w;x,s)-\wt S(\wt w_c;x,s)-
d\wt w^2|
\leq |d|\fc^2/100,
\quad
d:=\wt S''(0;x_0,s_0)/2.
\end{align}
with $|d|\asymp 1$. The hypotheses of
\Cref{prop:exit} hold with $(m, \delta, r)$ taken to be $(2, |d|\fc/100, \fc)$.

From \eqref{e:S_copy}, we also have
\begin{align}\begin{split}\label{e:h_tangent_sdiff}
\wt S(\wt w; x,s) &= \wt S(\wt w; x_0,s_0)+(s\ln s-(x-x_0+1/\wt w)\ln (x-x_0+1/\wt w)\\
&-(s-x+x_0-1/\wt w)\ln (s-x+x_0-1/\wt w))\\
&-(s_0\ln s_0-(1/\wt w)\ln (1/\wt w)-(s_0-1/\wt w)\ln (s_0-1/\wt w)).
\end{split}\end{align}
Notice that from the above expression, the critical point $0$ is a singular point. For $|\wt w| \gtrsim |s-s_0|$, we have the following estimates for the derivatives of \eqref{e:h_tangent_sdiff}
%\begin{align}
%S'(w; x,s) &= S'(w; x_0,s_0)+\ln \frac{x-w}{w-(x-s)}-\ln \frac{x_0-w}{w-(x_0-s_0)} 
%\end{align}
\begin{align}\begin{split}\label{e:tangent_h_der}
\del_{\wt w}^k \wt S(\wt w;x,s)
&=\del_{\wt w}^k \wt S(\wt w;x_0,s_0)+\del_{\wt w}^{k-1} \frac{1}{\wt w^2} \ln\left(\frac{(1+(x-x_0)\wt w)(1-s_0 \wt w)}{1+(x-s-x_0)\wt w}\right)\\
&=\del_{\wt w}^{k-1} \frac{s-s_0}{\wt w}+\OO(1)\lesssim \frac{|s-s_0|}{|\wt w|^k}+1,
\end{split}\end{align}
where we used that from \eqref{e:hor_tangent}, $\wt S(\wt w;x_0, s_0)$ is holomorphic in a small neighborhood of $0$

We recall from \eqref{e:tangent_h_der} when $k=1$, $\wt S'(\wt w;x,s)$ contains the following non-analytic term
\begin{align}
\frac{s-s_0}{\wt w}.
\end{align}
If $s>s_0$,  $0$ is a sink:  the steepest–descent paths point inward; 
If $s<s_0$  $0$ is a source, the steepest–descent paths point outward; and if $s=s_0$ it vanishes and $\wt S(\wt w;x,s)$ is analytic locally around $0$. 

By \eqref{e:limitPI}  if $s\geq s_0$ the integrand in the first integral of \eqref{e:twoint} behaves like $\OO(w^{-n(s-s_0)-2})$ as $w\rightarrow \infty$. Thus it has no poles and is holomorphic at $\infty$. After the change of variable \eqref{e:change_variable}, we have that the integrand is holomorphic in a neighborhood of $0$.

We begin with \(d=\wt S''(0;x_0,s_0)<0\). This case corresponds to Panel (C) and (D) in \Cref{f:horizontal_tangent}, and $\mathsf C^{\rm d}(0)$ is a contour surrounding $w_0$.

When $s\neq s_0$ so $(x,s)$ is not on the extended-side, the configurations are the same as in the vertical tangent case, by shrinking the interval $[x_0, x]$ to a point $0$, we refer to Panel (A)--(F) in \Cref{f:h_tangent}. In these cases, we can deform $\wt \sfC^{\rm d}(0)$ to the union of $\wt \sfD^{\rm d}(\wt w_c)$. 

In the following we discuss the cases when $s=s_0$ so $(x,s)$ is on the extended-side

\begin{enumerate}

\item \textbf{Tangent Location}
If $(x,s) = (x_0,s_0)$, there is a single critical point $w_c =0$; see panel~(G) of \Cref{f:h_tangent}. Recall from \eqref{e:hor_tangent_copy} that, at $0$, the expansion of $\wt S(\cdot; x,s)$ vanishes up to the quadratic term. Hence, from $0$ there emanate two steepest–descent and two steepest–ascent branches, alternating in angle. By Proposition~\ref{prop:exit}, these four trajectories intersect the circle $|\wt w| = \fc$ in the arcs $A^{(4)}_0(0,\fc), A^{(4)}_1(0,\fc), A^{(4)}_2(0,\fc), A^{(4)}_3(0,\fc)$ in alternating order. In particular, when $d = \wt S''(0; x_0,s_0) < 0$, the two descent branches lie along the real axis from $0$ to $ \pm \fc$, and $\wt \sfD^{\rm d}(\wt w_c)=\emptyset$.

After the change of variable \eqref{e:change_variable}, when $s\geq s_0$, the integrand of the first integral in \eqref{e:twoint} is holomorphic in a small neighborhood of $0$, so the circular contour 
  \(\wt\sfC^{\rm d}(0)\) deforms to the empty set, which is also $\wt \sfD^{\rm d}(\wt w_c)=\emptyset$.

\item \textbf{Extended-side}
If $s=s_0$ and $x\neq x_0$, $\wt S(\wt w;x,s)$ is analytic locally around $0$. Thus $0$ is a spurious critical point, and there is a critical point  $\wt w_c\neq 0$. If $x>x_0$, we have $\wt w_c<0$; see Panel (H) of \Cref{f:h_tangent};  If $x<x_0$, we have $\wt w_c>0$; see Panel (I) of \Cref{f:h_tangent}.  

In both cases, at $\wt w_c$, the expansion of $\wt S(\cdot; x,s)$ vanishes up to the quadratic term. Hence, from $\wt w_c$ emanate two steepest–descent and two steepest–ascent branches, alternating in angle. By Proposition~\ref{prop:exit}, we have that these four trajectories intersect the circle $|\wt w| = \fc$ in the arcs $A^{(4)}_0(0,\fc), A^{(4)}_1(0,\fc), A^{(4)}_2(0,\fc), A^{(4)}_3(0,\fc)$ in alternating order. In particular, when $d = \wt S''(0; x_0,s_0) < 0$, the two descent branches lie along the real axis from $\wt w_c$ to $\pm \fc$, $\wt w_c$ is an ascent critical point, and $\wt\sfD^{\rm d}(\wt w_c)=\emptyset$.

The same as in the first case, the circular contour 
  \(\wt \sfC^{\rm d}(0)\) deforms to the empty set, which is also $\wt\sfD^{\rm d}(\wt w_c)=\emptyset$.
\end{enumerate}

Assume that \(d=S''(0;x_0,s_0)>0\). In the following we also only discuss the cases when $s=s_0$, so $(x,s)$ is on the extended-side. 

The steepest--descent/ascent configurations are the same as before, but all
gradient directions are reversed; see \Cref{f:h_tangent2}.  

In both cases as in Panel (A) and (B) of \Cref{f:horizontal_tangent}, \(\wt \sfC^{\rm d}(0)\) consists of two pieces: a path from
\(0-\ri\fc\) to \(0+\ri\fc\) passing to the left of \(0\), and a circle
surrounding \(0\). When $s\geq s_0$,  the integrand of the first integral in \eqref{e:twoint} is holomorphic in a small neighborhood of $0$, so the circular piece of
  \(\wt \sfC^{\rm d}(0)\) deforms to the empty set.

In the configuration shown in panel~(A), where $(x,s)=(x_0, s_0)$, and $0$ is a descent critical point. We can deform the 
path from \(0-\ri\fc\) to \(0+\ri\fc\) in \(\wt \sfC^{\rm d}(0)\) to $\wt\sfD^{\rm d}(0)$.

In the configurations shown in panel~(B)--(C), where $s=s_0$ but $x\neq x_0$. $0$ is a spurious critical point, and there exists an descent critical point $\wt w_c\neq 0$. We can deform the 
path from \(0-\ri\fc\) to \(0+\ri\fc\) in \(\wt \sfC^{\rm d}(0)\) to $\wt \sfD^{\rm d}(\wt w_c)$.
\end{proof}

\begin{figure}
			\begin{subfigure}[t]{0.2\textwidth}
	
			\centering
		
			%
				
		\caption{}
			\end{subfigure}

				\caption{
		Gradient-flow structure near a critical point $w_0$ corresponding to horizontal tangent region, with \(d=S''(w_0;x_0,s_0)/2<0\).}
	\label{f:h_tangent2}
\end{figure}

\begin{figure}
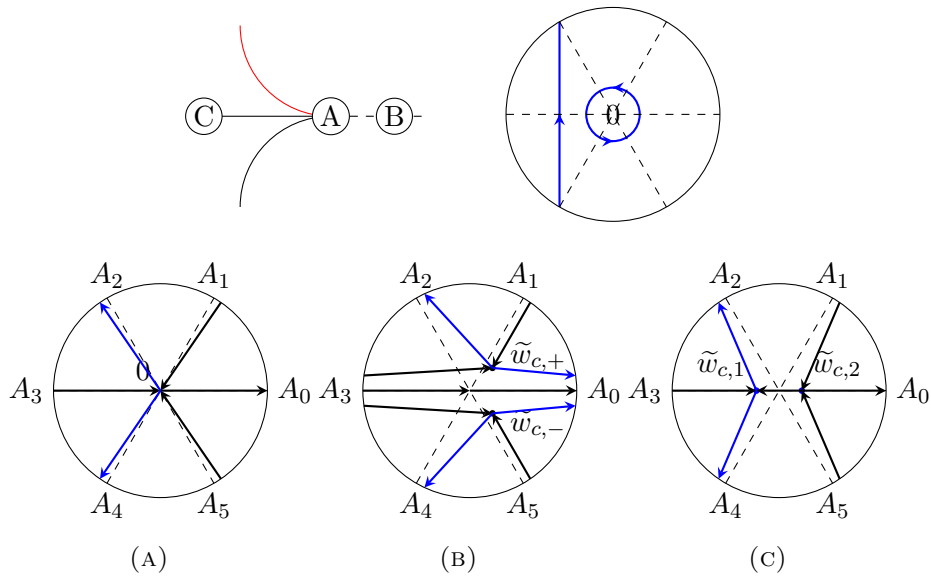

			\begin{subfigure}[t]{0.24\textwidth}
	
			\centering
		
			%
		
		\caption{}
	\end{subfigure}

				\caption{
		Gradient-flow structure near a critical point $0$ corresponding to horizontal cusp location, with \(d=\wt S'''(0;x_0,s_0)/6<0\).}
	\label{f:h_cusp}
\end{figure}

\begin{proof}[Proof of the second statement in \Cref{c:horizontal_tangent_steepest}]
By \Cref{c:cusp_turning_critical3}, uniformly for $|\wt w|=\fc$,
\begin{align}
|\wt S(\wt w;x,s)-\wt S(\wt w_c;x,s)-
d\wt w^3|
\leq |d|\fc^3/100,
\quad
d:=\wt S'''(0;x_0,s_0)/6.
\end{align}
with $|d|\asymp 1$. The hypotheses of
\Cref{prop:exit} hold with $(m, \delta, r)$ taken to be $(3, |d|\fc^3/100, \fc)$.

First suppose that $s\neq s_0$. Then $(x,s)$ does not lie on the
extended side. In this case the local configurations are the same as in
the vertical tangent case, with the interval $[x_0,x]$ collapsed, in the
local coordinate, to the point $0$. Thus we may use the six
configurations shown in Panels (A)--(F) of \Cref{f:vertical_cusp1} and  \Cref{f:vertical_cusp2}. The same
contour-deformation argument as in the proof of \Cref{l:vertical_tangent_steepest}
then shows that $\wt \sfC^{\rm d}(0)$ can be deformed to the union of the
corresponding local steepest-descent sets $\wt \sfD^{\rm d}(w_c)$, up to the
same exponentially small contributions from the auxiliary arcs.
This proves the claim when $s\neq s_0$.

It remains to consider the case $s=s_0$, so that $(x,s)$ lies on the
extended side. In the configurations of Panels (A)--(C) of
\Cref{f:h_cusp}, the point $0$ is a spurious critical point of multiplicity one. 

After
the change of variables \eqref{e:change_variable}, the integrand of the
first integral in \eqref{e:twoint} is holomorphic in a small
neighborhood of $0$ at $s=s_0$; indeed, the same holomorphicity holds in
the regime $s\geq s_0$. 

In the following we begin with \(d=\wt S'''(0;x_0,s_0)/6<0\).

\begin{enumerate}

\item \textbf{Cusp-Turning Location}
If $(x,s) = (x_0,s_0)$, there is a single critical point $\wt w_c = \wt w_0 =0$; see panel~(A) of \Cref{f:h_cusp}. Recall from \eqref{e:hor_cusp} that, at $0$, the expansion of $\wt S(\cdot; x,s)$ vanishes up to the cubic term. Hence, from $0$ there emanate three steepest–descent and three steepest–ascent branches, alternating in angle. By Proposition~\ref{prop:exit}, these six trajectories intersect the circle $|\wt w| = \fc$ in the arcs $A^{(6)}_0(0,\fc), A^{(6)}_1(0,\fc), \cdots, A^{(6)}_5(0,\fc)$ in alternating order. In particular, when $d = \wt S'''(0; x_0,s_0)/6 < 0$, the one ascent and one descent branches lie along the real axis from $-\fc$ to $0$ and $0$ to $\fc$. 

By Cauchy’s theorem, we may deform $\wt\sfC^{\rm d}(0)$ to $\wt \sfD^{\rm d}(w_c)$ together with (sub)arcs of $A^{(6)}_2(0;\fc)$ and $A^{(6)}_4(0;\fc)$. Thanks to \eqref{e:final_copy}, for any $\wt w$ on $A^{(6)}_2(0;\fc)$ and $A^{(6)}_4(0;\fc)$ we have
\begin{align}\label{e:h_cusp_diffS}
n\,\Re\!\big(\wt S(\wt w;x,s)-\wt S(\wt w_c;x,s)\big)
\leq -\,n\Bigl(|d|\,\fc^{3}-3|d|\fc^3/100\Bigr)
\le -\,n\fc'.
\end{align}

\item \textbf{Liquid region}
 If $s=s_0$,  there are only a conjugate pair critical points \(\wt w_{c,\pm}\); see panel (B) of \Cref{f:h_cusp}. Along the real axis, the flow runs from \(0-\fc\) to
\(0+\fc\). By real-analyticity, the configuration is symmetric with respect
to the real axis, so it suffices to describe the upper half-plane. From
\(w_{c,+}\) emanate two descent and two ascent branches, forming the usual
alternating pattern; these four branches stay in the upper half plane (flows cannot cross the real
axis). By Proposition~\ref{prop:exit}, they
meet \(|\wt w|=\fc\) through \(A^{(6)}_0(0, \fc), A^{(6)}_1(0, \fc),A^{(6)}_2(0, \fc),A^{(6)}_3(0, \fc)\) in alternating order. The flows from \(\wt w_{c,-}\) are the complex conjugates
of those from \(\wt w_{c,+}\). 

By Cauchy’s theorem, we may deform $\wt\sfC^{\rm d}(0)$ to $\wt \sfD^{\rm d}(\wt w_{c,+})\cup \wt \sfD^{\rm d}(\wt w_{c,-})$ together with (sub)arcs of $A^{(6)}_2(0;\fc)$ and $A^{(6)}_4(0;\fc)$, and \eqref{e:h_cusp_diffS} holds.

\item \textbf{Frozen region}
 If $s=s_0$, there are only two real critical points $\wt w_{c,1}<0<\wt  w_{c,2}$; see panel (C) of \Cref{f:h_cusp}. Each has two descent and two ascent branches, forming the
usual cross pattern of a simple saddle. Along the real axis, 
the descent branches goes from $-\fc$ to $\wt w_{c,1}$ and  $\wt w_{c,2}$ to $\wt w_{c,1}$ and $\fc$. The remaining branches are off the real axis. Since gradient flows
do not intersect, the two non-real descent branches from \(\wt  w_{c,1}\) exit the
circle through arcs \(A^{(6)}_2(0, \fc)\) and \(A^{(6)}_4(0, \fc)\); the two non-real ascent branches from \(\wt w_{c,2}\) exit the
circle through arcs \(A^{(6)}_1(0, \fc)\) and \(A^{(6)}_5(0, \fc)\), as shown in
panel (C) of \Cref{f:h_cusp}.

In this case $\wt w_{c,2}$ is an ascent critical point, and $\wt \sfD^{\rm d}(\wt w_{c,2})=\emptyset$. By Cauchy’s theorem, we may deform $\wt \sfC^{\rm d}(0)$ to $\wt \sfD^{\rm d}(\wt w_{c,1})$ together with (sub)arcs of $A^{(6)}_2(0;\fc)$ and $A^{(6)}_4(0;\fc)$, and \eqref{e:h_cusp_diffS} holds.

\end{enumerate}

The case  \(d =\wt S'''(0; x_0, s_0)/6 > 0\), the steepest–descent/ascent configurations are the same as before, but all gradient directions are reversed. Hence the same argument proves the desired contour deformation in the case $d>0$ as well. This completes the proof.

\end{proof}

\subsection{Horizontal tangent frozen chart}
\label{s:horizontal_frozen_chart_proof}
\begin{proof}[Proof of \Cref{c:horizontal_frozen_critical1}]

The tangent line through $(x_0, s_0)$ is horizontal. We denote the tangent location as $(x_0', s_0)$. As in
\eqref{e:frozen_change_coordinate}, we change coordinates by setting
\begin{align}\label{e:change_variable}
\wt w=\frac{1}{x'_0-w},
\qquad
\wt S(\wt w;x_0,s_0):=S(w;x_0,s_0),
\qquad
\wt\chi(\wt w):=\chi(w)
=
\frac{1}{s_0\wt w}+\sum_{i\geq 1}c_i\wt w^i .
\end{align}

Fix $\|(x,s)-(x_0,s_0)\|_2\leq \delta$. We study $\wt S(\wt w;x,s)$ in a small
neighborhood of $0$. From \eqref{e:critical_point2}, the critical points of
$\wt S(\wt w;x,s)$ satisfy
\begin{align}\label{e:x_exp3}
0
=
s\wt w\wt\chi(\wt w)-(x-x_0')\wt w-1
=
\frac{s}{s_0}
+\sum_{i\geq 1}sc_i\wt w^{i+1}
-(x-x_0')\wt w
-1 .
\end{align}
Equivalently,
\begin{align}
-\wt w-\frac{s-s_0}{(x_0'-x)s_0}
+\sum_{i\geq 1}\frac{sc_i}{(x-x_0')}\wt w^{i+1}
=0.
\end{align}
Thus the critical point equation has the form
\begin{equation}
F(\wt w;a,b):=-\wt w-b+\sum_{k\ge3} c_k \wt w^k=0,\quad b=\frac{s-s_0}{(x_0'-x)s_0}
\end{equation}
By \Cref{lem:bulk-perturbation}, we have 
\begin{align}\label{e:wc_critical_distance}
\wt w_c=-b+\OO(b^2)\asymp |s-s_0|.
\end{align}

We recall the expansion of $\wt S(\wt w;x_0, s_0)$ from \eqref{e:h_S_extend},
\begin{align}
d:= \wt S'(0;x_0,s_0)=s_0(x_0'-x_0),\quad 
|d|\asymp 1. 
\end{align}
The estimate \eqref{e:wc_critical_distance}, together with \eqref{e:Scritical2}, we conclude that for $|w-w_0|=r$,
\begin{align}
\wt S(\wt w; x,s) - \wt S(\wt w_c; x,s)=d\wt w+\OO(\delta+r^{2}+\delta  \ln(1/\delta )),\quad d:=\wt S'(0;x_0, s_0).
\end{align}
Taking $r=\fc$ and then choosing $\fc$ and $\delta$ sufficiently small gives
\eqref{e:hor_tangent_frozen_err}.

\end{proof}

\begin{figure}
		\begin{subfigure}[t]{0.23\textwidth}

				
			\caption{}
			\end{subfigure}
		\begin{subfigure}[t]{0.24\textwidth}
			
			%
				
			\caption{}
			\end{subfigure}
	
	\caption{
	Gradient flow structure near a nondegenerate saddle (\(m=2\) in
\Cref{prop:exit}), which is correspond to a vertical tangent location.}
	\label{f:side2}
\end{figure}

\begin{proof}[Proof of \Cref{l:local_descent_deformation3}]
By \Cref{c:horizontal_frozen_critical1}, uniformly for $|\wt w|=\fc$,
\begin{align}
|\wt S(\wt w;x,s)-\wt S(\wt w_c;x,s)-
d\wt w|
\leq |d|\fc/100,
\quad
d:=\wt S'(0;x_0,s_0).
\end{align}
with $|d|\asymp 1$. The hypotheses of
\Cref{prop:exit} hold with $(m, \delta, r)$ taken to be $(1, |d|\fc/100, \fc)$.

From \eqref{e:S_copy}, we also have
\begin{align}\begin{split}\label{e:h_tangent_sdiff}
&\wt S(\wt w; x,s) = \wt S(\wt w; x_0,s_0)-(s_0\ln s_0-(x_0-x'_0+1/\wt w)\ln (x_0-x'_0+1/\wt w)\\
&\phantom{\wt S(\wt w; x,s) = \wt S(\wt w; x_0,s_0)-(s_0\ln s_0}-(s_0-x_0+x'_0-1/\wt w)\ln (s_0-x_0+x'_0-1/\wt w))\\
&+(s\ln s-(x-x'_0+1/\wt w)\ln (x-x'_0+1/\wt w)-(s-x+x'_0-1/\wt w)\ln (s-x+x'_0-1/\wt w))
\end{split}\end{align}
Notice that from the above expression, the critical point $0$ is a singular point. For $|\wt w| \gtrsim |s-s_0|$, we have the following estimates for the derivatives of \eqref{e:h_tangent_sdiff}
\begin{align}\begin{split}\label{e:tangent_h_der}
\del_{\wt w}^k \wt S(\wt w;x,s)
&=\del_{\wt w}^k \wt S(\wt w;x_0,s_0)+\del_{\wt w}^{k-1} \frac{1}{\wt w^2} \ln\left(\frac{(1+(x-x'_0)\wt w)(1+(x_0-s_0-x_0') \wt w)}{(1+(x_0-x_0')\wt w)(1+(x-s-x'_0)\wt w)}\right)\\
&=\del_{\wt w}^{k-1} \frac{s-s_0}{\wt w}+\OO(1)\lesssim \frac{|s-s_0|}{|\wt w|^k}+1,
\end{split}\end{align}
where we used that from \eqref{e:h_S_extend}, $\wt S(\wt w;x_0, s_0)$ is holomorphic in a small neighborhood of $0$

\noindent\textbf{Replacing \(\wt \sfC^{\rm d}(w_0)\) by \(\wt \sfD^{\rm d}(w_c)\).}

First suppose that $s\neq s_0$. Then $(x,s)$ does not lie on the
extended side. In this case the local configurations are the same as in
the vertical tangent case, with the interval $[x_0,x]$ collapsed, in the
local coordinate, to the point $0$. Thus we may use the four
configurations shown in Panels (A)--(D) of \Cref{f:side}. The same
contour-deformation argument as in the proof of \Cref{l:local_descent_deformation}
then shows that $\wt \sfC^{\rm d}(w_0)$ can be deformed to the union of the
corresponding local steepest-descent sets $\wt \sfD^{\rm d}(w_c)$, up to the
same exponentially small contributions from the auxiliary $A^{(2)}_0$-arcs.
This proves the claim when $s\neq s_0$.

It remains to consider the case $s=s_0$, so that $(x,s)$ lies on the
extended side. In the configurations of Panels (A) and (B) of
\Cref{f:side2}, the point $0$ is a spurious critical point. There is no critical point.

We now show that the original contour also gives no contribution. By \eqref{e:limitPI}  if $s\geq s_0$ the integrand in the first integral of \eqref{e:twoint} behaves like $\OO(w^{-n(s-s_0)-2})$ as $w\rightarrow \infty$.  After
the change of variables \eqref{e:change_variable}, the integrand of the
first integral in \eqref{e:twoint} is holomorphic in a small
neighborhood of $0$ at $s=s_0$; indeed, the same holomorphicity holds in
the regime $s\geq s_0$. Therefore there are only two possibilities:
either
$
    \wt \sfC^{\rm d}(w_0)=\emptyset$,
or $\wt \sfC^{\rm d}(w_0)$ is a small circular contour surrounding $0$. In
the first case there is nothing to prove. In the second case, since the
integrand is holomorphic in a neighborhood of $0$, Cauchy's theorem
allows us to deform the circular contour to the empty contour. This completes the proof.

\noindent\textbf{Replacing \(\wt \sfD^{\rm d}(w_c)\) by \(\wt \sfS^{\rm d}(w_c)\).}

We introduce the shifted coordinate
\begin{align}\label{e:shift_coordinate}
  \sfb = s - s_0, \quad
\end{align}

We discuss two cases separately
\begin{enumerate}

\item \text{Close to the tangent line} $|\sfb|\leq (\ln n)^3/n$.
Let
$r_n:=(\ln n)^5 /n$. 
We define the local steepest–descent set $\wt \sfS(w_c)$ at \(w_c\) to be the portions of
steepest–descent trajectories starting at \(w_c\) up to their first exit from
the disk $\{\wt w: |\wt w|\leq r_n\}$. 

On the circle $\{|\wt w|=r_n\}$, we have
%\begin{align}\begin{split}\label{e:tangent_S_exp}
%&\phantom{{}={}}\wt S(w;x,s)-\wt S(w_c;x,s)\\
%&=\wt S(w;x_0,s_0)-\wt S(w_c;x_0,s_0)+\int_{w_c}^w \ln\left(\frac{u-x}{x-s-u}\right)-\ln\left(\frac{u-x_0}{x_0-s_0-u}\right)\rd u\\
%&=d w^2 +\int_{w_c}^w \ln \frac{u-x}{u-x_0}\rd u+ \OO(|w_c-w_0|^2+|w-w_0|^3+|w_c-w_0|^3+|w-w_c|(|x-x_0|+|s-s_0|)),
%\end{split}\end{align}
%From \eqref{e:wS_der} and the Taylor expansion \eqref{e:hor_tangent}
\begin{align}\begin{split}\label{e:v_tangent_S_exp}
&\phantom{{}={}}\wt S(\wt w;x,s)-\wt S(\wt w_c;x,s)\\
&=\wt S(\wt w;x_0,s_0)-\wt S(\wt w_c;x_0,s_0)+\int_{\wt w_c}^{\wt w} \frac{1}{u^2}\ln\left(\frac{(1+(x-x'_0)u)(1+(x_0-s_0-x_0') u)}{(1+(x_0-x_0')u)(1+(x-s-x'_0)u)}\right)\rd u\\
&=dw +\int_{\wt w_c}^{\wt w} \frac{1}{u^2}\ln\left(\frac{(1+(x-x'_0)u)(1+(x_0-s_0-x_0') u)}{(1+(x_0-x_0')u)(1+(x-s-x'_0)u)}\right)\rd u+ \OO(|\wt w_c|+|\wt w|^2),
\end{split}\end{align}
where the first statement is from \eqref{e:tangent_h_der}, and the second statement is from \eqref{e:h_S_extend}. For the integral in \eqref{e:v_tangent_S_exp}, we can bound it as
\begin{align}\begin{split}\label{e:v_integral_bound}
&\phantom{{}={}}\int_{w_c}^w \frac{1}{u^2}\ln\left(\frac{(1+(x-x'_0)u)(1+(x_0-s_0-x_0') u)}{(1+(x_0-x_0')u)(1+(x-s-x'_0)u)}\right)\rd u\\
&=\int_{\wt w_c}^{\wt w} \frac{s-s_0}{u}+\OO(|s-s_0|+|x-x_0|)\rd u\\
&=(s-s_0)\ln\frac{\wt w}{\wt w_c}+\OO(r_n(|s-s_0|+|x-x_0|))
\lesssim \sfb\ln(1/|\sfb|)+ r_n=\OO((\ln n)^4/n),
\end{split}\end{align}
where the first statement is from direct Taylor expansion; the second statement is from integrating; the last two statements follow from  $|\wt w|=r_n\ll 1$ and $|\wt w_c|\asymp |s-s_0|=|\sfb|$ from \eqref{e:twc_bound}.

 By plugging \eqref{e:v_integral_bound} into \eqref{e:v_tangent_S_exp}, we conclude 
\begin{align}\begin{split}\label{e:v_tangent_S_exp2}
\wt S(\wt w;x,s)-\wt S(\wt w_c;x,s)=d\wt w +\OO((\ln n)^4/n).
\end{split}\end{align}

The estimate \eqref{e:v_tangent_S_exp2} verifies the assumptions in \Cref{prop:exit} with $(m, \delta, r)$ taken to be $(1, (\ln n)^4/n, r_n)$. 
It follows that the two steepest–descent branches from \(\wt w_c\) exit the circle
\(\{\wt w: |\wt w|=r_n\}\) at points on \(A^{(2)}_1(0;r_n)\).  At any such
point,
\begin{equation}\label{eq:sd_high_drop}
n\,\Re\!\big(\wt S(\wt w; x,s)-\wt S(\wt w_c; x,s)\big)
= -\,n|d|\,r_n+\OO((\ln n)^4)
\le -\,\fc'\,(\ln n)^{5}.
\end{equation}
Moreover, since $\Re[\wt S(\wt w; x,s)]$ decreases along the steepest–descent, \eqref{eq:sd_high_drop}
holds for all \(w\in\wt \sfD(w_c)\setminus \wt \sfS(w_c)\).
Consequently, for  \(w\in\wt  \sfD(w_c)\setminus \wt \sfS(w_c)\) we have
\[
e^{n\Re[\wt S(\wt w;x,s)]}
\leq
e^{n\Re[\wt S(\wt w_c;x,s)]}e^{-\fc'(\ln n)^5}.
\]

\item \textbf{Other Frozen} $|\sfb|\geq (\ln n)^3/n$.

Locally around $\wt w_c$, $\wt S(\wt w_c; (x,s))$ is holomorphic and we have the Taylor expansion 
\begin{align}\begin{split}\label{e:hor_side_S0}
\wt S(\wt w;x,s)-\wt S(\wt w_{c};x,s)
&=d (\wt w-\wt w_{c})^2 +\OO\left(\frac{|\sfb||\wt w-\wt w_c|^3}{|\wt w_c|^3}\right)
\end{split}\end{align}
where the error term is from \eqref{e:tangent_h_der}, and $d=\wt S''(\wt w_c; (x,s))/2$. 
Thanks to \eqref{e:dertS2}, \eqref{e:wtchi} and \eqref{e:twc_bound}
\begin{align}\begin{split}\label{e:hor_side_S1}
&|d|\asymp|\wt S''(\wt w_{c};x,s)|=\left|\frac{1}{{\wt w}_c^3}\,
\frac{\del_{\wt w}\bigl({\wt w}\wt \chi({\wt w})\bigr)\big|_{{\wt w}={\wt w}_c}-(x-x'_0)/s}{\wt \chi({\wt w}_c)\bigl(1-\wt \chi({\wt w}_c)\bigr)}\right|\asymp \frac{1}{|\wt w_c|}\asymp\frac{1}{|\sfb|},
\end{split}
\end{align}

By plugging \eqref{e:hor_side_S1} into \eqref{e:hor_side_S0}, we get
\begin{align}\begin{split}\label{e:hor_side_S2}
\wt S(\wt w;x,s)-\wt S(\wt w_{c};x,s)
&=\Theta\left(1/|\sfb|\right)(\wt w-\wt w_{c})^2+\OO\left(\frac{|\wt w-\wt w_c|^3}{|\sfb|^2}\right),
\end{split}\end{align}

Let
$r_n:=(\ln n)^2 (|\sfb|/n)^{1/2}$. 
We define the local steepest–descent set $\wt \sfS(w_c)$ at \(w_c\) to be the portions of
steepest–descent trajectories starting at \(w_c\) up to their first exit from
the disk $\{w: |w-w_c|\leq r_n\}$.  
For any $w=w_c+r_n e^{\ri\theta}$, we have
\begin{align}\begin{split}
&\phantom{{}={}}\wt S(\wt w;x,s)-\wt S(\wt w_c;x,s)=\Theta(1/|\sfb|)(\wt w-\wt w_c)^2+\OO(|\wt w-\wt w_c|^3/|\sfb|^2)\\
&=\Theta(1/|\sfb|)(\wt w-\wt w_c)^2+\OO((\ln n)^3/(|\sfb|^{1/2}n^{3/2}))
=\Theta(1/|\sfb|)(\wt w-\wt w_c)^2+\OO((\ln n)^{3/2}/n)
\end{split}\end{align}
The above estimate verifies the assumptions in \Cref{prop:exit} with $(m, \delta, r)$ taken to be $(2, (\ln n)^{3/2}/n, r_n)$. We can therefore replace $\wt \sfD^{\rm d}(w_c)$ by $\wt \sfS^{\rm d}(w_c)$ as in the
close to tangent case.
\end{enumerate}

\end{proof}

\subsection{Proof of perturbation estimates}
\label{s:perturb}

\begin{lemma}\label{e:power_series}
Let $g(z)=\sum_{k\geq m} c_k z^k$ be analytic on the open disk $\{z\in \bC:|z|<\fc_0\}$,
then there exists a small $0<\fc\leq \fc_0$ and large $M\leq 1/(100\fc)$, such that the following holds
\begin{enumerate}
\item For any $|z|\leq \fc$ and $0\leq \ell\leq m$, $| \del_z^{\ell}g(z)|\leq M |z|^{m-\ell}$.
\item  For every $|z'|\leq \fc$, we can Taylor expand $g(z)$ around $z'$
\begin{align}
g(z)=\sum_{k\geq 0} c_k'(z-z')^k,\quad c_k'=\frac{\del_z^k g(z')}{k!},
\end{align}
which converges absolutely, and 
\begin{align}
\left|\sum_{k\geq m} c_k'(z-z')^k\right|\leq M|z-z'|^m.
\end{align}
\end{enumerate}
\end{lemma}
\begin{proof}
We can write
\[
g(z)=z^{m}h(z),\quad h(z):=\sum_{k\geq m} c_k z^{k-m},
\]
so $h$ is analytic on $\{|z|<\fc_0\}$. For each $0\leq \ell\leq m$, differentiating termwise gives
\[
\del_z^{\ell} g(z)=\sum_{k\geq m} k(k-1)\cdots(k-\ell+1)\,c_k\, z^{k-\ell}
= z^{m-\ell}\,h_\ell(z),
\]
where
\[
h_\ell(z):=\sum_{k\ge m} k(k-1)\cdots(k-\ell+1)\,c_k\, z^{k-m}
\]
is analytic on $\{|z|<\fc_0\}$.

We set
\begin{align}\label{e:der_bound}
M:=\max_{0\leq \ell\leq m}\ \sup_{|z|\leq\fc_0/2} |h_\ell(z)|<\infty, \quad
\fc:=\min\Big\{\fc_0/10,\ \frac{1}{100M}\Big\},
\end{align}
with $1/(100M)=+\infty$ if $M=0$.
Then $\fc\leq \fc_0$ and $M\leq 1/(100\fc)$ as desired.

 For $|z|\leq \fc\leq \fc_0/2$ and $0\leq \ell\leq m$, we have
\[
|\del_z^{\ell} g(z)|=|z|^{m-\ell}\,|h_\ell(z)|\leq M\,|z|^{m-\ell}.
\]
This gives the first statement. 

 Fix $z'$ with $|z'|\leq \fc$. Since $g$ is analytic on $\{|z|<\fc_0\}$, the Taylor
expansion at $z'$,
\[
g(z)=\sum_{k\geq 0} c_k'(z-z')^k,\quad c_k'=\frac{\del_z^k g(z')}{k!},
\]
converges absolutely for $|z-z'|<\fc_0-|z'|$.
For the remainder of order $m-1$, use the integral form of Taylor’s theorem:
for any $z$ with $|z-z'|<\fc_0-|z'|$,
\[
\sum_{k\ge m} c_k'(z-z')^k
=\frac{(z-z')^m}{(m-1)!}\int_0^1 (1-t)^{m-1}\,\del_z^{m}g\big(z'+t(z-z')\big)\,\rd t .
\]
For $|z|, |z'|\leq \fc$, along the segment $z'+t(z-z')$ we have
\(
\big|z'+t(z-z')\big|\leq |z'|+t(|z|+|z'|)\leq 3\fc\leq \fc_0/2,
\)
so by \eqref{e:der_bound} with $\ell=m$,
\(
|\del_z^{m}g\big(z'+t(z-z')\big)|\leq M.
\)
Hence
\[
\Big|\sum_{k\geq m} c_k'(z-z')^k\Big|
\leq \frac{|z-z'|^m}{(m-1)!}\int_0^1 (1-t)^{m-1} M\,\rd t
= \frac{M}{m!}\,|z-z'|^m \leq M\,|z-z'|^m.
\]
\end{proof}

\begin{proof}[Proof of \Cref{lem:bulk-perturbation}]
The power series $g(z)=\sum_{k\ge 2} c_k z^k$ satisfies the assumptions of \Cref{e:power_series} with $m=2$. Thus there exist constants $M,\fc$ with $M\le 1/(100\fc)$ such that for any $|z|\le \fc$ and $0\le \ell\le 3$, 
\begin{align}\label{e:bulk_gbound}
|g(z)|\le M |z|^{2}.
\end{align}

Let $G(z)=-z$. Then, on the circle $|z|=\fc$, we have
\begin{align}\begin{split}
\big|F(z;b)-G(z)\big|
&\le |b|+M |z|^2\le |b|+M\fc^2\le \fc/2+\fc/4< \fc=|G(z)|,
\end{split}\end{align}
where the first inequality follows from \eqref{e:gbound}; the second uses $|z|=\fc$; and the third holds provided $\fc M\le 1/4$.
By Rouch\'e's theorem, $F(\,\cdot\,;b)$ and $G$ have the same number of zeros inside $|z|<\fc$. We conclude that $F$ has exactly one zero in $|z|<\fc$. We can repeat this argument on the circle $|z|=2|b|\leq \fc$ 
\begin{align}\begin{split}
\big|F(z;b)-G(z)\big|
&\le |b|+M |z|^2\le |b|+M(2|b|)^2\leq  2|b|=|G(z)|,
\end{split}\end{align}
where the first inequality follows from \eqref{e:bulk_gbound}; the second uses $|z|=2|b|$; and the third holds provided $\fc M\le 1/4$. We conclude that $F$ has exactly one zero in $|z|\leq 2|b|$.
The zero is given by
\begin{align}
z=-b+\OO\left(M|z|^2\right)=-b+\OO\left(|b|^2\right).
\end{align}

\end{proof}

\begin{proof}[Proof of \Cref{lem:edge-perturbation}]
The power series $g(z)=\sum_{k\ge 3} c_k z^k$ satisfies the assumptions of \eqref{e:power_series} with $m=3$. Thus there exist constants $M,\fc$ with $M\le 1/(100\fc)$ such that
\begin{enumerate}
\item For any $|z|\le \fc$ and $0\le \ell\le 2$, 
\begin{align}\label{e:gbound}
|\partial_z^{\ell}g(z)|\le M |z|^{3-\ell}.
\end{align}
\item For every $|z'|\le \fc$, 
\begin{align}\label{e:gbound2}
\left|\sum_{k\ge 3} c_k'(z-z')^k\right|\le M|z-z'|^{3},\qquad c_k'=\frac{\partial_z^k g(z')}{k!}.
\end{align}
\end{enumerate}

%We first prove \eqref{e:small_real_z0}. Let $G(z)=-az$. Then, on the circle $|z|=\fc$, we have
%\begin{align}\begin{split}\label{e:edge_F-Gdiff}
%\big|F(z;a,b)-G(z)\big|
%&\le |b|+|z|^2+M |z|^3\le |b|+\fc^2+ M\fc^3\le 3\fc^2< |a|\fc=|G(z)|,
%\end{split}\end{align}
%where the first inequality follows from \eqref{e:gbound} with $\ell=0$; the second uses $|z|=\fc$; and the third holds provided $\fc M\le 1$, $|a|\geq 4\fc$ and $|b|\leq \fc^2$.
%By Rouch\'e's theorem, $F(\,\cdot\,;a,b)$ and $G$ have the same number of zeros inside $|z|<\fc$. We conclude that $F$ has exactly one zero in $|z|<\fc$. We can repeat this argument on the circle $|z|=4|b|/|a|\leq \fc$ 
%\begin{align}\begin{split}\label{e:edge_F-Gdiff}
%\big|F(z;a,b)-G(z)\big|
%&\le |b|+|z|^2+M |z|^3\le |b|+16|b|^2/|a|^2+ M(|b|/|a|)^3\leq  4|b|=|G(z)|,
%\end{split}\end{align}
%where the first inequality follows from \eqref{e:gbound} with $\ell=0$; the second uses $|z|=4|b|/|a|$; and the third holds provided $\fc M\le 1$, $|b|\geq 16|b|^2/|a|^2$. We conclude that $F$ has exactly one zero in $|z|\leq 4|b|/|a|$.
%The zero is given by
%\begin{align}
%z=-\frac{b}{a}+\OO\left(\frac{|z|^2+M|z|^3}{|a|}\right)=-\frac{b}{a}+\OO\left(\frac{b^2}{|a|^3}\right).
%\end{align}
%
%
%
%Next we study the case that $\rho_0:=\max\{|a|/4,\,|b|^{1/2}\}\le \fc$. 

Let $G(z)=z^2$. Then, on the circle $|z|=8\fc$, we have
\begin{align}\begin{split}\label{e:edge_F-Gdiff}
\big|F(z;a,b)-G(z)\big|
&\le |a||z|+|b|+M |z|^3\le 8|a|\fc+|b|+M(8\fc)^3< 64\fc^2=|G(z)|,
\end{split}\end{align}
where the first inequality follows from \eqref{e:gbound} with $\ell=0$; the second uses $|z|=8\fc$; and the third holds provided $100\fc M\le 1$,  $|a|\leq 4\fc$ and $|b|\le \fc^2$.
By Rouch\'e's theorem, $F(\,\cdot\,;a,b)$ and $G$ have the same number of zeros (counting multiplicity) inside $|z|<8\fc$. Since $G$ has exactly two zeros at $z=0$ (counting multiplicity), $F$ also has exactly two zeros in $|z|<8\fc$. We can repeat this argument on the circle $|z|=4\rho_0$ (as in \eqref{e:edge_F-Gdiff}) and conclude that $F$ has exactly two zeros in $|z|<8\rho_0$ as well. Finally, since $F(z;a,b)$ is real-analytic, i.e., $F(\overline{z};a,b)=\overline{F(z;a,b)}$, complex roots occur in conjugate pairs.

If $|b|^{1/2}> |a|/4$, then, for any $|z|\leq |b|^{1/2}/8$,
\begin{align}\begin{split}
F(z;a,b)\ge |b|-|b|/64-|a| (|b|^{1/2}/8) -M(|b|^{1/2}/8)^{3}
\ge \frac{1}{2}\Big(|b|-2|a|(|b|^{1/2}/8)\Big)>0,
\end{split}\end{align}
where the first inequality follows by substituting $|z|\le |b|^{1/2}/2$; the second holds provided  $M|b|^{1/2}\le M\fc\le  1$; and the last inequality follows from $|b|^{1/2}>|a|/4$. In this setting we also have that both roots are bounded by $4\rho_0\asymp |b|^{1/2}$. Thus both roots satisfy:
\begin{align}\label{e:edge_bothroot} 
 |b|^{1/2}/8\leq |z|\leq 8\rho_0= |b|^{1/2}.
\end{align}

Next, we prove that
\begin{align}\label{e:zlower}
|z|\geq \frac18\min\{|b/a|,\rho_0\}.
\end{align}
The claim follows from \eqref{e:edge_bothroot} when
$
|b|^{1/2}>{|a|}/{4}.
$
Suppose instead that
$
|b|^{1/2}\leq {|a|}/{4}.
$
If
$
|z|
\leq
 \min\{|b/a|,\rho_0\}/8
\leq
{|b|}/({8|a|})
\leq
{|b|^{1/2}}/{32},
$
then
\begin{align}
\begin{split}
|F(z;a,b)|
&\geq
|b|-{|b|}/{8}
-\left({|b|^{1/2}}/{32}\right)^2
-M\left({|b|^{1/2}}/{32}\right)^3\geq
|b|(1/4-M\fc)>0,
\end{split}
\end{align}
where we used \(|b|\leq\fc^2\) and \(100\fc M\leq1\). Thus \eqref{e:zlower} holds 
as claimed.

A double root occurs when $F(z;a,b)=\partial_z F(z;a,b)=0$:
\begin{align}\label{e:edge_derF}
\partial_z F(z;a,b)=2z-a+\sum_{k\ge 3} k c_k z^{k-1}=0.
\end{align}
By Rouch\'e's theorem (using the same argument as in \eqref{e:edge_F-Gdiff}), the equation \eqref{e:edge_derF} has one solution
\begin{align}\label{e:defzc}
z'=\frac{a}{2}\big(1+\OO(|a|)\big),
\end{align}
where we bound the higher-order terms in \eqref{e:edge_derF} by $\OO(|z|^2)$ using \eqref{e:gbound}.

Thus $F(z;a,b)$ is decreasing on the interval $[-4\rho_0, z']$ and increasing on the interval $[z', 4\rho_0]$. Setting $F(z';a,b')=0$ yields
\begin{align}\label{e:defbc}
b':=b'(a):=(z')^2-a z'+\sum_{k\ge 3} c_k (z')^k
=-(z')^2 +\sum_{k\ge 3} (1-k)c_k (z')^k= -\frac{a^2}{4}\big(1+\OO(|a|)\big),
\end{align}
where in the last step we bound the higher-order terms by $\OO(|z'|^3)=\OO(|a|^3)$ using \eqref{e:gbound}.
We conclude that for $b\ge b'$ the two roots are real, one on each side of $z'$; denote them by $z_0\le z'\le z_1$. If $b<b'$, there is a complex conjugate pair $z_\pm$ of roots.

We can Taylor expand around $z'$, noting that $F(z';a,b')=\partial_z F(z';a,b')=0$:
\begin{align}\label{e:bbc}
b-b'=\frac{\partial_z^2 F (z';a,b)}{2}\,(z-z')^2+\OO\big(|z-z'|^3\big),
\end{align}
where the error term is bounded by $M|z-z'|^3$ by \eqref{e:gbound2} (with $z'=(a/2)(1+\OO(|a|))$), and, by \eqref{e:gbound},
\begin{align}
\partial_z^2 F (z';a,b)=2+\sum_{k\ge 3}k(k-1) c_k (z')^{k-2}=2\big(1+\OO(z')\big)
=2 \big(1+\OO(|a|)\big). 
\end{align}

If $|b-b'|\le a^{2}/(100M^2)$, then by Rouch\'e's theorem (the same argument as in \eqref{e:edge_F-Gdiff}), \eqref{e:bbc} has exactly two solutions in the disk $|z-z'|\le 2\sqrt{|b-b'|}$. They satisfy
\begin{align}\label{e:edge_imz}
z-z'=\pm \big(1+\OO(|a|)\big)\sqrt{b-b'}.
\end{align}
Thus, for $b\ge b'$ there are three cases:
\begin{enumerate}
\item If $|b-b'|\le a^{2}/(100M^2)$, the claim \eqref{e:edge_realroot} follows from \eqref{e:edge_imz}. 
\item If $b-b'\asymp a^{2}$, note that as $b$ increases, the distances $z_1-z'$ and $z'-z_0$ also increase. We conclude that $z_1-z',\, z'-z_0\gtrsim |a|\asymp \sqrt{b-b'}$. In this setting $\rho_0\asymp |a|$, so we also have $z_1-z',\, z'-z_0\lesssim \rho_0\asymp |a|\asymp \sqrt{b-b'}$, and the claim \eqref{e:edge_realroot} follows. 
\item If $b-b'\ge (10a)^{2}$, we recall from \eqref{e:defzc} and \eqref{e:defbc}, $z'=(a/2)(1+\OO(|a|))$ and $b'=-(a^2/4)(1+\OO(|a|))$. 
Then $b\ge (8a)^2$ and $|b|\asymp \sqrt{b-b'}$. We also recall from \eqref{e:edge_bothroot} that in this case  $|a|\leq |b|^{1/2}/8\leq |z_0|, |z_1|\leq 8 |b|^{1/2}$. We then conclude that $z_1-z',\, z'-z_0\asymp |b|^{1/2}\asymp \sqrt{b-b'}$. 
\end{enumerate}

For $b<b'$ there are also two cases:
\begin{enumerate}
\item If $|b-b'|\le a^{2}/(100M^2)$, \eqref{e:edge_imz} implies that $|z_\pm|\asymp |a|$ and $|z_\pm-z'|\asymp |\Im z_\pm|\asymp \sqrt{b'-b}$. 

\item If $b'-b\gtrsim a^2$, write $z_+ =x+\ri y$ with $y>0$. Taking the real and imaginary parts of $F(z;a,b)=0$ gives
\begin{align}\label{e:edge_emF}
0=\Im[F(x+\ri y;a,b)]=2xy-ay+y\sum_{k\ge 3}\frac{c_k \Im[(x+\ri y)^k]}{y}.
\end{align}
We note that $\Im[(x+\ri y)^k]\le k |x+\ri y|^{k-1}y$, and, after dividing \eqref{e:edge_emF} by $y$ and using \eqref{e:gbound}, we obtain
\begin{align}\label{e:edge_Imbound}
x=\frac{a}{2}+\OO\big(|x+\ri y|^2\big),\qquad |x|\lesssim |a|+y.
\end{align}
Similarly, taking the real part yields
\begin{align}\label{e:edge_emrelation}
0=\Re[F(x+\ri y;a,b)]=-b-ax +(x^2-y^2) +\OO\big(|x+\ri y|^3\big).
\end{align}
Combining \eqref{e:edge_Imbound} and \eqref{e:edge_emrelation}, we get
\begin{align}\begin{split}\label{e:edge_ylowb}
b&=-ax+(x^2-y^2)+\OO\big(|x+\ri y|^3\big)\\
&=-\frac{a^2}{2} +\frac{a^2}{4}-y^2+\OO(|a|^3+y^3)
=-\frac{a^2}{4} -y^2+\OO(|a|^3+y^3).
\end{split}\end{align}
Recalling from \eqref{e:defbc} that $b'=-(a^2/4)\big(1+\OO(|a|)\big)$, if $b'-b\gtrsim a^2$, then \eqref{e:edge_ylowb} gives
\begin{align}
|\Im z_\pm|^2=y^2=-\frac{a^2}{4}-b-\OO(|a|^3+y^3)
=b'-b -\OO(|a|^3+y^3)\asymp b'-b,
\end{align}
and the claim \eqref{e:edge_complexroot} follows.
\end{enumerate}
\end{proof}

\begin{proof}[Proof of \Cref{lem:cubic-perturbation}]
The power series $g(z)=\sum_{k\ge 4} c_k z^k$ satisfies the assumptions of \eqref{e:power_series} with $m=4$. Thus there exist constants $M,\fc$ with $M\le 1/(100\fc)$ such that
\begin{enumerate}
\item For any $|z|\le \fc$ and $0\le \ell\le 4$, 
\begin{align}\label{e:cusp_gbound}
|\partial_z^{\ell}g(z)|\le M |z|^{4-\ell}.
\end{align}
\item For every $|z'|\le \fc$, 
\begin{align}\label{e:cusp_gbound2}
\left|\sum_{k\ge 4} c_k'(z-z')^k\right|\le M|z-z'|^{4},\qquad c_k'=\frac{\partial_z^k g(z')}{k!}.
\end{align}
\end{enumerate}

Let $G(z)=z^3$. Then, on the circle $|z|=\fc$, we have
\begin{align}\begin{split}\label{e:cusp_F-Gdiff}
\big|F(z;a,b)-G(z)\big|
&\le |a||z|+|b|+M |z|^4\le |a|\fc+|b|+M\fc^4\le \fc^3/2=|G(z)|/2,
\end{split}\end{align}
where the first inequality follows from \eqref{e:cusp_gbound} with $\ell=0$; the second uses $|z|=\fc$; and the third holds provided $\fc M\le 1/8$ and $\rho_0\le \rho\le \fc/4$.
By Rouch\'e's theorem, $F(\,\cdot\,;a,b)$ and $G$ have the same number of zeros (counting multiplicity) inside $|z|<\fc$. Since $G$ has exactly three zeros at $z=0$ (counting multiplicity), $F$ also has exactly three zeros in $|z|<\fc$. We can repeat this argument on the circle $|z|=4\rho_0$ (as in \eqref{e:cusp_F-Gdiff}) and conclude that $F$ has exactly three zeros in $|z|<4\rho_0$ as well. Finally, since $F(z;a,b)$ is real-analytic, i.e., $F(\overline{z};a,b)=\overline{F(z;a,b)}$, complex roots occur in conjugate pairs. Thus, among these three solutions, at least one is real.

Next, we prove that
\begin{align}\label{e:zlower2}
|z|
\geq
\frac18\min\{|b/a|,\rho_0\}.
\end{align}
Suppose first that
$
|b|^{2/3}>{|a|}/{4}.
$
For any \(|z|\leq |b|^{1/3}/8\), we have
\begin{align}
\begin{split}
|F(z;a,b)|
\geq
|b|-{|b|}/{8^3}
-{|a||b|^{1/3}}/{8}
-M\left({|b|^{1/3}}/{8}\right)^4
\geq
({|b|}/{4})
\left(1-{M|b|^{1/3}}/{64}\right)>0,
\end{split}
\end{align}
where the last two inequalities follow from
$
|a|<4|b|^{2/3},
$ and $
M|b|^{1/3}\leq M\fc\leq1.
$ 
The claim \eqref{e:zlower2} follows in this case.

Suppose next that
$
|b|^{2/3}\leq{|a|}/{4}.
$
 If
$
|z|
\leq
\min\{|b/a|,\rho_0\}/8
\leq {|b|}/({8|a|})
\leq
{|b|^{1/3}}/{32},
$
then
\begin{align}
\begin{split}
|F(z;a,b)|
\geq
|b|-{|b|}/{8}
-\left({|b|^{1/3}}/{32}\right)^3
-M\left({|b|^{1/3}}/{32}\right)^4\geq
|b|(1/4-M\fc)>0,
\end{split}
\end{align}
where we used \(|b|\leq\fc^3\) and \(\fc M\leq1/8\). Thus,
\eqref{e:zlower2} also holds.

Without loss of generality, assume that $b\ge 0$. The case $b<0$ follows by the same symmetry argument, so we omit it. 

For $a<0$, the derivative $\partial_z F(z;a,b)=3z^2-a+\sum_{k\ge4} k c_k z^{k-1}\ge 3z^2-a-M |z|^3>0$ for $z\in [-4\rho_0, 4\rho_0]$, provided $M\rho_0\le M\rho\le 2$. Thus $F(z;a,b)$ is increasing on $[-4\rho_0, 4\rho_0]$, and $F(z;a,b)$ has exactly one real solution. 

If $|a|\gtrsim |b|^{2/3}$, then $-b/a\lesssim b^{1/3}\le \rho_0$. Take $C\ge 2$ large enough such that $|a|\ge 2 b^{2/3}/C$. Then
\begin{align}
&F(-b/a;a,b)\ge (-b/a)^3 -M|b/a|^4\ge 0,\\
&F(-b/(Ca);a,b)\le -b+b/C+|b/(Ca)|^3 +M|b/(Ca)|^4\le -\frac{1}{2}\left(b-4|b/(Ca)|^3\right)\le 0,
\end{align}
provided that $M|b/a|\le 1$.
We conclude that the real root $z_0$ of $F(z;a,b)$ lies in $[ -b/a,-b/(Ca)]$. This yields \eqref{e:realroot1}.
%In this case we also have that
%\begin{align}
%\del_z F(z_0;a,b)\geq 3z_0^2-a-M|z_0|^3\geq -a, \quad \del_z F(z_0;a,b)\leq 3z_0^2-a+M|z_0|^3\lesssim 
%z_0^2-a\lesssim -a,
%\end{align}
%where we used that $M|z_0|\leq M|b/a|\leq 1$, and $|z_0|^2\leq |b/a|^2\lesssim |a|$. We conclude that $|\del_z F(z_0;a,b|\asymp |a|$. 

If $|b|^{2/3}\gtrsim |a|$, take $C\ge 4$ large enough such that $|a|\leq C^{1/3}b^{2/3}/4$. Then
\begin{align}\begin{split}\label{e:zbound}
&F((2b)^{1/3};a,b)\ge b -M|(2 b)^{1/3}|^4\ge 0,\\
&F((b/C)^{1/3};a,b)\le -b+b/C+|a| (b/C)^{1/3} +M|(b/C)|^{4/3}\le -\frac{1}{4}\left(b-4|a| (b/C)^{1/3}\right)\le 0,
\end{split}\end{align}
provided that $4M|b|^{1/3}\le 4 M\rho\le 1$.
We conclude that the real root $z_0$ of $F(z;a,b)$ lies in $[(b/C)^{1/3}, (2b)^{1/3}]$. 
%In this case we also have that
%\begin{align}
%\del_z F(z_0;a,b)\geq 3z_0^2-a-M|z_0|^3\geq z_0^2\geq |b/C|^{1/3}, \quad \del_z F(z_0;a,b)\leq 3z_0^2-a+M|z_0|^3\lesssim 
%z_0^2-a\lesssim |b|^{2/3},
%\end{align}
%where we used that $M|z_0|\leq M(2b)^{1/3}\leq 1$ and $(b/C)^{2/3}\leq |z_0|^2\leq (2b)^{2/3}$. We conclude that $|\del_z F(z_0;a,b|\asymp |b|^{2/3}$. 

Next, we study the complex conjugate pair of roots $z_\pm$. Let $z_\pm =x+\ri y$. Then 
\begin{align}\label{e:emF}
0=\Im[F(x+\ri y;a,b)]=-ay +(3x^2-y^2)y+y\sum_{k\ge 4}\frac{c_k\Im[(x+\ri y)^k]}{y}.
\end{align}
We note that $\Im[(x+\ri y)^k]\le k |x+\ri y|^{k-1}y$, and, dividing both sides of \eqref{e:emF} by $y$ and using \eqref{e:cusp_gbound}, we obtain
\begin{align}\label{e:Imbound}
|(3x^2-y^2)-a|\le M|x+\ri y|^3\leq (x^2+y^2)/8\Rightarrow y^2\geq \frac{1}{2}(x^2-a),
\end{align}
where we used that $M|x+\ri y|\leq \fc M\leq 1/2$.
In particular, 
\begin{align}\label{e:Imz_lowbb}
|\Im z_\pm|=y\geq \sqrt{|a|/2}
\end{align}

There are two cases:
\begin{enumerate}
\item If $|a|\gtrsim |b|^{2/3}$, we have already shown that $|z_\pm|\le 2\rho_0\lesssim \sqrt{|a|}$. Together with \eqref{e:Imz_lowbb},  we conclude \eqref{e:complexroot1}. Also in this case $|F'(z_\pm; a,b)|\asymp |a|$.

\item If $ |b|^{2/3}\gtrsim |a|$, we have already shown that $|z_\pm|\le 4\rho_0\lesssim b^{1/3}$. By the same argument as in \eqref{e:zbound}, if $|z|\le (b/C)^{1/3}$, then
\begin{align}\label{e:lowerbound}
|F(z;a,b)|\ge b-b/C-|a| (b/C)^{1/3} -M|(b/C)|^4>0.
\end{align}
Thus we conclude $|z_\pm|\asymp b^{1/3}$. Together with \eqref{e:Imbound}, we conclude that $|\Im z_\pm|=y\gtrsim \sqrt{x^2+y^2}=\sqrt{|z_\pm|}\gtrsim b^{1/3}$. This finishes the proof of \eqref{e:realroot2}. 
%Also in this case $|F'(z_\pm; a,b)|\asymp |b|^{2/3}$.
\end{enumerate}

Next, we study the case $a\ge 0$. A double root occurs when $F(z;a,b)=\partial_z F(z;a,b)=0$: 
\begin{align}\label{e:derF}
\partial_z F(z;a,b)=3z^2-a+\sum_{k\ge 4} k c_k z^{k-1}=0.
\end{align}
By Rouch\'e's theorem (using the same argument as in \eqref{e:cusp_F-Gdiff}), the equation \eqref{e:derF} has two solutions
\begin{align}\label{e:cusp_defzc}
z'_{\pm}=\pm \sqrt{a/3}\big(1+\OO(|a|)\big),
\end{align}
where we bound the higher-order terms in \eqref{e:derF} by $\OO(|z|^3)$ using \eqref{e:cusp_gbound}. 
Thus $F(z;a,b)$ is increasing on the intervals $[-4\rho_0, z'_-]\cup [z'_+, 4\rho_0]$ and decreasing on the interval $[z'_-,z'_+]$. Setting $F(z'_\pm;a,b'_\mp)=0$ yields
\begin{align}\label{e:cusp_defbc}
b'_\mp(a)=(z'_\pm)^3-a z'_\pm+\sum_{k\ge4} c_k (z'_\pm)^k
=-2(z'_\pm)^3+\sum_{k\ge4} (1-k)c_k (z'_\pm)^k= \mp 2(a/3)^{3/2}\big(1+\OO(\sqrt a)\big). 
\end{align}
In the last step we bound the higher-order terms by $\OO(|z'_\pm|^4)=\OO(|a|^2)$ using \eqref{e:cusp_gbound}.
We conclude that for $b'_-\le b\le b'_+$ the three roots are real, one in each of the intervals $[-4\rho_0, z'_-]$, $[z'_-,z'_+]$, and $[z'_+, 4\rho_0]$; denote them by $z_0\le z_1\le z_2$.
If $b>b'_+$ there is one real root $z_0\in [z'_+, 4\rho_0]$ and a complex pair $z_\pm$; if $b<b'_-$ there is one real root $z_0\in [-4\rho_0, z'_-]$ and a complex pair $z_\pm$.

In the following, we discuss the case $b\ge 0$. The case $b<0$ follows by the same symmetry argument, so we omit it.

If $0\le b\le b'_+$, then $b^{2/3}\lesssim a$ and $b/a\lesssim a^{1/2}\le \rho_0$. All zeros are bounded by $4\rho_0\lesssim \sqrt{a}$. The real roots $z_0\in [-4\rho_0, z'_-]$ and $z_2\in [z'_+, 4\rho_0]$ satisfy $z_0,z_2 \asymp \pm \sqrt{a}$. Next, we study the real root $z_1\in [z'_-, z'_+]$.
Take $C\ge 2$ large enough such that $a\ge 2 b^{2/3}/C$. Then
\begin{align}
F(-b/(Ca);a,b)\le -b+b/C+|b/(Ca)|^3 +M|b/(Ca)|^4\le -\frac{1}{2}\left(b-2|b/(Ca)|^3\right)\le 0.
\end{align}
We conclude that $z_1$ lies in $[z'_-, -b/(Ca)]$. If we further have $a\le 4 b^{2/3}$, then $a\asymp b^{2/3}$ and $z'_-\asymp -\sqrt{a}\asymp -b/a$. We conclude that $z_1\asymp -b/a$, which gives \eqref{e:realroot3}. Otherwise, if $a\ge 4 b^{2/3}$, then
\begin{align}
F(-2b/a;a,b)\ge b -|2b/a|^3-M|2b/a|^4\ge b -2|2b/a|^3\ge 0,
\end{align}
where we used that $2Mb/a\le M b^{1/3}\le M\rho\le 1$.
Thus $z_1$ lies in $[\max\{z'_-, -2b/a\}, -b/(Ca)]$. In this case as well, we have $z_1\asymp -b/a$, which gives \eqref{e:realroot3}.

If $b> b'_+$, then $b\gtrsim a^{3/2}$ and $b/a\lesssim a^{1/2}\le \rho_0$. All zeros are bounded by $4\rho_0\lesssim b^{1/3}$. There is one real root $z_0\in [z'_+, 4\rho_0]$ and two complex conjugate roots $z_\pm$. By the same argument as in \eqref{e:lowerbound}, all roots (in absolute value) are bounded below by $(b/C)^{1/3}$. Thus $z_0\asymp b^{1/3}$ and $|z_\pm|\asymp b^{1/3}$.

\end{proof}

\bibliography{TilingsConcentration.bib}
\bibliographystyle{alpha}

\end{document}